\documentclass[a4paper,11pt]{report}
\usepackage{amsmath,amssymb,amsthm,upref}
\usepackage[ansinew]{inputenc}
\usepackage[dvips]{graphics}
\usepackage{dsfont}
\usepackage{mathrsfs}
\usepackage{makeidx}
\usepackage[all]{xy}
\usepackage{multind}
\usepackage{graphicx}
\usepackage{hyperref}

 \makeindex{Index}
 \makeindex{Symbols}

\newcommand{\filt}{\mrm{CF}^+\mc{A}}
\newcommand{\Hom}[1]{\mathrm{Hom}_{\mbox{\tiny{$#1$}}}}

\newcommand{\cdc}{{Ch_\ast(\mathcal{A})}}
\newcommand{\cdco}{{Ch^\ast(\mathcal{A})}}
\newcommand{\cont}{{\mathcal C}}

\newcommand{\Dl}{\ensuremath{\Delta} }
\newcommand{\simp}{\Delta^\comp}

\newcommand{\ds}{\displaystyle}

\newcommand{\mc}[1]{\mathcal{#1}}
\newcommand{\mrm}[1]{\mathrm{#1}}
\newcommand{\mbf}[1]{\mathbf{#1}}

\newcommand{\comp}{{\scriptscriptstyle \ensuremath{\circ}}}

\DeclareMathSymbol{\minus}{\mathbin}{AMSb}{"72}
\DeclareMathSymbol{\downharpoonleft}
    {\mathrel}{AMSa}{"19}
\DeclareMathSymbol{\upharpoonright}
    {\mathrel}{AMSa}{"16}
\DeclareMathSymbol{\nwesq}
    {\mathrel}{AMSa}{"70}
\DeclareMathSymbol{\seesq}
    {\mathrel}{AMSa}{"79}

\newtheorem{thm}{\bf T{\footnotesize HEOREM}}[section]
\newtheorem{lema}[thm]{\bf L{\footnotesize EMMA}}
\newtheorem{prop}[thm]{\bf P{\footnotesize ROPOSITION}}
\newtheorem{cor}[thm]{\bf C{\footnotesize OROLLARY}}
\theoremstyle{definition}
\newtheorem{def{i}}[thm]{\bf D{\footnotesize EFINITION}}
\newtheorem{obs}[thm]{\bf R{\footnotesize EMARK}}
\newtheorem{ej}[thm]{\bf E{\footnotesize XAMPLE}}

\newenvironment{num}{\medskip
\refstepcounter{thm}\noindent {\bf (\thethm)}}{\vspace{1ex}\par}

\begin{document}

\thispagestyle{empty}



\mbox{}\\

\vfill

{\Huge \begin{center}\textbf{\textsc{SIMPLICIAL}}
\end{center}} {\Huge
\begin{center}\textbf{\textsc{DESCENT CATEGORIES}}\end{center}}

\vfill

\begin{center} BY\end{center}
\begin{center} BEATRIZ RODRIGUEZ GONZALEZ\footnote{Email: rgbea@algebra.us.es}\end{center}

\vspace{2cm}
\begin{center}
\begin{Large} DISSERTATION\end{Large}
\end{center}

\vspace{1cm}

\begin{flushright}
\begin{tabular}{l}ADVISORS:\\
LUIS NARV{\'A}EZ MACARRO\\
VICENTE NAVARRO AZNAR
\end{tabular}
\end{flushright}
 \vspace{2cm}

\begin{Large}
\begin{center}(Translated and revised version of the) thesis submitted \end{center}
\begin{center}for the degree of Doctor of Philosophy in Mathematics\end{center}
\begin{center} University of Seville, December 2007 \end{center}
\end{Large}
\vspace{1cm} \mbox{}


\newpage
\thispagestyle{empty}

\mbox{}

\vfill {\large

\begin{flushright}\textit{To my grandfather Mariano and my grandmother Carmen.}\end{flushright}

\vfill

{\large
\begin{center}\textbf{\textsc{Acknowledgements}}\end{center}}

First and foremost, I would like to thank my advisors Luis Narv{\'a}ez
Macarro and Vicente Navarro Aznar. Their support, dedication,
patience and expert guidance have made this dissertation possible.

On the other hand, I am indebted to Antonio Quintero Toscano for
our discussions concerning the (non-)exactness of geometric
realization. They were essential for the development of section
5.5.\\
I really appreciate the references and material that Francisco
Guill{\'e}n kindly gave to me. Also, it is my pleasure to thank Agust{\'\i}
Roig for his comments on the contents of section 6.2 and
subsection 6.2.1.

In addition, there are many people that I would like to thank for
their interest in my work and for their useful suggestions and
comments. Among them, I would like to
mention  I. Moerdijk, D. Cisinski, M. Vaqui{\'e}, L. Alonso Tarr{\'\i}o, B. Keller and G. Corti{\~n}as.\\

Finally, I would also like to gratefully acknowledge the support
provided by a Spanish FPU PhD Grant (ref. AP20033674), as well as
the partial support provided by the research projects
MTM-2004-07203-C02-01 and FQM-218:
\href{http://www.grupo.us.es/gfqm218/}{``Geometr{\'\i}a Algebraica,
Sistemas Diferenciales y Singularidades''}.

\vfill

\mbox{}


\renewcommand{\thepage}{\roman{page}}
\setcounter{page}{1}

\tableofcontents


\newpage
 \renewcommand{\thepage}{\arabic{page}}
\setcounter{page}{1}

\section*{Introduction}\addcontentsline{toc}{chapter}{Introduction}

In the f{i}eld of Algebraic Geometry, Grothendieck, at the
beginning of the sixties, glimpsed and impelled the introduction
of {\em derived categories} as the appropriate framework to handle
the general formulation of duality theorems, either
 ``continuous'' --on (quasi)coherent objects-- or ``discrete'' --on the
topological analogues motivated by $\ell$-adic cohomology-- (see
for instance the introduction in \cite{Hart}).\\
Grothendieck's program culminated in 1963 with Verdier's thesis,
where it is showed up, among other things, the importance of the
structure of {\em triangulated category}. This notion was also
related to ideas previously studied by Puppe in the f{i}eld of
Algebraic Topology \cite{V}.

Both notions, derived categories and triangulated categories, were
fundamental in the important developments in Algebraic Geometry
achieved in the period 1960-1975. Nevertheless, at the same time
these tools were considered quite sophisticated and not strictly
necessary for most questions. Consequently, the use of this
theory was not so widespread.\\
\indent The situation changed dramatically in the last three
decades. This was due to several reasons. Among them, the {\em
Riemann-Hilbert correspondence} and the discovering of {\em
perverse sheaves} in Algebraic Geometry. Later, it took place its
gradual introduction in Algebraic Topology, Representation Theory,
Mathematical Physics and Algebra in general, as well as in -as a
feedback- the development of the {\em theory of motives}.

Nowadays we can see a wide dif{f}usion of both notions, that have
become basic tools in Homological Algebra.
However, the notion of triangulated category does not seem to be
totally satisfactory. For instance, in \cite{GM} the non existence
of a functorial cone is remarked (see also \cite{Ne}).

Independently, and also in the sixties, Quillen introduced the
notion of {\em model category} \cite{Q}, establishing a general
abstract framework to study homotopy categories and
even {\em derived functors}.\\

On the other hand, at the middle of the twentieth century the
notion of \textit{simplicial object} arose to def{i}ne the
singular homology of topological spaces. Since then,
(co)simplicial objects have been present in the development of
homological and homotopical theories in Algebraic Topology and
Algebraic Geometry.

\indent Simplicial sets, and more generally simplicial techniques,
are also useful in the framework of model categories, for instance
through the natural notion of {\em simplicial model category}.
Another instance appeared later through the natural action of the
homotopy category of simplicial sets on the homotopy category of
any model category, which is in fact a key ingredient in the
triangulated category structure on the
latter one in the stable case (cf. \cite{Ho}).\\

\indent Model categories are a useful tool in the study of
localized categories arising in the {\em theory of motives} (cf.
\cite{FSV}, \cite{DLORV}) and more generally in the frameworks of
Algebraic Geometry and Homological Algebra, where they are a
complement of the notion of triangulated category. Nevertheless,
model categories do not always fulf{i}ll satisfactorily some
common situations. For instance, it is not easy to induce a model
category structure on the category of diagrams of a f{i}xed model
category. There is also some dif{f}iculty in handling f{i}ltered
structures that often appear in cohomological theories of
Algebraic Geometry.

On the other hand, simplicial structures --through the theory of
sheaves-- have been a relevant tool to deal with certain
(multiplicative) constructions, coming from Algebraic
Topology, in the framework of Algebraic Geometry \cite{God}.\\

In this work we introduce and develop the notion of {\em
$($co$)$simplicial descent category}, that is an evolution of the
corresponding ``cubic'' notion of Guill{\'e}n-Navarro \cite{GN}. It is
presented as an alternative or complementary instrument to be used
in the study of localized (or homotopic) categories
arising in Algebraic Geometry.\\

Let $\mc{D}$ be a category with f{i}nite coproducts, endowed with
a class of equivalences $\mrm{E}$. The category of simplicial
objects in $\mc{D}$, $\simp\mc{D}$, has an extremely rich
structure. Our aim is to transfer this richness to $\mc{D}$ with
the help of a ``simple'' functor
$\mbf{s}:\simp\mc{D}\rightarrow\mc{D}$. To this end we need to
impose some natural compatibility conditions between $\mbf{s}$ and
$\mrm{E}$. Thus, a {\em simplicial descent category} is the data
$(\mc{D},\mbf{s})$ satisfying certain
axioms, as\\
\noindent\textbf{Normalization:} the simple of the constant simplicial object associated with $X$ is equivalent to $X$\\
\textbf{Exactness:} $\mbf{s}(\Dl\mrm{E})\subset \mrm{E}$\\
\textbf{Factorization:} abstraction of Eilenberg-Zilber-Cartier's theorem appearing in \cite{DP}\\
\textbf{Acyclicity:} the image under $\mbf{s}$ of the simplicial
cone existing in $\simp\mc{D}$ must be a cone object in
$\mc{D}$ with respect to the equivalence class $\mrm{E}$.\\
Using the simple functor, we obtain {\em cone} and {\em cylinder}
\underline{functors} in $\mc{D}$ satisfying the ``usual'' properties.\\
\indent The notion of {\em cosimplicial descent category} turns
out to be the dual notion, that is, the opposite category of a
simplicial descent category.

\subsection*{Background}

To use a ``simple'' functor in order to transfer a structure is
not a new idea, and it has appeared since the beginning of
Topology, for instance when $\mbf{s}=$geometric realization, fact
that is emphasized in \cite{May} $\S$ 11.\\[0.3cm]
Grothendieck and his school introduced ``geometric'' simplicial
methods in Algebraic Geometry through the so-called
\textit{simplicial hypercovers}. They are an essential tool used
by Deligne to def{i}ne a {\em mixed Hodge structure} on the
cohomology of any complex algebraic variety $S$ (not
necessarily smooth).\\
Technically, the key point is the existence of a suitable
``simple'' functor
$$\left\{\begin{array}{l}\mbox{\footnotesize Cosimplicial Mixed}\\
                         \mbox{\footnotesize Hodge Complexes}\end{array}\right\}\longrightarrow
                        \left\{\begin{array}{l}\mbox{\footnotesize Mixed Hodge}\\
                                                        \mbox{\footnotesize complexes}\end{array}\right\}$$
that induces a mixed Hodge structure on the cohomology
of $S$ through a simplicial hypercover $X$ of $S$.\\
A similar procedure is followed in \cite{DB} to construct a
f{i}ltered De Rham complex on the cohomology of a singular variety.\\[0.3cm]

Another instance of simple functor motivating this work appears in
\cite{N}. Let $\textbf{Adgc}$ be the category of commutative
dif{f}erential graded algebras over a f{i}eld of characteristic 0.
In loc. cit. a ``simple'' functor
$$ \{\mbox{Cosimplicial objects in } \mathrm{Adgc}\}\rightarrow\{\mathrm{Adgc}\}$$
is introduced. This functor is known as the {\em Thom-Whitney
simple}. It is used mainly to transfer Sullivan techniques from
Algebraic Topology to Algebraic Geometry. As an application, the
author endows the rational homotopy spaces of an algebraic variety
with a mixed Hodge structure.\\[0.3cm]
In \cite{GN}, F. Guill{\'e}n and V. Navarro Aznar give an axiomatic
and abstract notion of ``simple'' functor, inspired in Deligne's
and Thom-Whitney simples, but formulated in the ``cubical''
framework \cite{GNPP}.
To this end, they introduce the notion of (cubical) {\it
cohomological descent category} and develop an extension criterion
of a functor $F:\{$smooth schemes$\}\rightarrow \mc{D}$ to the
category of non smooth schemes, provided that $\mc{D}$ is the
localized category of a (cubical) cohomological descent category.\\[0.3cm]
In this work we develop the notion of (co)simplicial descent
category, that is widely based in the previous notion of
(co)homological descent category, where the basic objects are
diagrams of cubical shape instead of simplicial objects.\\
\indent Both notions share the same philosophy, but there are
important dif{f}erences between them:\\
\indent -cubical diagrams are f{i}nite whereas simplicial objects
are inf{i}nite.\\
\indent -in the cubical case, the factorization axiom is not so
strong, in fact it is usually an automatic consequence of the
Fubini theorem on the index swapping in a double coend.
However, the simplicial factorization axiom is much stronger,
because it involves the diagonal object associated with a
bisimplicial object, and it is in fact an abstraction of the
Eilenberg-Zilber-Cartier theorem given in \cite{DP}. This theorem
uses the degeneracy maps of simplicial objects. Hence, strict
simplicial objects (with no degeneracy maps) are not enough for
our purposes.
Nevertheless, a cubical diagram does not have degeneracy maps.\\

On the other hand, working in the simplicial framework has other
advantages. For instance, one can induce a natural action of
$\simp Set$ on $\mc{D}$ (def{i}ned through the simple functor from
the action of $\simp Set$ on $\simp\mc{D}$).\\
We can also exploit the homotopy structure of $\simp\mc{D}$ when
$\mc{D}$ is a simplicial descent category. It turns out that
homotopic morphisms between simplicial objects (in the classical
sense of simplicial homotopy) are mapped by $\mbf{s}$ into
identical morphisms in the localized category
$\mc{D}[\mrm{E}^{-1}]$. In particular, (simplicial) homotopy
equivalencies are mapped by $\mbf{s}$ into equivalences. This
applies to augmentations with an ``extra degeneracy''.
\subsection*{Main results}
\indent \textbf{a)} We establish a set of axioms for the
\textit{(co)simplicial descent categories}. These axioms unify the
properties satisf{i}ed by a significant number of examples in the
frameworks of Algebraic Topology and Algebraic Geometry.\\
To be precise, a category $\mc{D}$ with f{i}nite coproducts and
f{i}nal object $\ast$ together with a saturated class of
equivalences $\mrm{E}$ and a ``simple'' functor
$\mbf{s}:\simp\mc{D}\rightarrow\mc{D}$ is
a \emph{simplicial descent category} if the following axioms are satisf{i}ed.\\
\textbf{Additivity:} The canonical morphism
$\mbf{s}X\sqcup\mbf{s}Y\rightarrow\mbf{s}(X\sqcup Y)$ is an
equivalence for any $X,Y\in\simp\mc{D}$.
Also $\mrm{E}\sqcup\mrm{E}\subseteq\mrm{E}$.\\
\textbf{Factorization:} There exists a natural equivalence
$\mu_Z:\mbf{s}\mrm{D}Z\rightarrow\mbf{s}\simp\mbf{s}Z$, where
$\mrm{D}Z$ is the diagonal object associated with the bisimplicial
object $Z$ and $\mbf{s}\simp\mbf{s}Z$ is its iterated simple.\\
\textbf{Normalization:} There exists a natural equivalence
$\lambda_X:\mbf{s}(X\times\Dl)\rightarrow X$, compatible with
$\mu$, relating an object $X$ in $\mc{D}$ to the simple of its associated constant simplicial object $X\times\Dl$.\\
\textbf{Exactness:} If $f$ is a morphism in $\simp\mc{D}$ such
that $f_n\in\mrm{E}$ for all $n$ then $\mbf{s}f\in\mrm{E}$.\\
\textbf{Acyclicity:} The simple of a morphism $f$ in $\simp\mc{D}$
is an equivalence if and only if the simple of its simplicial cone, $\mbf{s}Cf$, is acyclic.\\
\textbf{Symmetry:} The class $\{f:X\rightarrow Y\, | \,
\mbf{s}f\in\mrm{E}\}$ is invariant by the operation of inverting
the order of the face and degeneracy maps of $X$ and $Y$.\\
\indent \textbf{ b)} We introduce the following
`transfer lemma', that will be widely used to produce examples of (co)simplicial descent categories\\
\textit{Assume that $(\mc{D}',\mbf{s}',\mrm{E}',\lambda',\mu')$ is
a simplicial descent category and $\mc{D}$ is a category with a
simple functor $\mbf{s}$ and ``compatible'' natural
transformations $\lambda$ and $\mu$. Moreover, let
$\psi:\mc{D}\rightarrow \mc{D}'$ be a functor such that the
following diagram
$$\xymatrix@H=4pt@C=20pt@R=20pt{\simp\mc{D} \ar[r]^{\psi}\ar[d]_{\mbf{s}} & \simp\mc{D}' \ar[d]^{\mbf{s'}} & \\
                                \mc{D}\ar[r]^{\psi}         &                                \mc{D}'       }$$
commutes up to natural equivalence, ``compatible'' with
transformations $\lambda$, $\lambda'$ and $\mu$, $\mu'$. Then
$(\mc{D},\mrm{E}=\psi^{-1}\mrm{E}')$ is also a simplicial descent
category.}

\indent \textbf{c)} In any simplicial descent category $\mc{D}$
we can induce cone and cylinder functors in the following way.\\
Given a morphism $f:X\rightarrow Y$ in $\mc{D}$, consider it in
$\simp\mc{D}$ as a constant simplicial map
$f\times\Dl:X\times\Dl\rightarrow Y\times\Dl$. Then, the cone of
$f$ is by def{i}nition $\mbf{s}C(f\times\Dl)$, where $C$ is the simplicial cone associated with $f\times\Dl$.\\
Similarly, the ``cylinder'' of two morphisms
$A\stackrel{f}{\leftarrow} B\stackrel{g}{\rightarrow} C$ in
$\mc{D}$ is $\mbf{s}Cyl(f\times\Dl,g\times\Dl)$, where $Cyl$ is
the simplicial cylinder associated with $(f\times\Dl,g\times\Dl)$.\\
When $X$ is a simplicial set, then the classical cylinder
associated with $X$  is just our cylinder of $X=X=X$, and the one
associated with a morphism $f$ is our cylinder of
$X=X\stackrel{f}{\rightarrow} Y$.

\indent \textbf{d)} We provide a ``reasonable'' description of the
morphisms in $Ho\mc{D}=\mc{D}[\mrm{E}^{-1}]$, the homotopy
category of $\mc{D}$. In general the class of equivalences
$\mrm{E}$ does not has calculus of fractions. The key point to
obtain this description is the cylinder functor and its properties.\\
More specif{i}cally, consider the functor
$\mrm{R}:\mc{D}\rightarrow\mc{D}$ def{i}ned by
$\mrm{R}X:=\mbf{s}(X\times\Dl)$. Note that the natural
transformation $\lambda:\mrm{R}\rightarrow Id$ is a pointwise
equivalence. Then a morphism $f:X\rightarrow Y$ in $Ho\mc{D}$ is
represented by a sequence
$$\xymatrix@M=4pt@H=4pt@R=10pt{  X & \mrm{R}X \ar[l]_-{\lambda_X}\ar[r]^-{{f}'} & T & \mrm{R}Y \ar[l]_{w}\ar[r]^-{\lambda_{Y}}& Y\ .}$$
where all arrows except $f'$ are equivalences. Using this
description we can prove, for instance, that $Ho\mc{D}$
is additive when $\mc{D}$ is so.\\
\indent \textbf{e)} The shift $[1]:\mc{D}\rightarrow\mc{D}$ is
defined by $X[1]=c(X\rightarrow \ast)$, where $c$ is the cone
functor given above. We consider the class of distinguished
triangles in $Ho\mc{D}$ consisting of those triangles isomorphic,
for some $f$, to
$$\xymatrix@M=4pt@H=4pt{  X\ar[r]^f &  Y \ar[r] & c(f) \ar[r] & X[1]\ .}$$
Distinguished triangles satisfy all axioms of triangulated
category except the second axiom TR 2 (that is, the one involving
the shift of distinguished triangles), with no extra assumptions
(neither additivity).\\
Moreover, in the additive case $Ho\mc{D}$ is a ``suspended'' (or
right triangulated) category \cite{KV} in the simplicial case, and
it is ``cosuspended'' (or left triangulated) in the cosimplicial
case. In particular, if in addition the shift is an automorphism
of $Ho\mc{D}$ then $Ho\mc{D}$ is a triangulated category.\\
\indent \textbf{f)} In order to study the properties of the
simplicial cone and cylinder functors, we develop a much more
general construction, the ``\textit{total simplicial object}''
associated with a ``biaugmented bisimplicial object'' (or, more
generally, to a ``$n$-augmented $n$-simplicial object''). More
concretely, consider the bisimplicial object $Z$ given by the
picture
$$\xymatrix@M=4pt@H=4pt@C=25pt{                                   &                                                                                             &                                                                                                                                              &                                                                                                                                                          &                                                                                                                             & \\
   Z_{-1,2}\ar@{}[u]|{\vdots}\ar@<0ex>[d]\ar@<1ex>[d] \ar@<-1ex>[d]                        & {Z_{0,2}}\ar[l]\ar@<0ex>[d]\ar@<1ex>[d] \ar@<-1ex>[d]\ar@{}[u]|{\vdots}\ar@/_0.75pc/[r]     & {Z_{1,2}} \ar@{}[u]|{\vdots}\ar@<0.5ex>[l] \ar@<-0.5ex>[l]\ar@<0ex>[d]\ar@<1ex>[d] \ar@<-1ex>[d]\ar@/_1pc/[r]\ar@/_0.75pc/[r]                &Z_{2,2} \ar@{}[u]|{\vdots}\ar@<0ex>[l]\ar@<1ex>[l] \ar@<-1ex>[l]\ar@<0ex>[d]\ar@<1ex>[d] \ar@<-1ex>[d]\ar@/_1pc/[r]\ar@/_0.75pc/[r]\ar@{-}@/_1.25pc/[r]   & Z_{3,2}\ar@{}[u]|{\vdots}\ar@<0.33ex>[l]\ar@<-0.33ex>[l]\ar@<1ex>[l]\ar@<-1ex>[l]\ar@<0ex>[d]\ar@<1ex>[d] \ar@<-1ex>[d]     & \cdots \\
   Z_{-1,1}\ar@<0.5ex>[d] \ar@<-0.5ex>[d] \ar@/^1pc/[u]\ar@/^0.75pc/[u]                    & {Z_{0,1}}\ar[l]\ar@<0.5ex>[d] \ar@<-0.5ex>[d] \ar@/^1pc/[u]\ar@/^0.75pc/[u]\ar@/_0.75pc/[r] &  Z_{1,1} \ar@/^1pc/[u]\ar@/^0.75pc/[u]\ar@<0.5ex>[l]\ar@<-0.5ex>[l]\ar@<0.5ex>[d] \ar@<-0.5ex>[d]\ar@/_1pc/[r]\ar@/_0.75pc/[r]               &Z_{2,1} \ar@/^1pc/[u]\ar@/^0.75pc/[u] \ar@<0ex>[l]\ar@<1ex>[l] \ar@<-1ex>[l] \ar@<0.5ex>[d] \ar@<-0.5ex>[d]\ar@/_1pc/[r]\ar@/_0.75pc/[r]\ar@/_1.25pc/[r]  & Z_{3,1}\ar@<0.33ex>[l]\ar@<-0.33ex>[l]\ar@<1ex>[l]\ar@<-1ex>[l] \ar@<0.5ex>[d] \ar@<-0.5ex>[d]\ar@/^1pc/[u]\ar@/^0.75pc/[u] & \cdots \\
   Z_{-1,0} \ar@/^0.75pc/[u]                                                               &  Z_{0,0}\ar[d] \ar[l]\ar@/_0.75pc/[r] \ar@/^0.75pc/[u]                                      & {Z_{1,0}} \ar[d] \ar@<0.5ex>[l] \ar@<-0.5ex>[l]  \ar@/^0.75pc/[u]  \ar@/_1pc/[r]\ar@/_0.75pc/[r]                                             &Z_{2,0} \ar[d] \ar@<0ex>[l]\ar@<1ex>[l] \ar@<-1ex>[l]\ar@/^0.75pc/[u]\ar@/_1pc/[r]\ar@/_0.75pc/[r]\ar@/_1.25pc/[r]                                        & Z_{3,0}\ar[d] \ar@<0.33ex>[l]\ar@<-0.33ex>[l]\ar@<1ex>[l]\ar@<-1ex>[l] \ar@/^0.75pc/[u]                                     & \cdots \\
                                                                                           &  Z_{0,-1}\ar@/_0.75pc/[r]                                                             & Z_{1,-1}  \ar@<0.5ex>[l]\ar@<-0.5ex>[l]\ar@/_1pc/[r]\ar@/_0.75pc/[r]                                                                         &Z_{2,-1} \ar@<0ex>[l]\ar@<1ex>[l] \ar@<-1ex>[l]\ar@/_1pc/[r]\ar@/_0.75pc/[r]\ar@/_1.25pc/[r]                                                              & Z_{3,-1}\ar@<0.33ex>[l]\ar@<-0.33ex>[l]\ar@<1ex>[l]\ar@<-1ex>[l]                                                            & \cdots}$$
Then the total simplicial object associated with a biaugmented is
in degree $k$ the coproduct of the $k$-th diagonal of $Z$. The
face and degeneracy maps are def{i}ned respectively as coproducts
of those of $Z$.\\
It turns out that this total functor is left adjoint to the
``total decalage'' given in \cite{IlII} p.7. We can also consider
the total functor as the simplicial analogue to the total chain
complex associated with a double complex.\\
\indent \textbf{g)} On one hand, we have checked that all examples
of (cubical) (co)ho\-mo\-lo\-gi\-cal descent categories \cite{GN}
are simplicial descent categories. Among them we can mention
(filtered) cochain complexes, topological spaces or commutative dif{f}erential graded algebras.\\
On the other hand, we provide other examples than the cubical
ones. For instance, we consider a category of mixed Hodge
complexes, and we endow it with a cosimplicial descent category
structure. In this structure the simple functor is just the one
developed in \cite{DeIII}. As a corollary we obtain a triangulated
structure on its homotopy category
(similar to the one obtained in \cite{Be}).\\
Also the category of DG-modules over a f{i}xed DG-category is a
cosimplicial descent category, and we deduce the usual
triangulated structure existing in its homotopy category \cite{K}.\\
A possibly less known example is the category of cochain complexes
together with a biregular f{i}ltration, where the class
$\mrm{E}_2$ of equivalences consists of those morphisms which
induce isomorphism in the second term of the respective spectral sequences.\\
Now we get a triangulated structure on the
 category of bounded-below f{i}ltered complexes localized with
 respect to the class $\mrm{E}_2$. Moreover, the ``decalage''
 functor of a f{i}ltration \cite{DeII} I.3.3 is a triangulated functor
with values in the (usual) f{i}ltered derived category.

\subsection*{Contents}

\textsc{Chapter 1}: The f{i}rst chapter contains the
simplicial/combinatorial preliminaries.\\
We study the classical cone and cylinder functors in $\simp\mc{D}$
as particular cases of the {\it total functor} of biaugmented
bisimplicial objects, that was mentioned before. As far as the
author knows, this total functor has not been previously studied.
Particular and related cases can be found in \cite{EP} and \cite{AM}).\\
The total functor satisf{i}es interesting properties. For
instance, the iteration of totals of $n$-augmented $n$-simplicial
objects does not depend on the order in which we compute it
(analogously to the property of
the total complex associated with a multiple chain complex).\\[0.2cm]
\textsc{Chapter 2}: This chapter contains the def{i}nition of
(co)\-sim\-pli\-cial descent category, as well as some of their
properties, mainly those of homotopical type, related to the cone
and cylinder functors. The axioms of (co)\-sim\-pli\-cial descent
category are natural in the following sense: If $I$ is a small
category and $\mc{D}$ is a (co)\-sim\-pli\-cial descent category
then the category of functors from $I$ to $\mc{D}$ (endowed with
the pointwise simple and the pointwise equivalences) is again a
(co)\-sim\-pli\-cial
descent category.\\
In section 2.3 the ``factorization'' property of the cylinder is
established. In terms of the cone functor this property means the
following: Consider a commutative diagram $C$ in a simplicial
descent category $\mc{D}$
$$\xymatrix@M=4pt@H=4pt@C=15pt@R=15pt{
  X \ar[r]^f\ar[d]^g & Y\ar[d]^{g'} \\
  X'\ar[r]^{f'} &Y' \ .}$$
If we apply the cone functor by rows and columns respectively we
get
$$\xymatrix@M=4pt@H=4pt@C=15pt@R=15pt{
  X \ar[r]^f\ar[d]^g & Y\ar[d]^{g'} & c(f) \ar[d]^\alpha \\
  X'\ar[r]^{f'}      & Y'           & c(f') \\
  c(g) \ar[r]^\beta  & c(g') \ .}$$
Then the cone of $\alpha$ and the cone of $\beta$ are equivalent
in a natural way. This fact will play an important role in chapter
\ref{CapitDescTriang}, since it is the key point in the proofs of
the octahedron axiom and of the second axiom of triangulated categories.\\
Section 2.4 is devoted to the study of the properties of the
square
$$\xymatrix@M=4pt@H=4pt{
  \mbf{s}X \ar[r]^-{\mbf{s}f}\ar[d]_{\mbf{s}\epsilon} & \mbf{s}Y\ar[d]_{I_Y} \\
  \mbf{s}X_{-1}\ar[r]^-{I_{X_{-1}}}   & \mbf{s}Cyl(f,\epsilon) }$$
obtained by applying $\mbf{s}$ to the respective square in
$\simp\mc{D}$ induced by the simplicial cylinder.
One can check that this square ``commutes up to equivalence'', and
that $I_{X_{-1}}$ is an equivalence provided that $\mbf{s}f$ is
so. The reciprocal assertion also holds under some extra
assumptions, and it will be needed to prove the ``transfer
lemma''.\\
In section 2.5 we introduce the notion of functor of simplicial
descent categories, and we prove the ``transfer lemma''.\\[0.2cm]
\textsc{Chapter 3}: In this chapter we give a ``reasonable''
description of the morphisms of $Ho\mc{D}$. We use this
description to prove that $Ho\mc{D}$ is additive if $\mc{D}$ is so.\\[0.2cm]
\textsc{Chapter 4}: We prove here that the class of distinguished
triangles def{i}ned through the cone functor
$$\xymatrix@M=4pt@H=4pt{  X\ar[r]^f &  Y \ar[r] & c(f) \ar[r] & X[1]}$$
satisfies the axioms TR1, TR3 and TR4 of triangulated categories
(with no extra assumptions). In the additive case, the right
implication of TR2 holds, so $\mc{D}$ is a ``suspended'' category
\cite{KV}. Moreover, if the shift is an isomorphism of
categories then $Ho\mc{D}$ is a triangulated category.\\[0.2cm]
\textsc{Chapter 5}: In this chapter we exhibit examples of
simplicial descent categories. The f{i}rst one is the category of
chain complexes in an additive or abelian category, taking as
equivalences $\mrm{E}=$homotopy equivalences or
$\mrm{E}=$quasi-isomorphisms. In this example the axioms of
simplicial descent category are checked ``by hand'', whereas in
the remaining simplicial
examples they are checked by means of the transfer lemma.\\
The following picture contains the main examples of simplicial
descent categories included in this chapter as well as the
functors of simplicial descent categories between them
$$\xymatrix@M=4pt@H=4pt{ Top \ar@<0.5ex>[r]^{S} & \simp Set\ar@<0.5ex>[l]^{|\cdot|} \ar[r] & \simp Ab \ar[r] & Ch_\ast Ab }\ .$$
\textsc{Chapter 6}: Examples of cosimplicial descent categories
are provided in this chapter. Cochain complexes are obtained just
as the dual case of chain complexes. In section \ref{SecADGC} the
category of commutative dif{f}erential graded algebras (over a
f{i}eld of characteristic 0) is considered. The simple is just the
`Thom-Whitney simple' given in \cite{N}.\\
In section \ref{SecDGmod} we endow the category of DG-modules
\cite{K} with a cosimplicial descent category structure.\\
In the next section we prove that the category of (positive)
complexes together with a (biregular) f{i}ltration,
$\mrm{CF}^+\mc{A}$ has two dif{f}erent cosimplicial descent
structures. In the f{i}rst one,
(${}_1\mrm{CF}^+\mc{A}$,$\mrm{E}$), the equivalences $\mrm{E}$ are
the f{i}ltered quasi-isomorphisms. In the second one,
(${}_2\mrm{CF}^+\mc{A}$,$\mrm{E}_2$), the class of equivalences
$\mrm{E}_2$ consists of those morphisms inducing isomorphism in
the second term of the spectral sequence. These statements follow
from the transfer lemma. It is applied twice to the functors
contained in the following diagram
$$\xymatrix@M=4pt@H=4pt{ {}_2\mrm{CF}^+\mc{A} \ar[r]^{Dec} & {}_1\mrm{CF}^+\mc{A} \ar[r]^{Gr} & (Ch^\ast\mc{A})^{\mathbb{Z}} }\ .$$
The functor $Dec$ is the decalage functor given in \cite{DeII}
I.3.3, and $Gr$ is the graded functor, with values in the category
of $\mathbb{Z}$-graded cochain complexes, endowed with the degreewise descent category structure.\\
Both structures are used to induce a cosimplicial descent category
structure on ``the'' category of mixed Hodge complexes.\\

To finish, we include an appendix containing the
Eilenberg-Zilber-Cartier theorem \cite{DP}, and some extra
properties which are not easy to find in the existing literature.

\subsection*{Further research/Open problems}

Next we list some questions and problems related to this work.
Some of them are natural questions and others are further
applications and complements.

\textbf{I.} The category $Op(\mc{D})$ of operads over a symmetric
monoidal descent category $\mc{D}$ has a natural structure of
descent category (joint work with A. Roig). A related open problem
is to endow the category of operadic algebras over a f{i}xed
operad $\mc{P}\in Op(\mc{D})$ with a structure of cosimplicial
descent category.

\textbf{II.} Every cubical diagram $X$ in a f{i}xed category
$\mc{D}$ gives rise in a natural way to a simplicial object in
$\mc{D}$, $\tau X$ (\cite{N}, 12.1). If $\mc{D}$ is a simplicial
descent category, we can compose $\tau$ with the simple functor
$\mbf{s}:\simp\mc{D}\rightarrow \mc{D}$, obtaining in this way a
``cubical simple functor'' $\{\mbox{Cubical diagrams in
}\mc{D}\}\rightarrow\mc{D}$ . The following natural question
arises:

Is every (co)simplicial descent category a ``cubical''
(co)homological descent category in the sense of \cite{GN}?
If the answer is af{f}irmative, then the ``extension of functors''
theorem, given in --loc. cit.--, will be also valid for functors
with values in the localized category of a cosimplicial descent
category.

\textbf{III.} The localized category $Ho\mc{D}$ of a descent
category $\mc{D}$ has a translation functor
$\mrm{T}:Ho\mc{D}\rightarrow Ho\mc{D}$ as well as a class of
distinguished triangles (coming from the cone functor in $\mc{D}$).\\
The following question is motivated by model category theory:

Is $Ho\mc{D}$ an additive category provided that $\mrm{T}$ is an
automorphism in $Ho\mc{D}$ (stable case)?\\
If the translation functor $\mrm{T}:Ho\mc{D}\rightarrow Ho\mc{D}$
is not an automorphism, it would be interesting to provide an
abstract process of `stabilization', similar to the construction
of the category of spectra from the category of topological
spaces.

\textbf{IV.} Study the properties of the action of $\simp Set$
over a descent category $\mc{D}$ inherited from the natural action
of $\simp Set$ on $\simp\mc{D}$ through the simple functor.
Furthermore, possibility of carrying this action to the level of
derived categories. That is, to check whether this action gives
rise to another one from $Ho\simp Set$ on $Ho\mc{D}$ or not.

\textbf{V.} Def{i}ne sheaf cohomology with values in a descent
category $\mathcal{D}$, in a similar way to that given in \cite{N}
for the case $\mc{D}$=commutative dif{f}erential graded algebras.
In a more general sense, to tackle the def{i}nition of derived
functors in descent categories, following the recent work
\cite{GNPR}.

\textbf{VI.} Relationship between simplicial descent categories
and model categories. In the ``cubical'' case, it is known that
the subcategory of f{i}brant objects in a simplicial model
category is a cohomological descent category \cite{R}.

\textbf{VII.} Extension of the notion of descent category to the
context of f{i}bred categories and stacks.

\textbf{VIII.} Study the relationship with the recent work
\cite{Vo}, where it is also used the simplicial cylinder
$\widetilde{Cyl}$ introduced in the f{i}rst chapter of this work.
For instance, if $\mc{D}$ is a simplicial descent category, it
holds that the class $\{f\ | \ \mbf{s}f\in\mrm{E}\}$ is a
$\Delta$-closed class in the sense of loc. cit..

%
%
%
%
%
%

%


\setcounter{chapter}{0}
\chapter{Preliminaries}
\section{Simplicial objects}

In this section we will remind the def{i}nition and some basic
properties of simplicial objects in a f{i}xed category $\cont$.
For a more detailed exposition see \cite{May}, \cite{GZ} or
\cite{GJ}.

\begin{def{i}}[The Simplicial Category $\Dl$]\mbox{}\\
The simplicial category\index{Index}{simplicial!category}
$\Delta$\index{Symbols}{$\Dl$} has as objects the ordered sets
$[n]\equiv \{0,\ldots , n\}$ with $0<1<\cdots <n$, and as
morphisms the (weak) monotone functions, that is
$$ \Hom{\Dl }([m],[n])=\{ f:[m]\longrightarrow [n] \mbox{ with } f(i)\leq f(j)\mbox{ if } i\leq j \}.$$
There exists two kind of relevant morphisms in the category
$\Delta$.\\
The face morphisms $\partial_i =
\partial_i^n: [n-1]\rightarrow [n]$ are just those monotone functions
such that $\partial_i (\{0,\ldots,n-1\})=\{0,\ldots,i-1,i+1,\ldots
,n\}$, for all $i=0,\ldots,n$.\\
The degeneracy morphisms $\sigma_i^n=\sigma_i : [n+1]\rightarrow
[n]$ are characterized by $\sigma_i (i)=\sigma_i (i+1)=i$, for all
$i=0,\ldots,n $.

More specif{i}cally, $\partial_i(l)= l${ if }$l\leq i-1$ and
$\partial_i(l)=l+1${ if }$l\geq i$, whereas $\sigma_i(l)= l${ if
}$l\leq i$ and $\sigma_i(l)=l-1${ if }$l>i$.

These morphisms satisfy the so called ``simplicial
identities''\index{Index}{simplicial identities}, that are the
following equalities
\begin{equation}\label{igsimp} \begin{array}{rcl} \partial^{n+1}_j \partial^n_i           & = & \partial^{n+1}_i \partial^{n}_{j-1} \mbox{ if } i<j\\
                      \sigma^n_j \sigma^{n+1}_i & = & \sigma^n_i \sigma^{n+1}_{j+1}\,\;\; \mbox{ if }i\leq j\\
                      \sigma^{n-1}_j \partial^{n}_i & = & \left\{ \begin{array}{lcl} \partial^n_i \sigma^{n-2}_{j-1}& \mbox{ if } & i < j \\
                                                                              Id_{[n-1]} & \mbox{ if }& i=j \mbox{ or } i=j+1 \\
                                                                              \partial^{n-1}_{i-1} \sigma^{n-2}_{j} & \mbox{ if }& i> j+1 .\end{array}\right.\\
                      \end{array}\end{equation}
\end{def{i}}

The category $\Delta$ is generated by the face and degeneracy
morphisms, as described in the following proposition \cite{May}.

\begin{prop}\label{PropFactMorDl}
Let $f:[n]\rightarrow [m]$ be a morphism in $\Delta$ dif{f}erent
from the identity. Denote by $i_1> i_2\cdots >i_s$ those elements
of $[m]$ that do not belong to the image of $f$, and by
$j_1<j_2\cdots <j_l$ those elements of $[n]$ such that
$f(j_k)=f(j_k +1)$. Then
\begin{equation}\label{gendelta} f=\partial_{i_1}\cdots \partial_{i_s}\sigma_{j_1}\cdots\sigma_{j_l}\, ,\end{equation}
Moreover, the factorization of $f$ in this way is unique.
\end{prop}

\begin{obs}\label{catDeltagorro}
Let
$\widehat{\Dl}$\index{Symbols}{$\widehat{\Dl}$}\index{Index}{categoy!of
ordered sets} be the category whose object are all the (non empty)
f{i}nite ordered sets, and whose morphisms are the monotone maps.\\
If $E=\{e_0< e_1 < \cdots < e_n\}$ is an object of
$\widehat{\Dl}$, then $E$ is canonically isomorphic to $[n]$, and
$n$ is the cardinal of $E$ minus 1.\\
Then each object of $\widehat{\Dl}$ is isomorphic to a unique
object of $\Dl$, or equivalently, $\Dl$ is a skeletal subcategory
of $\widehat{\Dl}$. Then it follows the existence of a functor
$p:\Dl\rightarrow \widehat{\Dl}$ quasi-inverse of the inclusion
$i:\Dl\rightarrow
\widehat{\Dl}$.\\
The intrinsic meaning of some def{i}nitions and properties given
in terms of $\Dl$ become clarif{i}ed when expressed in terms of
 $\widehat{\Dl}$, as we will check along this chapter.
\end{obs}

We will use the following operations relative to $\widehat{\Dl}$.

\begin{num}\label{sumaOrdenada}\index{Index}{ordered sum}
Denote by
$$+:\widehat{\Dl}\times\widehat{\Dl}\rightarrow
\widehat{\Dl}$$ the ``ordered sum'' of ordered sets.\\
That is,{ if }$E$ and $F$ are objects of $\widehat{\Dl}$, then
$E+F$\index{Symbols}{$E+F$} is $E\sqcup F$ as a set. The order in
$E+F$ is the one compatible with those of $E$ and $F$, and such
that $e<f${ if }$e\in E$ and $f\in F$.\\
Analogously, $+:\Dl\times\Dl\rightarrow \Dl$ is such that
$[n]+[m]=p(i([n])+i([m]))=[n+m+1]$, where $[n]$ is identif{i}ed
with $\{0,\ldots,n\}\subset [n+m+1]$ and $[m]$ with
$\{n+1,\ldots,n+m+1\}\subset [n+m+1]$.
\end{num}

\begin{num}\label{OrdenOpuesto}\index{Symbols}{$op$}
Denote by
$$\widehat{op}: \widehat{\Dl}\rightarrow\widehat{\Dl}$$ the
functor which consists of taking the opposite order. That is,
$op(E)$ is equal to $E$ as a set, but it has the inverse order of
$E$.

\noindent Analogously, $op=p\comp\widehat{op}\comp
i:\Dl\rightarrow \Dl$ is the functor given by
$$op([n])=[n] \ \mbox{ and if }\  \theta:[n]\rightarrow [m],\ (op(\theta))(i)=m-\theta(n-i) \ ,$$
and then
$$op(\partial_i^n)=\partial_{n-i}^n:[n-1]\rightarrow [n]\mbox{ and } op(\sigma_j^n)=\sigma_{n-j}^n:[n+1]\rightarrow [n]\ .$$
\end{num}

\begin{def{i}}[The Category $\Dl_e$]\label{Deltaestricto}\mbox{}\\
Let $\Dl_e$\index{Symbols}{$\Dl_e$} be the strict simplicial
category\index{Index}{strict simplicial category}, that is the
subcategory of $\Dl$ with the same objects, but whose morphisms
are the injective monotone functions.\\
Analogously, $\Dl_e$ is generated by the face morphisms, and it is
a skeletal subcategory of the corresponding category
$\widehat{\Dl}_e$.
\end{def{i}}

\begin{def{i}}[Simplicial Objects]\mbox{}\\
A simplicial object $X$\index{Index}{simplicial!object} in a
category $\cont$ is a contravariant functor from the simplicial
category to $\cont$, that is, $X: \Delta ^{\comp}\rightarrow
\cont$.
\end{def{i}}

\begin{num}
As a corollary of (\ref{gendelta}), $X$ is characterized by the
data
$$\begin{array}{ccc} X_p=X([p]) & d_i=X(\partial_i) & s_j=X(\sigma_j) \end{array}$$
where the face and degeneracy maps $d_i$ and $s_j$ of $X$ satisfy
the following equalities, also called simplicial identities
\index{Index}{simplicial!identities}:
\begin{equation}\label{igobsim} \begin{array}{rcl} d^n_i d^{n+1}_j           & = & d^{n}_{j-1}d^{n+1}_i  \mbox{ if } i<j\\
                      s^{n+1}_i s^n_j  & = &s^{n+1}_{j+1} s^n_i \,\;\; \mbox{ if }i\leq j\\
                      d^{n}_i s^{n-1}_j  & = & \left\{ \begin{array}{lcl} s^{n-2}_{j-1} d^n_i & \mbox{ if } & i < j \\
                                                                              Id_{[n-1]} & \mbox{ if }& i=j \mbox{ or } i=j+1 \\
                                                                              s^{n-2}_{j} d^{n-1}_{i-1} & \mbox{ if }& i> j+1 .\end{array}\right.\\
                      \end{array}\end{equation}
Then, a simplicial object $X=\{ X_n , d_i , s_j \}$ can be
represented as follows
$$\xymatrix@1{ X_0 \; \ar@/_2pc/[rr]_{s_0} && {\;} X_1 \;
\ar@<0.5ex>[ll]^-{d_0} \ar@<-0.5ex>[ll]_-{d_1} \ar@/_2pc/[rr]
\ar@/_1pc/[rr] && \; X_2 \; \ar@<0ex>[ll] \ar@<1ex>[ll]
\ar@<-1ex>[ll] \ar@/_1.5pc/[rr] \ar@/_2pc/[rr] \ar@/_1pc/[rr] &&
{\;} X_3
\ar@<0.33ex>[ll]\ar@<-0.33ex>[ll]\ar@<1ex>[ll]\ar@<-1ex>[ll]&
\cdots\cdots }.$$
\end{num}

\begin{num}\label{ObjetSimplEstricto}
Analogously, a strict simplicial object\index{Index}{stric
simplicial object} is $X:\Dl_e^\comp \rightarrow \cont$, that is
given by
$$\xymatrix@1{ X_0 \;  && {\;} X_1 \;
\ar@<0.5ex>[ll]^-{d_0} \ar@<-0.5ex>[ll]_-{d_1}  && \; X_2 \;
\ar@<0ex>[ll] \ar@<1ex>[ll] \ar@<-1ex>[ll]  && {\;} X_3
\ar@<0.33ex>[ll]\ar@<-0.33ex>[ll]\ar@<1ex>[ll]\ar@<-1ex>[ll]&
\cdots\cdots }.$$
\end{num}

\begin{num} Dually, a cosimplicial object in $\cont$ is a functor
$X:\Dl\rightarrow\cont$, or equivalently, a simplicial object in
$\cont^\comp$. The strict cosimplicial objects in $\cont$ are
def{i}ned in the same way.
\end{num}

\begin{obs} From now on, we will use the notation $\mc{I}\cont$
for the category of functors from $\mc{I}$ to $\cont$.
\end{obs}

\begin{def{i}}
The simplicial objects in $\cont$ give rise to the
category\index{Index}{category!of simplicial objects} $\simp
\cont$\index{Symbols}{$\simp\cont$}, whose  morphisms are
the natural transformations between functors.\\
A morphism $\rho : X \rightarrow X'$ in $\simp \cont$ is a set of
morphisms $\rho_n : X_n \rightarrow X'_n$ in $\cont$ commuting
with the face and degeneracy maps, that is
$$ \rho_n d_X = d_{X'} \rho_{n+1} \;\;\; \rho_{n+1} s_X = s_{X'} \rho_n
\ .$$
\end{def{i}}

\begin{num}\label{equivDlyDlGorro} Having into account \ref{catDeltagorro}, $\simp \cont$ is
canonically equivalent to $\widehat{\Dl}^\comp \cont$. Denote by
$\mathcal{I}:\widehat{\Dl}^\comp \cont\rightarrow\simp \cont $
(resp. $\mathcal{P}:\simp \cont\rightarrow \widehat{\Dl}^\comp
\cont$) the functor induced by composition with $i:\Dl\rightarrow
\widehat{\Dl}$ (resp. $p:\widehat{\Dl}\rightarrow \Dl$).\\
Since $p\comp i=Id_{\Dl}$, we have that
$\mc{I}\comp\mc{P}=Id_{\simp\cont}$.
\end{num}
\begin{num}\label{CatObjSimplEstricCosEtc}
The categories of cosimplicial, strict simplicial and strict
cosimplicial objects in
$\cont$\index{Symbols}{$\simp_e\cont$}\index{Symbols}{$\Dl\cont$}\index{Symbols}{$\Dl_e\cont$}
are denoted respectively by $\Dl\cont$\index{Symbols}{$\Dl\cont$},
$\simp_e \cont$ and $\Dl_e\cont$\index{Symbols}{$\Dl_e\cont$}, and are def{i}ned in the same way as $\simp\cont$.\\[0.2cm]
Composing with the inclusion $\Dl_e\rightarrow \Dl$ we obtain the
forgetful functor
$$ \mrm{U}: \simp\cont\rightarrow \Dl_e^{\comp} \cont $$
consisting in forgetting the degeneracy maps of a simplicial
object.
\end{num}

\begin{obs} If $\cont$ has colimits, the forgetful functor $\mrm{U}: \simp\cont\rightarrow \Dl_e^{\comp} \cont$ has a left adjoint $\pi$
(called ``Dold-Puppe transformation'') constructed as usual.
However, it follows from the property \ref{PropFactMorDl} that
$\pi$ exists just assuming the existence of coproducts in $\cont$.
Next we remind the def{i}nition of the Dold-Puppe transformation
(cf. \cite{G} 1.2).
\end{obs}

\begin{prop}\label{Dold-Puppe}
If $\cont$ is a category with f{i}nite coproducts, the forgetful
functor $\mrm{U}:\simp\cont\rightarrow \Dl_e^{\comp} \cont$ admits
a left adjoint $\pi:\Dl_e^{\comp} \cont\rightarrow \simp\cont$. If
$A$ is a strict simplicial object in $\cont$, then $\pi A$ is
def{i}ned as
$$(\pi A)_n = \coprod_{\theta:[n]\twoheadrightarrow [m]} A_{m}^{\theta}$$
where the coproduct is indexed over the set of surjective
morphisms $\theta:[n]\twoheadrightarrow [m]$ in $\Dl$ and
$A_{m}^{\theta}=A_m$.
\end{prop}

\begin{num} Let us see how the action or $\pi A$ over the morphisms in $\Dl$ is
def{i}ned.\\
Let $f:[n']\rightarrow [n]$ be a morphism in $\Dl$. The morphism
in $\cont$
$$ (\pi A )(f): (\pi A)_n= \coprod_{\theta:[n]\twoheadrightarrow [m]} A_{m}^{\theta} \longrightarrow (\pi A)_{n'}=\coprod_{\rho:[n']\twoheadrightarrow [m']} A_{m'}^{\rho}$$
is def{i}ned as follows.\\
Given a surjective $\theta:[n]\twoheadrightarrow [m]$, it follows
from \ref{PropFactMorDl} that there exists a unique factorization
of $\theta\comp f:[n']\rightarrow [m]$ as
 $$ [n']\stackrel{\alpha}{\twoheadrightarrow}  [l]\stackrel{\beta}{\hookrightarrow}   [m]$$
where $\alpha$ is surjective and $\beta$ is injective.\\
Then the restriction of $(\pi A )(f)$ to $ A_{m}^{\theta}$ is just
$$ A(\beta): A_m^{\theta} \rightarrow A_l^{\alpha} \ .$$
\end{num}

\begin{def{i}}[Bisimplicial objects]\mbox{}\\
A bisimplicial object\index{Index}{bisimplicial object} in $\cont$
is by def{i}nition a simplicial object in $\simp\cont$. Then, the
category of bisimplicial objects in $\cont$ is
$\simp\simp\cont\simeq
(\Dl\times\Dl)^\comp\cont$\index{Symbols}{$\simp\simp\cont$}.\\
Given $Z\in\simp\simp\cont$ we will denote
$$\begin{array}{l} d_i^{(1)}=Z(\partial_i,Id):Z_{n,m}\rightarrow Z_{n-1,m}\ \ d_i^{(2)}=Z(Id,\partial_i):Z_{n,m}\rightarrow Z_{n,m-1}\\
                   s_j^{(1)}=Z(\sigma_j,Id):Z_{n,m}\rightarrow Z_{n+1,m}\ \ s_j^{(2)}:Z_{n,m}\rightarrow Z_{n,m+1}\
                   .\end{array}$$
\end{def{i}}
Now we introduce some remarks that will be useful along these
notes, as well as their dual constructions in the cosimplicial
case.

\begin{num} The diagonal functor $\mrm{D}:\simp\simp\cont\rightarrow\simp\cont$
\index{Symbols}{$\mrm{D}$}\index{Index}{diagonal of a bisimp.
object} is the functor induced by composition with $\Dl\rightarrow
\Dl\times\Dl$, $[n]\rightarrow ([n],[n])$.\\
Then,{ given }$Z\in\simp\simp\cont$, $(\mrm{D}Z)_n=Z_{n,n}$ and
$(\mrm{D}Z)(\theta)=Z(\theta,\theta)$, $\forall n\geq 0$ and
$\theta$ in $\Dl$.
\end{num}

\begin{num} The index swapping in $\simp\simp\cont$ gives rise to a canonical functor $\Gamma:\simp\simp\cont\rightarrow
\simp\simp\cont$,\index{Symbols}{$\Gamma$} with
$$(\Gamma Z)_{n,m}=Z_{m,n}\ \ (\Gamma Z)(\alpha,\beta)=Z(\beta,\alpha) \ ,$$
if $n,m\geq 0$ and $\alpha$, $\beta$ are morphisms in $\Dl$.\\
It holds that $\mrm{D}\comp \Gamma = \mrm{D}$.
\end{num}

\begin{num}\label{simpFuntor}
Each functor $F:\mc{C}\rightarrow\cont'$ induces by composition a
functor between the respective categories of simplicial objects of
$\mc{C}$ and $\cont'$, that will be denoted by $\simp
F:\simp\mc{C}\rightarrow\simp\cont$. Then $(\simp F(X))_n=F(X_n)$
$\forall n\geq 0$.
\end{num}
%

\begin{num}
Let $$-\times\Dl:\cont\longrightarrow \simp \cont
$$\index{Index}{constant simplicial object
}\index{Symbols}{$-\times\Dl$}be the simplicial constant functor.
More concretely, $X\times\Dl: \simp \rightarrow \cont$ is just the
constant functor equal to $X$
$$\begin{array}{lll} (X\times\Dl)_n= X \ \ \forall n\geq 0 & ; &
(X\times\Dl)(f)=Id_{X}\ \ \forall \mbox{ morphism $f$ of
}\Delta\end{array}$$
Note that{ if }$\ast$ is a f{i}nal object (resp. initial) in
$\cont$, so is $\ast\times\Dl$ in $\simp\cont$.
\end{num}

\begin{num}\label{ctebisimplicial}
The functor $-\times\Dl:\cont\rightarrow \simp\cont$ induces the
functors
$$ -\times\Dl\ ,\ \ \Dl\times- : \simp\cont\rightarrow \simp\simp\cont\ .$$
Given $X$ in $\simp\cont$ then
$$(X\times\Dl)_{n,m}=X_n \mbox{ and } (\Dl\times X)_{n,m}=X_m \ \ \forall n,m\geq 0 \ .$$
\end{num}

\begin{obs}
The category $\simp\cont$ inherits most of the properties
satisf{i}ed by $\cont$. For instance,{ if }$\cont$ has f{i}nite
coproducts the same holds for $\simp \cont$ through:
$$\begin{array}{rcl}(X\coprod Y)_n = X_n \coprod Y_n &  ; &(X\coprod Y)(f)=X(f)\coprod Y(f)\end{array}.$$
Analogously, $\simp\cont$ is additive (resp. monoidal, abelian,
complete, cocomplete, etc){ if }$\cont$ is so.
\end{obs}

%
%
%
%
%
%

\section{Augmented simplicial objects}

\begin{def{i}}
Let $\Dl_{+}$\index{Symbols}{$\Dl_{+}$} be the category whose
objects are the ordinal numbers $[n]=\{0,\ldots,n\}$, $n\geq -1$,
where $[-1]=\emptyset$, and whose morphisms are the (weak)
monotone functions.\\
Then $\Dl_{+}$ contains $\Dl$ as a full subcategory.\\
Denote by $\partial_0: [-1]\rightarrow [0]$ the trivial morphism,
that will be considered as a face morphism. As in the case of
$\Dl$, it also holds that every morphism has a unique
factorization in terms of the face and degeneracy maps. These maps
satisfy in the same way the simplicial identities.
\end{def{i}}

\begin{obs}\label{DeltagorroMasMonoidal} Analogously to \ref{catDeltagorro}, $\Dl_{+}$ is
a skeletal subcategory of
$\widehat{\Dl}_{+}$\index{Symbols}{$\widehat{\Dl}_{+}$}, that is
the category of (eventually empty) ordered sets and monotone
functions. We will also denote by
$i:\Dl_+\rightarrow\widehat{\Dl}_+$ and
$p:\widehat{\Dl}_+\rightarrow \Dl_+$ the inclusion and its
quasi-inverse.\\
The operation
$+:\widehat{\Dl}_{+}\times\widehat{\Dl}_{+}\rightarrow\widehat{\Dl}_{+}$\index{Index}{ordered
sum}\index{Symbols}{$E+F$} is def{i}ned as (\ref{sumaOrdenada}),
and it makes $\widehat{\Dl}_{+}$ (then also $\Dl_+$) into a
(strong) monoidal category (cf. \cite{ML} p.171).
\end{obs}
\begin{def{i}}

An augmented simplicial object\index{Index}{simplicial
object!augmented} in a category $\cont$ is a functor $X^+
:(\Dl_+)^{\comp}\rightarrow \cont$.

Denote by $\Dl_+^{\comp}\cont$ the category of augmented
simplicial objects in $\cont$.

Therefore, an augmented simplicial object $X^+ \in
\Dl_+^{\comp}\cont$ is characterized by the data $\{X_n , d_i ,
s_j\}$, where $X_n=X([n])$, $d_i= X(\partial_i)$ and $s_j=
X(\sigma_j)$. That is, $X$ is represented by the following diagram
$$\xymatrix@1{ X_{-1}{\,} & {\;} X_0 \; \ar[l]_{d_0} \ar@/_2pc/[rr]_{s_0} && {\;} X_1 \;
\ar@<0.5ex>[ll]^-{d_0} \ar@<-0.5ex>[ll]_-{d_1} \ar@/_2pc/[rr]
\ar@/_1pc/[rr] && \; X_2 \; \ar@<0ex>[ll] \ar@<1ex>[ll]
\ar@<-1ex>[ll] \ar@/_1.5pc/[rr] \ar@/_2pc/[rr] \ar@/_1pc/[rr] &&
{\;} X_3
\ar@<0.33ex>[ll]\ar@<-0.33ex>[ll]\ar@<1ex>[ll]\ar@<-1ex>[ll]&
\cdots\cdots }.$$
\end{def{i}}
\begin{def{i}}
Given a simplicial object $X\in\simp\cont$, an
augmentation\index{Index}{augmentation} of $X$ is an augmented
simplicial object $X^+ \in \Dl_+^{\comp}\cont$ such that its image
under the forgetful functor $\mrm{U}:\Dl_+^{\comp}\cont
\rightarrow\simp\cont$ induced by restriction is just $X$.
\end{def{i}}

\begin{obs}\label{aumentacion}
Given $X\in \simp\cont$, an augmentation of $X$ is a pair
$(X_{-1},d_0)$ where $X_{-1}$ is an object in $\cont$, and $d_0 :
X_{0}\rightarrow X_{-1}$ is such that $d_0 d_0 =d_1
d_0:X_{1}\rightarrow X_{-1}$.

If $\cont$ has a f{i}nal object 1, then every simplicial object
$X$ has a trivial augmentation, taking $X_{-1}=1$.
Hence, we have the functor $\simp\cont\rightarrow\simp_+ \cont$,
right adjoint of $\mrm{U}:\simp_+ \cont\rightarrow\simp\cont$.
\end{obs}

\begin{prop}\label{aumentvisionequivalente}\mbox{}\\
i) Given an object $S$ in $\cont$ and a morphism $\epsilon
:X\rightarrow S\times\Dl$ in $\simp\cont$, then
$\epsilon_0:X_0\rightarrow S$ is an augmentation of $X$, and the
correspondence $\epsilon\rightarrow \epsilon_0$ is a bijection
$$\mathrm{Hom}_{\simp\cont}(X,S\times\Dl)\simeq\{X^+\in(\Dl_+)^{\comp}\cont\mbox{ with }X_{-1}=S\mbox{ and }\mrm{U}X^+=X\}$$
ii) The functor $-\times\Dl:\cont\rightarrow\simp\cont$ is left
adjoint to $\simp\cont\rightarrow\cont$, $X_\cdot\rightarrow X_0$.
That is,
$$\mathrm{Hom}_{\simp\cont}(S\times\Dl,X)\simeq\mrm{Hom}_{\cont}(S,X_0)\ .$$
\end{prop}

\begin{proof}
\textit{i)} Clearly,{ if }$\epsilon : X \rightarrow S\times\Dl$,
then $d_0=\epsilon_0:X_{0}\rightarrow S$ is an augmentation of
$X$. Conversely,{ if }$d_0:X_0\rightarrow S$ is an augmentation of
$X$, then $\epsilon_n =(d_0)^{n+1} : X_n \rightarrow S$ def{i}nes
a morphism $\epsilon : X \rightarrow S\times\Dl$.
Moreover, it follows by induction that $\epsilon_n=(d_0)^{n+1}$.\\
To check \textit{ii)}, it is enough to note that given
$\epsilon:S\times\Dl\rightarrow X$ then
$\epsilon_n=(s_0)^n\epsilon_0$.
\end{proof}

Two relevant examples of augmented simplicial object are the
``decalage'' objects associated with any simplicial object (see
\cite{IlII}).

\begin{num}
Consider the functor $F:\widehat{\Dl}\rightarrow\widehat{\Dl}$
(resp. $G:\widehat{\Dl}\rightarrow\widehat{\Dl}$) given by
$F(E)=[0]+E$ (resp. $G(E)=E+[0]$), that is, the ordered set
obtained from
$E$ by adding a smallest (resp. greatest) element.\\
If we work in the category $\Dl$, the functor
$F:\Dl\rightarrow\Dl$ maps $[n]$ into $[n+1]$ and{ if }$\theta$ is
a morphism $\Dl$, $F(\theta)(0)=0$
and $F(\theta)(i)=\theta(i-1)+1$,{ if }$i>0$.\\
On the other hand, the functor $G:\Dl\rightarrow\Dl$ is such that
$G([n])=[n+1]$ and{ if }$\theta:[n]\rightarrow [m]$ then
$G(\theta)(i)=i${ if }$i<n+1$, $G(n+1)=m+1$.\\
\end{num}

\begin{def{i}}[``Decalage'' objects]\label{ejempDecalado}\mbox{}\\
The ``lower decalage functor'',\index{Index}{decalage!lower}
$dec_1:\simp\cont\rightarrow\simp\cont$\index{Symbols}{$dec_1$},
is def{i}ned by composing with $F$. If $X:\Dl^\comp\rightarrow
\cont$, then  $X\comp F$ is obtained by ``forgetting the f{i}rst
face and degeneracy maps'' of $X$
$$(dec_1(X))_n=X_{n+1} \ \ (dec_1(X))(d_i)=d_{i+1}\ \ (dec_1(X))(s_j)=s_{j+1}\ .$$
The morphism $d_1:X_1\rightarrow X_0$ gives rise to the
augmentation $dec_1(X)\rightarrow X_0\times\Dl$ given by
$$\xymatrix@1{ X_{0}{\,} & {\;} X_1 \; \ar[l]_{d_1} \ar@/_2pc/[rr]_{s_1} && {\;} X_{2} \; \ar@<0.5ex>[ll]^-{d_1} \ar@<-0.5ex>[ll]_-{d_{2}} \ar@/_2.5pc/[rr]_{s_1} \ar@/_1.5pc/[rr]_{s_{2}} && \; X_{3} \; \ar@<0ex>[ll]_{d_{2}} \ar@<2ex>[ll]_-{d_{1}} \ar@<-2.25ex>[ll]_-{d_{3}} \ar@/_1.5pc/[rr] \ar@/_2pc/[rr] \ar@/_1pc/[rr] && {\;} X_{4}\ar@<0.33ex>[ll]\ar@<-0.33ex>[ll]\ar@<1ex>[ll]\ar@<-1ex>[ll]&\cdots\cdots }\ .$$
In the same way, by composition with $G$ we obtain the functor
``upper decalage'',\index{Index}{decalage!upper}
$dec^1(X):\simp\cont\rightarrow\simp\cont$\index{Symbols}{$dec^1$}
that consists of ``forgetting the last face and degeneracy maps''
$$(dec^1(X))_n=X_{n+1} \ \ (dec^1(X))(d_i)=d_{i}\ \ (dec^1(X))(s_j)=s_{j}\ .$$
Therefore, $d_0:X_1\rightarrow X_0$ produces the augmentation
$dec^1(X)\rightarrow X_0\times\Dl$
$$\xymatrix@1{ X_{0}{\,} & {\;} X_1 \; \ar[l]_{d_0} \ar@/_2pc/[rr]_{s_0} && {\;} X_{2} \; \ar@<0.5ex>[ll]^-{d_0} \ar@<-0.5ex>[ll]_-{d_1} \ar@/_2.25pc/[rr]_{s_0} \ar@/_1.5pc/[rr]_{s_1} && \; X_{3} \; \ar@<0ex>[ll]_{d_1} \ar@<2ex>[ll]_-{d_0} \ar@<-2ex>[ll]_-{d_2} \ar@/_1.5pc/[rr] \ar@/_2pc/[rr] \ar@/_1pc/[rr] && {\;} X_{4} \ar@<0.33ex>[ll]\ar@<-0.33ex>[ll]\ar@<1ex>[ll]\ar@<-1ex>[ll]&\cdots\cdots}\ .$$
\end{def{i}}

\subsection{Simplicial homotopy and extra degeneracy}

We will remind the relationship between an augmentation with an
extra degeneracy and the notion of homotopy between simplicial
morphisms (cf. \cite{B}, p. 78).

F{i}rst we give the combinatorial def{i}nition of homotopic
morphisms in (cf. \cite{May} $\S$5).

\begin{def{i}}[Simplicial homotopy]\label{def{i}homotSimplicial}\mbox{}\\
The morphism $f:X\rightarrow Y$ in $\simp\cont$ is homotopic to
$g:X\rightarrow Y$\index{Index}{homotopic morphisms}, $f\sim g$,{
if }there exists morphisms $h_i:X_n\rightarrow Y_{n+1}$,
$i=0,\ldots,n$ satisfying the following identities
$$
\begin{array}{ll}
 i)  & d_0h_0=f\, ,\; d_{n+1}h_n=g\\[0.1cm]
 ii) & d_ih_j=\begin{cases} h_{j-1}d_i &\mbox{ if } i<j\\
                            d_{j}h_{i-1} &\mbox{ if } i=j\geq 1\\
                            h_jd_{i-1} &\mbox{ if }i>j+1\end{cases}\\[0.1cm]
 iii)& s_ih_j=\begin{cases} h_{j+1}s_i &\mbox{ if }i\leq j\\
                            h_js_{i-1} &\mbox{ if }i>j .\end{cases}
\end{array}$$
Note that the relation $\sim$ is not symmetric.
\end{def{i}}%

\begin{def{i}}\label{def{i}degmas}\mbox{}\\
\textbf{1.} An augmentation $X\rightarrow X_{-1}\times\Dl$ has a
``lower'' extra degeneracy\index{Index}{extra degeneracy}{ if
}there exists morphisms $s_{-1}^n=s_{-1}:X_{n}\rightarrow X_{n+1}$
in $\cont$ for all $n\geq -1$ such that the following simplicial
identities hold
\begin{equation}\label{extcono}
d_0 s_{-1}=Id  \ \  d_{i+1}s_{-1}=s_{-1}d_i \ \
s_{j}s_{-1}=s_{-1}s_{j-1} \  \forall i\geq 0 ,\  j\geq 0 \ ,
\end{equation}
where $\epsilon_0=d_0:X_0\rightarrow X_{-1}$.\\
\textbf{2.} Dually, an augmentation $X\rightarrow X_{-1}\times\Dl$
has an ``upper'' extra degeneracy{ if }there exists morphisms
$s_{n+1}^n=s_{n+1}:X_{n}\rightarrow X_{n+1}$ in $\cont$ for all
$n\geq -1$ such that
$$d_{n+1} s_{n+1}=Id \ \ \  d_{i}s_{n+1}=s_n d_i  \ \ \  s_{j}s_{n+1}=s_{n+2} s_{j}\ \  \forall i \leq n ,\ \  j\leq n+1 .$$
\end{def{i}}

\begin{ej}
Following the notations introduced in \ref{ejempDecalado}, it is
clear that the augmentation $dec_1(X)\rightarrow X_0\times\Dl$ has
a ``lower'' extra degeneracy, consisting
 of the forgotten degeneracy $s_0:X_n\rightarrow X_{n+1}$.\\
Analogously, $dec^1(X)\rightarrow X_0\times\Dl$ has an ``upper''
extra degeneracy $s_{n}:X_n\rightarrow X_{n+1}$.
\end{ej}

\begin{prop}[\cite{B} cap. 3, 3.2]\label{degmasHomot}\mbox{}\\
An augmentation $\epsilon: X\rightarrow X_{-1}\times\Dl$ has a
lower extra degeneracy{ if }and only{ if }there exists
$\zeta:X_{-1}\times\Dl\rightarrow X$ such that $Id_X \sim
\zeta\epsilon$ and $\epsilon\zeta=Id_{X_{-1}\times\Dl}$.\\
Dually, $\epsilon$ has an upper extra degeneracy{ if }and only{ if
}there exists $\zeta:X_{-1}\times\Dl\rightarrow X$ such that
$\zeta\epsilon \sim Id_X$ and $\epsilon\zeta=Id_{X_{-1}\times\Dl}$
\end{prop}


\section{Total object of a biaugmented bisimplicial object}

In this section $\mc{D}$ denotes a category with f{i}nite
coproducts.

In this case, one can consider the classical cone object
associated with a morphism $f:X\rightarrow Y$ in $\simp\mc{D}$, as
well as the classical cylinder object associated with $X$. It
turns
out that both classical objects are particular cases of the ``total'' functor developed in this section.\\

In fact, the total functor can be seen as the simplicial analogue
of the total chain complex associated with a double chain complex
in
an additive category. Even though this total functor is an extremely natural construction, the author could not find it in the literature.\\

Next we introduce this combinatorial construction associated with
any biaugmented bisimplicial object $Z$, although
we will only use some particular cases of $Z$ in this work.\\

In the cosimplicial setting, all dual constructions and properties
can be established.

\begin{def{i}}[biaugmented bisimplicial object]\mbox{}\\
Let $2-\Dl_{+}$\index{Symbols}{$2-\Dl_{+}$} be the full
subcategory of $\Dl_{+}\times\Dl_{+}$ whose objects are the pairs
$([n],[m])\in\Dl_{+}\times\Dl_{+}$ such that $[n]$ and $[m]$ are
not both empty.

A biaugmented bisimplicial object $Z$ is a functor
$Z:2-\Dl_{+}^\comp\rightarrow \mc{D}$, or equivalently a diagram
$Z_{-1,\cdot} \leftarrow Z^+_{\cdot,\cdot}\rightarrow
Z_{\cdot,-1}$, where $Z_{-1,\cdot}$, $Z^+_{\cdot,\cdot}$ and
$Z_{\cdot,-1}$ are the respective restrictions of $Z$ to
$[-1]\times\Dl$, $\Dl\times\Dl$ and $\Dl\times [-1]$.

Hence $Z$ is a diagram of $\mc{D}$ of the form
\begin{equation}\label{diagrObjBisimpBiaum}\xymatrix@M=4pt@H=4pt@C=25pt{                                   &                                                                                             &                                                                                                                                              &                                                                                                                                                          &                                                                                                                             & \\
                                                                        Z_{-1,2}\ar@{}[u]|{\vdots}\ar@<0ex>[d]\ar@<1ex>[d] \ar@<-1ex>[d]                        & {Z_{0,2}}\ar[l]\ar@<0ex>[d]\ar@<1ex>[d] \ar@<-1ex>[d]\ar@{}[u]|{\vdots}\ar@/_0.75pc/[r]     & {Z_{1,2}} \ar@{}[u]|{\vdots}\ar@<0.5ex>[l] \ar@<-0.5ex>[l]\ar@<0ex>[d]\ar@<1ex>[d] \ar@<-1ex>[d]\ar@/_1pc/[r]\ar@/_0.75pc/[r]                &Z_{2,2} \ar@{}[u]|{\vdots}\ar@<0ex>[l]\ar@<1ex>[l] \ar@<-1ex>[l]\ar@<0ex>[d]\ar@<1ex>[d] \ar@<-1ex>[d]\ar@/_1pc/[r]\ar@/_0.75pc/[r]\ar@{-}@/_1.25pc/[r]   & Z_{3,2}\ar@{}[u]|{\vdots}\ar@<0.33ex>[l]\ar@<-0.33ex>[l]\ar@<1ex>[l]\ar@<-1ex>[l]\ar@<0ex>[d]\ar@<1ex>[d] \ar@<-1ex>[d]     & \cdots \\
                                                                        Z_{-1,1}\ar@<0.5ex>[d] \ar@<-0.5ex>[d] \ar@/^1pc/[u]\ar@/^0.75pc/[u]                    & {Z_{0,1}}\ar[l]\ar@<0.5ex>[d] \ar@<-0.5ex>[d] \ar@/^1pc/[u]\ar@/^0.75pc/[u]\ar@/_0.75pc/[r] &  Z_{1,1} \ar@/^1pc/[u]\ar@/^0.75pc/[u]\ar@<0.5ex>[l]\ar@<-0.5ex>[l]\ar@<0.5ex>[d] \ar@<-0.5ex>[d]\ar@/_1pc/[r]\ar@/_0.75pc/[r]               &Z_{2,1} \ar@/^1pc/[u]\ar@/^0.75pc/[u] \ar@<0ex>[l]\ar@<1ex>[l] \ar@<-1ex>[l] \ar@<0.5ex>[d] \ar@<-0.5ex>[d]\ar@/_1pc/[r]\ar@/_0.75pc/[r]\ar@/_1.25pc/[r]  & Z_{3,1}\ar@<0.33ex>[l]\ar@<-0.33ex>[l]\ar@<1ex>[l]\ar@<-1ex>[l] \ar@<0.5ex>[d] \ar@<-0.5ex>[d]\ar@/^1pc/[u]\ar@/^0.75pc/[u] & \cdots \\
                                                                        Z_{-1,0} \ar@/^0.75pc/[u]                                                               &  Z_{0,0}\ar[d] \ar[l]\ar@/_0.75pc/[r] \ar@/^0.75pc/[u]                                      & {Z_{1,0}} \ar[d] \ar@<0.5ex>[l] \ar@<-0.5ex>[l]  \ar@/^0.75pc/[u]  \ar@/_1pc/[r]\ar@/_0.75pc/[r]                                             &Z_{2,0} \ar[d] \ar@<0ex>[l]\ar@<1ex>[l] \ar@<-1ex>[l]\ar@/^0.75pc/[u]\ar@/_1pc/[r]\ar@/_0.75pc/[r]\ar@/_1.25pc/[r]                                        & Z_{3,0}\ar[d] \ar@<0.33ex>[l]\ar@<-0.33ex>[l]\ar@<1ex>[l]\ar@<-1ex>[l] \ar@/^0.75pc/[u]                                     & \cdots \\
                                                                                                                                                                &  Z_{0,-1}\ar@/_0.75pc/[r]                                                             & Z_{1,-1}  \ar@<0.5ex>[l]\ar@<-0.5ex>[l]\ar@/_1pc/[r]\ar@/_0.75pc/[r]                                                                         &Z_{2,-1} \ar@<0ex>[l]\ar@<1ex>[l] \ar@<-1ex>[l]\ar@/_1pc/[r]\ar@/_0.75pc/[r]\ar@/_1.25pc/[r]                                                              & Z_{3,-1}\ar@<0.33ex>[l]\ar@<-0.33ex>[l]\ar@<1ex>[l]\ar@<-1ex>[l]                                                            & \cdots
                                                                                                                              }\end{equation}
\mbox{}\\
We will denote by $2-\Dl_{+}^{\comp}\mc{D}$ the category of
biaugmented bisimplicial objects, whose morphisms are natural
transformations between functors.
\end{def{i}}

\begin{num}\label{equiv2DlyDlGorro} Again, $2-\Dl_{+}^{\comp}\mc{D}$ is canonically equivalent to the category
$2-\widehat{\Dl}_+^\comp\mc{D}$, where
$2-\widehat{\Dl}_+$\index{Symbols}{$2-\widehat{\Dl}_+$} is the
full subcategory of
 $\widehat{\Dl}_+\times\widehat{\Dl}_+$ having as morphisms the pairs
$(E,F)$ of ordered sets $E$ and $F$
dif{f}erent from $(\emptyset,\emptyset)$.\\
The functors giving rise to this equivalence will be denoted by
$\mc{P}:2-\Dl_{+}^{\comp}\mc{D}\rightarrow
2-\widehat{\Dl}_+^\comp\mc{D}$ and
$\mc{I}:2-\widehat{\Dl}_+^\comp\mc{D}\rightarrow
2-\Dl_{+}^{\comp}\mc{D}$.
\end{num}

Next we will exhibit two examples of biaugmented bisimplicial
object. F{i}rst, we introduce the ``total decalage'' object
associated with a simplicial object (cf. \cite{IlII}, p.7).

\begin{num}
The ordered sum of sets (see \ref{DeltagorroMasMonoidal}),
$+:\widehat{\Dl}_+\times\widehat{\Dl}_+\rightarrow\widehat{\Dl}_+$,
can be restricted to
$$\wp:2-\Dl_+\rightarrow \Dl \ ,$$
with $\wp([n],[m])=[n]+[m]=[n+m+1]$, where $[n]$ is identif{i}ed
with $\{0,\ldots,n\}\subset [n+m+1]$ and $[m]$ with
$\{n+1,\ldots,n+m+1\}\subset [n+m+1]$.
\end{num}
\begin{ej}\label{ejDecaladoTotal}
The functor \textit{total decalage}\index{Index}{total decalage}
$Dec: \simp\mc{D}\rightarrow
2-\Dl_+^\comp\mc{D}$\index{Symbols}{$Dec$} is def{i}ned by
composition with
$\wp:2-\Dl_+\rightarrow \Dl$.\\
If $X$ is a simplicial object, the biaugmented bisimplicial object
$Dec(X)$ consists of
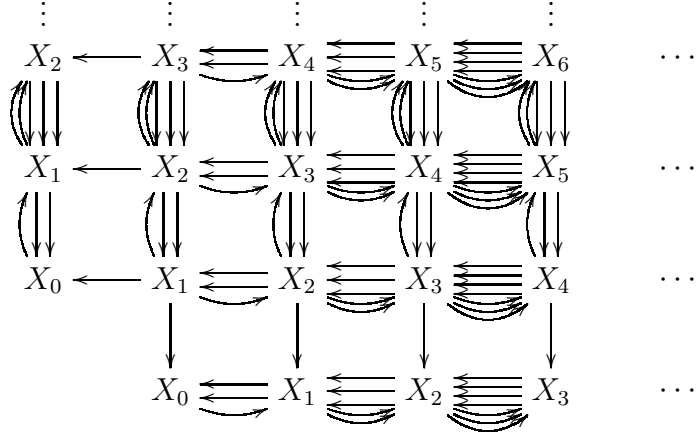
\vspace{-0.5cm}
$$\xymatrix@M=4pt@H=4pt@C=25pt{                                   &                                                                                             &                                                                                                                                              &                                                                                                                                                          &                                                                                                                             & \\
                                      X_{2}\ar@{}[u]|{\vdots}\ar@<0ex>[d]\ar@<1ex>[d] \ar@<-1ex>[d]      & {X_{3}}\ar[l]\ar@<0ex>[d]\ar@<1ex>[d] \ar@<-1ex>[d]\ar@{}[u]|{\vdots}\ar@/_0.75pc/[r]     & {X_{4}} \ar@{}[u]|{\vdots}\ar@<0.5ex>[l] \ar@<-0.5ex>[l]\ar@<0ex>[d]\ar@<1ex>[d] \ar@<-1ex>[d]\ar@/_1pc/[r]\ar@/_0.75pc/[r]                &X_{5} \ar@{}[u]|{\vdots}\ar@<0ex>[l]\ar@<1ex>[l] \ar@<-1ex>[l]\ar@<0ex>[d]\ar@<1ex>[d] \ar@<-1ex>[d]\ar@/_1pc/[r]\ar@/_0.75pc/[r]\ar@{-}@/_1.25pc/[r]   & X_{6}\ar@{}[u]|{\vdots}\ar@<0.33ex>[l]\ar@<-0.33ex>[l]\ar@<1ex>[l]\ar@<-1ex>[l]\ar@<0ex>[d]\ar@<1ex>[d] \ar@<-1ex>[d]     & \cdots \\
                                      X_{1}\ar@<0.5ex>[d] \ar@<-0.5ex>[d] \ar@/^1pc/[u]\ar@/^0.75pc/[u]  & {X_{2}}\ar[l]\ar@<0.5ex>[d] \ar@<-0.5ex>[d] \ar@/^1pc/[u]\ar@/^0.75pc/[u]\ar@/_0.75pc/[r] &  X_{3} \ar@/^1pc/[u]\ar@/^0.75pc/[u]\ar@<0.5ex>[l]\ar@<-0.5ex>[l]\ar@<0.5ex>[d] \ar@<-0.5ex>[d]\ar@/_1pc/[r]\ar@/_0.75pc/[r]               &X_{4} \ar@/^1pc/[u]\ar@/^0.75pc/[u] \ar@<0ex>[l]\ar@<1ex>[l] \ar@<-1ex>[l] \ar@<0.5ex>[d] \ar@<-0.5ex>[d]\ar@/_1pc/[r]\ar@/_0.75pc/[r]\ar@/_1.25pc/[r]  & X_{5}\ar@<0.33ex>[l]\ar@<-0.33ex>[l]\ar@<1ex>[l]\ar@<-1ex>[l] \ar@<0.5ex>[d] \ar@<-0.5ex>[d]\ar@/^1pc/[u]\ar@/^0.75pc/[u] & \cdots \\
                                      X_{0}  \ar@/^0.75pc/[u]                                            &  X_{1}\ar[d] \ar[l]\ar@/_0.75pc/[r] \ar@/^0.75pc/[u]                                      & {X_{2}} \ar[d] \ar@<0.5ex>[l] \ar@<-0.5ex>[l]  \ar@/^0.75pc/[u]  \ar@/_1pc/[r]\ar@/_0.75pc/[r]                                             &X_{3} \ar[d] \ar@<0ex>[l]\ar@<1ex>[l] \ar@<-1ex>[l]\ar@/^0.75pc/[u]\ar@/_1pc/[r]\ar@/_0.75pc/[r]\ar@/_1.25pc/[r]                                        & X_{4}\ar[d] \ar@<0.33ex>[l]\ar@<-0.33ex>[l]\ar@<1ex>[l]\ar@<-1ex>[l] \ar@/^0.75pc/[u]                                     & \cdots \\
                                                                                                         &  X_0 \ar@/_0.75pc/[r]                                                                        & X_1  \ar@<0.5ex>[l]\ar@<-0.5ex>[l]\ar@/_1pc/[r]\ar@/_0.75pc/[r]                                                                         &X_2 \ar@<0ex>[l]\ar@<1ex>[l] \ar@<-1ex>[l]\ar@/_1pc/[r]\ar@/_0.75pc/[r]\ar@/_1.25pc/[r]                                                                   & X_3\ar@<0.33ex>[l]\ar@<-0.33ex>[l]\ar@<1ex>[l]\ar@<-1ex>[l]                                                                 & \cdots }$$
\mbox{}\\
Following the notations introduced in \ref{ejempDecalado}, the
rows of the diagram are the augmented simplicial objects
$dec^k(X)\rightarrow X_{k-1}$, $k\geq 1$, where
$dec^k(X)=dec^1(\stackrel{k)}{\cdots}dec^1(X))$ is obtained from
$X$ by forgetting the last $k$ face and degeneracy maps, that is
$$
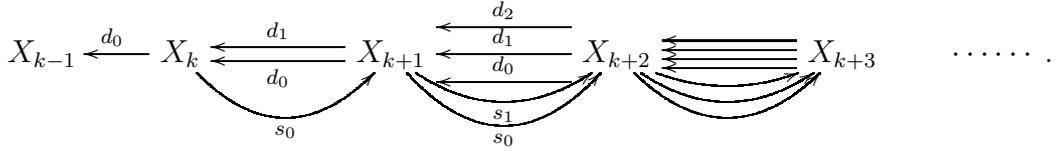
\xymatrix@1{ X_{k-1}{\,} & {\;} X_k \; \ar[l]_{d_0} \ar@/_2pc/[rr]_{s_0} && {\;} X_{k+1} \; \ar@<0.5ex>[ll]^-{d_0} \ar@<-0.5ex>[ll]_-{d_1} \ar@/_2.25pc/[rr]_{s_0} \ar@/_1.5pc/[rr]_{s_1} && \; X_{k+2} \; \ar@<0ex>[ll]_{d_1} \ar@<2ex>[ll]_-{d_0} \ar@<-2ex>[ll]_-{d_2} \ar@/_1.5pc/[rr] \ar@/_2pc/[rr] \ar@/_1pc/[rr] && {\;} X_{k+3} \ar@<0.33ex>[ll]\ar@<-0.33ex>[ll]\ar@<1ex>[ll]\ar@<-1ex>[ll]&\cdots\cdots}\ .$$
In the same way, the columns are the augmented simplicial objects
$dec_k(X)\rightarrow X_{k-1}$, this time obtained by forgetting
the f{i}rst face and degeneracy maps of $X$, that is
$$
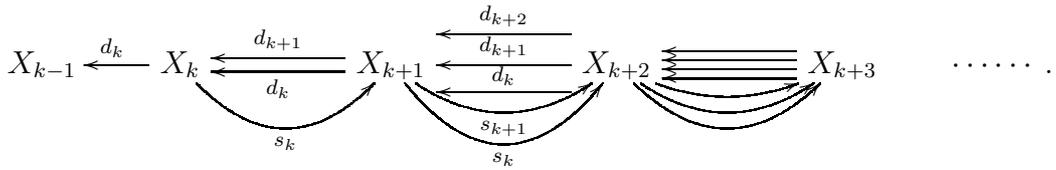
\xymatrix@1{ X_{k-1}{\,} & {\;} X_k \; \ar[l]_{d_k} \ar@/_2pc/[rr]_{s_k} && {\;} X_{k+1} \; \ar@<0.5ex>[ll]^-{d_k} \ar@<-0.5ex>[ll]_-{d_{k+1}} \ar@/_2.5pc/[rr]_{s_k} \ar@/_1.5pc/[rr]_{s_{k+1}} && \; X_{k+2} \; \ar@<0ex>[ll]_{d_{k+1}} \ar@<2ex>[ll]_-{d_{k}} \ar@<-2.25ex>[ll]_-{d_{k+2}} \ar@/_1.5pc/[rr] \ar@/_2pc/[rr] \ar@/_1pc/[rr] && {\;} X_{k+3}\ar@<0.33ex>[ll]\ar@<-0.33ex>[ll]\ar@<1ex>[ll]\ar@<-1ex>[ll]&\cdots\cdots }\ .$$
Analogously, $\widehat{Dec}:\widehat{\Dl}^\comp\mc{D}\rightarrow
2-\widehat{\Dl}^\comp\mc{D}$ is def{i}ned as
$[\widehat{Dec}(X)](E,F)=X(E+F)$, and it holds that
$Dec=\mc{I}\widehat{Dec}\mc{P}$ (see \ref{equiv2DlyDlGorro},
\ref{equivDlyDlGorro}).
\end{ej}

\begin{ej}\label{ejObjBisimpfEpsilon}
Let $f:X\rightarrow Y$ be a morphism in $\simp\mc{D}$ and
$\epsilon:X\rightarrow X_{-1}\times\Dl$ an augmentation of $X$.\\
Hence, $X$ gives rise to the bisimplicial object $Z^{+}=\Dl\times
X$, that is, $Z^+_{i,j}=X_j$. Moreover, $f$ and $\epsilon$
def{i}ne the augmentations $Z^+_{i,j}\rightarrow Z_{-1,j}=Y_j$ and
$Z^+_{i,j}\rightarrow Z_{i,-1}=X_{-1}$.\\
Therefore the diagram $Z=Z(f,\epsilon):\ Z_{-1,\cdot} \leftarrow
Z^+_{\cdot,\cdot}\rightarrow Z_{\cdot,-1}$ is a biaugmented
bisimplicial object.
Consequently, given $([n],[m])$, $([n'],[m'])$ in $2-\Dl_+$ and a
morphism $(\theta,\theta'):([n'],[m'])\rightarrow ([n],[m])$, we
have that
$$ Z_{n,m}=\begin{cases}
           Y_m & \mbox{ if } n=-1 \\
           X_m & \mbox{ if } n\neq -1
           \end{cases}\ \ Z(\theta,\theta')=\begin{cases}
                                            Y(\theta') & \mbox{{ if }} n=-1  \, (\mbox{hence }n'=-1)\\
                                            X(\theta') & \mbox{{ if }} n'\neq -1 \, (\mbox{hence }n\neq -1)\\
                                            f_{m'}X(\theta') & \mbox{{ if }} n'= -1 \mbox{ and }n\neq -1
                                            \end{cases}$$
Visually, $Z(f,\epsilon)$ is the following diagram
\begin{equation}\label{diagramaCilindro}\xymatrix@M=4pt@H=4pt@C=25pt{                                   &                                                                                             &                                                                                                                                              &                                                                                                                                                          &                                                                                                                             & \\
                                      Y_{2}\ar@{}[u]|{\vdots}\ar@<0ex>[d]\ar@<1ex>[d] \ar@<-1ex>[d]      & {X_{2}}\ar[l]\ar@<0ex>[d]\ar@<1ex>[d] \ar@<-1ex>[d]\ar@{}[u]|{\vdots}\ar@/_0.75pc/[r]     & {X_{2}} \ar@{}[u]|{\vdots}\ar@<0.5ex>[l] \ar@<-0.5ex>[l]\ar@<0ex>[d]\ar@<1ex>[d] \ar@<-1ex>[d]\ar@/_1pc/[r]\ar@/_0.75pc/[r]                &X_{2} \ar@{}[u]|{\vdots}\ar@<0ex>[l]\ar@<1ex>[l] \ar@<-1ex>[l]\ar@<0ex>[d]\ar@<1ex>[d] \ar@<-1ex>[d]\ar@/_1pc/[r]\ar@/_0.75pc/[r]\ar@{-}@/_1.25pc/[r]   & X_{2}\ar@{}[u]|{\vdots}\ar@<0.33ex>[l]\ar@<-0.33ex>[l]\ar@<1ex>[l]\ar@<-1ex>[l]\ar@<0ex>[d]\ar@<1ex>[d] \ar@<-1ex>[d]     & \cdots \\
                                      Y_{1}\ar@<0.5ex>[d] \ar@<-0.5ex>[d] \ar@/^1pc/[u]\ar@/^0.75pc/[u]  & {X_{1}}\ar[l]\ar@<0.5ex>[d] \ar@<-0.5ex>[d] \ar@/^1pc/[u]\ar@/^0.75pc/[u]\ar@/_0.75pc/[r] &  X_{1} \ar@/^1pc/[u]\ar@/^0.75pc/[u]\ar@<0.5ex>[l]\ar@<-0.5ex>[l]\ar@<0.5ex>[d] \ar@<-0.5ex>[d]\ar@/_1pc/[r]\ar@/_0.75pc/[r]               &X_{1} \ar@/^1pc/[u]\ar@/^0.75pc/[u] \ar@<0ex>[l]\ar@<1ex>[l] \ar@<-1ex>[l] \ar@<0.5ex>[d] \ar@<-0.5ex>[d]\ar@/_1pc/[r]\ar@/_0.75pc/[r]\ar@/_1.25pc/[r]  & X_{1}\ar@<0.33ex>[l]\ar@<-0.33ex>[l]\ar@<1ex>[l]\ar@<-1ex>[l] \ar@<0.5ex>[d] \ar@<-0.5ex>[d]\ar@/^1pc/[u]\ar@/^0.75pc/[u] & \cdots \\
                                      Y_{0}  \ar@/^0.75pc/[u]                                            &  X_{0}\ar[d] \ar[l]\ar@/_0.75pc/[r] \ar@/^0.75pc/[u]                                      & {X_{0}} \ar[d] \ar@<0.5ex>[l] \ar@<-0.5ex>[l]  \ar@/^0.75pc/[u]  \ar@/_1pc/[r]\ar@/_0.75pc/[r]                                             &X_{0} \ar[d] \ar@<0ex>[l]\ar@<1ex>[l] \ar@<-1ex>[l]\ar@/^0.75pc/[u]\ar@/_1pc/[r]\ar@/_0.75pc/[r]\ar@/_1.25pc/[r]                                        & X_{0}\ar[d] \ar@<0.33ex>[l]\ar@<-0.33ex>[l]\ar@<1ex>[l]\ar@<-1ex>[l] \ar@/^0.75pc/[u]                                     & \cdots \\
                                                                                                         &  X_{-1} \ar@/_0.75pc/[r]                                                                        & X_{-1} \ar@<0.5ex>[l]\ar@<-0.5ex>[l]\ar@/_1pc/[r]\ar@/_0.75pc/[r]                                                                             &X_{-1} \ar@<0ex>[l]\ar@<1ex>[l] \ar@<-1ex>[l]\ar@/_1pc/[r]\ar@/_0.75pc/[r]\ar@/_1.25pc/[r]                                                                   & X_{-1}\ar@<0.33ex>[l]\ar@<-0.33ex>[l]\ar@<1ex>[l]\ar@<-1ex>[l]                                                                 & \cdots
                                                                                                         }\end{equation}
\mbox{}\\
where the horizontal morphisms are either $f_n:X_n\rightarrow Y_n$
or the identity $X_n\rightarrow X_n$, whereas the vertical
morphisms are either the face and degeneracy maps of $X$
or $\epsilon_0=d_0:X_0\rightarrow X_{-1}$.\\
Clearly $Z(f,\epsilon)$ is natural with respect to $(f,\epsilon)$.
In addition, since $\epsilon:X\rightarrow X_{-1}\times\Dl$ is
characterized by $\epsilon_0$ (see \ref{aumentvisionequivalente}),
then $Z(f,\epsilon)$ preserves all the information contained in
$(f,\epsilon)$.
\end{ej}

\begin{def{i}}[Total object of a biaugmented bisimplicial object]\label{def{i}TotCombinatorio}\index{Index}{total!of a biaug. bisimpl. object}\index{Symbols}{$Tot(Z)$}\mbox{}\\
The functor $Tot:2-\Dl_{+}^{\comp}\mc{D}\rightarrow \simp\mc{D}$
is defined as follows.
If $Z$ is a biaugmented bisimplicial object, then $Tot(Z)$ is in
degree $n$ the coproduct of the $n$-th antidiagonal of
(\ref{diagrObjBisimpBiaum}), that is
$$Tot(Z)_n=\ds\coprod_{i+j=n-1}Z_{i,j}\ .$$
The face maps in $Tot(Z)$ are def{i}ned as coproducts of the
morphisms in (\ref{diagrObjBisimpBiaum}) going from the $n$-th
antidiagonal to the
 $n-1$-th one, and analogously for the degeneracy maps.\\
More specif{i}cally, set $\theta^{(1)}=Z(\theta,
Id):Z_{n,k}\rightarrow Z_{m,k}$ and
$\theta^{(2)}=Z(Id,\theta):Z_{k,n}\rightarrow Z_{k,m}$
for every $\theta:[m]\rightarrow [n]$ and $k\geq -1$.\\
The face maps $d_k:Tot(Z)_n\rightarrow Tot(Z)_{n-1}$ are given by
$$ d_k |_{Z_{i,j}} = \left\{ \begin{array}{ll} d_{k}^{(1)}:Z_{i,j}\rightarrow Z_{i-1,j} & \mbox{ if } i\geq k \\
                                               d^{(2)}_{k-i-1}:Z_{i,j}\rightarrow Z_{i,j-1} & \mbox{ if } i < k \, , \end{array}\right. \ .$$
The degeneracy maps $s_k:Tot(Z)_n\rightarrow Tot(Z)_{n+1}$ are
$$ s_k |_{Z_{i,j}} = \left\{ \begin{array}{ll} s_{k}^{(1)}:Z_{i,j}\rightarrow Z_{i+1,j} & \mbox{ if } i\geq k \\
                                              s^{(2)}_{k-i-1}:Z_{i,j}\rightarrow Z_{i,j+1} & \mbox{ if } i < k \, , \end{array}\right. \ .$$
\end{def{i}}
\begin{num}\label{inclusionesTot}
Consider $Z\in (2-\Dl_{+})^\comp\mc{D}$. The canonical maps
$$Z_{-1,n}\rightarrow \ds\coprod_{i+j=n-1}Z_{i,j}\ \mbox{ and }\ Z_{n,-1}\rightarrow \ds\coprod_{i+j=n-1}Z_{i,j}$$
give rise to the following canonical morphisms in $\simp\mc{D}$
$$Z_{-1,\cdot}\rightarrow Tot(Z)\leftarrow Z_{\cdot,-1}\ ,$$
that are natural in $Z$.
\end{num}
\begin{num}\label{def{i}TotAumentado}
The functor $Tot$ can be extended to
$Tot_+:(\Dl_+\times\Dl_+)^{\comp}\mc{D}\rightarrow(\Dl_+)^{\comp}\mc{D}$ as follows.\\
If ${}_+Z\in (\Dl_+\times\Dl_+)^{\comp}\mc{D}$ and $Z${ if }the
restriction of ${}_+Z$ to $2-\Dl_{+}$ then the morphism
$$d_0=d_0^{(2)}\sqcup d_0^{(1)}:Z_{-1,0}\sqcup
Z_{0,-1}\rightarrow {}_+Z_{-1,-1}$${ is }an augmentation of
$Tot(Z)$.
\end{num}

\begin{obs}
As far as the author knows, the functors $Tot$ and $Tot_+$ do not
appear in the literature.\\
However, a particular case of $Tot_+$ is introduced in \cite{EP},
that is called the ``join'' of two augmented simplicial sets. If
$X\rightarrow X_{-1}$ and $Y\rightarrow Y_{-1}$ are augmented
simplicial sets then their join is just the image under the total
functor of ${}_+Z=\{X_n \times Y_m\}_{n,m\geq -1}$. In loc. cit.
the associativity of this join is stated. This associativity can
be deduced as well from (the augmented version) of
proposition \ref{TotGorroIterado}.\\
Hence, the notion of join is completely similar to the one of
tensor product of two chain complexes $A$ and $B$ of modules over
a commutative ring $R$, since the tensor product of $A$ and $B$ is
just the total chain complex associated with the double complex
$\{A_n\otimes_{R} B_m\}_{n,m}$.
\end{obs}

The def{i}nition given above for $Tot$ is purely combinatorial, we
introduced it in this way because we will need the formulae for
computations. However, this construction can be understood in a
more intuitive way{ if }we consider the following equivalent
def{i}nitions.

\begin{num}\label{def{i}Tot} The functor $Tot:2-\Dl_{+}^{\comp}\mc{D}\rightarrow
\simp\mc{D}$ can be described as follows.

\begin{itemize}
\item Consider the category $\Dl/[1]$\index{Symbols}{$\Dl/[1]$}
with objects the pairs $([n],\sigma)$, where
$\sigma:[n]\rightarrow [1]$ is a morphism in $\Dl$. A morphism
$\theta:([n],\sigma)\rightarrow ([m],\rho)$ is
$\theta:[n]\rightarrow [m]$ in $\Dl$ such that
$\rho\theta=\sigma$. We will denote $([n],\sigma)$ by $\sigma$ if
$[n]$ is understood.

\item If $\sigma:[n]\rightarrow [1]$, let $i_\sigma\in\{0,\ldots
,n+1\}$ be such that $\{i_\sigma,\ldots,n\}=\sigma^{-1}(1)${ if
}it is non empty, and $i_\sigma=n+1${ if
}$\sigma^{-1}(1)=\emptyset$. Note that the correspondence
$\sigma\rightarrow i_\sigma$ is a bijection. Moreover, we will
identify $[i_\sigma-1]$ with $\sigma^{-1}(0)$, as well as
$[n-i_{\sigma}]$ with $\sigma^{-1}(1)$ (after relabelling in a
suitable way).

\item Let $\Psi:\Dl/[1]\rightarrow 2-\Dl_{+}$ be the functor def{i}ned by $\Psi(\sigma)=[i_\sigma-1]\times[n-i_\sigma]$.\\
Given $\theta:([n],\sigma)\rightarrow ([m],\rho)$, since
$\theta(\sigma^{-1}(j))\subseteq \rho^{-1}(j)$ for $j=0,1$ then
$\theta|_{\sigma^{-1}(0)}:[i_\sigma-1]\rightarrow [i_{\rho}-1]$,
and $\theta|_{\sigma^{-1}(1)}:[n-i_{\sigma}]\rightarrow
[m-i_{\rho}]$
are morphisms in $\Dl_+$.\\
Hence, def{i}ne
$\Psi(\theta)=\theta|_{\sigma^{-1}(0)}\times\theta|_{\sigma^{-1}(1)}:[i_\sigma-1]\times[n-i_\sigma]\rightarrow
[i_\rho-1]\times[m-i_\rho]$.\\

\item If $Z$ is a biaugmented bisimplicial object, $Tot(Z)$ is the
simplicial object def{i}ned in degree $n$ as
$$ Tot(Z)_n=\ds\coprod_{\sigma:[n]\rightarrow [1]} Z_{i_\sigma-1,n-i_\sigma}$$
If $\rho: [n]\rightarrow [1]$ and $\theta:[n]\rightarrow [m]$,
then $\theta:([n],\rho\theta)\rightarrow ([m],\rho)$ is in
$\Dl/[1]$, and the restriction of $Tot(Z)(\theta)$ to
$Z_{i_\rho-1,m-i_\rho}$ is
$$Z(\Psi(\theta)):Z_{i_\rho-1,m-i_\rho}\rightarrow Z_{i_{\rho\theta}-1,n-i_{\rho\theta}}$$
\end{itemize}
\end{num}
\begin{num}\label{def{i}TotGorro} In terms of $\widehat{\Dl}_+$,
$\widehat{Tot}:2-\widehat{\Dl}_{+}^{\comp}\mc{D}\rightarrow
\widehat{\Dl}^{\comp}\mc{D}$ is def{i}ned as
$$ [\widehat{Tot} (Z)](E)= \ds\coprod_{E= E_0 + E_1} Z(E_0,E_1) \ ,$$
where the coproduct is indexed over the set of pairs $(E_0,E_1)\in
2-\widehat{\Dl}_{+}$ such that $E$ is the ordered sum of $E_0$ and
$E_1$ (see \ref{DeltagorroMasMonoidal}).

Let $f:E'\rightarrow E$ be a morphism in $\widehat{\Dl}$ and
$(E_0,E_1)$ an object in $2-\widehat{\Dl}_{+}$ such that
$E_0+E_1=E$.\\
Denote by $E'_i$ the set $f^{-1}(E_i)$ with the order induced by
the one of $E'$, for $i=0,1$. It is clear that $E'=E'_0+E'_1$, as
well as $f=f_0+f_1$, where $f_i=f|_{E'_i}:E'_i\rightarrow
E_i$ for $i=0,1$.\\
Then $[\widehat{Tot} (Z)](f):[\widehat{Tot} (Z)](E)\rightarrow
[\widehat{Tot} (Z)](E')$ is def{i}ned as
$$[\widehat{Tot} (Z)](f)|_{Z(E_0,E_1)}= Z(f_0,f_1):Z(E_0,E_1)\rightarrow Z(E'_0,E'_1) \ .$$
\end{num}

\begin{num}\label{def{i}TotGorroAumentado}
The functor $\widehat{Tot}$ can be also extended in a natural way
to
$$\widehat{Tot}^+:(\widehat{\Dl}_+\times \widehat{\Dl}_+)^\comp\mc{D}\rightarrow \widehat{\Dl}_+^\comp\mc{D}\ ,$$
with $[\widehat{Tot}^+(Z)](\emptyset)=Z(\emptyset,\emptyset)$.
\end{num}

\begin{prop}
The two def{i}nitions given for $Tot$ coincide, and do def{i}ne a
functor
$Tot:2-\Dl_{+}^{\comp}\mc{D}\rightarrow \simp\mc{D}$.\\
In addition, these constructions correspond to
$\widehat{Tot}:2-\widehat{\Dl}_{+}^{\comp}\mc{D}\rightarrow
\widehat{\Dl}^{\comp}\mc{D}$ under the canonical equivalence of
the categories $\Dl_+$ and $\widehat{\Dl}_+$. That is, under the
notations of \ref{equiv2DlyDlGorro} and \ref{equivDlyDlGorro} we
have that
$$ Tot= \mc{I}\comp\widehat{Tot}\comp \mc{P}\ .$$
\end{prop}

\begin{proof}
Let us see first the second statement.\\
Define $\widehat{\Dl}/[1]$ as \ref{def{i}Tot}, as well as
$\widehat{\Psi}:\widehat{\Dl}/[1]\rightarrow 2-\widehat{\Dl}_{+}$.
More specif{i}cally,{ if }$\sigma:E\rightarrow [1]$, then
$\widehat{\Psi}(\sigma)=\sigma^{-1}(0)\times\sigma^{-1}(1)$.\\
Any decomposition $E=E_0+E_1$ in $2-\widehat{\Dl}_+$ determines in
a unique way an object $\sigma:E\rightarrow [1]$ in
$\widehat{\Dl}/[1]$ (just take $\sigma(E_i)=i$, $i=0,1$).\\
Hence \ref{def{i}TotGorro} can be rewritten as
$$[\widehat{Tot} (Z)](E)= \ds\coprod_{\sigma:E\rightarrow [1]} Z(\sigma^{-1}(0),\sigma^{-1}(1))\ .$$
Moreover,{ if }$f:E'\rightarrow E$, then the restriction of
$[Tot(Z)](f)$ to $Z(\sigma^{-1}(0),\sigma^{-1}(1))$ coincides with
$Z(f|_{(f\sigma)^{-1}(0)},f|_{(f\sigma)^{-1}(1)})$.
Consequently
$$[\mc{I}\comp\widehat{Tot}\comp \mc{P}(Z)]([n])=\ds\coprod_{\sigma:[n]\rightarrow [1]} (\mc{P}Z)(\sigma^{-1}(0),\sigma^{-1}(1))=\ds\coprod_{\sigma:[n]\rightarrow [1]} Z(p(\sigma^{-1}(0)),p(\sigma^{-1}(1)))\ ,$$
that is equal to $Tot(Z)_n $, and clearly the action of
$\mc{I}\comp\widehat{Tot}\comp \mc{P}(Z)$ and of $Tot(Z)$ over the
morphisms of $\Dl$ coincides.

Secondly, let us check that the two def{i}nitions given for the
total functor coincide. Note that any monotone function
$\sigma:[n]\rightarrow [1]$ is characterized by the integer
$i_\sigma\in\{0,n+1\}$ such that $\sigma^{-1}(1)=\{i_\sigma,\ldots
,n\}$.\\ It follows that the correspondence $\sigma\rightarrow
i_{\sigma}-1$ is one-to-one between the sets
$$\{(i,j)\ | \ i+j=n-1,\ i,j\geq -1\}\mbox{ and } \{(i_\sigma-1,n-i_\sigma)\ | \ \sigma:[n]\rightarrow [1]\}\ .$$
Let $\theta=\partial_k:[n-1]\rightarrow [n]$ and
$\rho:[n]\rightarrow [1]$ be morphisms in $\Dl$. We will compute
$$\Psi(\partial_k):[i_{\rho\partial_k}-1]\times[n-1-i_{\rho\partial_k}]\rightarrow [i_\rho-1]\times[n-i_\rho]\ .$$
\indent If $i_\rho\leq k$, then $i_{\rho\partial_k}=i_\rho$ and it
holds that
$$\partial_k|_{{(\rho\partial_k)}^{-1}(0)}=Id:[i_\rho-1]\rightarrow [i_\rho-1]\mbox{ ; }\partial_k|_{{(\rho\partial_k)}^{-1}(1)}=\partial_{k-i_\rho}:[n-i_\rho-1]\rightarrow [n-i_\rho] \ .$$
In this case
$[Tot(Z)](\partial_k)|_{Z_{i_\rho-1,n-i_\rho}}=d^{(2)}_{k-i_\rho}:Z_{i_\rho-1,n-i_\rho}\rightarrow
Z_{i_{\rho}-1,n-i_{\rho}-1}$.\\
Then, taking $i=i_\rho-1$,
$[Tot(Z)](\partial_k)|_{Z_{i_\rho-1,n-i_\rho}}$ corresponds to
$d_{k-i-1}^{(2)}:Z_{i,j}\rightarrow Z_{i,j-1}$.

On the other hand,{ if }$i_\rho >k$, then
$i_{\rho\partial_k}=i_\rho-1$ and it can be checked analogously
that
$[Tot(Z)](\partial_k)|_{Z_{i_\rho-1,n-i_\rho}}=Z(\partial_k,Id):Z_{i_\rho-1,n-i_\rho}\rightarrow
Z_{i_{\rho}-2,n-i_{\rho}}$.\\
Hence, setting $i=i_\rho-1$ we have that
$[Tot(Z)](\partial_k)|_{Z_{i_\rho-1,n-i_\rho}}$ is
$d_{k}^{(1)}:Z_{i,j}\rightarrow Z_{i-1,j}$.\\
The analogous equality involving the degeneracy maps
$\sigma_k:[n+1]\rightarrow [n]$ can be checked in a similar way.

It remains to see that $Tot(Z)\in\simp\mc{D}$, that is, that
$Tot(Z)(\theta\theta')=Tot(Z)(\theta')Tot(Z)(\theta)$ for every
composable morphisms $\theta$ and $\theta'$ in $\Dl$. The equality
$Tot(Z)(Id)=Id$ follows from the naturality of $\Psi$.

F{i}nally, given $\psi:Z\rightarrow S$ then
$Tot(\psi):Tot(Z)\rightarrow Tot(S)$ is def{i}ned in degree $n$ by
$Tot(\psi)|_{Z_{i,j}}=\psi_{i,j}:Z_{i,j}\rightarrow S_{i,j}$.
Clearly, it is a morphism between simplicial objects, natural in
$\psi$.
\end{proof}

\begin{obs}
Let $\mc{I}:\widehat{\Dl}_+^\comp\mc{D}\rightarrow
{\Dl}_+^\comp\mc{D}$ and $\mc{P}:({\Dl}_+\times
{\Dl}_+)^\comp\mc{D}\rightarrow
(\widehat{\Dl}_+\times\widehat{\Dl}_+)^\comp\mc{D}$ be the
equivalences of categories induced by
$i:\Dl_+\rightarrow\widehat{\Dl}_+$ and
$p:\widehat{\Dl}_+\rightarrow
\Dl_+$.\\
Hence, it holds that $ Tot^+= \mc{I}\comp\widehat{Tot}^+\comp
\mc{P}$.
\end{obs}

\begin{prop}\label{relacionTotyTotGorro}
The functor
$\widehat{Tot}:2-\widehat{\Dl}_{+}^{\comp}\mc{D}\rightarrow
\widehat{\Dl}^\comp\mc{D}$ is left adjoint to the total decalage
functor $\widehat{Dec}:\simp\mc{D}\rightarrow
2-\Dl_{+}^{\comp}\mc{D}$ introduced in \ref{ejDecaladoTotal}.\\
Hence $(Tot,Dec)$ is also an adjoint pair of functors.\\
Consequently, the functors $\widehat{Tot}$ and $Tot$ commutes with
colimits and in particular with coproducts.
\end{prop}

\begin{proof}
Since $\mc{P}$ and $\mc{I}$ are quasi-inverse equivalences of
categories, and since $Tot=\mc{I}\comp\widehat{Tot}\comp \mc{P}$
and $Dec=\mc{I}\comp\widehat{Dec}\comp \mc{P}$, it follows from
the adjunction $(\widehat{Tot},\widehat{Dec})$ that
$(Tot,Dec)$ is also an adjoint pair.\\
Consider $Z\in 2-\widehat{\Dl}^\comp_+\mc{D}$ and $Y\in
\widehat{\Dl}^\comp\mc{D}$. Let us check that there is a canonical
bijection
$$\mrm{Hom}_{\widehat{\Dl}^\comp\mc{D}}\left( \widehat{Tot}(Z), Y \right) \simeq \mrm{Hom}_{2-\widehat{\Dl}^\comp_+\mc{D}}\left(Z,\widehat{Dec}(Y)\right)\ .$$
A morphism $F:\widehat{Tot}(Z)\rightarrow Y$ consists of a
collection of morphisms in $\mc{D}$
$$ G(E_0,E_1)=F(E)|_{Z(E_0,E_1)}:Z(E_0,E_1)\rightarrow Y(E)=[\widehat{Dec}(Y)](E_0,E_1)$$
for every equality of ordered sets $E=E_0+E_1$.\\
Moreover, $G$ is natural in $(E_0,E_1)$ since $F$ is natural in
$E$.
To see this, given an expression $E=E_0+E_1$, it is enough to note
that there exists a bijection between the set of morphisms
$(f_0,f_1):(E'_0,E'_1)\rightarrow (E_0,E_1)$ in
$2-\widehat{\Dl}_+$ and the morphisms $f:E'\rightarrow E$. To see
this, set $E'=E'_0+E'_1$; $f=f_0+f_1:E'\rightarrow E$ on one hand,
and $E'_i=f^{-1}(E_i)$; $f_i=f|_{E'_i}$, $i=0,1$ on the other hand.\\
The naturality of $F$ means that for every ordered set $E$ the
following equality holds
$$ F(E')\comp \widehat{Tot}(Z)(f) = Y(f)\comp F(E) :[\widehat{Tot}(Z)](E)\rightarrow Y(E') \ ,$$
and this happens{ if }and only{ if }this equality holds on each
component $Z(E_0,E_1)$ of $\widehat{Tot}(Z)$. That is to say,{ if }and only{ if }%
\begin{eqnarray} G(E'_0,E'_1) \comp Z(f_0,f_1) = (F(E')\comp \widehat{Tot}(Z)(f))|_{Z(E_0,E_1)}=\nonumber\\
                     = (Y(f)\comp F(E))|_{Z(E_0,E_1)}=[\widehat{Dec}(Y)](f_0,f_1)\comp G(E_0,E_1)\nonumber
\end{eqnarray}
and this happens{ if }and only{ if }$G$ is natural with respect to
each morphism $(f_0,f_1)$.
\end{proof}

\begin{obs} In \cite{CR} the authors consider a ``decalage'' functor $dec:\simp Set\rightarrow\simp\simp Set$, by forgetting the
biaugmentation  of the biaugmented bisimplicial object $Dec(X)$
associated with $X\in\simp Set$. In loc. cit. another
``combinatorial'' functor $\simp\simp Set\rightarrow Set$ is also
introduced. This functor is right adjoint to $dec$ and is
def{i}ned using f{i}ber products in $Set$ instead of coproducts.
\end{obs}

Next we will see that the classical simplicial cone and cylinder
objects are particular cases of the functor $Tot$.
\begin{ej}[Simplicial cone]\label{conoSimpl}\mbox{}\\
Assume that $\mc{D}$ has a f{i}nal object $1$ and $f:X\rightarrow
Y$ is a morphism in $\simp\mc{D}$. The classical def{i}nition of
cone object associated with $f$ \cite{DeIII} is
$Cf\in\simp\mc{D}$\index{Index}{cono!simplicial}\index{Symbols}{$Cf$},
with
$$(Cf)_n=Y_n\sqcup X_{n-1}\sqcup\cdots\sqcup X_{0}\sqcup 1\ .$$
Let $Z$ be the biaugmented bisimplicial object associated with
$f:X\rightarrow Y$ and to the trivial augmentation $X\rightarrow
1\times\Dl$ (see \ref{ejObjBisimpfEpsilon}).\\
Hence $Cf$ is just the total simplicial object of $Z$.
\end{ej}

\begin{ej}['Cubical' cylinder]\label{ejCilindroCub}\mbox{}\\
Given a simplicial object $X$ in $\mc{D}$, let us remind the
classical notion of cylinder associated with $X$, that will be
denoted by
$\widetilde{Cyl}(X)$\index{Symbols}{$\widetilde{Cyl}(X)$}. We will
say that $\widetilde{Cyl}(X)$ is the ``cubical'' cylinder object
of $X$. It is def{i}ned in degree $n$ as
$$\widetilde{Cyl}(X)_n=\ds\coprod_{\sigma:[n]\rightarrow [1]}X_n^\sigma\ ,$$
where $X_n^\sigma=X_n$ $\forall \sigma$. Given
$\theta:[m]\rightarrow [n]$ in $\Dl$,
$\widetilde{Cyl}(\theta):\widetilde{Cyl}(X)_n\rightarrow
\widetilde{Cyl}(X)_m$ is
$$\widetilde{Cyl}(\theta)|_{X^{\sigma}_n}=X(\theta):X^{\sigma}_n\rightarrow
X^{\sigma\theta}_m \ .$$
It holds that $\widetilde{Cyl}(X)$ is the total object of
$Dec(X)\in (2-\Dl_+)^\comp\mc{D}$ (given in
\ref{ejDecaladoTotal}).
\end{ej}
\begin{obs}\label{HomotbyCylTilde} We recall that the cubical cylinder object
characterizes simplicial homotopies [\cite{May} prop. 6.2].\\
In other words, let $u_0,u_1:[n]\rightarrow[1]$ be the morphisms
such that $u_i([n])=i$, $i=0,1$. If $X$ is a simplicial object,
def{i}ne $I,J:X\rightarrow \widetilde{Cyl}(X)$ as
$$I_n=Id:X_n\rightarrow X_n^{u_1}\mbox{ and } J_n=Id:X_n\rightarrow X_n^{u_0}\ .$$
Given $f,g:X\rightarrow Y$ in $\simp\mc{D}$, we have that $f\sim
g$ (that is, $f$ is homotopic to $g$){ if }and only{ if }there
exists $H:\widetilde{Cyl}(X)\rightarrow X$ such that $H\comp I=f$
and $H\comp J=g$.
\end{obs}
%


\section{Total object of \textit{n}-augmented \textit{n}-simplicial objects}

In this section we will generalize in a natural way the functors
$Tot$ and $\widehat{Tot}$ to the categories $n-\Dl_+^\comp\mc{D}$
and $n-\widehat{\Dl}_+^\comp\mc{D}$, respectively. In addition,
the functor $Tot_n$ (resp. $\widehat{Tot}_n$) can be obtain as
well as iterations of $Tot$ (resp. $\widehat{Tot}$).

\begin{def{i}}
Let $n-\Dl_+$\index{Symbols}{$n-\Dl_+$} be the full subcategory of
$\Dl_+\times\stackrel{n)}{\cdots}\times\Dl_+$ whose objects are
$([m_0],\ldots,[m_{n-1}])$ with $\sum m_k\neq -n$. Analogously,
let $n-\widehat{\Dl}_+$ be the category def{i}ned using
$\widehat{\Dl}_+$ instead of ${\Dl}_+$.\\
Again, we will refer to the elements of the category
$n-\Dl_+^\comp\mc{D}$ of contravariant functors from $n-\Dl_+$ to
$\mc{D}$ as \textit{n}-augmented \textit{n}-simplicial objects.
\end{def{i}}

\begin{num}\label{objeto3Delta}
Similarly to the case $n=2$, an object $T$ of
$3-\Dl_+^\comp\mc{D}$ can be visualized as
$$\xymatrix@R=33pt@C=33pt@H=4pt@M=4pt@!0{                                           & T_{i,-1,k} \ar[rr]\ar[ld]           &                   & T_{i,-1,-1} \\
                                             T_{-1,-1,k}                            &                                               &                   &\\
                                                                                    & T_{i,j,k} \ar[uu] \ar[rr]\ar[ld]&                   & T_{i,j,-1}\ ,\ar[uu]\ar[ld]\\
                                             T_{-1,j,k} \ar[uu]   \ar[rr]           &                                               & T_{-1,j,-1} & }$$
where the indexes $i,j,k$ are positive integers.
\end{num}

\begin{def{i}} Def{i}ne
$Tot_n:n-\Dl_+^\comp\mc{D}\rightarrow\simp\mc{D}$\index{Symbols}{$Tot_n$}\index{Index}{total object!of \textit{n}-aug. \textit{n}-simpl. objects} as follows.\\
Let $\Dl/[n-1]$ be the category def{i}ned as $\Dl/[1]$ in
\ref{def{i}Tot}.\\
Given $\sigma:[m]\rightarrow [n-1]$ and $k\in [n-1]$, note that
each ordered subset $\sigma^{-1}(k)$ of $[m]$ is isomorphic to a
unique object
$[m_k(\sigma)]$ of $\Dl_+$.\\
If $Z\in n-\Dl_+^\comp\mc{D}$, set
$$ [Tot_n(Z)]_m=\ds\coprod_{\sigma:[m]\rightarrow [n-1]} Z_{m_0(\sigma),\ldots,m_{n-1}(\sigma)}\ .$$
Moreover, if $\theta:[m']\rightarrow [m]$, then
$\theta|_{(\sigma\theta)^{-1}(k)}:(\sigma\theta)^{-1}(k)\rightarrow
(\sigma)^{-1}(k)$ induces a monotone function
$\theta_k:[m'_k(\sigma\theta)]\rightarrow [m_k(\sigma)]$, and the
restriction of $[Tot_n(Z)](\theta)$ to
$Z_{m_0(\sigma),\ldots,m_{n-1}(\sigma)}$ is
$$Z(\theta_0,\ldots,\theta_{m-1}):Z_{m_0(\sigma),\ldots,m_{n-1}(\sigma)}\rightarrow Z_{m'_0(\sigma\theta),\ldots,m'_{n-1}(\sigma\theta)}\ .$$
\end{def{i}}

Next we provide a description of the total functor in terms of
$\widehat{\Dl}_+$.

\begin{num}
The functor
$\widehat{Tot}_{n}:n-\widehat{\Dl}_+^\comp\mc{D}\rightarrow
\widehat{\Dl}^\comp\mc{D}$ is given by
$$ [\widehat{Tot}_{n}(Z)](E)= \ds\coprod_{E= E_0 +\cdots +E_{n-1}} Z(E_0,\ldots,E_{n-1}) \ .$$
If $f:E'\rightarrow E$ is a morphism in $\widehat{\Dl}$ and $E=E_0
+\cdots +E_{n-1}$ then
$$E'=E'_0 +\cdots + E'_{n-1}\ \mbox{ and }\ f=f_0+\cdots + f_{n-1} \ , $$
where $E'_i=\theta^{-1}(E_i)$ and $f_i=f|_{E'_i}:E'_i\rightarrow
E_i$, for $i=0,\ldots ,n-1$.\\[0.1cm]
Hence $[\widehat{Tot}_n (Z)](f):[\widehat{Tot} (Z)](E)\rightarrow
[\widehat{Tot}_n (Z)](E')$ is
$$[\widehat{Tot}_n (Z)](f)|_{Z(E_0,\ldots,E_{n-1})}= Z(f_0,\ldots ,f_{n-1}):Z(E_0,\ldots,E_{n-1})\rightarrow Z(E'_0,\ldots,E'_{n-1}) \ .$$
\end{num}

\begin{obs}
As in \ref{equiv2DlyDlGorro}, the equivalences of categories
$i:\Dl\rightarrow \widehat{\Dl}$, $p:\widehat{\Dl}\rightarrow \Dl$
induce the quasi-inverse equivalences
$$\mc{I}:n-\widehat{\Dl}^\comp\mc{D}\rightarrow n-\simp\mc{D}\ \ \ \mc{P}:n-\simp\mc{D}\rightarrow n-\widehat{\Dl}^\comp\mc{D}\ .$$
Then, the following analogue of \ref{relacionTotyTotGorro} also
holds
$${Tot}_n = \mc{I}\comp \widehat{Tot}_n \comp \mc{P}  \ .$$
\end{obs}

\begin{obs}
Since $+:2-\widehat{\Dl}_+\rightarrow \widehat{\Dl}$ is
associative, we can consider the sum $E_0 +\cdots +E_{n-1}$ with
no need of fixing the order in which the sums are made.
Consequently, the following property holds.\\
Consider a non empty ordered set $E$. It is clear that there is a
bijection between the decompositions of $E$ as an ordered sum of
$n$ of its subsets, $E=E_0+\cdots +E_{n-1}$, and the
decompositions $E=E'_0 + E'_1$ together with two more
decompositions $E'_0=E_0+\cdots +E_{k-1}$ and
$E'_1=E_{k}+\cdots +E_{n-1}$, for a f{i}xed $0\leq k \leq n-1$.\\
Iterating this procedure, it follows that $\widehat{Tot}_n$ (resp.
$Tot_n$) agrees with consecutive iterations of
$\widehat{Tot}=\widehat{Tot}_2$ (resp. $Tot=Tot_2$).
Let us see the case $n=3$ in more detail.
\end{obs}

\begin{num}
Consider $T\in 3-\widehat{\Dl}_+^\comp\mc{D}$ and an ordered set
$E$.\\
Roughly speaking, $\coprod_{F=G+H} T(E,G,H)$ is the total object
of $T(E,-,-)$ evaluated at $F$.\\
However,{ when }$E\neq\emptyset$ it turns out that $T(E,-,-)\in
(\widehat{\Dl}_+\times \widehat{\Dl}_+)^\comp \mc{D}$, whereas
$T(\emptyset,-,-)\in 2-\widehat{\Dl}_+^\comp\mc{D}$. Therefore, we
introduce the following notations.

\noindent Given $E\in\widehat{\Dl}_+$, we define
$$\widehat{Tot}^\ast (T(E,-,-))=\begin{cases}
                                     \widehat{Tot}^+(T(E,-,-))       & \mbox{{ if }} E\neq\emptyset\\
                                     \widehat{Tot}(T(\emptyset,-,-)) & \mbox{{ if }} E=\emptyset \ ,
                                     \end{cases}$$
where $\widehat{Tot}^+:(\widehat{\Dl}_+\times
\widehat{\Dl}_+)^\comp \mc{D}\rightarrow \widehat{\Dl}_+^\comp
\mc{D}$ is def{i}ned as $Tot^+$ in \ref{def{i}TotGorroAumentado}.\\
We define also
$$\widehat{Tot}^\ast (T(-,-,E))=\begin{cases}
                        \widehat{Tot}^+(T(-,-,E))       & \mbox{{ if }} E\neq\emptyset\\
                        \widehat{Tot}(T(-,-,\emptyset)) & \mbox{{ if }} E=\emptyset \ .
                        \end{cases}$$
\end{num}

\begin{lema}
There exists functors
$\widehat{Tot}_{(0)},\widehat{Tot}_{(2)}:3-\widehat{\Dl}_+^\comp\mc{D}\rightarrow
2-\widehat{\Dl}_+^\comp\mc{D}$, such that
$$\begin{array}{l}
 {[\widehat{Tot}_{(0)}(T)]}(E,F)=[\widehat{Tot}^\ast (T(E,-,-))](F)\ \mbox{ and}\\
 {[\widehat{Tot}_{(2)}(T)]}(E,F)=[\widehat{Tot}^\ast (T(-,-,F))](E)
 \end{array}$$
for any $T$ in $3-\widehat{\Dl}_+^\comp\mc{D}$.
\end{lema}

\begin{proof}
Firstly, it suf{f}{i}ces to prove the statement for
$Tot_{(0)}(T)$. In this case,{ if }$R\in
3-\widehat{\Dl}_+^\comp\mc{D}$ is given by $R(E,F,G)=T(G,F,E)$, we
have that $Z=\widehat{Tot}_{(0)}(R)$ is a biaugmented bisimplicial
object. Note that $\widehat{Tot}_{(2)}(T)$ is just the object
in $2-\widehat{\Dl}_+^\comp\mc{D}$ obtained by interchanging the indexes $E$ and $F$ in $Z$.\\
Let us check that $\widehat{Tot}_{(0)}(T)$ is in fact in
$2-\widehat{\Dl}_+^\comp\mc{D}$. Given $(f,g):(E',F')\rightarrow
(E,F)$, then
$$[\widehat{Tot}_{(0)}(T)](f,g):[\widehat{Tot}_{(0)}(T)](E,F)\rightarrow
[\widehat{Tot}_{(0)}(T)](E',F')$$ is induced by $T$ in a natural way.\\
To see this, we have that
$$ [\widehat{Tot}_{(0)}(T)](E,F)=\!\!\!\ds\coprod_{F=F_0+F_1}\!\!\! T(E,F_0,F_1)\ ; \ [\widehat{Tot}_{(0)}(T)](E',F')=\!\!\!\ds\coprod_{F'=F'_0+F'_1}\!\!\! T(E',F'_0,F'_1)\ . $$
Given $F_0$ and $F_1$ with $F=F_0+F_1$, set $F'_i=g^{-1}(F_i)$,
$g_i=g|_{F'_i}:F'_i\rightarrow F_i$ for $i=0,1$.\\
Then
$$[\widehat{Tot}_{(0)}(T)](f,g)|_{T(E,F_0,F_1)}=T(f,g_0,g_1):T(E,F_0,F_1)\rightarrow
T(E',F'_0,F'_1)\ ,$$ and clearly $\widehat{Tot}_{(0)}(T)$ is
functorial in $(f,g)$.
\end{proof}

\begin{prop}\label{TotGorroIterado}
The functors $\widehat{Tot}_3$ ,
$\widehat{Tot}\comp\widehat{Tot}_{(0)}$ and
$\widehat{Tot}\comp\widehat{Tot}_{(2)}$ are isomorphic.\\
In other words, given $T\in 3-\widehat{\Dl}^\comp\mc{D}$, there
are canonical and functorial isomorphisms
$$ \widehat{Tot}_3(T)\simeq\widehat{Tot}(\widehat{Tot}^\ast_{(0)}(T))\ \ \widehat{Tot}_3(T)\simeq\widehat{Tot}(\widehat{Tot}^\ast_{(2)}(T))\ .$$
\end{prop}

\begin{proof}
By definition
$$ [\widehat{Tot}_{3}(T)](E)= \ds\coprod_{E= E_0 +E_1 +E_{2}} T(E_0,E_1,E_{2})\ .$$
Clearly, the sets
$$\begin{array}{l}
 \{(E_0,E_1,E_2)\in 3-\widehat{\Dl}\ |\ E= E_0 +E_1 +E_{2}\}\nonumber \\
 \{(F,G)\times(G_0,G_1)\in (2-\widehat{\Dl})\times(2-\widehat{\Dl}) \ |\ E= F +G\mbox{ and }G=G_0+G_1 \} \nonumber \\
 \{(E,F)\times(E_0,E_1)\in (2-\widehat{\Dl})\times(2-\widehat{\Dl}) \ |\ E= F +G\mbox{ and }F=F_0+F_1 \} \nonumber
 \end{array}$$
are bijective.  Hence, after reordering the terms in the coproduct
def{i}ning $[\widehat{Tot}_{3}(Z)](E)$ we obtain canonical
isomorphisms
\begin{eqnarray}
 {[\widehat{Tot}_{3}(T)]}(E)\!\stackrel{\sigma_E}{\simeq}\!\!\!\!\! \ds\coprod_{E= F +G}\!\! \left(\ds\coprod_{G=G_0+G_1}\!\!\!\!\! T(F,G_0,G_1)\!\!\right)\! =\!\!\!\! \ds\coprod_{E= F +G}\!\!\!\! \widehat{Tot}^\ast_{(0)}(T)(F,G)\! =\! [\widehat{Tot}(\widehat{Tot}^\ast_{(0)}(T))](E)\nonumber \\
 {[\widehat{Tot}_{3}(T)]}(E)\!\stackrel{\rho_E}{\simeq}\!\!\!\!\! \ds\coprod_{E= F +G}\!\! \left(\ds\coprod_{F=F_0+F_1}\!\!\!\!\! T(F_0,F_1,G)\!\!\right)\! =\!\!\!\! \ds\coprod_{E= F +G}\!\!\!\! \widehat{Tot}^\ast_{(2)}(T)(F,G)\! =\! [\widehat{Tot}(\widehat{Tot}^\ast_{(2)}(T))](E)\nonumber
\end{eqnarray}
Therefore, for every $f:E'\rightarrow E$ it holds that
$[\widehat{Tot}(\widehat{Tot}^\ast_{(0)}(T))](f)\comp \sigma_E =
\sigma_{E'}\comp [\widehat{Tot}_{3}(T)](f)$ (and similarly for
$\rho$).

\noindent Consider $E_i$, $i=0,1,2$ with $E=E_0+E_1+E_2$.\\
Then
$[\widehat{Tot}_{3}(T)](f)|_{T(E_0,E_1,E_2)}=T(f_0,f_1,f_2):T(E_0,E_1,E_2)\rightarrow
T(E'_0,E'_1,E'_2)$, where
$E'_i=f^{-1}(E_i)$ and $f_i=f|_{E'_i}$ for $i=0,1,2$.\\
In addition, the restriction of $\sigma_{E'}$ to
$T(E'_0,E'_1,E'_2)$ is the identity on
 the same component of
$\widehat{Tot}^\ast_{(0)}(E'_0,E'_1+E'_2)$.

On the other hand, the restriction of $\sigma_E$ to
$T(E_0,E_1,E_2)$ is the identity on the same component of
$\widehat{Tot}^\ast_{(0)}(E_0,E_1+E_2)$, whereas the restriction
of $[\widehat{Tot}(\widehat{Tot}^\ast_{(0)}(T))](f)$ to
$\widehat{Tot}^\ast_{(0)}(E_0,E_1+E_2)$ is
$\widehat{Tot}^\ast_{(0)}(f_0,f_1+f_2)$, since
$f|_{f^{-1}(E_0+E_1)}=f_0+f_1$.\\
F{i}nally,
$\widehat{Tot}^\ast_{(0)}(f_0,f_1+f_2)|_{T(E_0,E_1,E_2)}$
coincides with $T(f_0,f_1,f_2)$ by def{i}nition, and the proof is
concluded.
\end{proof}

The above proposition also holds for the functor $Tot$.

\begin{num}
Consider $T\in3-\simp\mc{D}$ and $n\geq -1$. We define
$${Tot}^\ast (T([n],-,-))=\begin{cases}
                                     {Tot}^+(T([n],-,-))       & \mbox{{ if }} n > -1\\
                                     {Tot}(T([-1],-,-)) & \mbox{{ if }} n=-1 \ .
                                     \end{cases}$$
Similarly,
$${Tot}^\ast (T(-,-,[n]))=\begin{cases}
                        {Tot}^+(T(-,-,[n]))       & \mbox{{ if }} \mbox{ if } n > -1\\
                        {Tot}(T(-,-,[-1])) & \mbox{{ if }} n=-1 \ .
                        \end{cases}$$
\end{num}

\begin{cor}\label{TotIterado}
There exists functors ${Tot}_{(0)}$,
${Tot}_{(2)}:3-\simp\mc{D}\rightarrow 2-\simp\mc{D}$ such that
$$\begin{array}{l}
 {[{Tot}_{(0)}(T)]}([n],[m])=[\widehat{Tot}^\ast (T([n],-,-))]([m])\ \mbox{ and }\\
 {[{Tot}_{(2)}(T)]}([n],[m])=[\widehat{Tot}^\ast (T(-,-,[m]))]([n])\ .
\end{array}$$
Moreover, the functors $Tot_3$, $Tot\comp{Tot}_{(0)}$ and
$Tot\comp {Tot}_{(2)}:3-\simp\mc{D}\rightarrow \simp\mc{D}$ are
isomorphic.
\end{cor}

\begin{proof}
We have that ${Tot}_{(0)}=
\mc{I}\comp\widehat{Tot}_{(0)}\comp\mc{P}$, and the same holds for
${Tot}_{(2)}$. More specif{i}cally, given $T\in 3-\simp\mc{D}$
then
$$[\mc{I}\comp\widehat{Tot}_{(0)}\comp\mc{P}(T)]([n],[m])=[\widehat{Tot}_{(0)}\comp\mc{P}(T)]([n],[m])=[\widehat{Tot}^\ast ((\mc{P}T)([n],-,-))]([m])\ ,$$
that agrees with $[{Tot}^\ast (T([n],-,-))]([m])$ since $Tot$ and
$Tot^+$ are obtained from $\widehat{Tot}$ and
$\widehat{Tot}^+$ by composition with $\mc{I}$ and $\mc{P}$.\\
Consequently ${Tot}_{(0)}$ and $Tot_{(1)}$ are functors, and the
statement follows from \ref{TotGorroIterado} together with
$\mc{P}\mc{I}\simeq Id$.
\end{proof}


\section{Simplicial cylinder object}\label{defCilindroSimp}

In this section we introduce the \emph{simplicial cylinder
object}, that is a generalization of the simplicial cone object
$Cf$ associated with a morphism $f:X\rightarrow Y$ between
simplicial objects in $\mc{D}$, \ref{conoSimpl}.

\begin{def{i}}\label{Def{i}Omega}
Let $\Omega(\mc{D})$\index{Symbols}{$\Omega$} be the category of
pairs $(f,\epsilon)$ consisting of diagrams in $\simp\mc{D}$
$$\xymatrix@M=4pt@H=4pt{ X \ar[r]^{f} \ar[d]_{\epsilon} &  Y \\
                         X_{-1}\times\Dl . &}$$
A morphism in $\Omega(\mc{D})$ is a triple
$(\alpha,\beta,\gamma):(f,\epsilon)\rightarrow(f',\epsilon')$ such that $$\xymatrix@M=4pt@H=4pt{ X_{-1}\times\Dl \ar[d]_-{\alpha} & X \ar[r]^{f} \ar[l]_-{\;\epsilon}\ar[d]_{\beta} &  Y \ar[d]^{\gamma}\\
                                                                                                     X'_{-1}\times\Dl & X' \ar[r]^{f'} \ar[l]_-{\;\epsilon'} &  Y'   }$$
commutes.\\
In a similar way, we will denote by
$\mathrm{Co}\Omega(\mc{D})$\index{Symbols}{$\mathrm{Co}\Omega$}
the category whose objects are the diagrams in $\simp\mc{D}$
$$\xymatrix@M=4pt@H=4pt{ X  &  Y \ar[l]_{f}\\
                         X_{-1}\times\Dl \ar[u]^{\varepsilon} . &}$$
\end{def{i}}

\begin{num}\label{def{i}PsiCilindroSimp}
Let $\psi:\Omega(\mc{D})\rightarrow 2-\Dl_{+}^\comp\mc{D}$ be the
functor that maps the pair $(f,\epsilon)$ into the biaugmented
bisimplicial (\ref{diagramaCilindro}) of
\ref{ejObjBisimpfEpsilon}.
\end{num}

\begin{def{i}}[Simplicial cylinder object]\mbox{}\\
The simplicial cylinder functor
$Cyl:\Omega(\mc{D})\rightarrow\simp\mc{D}$\index{Symbols}{$Cyl$}\index{Index}{simplicial!cylinder}
is the composition
$$\xymatrix@H=4pt@M=4pt{\Omega(\mc{D})\ar[r]^-{\psi} & 2-\Dl_{+}^\comp\mc{D}\ar[r]^-{Tot}& \simp\mc{D}}\ .$$
In other words, $Cyl(f,\epsilon)=Tot(\psi(f,\epsilon))$.
\end{def{i}}

\begin{num}\label{def{i}CylenTerminosdeSigma} Having in mind the description \ref{def{i}Tot} of $Tot$, the functor $Cyl$
can be described as follows.\\
Denote by $u_i:[n]\rightarrow [1]$ the morphisms in $\Dl$ with
$u_i([n])=i$,{ if }$i=0,1$ and
$$\Lambda_n=\{\sigma:[n]\rightarrow [1]\ |\ \sigma\neq u_0 , \ u_1 \}\ .$$
Given $\sigma:[n]\rightarrow [1]$, we will identify the ordered
set $\sigma^{-1}(i)$ with the corresponding object in $\Dl_+$ for
$i=0,1$.
Then
$$[Cyl(f,\epsilon)]([n])= \ds\coprod_{\sigma:[n]\rightarrow[1]} Cyl(f,\epsilon)^{(\sigma)}\ ,$$
where
$$[Cyl(f,\epsilon)]^{(\sigma)}=\begin{cases}
                                             Y([n]) & \mbox{{ if }} \sigma=u_1\\
                                             X_{-1} & \mbox{{ if }} \sigma=u_0\\
                                             X(\sigma^{-1}(1)) & \mbox{{ if }}\sigma\in\Lambda_n\ .
                                           \end{cases}$$
If $\theta:[m]\rightarrow [n]$, the restriction of
$\Theta=[Cyl(f,\epsilon)](\theta)$ to the component indexed by
$(\sigma)$ is
\begin{equation}\label{def{i}CylAbstracto}\Theta|_{(\sigma)}=
\begin{cases}
 X(\theta|_{(\sigma\theta)^{-1}(1)}):X(\sigma^{-1}(1))\rightarrow X((\sigma\theta)^{-1}(1))& \mbox{if } \sigma\theta\in\Lambda_m\\
 Y(\theta):Y([n])\rightarrow Y([m])& \mbox{if } \sigma=u_1\\
 f([m])\comp X(\theta|_{(\sigma\theta)^{-1}(1)}):X(\sigma^{-1}(1))\rightarrow Y([m])& \mbox{if } \sigma\in\Lambda_n\mbox{ and }\sigma\theta=u_1\\
 Id:X_{-1}\rightarrow X_{-1}& \mbox{if } \sigma=u_0\\
 \epsilon(\sigma^{-1}(1)):X(\sigma^{-1}(1))\rightarrow X_{-1}& \mbox{if }\sigma\in\Lambda_n\mbox{ and }\sigma\theta=u_0\ .
\end{cases}\end{equation}
\end{num}

\begin{num}
More specif{i}cally, $Cyl(f,\epsilon)$ is in degree $n$
$$ Cyl(f,\epsilon)_n = Y_n \sqcup X_{n-1} \sqcup X_{n-2} \sqcup \cdots \sqcup X_0 \sqcup X_{-1} \ .$$
The face maps $d_i : Cyl(f,\epsilon)_{n} \rightarrow
Cyl(f,\epsilon)_{n-1}$ are given componentwise by $d_i | _{Y_{n}}
= d^{Y}_{i}$, and{ if }$1\leq k \leq n+1$ then
$$
 d_i |_{X_{n-k}} = \left\{ \begin{array}{ll} d_{i-k}^{X} & \mbox{{ if }} i \geq k \\
                                                  Id_{X_{n-k}} & \mbox{{ if }} i < k \; \mbox{ and } (k,i)\neq (1,0)\\
                                                  f_{n-1} & \mbox{{ if }} (k,i)=(1,0)\end{array}\right.$$
where $d^X_0=\epsilon_0: X_0\rightarrow X_{-1}$. Visually,{ if
}$1\leq i \leq n$, then $d_i$ is
$$\xymatrix@1@C=3pt{ Y_{n}\ar[rrd]_(.35){d_i^Y}  & \sqcup &  X_{n-1} \ar[rrd]_(.35){d_{i-1}^X} & \sqcup & X_{n-2} \ar[rrd]_(.35){d_{i-2}^X}& \sqcup & \cdots  &         & \sqcup  & X_{n-i}\ar[rrd]_(.35){d_{0}^X} & \sqcup & X_{n-i-1}\ar[d]^(.5){Id} & \sqcup & \cdots  & \sqcup &  X_0 \ar[d]^(.5){Id}   & \sqcup & X_{-1} \ar[d]^(.5){Id }& \\
                                                          &        &  Y_{n-1}                                    & \sqcup & X_{n-2}                                 & \sqcup & X_{n-3} &  \sqcup &         & \cdots                                  & \sqcup & X_{n-i-1}                          & \sqcup & \cdots  & \sqcup &  X_0                             & \sqcup & X_{-1}                & ,}$$
whereas $d_0$ is
$$\xymatrix@1@C=4pt{ Y_{n}\ar[rrd]_(.35){d_{0}^{Y}}  & \sqcup &  X_{n-1} \ar[d]^(.5){f_{n-1}} & \sqcup & X_{n-2} \ar[d]^(.5){Id}& \sqcup & X_{n-3} \ar[d]^(.5){Id} &  \sqcup & \cdots  & \sqcup &  X_0 \ar[d]^(.5){Id}   & \sqcup & X_{-1} \ar[d]^(.5){Id}& \\
                                                              &         &  Y_{n-1}                    & \sqcup & X_{n-2}                           & \sqcup & X_{n-3}                            &  \sqcup & \cdots  & \sqcup &  X_0           & \sqcup & X_{-1} & .}$$

The degeneracy maps $s_j: Cyl(f,\epsilon)_n \rightarrow
Cyl(f,\epsilon)_{n+1}$ are def{i}ned as $s_j {|}_{Y_n}= s_j^{Y}$
and given $1 \leq k \leq n+1$ then
$$ s_j |_{X_{n-k}} = \left\{ \begin{array}{ll} s_{j-k}^{X} & \mbox{{ if }} j\geq k \\
                                           Id_{X_{n-k}} & \mbox{{ if }} j < k \, , \end{array}\right.$$
that is to say
$$\xymatrix@1@C=4pt{        &       & Y_{n}\ar[lld]_(.65){s_j^Y}  & \sqcup & X_{n-1} \ar[lld]_(.65){s_{j-1}^X} & \sqcup & X_{n-2} \ar[lld]_(.65){s_{j-2}^X}& \sqcup & \cdots    &  \sqcup & X_{n-j}\ar[lld]_(.65){s_{0}^X} & \sqcup & X_{n-j-1}\ar[d]^(.5){Id} & \sqcup & \cdots  & \sqcup  & X_{-1} \ar[d]^(.5){Id }& \\
                    Y_{n+1} &\sqcup & X_n                         & \sqcup & X_{n-1}                           & \sqcup & \cdots                           &        & X_{n-j+1} &  \sqcup & X_{n-j}                        & \sqcup & X_{n-j-1}                            & \sqcup & \cdots  & \sqcup  & X_{-1}                & .}$$
\end{num}

\begin{num}\label{inclusionesCyl}
Given $\xymatrix@M=4pt@H=4pt@C=15pt{D: X_{-1}\times\Dl & X
\ar[r]^-{f}\ar[l]_-{\epsilon} &  Y }$ in $\Omega(\mc{D})$, it
follows from  \ref{inclusionesTot} that the canonical morphisms
$$i_{Y_n}:Y_n\rightarrow Cyl(D)_n\mbox{ and }i_{X_{-1}}X_{-1}\rightarrow Cyl(D)_n $$
induce the following diagram, natural in $(f,\epsilon)$,
\begin{equation}\label{CuadradoCyl}\xymatrix@M=4pt@H=4pt{ X \ar[r]^-{f}   \ar[d]_{\epsilon}           & Y \ar[d]^{i_Y}\\
                       X_{-1}\times\Dl \ar[r]^-{i_{X_{-1}}}&  Cyl(D) \ .}
                       \end{equation}
Note that{ if }$0$ is an initial object in $\mc{D}$, the morphisms
$i_{X_{-1}}$ and $i_Y$ are just the image under the functor $Cyl$
of the maps
$$\xymatrix@M=4pt@H=4pt{ X_{-1}\times\Dl \ar[d]_-{Id} & 0 \ar[r] \ar[l] \ar[d] & 0 \ar[d] &\ar@{}[d]|{;}&  0 \ar[d] & 0 \ar[r] \ar[l] \ar[d]               & Y \ar[d]_-{Id}  \\
                         X_{-1}\times\Dl & X \ar[r]^{f} \ar[l]_-{\;\epsilon} &  Y            &             &  X_{-1}\times\Dl      & X \ar[r]^{f} \ar[l]_-{\;\epsilon} & Y,   }$$
since $Cyl(X_{-1}\times\Dl\leftarrow 0\rightarrow 0)=
X_{-1}\times\Dl$ and $Cyl(0\leftarrow 0\rightarrow Y)= Y$.\\
\end{num}
\begin{obs}\label{def{i}ConoSimplicial}
Let $Fl(\simp\mc{D})$\index{Symbols}{$Fl(\simp\mc{D})$} be the category of morphisms in $\simp\mc{D}$.\\
If $\mc{D}$ has a f{i}nal object $1$, consider the inclusion
$\mc{I}:Fl(\simp\mc{D})\rightarrow\Omega(\mc{D})$ that maps the
diagram $1\longleftarrow
X\stackrel{f}{\longrightarrow}Y$ into the morphism $f:X\rightarrow Y$.\\
Hence $\mc{I}$ is right adjoint to the forgetful functor
$\mrm{U}:\Omega(\mc{D})\rightarrow Fl(\simp\mc{D})$,
$\mrm{U}(X,f,\epsilon)=f$, and the  \textit{simplicial cone}
functor $C:Fl(\simp\mc{D})\rightarrow\simp\mc{D}$ is
$Cyl\comp\mc{I}$.\\
\indent Given $X\in\simp\mc{D}$, $CX$\index{Symbols}{$CX$} will
mean $C(Id_X)$, and{ if }$S$ is an object in $\mc{D}$,
$Cyl(S)$\index{Symbols}{$Cyl(S)$} will denote
$Cyl(S\times\Dl,Id_{S\times\Dl},Id_{S\times\Dl})$.\\
\end{obs}
\noindent Next we study some of the properties of the functor
$Cyl$.

\begin{prop}\label{CylCoproducto}
The functor $Cyl: \Omega(\mc{D})\rightarrow\simp\mc{D}$ commutes
with coproducts, that is
$$Cyl(f,\epsilon)\sqcup Cyl(f',\epsilon')\simeq Cyl(f\sqcup f',\epsilon\sqcup\epsilon') \ .$$
\end{prop}

\begin{proof}
The statement follows directly from the commutativity of the
functors $\Psi$ and $Tot$ with coproducts (see
\ref{relacionTotyTotGorro}).
\end{proof}

The proof of the following proposition will be given later in
\ref{HomotCuadradoCyl}.

\begin{prop}\label{CylHomotopia}
For every $\xymatrix@M=4pt@H=4pt@C=15pt{D: X_{-1}\times\Dl & X
\ar[r]^-{f}\ar[l]_-{\epsilon} &  Y }$ in $\Omega(\mc{D})$, the
diagram $($\ref{CuadradoCyl}$)$ commutes up to simplicial
homotopy, natural in $D$.
\end{prop}

\begin{prop}\label{Cylaumentado}
Given a commutative diagram in $\simp\mc{D}$
\begin{equation}\label{CuadradConmSimpl}\xymatrix@M=4pt@H=4pt{ X \ar[r]^-{f}   \ar[d]_{\epsilon}           & Y \ar[d]^{\rho'}\\
                       X_{-1}\times\Dl \ar[r]^-{\rho}&  T\times\Dl \ ,}
                       \end{equation}
there exists a unique $H:Cyl(f,\epsilon)\rightarrow T\times\Dl$
such that $H\comp i_{X_{-1}}=\rho$ and $H\comp i_Y=\rho'$.
Equivalently,
$$\xymatrix@M=4pt@H=4pt{ X \ar[r]^f   \ar[d]_{\epsilon}          & Y \ar@/^1pc/[rdd]^{\rho'} \ar[d]^{i_Y}                 &\\
                        X_{-1}\times\Dl \ar@/_1pc/[drr]_{\rho}  \ar[r]^{i_{X_{-1}}}  &  Cyl(f,\epsilon)\ar@{}[d]|{\sharp}  \ar@{}[r]|{\sharp}\ar[dr]^{H}   &\\
                                                                 &                                                 & T\times\Dl \ .} $$
In addition, $H$ is natural in $($\ref{CuadradConmSimpl}$)$.
\end{prop}

\begin{proof}
The data $T$, $\rho$ and $\rho'$ allows to construct ${}_+Z\in
(\Dl_+\times\Dl_+)^\comp\mc{D}$ from $\Psi(f,\epsilon)$. The
restriction of ${}_+Z$ to $2-\Dl_+^\comp\mc{D}$ is $\Psi(f,\epsilon)$ and ${}_+Z_{-1,-1}=T$.\\
Hence, the morphism $H_0 :Y_0\sqcup X_{-1}\rightarrow T$ with
$H_0\mid_{Y_{0}}=\rho'_0$,
$H_0\mid_{X_{-1}}=\rho$ is an augmentation of $Cyl(f,\epsilon)$ by \ref{def{i}TotAumentado}.\\
Moreover $H_0$ is the unique morphism such that $H_0\comp
i_{X_{-1}}=\rho_0$ and $H_0\comp i_{Y_0}=\rho'_0$.
We deduce from \ref{aumentvisionequivalente} that $H_0$ induces
$H:Cyl(f,\epsilon)\rightarrow T\times\Dl$ with $H_n=H_0 \comp
(d_0)^n $, and such that $H\comp
i_{X_{-1}}=\rho:X_{-1}\times\Dl\rightarrow T\times\Dl$ and $H\comp
i_Y=\rho':Y\rightarrow T\times\Dl$, because these morphisms agree
in degree 0.\\
Since $H$ is determined by $H_0$, it follows that $H$ is the
unique morphism satisfying the required equality, and $H$ is
natural in $($\ref{CuadradConmSimpl}$)$ because $H_0$ is so.
\end{proof}

Now we develop a property of ${Cyl}$ that will be needed later for
the study of the relationship between simplicial descent
categories and triangulated categories.

\begin{num}
We have that $\simp\mc{D}$ has coproducts because $\mc{D}$ has.
Then we can consider the cylinder functor of a morphism
$f:X\rightarrow Y$ between bisimplicial objects, where $X$ has an
augmentation $\epsilon$. This can be an augmentation with respect
to any of the two simplicial indexes of $X$. Therefore we need to
introduce the following notations.
\end{num}

\begin{def{i}}\label{cylBisimpl}
Consider the category
$\Omega^{(1)}(\simp\mc{D})$\index{Symbols}{$\Omega^{(1)}$} whose
objects are the diagrams
$$\xymatrix@M=4pt@H=4pt{
{\Dl\times Z_{-1}} & Z\ar[l]_-{\epsilon} \ar[r]^f & T }\ ,$$
that in degree $n,m$ is $\xymatrix@M=4pt@H=4pt{Z_{-1,m} &
Z_{n,m}\ar[l]_{\epsilon_{n,m}} \ar[r]^{f_{n,m}} & T_{n,m} }$.\\
Hence, the functor
$Cyl^{(1)}_{\simp\mc{D}}:\Omega^{(1)}(\simp\mc{D})\rightarrow
\simp\simp\mc{D}$\index{Symbols}{$Cyl^{(1)}_{\simp\mc{D}}$} is
$$Cyl^{(1)}_{\simp\mc{D}}(f,\epsilon)_{n,m}=T_{n,m} \sqcup Z_{n-1,m} \sqcup  {\cdots} \sqcup Z_{0,m} \sqcup Z_{-1,m} \ .$$
If $\alpha:[m']\rightarrow [m]$ then
$[Cyl^{(1)}_{\simp\mc{D}}(f,\epsilon)](Id,\alpha):Cyl^{(1)}_{\simp\mc{D}}(f,\epsilon)_{n,m}\rightarrow
Cyl^{(1)}_{\simp\mc{D}}(f,\epsilon)_{n,m'}$ is
$${T}(Id,\alpha)\sqcup Z(Id,\alpha) \sqcup \stackrel{n)}{\cdots} \sqcup Z(Id,\alpha)\sqcup {Z_{-1}}(\alpha)\ ,$$
whereas{ if }$\beta:[n']\rightarrow [n]$, we def{i}ne
$[Cyl^{(1)}_{\simp\mc{D}}(f,\epsilon)](\beta,Id):Cyl^{(1)}_{\simp\mc{D}}(f,\epsilon)_{n,m}\rightarrow
Cyl^{(1)}_{\simp\mc{D}}(f,\epsilon)_{n',m}$ using the formulae
(\ref{def{i}CylAbstracto}), by forgetting the index $m$.

The category
$\Omega^{(2)}(\simp\mc{D})$\index{Symbols}{$\Omega^{(2)}$} is
def{i}ned in the same way, but considering the diagrams
$$\xymatrix@M=4pt@H=4pt{D: Z_{-1}\times\Dl & Z\ar[l]_-{\epsilon} \ar[r]^f & T }\ .$$
Similarly, we def{i}ne the diagram
$Cyl^{(2)}_{\simp\mc{D}}:\Omega^{(2)}(\simp\mc{D})\rightarrow
\simp\simp\mc{D}$\index{Symbols}{$Cyl^{(2)}_{\simp\mc{D}}$} by
applying $Cyl$ to the second index of the diagram $D$. Then, we
obtain the following square of functors
\begin{equation}\label{diagGammaOmega}\xymatrix@H=4pt@M=4pt@C=37pt{ \Omega^{(1)}(\simp\mc{D})\ar[d]_{\Gamma} \ar[r]^{Cyl_{\simp\mc{D}}^{(1)}} & \simp\simp\mc{D}\ar[d]_{\Gamma}\\
                         \Omega^{(2)}(\simp\mc{D}) \ar[r]^{Cyl_{\simp\mc{D}}^{(2)}} & \simp\simp\mc{D} }\ .\end{equation}
\end{def{i}}

\begin{num}
As happens with $Cyl$, we have the canonical inclusions of
$\Dl\times Z_{-1}$ (resp. $Z_{-1}\times\Dl$) and $T$ in
$Cyl_{\simp\mc{D}}^{(1)}$ (resp. $Cyl_{\simp\mc{D}}^{(2)}$).
\end{num}

\begin{num}\label{NotacionesFactCyl}
Assume that the following diagram commutes in $\mc{D}$
\begin{equation}\label{diagrTresporTres}\xymatrix@M=4pt@H=4pt{Z'  & X'\ar[l]_{g'}\ar[r]^{f'}   & Y' \\
                                                              Z  \ar[u]^{\alpha} \ar[d]_{\alpha'}& X  \ar[l]_{g}\ar[r]^{f} \ar[u]^{\beta} \ar[d]_{\beta'} & Y \ar[u]^{\gamma} \ar[d]_{\gamma'} \\
                                                              Z'' & X''\ar[l]_{g''}\ar[r]^{f''} & Y''.}\end{equation}
Consider $(\ref{diagrTresporTres})$ in $\simp\mc{D}$ through the
functor $-\times\Dl$. Applying ${Cyl}$ by rows and columns we
obtain
$$\xymatrix@M=4pt@H=4pt@R=8pt{{Cyl}(f',g')& {Cyl}(f,g)\ar[l]_{\rho}\ar[r]^{\rho'}& {Cyl}(f'',g'')\\
                              {Cyl}(\alpha',\alpha)& {Cyl}(\beta',\beta)\ar[l]_{G}\ar[r]^{F}& {Cyl}(\gamma',\gamma).}$$
Then the diagrams of $\simp\mc{D}$
$$\xymatrix@M=4pt@H=4pt@C=17pt{Y'\ar[d]_{i}& Y\ar[l]_{\gamma}\ar[r]^{\gamma'}\ar[d]_{i} & Y''\ar[d]_{i}\ar@{}[rd]|{;} & Z''\ar[d]_{i}& X''\ar[l]_{g''}\ar[r]^{f''}\ar[d]_{i} & Y''\ar[d]_{i}\\
                         Cyl(f',g')& Cyl(f,g)\ar[l]_{\rho}\ar[r]^{\rho'}& Cyl(f'',g'')                                          & Cyl(\alpha',\alpha)& Cyl(\beta',\beta)\ar[l]_{G}\ar[r]^{F}& Cyl(\gamma',\gamma).}$$
give rise to the morphisms between bisimplicial objects
$$\begin{array}{l}\varphi:\Dl\times Cyl(\gamma',\gamma)\rightarrow Cyl_{\simp\mc{D}}^{(2)}(\rho'\times\Dl,\rho\times\Dl)\\
                  \phi:Cyl(f'',g'')\times\Dl\rightarrow Cyl_{\simp\mc{D}}^{(1)}(\Dl\times F,\Dl\times G)\ .\end{array}$$
\end{num}

\begin{lema}\label{factorizCylTilde}
Under the above notations, there exists a canonical morphism\\
$\Theta:{Cyl}_{\simp\mc{D}}^{(1)}(\Dl\times F,\Dl\times
G)\rightarrow
{Cyl}_{\simp\mc{D}}^{(2)}(\rho'\times\Dl,\rho\times\Dl)$ in
$\simp\simp\mc{D}$, such that the following diagram commutes
$$\xymatrix@R=10pt@C=25pt{                                                     &{Cyl}_{\simp\mc{D}}^{(1)}(\Dl\times F,\Dl\times G) \ar[dd]_{\Theta}&                                          \\
                               Cyl(f'',g'')\times\Dl\ar[dr]_{i}\ar[ur]^{\phi}  &                                                                  &\Dl\times Cyl(\gamma',\gamma) \ar[lu]_{i}\ar[ld]^{\varphi} \ .\\
                                                                               &{Cyl}_{\simp\mc{D}}^{(2)}(\rho'\times\Dl,\rho\times\Dl)           &                               }$$
\end{lema}

\begin{proof} Set
$A_{\cdot\cdot}=A\times\Dl\times\Dl\in\simp\simp\mc{D}${ if
}$A\in\mc{D}$, as well as $h_{\cdot\cdot}=h\times\Dl\times\Dl${ if
}$h$ is a morphism in $\mc{D}$.
The diagram of $\simp\simp\mc{D}$
\begin{equation}\label{diagrConmCuatroporCuatro}\xymatrix@M=4pt@H=4pt{   Z'_{\cdot\cdot} & X'_{\cdot\cdot}\ar[l]_{g'_{\cdot\cdot}}\ar[r]^{f'_{\cdot\cdot}}   & Y'_{\cdot\cdot} \ar[r]^-{i} & Cyl(f',g')\times\Dl\\
                                                                         Z_{\cdot\cdot} \ar[u]^{\alpha_{\cdot\cdot}} \ar[d]_{\alpha'_{\cdot\cdot}}& X_{\cdot\cdot}  \ar[l]_{g_{\cdot\cdot}}\ar[r]^{f_{\cdot\cdot}} \ar[u]^{\beta_{\cdot\cdot}} \ar[d]_{\beta'_{\cdot\cdot}} & Y_{\cdot\cdot} \ar[u]^{\gamma_{\cdot\cdot}} \ar[d]_{\gamma'_{\cdot\cdot}}\ar[r]^-{i}& Cyl(f,g)\times\Dl \ar[u]^{{\rho}\times\Dl}\ar[d]_{{\rho'}\times\Dl} \\
                                                                         Z''_{\cdot\cdot} \ar[d]_{i} & X''_{\cdot\cdot}\ar[d]_{i}\ar[l]_{g''_{\cdot\cdot}}\ar[r]^{f''_{\cdot\cdot}} & Y''_{\cdot\cdot}\ar[d]_{i}\ar[r]^-{i}& Cyl(f'',g'')\times\Dl\ar[d]_{\phi}\\
                                                                         \Dl\!\times\! Cyl(\alpha',\alpha)& \Dl\!\times\!  Cyl(\beta',\beta)\ar[l]_{\Dl\times G}\ar[r]^{\Dl\times F}& \Dl\!\times\!   Cyl(\gamma',\gamma)\ar[r]^-{i}& Cyl_{\simp\mc{D}}^{(1)}(\Dl\!\times\! F,\Dl\!\times\!G),}\end{equation}
where each $i$ is degreewise the canonical inclusion given by the
coproduct, is commutative. (\ref{diagrConmCuatroporCuatro}) is in
degrees $n,m$
$$\xymatrix@M=4pt@H=4pt@R=20pt@C=15pt{Z'                                     & X'\ar[l]_{g'}\ar[r]^{f'}                               & Y'\ar[r]^-{i_n}                                 & Y'\sqcup \coprod^{n}X'\sqcup Z'\\
                                      Z\ar[u]^{\alpha} \ar[d]_{\alpha'}      & X  \ar[l]_{g}\ar[r]^{f} \ar[u]^{\beta} \ar[d]_{\beta'} & Y \ar[u]^{\gamma} \ar[d]_{\gamma'}\ar[r]^-{i_n} & Y\sqcup \coprod^{n}X\sqcup Z\ar[u]^{\rho_n}\ar[d]_{\rho'_n} \\
                                      Z''\ar[d]_{i_m}                          & X''\ar[d]_{i_m}\ar[l]_{g''}\ar[r]^{f''}                   & Y''\ar[d]_{i_m}\ar[r]^-{i_n}                      & Y''\sqcup\coprod^{n}X''\sqcup Z''\ar[d]_{\phi_{n,m}}\\
                                      Z''\sqcup \coprod^{m} Z\sqcup Z'       & X''\sqcup \coprod^{m} X\sqcup X' \ar[l]_{G_m}\ar[r]^{F_m}          & Y''\sqcup \coprod^{m} Y\sqcup Y' \ar[r]^-{i}  & T_{n,m} }$$
where $T_{n,m}=(Y''\sqcup \coprod^{m} Y\sqcup Y')\sqcup$
$\ds{\coprod^n}$ $(X''\sqcup \coprod^{m} X\sqcup X')$ $\sqcup$
$(Z''\sqcup \coprod^{m} Z\sqcup Z')$.\\
On the other hand,{ if }$R={Cyl}_{\simp\mc{D}}^{(2)}(\rho'\times\Dl,\rho\times\Dl)$, we have that\\
$R_{n,m}= (Y''\sqcup\coprod^{n}X''\sqcup Z'')\sqcup$
$\ds{\coprod^m}$ $(Y\sqcup \coprod^{n}X\sqcup Z)$ $\sqcup$
$(Y'\sqcup \coprod^{n}X'\sqcup Z')$, that is obtained by
reordering the coproduct in $T_{n,m}$.\\
Therefore, let $\Theta_{n,m}:T_{n,m}\rightarrow R_{n,m}$ be the
canonical isomorphism that reorders the coproduct.\\
It is clear that
$$\xymatrix@R=10pt{                                                                                                   &(Y''\sqcup \coprod^{m} Y\sqcup Y')\sqcup\ds{\coprod^n}(X''\sqcup \coprod^{m} X\sqcup X')\sqcup(Z''\sqcup \coprod^{m} Z\sqcup Z') \ar[dd]_{\Theta_{n,m}}\\
                               Y''\sqcup\coprod^{n}X''\sqcup Z'' \ar@/_1.5pc/[dr]_-{i_n} \ar@/^1.5pc/[ur]^-{\phi_{n,m}}  &                                                                  \\
                                                                                                                      &(Y''\sqcup\coprod^{n}X''\sqcup Z'')\sqcup\ds{\coprod^m}(Y\sqcup \coprod^{n}X\sqcup Z)\sqcup(Y'\sqcup\coprod^{n}X'\sqcup Z')}$$
is commutative, and similarly $\Theta_{n,m}\comp i_m=\varphi_{n,m}$.\\
Hence, it remains to show that $\Theta=\{\Theta_{n,m}\}_{n,m}$ is
a morphism of bisimplicial objects.

\noindent Following the terminology in
\ref{def{i}CylenTerminosdeSigma}, write
$T_{n,m}={Cyl}_{\simp\mc{D}}^{(1)}(\Dl\times F,\Dl\times G)_{n,m}$
as
$$Cyl(\gamma',\gamma)_m^{u_1}\sqcup \ds \coprod_{\rho\in\Lambda_n} Cyl(\beta',\beta)_m^{\rho} \sqcup Cyl(\alpha',\alpha)_m^{u_0}=$$
$$ ({Y''}^{u_1,u_1}\sqcup\!\!\ds\coprod_{\sigma\in\Lambda_m}\!\! Y^{u_1,\sigma}\sqcup {Y'}^{u_1,u_0}) \sqcup
 \!\!\ds\coprod_{\rho\in\Lambda_n}\!\!({X''}^{\rho,u_1}\sqcup\!\!\ds\coprod_{\sigma\in\Lambda_m}\!\! X^{\rho,\sigma} \sqcup {X'}^{\rho,u_0})\sqcup
 ({Z''}^{u_0,u_1}\sqcup\!\!\ds\coprod_{\sigma\in\Lambda_m}\!\! Z^{u_0,\sigma}\sqcup{Z'}^{u_0,u_0})$$
On the other hand,
$R_{n,m}={Cyl}_{\simp\mc{D}}^{(2)}(\rho'\times\Dl,\rho\times\Dl)_{n,m}$
is
$$Cyl(f'',g'')_n^{u_1}\sqcup \ds \coprod_{\rho\in\Lambda_m} Cyl(f,g)_n^{\rho} \sqcup Cyl(f',g')_n^{u_0} =$$
$$({Y''}^{u_1,u_1} \sqcup\!\!\ds\coprod_{\sigma\in\Lambda_n}\!\! {X''}^{u_1,\sigma} \sqcup {Z''}^{u_1,u_0}) \sqcup\!\! \ds\coprod_{\rho\in\Lambda_m}\!\!( {Y}^{\rho,u_1} \sqcup\ds\coprod_{\sigma\in\Lambda_n} X^{\rho,\sigma} \sqcup {Z}^{\rho,u_0})\sqcup ({Y'}^{u_0,u_1} \sqcup\!\!\ds\coprod_{\sigma\in\Lambda_n}\!\! {X'}^{u_0,\sigma} \sqcup{Z'}^{u_0,u_0})$$
Then, $\Theta_{n,m}$ maps the component $(\rho,\sigma)$ of
$T_{n,m}$ into the component $(\sigma,\rho)$ of $R_{n,m}$.\\
If $\theta:[n']\rightarrow [n]$, the verif{i}cation of the
equalities $\Theta_{n',m}\comp T(\theta,Id)=
R(\theta,Id)\comp\Theta_{n,m}$, and $\Theta_{m,n'}\comp
T(Id,\theta)=
R(Id,\theta)\comp\Theta_{m,n}$ is a straightforward computation.\\
Let us see, for instance, the f{i}rst equality, because the second
one is totally similar. We have that
$$T(\theta,Id)|_{(\rho,\sigma)}=
\begin{cases}
 Id:Cyl(\beta',\beta)_m^{\rho}\rightarrow Cyl(\beta',\beta)_m^{\rho\theta}&\mbox{ if } \rho\theta\in\Lambda_{n'}\\
 F_m:Cyl(\beta',\beta)_m^{\rho}\rightarrow Cyl(\gamma',\gamma)_m^{u_1}&\mbox{ if } \rho\theta=u_1\mbox{ and }\rho\in\Lambda_n\\
 G_m:Cyl(\beta',\beta)_m^{\rho}\rightarrow Cyl(\alpha',\alpha)_m^{u_0}&\mbox{ if } \rho\theta=u_0\mbox{ and }\rho\in\Lambda_n\\
 Id:Cyl(\gamma',\gamma)_m^{u_1}\rightarrow Cyl(\gamma',\gamma)_m^{u_1}&\mbox{ if } \rho=u_1\\
 Id:Cyl(\alpha',\alpha)_m^{u_0}\rightarrow Cyl(\alpha',\alpha)_m^{u_0}&\mbox{ if } \rho=u_0\ .
\end{cases}$$
Note also that the restriction of $R(\theta,Id)\comp\Theta_{n,m}$
to the component $(\rho,\sigma)$ agrees with the restriction of
$R(\theta,Id)$ to the component $(\sigma,\rho)$, that is by
def{i}nition
$$R(\theta,Id)|_{(\sigma,\rho)}=
\begin{cases}
 Cyl(f,g)(\theta)|_{(\sigma,\rho)}:Cyl(f,g)_n^{\sigma}\rightarrow Cyl(f,g)_{n'}^{\sigma}&\mbox{ if } \sigma\in\Lambda_{m}\\
 Cyl(f'',g'')(\theta)|_{(\sigma,\rho)}:Cyl(f'',g'')_n^{u_1}\rightarrow Cyl(f'',g'')_{n'}^{u_1}&\mbox{ if } \sigma=u_1\\
 Cyl(f',g')(\theta)|_{(\sigma,\rho)}:Cyl(f',g')_n^{u_0}\rightarrow Cyl(f',g')_{n'}^{u_0}&\mbox{ if } \sigma=u_0\ .
\end{cases}$$
We remind that
$$Cyl(f,g)(\theta)|_{(\sigma,\rho)}=
\begin{cases}
 Id:X^{\sigma,\rho}\rightarrow X^{\sigma,\rho\theta}&\mbox{ if } \rho\theta\in\Lambda_{n'}\\
 Id:Y^{\sigma,u_1}\rightarrow X^{\sigma,u_1}&\mbox{ if } \rho=u_1\\
 Id:Z^{\sigma,u_0}\rightarrow X^{\sigma,u_0}&\mbox{ if } \rho=u_0\\
 f:X^{\sigma,\rho}\rightarrow Y^{\sigma,u_1}&\mbox{ if } \rho\in\Lambda_{n}\mbox{ and }\rho\theta=u_1\\
 g:X^{\sigma,\rho}\rightarrow Z^{\sigma,u_0}&\mbox{ if } \rho\in\Lambda_{n}\mbox{ and }\rho\theta=u_0\ ,
\end{cases}$$
and analogously for $ Cyl(f'',g'')(\theta)|_{(\sigma,\rho)}$ and
$Cyl(f',g')(\theta)|_{(\sigma,\rho)}$.\\
Assume that $\rho\theta\in \Lambda_{n'}$. Then
$$T(\theta,Id)|_{(\rho,\sigma)}=
\begin{cases}
  Id:X''^{\rho,u_1}\rightarrow X''^{\rho\theta,u_1} &\mbox{{ if }} \sigma=u_1\\
  Id:X^{\rho,\sigma}\rightarrow {X}^{\rho\theta,\sigma} &\mbox{{ if }} \sigma\in\Lambda_m \\
  Id:X'^{\rho,u_0}\rightarrow {X'}^{\rho\theta,u_0} &\mbox{{ if }} \sigma=u_0 \ .
\end{cases}$$
Hence, the result of interchanging the indexes in the above
formula is
$$ \Theta_{n,m}\comp T(\theta,Id)|_{(\rho,\sigma)}= \left\{
\begin{array}{ll}
  Id:{X''}^{u_1,\rho}\rightarrow {X''}^{u_1,\rho\theta}     &\mbox{ if } \sigma=u_1\\
  Id:X^{\sigma,\rho}\rightarrow {X}^{\sigma,\rho\theta} &\mbox{ if } \sigma\in\Lambda_m \\
  Id:{X'}^{u_0,\rho}\rightarrow {X'}^{u_0,\rho\theta}     &\mbox{ if } \sigma=u_0
\end{array}
\right\} = R(\theta,Id)|_{(\rho,\sigma)} \ ,$$
and the equality holds as in the remaining cases.

\end{proof}

\begin{def{i}}[Simplicial path object]\label{DefiPath}\mbox{}\\
The simplicial path functor\index{Index}{path object}
$Path:Co\Omega(\mc{D})\rightarrow
\Dl\mc{D}$\index{Symbols}{$Path$} is just the dual notion of
$Cyl$, consequently it satisfies the dual properties included in
this section.
\end{def{i}}


\section{Symmetric notions of cylinder and cone}\label{condeg+bis}

The biaugmented bisimplicial object $Z$ associated with an object
$(f,\epsilon)$ in $\Omega(\mc{D})$ is clearly asymmetric, and we
could consider as well the object $Z'$ obtained from $Z$ by
interchanging the indexes. The total simplicial object of $Z'$ is
another cylinder object associated with $(f,\epsilon)$ that will
be studied in this section.

\begin{num}\label{Upsilon}
Let $op:\Dl\rightarrow\Dl$ be the isomorphism of categories that
`reverses the order', introduced in \ref{OrdenOpuesto}.
Denote by
$\Upsilon:\simp\mc{D}\rightarrow\simp\mc{D}$\index{Symbols}{$\Upsilon$}
the functor obtained by composition with $op$.\\
Therefore
$$(\Upsilon X)_n=X_n \;\;\; d_i^{\Upsilon X}=d_{n-i}^X: X_n\rightarrow X_{n-1} \;\;\; s_j^{\Upsilon X}=s_{n-j}^X:X_n\rightarrow X_{n+1} .$$
Let $\Upsilon:\Omega(\mc{D})\rightarrow\Omega(\mc{D})$ be the
induced functor, that is also an isomorphism, and is given by
$\Upsilon(f,\epsilon)=(\Upsilon(f),\Upsilon(\epsilon))$.
\end{num}

The result of ``conjugate'' the cylinder functor $Cyl$ with
respect to $\Upsilon$ is the following alternative def{i}nition of
cylinder.

\begin{def{i}}\label{defCilPrima}
Set $Cyl'=\Upsilon\comp
Cyl\comp\Upsilon:\Omega(\mc{D})\rightarrow \simp\mc{D}$.\\
Given $\xymatrix@M=4pt@H=4pt@C=20pt{ X_{-1}\times\Dl & X
\ar[r]^{f}\ar[l]_-{\epsilon} &  Y }$ in $\Omega(\mc{D})$, then
$Cyl'(D)$\index{Symbols}{$Cyl'$} is in degree $n$
$$Cyl'(D)_n=Y_n\sqcup X_{n-1}\sqcup\cdots\sqcup X_0\sqcup X_{-1}.$$
The face morphisms $d_i^{Cyl'(D)}:Cyl'(D)_n\rightarrow
Cyl'(D)_{n-1}$ are def{i}ned as
$$ d_i^{Cyl'(D)}\mid_{Y_n}=d_i^Y\, ,\;\; d_i^{Cyl'(D)}\mid_{X_k}=\left\{\begin{array}{ll} d_i^{X}\! &\! i\leq k \\
                                                                                                                            Id    \!  & \!i>k\, (i,k)\!\neq\!(n,n\!-\!1)\\
                                                                                                                            f_{n-1}\! & \!(i,k)\neq (n,n-1) \end{array}\right.$$
The degeneracy maps $s_j^{Cyl'(D)}:Cyl'(D)_n\rightarrow
Cyl'(D)_{n+1}$ are
$$s_j^{Cyl'(D)}\mid_{Y_n}=s_j^Y\, ,\;\; s_j^{Cyl'(D)}\mid_{X_k}=\left\{\begin{array}{ll} s_j^{X} & j\leq k \\
                                                                                                                            Id        & j>k  \end{array}\right. .$$
Visually, for $0\leq i < n$, $d_i^{Cyl'(D)}$ is
$$\xymatrix@1@C=3pt{ Y_{n}\ar[rrd]_(.35){\partial_i^Y}  & \sqcup &  X_{n-1} \ar[rrd]_(.35){\partial_{i}^X} & \sqcup & X_{n-2} \ar[rrd]_(.35){\partial_{i}^X}& \sqcup & \cdots  &         & \sqcup  & X_{i}\ar[rrd]_(.35){\partial_{i}^X} & \sqcup & X_{i-1}\ar[d]^(.5){Id} & \sqcup & \cdots  & \sqcup &  X_0 \ar[d]^(.5){Id}   & \sqcup & X_{-1} \ar[d]^(.5){Id}& \\
                                                          &        &  Y_{n-1}                                    & \sqcup & X_{n-2}                                 & \sqcup & X_{n-3} &  \sqcup &         & \cdots                                  & \sqcup & X_{i-1}                          & \sqcup & \cdots  & \sqcup &  X_0                             & \sqcup & X_{-1}                & .}$$
The case $i=n$ is
$$\xymatrix@1@C=4pt{ Y_{n}\ar[rrd]_(.35){\partial_{n}^{Y}}  & \sqcup &  X_{n-1} \ar[d]^(.5){f_n} & \sqcup & X_{n-2} \ar[d]^(.5){Id}& \sqcup & X_{n-2} \ar[d]^(.5){Id} &  \sqcup & \cdots  & \sqcup &  X_0 \ar[d]^(.5){Id}   & \sqcup & X_{-1} \ar[d]^(.5){Id }& \\
                                                              &         &  Y_{n-1}                    & \sqcup & X_{n-2}                           & \sqcup & X_{n-2}                            &  \sqcup & \cdots  & \sqcup &  X_0                           & \sqcup & X_{-1} & .}$$
F{i}nally $s_j^{Cyl'(D)}$ is expressed as
$$\xymatrix@1@C=4pt{        &       & Y_{n}\ar[lld]_(.65){s_j^Y}  & \sqcup & X_{n-1} \ar[lld]_(.65){s_{j}^X} & \sqcup & X_{n-2} \ar[lld]_(.65){s_{j}^X}& \sqcup & \cdots    &  \sqcup & X_{j}\ar[lld]_(.65){s_{j}^X} & \sqcup & X_{j-1}\ar[d]^(.5){Id} & \sqcup & \cdots  & \sqcup & X_0 \ar[d]^(.5){Id}   & \sqcup & X_{-1} \ar[d]^(.5){Id }& \\
                    Y_{n+1} &\sqcup & X_n                         & \sqcup & X_{n-1}                           & \sqcup & \cdots                           &        & X_{j+1} &  \sqcup & X_{j}                        & \sqcup & X_{j-1}                            & \sqcup & \cdots  & \sqcup & X_0                           & \sqcup & X_{-1}                & .}$$

\end{def{i}}

\begin{num} Again, it follows from the properties of $Tot$ the
existence of canonical inclusions
$$ X_{-1}\times\Dl \rightarrow Cyl'(f,\epsilon)\ \ \ Y\rightarrow Cyl'(f,\epsilon)\ . $$
\end{num}

The next result is a consequence of the def{i}nitions of
 $Tot$ and $Cyl'$. We denote by
$$\Gamma:2-\Dl_+\longrightarrow 2-\Dl_+$$ the functor that interchanges the indexes of a biaugmented bisimplicial object.

\begin{prop} The functor $Cyl':\Omega(\mc{D})\rightarrow \simp\mc{D}$
agrees with the composition
$$ \Omega(\mc{D})\stackrel{\psi}{\longrightarrow} 2-\Dl_+ \stackrel{\Gamma}{\longrightarrow} 2-\Dl_+ \stackrel{Tot}{\longrightarrow} \simp\mc{D}\ ,$$
where $\psi$ is the functor given in $($\ref{def{i}PsiCilindroSimp}$)$.\\
In other words, the simplicial object $Cyl'(f,\epsilon)$ coincides
with the total object of the biaugmented bisimplicial object
obtained by interchanging the indexes in
$($\ref{diagramaCilindro}$)$, that is to say
$$\xymatrix@M=4pt@H=4pt@C=25pt{                                   &                                                                                             &                                                                                                                                              &                                                                                                                                                          &                                                                                                                             & \\
                                      X_{-1}\ar@{}[u]|{\vdots}\ar@<0ex>[d]\ar@<1ex>[d] \ar@<-1ex>[d]      & {X_{0}}\ar[l]\ar@<0ex>[d]\ar@<1ex>[d] \ar@<-1ex>[d]\ar@{}[u]|{\vdots}\ar@/_0.75pc/[r]     & {X_{1}} \ar@{}[u]|{\vdots}\ar@<0.5ex>[l] \ar@<-0.5ex>[l]\ar@<0ex>[d]\ar@<1ex>[d] \ar@<-1ex>[d]\ar@/_1pc/[r]\ar@/_0.75pc/[r]                &X_{2} \ar@{}[u]|{\vdots}\ar@<0ex>[l]\ar@<1ex>[l] \ar@<-1ex>[l]\ar@<0ex>[d]\ar@<1ex>[d] \ar@<-1ex>[d]\ar@/_1pc/[r]\ar@/_0.75pc/[r]\ar@{-}@/_1.25pc/[r]   & X_{3}\ar@{}[u]|{\vdots}\ar@<0.33ex>[l]\ar@<-0.33ex>[l]\ar@<1ex>[l]\ar@<-1ex>[l]\ar@<0ex>[d]\ar@<1ex>[d] \ar@<-1ex>[d]     & \cdots \\
                                      X_{-1}\ar@<0.5ex>[d] \ar@<-0.5ex>[d] \ar@/^1pc/[u]\ar@/^0.75pc/[u]  & {X_{0}}\ar[l]\ar@<0.5ex>[d] \ar@<-0.5ex>[d] \ar@/^1pc/[u]\ar@/^0.75pc/[u]\ar@/_0.75pc/[r] &  X_{1} \ar@/^1pc/[u]\ar@/^0.75pc/[u]\ar@<0.5ex>[l]\ar@<-0.5ex>[l]\ar@<0.5ex>[d] \ar@<-0.5ex>[d]\ar@/_1pc/[r]\ar@/_0.75pc/[r]               &X_{2} \ar@/^1pc/[u]\ar@/^0.75pc/[u] \ar@<0ex>[l]\ar@<1ex>[l] \ar@<-1ex>[l] \ar@<0.5ex>[d] \ar@<-0.5ex>[d]\ar@/_1pc/[r]\ar@/_0.75pc/[r]\ar@/_1.25pc/[r]  & X_{3}\ar@<0.33ex>[l]\ar@<-0.33ex>[l]\ar@<1ex>[l]\ar@<-1ex>[l] \ar@<0.5ex>[d] \ar@<-0.5ex>[d]\ar@/^1pc/[u]\ar@/^0.75pc/[u] & \cdots \\
                                      X_{-1}  \ar@/^0.75pc/[u]                                            &  X_{0}\ar[d] \ar[l]\ar@/_0.75pc/[r] \ar@/^0.75pc/[u]                                      & {X_{1}} \ar[d] \ar@<0.5ex>[l] \ar@<-0.5ex>[l]  \ar@/^0.75pc/[u]  \ar@/_1pc/[r]\ar@/_0.75pc/[r]                                             &X_{2} \ar[d] \ar@<0ex>[l]\ar@<1ex>[l] \ar@<-1ex>[l]\ar@/^0.75pc/[u]\ar@/_1pc/[r]\ar@/_0.75pc/[r]\ar@/_1.25pc/[r]                                        & X_{3}\ar[d] \ar@<0.33ex>[l]\ar@<-0.33ex>[l]\ar@<1ex>[l]\ar@<-1ex>[l] \ar@/^0.75pc/[u]                                     & \cdots \\
                                                                                                         &  Y_{0} \ar@/_0.75pc/[r]                                                                        & Y_{1} \ar@<0.5ex>[l]\ar@<-0.5ex>[l]\ar@/_1pc/[r]\ar@/_0.75pc/[r]                                                                             &Y_{2} \ar@<0ex>[l]\ar@<1ex>[l] \ar@<-1ex>[l]\ar@/_1pc/[r]\ar@/_0.75pc/[r]\ar@/_1.25pc/[r]                                                                   & Y_{3}\ar@<0.33ex>[l]\ar@<-0.33ex>[l]\ar@<1ex>[l]\ar@<-1ex>[l]                                                                 & \cdots \ .}$$
\end{prop}

\begin{obs}
The functor $Cyl'$ satisf{i}es as well the analogous propositions
to \ref{CylCoproducto}, \ref{CylHomotopia} and \ref{Cylaumentado},
that we will labelled as \ref{CylCoproducto}', \ref{CylHomotopia}'
and \ref{Cylaumentado}'.
\end{obs}

\begin{def{i}}\label{def{i}ConoSimetrico}
If $\mc{D}$ has a f{i}nal object, the following alternative notion
of simplicial cone\index{Symbols}{$C'$} is induced by $Cyl'$
$$C'=\Upsilon\comp C\comp \Upsilon:Fl(\simp\mc{D})\rightarrow\simp\mc{D}\ .$$
\end{def{i}}

The following proposition will be useful in the next chapter, and
it will be a key point in the study of the relationship between
the cone and cylinder axioms.

\begin{num}
Consider a commutative diagram in $\simp\mc{D}$
\begin{equation}\label{diagrCylIterado}\xymatrix@R=33pt@C=33pt@H=4pt@M=4pt@!0{                                 & X_{-1}\times\Dl  \ar[rr]^-{h}\ar[ld]_{r}   &                   & U_{-1}\times\Dl \\
                                                                   Y_{-1}\times\Dl                             &                                            &                   &\\
                                                                                                               & X \ar[uu]_{\beta} \ar[rr]^-{g} \ar[ld]_{q} &                   & U \ ,\ar[uu]_{\gamma}\ar[ld]^{p}\\
                                                                              Y \ar[uu]^{\alpha}   \ar[rr]_{f}   &                                            &            V      & }\end{equation}
We will denote by
$$\xymatrix@M=4pt@H=4pt{ U_{-1}\times\Dl  & Cyl'(g,\beta)\ar[r]^-{t} \ar[l]_-{\delta} & Cyl'(f,\alpha)} $$
the object of $\Omega(\mc{D})$ obtained by applying $Cyl'$ in one
direction, and by
$$\xymatrix@M=4pt@H=4pt{ Y_{-1}\times\Dl & Cyl(q,\beta)\ar[r]^-{u} \ar[l]_-{\zeta} &  Cyl(p,\gamma)} \ .$$
the result of applying $Cyl$ to (\ref{diagrCylIterado}) in the
other sense.
\end{num}

\begin{prop}\label{lemaDiagrCubico}
There exists a natural isomorphism in $($\ref{diagrCylIterado}$)$
$$Cyl'(u,\zeta)\simeq Cyl(t,\delta) \ .$$
\end{prop}

\begin{proof}
We can add in a suitable way two new simplicial indexes to each
simplicial object in (\ref{diagrCylIterado}) in order to obtain a
$3$-augmented $3$-simplicial object $T$ (see \ref{objeto3Delta}).
In other words, for $i,j,k\geq 0$ def{i}ne
$$\begin{array}{cccc} T_{i,j,k}=X_j &  T_{i,j,-1}=U_j &   T_{-1,j,k}=Y_j & T_{-1,j,-1}=V_j \\
                      T_{i,-1,k}=X_{-1} & T_{i,-1-1}=U_{-1} & T_{-1,-1,k}=Y_{-1} &
\end{array}$$
It follows from \ref{TotIterado} that $Tot\comp {Tot}_{(0)}(T)\simeq Tot\comp {Tot}_{(2)}(T)$.\\
If $i\geq 0$, we can f{i}x as $i$ the f{i}rst index of $T$ and
apply $Tot^+$. The result is the augmented simplicial object
$Cyl'(g,\beta)\rightarrow U_{-1}\times\Dl$, whereas{ if }$i=-1$ we
obtain $Cyl'(f,\alpha)$. Hence
$${Tot}_{(0)}(T)([n],[m])=\begin{cases}
                            Cyl'(g,\beta)_m & \mbox{if } n,m\geq 0\\
                            U_{-1} & \mbox{if } m=-1\\
                            Cyl'(f,\alpha) & \mbox{if } n=-1 \ .
                           \end{cases}$$
Then, ${Tot}_{(0)}(T)$ is the biaugmented bisimplicial object
associated with
$$\xymatrix@M=4pt@H=4pt{ Cyl'(f,\alpha) &
Cyl'(g,\beta)\ar[l]_-{t} \ar[r]^-{\delta} & U_{-1}\times\Dl}\ ,$$
in (\ref{ejObjBisimpfEpsilon}), so $Tot\comp
{Tot}_{(0)}(T)=Cyl(t,\delta)$.

On the other hand,{ if }we f{i}x $k\geq 0$ in $T$ and apply
$Tot^+$ we get $Cyl(q,\beta)\rightarrow Y_{-1}\times\Dl$, whereas
setting $k=-1$ and applying $Tot$ we obtain $Cyl(p,\gamma)$.
Consequently,
$${Tot}_{(2)}(T)([n],[m])=\begin{cases}
                            Cyl(q,\beta)_n & \mbox{if } n,m\geq 0\\
                            Y_{-1} & \mbox{if } n=-1\\
                            Cyl(p,\gamma)_n & \mbox{if } m=-1 \ .
                           \end{cases}$$
Therefore $Tot\comp {Tot}_{(2)}(T)=Cyl'(u,\zeta)$, and we are
done.
\end{proof}

\begin{cor}\label{relacionCylCylPrima}
Let $f:A\rightarrow B$ and $g:A\rightarrow C$ morphisms in
$\mc{D}$. Then
$$Cyl(f\times\Dl,g\times\Dl)\simeq Cyl'(g\times\Dl,f\times\Dl)\ .$$
Moreover, this isomorphism is compatible with the respective
inclusion of $B\times\Dl$ and $C\times\Dl$ into both cylinder
objects.
\end{cor}

\begin{proof}
Let 0 be the initial object of $\mc{D}$, and $D_\cdot$ be the
constant simplicial object $D\times\Dl$, for every object $D$ in
$\mc{D}$. It is enough to apply the above proposition to
$$\xymatrix@R=30pt@C=30pt@H=4pt@M=4pt@!0{                       & A_\cdot  \ar[rr]^-{g_\cdot}\ar[ld]_{f_\cdot}   &                   & C_\cdot \\
                                           B_\cdot           &                                            &                   &\\
                                                                                                               & 0_\cdot \ar[uu] \ar[rr] \ar[ld] &                   & 0_\cdot \ ,\ar[uu]\ar[ld]\\
                                           0_\cdot \ar[uu]   \ar[rr]  &                                            &            0_\cdot      & }$$
and note that $Cyl'(D_\cdot \leftarrow 0_\cdot \rightarrow
0_\cdot)=Cyl(D_\cdot \leftarrow 0_\cdot \rightarrow
0_\cdot)=D_\cdot$ for any $D$ in $\mc{D}$. In addition, since the
isomorphism obtained in this way is just to reorder a coproduct,
it is clear that the canonical inclusions of $B\times\Dl$ and
$C\times\Dl$ are preserved.
\end{proof}


\section{Cubical cylinder object}

The construction developed in this section is just a
generalization of the cubical cylinder object $\widetilde{Cyl}(X)$
associated with a simplicial object $X$ in $\mc{D}$,
\ref{ejCilindroCub}.

\begin{def{i}}\label{catSquareCont}
Let $\square_1$\index{Symbols}{$\square_1$} be the category
$$\xymatrix@H=4pt@M=4pt{\bullet\ar[r] &\bullet &\ar[l] \bullet}\ .$$
Then,{ if }we f{i}x a category $\cont$, the category
$\square_1^\comp\cont$ has as objects the diagrams in $\cont$
$$\xymatrix@H=4pt@M=4pt{Z  &X \ar[r]^f\ar[l]_g & Y}\ ,$$
that will be represented by $(f,g)$.\\
The morphisms in $\square_1^\comp\cont$ are commutative diagrams
$$\xymatrix@M=4pt@H=4pt{ Z \ar[d] & X \ar[r] \ar[l]\ar[d] &  Y \ar[d]\\
                         Z' & X' \ar[r] \ar[l] &  Y'  . }$$
Similarly, $\square_1\cont$ is the category whose objects are
diagrams
$$\xymatrix@H=4pt@M=4pt{Z\ar[r]^g &X &\ar[l]_f Y}\ .$$
\end{def{i}}

\begin{def{i}}
Def{i}ne the functor $\Phi:\square_1^\comp\simp\mc{D}\rightarrow
2-\Dl_{+}^\comp\mc{D}$ as follows.\\
Given $(f,g)\in \square_1^\comp\simp\mc{D}$, the biaugmented
bisimplicial object $\Phi(f,g)$ is
$$\xymatrix@H=4pt@M=4pt{{T}:\ar@{}[r] & T_{-1,\cdot} &T^+_{\cdot,\cdot} \ar[r]^-{\zeta}\ar[l]_-{\epsilon} & T_{\cdot,-1}}\ ,$$
where $T^+_{\cdot,\cdot}=Dec(X)|_{\Dl\times\Dl}$ (see
\ref{ejCilindroCub}), $T_{-1,\cdot}=Y_\cdot$ and
$T_{\cdot,-1}=Z_\cdot$. In other words
$$ T_{i,j} = \left\{ \begin{array}{ll}  X_{i+j+1} & \mbox{ if } i,j\geq 0 \\
                                          Y_{j} & \mbox{ if } i = -1\\
                                          Z_{i}  & \mbox{ if } j = -1\ . \end{array}\right. $$
Visually, $T$ can be visualized as
\vspace{-0.5cm}
$$\xymatrix@M=4pt@H=4pt@C=25pt{                                   &                                                                                             &                                                                                                                                              &                                                                                                                                                          &                                                                                                                             & \\
                                      Y_{2}\ar@{}[u]|{\vdots}\ar@<0ex>[d]\ar@<1ex>[d] \ar@<-1ex>[d]      & {X_{3}}\ar[l]\ar@<0ex>[d]\ar@<1ex>[d] \ar@<-1ex>[d]\ar@{}[u]|{\vdots}\ar@/_0.75pc/[r]     & {X_{4}} \ar@{}[u]|{\vdots}\ar@<0.5ex>[l] \ar@<-0.5ex>[l]\ar@<0ex>[d]\ar@<1ex>[d] \ar@<-1ex>[d]\ar@/_1pc/[r]\ar@/_0.75pc/[r]                &X_{5} \ar@{}[u]|{\vdots}\ar@<0ex>[l]\ar@<1ex>[l] \ar@<-1ex>[l]\ar@<0ex>[d]\ar@<1ex>[d] \ar@<-1ex>[d]\ar@/_1pc/[r]\ar@/_0.75pc/[r]\ar@{-}@/_1.25pc/[r]   & X_{6}\ar@{}[u]|{\vdots}\ar@<0.33ex>[l]\ar@<-0.33ex>[l]\ar@<1ex>[l]\ar@<-1ex>[l]\ar@<0ex>[d]\ar@<1ex>[d] \ar@<-1ex>[d]     & \cdots \\
                                      Y_{1}\ar@<0.5ex>[d] \ar@<-0.5ex>[d] \ar@/^1pc/[u]\ar@/^0.75pc/[u]  & {X_{2}}\ar[l]\ar@<0.5ex>[d] \ar@<-0.5ex>[d] \ar@/^1pc/[u]\ar@/^0.75pc/[u]\ar@/_0.75pc/[r] &  X_{3} \ar@/^1pc/[u]\ar@/^0.75pc/[u]\ar@<0.5ex>[l]\ar@<-0.5ex>[l]\ar@<0.5ex>[d] \ar@<-0.5ex>[d]\ar@/_1pc/[r]\ar@/_0.75pc/[r]               &X_{4} \ar@/^1pc/[u]\ar@/^0.75pc/[u] \ar@<0ex>[l]\ar@<1ex>[l] \ar@<-1ex>[l] \ar@<0.5ex>[d] \ar@<-0.5ex>[d]\ar@/_1pc/[r]\ar@/_0.75pc/[r]\ar@/_1.25pc/[r]  & X_{5}\ar@<0.33ex>[l]\ar@<-0.33ex>[l]\ar@<1ex>[l]\ar@<-1ex>[l] \ar@<0.5ex>[d] \ar@<-0.5ex>[d]\ar@/^1pc/[u]\ar@/^0.75pc/[u] & \cdots \\
                                      Y_{0}  \ar@/^0.75pc/[u]                                            &  X_{1}\ar[d] \ar[l]\ar@/_0.75pc/[r] \ar@/^0.75pc/[u]                                      & {X_{2}} \ar[d] \ar@<0.5ex>[l] \ar@<-0.5ex>[l]  \ar@/^0.75pc/[u]  \ar@/_1pc/[r]\ar@/_0.75pc/[r]                                             &X_{3} \ar[d] \ar@<0ex>[l]\ar@<1ex>[l] \ar@<-1ex>[l]\ar@/^0.75pc/[u]\ar@/_1pc/[r]\ar@/_0.75pc/[r]\ar@/_1.25pc/[r]                                        & X_{4}\ar[d] \ar@<0.33ex>[l]\ar@<-0.33ex>[l]\ar@<1ex>[l]\ar@<-1ex>[l] \ar@/^0.75pc/[u]                                     & \cdots \\
                                                                                                         &  Z_0 \ar@/_0.75pc/[r]                                                                        & Z_1  \ar@<0.5ex>[l]\ar@<-0.5ex>[l]\ar@/_1pc/[r]\ar@/_0.75pc/[r]                                                                         &Z_2 \ar@<0ex>[l]\ar@<1ex>[l] \ar@<-1ex>[l]\ar@/_1pc/[r]\ar@/_0.75pc/[r]\ar@/_1.25pc/[r]                                                                   & Z_3\ar@<0.33ex>[l]\ar@<-0.33ex>[l]\ar@<1ex>[l]\ar@<-1ex>[l]                                                                 & \cdots }$$

\noindent The horizontal augmentations, $\epsilon_n:X_n\rightarrow
Y_{n-1}$, are equal to $f_{n-1}d_0$, whereas the vertical ones,
$\zeta_n:X_n\rightarrow Z_{n-1}$, are $g_{n-1}d_n$.
\end{def{i}}
\begin{def{i}}[Cubical cylinder object]\label{defCilCubico}\mbox{}\\
We def{i}ne the cubical cylinder functor\index{Index}{cubical
cylinder}
$\widetilde{Cyl}:\square_1^\comp\simp\mc{D}\rightarrow\simp\mc{D}$\index{Symbols}{$\widetilde{Cyl}$}
as the composition
$$\xymatrix@H=4pt@M=4pt{\square_1^\comp\simp\mc{D}\ar[r]^-{\Phi} & 2-\Dl_{+}^\comp\mc{D}\ar[r]^-{Tot}& \simp\mc{D}}\ .$$
In other words, $\widetilde{Cyl}(f,g)=Tot(\Phi(f,g))$, that is in
degree $n$
$$ \widetilde{Cyl}(f,g)_n = Y_n \sqcup X_{n} \sqcup  \stackrel{n)}{\cdots} \sqcup X_n \sqcup Z_{n} \ .$$
\end{def{i}}

\begin{num}\label{def{i}CilTildeAbstracto}
Equivalently, consider the maps $u_i:[n]\rightarrow [1]$ given by
$u_i([n])=i$, $i=0,1$, and let $\Lambda$\index{Symbols}{$\Lambda$}
be the set of morphisms $\sigma:[n]\rightarrow [1]$ dif{f}erent
from $u_0$ and $u_1$. Then
$$ \widetilde{Cyl}(f,g)_n = \widetilde{Cyl}(f,g)_n^{u_1} \sqcup\ds\coprod_{\sigma\in\Lambda} \widetilde{Cyl}(f,g)_n^\sigma \sqcup \widetilde{Cyl}(f,g)_n^{u_0} \ ,$$
where $\widetilde{Cyl}(f,g)_n^{u_1}=Y_n$, $\widetilde{Cyl}(f,g)_{n}^{u_0}=Z_n$ and $\widetilde{Cyl}(f,g)_n^\sigma=X_n$ $\forall\sigma\in\Lambda$.\\
If $\theta:[m]\rightarrow [n]$ is a morphism in $\Dl$, it follows
that the restriction of $\widetilde{Cyl}(f,g)(\theta):
\widetilde{Cyl}(f,g)_n \rightarrow \widetilde{Cyl}(f,g)_m$ to the
component $\sigma\in\Dl/[1]$ is
$$
\widetilde{Cyl}(f,g)(\theta)|_{\widetilde{Cyl}(f,g)_n^\sigma}=\begin{cases}
                                                                 X(\theta):X_n^\sigma\rightarrow X_m^{\sigma\theta} &\mbox{ if } \sigma\theta\in\Lambda\\
                                                                     f_m X(\theta):X_n^\sigma\rightarrow Y_m^{u_1} &\mbox{ if } \sigma\in\Lambda, \sigma\theta=u_1\\
                                                                      g_m X(\theta):X_n^\sigma\rightarrow Z_m^{u_1} &\mbox{ if } \sigma\in\Lambda, \sigma\theta=u_0\\
                                                                      Y(\theta):Y_n^{u_1}\rightarrow Y_m^{u_1} &\mbox{ if } \sigma=u_1\\
                                                                      Z(\theta):Z_n^{u_1}\rightarrow Z_m^{u_1} &\mbox{ if } \sigma=u_0\ .
                                                            \end{cases}
                                                            $$
\end{num}

\begin{num}\label{inclusionesCylTilde}
Consider $\xymatrix@M=4pt@H=4pt@C=15pt{Z & X
\ar[r]^-{f}\ar[l]_-{g} & Y }$ in $\square_1^\comp\simp\mc{D}$. The
canonical morphisms
$$j_{Y_n}:Y_n\rightarrow \widetilde{Cyl}(f,g)_n\mbox{ and }j_{Z_n}:Z_n\rightarrow\widetilde{ Cyl}(f,g)_n$$
give rise to the following diagram, natural in $(f,g)$
\begin{equation}\label{CuadradoCylTilde}\xymatrix@M=4pt@H=4pt{ X \ar[r]^{f}   \ar[d]_{g}           & Y \ar[d]^{j_Y}\\
                       Z \ar[r]^-{j_Z}&  \widetilde{Cyl}(f,g)\ . }
                       \end{equation}
\end{num}

\begin{obs} If $X$ is a simplicial object in $\mc{D}$,
then $\widetilde{Cyl}(Id_X,Id_X)$ is just $\widetilde{Cyl}(X)$,
the cubical cylinder object associated with $X$, that was
introduced in \ref{ejCilindroCub}.
\end{obs}

\begin{obs}
We will explain here why the cylinder considered in this section
is called ``cubical''.\\
Firstly, the category $\square_1^\comp\cont$ is a subcategory of
the category of (all) ``cubical diagrams'' in $\cont$ introduced
in \cite{GN}.
If $\cont$ has coproducts, there exists a functor
$$\square_1^\comp\cont \rightarrow \Dl_e^\comp\cont$$
that assigns to $\xymatrix@H=4pt@M=4pt{Z  &X \ar[r]^f\ar[l]_g &
Y}$ in $\cont$ the strict simplicial object $E(f,g)$ given by
$$\xymatrix@1{ Y\sqcup Z \;  && {\;} X \;
\ar@<0.5ex>[ll]^-{f} \ar@<-0.5ex>[ll]_-{g}  && \; 0 \;
\ar@<0ex>[ll] \ar@<1ex>[ll] \ar@<-1ex>[ll]  && {\;} 0
\ar@<0.33ex>[ll]\ar@<-0.33ex>[ll]\ar@<1ex>[ll]\ar@<-1ex>[ll]&
\cdots\cdots }\ ,$$ where $0${ if }the initial object in
$\cont$.\\
On the other hand, we have the Dold-Puppe transform
$\pi:\Dl_e^\comp\cont\rightarrow \simp\cont$ (see \ref{Dold-Puppe}).\\
Setting $\cont=\simp\mc{D}$, then
$\pi(E(f,g))\in\simp\simp\mc{D}$, and its diagonal is just
$\widetilde{Cyl}(f,g)$.
\end{obs}
\begin{prop}
The functor $\widetilde{Cyl}:
\square_1^\comp\simp\mc{D}\rightarrow\simp\mc{D}$ commutes with
coproducts, that is
$$\widetilde{Cyl}(f,g)\sqcup \widetilde{Cyl}(f',g')\simeq \widetilde{Cyl}(f\sqcup f',g\sqcup g')$$
\end{prop}

\begin{proof} The statement is a consequence of the commutativity
of $\Phi$ and $Tot$ with coproducts.
\end{proof}

\begin{prop}\label{HomotCuadrCylTilde}
Given $\xymatrix@M=4pt@H=4pt@C=15pt{Z & X \ar[r]^-{f}\ar[l]_-{g} &
Y }$ in $\square_1^\comp\simp\mc{D}$, the diagram
$($\ref{CuadradoCylTilde} $)$ commutes up to simplicial homotopy
equivalence, natural in $(f,g)$.
\end{prop}

\begin{proof}
Applying $\widetilde{Cyl}$ to the following morphism of
$\square_1^\comp\simp\mc{D}$
$$\xymatrix@M=4pt@H=4pt@C=25pt{X\ar[d]_g & X \ar[r]^-{Id}\ar[l]_-{Id}\ar[d]_{Id} &  X\ar[d]_f\\
                               Z & X \ar[r]^-{f}\ar[l]_-{g} &  Y \ ,}$$
we obtain $H:\widetilde{Cyl}(X)\rightarrow
\widetilde{Cyl}(f,g)$.\\
Following the notations introduced in \ref{HomotbyCylTilde}, from
the naturality of (\ref{CuadradoCylTilde}) we deduce that the
diagram
$$\xymatrix@M=4pt@H=4pt@C=20pt{X\ar[d]_g\ar[r]^-{J} & \widetilde{Cyl}(X) \ar[d]_{H} &  X\ar[d]_f\ar[l]_-{I}\\
                               Z\ar[r]^-{j_Z} & \widetilde{Cyl}(f,g) &  Y\ar[l]_-{j_Y} \ ,}$$
commutes. Then $H\comp I=j_Y\comp f$ and $H\comp J= j_Z \comp g$,
therefore $j_Y\comp f\sim j_Z \comp g$.
\end{proof}

\begin{obs} The functor $\widetilde{Cyl}$ also satisf{i}es an analogue of \ref{Cylaumentado}, that will not be used in this work.\\
Given a commutative diagram in $\simp\mc{D}$
\begin{equation}\label{CuadradConmdeObfSimpl}\xymatrix@M=4pt@H=4pt{ X \ar[r]^f   \ar[d]_{g}           & Y \ar[d]^{\rho'}\\
                       Z \ar[r]^{\rho}&  T ,}
                       \end{equation}
there exists a morphism, natural in (\ref{CuadradConmdeObfSimpl}),
$H:\widetilde{Cyl}(f,g)\rightarrow T$ such that $H\comp
j_{Z}=\rho$ and $H\comp j_Y=\rho'$, that is
$$\xymatrix@M=4pt@H=4pt{ X \ar[r]^f   \ar[d]_{g}          & Y \ar@/^1pc/[rdd]^{\rho'} \ar[d]^{j_Y}                 &\\
                       Z \ar@/_1pc/[drr]_{\rho}  \ar[r]^-{j_{Z}}  &  \widetilde{Cyl}(f,g)\ar@{}[d]|{\sharp}  \ar@{}[r]|{\sharp}\ar[dr]^{H}   &\\
                                                                 &                                                 & T\ .} $$
Indeed, it is enough to consider $H$ such that
$H_n|_{X_n}=\rho'_nf_n:X_n\rightarrow T_n$, where $X_n$ denotes a
component of $\widetilde{Cyl}(f,g)_n$.
\end{obs}

\begin{obs} A property of ``factorization'' (similar to \ref{factorizCylTilde}) also holds for $\widetilde{Cyl}$, with respect
to a diagram (\ref{diagrTresporTres}) of simplicial objects
(introduced in \ref{NotacionesFactCyl}). This property will not be
used in this work.
\end{obs}

Next we study the relationship between $Cyl$ and $\widetilde{Cyl}$
(in those cases in which they are comparable).

The following result is a direct consequence of the def{i}nitions
of $Cyl$ and $\widetilde{Cyl}$.

\begin{prop}\label{relacionCylCylTildeCtes}
If $\xymatrix@M=4pt@H=4pt@C=15pt{C & A \ar[r]^-{f}\ar[l]_-{g} &
B}$ is a diagram in $\mc{D}$, then
$$ Cyl(f\times\Dl,g\times\Dl) = \widetilde{Cyl}(f\times\Dl,g\times\Dl)\ .$$
\end{prop}

\begin{prop}\label{relacionCylCylTildeDiag}
Let $\xymatrix@M=4pt@H=4pt@C=15pt{Z & X \ar[r]^-{f}\ar[l]_-{g} &
Y}$ be a diagram in $\simp\mc{D}$. Then, following the notations
given in \ref{cylBisimpl}, the diagonal simplicial object of the
bisimplicial object $Cyl^{(1)}_{\simp\mc{D}}(\Dl\times f,\Dl\times
g)$ is equal to $\widetilde{Cyl}(f,g)$.\\
Similarly, $\mrm{D}Cyl^{(2)}_{\simp\mc{D}}(f\times \Dl,g\times
\Dl)=\widetilde{Cyl}(f,g)$.
\end{prop}

\begin{proof}
The bisimplicial object $T=Cyl_{\simp\mc{D}}^{(1)}(\Dl\times
f,\Dl\times g)$ is
$$T_{n,m}=Y_m \sqcup X_{m}^{(n-1)} \sqcup  {\cdots} \sqcup X_m^{(0)} \sqcup Z_{m} $$
where $X_m^{(n-i)}=X_m$ $\forall i=1,\ldots ,n$.
The face maps $d_i^{(2)}:T_{n,m}\rightarrow T_{n,m-1}$ respect to
the second index are $d_i^{(2)}=d_i^Y \sqcup d_i^X \sqcup
\stackrel{n)}{\cdots} \sqcup d_i^X  \sqcup d_i^Z$, and
analogously for the degeneracy maps $s_k^{(2)}:T_{n,m}\rightarrow T_{n,m+1}$.\\
On the other hand $d_i^{(1)}:T_{n,m}\rightarrow T_{n-1,m}$ is
$d_i^{(1)} |_{Y_{m}}=Id$, $d_i^{(1)}|_{Z_{m}}=Id$ and
$$ d_i^{(1)} |_{X_m^{(n-k)}} = \left\{ \begin{array}{ll} Id:X_m^{(n-k)}\rightarrow X_m^{(n-k-1)} &\mbox{ if } i \geq k \; \mbox{ and } (k,i)\neq (1,0) \\
                                                  Id:X_m^{(n-k)}\rightarrow X_m^{(n-k)} & \mbox{ if } i > k \; \mbox{ and } (k,i)\neq (n,n)\\
                                                  f_{m}:X_m^{(n-1)}\rightarrow Y_m & \mbox{ if } (k,i)=(1,0)\\
                                                  g_{m}:X_m^{(0)}\rightarrow Z_m  & \mbox{ if } (k,i)=(n,n)
                                                  \end{array}\right.\ .$$
The degeneracy maps are built in a similar way using the
def{i}nition of $Cyl$.\\
Clearly, the diagonal of $T$, $\mrm{D}T$, coincides with
$\widetilde{Cyl}(f,g)$.
The last statement follows from the commutativity of diagram
(\ref{diagGammaOmega}) and from the fact $\mrm{D}\Gamma=\mrm{D}$.
\end{proof}

\begin{prop}\label{CylRetractoCylTilde}
If $\xymatrix@M=4pt@H=4pt@C=15pt{ X_{-1}\times\Dl & X
\ar[r]^-{f}\ar[l]_-{\epsilon} &  Y }$ is a diagram in
$\simp\mc{D}$, then $Cyl(f,\epsilon)$ is a retract of
$\widetilde{Cyl}(f,\epsilon)$.\\
In other words, there exists morphisms
$\alpha:{Cyl}(f,\epsilon)\rightarrow\widetilde{Cyl}(f,\epsilon)$
and $\beta:\widetilde{Cyl}(f,\epsilon)\rightarrow
{Cyl}(f,\epsilon)$ such that $\beta\alpha=Id$.\\
In addition, $\alpha$ and $\beta$ are natural in $(f,\epsilon)$
and commute with the inclusions of $X_{-1}$ and $Y$ into the
respective cylinders.
\end{prop}

\begin{proof}
We have that
$$\alpha_n:Y_n \sqcup X_{n-1} \sqcup  {\cdots} \sqcup X_0 \sqcup X_{-1}\longrightarrow Y_n \sqcup X_{n}^{(n-1)} \sqcup \stackrel{n)}{\cdots} \sqcup X_n^{(0)} \sqcup X_{-1}$$
is def{i}ned as the identity on $Y_n$ and $X_{-1}$, and on
$X_{n-k}$ is $(s_0)^k:X_{n-k}\rightarrow X_n^{(n-k)}$.\\
It holds that $\alpha$ commutes with the face morphisms, since
$$\alpha_{n-1} d_i|_{X_{n-k}}=\left\{ \begin{array}{ll} (s_0)^k d_{i-k}^{X} & \mbox{ if } i \geq k \\
                                                        (s_0)^{k-1} & \mbox{ if } i < k \; \mbox{ and } (k,i)\neq (1,0)\\
                                                        f_{n-1} & \mbox{ if } (k,i)=(1,0)\end{array}\right.\ ,$$
whereas
$$ d_i\alpha_{n}|_{X_{n-k}}=\left\{ \begin{array}{ll}   d_{i}^{X}(s_0)^k  & \mbox{ if } i \geq k \\
                                                        d_i(s_0)^k & \mbox{ if } i < k \; \mbox{ and } (k,i)\neq (1,0)\\
                                                        f_{n-1}d_0 s_0 & \mbox{ if } (k,i)=(1,0)\end{array}\right.\ .$$
The equality $(s_0)^k d_{i-k}=d_{i}^{X}(s_0)^k$ follows from the
iteration of the simplicial identity
$d_{j+1}s_0=s_0 d_j${ if }$j> 1$.\\
In addition, since $(s_0)^l=s_{l-1}(s_{0})^{l-1}$, then
$d_i(s_0)^k=d_i(s_0)^i(s_0)^{k-i}=d_is_{i-1}(s_0)^{k-1}=(s_0)^{k-1}$.
One can check similarly that $\alpha$ commutes with the degeneracy
maps.\\
On the other hand, $\beta_n:Y_n \sqcup X_{n}^{(n-1)} \sqcup
\stackrel{n)}{\cdots} \sqcup X_n^{(0)} \sqcup X_{-1}\rightarrow
Y_n \sqcup X_{n-1} \sqcup  {\cdots} \sqcup X_0 \sqcup X_{-1}$ is
the identity on $Y_n$ and $X_{-1}$, and on $X_{n}^{(n-k)}$ is
$(d_0)^k :X_{n}^{(n-k)}\rightarrow X_{n-k}$.\\
We deduce again from the simplicial identities that $\beta$ is in
fact a morphism in $\simp\mc{D}$, as well as the equality
 $\beta\alpha=Id$ holds, and it is clear that the inclusions of $Y$ and $X_{-1}$ commutes with both morphisms.
 \end{proof}

\begin{cor}\label{HomotCuadradoCyl} The diagram
$($\ref{CuadradoCyl}$)$ given in \ref{inclusionesCyl} commutes up
to simplicial homotopy.
\end{cor}

\begin{proof}
Let $\xymatrix@M=4pt@H=4pt@C=15pt{ X_{-1}\times\Dl & X
\ar[r]^-{f}\ar[l]_-{\epsilon} &  Y }$ be a diagram in
$\simp\mc{D}$. It follows from $\ref{HomotCuadrCylTilde}$ that
there exists $R:\widetilde{Cyl}(X) \rightarrow
\widetilde{Cyl}(f,\epsilon)$ such that $R\comp I=j_Y\comp f$ and
$R\comp J= j_{X_{-1}} \comp \epsilon$.\\
Therefore $H=\beta\comp R:\widetilde{Cyl}(X)\rightarrow
Cyl(f,\epsilon)$ is such that $H \comp I= i_Y\comp f$ and $H\comp
J= i_{X_{-1}}\comp\epsilon$.
\end{proof}

%
%

\chapter{Simplicial Descent Categories}

The notion of (co)simplicial descent category is widely based in
the one of ``(co)homological descent category'', introduced in
\cite{GN}. In loc. cit. the basic objects are diagrams of
``cubical'' nature instead of simplicial objects in a f{i}xed
category.\\

\section{Def{i}nition}\label{descensosimplicial}

\begin{def{i}}
Consider a category $\mc{D}$ and a class of morphisms $\mrm{E}$ in
$\mc{D}$. We will denote by $Ho\mc{D}$\index{Symbols}{$Ho\mc{D}$}
the localization of $\mc{D}$ with respect to $\mrm{E}$, and by
$\gamma:\mc{D}\rightarrow Ho\mc{D}$\index{Symbols}{$Ho\mc{D}$} the
canonical functor. The class $\mrm{E}$ is
saturated\index{Index}{saturated class} if
$$\mbox{a morphism } f\mbox{ is in }\mrm{E} \;\; \Longleftrightarrow\;\; \gamma(f)\mbox{ is an isomorphism in }Ho\mc{D} \ .$$
Equivalently, $\mrm{E}$ is saturated if
$\mrm{E}=\gamma^{-1}\{\mbox{ isomorpshisms of }Ho\mc{D}\,\}$.
\end{def{i}}
\begin{obs}\mbox{}\\
\textbf{i)} If $\mrm{E}$ is saturated, the ``2 out of 3'' property
holds for $\mrm{E}$\index{Index}{two-of-three property}. That is,
if two of the morphisms $f$, $g$ or $gf$ are in $\mrm{E}$ then so is the third.\\[0.3cm]
\textbf{ii)} An enough (and necessary) condition for $\mrm{E}$
being saturated is that the morphisms in $\mrm{E}$ are just those
that are mapped by a certain functor into isomorphisms. In other words:\\
If $F:\mc{D}\rightarrow\mc{C}$ is a functor and $\mrm{E}=\{f \;|\;
F(f)\mbox{ is an isomorphism }\}$ then $\mrm{E}$ is saturated.
\end{obs}

The following lemma will be needed later, whose (trivial) proof is
left to the reader.

\begin{lema}\label{DcerradporRetracto}
If $\mrm{E}$ is a saturated class of morphisms in a category
$\mc{D}$ with f{i}nal object $1$, then the class of acyclic
objects of $\mc{D}$ $($with respect to $\mrm{E})$
$$\mc{A} =\{\mbox{objects }A\mbox{ of }\mc{D}\ |\  A\rightarrow 1\mbox{ is an equivalence}\} $$
is closed under retracts.\\
More specif{i}cally, if $A\stackrel{r}{\rightarrow}
B\stackrel{p}{\rightarrow}A$ is such that $p\comp r=Id_A$ and
$B\rightarrow 1$ is an equivalence, then $A\rightarrow 1$ is also
an equivalence.
\end{lema}

\begin{proof}
By the 2 out of 3 property, it is enough to see that $r$ is an equivalence.\\
Let $\xi:B\rightarrow 1$ be the trivial morphism. Then $\xi\comp
r\comp p=\xi:B\rightarrow 1$, since $1$ is f{i}nal object. Therefore $\rho=r\comp p:B\rightarrow B$ is in $\mrm{E}$ .\\
Thus, the equalities $p\comp r=Id_A$ and $r\comp (p\comp
\rho^{-1})=Id_B$ hold in $Ho\mc{D}$. So $p$ is a right inverse of
$r$ and $p\rho^{-1}$ is a left inverse of $r$. Consequently $r$ is
an isomorphism in $Ho\mc{D}$ with $p\rho^{-1}=p$ as inverse. Since
$\mrm{E}$ is saturated, we deduce that $r\in\mrm{E}$.
\end{proof}

Before going into details with the notion of simplicial descent
categories, we introduce the following notations.

\begin{num}\label{FuntorMonoidalCasiEstricto}
Let $\cont$ and $\mc{D}$ be categories with f{i}nite coproducts
(in particular initial object $0$).\\
Note that every functor $\psi:\cont\rightarrow\mc{D}$ is (lax)
monoidal with respect to the coproduct. The K{\"u}nneth morphism is
the one given by the universal property of the coproduct
$$\sigma_{X,Y}: \psi(X)\sqcup \psi(Y)\rightarrow \psi(X\sqcup Y)\ \ \forall\ X,Y\in\cont\ \ ; \ \sigma_0: =0\rightarrow \psi(0)\ . $$
The morphism $\sigma_{X,Y}$ is the unique morphism such that the
diagram
$$\xymatrix@H=4pt@M=4pt@R=7pt{
 \psi(X)\ar[rd]_{i_{\psi(X)}} \ar@/^1pc/[rrd]^-{\psi(i_X)} & & \\
                                                          & \psi(X)\sqcup\psi(Y)\ar[r]^{\sigma_{X,Y}} & \psi(X\sqcup Y) \\
 \psi(Y)\ar[ru]^{i_{\psi(Y)}} \ar@/_1pc/[rru]_-{\psi(i_Y)} & &}$$
commutes. We will denote this natural transformation by
$\sigma_{\psi}$\index{Symbols}{$\sigma_\psi$},
or just $\sigma$ if $\psi$ is understood.\\
If $\mrm{E}$ is a class of morphisms of $\mc{D}$, the functor
$\psi$ is said to be \textit{quasi-strict
monoidal}\index{Index}{quasi-strict monoidal} (with respect to
$\mrm{E}$) if $\sigma_{X,Y}$ and $\sigma_0$ belongs to
$\mrm{E}$ for every objects $X$ and $Y$ in $\cont$.\\
Dually, if $\cont$ and $\mc{D}$ have f{i}nite products then every
functor $\psi:\cont\rightarrow \mc{D}$ is (lax) comonoidal with
respect to the product. This time the K{\"u}nneth morphism is the
canonical morphism $\sigma_\psi: \psi(X\times Y)\rightarrow
\psi(X)\times\psi(Y)$ given by the universal property of the
product.
\end{num}

\begin{num}\label{Def{i}Simpsimpl}
Assume that a functor $\mbf{s}:\simp\mc{D}\rightarrow \mc{D}$ is
given. Under the notations introduced in \ref{simpFuntor}, the
image under\index{Symbols}{$\simp\mbf{s}$}
$$\simp\mbf{s}:\simp\simp\mc{D}\rightarrow \simp\mc{D} $$
of a bisimplicial object $T$ in $\mc{D}$ is the simplicial object
$$ (\simp\mbf{s}(T))_n = \mbf{s}(T_{n,\cdot}) = \mbf{s}(m\rightarrow T_{n,m})\ .$$
\end{num}
\begin{def{i}}[Simplicial descent category]\label{catdescdebil}\mbox{}\\
A (\textit{simplicial}) \textit{descent category}
\index{Index}{category!simplicial descent} consists of the
data $(\mc{D},\mrm{E},\mbf{s},\mu,\lambda)$ where:\\[0.2cm]
$\mathbf{(SDC\; 1)}$ $\mc{D}$ is a category with f{i}nite
coproducts
(in particular with initial object $0$) and with f{i}nal object $1$.\\[0.2cm]
$\mathbf{(SDC\; 2)}$ $\mrm{E}$ is a saturated class of morphisms
in $\mc{D}$, stable by coproducts (that is
$\mrm{E}\sqcup\mrm{E}\subseteq \mrm{E}$). The morphisms in $\mrm{E}$ will be called \textit{equivalences}.\\[0.2cm]
$\mathbf{(SDC\; 3)} $\textbf{Additivity:} The \textit{simple}
functor $\mbf{s}:\simp\mc{D}\rightarrow \mc{D}$ commutes with
coproducts up to equivalence. In other words, the canonical
morphism
$\mbf{s}X \sqcup \mbf{s}Y\rightarrow \mbf{s}(X\sqcup Y)$ is in $\mrm{E}$ for all $X$, $Y$ in $\simp\mc{D}$.\\[0.2cm]
$\mathbf{(SDC\; 4)}$ \textbf{Factorization:} Let
$\mrm{D}:\simp\simp\mc{D}\rightarrow\simp\mc{D}$ be the diagonal
functor. Consider the functors $\mbf{s}(\simp\mathbf{s}),
\mbf{s}\mrm{D}:\simp\simp\mc{D}\rightarrow \mc{D}$. Then $\mu$ is
a natural transformation $\mu: \mbf{s}\mrm{D}\rightarrow
\mbf{s}(\simp\mathbf{s})$ such that $\mu_T\in \mrm{E}$ for all $T\in\simp\simp\mc{D}$.\\[0.2cm]
$\mathbf{(SDC\; 5)}$ \textbf{Normalization:}
$\lambda:\mathbf{s}(-\times \Dl)\rightarrow Id_{\mc{D}}$ is a
natural transformation compatible with
$\mu$ and such that $\lambda_X\in \mrm{E}$ for all $X\in\mc{D}$.\\[0.2cm]
$\mathbf{(SDC\; 6)}$ \textbf{Exactness:} If $f:X\rightarrow Y$ is
a morphism in $\simp\mc{D}$ with $f_n \in \mrm{E}$ $\forall n$
then $\mbf{s}(f)\in\mrm{E}$.\\[0.2cm]
$\mathbf{(SDC\; 7)}$ \textbf{Acyclicity:}
If $f:X\rightarrow Y$ is a morphism in $\simp\mc{D}$, then
$\mbf{s}f\in\mrm{E}$ if and only if the simple of its simplicial
cone is acyclic, that means that
$\mbf{s}(Cf)\rightarrow 1$ is an equivalence.\\[0.2cm]
$\mathbf{(SDC\; 8)}$ \textbf{Symmetry:} $\mbf {s}\Upsilon
f\in\mrm{E}$ if (and only if) $\mbf{s}f\in\mrm{E}$, where
$\Upsilon$ is the functor
$\Upsilon:\simp\mc{D}\rightarrow\simp\mc{D}$ that reverses the
order of the face and degeneracy maps in a simplicial object in
$\mc{D}$.
\end{def{i}}

\begin{num}[Compatibility between $\lambda$ and $\mu$]
Given $X\in\simp\mc{D}$, we have that
$$\begin{array}{l}\mbf{s}\comp\simp\mbf{s}(X\times\Dl)= \mbf{s}(n\rightarrow \mbf{s}(m\rightarrow X_n ))= \mbf{s}(n\rightarrow \mbf{s}(X_n \times\Dl)) \\
                  \mbf{s}\comp\simp\mbf{s}(\Dl\times X)= \mbf{s}(n\rightarrow \mbf{s}(m\rightarrow X_m))= \mbf{s}(\mbf{s}(X)\times\Dl)\ .
                  \end{array}$$
In the f{i}rst case, the morphisms
$\lambda_{X_n}:\mbf{s}(X_n\times\Dl)\rightarrow X_n$ give rise to
a morphism of simplicial objects.

The compatibility condition between $\lambda$ and $\mu$ means that
the following compositions must be equal to the identity in
$\mc{D}$:
\begin{equation}\label{compatibLambdaMu}\xymatrix@M=4pt@H=4pt@R=10pt{ \mbf{s}X \ar[r]^-{\mu_{\Dl\times X}} &  \mbf{s}((\mbf{s}X)\times\Dl)\ar[r]^-{\lambda_{\mbf{s}X}} & \mbf{s}X\\
                                                                      \mbf{s}X \ar[r]^-{\mu_{X\times \Dl}} &  \mbf{s}(n\rightarrow \mbf{s}(X_n \times\Dl))\ar[r]^-{\mbf{s}(\lambda_{X_n})} & \mbf{s}X\ .}\end{equation}
\end{num}

\begin{obs}[Comments on the symmetry axiom]\mbox{}\\
We will use later the following property of the image under the
simple functor of the simplicial cylinder associated with a
morphism $f:X\rightarrow Y$ and with an augmentation
$\epsilon:X\rightarrow X_{-1}\times\Dl$: \\
($\ast$) The simple of $f$ is an equivalence when the simple of
the
canonical inclusion $i_{X_{-1}}:X_{-1}\times\Dl\rightarrow Cyl(f,\epsilon)$ is so.\\
The converse property will be also needed to prove the ``transfer
lemma'' \ref{FDfuerte}, at least under some extra hypothesis.\\
The symmetry axiom is imposed in order to have the converse
statement of ($\ast$)
(see section \ref{SeccionCritAcCyl}).\\
However, other possibility is to impose the axiom: $sf\in\mrm{E}$
if and only if $\mbf{s}(i_{X_{-1}})\in\mrm{E}$, and
remove the symmetry axiom from the notion of simplicial descent category (in this case (SDC 7) holds setting $X_{-1}=1$).\\
We decide to do it in this way because the aim of this work is
just to establish a set of axioms ensuring the desired properties,
and we would like these axioms to be ``less restrictive as possible''.\\
An alternative to (SDC 8) is the existence of an isomorphism of
functors between $\mbf{s}$ and
$\mbf{s}\Upsilon:\simp\mc{D}\rightarrow \mc{D}$. But even this
property holds in many of our examples, it is not true for
$\mc{D}=Set$ and the diagonal $\mrm{D}:\simp\simp\mc{D}\rightarrow
\simp\mc{D}$ as simple functor.
\end{obs}

\begin{obs}
We consider in $\mc{D}$, $\simp\mc{D}$ and $\simp\simp\mc{D}$ the
trivial monoidal structures coming from the coproduct. Then we
have automatically that $\mbf{s}$ and $\simp\mbf{s}$ are (lax)
monoidal functors, with
$\sigma=\sigma_{\mbf{s}}$ and $\sigma_{\simp\mbf{s}}$ as respective K{\"u}nneth morphisms (see \ref{FuntorMonoidalCasiEstricto}).\\
In addition the natural transformations
$\lambda$ and $\mu$ are also monoidal.\\
That is to say, these transformations are compatible with $\sigma$
in the following sense.
Given objects $X$, $Y$ in $\mc{D}$, the diagram
\begin{equation}\label{ComptLambdaSigma}\xymatrix@H=4pt@M=4pt{ \mbf{s}(X\times\Dl) \sqcup \mbf{s} (Y\times\Dl) \ar[r]^-{\sigma}\ar[rd]_{\lambda_X \sqcup \lambda_Y} & \mbf{s}((X\sqcup Y)\times\Dl)\ar[d]^{\lambda_{X\sqcup Y}}\\
                                                                                                                             & X\sqcup Y}\end{equation}
commutes. On the other hand, let $Z$ and $T$ be bisimplicial
objects in $\mc{D}$. Then we have the following commutative
diagram
$$\xymatrix@H=4pt@M=4pt{ \mbf{s}\mrm{D}Z \sqcup \mbf{s}\mrm{D}T \ar[r]^{\sigma}\ar[d]_{\mu_Z \sqcup \mu_T}                     & \mbf{s}\mrm{D}(Z\sqcup T)\ar[d]_{\mu_{Z\sqcup T}}\\
                         \mbf{s}(\simp\mbf{s}Z) \sqcup \mbf{s}(\simp\mbf{s}T)\ar[r]^-{\sigma \simp\sigma} & \mbf{s}(\simp\mbf{s}(Z\sqcup Y)) .}$$
\end{obs}

\begin{obs}\label{obsfactorizacion}{\textbf{Factorization:}}\mbox{}\\
Recall the functor $\Gamma : \simp\simp
\mc{D}\rightarrow\simp\simp \mc{D}$ that swaps the indexes in a
bisimplicial object in $\mc{D}$. In the factorization axiom we may
also consider $\mbf{s}(\simp\mathbf{s})\Gamma
:\simp\simp\mc{D}\rightarrow \mc{D}$.\\
Assuming (SDC 4), since $\mrm{D}\Gamma=\mrm{D}$ we deduce the
existence of the natural transformation $\mu':\mbf{s}\mrm{D}
\rightarrow \mbf{s}(\simp\mathbf{s})\Gamma$ given by
$\mu'(Z)=\mu(\Gamma Z)$ and such that $\mu'(Z)\in \mrm{E}$,
$\forall Z\in\simp\simp\mc{D}$. Then
$$ \xymatrix@M=4pt@H=4pt{\mbf{s}(\simp\mathbf{s})Z & \mrm{D}Z \ar[r]^-{\mu'(Z)}\ar[l]_-{\mu(Z)} & \mbf{s}(\simp\mathbf{s})(\Gamma Z)  }\ .$$
\end{obs}

\begin{prop}\label{AxiomaExactitudPrima}
The axiom $\mathbf{(SDC\; 6)}$ in the notion of simplicial descent
category can be replaced by the following
alternative axiom\\
$\mathbf{(SDC\; 6)'}$ If $X$ is an object in $\simp\mc{D}$ such
that $X_n\rightarrow 1$ is an equivalence for every $n$, then
$\mbf{s}X\rightarrow 1$ is also in $\mrm{E}$.
\end{prop}

\begin{proof}
Assume that $\mc{D}$ satisf{i}es (SDC 6$)'$ instead of
(SDC 6), together with the remaining axioms of simplicial descent category.\\
Let $f:X\rightarrow Y$ be a morphism in $\simp\mc{D}$
with $f_n\in\mrm{E}$ for all $n$.\\
Then, for a f{i}xed $n\geq 0$ we have that
$\mbf{s}(f_n\times\Dl):\mbf{s}(X\times\Dl)\rightarrow\mbf{s}(Y\times\Dl)$
is an equivalence, since
$\lambda_Y\comp\mbf{s}(f_n\times\Dl)=f_n\comp\lambda_X$ and
the 2 out of 3 property holds for $\mrm{E}$.\\
Hence, it follows from the acyclicity axiom that $\mbf{s}C(f_n\times\Dl)\rightarrow 1$ is in $\mrm{E}$, for every $n$.\\
Consider now the bisimplicial object $T\in\simp\simp\mc{D}$
def{i}ned by
$$ T_{n,m} = C(f_n\times\Dl)_m = Y_n \sqcup \coprod^m X_n \sqcup 1 \ .$$
Equivalently, $T$ is the image under $Cyl^{(2)}$ of
$\xymatrix{1\times\Dl\times\Dl & X\times\Dl \ar[r]^{f\times\Dl}
\ar[l] & Y\times\Dl}$, (see \ref{cylBisimpl}).
Therefore
$$\mbf{s}(\simp\mbf{s} (T)) =\mbf{s}(n\rightarrow \mbf{s}(m\rightarrow C(f_n\times\Dl)_m))= \mbf{s}(n\rightarrow \mbf{s}(C(f_n\times\Dl)))\ ,$$
and from (SDC 6)$'$ we deduce that $\mbf{s}\simp\mbf{s} (T)\rightarrow 1$ is an equivalence.\\
Moreover, by the factorization axiom we have that
$\mu_T:\mbf{s}\mrm{D}T\rightarrow \mbf{s}(\simp\mbf{s} (T))$ is in
$\mrm{E}$, and again the 2 out of 3 property implies that
the trivial morphism $\rho:\mbf{s}\mrm{D}T\rightarrow 1$ is also in $\mrm{E}$.\\
On the other hand,
$\mrm{D}T=\mrm{D}Cyl^{(2)}(\xymatrix{1\times\Dl\times\Dl &
X\times\Dl \ar[r]^{f\times\Dl} \ar[l] & Y\times\Dl})$ that agrees
with $\widetilde{C}f$, the image under $\widetilde{Cyl}$ of
$\xymatrix{1\times\Dl & X \ar[r]^{f} \ar[l] &
Y}$, because of proposition \ref{relacionCylCylTildeDiag}.\\
From \ref{CylRetractoCylTilde} it follows that $\mbf{s}Cf$ is a
retract of $\mbf{s}\widetilde{C}f=\mbf{s}\mrm{D}T$, that is an
acyclic object.
Then we conclude by lemma \ref{DcerradporRetracto} that
$\mbf{s}Cf$ is acyclic, and using the acyclicity axiom we get that
$\mbf{s}f\in\mrm{E}$ .
\end{proof}

The following properties are direct consequences of the axioms.

\begin{prop}\label{AxiomaConoPrima}
If $f:X\rightarrow Y$ is a morphism between simplicial objects in
a simplicial descent category $\mc{D}$, then
\begin{center}$\mbf{s}f\in\mrm{E}$ if and only if $\mbf{s}(C'f)\rightarrow 1$ is an equivalence, \end{center}
where $C'$ is the ``symmetric'' notion of simplicial cone, given
in \ref{def{i}ConoSimetrico}.
\end{prop}

\begin{proof}
It follows from the symmetry axiom that
$\mbf{s}f\in\mrm{E}$ if and only if $\mbf{s}\Upsilon f\in\mrm{E}$.\\
By the acyclicity axiom, this happens if and only if
$\mbf{s}C(\Upsilon f)\rightarrow 1$ is an equivalence. If
$\tau:C(\Upsilon f)\rightarrow 1\times\Dl$ is the trivial
morphism, then $\mbf{s}C(\Upsilon f)\rightarrow 1$ is an
equivalence if and only if $\mbf{s}(\tau):\mbf{s}C(\Upsilon
f)\rightarrow \mbf{s}(1\times\Dl)$ is so, because the morphism
$\mbf{s}(1\times\Dl)\rightarrow 1$ is in $\mrm{E}$ by the normalization axiom.\\
Again this condition is equivalent to the fact that
$\mbf{s}(\Upsilon\tau):\mbf{s}(\Upsilon C(\Upsilon f))\rightarrow
\mbf{s}(1\times\Dl)$ is an equivalence, since $\Upsilon
(1\times\Dl) = 1\times\Dl$.
F{i}nally, by def{i}nition $C'f = \Upsilon C\Upsilon f$, and the
statement follows from the acyclicity of the object
$\mbf{s}(1\times\Dl)$.
\end{proof}

\begin{prop}\label{exactitudCyl}
Let $\mc{D}$ be a simplicial descent category.\\
Consider a morphism $(\alpha,\beta,\gamma):D\rightarrow D'$ in
$\Omega(\mc{D})$,
$$\xymatrix@M=4pt@H=4pt{ X_{-1}\times\Dl \ar[d]_-{\alpha} & X \ar[r]^-{f} \ar[l]_-{\;\epsilon}\ar[d]_{\beta} &  Y \ar[d]^{\gamma}\\
                         X'_{-1}\times\Dl                 & X' \ar[r]^-{f'} \ar[l]_-{\;\epsilon'}             &  Y'   ,}$$
such that $\alpha_n$, $\beta_n$ and $\gamma_n$ are in $\mrm{E}$
for all $n$. Then the induced morphism
$\mbf{s}(Cyl(\alpha,\beta,\gamma)):\mbf{s}(Cyl(D))\rightarrow
\mbf{s}(Cyl(D'))$ is also in $\mrm{E}$.
\end{prop}

\begin{proof}
By def{i}nition
$Cyl(\alpha,\beta,\gamma)_n=\gamma_n\sqcup\beta_{n-1}\sqcup\cdots\sqcup\beta_0\sqcup\alpha$.
So it follows from (SDC 2) that
$Cyl(\alpha,\beta,\gamma)_n\in\mrm{E}$ $\forall n$, and from (SDC
6) that $\mbf{s}(Cyl(\alpha,\beta,\gamma))\in\mrm{E}$.
\end{proof}

\begin{cor}
If $$\xymatrix@M=4pt@H=4pt{  X \ar[r]^-{f} \ar[d]_{\beta} &  Y \ar[d]^{\gamma}\\
                                                                              X' \ar[r]^-{f'}        &  Y'   ,}$$
is a morphism in $Fl(\simp\mc{D})$ such that $\beta_n$ and
$\gamma_n$ are equivalences for all $n$, then
$\mbf{s}C(\beta,\gamma):\mbf{s}(Cf)\rightarrow \mbf{s}(Cf')$ is
also in $\mrm{E}$.
\end{cor}

\begin{proof} Just set $X_{-1}=1$ and $\alpha=Id$ in the last proposition.
\end{proof}

\begin{prop}\label{CategFuntores}
If $I$ is a small category and
$(\mc{D},\mrm{E}_{\mc{D}},\mbf{s}_{\mc{D}},\lambda_{\mc{D}},\rho_{\mc{D}})$
is a simplicial descent category then the category of functors from $I$ to $\mc{D}$, $I\mc{D}$, has a natural structure of simplicial descent category.\\
Given $X:\simp\rightarrow I\mc{D}$, the image under the simple
functor in $I\mc{D}$, $\mbf{s}_{I\mc{D}}$, of a simplicial object
$X$ in $I\mc{D}$ is def{i}ned as
$$(\mbf{s}_{I\mc{D}}(X))(i)=\mbf{s}_{\mc{D}}(n\rightarrow X_n(i))$$
and $\mrm{E}_{I\mc{D}}=\{f\mbox{ such that
}f(i)\in\mrm{E}_{\mc{D}}\; \forall i\in I\}$.
\end{prop}

\begin{proof}
Def{i}ne also $\lambda_{I\mc{D}}$ and $\rho_{I\mc{D}}$ through the
identif{i}cation
\begin{equation}\label{IdentifCatPeq}\simp(I\mc{D})\equiv I{\simp\mc{D}}\ .\end{equation}
Then, (SDC 1) is clear, since the coproduct in $I\mc{D}$ is
def{i}ned degreewise.
The verif{i}cation of  the axioms (SDC 3)$,\ldots , $(SDC 6) and
(SDC 8) is straightforward.
To see (SDC 2), let us check that $\mrm{E}_{I\mc{D}}$ is a saturated class.\\
Let $\gamma_{I{\mc{D}}}:I{\mc{D}}\rightarrow
I{\mc{D}}[\mrm{E}_{I\mc{D}}^{-1}]$ and
$\gamma_{\mc{D}}:\mc{D}\rightarrow \mc{D}[\mrm{E}_{\mc{D}}^{-1}]$
be the localizations of $I{\mc{D}}$ and $\mc{D}$ with respect to $\mrm{E}_{I\mc{D}}$ and $\mrm{E}_{\mc{D}}$ respectively.\\
Given an object $j$ in $I$, consider the ``evaluation'' functor $\pi_j:I{\mc{D}}\rightarrow\mc{D}$, given by $\pi_j(P)=P(j)$.\\
Then $\gamma_{\mc{D}}\pi_j(\mrm{E}_{I\mc{D}})\subseteq$
$\{$isomorphisms of $\mc{D}[\mrm{E}_{\mc{D}}^{-1}]\}$, and the
composition $\gamma_{\mc{D}}\pi_j$ gives rise to the following
commutative diagram of functors,
$$
\xymatrix{I{\mc{D}}\ar[r]^{\pi_j} \ar[d]_{\gamma_{I{\mc{D}}}} & \mc{D}\ar[d]_{\gamma_\mc{D}} \\
            I{\mc{D}}[\mrm{E}_{I\mc{D}}^{-1}] \ar[r]^-{\overline{\pi}_j} & \mc{D}[\mrm{E}_{\mc{D}}^{-1}] \ .}
$$
Therefore, if $f$ is a morphism in $I\mc{D}$ such that
$\gamma_{I{\mc{D}}}(f)$
is an isomorphism, it follows that $\gamma_{\mc{D}}(f(j))$ is so for every $j\in I$.\\
Hence, $f(j)\in\mrm{E}_\mc{D}$ $\forall j$ and by def{i}nition $f\in\mrm{E}_{I\mc{D}}$.\\
To check (SDC 7) it is enough to note that, if we denote by
$C_{I\mc{D}}:Fl(\simp I\mc{D})\rightarrow \simp I\mc{D}$ and
$C_{\mc{D}}:Fl(\simp \mc{D})\rightarrow \simp \mc{D}$ the respective cone functors, then (by def{i}nition of the coproduct in  $I\mc{D}$), it holds that $[C_{I{\mc{D}}}(f)](j)=C_{\mc{D}}(f(j))$.\\
F{i}nally, (SDC 8) follows from the equality
$(\mbf{s}_{I\mc{D}}(\Upsilon f))(i)=\mbf{s}_{\mc{D}}(\Upsilon
f(i))$, for each $i\in {I}$.
\end{proof}

\begin{cor}\label{SimpDcatDescGamma}
If $\mc{D}$ is a $($simplicial$)$ descent category then
$\simp\mc{D}$ is so, where the simple functor
${\mathbf{\widetilde{s}}}$ is def{i}ned as
$$[\mbf{\widetilde{s}}(Z)]_n=\mbf{s}(m\rightarrow Z_{m,n})\mbox{ for all }Z\in\simp\simp\mc{D}$$
and where the class of equivalences is
$$\mrm{E}_{\simp\mc{D}}=\{f\mbox{ such that }f_n\in\mrm{E}_{\mc{D}}\; \forall n\}\ .$$
\end{cor}

\begin{obs} Following the notations in
(\ref{Def{i}Simpsimpl}) and \ref{obsfactorizacion}, the functor
$\widetilde{s}:\simp\simp\mc{D}\rightarrow \simp\mc{D}$ is the
composition $\simp\mbf{s}\comp\Gamma:\simp\simp\mc{D}\rightarrow
\simp\mc{D}$. This follows from the identif{i}cation
(\ref{IdentifCatPeq}), that now is just $\Gamma$.\\
A natural question is if it is also possible to consider
$\simp\mbf{s}$ as a simple. However, this time the transformation
$\lambda$ should be a morphism in $\simp\mc{D}$ relating
$(\mbf{s}X)\times\Dl$ to $X$, and in general this transformation
$\lambda$ does not exist.
\end{obs}

\subsection*{Cosimplicial Descent Categories}
\begin{def{i}}
A cosimplicial descent category\index{Index}{cosimplicial descent
category} consists of the data
$(\mc{D},\mrm{E},\mbf{s},\mu,\lambda)$ where $\mc{D}$ is a
category, $\mbf{s}:\Dl\mc{D}\rightarrow\mc{D}$ is a functor,
$\mrm{E}$ is a class of morphisms in $\mc{D}$, and
$\mu:\mbf{s}\simp\mbf{s}\rightarrow \mbf{s}\mrm{D} $ and
$\lambda:Id_{\mc{D}}\rightarrow \mbf{s}(-\times\Dl)$ are natural
transformations, such that $\mc{D}^\comp$, the opposite category
of $\mc{D}$, together with
$(\mrm{E}^\comp,\mbf{s}^\comp,\mu^\comp,\lambda^\comp)$, induced by $(\mrm{E},\mbf{s},\mu,\lambda)$ in $\mc{D}^\comp$, is a simplicial descent category.\\
More specif{i}cally, a  cosimplicial descent category is the data
$(\mc{D},\mrm{E},\mbf{s},\mu,\lambda)$ where:\\
$\mathbf{(CDC\; 1)}$ $\mc{D}$ is a category with f{i}nite products and initial object $0$.\\
$\mathbf{(CDC\; 2)}$ $\mrm{E}$ is a saturated class of morphisms
in $\mc{D}$, stable by products. That is, given $f,g
\in \mrm{E}$ then $f\sqcap g \in \mrm{E}$.\\
$\mathbf{(CDC\; 3)}$\textbf{Additivity:}
$\mbf{s}:\Dl\mc{D}\rightarrow \mc{D}$ is a quasi-strict comonoidal
functor with
respect to the product.\\
$\mathbf{(CDC\; 4)}$ \textbf{Factorization:} If
$\mrm{D}:\Dl\Dl\mc{D}\rightarrow\Dl\mc{D}$ is the diagonal
functor, consider $\mbf{s}(\Dl\mathbf{s}),
\mbf{s}\mrm{D}:\Dl\Dl\mc{D}\rightarrow \mc{D}$. Then
$\mu:\mbf{s}(\Dl\mathbf{s})\rightarrow\mbf{s}\mrm{D} $ is a
natural transformation such that $\mu(Z)\in \mrm{E}$ for every
$Z\in\Dl\Dl\mc{D}$.\\
$\mathbf{(CDC\; 5)}$ \textbf{Normalization:}
$\lambda:Id_{\mc{D}}\rightarrow\mathbf{s}(-\times \Dl) $ is a natural transformation, compatible with $\mu$, such that for every $X\in\mc{D}$, $\lambda(X)\in\mrm{E}$.\\
$\mathbf{(CDC\; 6)}$ \textbf{Exactness:} Given $f:X\rightarrow Y$
in $\Dl\mc{D}$ such that $f_n \in \mrm{E}$ $\forall n$
then $\mbf{s}(f)\in\mrm{E}$.\\
$\mathbf{(CDC\; 7)}$ \textbf{Acyclicity:}
If $f:X\rightarrow Y$ is a morphism in $\Dl\mc{D}$, then
$\mbf{s}f\in\mrm{E}$ if and only if the simple of its cosimplicial
path object is acyclic. That is, if and only if
$0\rightarrow\mbf{s}({Path}f)$ is an equivalence.\\
$\mathbf{(CDC\; 8)}$ \textbf{Symmetry:} it holds that $\mbf
{s}\Upsilon f\in\mrm{E}$ if (and only if) $\mbf{s}f\in\mrm{E}$.
\end{def{i}}
%


\section{Cone and Cylinder objects in a simplicial descent category}\label{pdadesDescenso}

From now on, $\mc{D}$ will denote a simplicial descent category.
The simplicial cone and cylinder functors in the category
$\simp\mc{D}$ (section \ref{defCilindroSimp}) induce cone and
cylinder functors in $\mc{D}$ through the constant and simple functors.\\
In addition, since $\mc{D}$ is a simplicial descent category,
these functors satisfy the ``usual'' properties (as in the chain
complex
case or topological case).\\
Of course, dual properties to those contained in this section remain valid in the cosimplicial setting.\\
Again, we mean by $Ho\mc{D}$ the localized category of $\mc{D}$
with respect to $\mrm{E}$.

\begin{def{i}}Let $\mrm{R}:\mc{D}\rightarrow\mc{D}$\index{Symbols}{$\mrm{R}$} be the functor def{i}ned as $\mrm{R}=\mbf{s}\comp
(-\times\Dl)$.\\
The normalization axiom provides the natural transformation
$\lambda:\mrm{R}\rightarrow Id_{\mc{D}}$, such that
$\lambda_X:\mrm{R}X\rightarrow X\in\mrm{E}$ $\forall X$ in
$\mc{D}$.\\
Therefore, given a morphism $f$ in $\mc{D}$ it follows from the 2
out of 3 property and from the naturality of $\lambda$ that
$f\in\mrm{E}$ if and only if $\mrm{R}f\in\mrm{E}$.
\end{def{i}}
\begin{num} Following the notations introduced in \ref{Def{i}Omega}
and \ref{catSquareCont}, the functor
$-\times\Dl:\mc{D}\rightarrow\simp\mc{D}$ induces
$-\times\Dl:\square_1^\comp \mc{D} \rightarrow \Omega(\mc{D})$, as
well as $\mbf{s}:\simp\mc{D}\rightarrow \mc{D}$ induces
$\mbf{s}: Co\Omega(\mc{D}) \rightarrow \square_1 \mc{D}$.\\
Again, we identify $Fl(\mc{D})$ with the full subcategory of
$\square_1^\comp \mc{D}$ whose objects are the diagrams
$1\leftarrow X\rightarrow Y$.
\end{num}

\begin{def{i}}[cone and cylinder functors]\label{ConoyCilindro}\mbox{}\\
We def{i}ne the cylinder functor $cyl:\square_1^\comp
\mc{D}\rightarrow \mc{D}$\index{Index}{cylinder!in a simpl.
descent cat.}\index{Symbols}{$cyl$} as the composition
$$\xymatrix@M=4pt@H=4pt{\square_1^\comp\mc{D} \ar[r]^-{-\times\Dl} &
\Omega(\mc{D}) \ar[r]^-{Cyl} &\simp\mc{D} \ar[r]^-{\mbf{s}}&
\mc{D}}\ .$$
More specif{i}cally, given morphisms $f:A\rightarrow B$ and
$g:A\rightarrow C$ in $\mc{D}$, the cylinder of $(f,g)\in
\square_1^\comp \mc{D}$ is the image under the simple functor of
the simplicial object $Cyl(f\times\Dl,g\times\Dl)$.\\
If 1 is a f{i}nal object in $\mc{D}$, the cone functor
$c:Fl(\mc{D})\rightarrow\mc{D}$\index{Index}{cone!in a simpl.
descent cat.}\index{Symbols}{$c$} is def{i}ned in a similar way as
$c=\mbf{s}\comp C\comp( - \times\Dl)$. Hence, the cone of
$f:A\rightarrow B$ is the simple of
$C(f\times\Dl)=Cyl(f\times\Dl,A\times\Dl\rightarrow 1\times\Dl)$.
\end{def{i}}

\begin{num}\label{inclusionescylenD} We deduce from \ref{inclusionesCyl} that if $\xymatrix@M=4pt@H=4pt{X &\ar[l]_{g}  Y\ar[r]^f & Z}$
is in $\square_1^\comp \mc{D}$, then applying $cyl$ we obtain an
object in $\square_1 \mc{D}$, natural in $(f,g)$, consisting of
$$\xymatrix@M=4pt@H=4pt{\mrm{R}X \ar[r]^-{I_X} & cyl(f,g)& \mrm{R}Z\ar[l]_-{I_Z}}\ .$$
\end{num}

In the f{i}rst chapter we have developed other notions of
simplicial cylinder dif{f}erent from $Cyl$, that are
$\widetilde{Cyl}$
and $Cyl'$ (given respectively in \ref{defCilCubico} and \ref{defCilPrima}).\\
However, the def{i}nition of $cyl$ ``does not depend'' on the
choice between $Cyl$ and $\widetilde{Cyl}$. On the other hand, if
we take $Cyl'$ instead of $Cyl$ in the def{i}nition of $cyl$ (that
is, its conjugate with respect to $\Upsilon$), the result is the
same as if we swap the variables $f$ and $g$ in $cyl$.

\begin{prop}\label{RelacionCilyCilTildeenD}\mbox{}\\
\textrm{a)} Given an object $\xymatrix@M=4pt@H=4pt{C & A \ar[r]^f
\ar[l]_g & B}$ in $\square_1^{\circ} \mc{D}$, it holds that
$$cyl(f,g)=\mbf{s}(\widetilde{Cyl}(f\times\Dl,g\times\Dl))\ .$$
\textrm{b)} There exists a natural isomorphism in $(f,g)$
$$cyl(f,g)\simeq cyl'(g,f)\ ,$$
where $cyl'(g,f)=\mbf{s}(Cyl'(g\times\Dl,f\times\Dl))$. Moreover
this isomorphism commutes with the canonical inclusions of
$\mrm{R}B$ and $\mrm{R}C$ into the respective cylinders.
\end{prop}

\begin{proof}
Part a) is a direct consequence of \ref{relacionCylCylTildeCtes},
whereas b) follows from \ref{relacionCylCylPrima}.
\end{proof}

Next we introduce some properties of the functors $cyl$ and $c$
that can be deduced trivially from the def{i}nitions.

\begin{prop}\label{exactitudcyl}
Consider a morphism  in $\square_1^\comp \mc{D}$
$$\xymatrix@M=4pt@H=4pt{ X \ar[d]_{\alpha} & Y \ar[r]^f \ar[l]_g \ar[d]_{\beta} &  Z \ar[d]_{\gamma}\\
                         X' & Y' \ar[r]^{f'} \ar[l]_{g'} &  Z'  , }$$
such that $\alpha,\beta,\gamma$ are equivalences. Then the induced
morphism $cyl(f,g)\rightarrow cyl(f',g')$ is also an equivalence.
\end{prop}

\begin{proof}
The statement follows trivially from the def{i}nition of $cyl$ and
from \ref{exactitudCyl}.
\end{proof}

\begin{cor}\label{exactitudcono}
Consider a commutative square in $\mc{D}$
$$\xymatrix@M=4pt@H=4pt{ X \ar[d]_{\alpha} \ar[r]^-f & Y \ar[d]_{\beta} \\
                         X' \ar[r]^-{f'}             & Y'  \ , }$$
such that $\alpha$ and $\beta$ are equivalences. Then the induced
morphism $c(f)\rightarrow c(f')$ is also an equivalence.
\end{cor}

\begin{prop}[Additivity of the cone and cylinder functors]\label{aditividadcylyc}\mbox{}\\
\textsc{a)} The cylinder functor is ``additive up to
equivalence''. In other words, given $(f,g),(f',g')\in
\square_1^\comp \mc{D}$, then the natural morphism
$$\sigma_{cyl} :  cyl(f,g)\sqcup cyl(f',g')\rightarrow cyl(f\sqcup f', g\sqcup g')\ ,$$
def{i}ned as in \ref{FuntorMonoidalCasiEstricto}, is an equivalence.\\
\textsc{b)} The morphism
$$\sigma_{c}: c(f)\sqcup c(f')\rightarrow c(f\sqcup f')$$
is an equivalence for any $f,f'\in Fl(\mc{D})$ if and only if
$1\sqcup 1 \rightarrow 1$ is an equivalence.
\end{prop}

\begin{proof} The f{i}rst statement follows from the additivity of the simplicial cylinder functor
$Cyl$ together with the axiom (SDC 3).\\
Indeed, given $(f,g),(f',g')\in \square_1^\comp \mc{D}$, we deduce
from proposition \ref{CylCoproducto} that
$$\sigma_{Cyl}:Cyl(f\times\Dl,g\times\Dl)\sqcup Cyl(f'\times\Dl,g'\times\Dl)\rightarrow Cyl((f\sqcup f')\times\Dl, (g\sqcup g')\times\Dl)$$
is an isomorphism. Therefore the following morphism is also an isomorphism\\
$$\mbf{s}\sigma_{Cyl}:\mbf{s}(Cyl(f\times\Dl,g\times\Dl)\sqcup
Cyl(f'\times\Dl,g'\times\Dl))\rightarrow cyl(f\sqcup f', g\sqcup
g')\ .$$
On the other hand, we have that
$$\sigma_{\mbf{s}}:cyl(f,g)\sqcup cyl(f',g')\rightarrow \mbf{s}(Cyl(f\times\Dl,g\times\Dl)\sqcup Cyl(f'\times\Dl,g'\times\Dl))$$
is an equivalence, and we are done since $\sigma_{cyl}=\mbf{s}\sigma_{Cyl}\comp\sigma_{\mbf{s}}$.\\[0.1cm]
\indent Let us prove b). The morphism $\mrm{R}1\rightarrow 1$ is
in $\mrm{E}$ because of (SDC 5), hence $1\sqcup 1\rightarrow 1$ is
an equivalence if and only if
$Id\sqcup Id: \mrm{R}1\sqcup \mrm{R}1\rightarrow \mrm{R}1$ is so.\\
F{i}rst, if in b) we set $f=f'=0 :0\rightarrow 0$ then
$\mrm{R}1\sqcup \mrm{R}1\rightarrow \mrm{R}1$ is just
$\sigma_{c}: c(0)\sqcup c(0)\rightarrow c(0\sqcup 0)=c(0)$.\\
\indent To see the remaining implication, assume that $1\sqcup
1\rightarrow 1$ is an equivalence.  If $D$ is an object in
$\mc{D}$, denote by $\rho_D$ the trivial morphism $D\rightarrow 1$.\\
Given morphisms $f:A\rightarrow B$ and $f':A'\rightarrow B'$ in
$\mc{D}$, it follows from part a) that
\begin{equation}\label{AuxSigmacyl}
 \sigma_{cyl}: c(f)\sqcup c(f')=cyl(f,\rho_A)\sqcup cyl(f',\rho_B)\rightarrow cyl(f\sqcup f',\rho_A \sqcup \rho_{A'})\end{equation}
is an equivalence. Also, by proposition \ref{exactitudcyl} we have
that the commutative diagram
$$\xymatrix@M=4pt@H=4pt@C=30pt{ B\sqcup B' \ar[d]_{Id} & A\sqcup A' \ar[r]^-{\rho_A\sqcup \rho_{A'}} \ar[l]_-{f\sqcup f'} \ar[d]_{Id} & 1\sqcup 1 \ar[d]\\
                                B\sqcup B'             & A\sqcup A' \ar[r]^-{\rho_{A\sqcup A'}} \ar[l]_-{f\sqcup f'} &  1  , }$$
gives rise to an equivalence $cyl(f\sqcup f',\rho_A \sqcup
\rho_{A'})\rightarrow c(f\sqcup f')$ such that composed with
(\ref{AuxSigmacyl}) is just $\sigma_c$.
\end{proof}

\begin{prop}\label{cilaumentado}
Consider the following commutative diagram in $\mc{D}$
\begin{equation}\label{diagrConmAumentacion}\xymatrix@H=4pt@M=4pt{X \ar[r]^f\ar[d]_{g} & Y\ar[d]_\beta\\
                                                                  Z\ar[r]^\alpha       & T\ .}\end{equation}
There exists $\rho:cyl(f,g)\rightarrow \mrm{R}T$, natural in
$($\ref{diagrConmAumentacion}$)$, such that $\rho\,
I_X=\mrm{R}\alpha$ and $\rho\, I_Z=\mrm{R}\beta$. Visually
$$\xymatrix@M=4pt@H=4pt{ \mrm{R}Y \ar[d]_{\mrm{R}g}\ar[r]^{\mrm{R}f}         & \mrm{R}Z \ar[d]^{I_Z}\ar@/^1pc/[rdd]^{\mrm{R}\beta} &\\
                         \mrm{R}X \ar@/_1pc/[drr]_{\mrm{R}\alpha}  \ar[r]^{I_X}  &  cyl(f,g)\ar@{}[d]|{\sharp}  \ar@{}[r]|{\sharp}\ar[dr]^{\rho}   &\\
                                                                 &                                                 & \mrm{R}T\,.} $$
\end{prop}

\begin{proof}
Just consider the diagram (\ref{diagrConmAumentacion}) in
$\simp\mc{D}$ through $-\times\Dl$, and apply $\mbf{s}$ to the
morphism $H:Cyl(f\times\Dl,g\times\Dl)\rightarrow T\times\Dl$
given in proposition \ref{Cylaumentado}.
\end{proof}

\begin{prop}\label{axiomaconoenD}
Given $f:X\rightarrow Y$ in $\mc{D}$ then $f\in \mrm{E}$ if and
only if its cone is acyclic, that is, if and only if
$c(f)\rightarrow 1$ is in $\mrm{E}$.
\end{prop}

\begin{proof}
It is enough to apply (SDC 7) to $f\times\Dl$, since $f\in\mrm{E}$
if and only if $\mrm{R}f\in\mrm{E}$.
\end{proof}
\begin{cor}\label{propRetracto}
The class $\mrm{E}$ is closed under retracts. In other words, if
$g:X'\rightarrow Y'$ is a morphism in $\mrm{E}$ and there exists a
commutative diagram in $\mc{D}$
\begin{equation}\label{retracto}\xymatrix@M=4pt@H=4pt{ X \ar[r]^r\ar[d]_{f} & X'\ar[r]^p\ar[d]_{g} & X\ar[d]_{f} \\
                         Y\ar[r]^{r'}         & Y'\ar[r]^{p'}& Y}\end{equation}
with $pr=Id_X$, $p'r'=Id_Y$ $($that is, $f$ is a retract of $g)$,
then $f\in\mrm{E}$.
\end{cor}

\begin{proof}
The image of (\ref{retracto}) under the cone functor is
$$\xymatrix@M=4pt@H=4pt{c(f)\ar[r]^{R}& c(g)\ar[r]^{P} & c(f)}$$
with $PR=Id_{c(f)}$ in $\mc{D}$. Since $g\in\mrm{E}$, then
$\xi:c(g)\rightarrow 1$ is an equivalence. Hence $c(f)$ is a
retract of an acyclic object, and from lemma
\ref{DcerradporRetracto} we deduce that $c(f)\rightarrow 1$ is in
$\mrm{E}$, so $f\in\mrm{E}$.
\end{proof}
%

%
%
%
%
%
%

\section{Factorization property of the cylinder functor}\label{SeccionFactCyl}

This section is devoted to the study of a relevant property of
``factorization'' satisf{i}ed by the functor $cyl$ in a simplicial
descent category $\mc{D}$. This property will be very useful in
the following section, in fact it is a key point in the
development of the relationship between the notions of simplicial
descent category and triangulated category.

\begin{num}
Assume given the following commutative diagram in $\mc{D}$
\begin{equation}\label{diagrConmTresporTres}\xymatrix@M=4pt@H=4pt{Z'  & X'\ar[l]_{g'}\ar[r]^{f'}   & Y' \\
                        Z  \ar[u]^{\alpha} \ar[d]_{\alpha'}& X  \ar[l]_{g}\ar[r]^{f} \ar[u]^{\beta} \ar[d]_{\beta'} & Y \ar[u]^{\gamma} \ar[d]_{\gamma'} \\
                        Z'' & X''\ar[l]_{g''}\ar[r]^{f''} & Y''.}\end{equation}
Applying the functor $cyl$ by rows and columns we obtain
$$\xymatrix@M=4pt@H=4pt@R=8pt{cyl(f',g')& cyl(f,g)\ar[l]_{\delta}\ar[r]^{\delta'}& cyl(f'',g'')\\
                        cyl(\alpha',\alpha)& cyl(\beta',\beta)\ar[l]_{\widehat{g}}\ar[r]^{\widehat{f}}& cyl(\gamma',\gamma).}$$
Denote by $\psi:cyl(\mrm{R}\gamma',\mrm{R}\gamma)\rightarrow
cyl(\delta',\delta)$ and
$\psi':cyl(\mrm{R}f'',\mrm{R}g'')\rightarrow
cyl(\widehat{f},\widehat{g})$ the respective morphisms obtained by
applying $cyl$ to the morphisms in $\square_1^{\circ}(\mc{D})$:
$$\xymatrix@M=4pt@H=4pt@C=17pt{\mrm{R}Y'\ar[d]_{I}& \mrm{R}Y\ar[l]_{\mrm{R}\gamma}\ar[r]^{\mrm{R}\gamma'}\ar[d]_{I} & \mrm{R}Y''\ar[d]_{I}\ar@{}[rd]|{;} & \mrm{R}Z''\ar[d]_{I}& \mrm{R}X''\ar[l]_{\mrm{R}g''}\ar[r]^{\mrm{R}f''}\ar[d]_{I} & \mrm{R}Y''\ar[d]_{I}\\
                         cyl(f',g')& cyl(f,g)\ar[l]_{\delta}\ar[r]^{\delta'}& cyl(f'',g'')                                          & cyl(\alpha',\alpha)& cyl(\beta',\beta)\ar[l]_{\widehat{g}}\ar[r]^{\widehat{f}}& cyl(\gamma',\gamma)\ ,}$$
where $I$ means the corresponding canonical inclusion.\\
In the same way, denote by
$\widehat{\lambda}:cyl(\mrm{R}\gamma',\mrm{R}\gamma)\rightarrow
cyl(\gamma',\gamma)$ and $\widetilde{\lambda}:
cyl(\mrm{R}f'',\mrm{R}g'')\rightarrow cyl(f'',g'')$ the
equivalences obtained from
$$\xymatrix@M=4pt@H=4pt@C=30pt{\mrm{R}Y'\ar[d]_{\lambda_{Y'}}& \mrm{R}Y\ar[l]_{\mrm{R}\gamma}\ar[r]^{\mrm{R}\gamma'}\ar[d]_{\lambda_{Y}} & \mrm{R}Y''\ar[d]_{\lambda_{Y''}}\ar@{}[rd]|{;} & \mrm{R}Z''\ar[d]_{\lambda_{Z''}}& \mrm{R}X''\ar[l]_{\mrm{R}g''}\ar[r]^{\mrm{R}f''}\ar[d]_{\lambda_{X''}} & \mrm{R}Y''\ar[d]_{\lambda_{Y''}}\\
                                   Y'& Y\ar[l]_{\gamma}\ar[r]^{\gamma'}& Y''                                                                                                              &Z''& X''\ar[l]_{g''}\ar[r]^{f''}& Y''.}$$
\end{num}

\begin{prop}\label{factorizacionCubica} Under the above notations, the cylinder objects $cyl(\delta',\delta)$
and $cyl(\widehat{f},\widehat{g})$ are naturally isomorphic in
$Ho\mc{D}$.\\
More concretely, let $\widetilde{T}$ be the simplicial object in
$\mc{D}$ obtained by applying $\widetilde{Cyl}$ to the diagram
${Cyl}(\alpha',\alpha)\leftarrow {Cyl}(\beta',\beta)\rightarrow
{Cyl}(\gamma',\gamma)$.\\
Then there exists isomorphisms
$\Psi:\mbf{s}(\widetilde{T})\rightarrow cyl(\delta',\delta)$,
$\Phi:\mbf{s}(\widetilde{T})\rightarrow
cyl(\widehat{f},\widehat{g})$ in $Ho\mc{D}$, natural in
$(\ref{diagrConmTresporTres})$, such that the diagram
\begin{equation}\label{diagDosCubos}\xymatrix@R=50pt@C=55pt@H=4pt@M=4pt@!0{                                                                   & \mrm{R}^2 Y'' \ar'[d][dd]_{I}\ar[ld]_{\mrm{R}I}\ar[rr]^{\lambda_{\mrm{R}Y''}}&                                              &  \mrm{R}Y''\ar'[d][dd]_{I}\ar[ld]_{I}&                                                          & \mrm{R}^2 Y'' \ar[dd]^{\mrm{R}I}\ar[ld]_{I} \ar[ll]_{\lambda_{\mrm{R}Y''}}  \\
                                             \mrm{R}cyl(f'',g'')\ar[dd]_{I}\ar[rr]^(0.6){\lambda}                &                                                                          & cyl(f'',g'')\ar[dd]                          &                                      &  cyl(\mrm{R}f'',\mrm{R}g'') \ar[dd]\ar[ll]_(0.4){\widetilde{\lambda}}\ar[dd]^(0.6){\psi'}  &                                               \\
                                                                                                            & cyl(\mrm{R}\gamma',\mrm{R}\gamma)\ar[ld]_{\psi}\ar'[r][rr]^(0.2){\widehat{\lambda}}&                                              & cyl(\gamma',\gamma)\ar[ld]           &                                                                     & \mrm{R}cyl(\gamma',\gamma)\, ,\ar[ld]^{I}\ar'[l][ll]_(0.2){\lambda}\\
                                              cyl(\delta',\delta)                                           &                                                                          & \mbf{s}\widetilde{T}  \ar[rr]^{\Phi}\ar[ll]_{\Psi}&                                      &  cyl(\widehat{f},\widehat{g})                            &}\end{equation}
commutes in $Ho\mc{D}$, and the same holds for $Z',Z''$ and $Y'$.
\end{prop}

\begin{proof} {F}irst of all, note that it is enough to prove
the commutativity in $Ho\mc{D}$ of
\begin{equation}\label{digrFactCil}\xymatrix@H=4pt@M=4pt{\mrm{R}cyl(f'',g'')\ar[d]_{I}\ar[r]^{\lambda}                      & cyl(f'',g'')\ar[d]                              & cyl(\mrm{R}f'',\mrm{R}g'') \ar[l]_{\widetilde{\lambda}}\ar[d]^{\psi'} \\
                        cyl(\delta',\delta)                                                     & \mbf{s}\widetilde{T}  \ar[r]^{\Phi}\ar[l]_{\Psi}& cyl(\widehat{f},\widehat{g})\\
                        cyl(\mrm{R}\gamma',\mrm{R}\gamma)\ar[u]_{\psi}\ar[r]^{\widehat{\lambda}}& cyl(\gamma',\gamma)\ar[u]                       & \mrm{R}cyl(\gamma',\gamma)\ar[u]^{I}\ar[l]_{\lambda}
                       }\end{equation}
since the remaining squares in diagram (\ref{diagDosCubos})
commutes because of the def{i}nitions of the arrows involved in
them, or because of the commutativity of the rest of the diagram,
together with the fact that the horizontal arrows are
isomorphisms in $Ho\mc{D}$.\\
Set $Cyl(h,t)=Cyl(h\times\Dl,t\times\Dl)$ if
$h,t$ are morphisms in $\mc{D}$.\\
Consider diagram  \ref{diagrConmTresporTres} in $\simp\mc{D}$
through the functor $-\times\Dl$ and denote by
$$\xymatrix@M=4pt@H=4pt@R=8pt{{Cyl}(f',g')& {Cyl}(f,g)\ar[l]_{\rho}\ar[r]^{\rho'}& {Cyl}(f'',g'')\\
                              {Cyl}(\alpha',\alpha)& {Cyl}(\beta',\beta)\ar[l]_{G}\ar[r]^{F}& {Cyl}(\gamma',\gamma)}$$
the result of applying ${Cyl}$ by rows and columns respectively.
We also follow the notations introduced in
\ref{NotacionesFactCyl} for $\phi$ and $\varphi$.\\
Let $\Theta:{Cyl}_{\simp\mc{D}}^{(1)}(\Dl\times F,\Dl\times
G)\rightarrow
{Cyl}_{\simp\mc{D}}^{(2)}(\rho'\times\Dl,\rho\times\Dl)$ be the
canonical isomorphism in $\simp\simp\mc{D}$ given in
\ref{factorizCylTilde}. It is such that the diagram
\begin{equation}\label{diagrRombo}\xymatrix@R=10pt@C=25pt{                                              &{Cyl}_{\simp\mc{D}}^{(1)}(\Dl\times F,\Dl\times G) \ar[dd]_{\Theta}&                                          \\
                                                          Cyl(f'',g'')\times\Dl\ar[dr]_{i}\ar[ur]^{\phi}&                                                                  &\Dl\times Cyl(\gamma',\gamma) \ar[lu]_{i}\ar[ld]^{\varphi} \\
                                                                                                        &{Cyl}_{\simp\mc{D}}^{(2)}(\rho'\times\Dl,\rho\times\Dl)           &}\end{equation}
commutes. We will compute the image of the above diagram under the
functors $\mbf{s}\comp\mrm{D}$, $\mbf{s}\comp\simp\mbf{s}$ and
$\mbf{s}\comp\simp\mbf{s}\Gamma$.\\
Set $T={Cyl}_{\simp\mc{D}}^{(1)}(\Dl\times F,\Dl\times G)$ and
$R={Cyl}_{\simp\mc{D}}^{(2)}(\rho'\times\Dl,\rho\times\Dl)$. We
deduce from \ref{relacionCylCylTildeDiag} that
$\mrm{D}T=\widetilde{Cyl}(F,G)=\widetilde{T}$ and
$\mrm{D}R=\widetilde{Cyl}(\rho',\rho)$.\\
Therefore, applying $\mbf{s}\comp\mrm{D}$ to (\ref{diagrRombo}) we
obtain the following commutative diagram in $\mc{D}$
\begin{equation}\label{diagrRomboBis}\xymatrix@R=10pt@C=25pt{                                                                     &\mbf{s}\widetilde{Cyl}( F,G) \ar[dd]_{\mbf{s}\mrm{D}\Theta}&                                          \\
                                                          cyl(f'',g'')\ar[dr]_{\mbf{s}\mrm{D}i}\ar[ur]^{\mbf{s}\mrm{D}\phi}&                                                                  &cyl(\gamma',\gamma) \ar[lu]_{\mbf{s}\mrm{D}i}\ar[ld]^{\mbf{s}\mrm{D}\varphi} \ .\\
                                                                                               &\mbf{s}\widetilde{Cyl}(\rho',\rho)&}\end{equation}
where $\mbf{s}\mrm{D}\Theta$ is an isomorphism.\\
Let us compute the image under $\simp\mbf{s}$ of
$$\xymatrix@H=4pt@M=4pt{Cyl(f'',g'')\times\Dl\ar[r]^-{\phi} & {Cyl}_{\simp\mc{D}}^{(1)}(\Dl\times F,\Dl\times G) & \Dl\times Cyl(\gamma',\gamma) \ar[l]_-{i}\ .}$$
It follows from the def{i}nitions that
$$(\simp\mbf{s}T)_n=\mbf{s}(m\rightarrow (Cyl(\gamma',\gamma))_m \sqcup \ds\coprod^n (Cyl(\beta',\beta))_{m} \sqcup (Cyl(\alpha',\alpha))_{m})\ .$$
In addition
$$
\begin{array}{l}
 (\simp\mbf{s}(Cyl(f'',g'')\times\Dl))_n=\mbf{s}(m\rightarrow Y''\sqcup \ds\coprod^n X'' \sqcup Z'')=\mrm{R}(Y''\sqcup \ds\coprod^n X'' \sqcup Z'')\\
 (\simp\mbf{s}(\Dl\times Cyl(\gamma',\gamma))_n=\mbf{s}(m\rightarrow Cyl(\gamma',\gamma)_m )= cyl(\gamma',\gamma)\ .\end{array}
$$
By the universal property of the coproduct we have the morphisms
$$\begin{array}{l}\xymatrix@H=4pt@M=4pt{cyl(\gamma',\gamma)\sqcup\coprod^n cyl(\beta',\beta)\sqcup cyl(\alpha',\alpha) \ar[r]^-{\sigma_n}  & (\simp\mbf{s}T)_n} \\
                  \xymatrix@H=4pt@M=4pt{\mrm{R} Y''\sqcup\coprod^n\mrm{R} X''\sqcup\mrm{R}Z'' \ar[r]^-{\sigma_n}                           & \mrm{R}(Y'' \sqcup \coprod^n X'' \sqcup Z'')}
                  \end{array}$$
in such a way that the following diagram commutes
$$\xymatrix@H=4pt@M=4pt@C=30pt{
 Cyl(\mrm{R}f'',\mrm{R}g'')_n \ar[d]_{\sigma_n} \ar[r]^-{\widetilde{\psi}} & Cyl(\mbf{s}F,\mbf{s}G)_n \ar[d]_{\sigma_n}& \\
 (\simp\mbf{s}(Cyl(f'',g'')\times\Dl))_n \ar[r]^-{(\simp\mbf{s}\phi)_n}    & (\simp\mbf{s}T)_n                         & cyl(\gamma',\gamma)\ar[l]_{(\simp\mbf{s}i)_n}\ar[lu]_{i} \ ,}
$$
where both morphisms $\sigma_n$ are equivalences because of (SDC
3), and where $\mbf{s}F=\widehat{f}$, $\mbf{s}G=\widehat{g}$ and
$\mbf{s}\widetilde{\psi}=\psi'$. Then, applying $\mbf{s}$, the
following diagram commutes
$$\xymatrix@H=4pt@M=4pt@C=30pt{
 cyl(\mrm{R}f'',\mrm{R}g'') \ar[d]_{\mbf{s}(\{\sigma_n\})} \ar[r]^-{\psi'}             & cyl(\widehat{f},\widehat{g}) \ar[d]_{\mbf{s}(\{\sigma_n\})}& \\
 \mbf{s}(\simp\mbf{s}(Cyl(f'',g'')\times\Dl)) \ar[r]^-{\mbf{s}(\simp\mbf{s}\phi)}         & \mbf{s}(\simp\mbf{s}T)                                     & \mrm{R}cyl(\gamma',\gamma)\ar[l]_{\mbf{s}\simp\mbf{s}i}\ar[lu]_{I}\ ,}
$$
where $\mbf{s}(\{\sigma_n\})$ is an equivalence (by the
exactness axiom).\\
On the other hand, the natural transformation
$\mu:\mbf{s}\comp\mrm{D}\rightarrow \mbf{s}\comp\simp\mbf{s}$
gives rise to
$$\xymatrix@H=4pt@M=4pt@C=30pt{
 \mbf{s}(\simp\mbf{s}(Cyl(f'',g'')\times\Dl)) \ar[r]^-{\mbf{s}(\simp\mbf{s}\phi)} & \mbf{s}(\simp\mbf{s}T)                    & \mrm{R}cyl(\gamma',\gamma)\ar[l]_-{\mbf{s}\simp\mbf{s}i}\\
 cyl(f'',g'') \ar[u]_{\mu_{Cyl(f'',g'')\times\Dl}} \ar[r]^-{\mbf{s}\mrm{D}\phi}     & \mbf{s}\widetilde{Cyl}(F,G)\ar[u]_{\mu_T}  &  cyl(\gamma',\gamma)\ar[u]_{\mu_{\Dl\times Cyl(\gamma',\gamma)}}\ar[l]_-{\mbf{s}\mrm{D}i}\ . }$$
From the equations (\ref{compatibLambdaMu}) describing the
compatibility between $\lambda$ and $\mu$ we deduce that the
following diagram commutes in $Ho\mc{D}$
$$\xymatrix@H=4pt@M=4pt@C=30pt{
 \mbf{s}(\simp\mbf{s}(Cyl(f'',g'')\times\Dl)) \ar[r]^-{\mbf{s}(\simp\mbf{s}i)}\ar[d]_{\mbf{s}(\{\lambda_n\})} & \mbf{s}(\simp\mbf{s}T)& \mrm{R}cyl(\gamma',\gamma)\ar[l]_-{\mbf{s}\simp\mbf{s}i}\ar[d]^{\lambda_{cyl(\gamma',\gamma)}}\\
 cyl(f'',g'') \ar[r]^-{\mbf{s}\mrm{D}\phi}     & \mbf{s}\widetilde{Cyl}(F,G)\ar[u]_{\mu_T}                                            &  cyl(\gamma',\gamma)\ar[l]_-{\mbf{s}\mrm{D}i}\ . }$$
If we join the two resulting diagrams, we obtain
$$\xymatrix@H=4pt@M=4pt@C=30pt{
 cyl(\mrm{R}f'',\mrm{R}g'') \ar[d]_{\mbf{s}(\{\sigma_n\})} \ar[r]^-{\psi'}                                    & cyl(\widehat{f},\widehat{g})\ar[d]_{\mbf{s}(\{\sigma_n\})}& \mrm{R}cyl(\gamma',\gamma)\ar[l]_-{I}\ar[dd]^{\lambda_{cyl(\gamma',\gamma)}}\\
 \mbf{s}(\simp\mbf{s}(Cyl(f'',g'')\times\Dl))\ar[d]_{\mbf{s}(\{\lambda_n\})}                                  & \mbf{s}(\simp\mbf{s}T)                                    & \\
 cyl(f'',g'') \ar[r]^-{\mbf{s}\mrm{D}\phi}                                                                    & \mbf{s}\widetilde{Cyl}(F,G)\ar[u]_{\mu_T}                 & cyl(\gamma',\gamma)\ar[l]_-{\mbf{s}\mrm{D}i}\ , }$$
that is just the right side of (\ref{digrFactCil}) taking
$\Phi=\mbf{s}(\{\sigma_n\})^{-1}\comp\mu_T$, and noting that
$\mbf{s}(\{\lambda_n\comp\sigma_n\})=\widetilde{\lambda}$.\\
Indeed, it follows from the compatibility (\ref{ComptLambdaSigma})
between $\lambda$ and $\sigma$ that the composition
$$ \mrm{R}Y''\sqcup\coprod^n \mrm{R}X''\sqcup \mrm{R}Z'' \stackrel{\sigma_n}{\longrightarrow} \mrm{R}(Y''\sqcup\coprod^n X''\sqcup Z'') \stackrel{\lambda_n}{\longrightarrow} Y''\sqcup\coprod^n X''\sqcup Z''$$
is equal to
$\lambda_{Y''}\sqcup\coprod^n\lambda_{X''}\sqcup\lambda_{Z''}$.

\noindent It remains to see the existence of the left side of
(\ref{digrFactCil}).
In order to do that, we argue in a similar way to compute the
image under $\simp\mbf{s}\Gamma$ of
$$\xymatrix@H=4pt@M=4pt{Cyl(f'',g'')\times\Dl\ar[r]^-{i} & {Cyl}_{\simp\mc{D}}^{(2)}(\rho'\times\Dl,\rho\times\Dl) & \Dl\times Cyl(\gamma',\gamma) \ar[l]_-{\varphi}\ .}$$
By def{i}nition
$$\begin{array}{l}
(\simp\mbf{s}\Gamma R)_n=\mbf{s}(m\rightarrow Cyl(f'',g'')_m \sqcup \ds\coprod^n Cyl(f,g)_{m} \sqcup Cyl(f',g')_{m})\\
 (\simp\mbf{s}\Gamma( Cyl(f'',g'')\times\Dl))_n= (\simp\mbf{s}(\Dl\times Cyl(f'',g'')))_n=\mbf{s}(m\rightarrow  Cyl(f'',g'')_m)\\
 (\simp\mbf{s}\Gamma (\Dl\times Cyl(\gamma',\gamma)))_n=(\simp\mbf{s} Cyl(\gamma',\gamma)\times\Dl)_n=\mbf{s}(m\rightarrow Y''\sqcup\coprod^n Y\sqcup Y') \ .
 \end{array}$$
Again by the universal property of the coproduct we have the
following commutative diagram
$$\xymatrix@H=4pt@M=4pt@C=35pt{
cyl(f'',g'')  \ar[r]^-{(\simp\mbf{s}\Gamma i)_n}\ar[rd]_{i} & (\simp\mbf{s}\Gamma R)_n                           & (\simp\mbf{s}(Cyl(\gamma',\gamma)\times\Dl))_n \ar[l]_-{(\simp\mbf{s}\Gamma\varphi)_n} \\
                                                            & Cyl(\mbf{s}\rho',\mbf{s}\rho)_n  \ar[u]_{\sigma_n} & cyl(\gamma',\gamma)\ar[u]_{\sigma_n}\ar[l]_{\overline{\psi}_n}  \ ,}
$$
where $\mbf{s}\rho'=\delta'$, $\mbf{s}\rho=\delta$ and
$\mbf{s}\overline{\psi}=\psi$. Hence applying $\mbf{s}$ we obtain
$$\xymatrix@H=4pt@M=4pt@C=35pt{
 \mrm{R}cyl(f'',g'')  \ar[r]^-{\mbf{s}(\simp\mbf{s}\Gamma i)}\ar[rd]_{I}        & \mbf{s}(\simp\mbf{s}\Gamma R)                          & \mbf{s}(\simp\mbf{s}(Cyl(\gamma',\gamma)\times\Dl)) \ar[l]_-{\mbf{s}(\simp\mbf{s}\Gamma\varphi)} \\
                                                                                & cyl(\delta',\delta)  \ar[u]_{\mbf{s}(\{\sigma_n\})}    & cyl(\gamma',\gamma)\ar[u]_{\mbf{s}(\{\sigma_n\})}\ar[l]_{\psi}  \ .}
$$
The natural transformation $\mu\Gamma:\mbf{s}\mrm{D}\rightarrow
\mbf{s}\comp\simp\mbf{s}\comp\Gamma$ gives rise to
$$\xymatrix@H=4pt@M=4pt@C=35pt{
 cyl(f'',g'')\ar[d]^{\mu_{\Dl\times Cyl(f'',g'')}} \ar[r]^-{I} & \widetilde{Cyl}(\rho',\rho)    \ar[d]^{\mu_{\Gamma R}} & cyl(\gamma',\gamma)\ar[l]_{\mbf{s}\mrm{D}\varphi} \ar[d]^{\mu_{Cyl(\gamma',\gamma)\times\Dl}} \\
 \mrm{R}cyl(f'',g'')  \ar[r]^-{\mbf{s}(\simp\mbf{s}\Gamma i)}  & \mbf{s}(\simp\mbf{s}\Gamma R)                          & \mbf{s}(\simp\mbf{s}(Cyl(\gamma',\gamma)\times\Dl))\ar[l]_-{\mbf{s}(\simp\mbf{s}\Gamma\varphi)} \ ,}$$
where we can replace $\mu_{\Dl\times Cyl(f'',g'')}$ by
$\lambda_{cyl(f'',g'')}$ and $\mu_{Cyl(\gamma',\gamma)\times\Dl}$
by $\mbf{s}(\{\lambda_n\})$.\\
Putting all together we get
$$\xymatrix@H=4pt@M=4pt@C=35pt{
  cyl(f'',g'')                                      \ar[r]^-{I}   & \widetilde{Cyl}(\rho',\rho)    \ar[d]^{\mu_{\Gamma R}} & cyl(\gamma',\gamma)\ar[l]_{\mbf{s}\mrm{D}\varphi}  \\
                                                                  & \mbf{s}(\simp\mbf{s}\Gamma R)                          & \mbf{s}(\simp\mbf{s}(Cyl(\gamma',\gamma)\times\Dl)) \ar[l]_-{\mbf{s}(\simp\mbf{s}\Gamma\varphi)} \ar[u]_{\mbf{s}(\{\lambda_n\})}\\
  \mrm{R}cyl(f'',g'') \ar[uu]_{\lambda_{cyl(f'',g'')}} \ar[r]^{I} & cyl(\delta',\delta)\ar[u]_{\mbf{s}(\{\sigma_n\})}      & cyl(\gamma',\gamma)\ar[u]_{\mbf{s}(\{\sigma_n\})}\ar[l]_{\psi}  \ .}
$$
Again,
$\mbf{s}(\{\lambda_n\}\comp\mbf{s}(\{\sigma_n\})=\widehat{\lambda}$.
Then, adjoining (\ref{diagrRomboBis}), the result is
$$\xymatrix@H=4pt@M=4pt@C=35pt{
  cyl(f'',g'')  \ar@{-}[d]^{Id}      \ar[r]^-{\mbf{s}\mrm{D}\phi} & \widetilde{Cyl}(F,G)    \ar[d]^{\mbf{s}\mrm{D}\Theta}  & cyl(\gamma',\gamma)\ar[l]_{\mbf{s}\mrm{D}i}\ar@{-}[d]^{Id}  \\
  cyl(f'',g'')                                      \ar[r]^-{I}   & \widetilde{Cyl}(\rho',\rho)    \ar[d]^{\mu_{\Gamma R}} & cyl(\gamma',\gamma)\ar[l]_{\mbf{s}\mrm{D}\varphi}  \\
                                                                  & \mbf{s}(\simp\mbf{s}\Gamma R)                          & \\
  \mrm{R}cyl(f'',g'') \ar[uu]_{\lambda_{cyl(f'',g'')}} \ar[r]^{I} & cyl(\delta',\delta)\ar[u]_{\mbf{s}(\{\sigma_n\})}      & cyl(\gamma',\gamma)\ar[uu]_{\widehat{\lambda}}\ar[l]_{\psi}  }
$$
To finish just take
$\Psi=\mbf{s}(\{\sigma_n\})^{-1}\comp\mu_{\Gamma
R}\comp\mbf{s}\mrm{D}\Theta$.
\end{proof}

In order to deduce from this result an analogous property for the
cone functor, we need the following lemma.

\begin{lema}\label{relCyltildeconCil}
Let $\xymatrix@M=4pt@H=4pt{Z&X\ar[l]_g \ar[r]^f &Y}$ be a diagram
of simplicial objects in $\mc{D}$. There exists an isomorphism
$\Phi:\mbf{s}(\widetilde{Cyl}(f,g))\rightarrow
cyl(\mbf{s}f,\mbf{s}g)$ in $Ho\mc{D}$, natural in $f$ and $g$,
that fits into the following commutative diagram of $Ho\mc{D}$
\begin{equation}\label{cohete}\xymatrix@M=4pt@H=4pt@R=13pt{  \mbf{s}Y \ar[rd]_{\mbf{s}j_Y}    \ar[r]^{\mu_{\Dl\times Y}} & \mrm{R}(\mbf{s}Y) \ar[rd]^{I_{\mbf{s}Y}}                        &     \\
                                                                                               & \mbf{s}(\widetilde{Cyl}(f,g)) \ar[r]^-{\Phi}   &  cyl(\mbf{s}f,\mbf{s}g)\ , \\
                                 \mbf{s}Z\ar[ru]^{\mbf{s}j_Z}     \ar[r]^{\mu_{\Dl\times Z}}  & \mrm{R}(\mbf{s}Z) \ar[ru]_{I_{\mbf{s}Z}}                        &
                                 }\end{equation}
where $j_Y$, $j_Z$ are the canonical inclusions given in
$\ref{inclusionesCylTilde}$.
\end{lema}

\begin{proof}
This result is a consequence of the factorization, additivity and
exactness axioms, together with \ref{relacionCylCylTildeDiag}.\\
Indeed, consider $(\Dl\times f,\Dl\times
g)\in\Omega^{(1)}(\simp\mc{D})$.\\
Then $T=Cyl_{\simp\mc{D}}^{(1)}(\Dl\times f,\Dl\times g)$ is a
bisimplicial object in $\mc{D}$ whose diagonal is just
$\widetilde{Cyl}(f,g)$.\\
By the factorization axiom, $\mu_T:\mbf{s}(\mrm{D}T)\rightarrow
\mbf{s}\comp\simp\mbf{s}(T)=\mbf{s}(n\rightarrow
\mbf{s}(m\rightarrow T_{n,m}))$ is an equivalence.
The simplicial object $\simp\mbf{s}(T)$ is given in degree $n$ by
$$(\simp\mbf{s}(T))_n =\mbf{s}(m\rightarrow Y_m \sqcup X_{m}^{(n-1)} \sqcup  {\cdots} \sqcup X_m^{(0)} \sqcup Z_{m})\ .$$
Hence by (SDC 3) we have an equivalence
$$\sigma_n : \mbf{s}Y \sqcup (\mbf{s}X)^{(n-1)} \sqcup  {\cdots} \sqcup (\mbf{s}X)^{(0)} \sqcup \mbf{s}Z\longrightarrow (\simp\mbf{s}(T))_n \ .$$
Moreover, $\mbf{s}Y \sqcup (\mbf{s}X)^{(n-1)} \sqcup  {\cdots}
\sqcup (\mbf{s}X)^{(0)} \sqcup \mbf{s}Z =
Cyl((\mbf{s}f)\times\Dl,(\mbf{s}g)\times\Dl)_n$ and since
$\sigma_n$ is obtained by the universal property of the coproduct,
we have that the following diagram commutes
$$\xymatrix@M=4pt@H=4pt@R=16pt{  \mbf{s}Y \ar[rd]_-{i_{\mbf{s}Y}}    \ar@/^1pc/[rrd]^{\mbf{s}(i_{Y})}  &                        &     \\
                                                                                                & Cyl((\mbf{s}f)\times\Dl,(\mbf{s}g)\times\Dl)_n \ar[r]^-{\sigma_n}   & (\simp\mbf{s}(T))_n\ . \\
                                 \mbf{s}Z\ar[ru]^-{i_{\mbf{s}Z}}      \ar@/_1pc/[rru]_{\mbf{s}(i_{Z})}  &                         &   }$$
By the naturality of $\sigma$, we get the morphism between
simplicial objects $\varrho=\{\sigma_n\}_n :
Cyl((\mbf{s}f)\times\Dl,(\mbf{s}g)\times\Dl)\rightarrow
\simp\mbf{s}(T)$.\\
Applying $\mbf{s}$ we obtain the commutative diagram
$$\xymatrix@M=4pt@H=4pt@R=13pt{  \mrm{R}(\mbf{s}Y) \ar[rd]_-{I_{\mbf{s}Y}}     \ar@/^1pc/[rrd]^{\mbf{s}(\mbf{s}(i_{Y}))}  &                        &     \\
                                                                                                & cyl(\mbf{s}f,\mbf{s}g) \ar[r]^-{\mbf{s}\varrho}   &  \mbf{s}\comp\simp\mbf{s}(T)\ , \\
                                  \mrm{R}(\mbf{s}Z)\ar[ru]^-{I_{\mbf{s}Z}}      \ar@/_1pc/[rru]_{\mbf{s}(\mbf{s}(i_{Z}))}  &                         &   }$$
where $\mbf{s}\varrho$ is an equivalence by the exactness axiom.\\
F{i}nally, $\mrm{R}(\mbf{s}Y)=\mbf{s}\comp\simp\mbf{s}(\Dl\times
Y)$, and from the naturality of $\mu$ follows that the following
diagram commutes
$$\xymatrix@M=4pt@H=4pt@R=13pt{        \mbf{s}Y   \ar[rd]_{\mbf{s}(j_{Y})}\ar[r]^{\mu_{\Dl\times Y}}  & \mrm{R}(\mbf{s}Y) \ar[rd]^-{I_{\mbf{s}Y}} &     \\
                                                                                   & \mbf{s}(\widetilde{Cyl}(f,g)) \ar[r]^-{\mu_T}   & \mbf{s}\comp\simp\mbf{s}(T) \ . \\
                                       \mbf{s}Z  \ar[ru]^{\mbf{s}(j_{Z})}\ar[r]^{\mu_{\Dl\times Z}}  & \mrm{R}(\mbf{s}Z)\ar[ru]_-{I_{\mbf{s}Z}}&   }$$
Therefore, it suf{f}{i}ces to take
$\Phi=(\mbf{s}\varrho)^{-1}\comp
\mu_T:\mbf{s}(\widetilde{Cyl}(f,g))\rightarrow
cyl(\mbf{s}f,\mbf{s}g)$.
\end{proof}
\begin{num}
Consider the following commutative square of $\mc{D}$
\begin{equation}\label{diagrConmLema}\xymatrix@M=4pt@H=4pt{  X \ar[r]^-{f} \ar[d]_{g} &  Y \ar[d]^{g'}\\
                                      X' \ar[r]^-{f'}        &  Y'   .}\end{equation}
Let $\widehat{g}:c(f)\rightarrow c(f')$ and
$\widehat{f}:c(g)\rightarrow c(g')$ be the morphisms deduced from
the functoriality of the cone, as well as
$\psi:c(\mrm{R}f')\rightarrow c(\widehat{f})$,
$\psi':c(\mrm{R}g')\rightarrow c(\widehat{g})$,
$\widetilde{\lambda}:c(\mrm{R}f')\rightarrow c(f)$ and
$\widehat{\lambda}:c(\mrm{R}g')\rightarrow c(g)$ those obtained
from the following commutative diagrams
$$\xymatrix@M=4pt@H=4pt{\mrm{R}X'\ar[r]^{\mrm{R}f'}\ar[d]_{I} & \mrm{R}Y'\ar[d]_{I}\ar@{}[rd]|{;} & \mrm{R}Y\ar[d]_{I}\ar[r]^{\mrm{R}g'} & \mrm{R}Y'\ar[d]_{I}\ar@{}[rd]|{;} & \mrm{R}X'\ar[r]^{\mrm{R}f'}\ar[d]_{\lambda_{X'}} & \mrm{R}Y'\ar[d]_{\lambda_{Y'}}\ar@{}[rd]|{;}& \mrm{R}Y\ar[r]^{\mrm{R}g'}\ar[d]_{\lambda_{Y}} & \mrm{R}Y'\ar[d]_{\lambda_{Y'}}\\
                                   c(g)          \ar[r]^{\widehat{f}}   & c(g')                             & c(f) \ar[r]^{\widehat{g}}                      & c(f')                             & X'\ar[r]^{f'}                                    & Y'                                          & Y\ar[r]^{g'}                                    & Y'\ ,}$$
where each $I$ denotes the corresponding canonical inclusion.
\end{num}

\begin{cor}\label{coherenciaCono}
Under the above notations, the cone objects $c(\widehat{f})$ and
$c(\widehat{g})$ are naturally isomorphic in $Ho\mc{D}$.\\
If $\widetilde{T}\in\simp\mc{D}$ is $\widetilde{Cyl}$ of the
diagram $1\times\Dl\leftarrow C(g)\rightarrow C(g')$, then there
exists isomorphisms $\Psi:\mbf{s}(\widetilde{T})\rightarrow
c(\widehat{g})$ and $\Phi:\mbf{s}(\widetilde{T})\rightarrow
c(\widehat{f})$ in $Ho\mc{D}$, natural in $(\ref{diagrConmLema})$,
such that the diagram
$$\xymatrix@R=50pt@C=50pt@H=4pt@M=4pt@!0{                                                                   & \mrm{R}^2 Y' \ar'[d][dd]_{I}\ar[ld]_{\mrm{R}I}\ar[rr]^{\lambda_{\mrm{R}Y'}}&                                              &  \mrm{R}Y'\ar'[d][dd]_{I}\ar[ld]_{I}&                                                          & \mrm{R}^2 Y' \ar[dd]^{\mrm{R}I}\ar[ld]_{I} \ar[ll]_{\lambda_{\mrm{R}Y'}}  \\
                                             \mrm{R}c(f')\ar[dd]_{I}\ar[rr]^(0.6){\lambda}                &                                                                          & c(f')\ar[dd]_(0.65){\eta}                          &                                      &  c(\mrm{R}f') \ar[dd]\ar[ll]_(0.4){\widetilde{\lambda}}\ar[dd]^(0.6){\psi}  &                                               \\
                                                                                                            & c(\mrm{R}g')\ar[ld]_{\psi'}\ar'[r][rr]^(0.2){\widehat{\lambda}}&                                              & c(g')\ar[ld]^{\eta'}           &                                                                     & \mrm{R}c(g')\, ,\ar[ld]^{I}\ar'[l][ll]_(0.2){\lambda}\\
                                              c(\widehat{g})                                           &                                                                          & \mbf{s}\widetilde{T}  \ar[rr]^{\Phi}\ar[ll]_{\Psi}&                                      &  c(\widehat{f})                            &  }$$
commutes in $Ho\mc{D}$. The map $\eta:c(f')\rightarrow
\mbf{s}\widetilde{T}$ is the simple of the morphism induced by the
inclusions of $Y'$ and $X'$ into $C(g)$ and ${C}(g')$
respectively, whereas $\eta':c(g')\rightarrow
\mbf{s}\widetilde{T}$ is just the simple of the inclusion of
$C(g')$ into $\widetilde{T}$.
\end{cor}

\begin{proof}
Again, it suf{f}{i}ces to prove the commutativity of
$$\xymatrix@H=4pt@M=4pt{\mrm{R}c(f')\ar[d]_{I}\ar[r]^{\lambda}                 & c(f')\ar[d]_{\eta}                                    & c(\mrm{R}f') \ar[l]_-{\widetilde{\lambda}}\ar[d]^{{\psi}} \\
                        c(\widehat{g})                                         & \mbf{s}\widetilde{T}  \ar[r]^-{\Phi}\ar[l]_-{\Psi}    & c(\widehat{f})\\
                        c(\mrm{R}g')\ar[u]^{{\psi}'}\ar[r]^-{\widehat{\lambda}} & c(g')\ar[u]^{\eta'}                                  &\mrm{R}c(g')\ar[u]^{I}\ar[l]_{\lambda}} \ .
$$
Let us see that this is true for the left side of the diagram,
since the commutativity of the right side can be checked
similarly.\\
We complete diagram (\ref{diagrConmLema}) to
$$\xymatrix@M=4pt@H=4pt{1  & 1\ar[l]\ar[r]   & 1 \\
                        1  \ar[u] \ar[d]& X  \ar[l]\ar[r]^{f} \ar[u] \ar[d]_{g} & Y \ar[u]_{\upsilon} \ar[d]^{g'} \\
                        1 & X'\ar[l]\ar[r]^{f'} & Y'.}$$
The image under $cyl$ of the rows of this diagram is
$$\xymatrix@M=4pt@H=4pt@R=8pt{cyl(1)& c(f)\ar[l]_{\varrho}\ar[r]^{\widehat{g}}& c(f')\ .}$$
Denote by
$\widetilde{\psi}:cyl(\mrm{R}g',\mrm{R}\upsilon)\rightarrow
cyl(\widehat{g},\varrho)$ and
$\overline{\lambda}:cyl(\mrm{R}g',\mrm{R}\upsilon)\rightarrow
c(g')$ the result of applying $cyl$ to
$$\xymatrix@M=4pt@H=4pt@C=17pt{\mrm{R}1\ar[d]_{I}& \mrm{R}Y\ar[l]_{\mrm{R}\upsilon} \ar[r]^{\mrm{R}g'}\ar[d]_{I} & \mrm{R}Y'\ar[d]_{I}\ar@{}[rd]|{;} & \mrm{R}1\ar[d]_{\lambda_{1}} & \mrm{R}Y\ar[l]_{\mrm{R}\upsilon}\ar[r]^{\mrm{R}g'}\ar[d]_{\lambda_{Y}} & \mrm{R}Y'\ar[d]_{\lambda_{Y'}} \\
                              cyl(1)             & c(f)\ar[l]_{\varrho}\ar[r]^{\widehat{g}}                      & c(f')                             & 1                        & Y\ar[l]_{\upsilon}\ar[r]^{g'}                                          & Y'      \ .}$$
By the last proposition we have that if $\widehat{T}$ is
$\widetilde{Cyl}$ of the diagram $Cyl(1\times\Dl)\leftarrow
C(g)\rightarrow C(g')$ then we have an isomorphism $\Psi'$ in
$Ho\mc{D}$ such that the following diagram commutes
$$\xymatrix@H=4pt@M=4pt{\mrm{R}c(f')\ar[r]^{I}\ar[d]_{\lambda} & cyl(\widehat{g},\varrho)              & cyl(\mrm{R}g',\mrm{R}\upsilon)\ar[l]_{\widetilde{\psi}}\ar[d]_-{\overline{\lambda}}\\
                        c(f')\ar[r]^{\widehat{\eta}}           & \mbf{s}\widehat{T}  \ar[u]^-{\Psi'}   & c(g')\ar[l]_{\widehat{\eta}'} \ ,}
$$
where $\widehat{\eta}'$ is the simple of the canonical inclusion
of $C(g')$ into $\widehat{T}$, whereas $\eta$ is the simple of the
morphism induced by the inclusions of  $Y'$, $X'$ and $1\times\Dl$
into $C(g)$, ${C}(g')$ and $Cyl(1\times\Dl)$ respectively.\\
If
$\mrm{R}1\stackrel{J_1}{\rightarrow}cyl(1)\stackrel{I_1}{\leftarrow}\mrm{R}1$
are the canonical inclusions obtained by applying $cyl$ to
$1\leftarrow 1\rightarrow 1$, by \ref{cilaumentado}, we deduce the
existence of a morphism $\rho:cyl(1)\rightarrow\mrm{R}1$
with $\rho I_1=\rho J_1=Id_{\mrm{R}1}$.\\
Since $Id:1\rightarrow 1$ is an equivalence, it follows from the
acyclicity axiom that $J_1$ is so, and because of the 2 out of 3
property we get that $\rho$ is in $\mrm{E}$.\\
Therefore $\rho'=\lambda_1\rho:cyl(1)\rightarrow 1$ is also an
equivalence, and the morphism of cubical diagrams
$$\xymatrix@M=4pt@H=4pt{ cyl(1)\ar[d]_{\rho'}& c(f)\ar[l]_{\varrho}\ar[r]^{\widehat{g}}\ar[d]_{Id} & c(f')\ar[d]_{Id}\\
                                   1                   & c(f)\ar[l]_{}\ar[r]^{\widehat{g}}            & c(f')}$$
gives rise to the equivalence
$\tau:cyl(\widehat{g},\varrho)\rightarrow c(\widehat{g})$ such
that the following diagram commutes
$$\xymatrix@M=4pt@H=4pt{ cyl(\widehat{g},\varrho)\ar[d]_{\tau} &   \mrm{R}c(f')\ar[d]_{Id}\ar[l]_I\\
                                 c(\widehat{g})                       &   \mrm{R}c(f')\ar[l]_{I}\ .}$$
On the other hand, the trivial morphism
$\lambda_1:\mrm{R}1\rightarrow 1$ induces
$$\xymatrix@M=4pt@H=4pt{ \mrm{R}1\ar[d]_{\lambda_1} & \mrm{R}Y \ar[l]_{\mrm{R}\upsilon}\ar[d]_{Id}\ar[r]^{\mrm{R}g'} & \mrm{R}Y'\ar[d]_{Id} \\
                                  1                   & \mrm{R}Y\ar[r]^{\mrm{R}g'}\ar[l]                             & \mrm{R}Y'}$$
and applying $cyl$, we get an equivalence
$\tau':cyl(\mrm{R}g',\mrm{R}\upsilon)\rightarrow c(\mrm{R}g')$.\\
Hence, the following diagram is also commutative
$$\xymatrix@H=4pt@M=4pt{
 \mrm{R}c(f')\ar[r]^{I}                            & c(\widehat{g})                        & c(\mrm{R}g')\ar[l]_{\psi'}\\
 \mrm{R}c(f')\ar[r]^{I}\ar[d]_{\lambda}\ar[u]^{Id} & cyl(\widehat{g},\varrho)\ar[u]^{\tau} & cyl(\mrm{R}g',\mrm{R}\upsilon)\ar[u]^{\tau'}\ar[l]_{\widetilde{\psi}}\ar[d]_-{\overline{\lambda}}\\
 c(f')\ar[r]^{\widehat{\eta}}                            & \mbf{s}\widehat{T}  \ar[u]^-{\Psi'}            & c(g')\ar[l]_{\widehat{\eta}'} \ .}
$$
Indeed, $\tau\comp\widetilde{\psi}=\psi'\comp\tau'$ since both
morphisms are respectively the image under $cyl$ of the morphism
in $\square^{\comp}_1\mc{D}$ given by these two compositions
$$\xymatrix@H=4pt@M=4pt{
\mrm{R}1\ar[d]_{I_1}  & \mrm{R}Y\ar[r]^{\mrm{R}g'}\ar[l]_{\mrm{R}\upsilon}\ar[d]_{I}  & \mrm{R}Y'\ar[d]_{I}           & \mrm{R}1\ar[d]_{}  & \mrm{R}Y\ar[r]^{\mrm{R}g'}\ar[l]_{\mrm{R}\upsilon}\ar[d]_{Id} & \mrm{R}Y'\ar[d]_{Id}  \\
 cyl(1)\ar[d]             & c(f)\ar[r]^{\widehat{g}}\ar[l]_{\varrho} \ar[d]_{Id}      & c(f')\ar[d]_{Id}\ar@{}[r]|{;} & 1\ar[d]            & \mrm{R}Y\ar[r]^{\mrm{R}g'}\ar[l]_{\mrm{R}\upsilon}\ar[d]_{I}  & \mrm{R}Y'\ar[d]_{I}  \\
 1                        & c(f)\ar[r]^{\widehat{g}} \ar[l]                           & c(f')                         & 1                  & c(f)\ar[r]^{\widehat{g}} \ar[l]                               & c(f')                     \ . }$$
Moreover, it holds in $Ho\mc{D}$ that
$\overline{\lambda}\comp(\tau')^{-1}=\widehat{\lambda}$ because in
$\mc{D}$, $\widehat{\lambda}\comp\tau'$ is $cyl$ of the
composition
$$\xymatrix@H=4pt@M=4pt{
 \mrm{R}1 \ar[d]_{\lambda_1} & \mrm{R}Y\ar[r]^{\mrm{R}g'}\ar[l]_{\mrm{R}\upsilon}\ar[d]_{Id} & \mrm{R}Y'\ar[d]_{Id}\\
 1 \ar[d]_{Id} & \mrm{R}Y\ar[r]^{\mrm{R}g'}\ar[l]\ar[d]_{\lambda_Y} & \mrm{R}Y'\ar[d]_{\lambda_{Y'}}\\
 1             & Y\ar[r]^{g'}\ar[l]& Y'\ .
 }$$
that agrees with $\overline{\lambda}$ by def{i}nition.\\
Hence, setting $\Psi''=\tau\comp\Psi'$ we have the following
commutative diagram
\begin{equation}\label{diagAuxiliarC2}
 \xymatrix@H=4pt@M=4pt{
 \mrm{R}c(f')\ar[r]^{I}\ar[d]_{\lambda} & c(\widehat{g})                                  & c(\mrm{R}g')\ar[l]_{\psi'}\ar[d]_{\widehat{\lambda}}\\
 c(f')\ar[r]^{\widehat{\eta}}           & \mbf{s}\widehat{T}  \ar[u]^-{\Psi''}            & c(g')\ar[l]_{\widehat{\eta}'} \ .}
\end{equation}
and it remains to show that
$$\widehat{T}=\widetilde{Cyl}(Cyl(1\times\Dl)\leftarrow C(g)\rightarrow C(g'))\mbox{ and }\widetilde{T}=\widetilde{Cyl}(1\times\Dl\leftarrow C(g)\rightarrow C(g'))$$
are such that $\mbf{s}(\widehat{T})$ and
$\mbf{s}(\widetilde{T})$ are naturally equivalent.\\
Indeed, if $\nu:Cyl(1\times\Dl)\rightarrow 1\times\Dl$ is the
morphisms deduced from \ref{Cylaumentado}, then the diagram
$$\xymatrix@M=4pt@H=4pt{Cyl(1\times\Dl) \ar[d]_{\nu} & C(g)\ar[l]\ar[r] \ar[d]_{Id}  & C(g') \ar[d]_{Id}\\
                        1\times\Dl       & C(g)\ar[l]\ar[r]   & C(g') }$$
gives rise (by applying $\widetilde{Cyl}$) to the morphism
$\vartheta:\widehat{T}\rightarrow\widetilde{T}$ between simplicial
objects. The image under $\mbf{s}$ of the above diagram produces
the morphism in $\square_1^{\circ}\mc{D}$ consisting of
$$\xymatrix@M=4pt@H=4pt{cyl(1)  \ar[d]_{\mbf{s}\nu}  & c(g)\ar[l]_-{\varrho'}\ar[r]^-{\widehat{f}} \ar[d]_{Id}  & c(g') \ar[d]_{Id}\\
                        \mrm{R}1                     & c(g)\ar[l]_-{\nu'}\ar[r]^-{\widehat{f}} & c(g') }$$
and such that $\mbf{s}\nu\in\mrm{E}$, since $cyl(1)=c(Id_1)$ and
$\mrm{R}1$ are equivalent to $1$. Then, by \ref{exactitudcyl} we
have that the induced morphism
$\vartheta':cyl({\widehat{f}},\mbf{s}\varrho')\rightarrow
cyl({\widehat{f}},\nu')$ is an equivalence.\\
F{i}nally, from lemma \ref{relCyltildeconCil} we deduce the
commutative diagram in $\mc{D}$
$$\xymatrix@M=4pt@H=4pt{\mbf{s}\widehat{T} \ar[d]_{\mbf{s}\vartheta}\ar[r]^-{\widehat{\Phi}} & cyl(\widehat{f},\varrho')  \ar[d]_{\vartheta'}\\
                        \mbf{s}\widetilde{T}                  \ar[r]^-{\widetilde{\Phi}}    & cyl({\widehat{f}},\nu') \ .}$$
It follows that
$\mbf{s}\vartheta:\mbf{s}\widehat{T}\rightarrow\mbf{s}\widetilde{T}$
is an equivalence. Consequently, by the def{i}nition of
$\vartheta$ it is clear that, after adjoining this morphism to
(\ref{diagAuxiliarC2}), we get the desired commutative diagram
$$ \xymatrix@H=4pt@M=4pt{
 \mrm{R}c(f')\ar[r]^{I}\ar[d]_{\lambda}      & c(\widehat{g})                                                & c(\mrm{R}g')\ar[l]_{\psi'}\ar[d]_{\widehat{\lambda}}\\
 c(f')\ar[r]^{\widehat{\eta}}\ar[rd]_{\eta}  & \mbf{s}\widehat{T}  \ar[u]^-{\Psi''}\ar[d]_{\mbf{s}\vartheta} & c(g')\ar[l]_{\widehat{\eta}'} \ar[ld]^{\eta'}\\
                                             & \mbf{s}\widetilde{T} \ .                                      & }
$$
\end{proof}

%
%
%
%
%
%

\section{Acyclicity criterion for the cylinder functor}\label{SeccionCritAcCyl}

In this section we develop a generalization of the acyclicity
axiom of the notion of simplicial descent category.

More concretely, the question is if (SDC 7) remains true when we
consider any augmentation $X\rightarrow X_{-1}\times\Dl$ instead
of $X\rightarrow 1\times\Dl$.\\
The answer is af{f}{i}rmative for the ``only if'' part of (SDC 7),
whereas the other part will be true under certain extra
hypothesis.

The acyclicity criterion is necessary in the next chapter in order
to establish the ``transfer lemma'', and it will play a crucial
role in the study of $Ho\mc{D}$.

\begin{prop}\label{axiomacylenD}
Consider morphisms $\xymatrix@M=4pt@H=4pt{ Z & X \ar[r]^{f}
\ar[l]_-{g} & Y }$ in $\mc{D}$. The functor $cyl$ gives rise to
the diagram
$$\xymatrix@H=4pt@M=4pt{
 \mrm{R}X\ar[d]_{\mrm{R}g}\ar[r]^{\mrm{R}f} & \mrm{R}Y\ar[d]_{I_Y} \\
 \mrm{R}Z  \ar[r]^{I_Z} & {cyl}(f,g)}$$
that satisfies the following properties\\
\textbf{a)} $f \in \mrm{E}$ if and only if $I_{Z}$ is in $\mrm{E}$.\\
\textbf{b)} $g \in \mrm{E}$ if and only if $I_{Y}$ is in $\mrm{E}$.\\
In addition, the functor $cyl'$ verifies the analogous properties.
\end{prop}

\begin{proof}
Consider the commutative diagram in $\mc{D}$
$$\xymatrix@M=4pt@H=4pt{1 & 0\ar[l]\ar[r] & 0\\
                        Z  \ar[d]_{Id} \ar[u]& 0  \ar[l]\ar[r] \ar[u] \ar[d]& 0 \ar[u] \ar[d] \\
                        Z  & X\ar[l]_g\ar[r]^f   & Y \ .
                        }$$
If we apply $cyl$ by rows and columns we get
$$\xymatrix@M=4pt@H=4pt@R=10pt{ \mrm{R}1 & \mrm{R}Z\ar[r]^{I_Z}\ar[l]_{\varrho}        & cyl(f,g)\\
                                c(Id_Z)  & \mrm{R}X \ar[r]^{\mrm{R}f}\ar[l]_{\varrho'} & \mrm{R}Y\ .}$$
By the factorization property of the cylinder functor,
\ref{factorizacionCubica}, we deduce that $cyl(I_Z,\varrho)$ and
$cyl(\mrm{R}f,\varrho')$ are isomorphic in $Ho\mc{D}$.\\
Since $Id_Z\in\mrm{E}$, then $c(Id_Z)\rightarrow 1$ is an
equivalence, and therefore the morphism
$cyl(\mrm{R}f,\varrho')\rightarrow c(\mrm{R}f)$ obtained by
applying $cyl$ to
$$\xymatrix@M=4pt@H=4pt{
 c(Id_Z)\ar[d]  & \mrm{R}X \ar[r]^{\mrm{R}f}\ar[l]_{\varrho'}\ar[d]_{Id} & \mrm{R}Y\ar[d]_{Id} \\
 1              & \mrm{R}X\ar[r]^{\mrm{R}f}\ar[l] & \mrm{R}Y}$$
is an equivalence.\\
On the other hand, it follows from the normalization axiom that
$\mrm{R}1\rightarrow 1$ is in $\mrm{E}$, and arguing as before we
obtain an equivalence $cyl(I_Z,\varrho)\rightarrow c(I_Z)$.\\
Therefore $c(\mrm{R}f)$ is isomorphic to $c(I_Z)$ in $Ho\mc{D}$,
so
\begin{center}$c(\mrm{R}f)\rightarrow 1$ is an equivalence if and only if $c(I_Z)\rightarrow 1$ is so. \end{center}
Then, by \ref{axiomaconoenD}, $I_Z$ is an equivalence if and only
if $\mrm{R}f$ is so, and by (SDC 4) this happens if and only if
$f\in\mrm{E}$.\\
The similar result for $cyl'$ follows from the symmetry axiom.
Indeed,
$$cyl'(f,g)=\mbf{s}Cyl'(f\times\Dl,g\times\Dl)=\mbf{s}\Upsilon Cyl \Upsilon (f\times\Dl,g\times\Dl)=\mbf{s}\Upsilon Cyl (f\times\Dl,g\times\Dl),$$
since $\Upsilon(h\times\Dl)=h$ for every morphism $h$ in $\mc{D}$.\\
Hence $I_{Z}=\mbf{s}(i_Z)\in\mrm{E}$ if and only if
$\mbf{s}(\Upsilon i_Z): \mrm{R}Z \rightarrow cyl'(f,g)$ is in
$\mrm{E}$, but this morphism is just the canonical inclusion of
$\mrm{R}Z$ into $cyl'(f,g)$, so we get a).\\
To see b), from the statement a) for $cyl'$ we get that
$g\in\mrm{E}$ if and only if the inclusion $\mrm{R}Y\rightarrow
cyl'(g,f)$ is so. From lemma \ref{relCyltildeconCil} we deduce the
existence of a morphism $\tau:cyl'(g,f)\rightarrow cyl(f,g)$ such
that the following diagram commutes
$$\xymatrix@H=4pt@M=4pt@R=10pt{
 \mrm{R}Y\ar[rd]\ar@/^1pc/[rrd] &                       & \\
                         & cyl'(g,f)\ar[r]^{\tau} & cyl(f,g)\ .\\
 \mrm{R}Z\ar[ru]\ar@/_1pc/[rru] &                       &
}$$
Then $\mrm{R}Y\rightarrow cyl'(g,f)$ is an equivalence if and only
if $\mrm{R}Y\rightarrow cyl(f,g)$ is so, that f{i}nish the proof of b).\\
The statement b) for $cyl'$ can be proved analogously using the
symmetry axiom.
\end{proof}

\begin{lema}\label{axiomaCilTilde}
Consider the morphisms $\xymatrix@M=4pt@H=4pt{ Z & X \ar[r]^{f}
\ar[l]_-{g} & Y }$ in $\simp\mc{D}$. The functors
$\widetilde{Cyl}$ and $\mbf{s}$ give rise to a diagram in $\mc{D}$
\begin{equation}\label{diagrCylTildeLema}\xymatrix@H=4pt@M=4pt{
 \mbf{s}X\ar[d]_{\mbf{s}g}\ar[r]^{\mbf{s}f} & \mbf{s}Y\ar[d]_{\mbf{s}j_Y} \\
 \mbf{s}Z  \ar[r]^-{\mbf{s}j_Z} & \mbf{s}\widetilde{Cyl}(f,g)}
 \end{equation}
such that\\
\textbf{a)} $\mbf{s}f \in \mrm{E}$ if and only if $\mbf{s}j_{Z}$ is in $\mrm{E}$.\\
\textbf{b)} $\mbf{s}g \in \mrm{E}$ if and only if $\mbf{s}j_{Y}$
is in $\mrm{E}$.
\end{lema}

\begin{proof}
F{i}rst, assume that $X$, $Y$ and $Z$ are objects in $\mc{D}$. In
this case, by \ref{RelacionCilyCilTildeenD} a),
$\mbf{s}\widetilde{Cyl}(f\times\Dl,g\times\Dl)=cyl(f,g)$  and the statement follows from proposition \ref{axiomacylenD}.\\
If $X$, $Y$ and $Z$ are simplicial objects, by diagram
(\ref{cohete}) of lemma \ref{relCyltildeconCil} we deduce that
$\mbf{s}j_{Z}\in\mrm{E}$ if and only if $I_{\mbf{s}Z}\in\mrm{E}$.
Then it follows from the constant case that this holds if and
only if $\mbf{s}f\in\mrm{E}$.\\
Similarly, $\mbf{s}j_Y$ if and only if $I_{\mbf{s}Y}$ is so, if
and only if $\mbf{s}g$ is so, by the constant case.
\end{proof}

\begin{thm}\label{AxCil} Let $f:X\rightarrow Y$ be a morphism in $\simp\mc{D}$ and $\epsilon:X\rightarrow
X_{-1}\times\Dl$ an augmentation. Then\\
\textbf{a)} if $\mbf{s}f$ is an equivalence then the simple of
$i_{X_{-1}}:X_{-1}\times\Dl\rightarrow Cyl(f,\epsilon)$ is so.\\
\textbf{b)} if $\mbf{s}\epsilon$ is an equivalence then the simple
of $i_{Y}:Y\rightarrow Cyl(f,\epsilon)$ is so.
\end{thm}

\begin{proof}
Assume that $\mbf{s}f\in\mrm{E}$. Then we deduce from
\ref{axiomaCilTilde} that
$\mbf{s}j_{X_{-1}}:\mrm{R}X_{-1}\rightarrow
\mbf{s}\widetilde{Cyl}(f,\epsilon)$ is in $\mrm{E}$.\\
On the other hand, by \ref{CylRetractoCylTilde} we deduce that the
following diagram commutes in $\mc{D}$
$$\xymatrix@M=4pt@H=4pt{  \mbf{s}(Cyl(f,\epsilon))\ar[r]^{\mbf{s}\alpha}                & \mbf{s}(\widetilde{Cyl}(f,\epsilon))\ar[r]^{\mbf{s}\beta} & \mbf{s}(Cyl(f,\epsilon)) \\
                         \mrm{R}X_{-1} \ar[r]^{Id}\ar[u]_{\mbf{s}(i_{X_{-1}})}         & \mrm{R}X_{-1}\ar[r]^{Id}\ar[u]_{\mbf{s}(j_{X_{-1}})}     & \mrm{R}X_{-1}\ar[u]_{\mbf{s}(i_{X_{-1}})}\ ,}$$
where $\mbf{s}\beta\comp\mbf{s}\alpha=Id$.\\
Therefore $\mbf{s}(i_{X_{-1}})$ is a retract of
$\mbf{s}(j_{X_{-1}})\in\mrm{E}$, and by \ref{propRetracto}
$\mbf{s}(i_{X_{-1}})\in\mrm{E}$.\\
To see b), one can argue in a similar way.
\end{proof}

The previous result remains valid using the symmetric notion of
cylinder, $Cyl'$, introduced in \ref{defCilPrima}.

\begin{cor}\label{AxCilPrima}
Let $f:X\rightarrow Y$ be a morphism in $\simp\mc{D}$ and
$\epsilon:X\rightarrow
X_{-1}\times\Dl$ be an augmentation. Then\\
\textbf{a)} if $\mbf{s}f$ is an equivalence, the simple of
$i_{X_{-1}}:X_{-1}\times\Dl\rightarrow Cyl'(f,\epsilon)$ is so.\\
\textbf{b)} if $\mbf{s}\epsilon$ is an equivalence, the simple of
$i_{Y}:Y\rightarrow Cyl'(f,\epsilon)$ is so.
\end{cor}

\begin{proof} Again, the statement follows from the previous proposition together with (SDC 8).
Let us see b), since a) can be proved analogously.\\
If $\mbf{s}\epsilon\in\mrm{E}$ then $\mbf{s}\Upsilon\epsilon$ is
also in $\mrm{E}$. By the previous proposition, the simple of
$\Upsilon Y \rightarrow Cyl(\Upsilon f,\Upsilon\epsilon)$ is an equivalence.\\
Therefore the simple of $Y \rightarrow \Upsilon Cyl(\Upsilon
f,\Upsilon\epsilon)=Cyl'(f,\epsilon)$ is in $\mrm{E}$, and we are
done.
\end{proof}

\begin{thm}\label{AxCil+retracto}
Let $\xymatrix@M=4pt@H=4pt{ X_{-1}\times\Dl & X \ar[r]^-{f}
\ar[l]_-{\epsilon} &  Y }$ be a diagram in $\simp\mc{D}$ such that
there exists $\epsilon':Y\rightarrow X_{-1}\times\Dl$ with
$\epsilon'f=\epsilon$. Then
\begin{center}$\mbf{s}f\in \mrm{E}$ if and only if $\mbf{s}i_{X_{-1}}:\mrm{R}X_{-1}\rightarrow \mbf{s}Cyl(f,\epsilon)$ is in $\mrm{E}$.\end{center}
\end{thm}

\begin{proof}
The ``only if'' part follows from \ref{AxCil}. Assume that
$\mbf{s}i_{X_{-1}}\in\mrm{E}$.
From the commutativity of the diagram of $\simp\mc{D}$
$$\xymatrix@M=4pt@H=4pt{
 X \ar[r]^-{f}\ar[d]_{\epsilon} &  Y\ar[d]_{\epsilon'}\\
 X_{-1}\times\Dl\ar[r]^{Id}     &  X_{-1}\times\Dl}$$
we get by \ref{Cylaumentado} the existence of a morphism
$H:Cyl(f,\epsilon)\rightarrow X_{-1}\times\Dl$ such that
$H\comp i_{X_{-1}}=Id_{X_{-1}\times\Dl}$ and $H\comp i_Y=\epsilon'$.\\
Since $\mbf{s}i_{X_{-1}}\in\mrm{E}$, it follows from the 2 out of 3 property that $\mbf{s}H$ is an equivalence.\\
On the other hand, applying \ref{lemaDiagrCubico} to the diagram
$$\xymatrix@R=33pt@C=33pt@H=4pt@M=4pt@!0{                                 & X_{-1}\times\Dl  \ar[rr]\ar[ld]_{Id}          &                   & 1\times\Dl \\
                                         X_{-1}\times\Dl                  &                                               &                   &\\
                                                                          & X \ar[uu]_{\epsilon}\ar[rr]^-{Id} \ar[ld]_{f} &                   & X \ ,\ar[uu]\ar[ld]^{f}\\
                                         Y \ar[uu]^{\epsilon'}\ar[rr]_{Id}&                                               &            Y      &
}$$
we get an isomorphism between the simplicial object obtained as
the image under $Cyl$ of
$$\xymatrix@H=4pt@M=4pt{ 1\times\Dl & Cyl'(Id_X,\epsilon)\ar[l]\ar[r]^-{F} & Cyl'(Id_Y,\epsilon')}$$
and the image under $Cyl'$ of
$$\xymatrix@H=4pt@M=4pt{ X_{-1}\times\Dl & Cyl(f,\epsilon)\ar[l]_-{H}\ar[r]^-{G} & C(f)}\ .$$
In other words, $C(F)\simeq Cyl'(G,H)$.\\
As $\mbf{s}(Id_X)=Id_{\mbf{s}X}\in\mrm{E}$ then by
\ref{AxCilPrima} we deduce that the simple of the canonical
inclusion $j_{X_{-1}}:X_{-1}\times\Dl\rightarrow
Cyl'(Id_X,\epsilon)$ is an equivalence. Similarly, the same holds
for the simple of $l_{X_{-1}}:
X_{-1}\times\Dl\rightarrow Cyl'(Id_Y,\epsilon')$.\\
Moreover, from the naturality of $Cyl'$ we get that $F\comp
j_{X_{-1}}=l_{X_{-1}}$, hence $\mbf{s}F\in\mrm{E}$.\\
Then, by the acyclicity axiom, $\mbf{s}C(F)\rightarrow 1$ is in
$\mrm{E}$, consequently $\mbf{s}Cyl'(G,H)\rightarrow 1$ is also an equivalence.\\
But $\mbf{s}H\in\mrm{E}$, and again it follows from
\ref{AxCilPrima} that the simple of $C(f)\rightarrow Cyl(G,H)$ is
in $\mrm{E}$. Therefore $\mbf{s}C(f)$ is acyclic, and from the
acyclicity criterion we deduce that $\mbf{s}f\in\mrm{E}$.
\end{proof}

\begin{prop}\label{CDSimplicadg+}\mbox{}\\
i) If $f,g$ are homotopic morphisms in $\simp\mc{D}$ (see
\ref{def{i}homotSimplicial}) then
$\mbf{s}(f)=\mbf{s}(g)$ in $Ho\mc{D}$.\\
ii) If $\epsilon:X\rightarrow X_{-1}\times\Dl$ has a (lower or upper) extra degeneracy (see \ref{def{i}degmas}) then $\mbf{s}(\epsilon)$ is an equivalence.\\
\end{prop}

\begin{proof}
Let $\widetilde{Cyl}(X)=\widetilde{Cyl}(Id_X,Id_X)$ be the cubical
cylinder object associated with $X$, given in \ref{ejCilindroCub}.
The morphisms f and g are homotopic in $\simp\mc{D}$, so there
exists a homotopy $H:\widetilde{Cyl}(X)\rightarrow X$ such that
the following diagram commutes in $\simp\mc{D}$
$$\xymatrix@M=4pt@H=4pt@R=15pt{ X \ar[rd]^I    \ar@/^1pc/[rrd]^{f}&                         &     \\
                                                           & \widetilde{Cyl}(X) \ar[r]^{H}   &  X \ .\\
                          X\ar[ru]_{J}  \ar@/_1pc/[rru]_{g}&                         & }$$
Then $\mbf{s}(f)=\mbf{s}(H)\comp \mbf{s}(I)$ and
$\mbf{s}(g)=\mbf{s}(H)\comp \mbf{s}(J)$ in $\mc{D}$, so it is
enough to check the equality $\mbf{s}(I)=\mbf{s}(J)$ in $Ho\mc{D}$.\\
If $cyl(\mbf{s}X)=cyl(Id_{\mbf{s}X},Id_{\mbf{s}X})$, by
\ref{relCyltildeconCil}, it suf{f}{i}ces to see that the
inclusions $I_{\mbf{s}X}$ and $J_{\mbf{s}X}:\mrm{R}X\rightarrow
cyl(\mbf{s}X)$ coincide in $Ho\mc{D}$.
Note that both morphisms $I_{\mbf{s}X}$, $J_{\mbf{s}X}$ are equivalences, because of \ref{axiomacylenD}.\\
On the other hand, it follows from \ref{cilaumentado} the
existence of $\rho:cyl(\mbf{s}(X))\rightarrow \mrm{R}\mbf{s}(X)$
such that $\rho\comp
I_{\mrm{R}s(X)}=\rho\comp J_{\mrm{R}s(X)}=Id$.\\
Hence, in $Ho\mc{D}$, $\rho$ is an isomorphism such that
$I_{\mrm{R}s(X)}=J_{\mrm{R}\mbf{s}X}=\rho^{-1}$, that f{i}nish the
proof.
\textit{ii)} follows from \textit{i)}, having into account
\ref{degmasHomot}.
\end{proof}
%
%

\begin{cor}\label{conmutcyl}
If $D$ is the object ${X_{-1}\times\Dl\leftarrow X\rightarrow Y}$
of $\Omega(\mc{D})$, the following diagram commutes in
$Ho{\mc{D}}$
$$\xymatrix@M=4pt@H=4pt{ \mbf{s}(X) \ar[r]^{\mbf{s}(f)}   \ar[d]_{\mbf{s}(\epsilon)}           & \mbf{s}(Y) \ar[d]^{\mbf{s}(i_Y)}\\
                       \mbf{s}(X_{-1}\times\Dl) \ar[r]^{\mbf{s}(i_{X_{-1}})}&  \mbf{s}(Cyl(D)) .} $$
\end{cor}
\begin{proof}
In the f{i}rst chapter we proved that $i_Y\comp f$ is homotopic to
 $i_{X_{-1}}\comp\epsilon$ (see
\ref{CylHomotopia}), therefore the statement is a consequence of
the above proposition.
\end{proof}

\begin{cor}\label{comptrivial}
If $f:X\rightarrow Y$  is a morphism between simplicial objects
then the composition $\mbf{s}(X)\stackrel{\mbf{s}(f)}{\rightarrow}
\mbf{s}(Y)\rightarrow \mbf{s}(Cf)$ is trivial in $Ho\mc{D}$, that
is, it factors through the f{i}nal object 1.
\end{cor}

\begin{proof}
By the previous proposition we have the following commutative
diagram in $Ho\mc{D}$
$$\xymatrix@M=4pt@H=4pt{ \mbf{s}(X) \ar[r]^{\mbf{s}(f)}   \ar[d]          & \mbf{s}(Y) \ar[d]^{\mbf{s}(i_Y)}\\
                       \mrm{R}1 \ar[r]&  \mbf{s}(C(f)) .} $$
The result follows from the equivalence existing between
$\mrm{R}1$ and $1$ (by (SDC 5)), so they are isomorphic in
$Ho\mc{D}$.
\end{proof}
%

\begin{cor}\label{pdadesCylenD}
Given an object $\xymatrix@M=4pt@H=4pt{Z &\ar[l]_{g}  X\ar[r]^f &
Y}$ in $\square_1^\comp \mc{D}$, then the following diagram
commutes in $Ho\mc{D}$
$$\xymatrix@M=4pt@H=4pt{\mrm{R}X \ar[d]_{\mrm{R}g}\ar[r]^{\mrm{R}f} & \mrm{R}Y \ar[d]^{I_Y}\\
                        \mrm{R}Z \ar[r]^{I_Z} & cyl(f,g) \ .}$$
Moreover, if $cyl(X)=cyl(Id_X,Id_X)$, there exists
$H:cyl(X)\rightarrow cyl(f,g)$ such that $H\comp
I_X=I_Y\comp\mrm{R}f$ and $H\comp J_X=I_Z\comp \mrm{R}g$, where
$I,J$ are the inclusions of $X$ into $cyl(X)$.\\
In particular, the composition $\xymatrix@M=4pt{\mrm{R}X
\ar[r]^{\mrm{R}f} & \mrm{R}Y \ar[r]^{I_Y} & c(f)}$ is trivial in
$Ho\mc{D}$.
\end{cor}
%

\begin{proof}
Having in mind the commutativity up to homotopy of the diagram of
simplicial objects
$$\xymatrix@M=4pt@H=4pt{ X\times\Dl \ar[r]^-{f\times\Dl}   \ar[d]          & Y\times\Dl \ar[d]^{i_Y}\\
                         Z\times\Dl \ar[r]^-{i_Z} & Cyl(f\times\Dl,g\times\Dl) \ ,} $$
there exists $L:\widetilde{Cyl}(X\times\Dl)\rightarrow
Cyl(f\times\Dl,g\times\Dl)$ such that $L\comp i_X =i_Y\comp
f\times\Dl$ and $L\comp j_X=i_Z\comp g\times\Dl$. In addition,
$\mbf{s}\widetilde{Cyl}(X\times\Dl)$ is equal to $cyl(X)$,
therefore it is enough to take $H=\mbf{s}L$.
\end{proof}


%
%
\begin{obs}\label{cubConmenD}
Consider the commutative diagram in $\mc{D}$
$$\xymatrix@M=4pt@H=4pt{ X \ar[d]_{\alpha} & Y \ar[r]^f \ar[l]_g \ar[d]_{\beta} &  Z \ar[d]_{\gamma}\\
                         X' & Y' \ar[r]^{f'} \ar[l]_{g'} &  Z'  . }$$
From the functoriality of $cyl$ we obtain
$\delta:cyl(f,g)\rightarrow cyl(f',g')$ such that in the diagram
$$\xymatrix@R=40pt@C=40pt@H=4pt@M=4pt@!0{                                                                &  \mrm{R}Y\ar[rr]^{\mrm{R}f}\ar'[d][dd]^{\beta}\ar[ld]_{\mrm{R}g}&                            & \mrm{R}Z \ar[dd]^{\gamma}\ar[ld]_{I_Z}\\
                                         \mrm{R}X \ar[rr]^{\;\;\;\;\;\;\;\;I_X}\ar[dd]_{\alpha} &                                                                                       & cyl(f,g)  \ar[dd]\\
                                                                                                         & \mrm{R}Y'  \ar'[r][rr]^{\mrm{R}f'}\ar[ld]_{\mrm{R}g'}                                &                            & \mrm{R}Z'\, ,\ar[ld]^{I_{Z'}}\\
                                             \mrm{R}X'  \ar[rr]^{I_{X'}}                                     &                                                                                       &  cyl(f',g') }$$
all the faces commute in $\mc{D}$, except the lower and upper
ones, that commute in $Ho\mc{D}$.
\end{obs}

To f{i}nish this section, we give the following result just for
completeness, because we will not use it in these notes. It is a
consequence of proposition \ref{axiomacylenD}.

\begin{prop} If $f,g$ are morphisms in $\mc{D}$ such that $f\sqcup
g$ is an equivalence, then $f$ and $g$ are so.
\end{prop}

\begin{proof} Let $f:A\rightarrow B$ and $g:A'\rightarrow B'$ be morphisms such that $f\sqcup g\in\mrm{E}$, and let us prove that $f\in\mrm{E}$.
We will check that the canonical inclusion $\mrm{R}1\rightarrow
c(f)$ is in $\mrm{E}$ (that is equivalent to check that $c(f)\rightarrow 1$ is in $\mrm{E}$).\\
Consider the trivial morphisms $\rho_A:A\rightarrow 1$ and
$\rho_{A'}:A'\rightarrow 1$, as well as the morphism
$\rho=\rho_A\sqcup\rho_{A'}:A\sqcup A'\rightarrow 1\sqcup 1$.
Since $f\sqcup g\in\mrm{E}$, we deduce from \ref{axiomacylenD}
that the inclusion
$$I:\mrm{R}(1\sqcup 1)\rightarrow cyl(f\sqcup g,\rho)$$
is an equivalence. If we denote by $I_f:\mrm{R}1\rightarrow c(f)$
and $I_g:\mrm{R}1\rightarrow c(g)$ the canonical inclusions, by
\ref{aditividadcylyc} we get an equivalence $\sigma_{cyl}$ such
that the diagram
$$\xymatrix@M=4pt@H=4pt@C=35pt{ \mrm{R}1\sqcup \mrm{R}1  \ar[d]_{\sigma_{\mrm{R}}} \ar[r]^-{I_f\sqcup I_g} & c(f)\sqcup c(g)  \ar[d]_{\sigma_{cyl}} \\
                         \mrm{R}(1\sqcup 1)          \ar[r]^-{I}            & cyl(f\sqcup g,\rho)  }$$
commutes. Moreover, from (SDC 3) we deduce that
$\sigma_{\mrm{R}}\in\mrm{E}$, and therefore $I_f\sqcup I_g$ is an
equivalence. F{i}nally, it is enough to see that $I_f$
is a retract of $I_f\sqcup I_g$, in such a case the proof would be concluded by \ref{propRetracto}.\\
\indent The ``zero'' morphism $\alpha:c(g)\rightarrow c(f)$ is
def{i}ned as follows. The morphism $C(g\times\Dl)\rightarrow
1\times\Dl$ gives rise to $c(g)\rightarrow \mrm{R}1$, and by
composing with $I_f$ we get the desired morphism $\alpha$.
Moreover, $\alpha\comp I_g=I_f$ since at the simplicial level
$1\times\Dl\rightarrow C(g\times\Dl)\rightarrow 1\times\Dl$
is the identity.\\
Then, we obtain the commutative diagram
$$\xymatrix@M=4pt@H=4pt@C=35pt{ c(f)    \ar[r]                   & c(f)\sqcup c(g) \ar[r]^-{Id\sqcup \alpha}                               & c(f)                  \\
                         \mrm{R}1   \ar[u]_{I_f} \ar[r]   & \mrm{R}1\sqcup \mrm{R}1 \ar[u]_{I_f\sqcup I_g}  \ar[r]^-{Id\sqcup Id}  & \mrm{R}1   \ar[u]_{I_f}   }$$
where the horizontal compositions are the identity, and it follows
that $I_f$ is in fact a retract of $I_f\sqcup I_g$.
\end{proof}

%


%

\section{Functors of simplicial descent categories}

The aim of this section is to state and prove the ``transfer
lemma'', that will allow us to transfer the simplicial descent
structure from a f{i}xed simplicial descent category to a new
category through a suitable functor.

%
%
%

%
\begin{def{i}}\label{Def{i}FuntorcatDesc}
Let $(\mc{D},\mbf{s},\mrm{E},\mu,\lambda)$ and
$(\mc{D'},\mrm{E'},\mbf{s'},\mu',\lambda')$ be simplicial descent
categories. We say that a functor $\psi:\mc{D}\rightarrow \mc{D}'$
is a \textit{functor of simplicial descent categories}\index{Index}{functor!of simpl. descent categories} if\\
$\mathbf{(FD\; 0)}$ $\psi$ preserves equivalences, that is, $\psi(\mrm{E})\subseteq \mrm{E}'$.\\
$\mathbf{(FD\; 1)}$ $\psi$ is a quasi-strict monoidal functor with respect to the coproduct (see \ref{FuntorMonoidalCasiEstricto}).\\
$\mathbf{(FD\; 2)}$ Consider the diagram
\begin{equation}\label{conmutSimples}\xymatrix@H=4pt@C=20pt@R=20pt{\simp\mc{D} \ar[r]^{\simp\psi}\ar[d]_{\mbf{s}} & \simp\mc{D}' \ar[d]^{\mbf{s'}} & \\
                                \mc{D}\ar[r]^{\psi}         &                                \mc{D}'        \ar[r]^-{\gamma}                        & Ho\mc{D}' \ .}\end{equation}
There exists a natural isomorphism of functors
$\Theta:\gamma\comp\psi\comp\mbf{s} \rightarrow
\gamma\comp\mbf{s}'\comp\simp\psi$ in such a way that $\Theta$
comes from a functorial ``zig-zag'' with values in $\mc{D}'$.
Moreover, $\Theta$ must be compatible with $\lambda$, $\lambda'$
and with $\mu$, $\mu'$.\\[0.2cm]
More concretely, there exists functors $A^0,\ldots ,
A^r:\simp\mc{D}\rightarrow \mc{D}'$ such that
$A^0=\psi\comp\mbf{s}$ and $A^r= \mbf{s}'\comp\simp\psi$, and they
are related by the natural transformations
\begin{equation}\label{RistraTheta}\xymatrix@M=4pt@H=4pt@C=25pt{
 \psi\comp\mbf{s}=A^0 \ar@{-}[r]^-{\Theta^0} & A^1 \ar@{-}[r]^-{\Theta^1} & \cdots \cdots\ar@{-}[r]^-{\Theta^{r-1}} & A^r=\mbf{s}'\comp\simp\psi}
\end{equation}
where either $\Theta^i:A^i\rightarrow A^{i+1}$ or
$\Theta^i:A^{i+1}\rightarrow A^{i}$, and such that
$\Theta^i_X\in\mrm{E}'$ for all $X$ in $\simp\mc{D}$ and for all
$i$. The natural transformation $\Theta$ must be the image under
$\gamma$ of the zig-zag (\ref{RistraTheta}).
\end{def{i}}

\begin{num}\label{compatFD}
Let us describe more specif{i}cally the compatibility condition
that is mentioned in (FD 2). We will denote also by $\psi$ the
induced morphisms $\simp\psi:\simp\mc{D}\rightarrow \simp\mc{D}'$
and $\simp\simp\psi:\simp\simp\mc{D}\rightarrow
\simp\simp\mc{D}'$.\\[0.2cm]
\textsc{i.}  Given an object $X$ in $\mc{D}$, the following
diagram must commute in $Ho\mc{D}'$
\begin{equation}\label{compatFDii}\xymatrix@M=4pt@H=4pt@C=50pt@R=40pt{\psi(\mbf{s}(X\times\Dl)) \ar[r]^-{\psi(\lambda_X)}   \ar[d]_-{\Theta_{X\times\Dl}}  & \psi(X).\\
                                      \mbf{s'}(\psi(X)\times\Dl) \ar[ru]_-{\lambda'_{\psi(X)}} &}\end{equation}
\textsc{ii.} If $Z\in\simp\simp\mc{D}$, then $\psi(\mrm{D}Z)=\mrm{D}(\psi(Z))$ and $\Theta_{\mrm{D}Z}:\gamma\psi(\mbf{s}(\mrm{D}Z))\rightarrow\gamma\mbf{s'}\mrm{D}(\psi(Z))$.\\
The following diagram must be commutative in $Ho\mc{D'}$
\begin{equation}\label{compatFDi}\xymatrix@M=4pt@H=4pt@C=40pt@R=40pt{\mbf{s'}\mrm{D}(\psi Z)\ar[r]^-{\mu'_{\psi(Z)}} &   \mbf{s'} \simp\mbf{s'}(\psi Z)   \\
                                                                     \psi\mbf{s}(\mrm{D}Z)\ar[u]^{\Theta_{\mrm{D}Z}}    \ar[r]^{\psi(\mu_{Z})}  & \psi\mbf{s}(\simp\mbf{s}Z) \ar[u]_{(\mbf{s}'\simp\Theta\comp\Theta)_Z}   ,}\end{equation}
where the natural transformation
$\mbf{s}'\simp\Theta\comp\Theta:\psi\comp\mbf{s}\comp\simp\mbf{s}\rightarrow\mbf{s}'\comp\simp\mbf{s}'\comp\psi$
is def{i}ned in $Ho\mc{D}'$ and is induced in a natural way by $\Theta$.\\
Concretely, on one hand we have that
$$\Theta_{\simp\mbf{s}Z}: \psi\mbf{s}(\simp\mbf{s}Z)\longrightarrow \mbf{s}'\psi(\simp\mbf{s}Z)\ .$$
On the other hand, f{i}xing as $n$ the f{i}rst index of $Z$ we get
$Z_{n,\cdot}\in\simp\mc{D}$, and the evaluation of
(\ref{RistraTheta}) in $Z_{n,\cdot}$ is the natural sequence of
morphisms in $\mc{D}'$
$$\xymatrix@M=4pt@H=4pt@C=30pt{
 \psi(\mbf{s}Z_{n,\cdot}) \ar@{-}[r]^-{\Theta^0_{Z_{n,\cdot}}} & A^1(Z_{n,\cdot}) \ar@{-}[r]^-{\Theta^1_{Z_{n,\cdot}}} & \cdots \cdots\ar@{-}[r]^-{\Theta^{r-1}_{Z_{n,\cdot}}} & \mbf{s}' (\psi Z_{n,\cdot})\ .}
$$
Therefore, we obtain the sequence in $\simp\mc{D}'$
$$\xymatrix@M=4pt@H=4pt@C=35pt{
 \psi(\simp\mbf{s}Z) \ar@{-}[r]^-{\simp\Theta^0_{Z}} & \simp A^1Z \ar@{-}[r]^-{\simp\Theta^1_{Z}} & \cdots \cdots\ar@{-}[r]^-{\simp\Theta^{r-1}_{Z}} & \simp\mbf{s}' (\psi Z)\ .}
$$
and, by the exactness axiom, the result of applying $\mbf{s}'$ is
the sequence of equivalences in $\mc{D}'$
$$\xymatrix@M=4pt@H=4pt@C=25pt{
 \mbf{s}'\psi(\simp\mbf{s}Z) \ar@{-}[r] & \mbf{s}'\simp A^1 Z \ar@{-}[r] & \cdots \cdots\ar@{-}[r] & \mbf{s}'\simp\mbf{s}' (\psi Z)\ .}
$$
that gives rise to the morphism in $Ho\mc{D}'$
$$(\mbf{s}'\simp\Theta)_Z: \mbf{s}'\simp(\psi\mbf{s})Z\longrightarrow \mbf{s}'\simp\mbf{s}' (\psi Z) $$
Then, $(\Theta\comp\mbf{s}'\simp\Theta)_Z$ is the composition
$$\xymatrix@H=4pt@M=4pt@C=40pt{\psi\mbf{s}(\simp\mbf{s}Z) \ar[r]^-{\Theta_{\simp\mbf{s}Z}} & \mbf{s}'\psi(\simp\mbf{s}Z)\ar[r]^-{\mbf{s}'(\simp\Theta)_Z} & \mbf{s}'\comp\simp\mbf{s}'(\psi Z)}\ .$$
\end{num}

\begin{obs}[Case $\Theta=Id$]
Assume that (\ref{conmutSimples}) is commutative,  that is,
$\Theta=Id: \psi\comp\mbf{s} \rightarrow \mbf{s}'\comp\simp\psi$
as functors $\simp\mc{D}\rightarrow \mc{D}'$. In this case the
compatibility condition between $\lambda$, $\lambda'$ and $\mu$,
$\mu'$ means that the following equalities hold in $Ho\mc{D}'$
$$\psi(\lambda(X))=\lambda'(\psi(X)) \;\; \forall X\in\mc{D}\;\;\; ; \;\;\; \psi(\mu(Z))=\mu'(\psi(Z))\;\; \forall Z\in\simp\simp\mc{D}\ .$$
\end{obs}

\begin{ej} If $F:\mc{D}\rightarrow \mc{D}'$ is a functor of simplicial descent categories and $I$
is a small category, then $F_\ast:I\mc{D}\rightarrow I\mc{D}'$ is
also a functor of simplicial descent categories, with the descent
structures introduced in \ref{CategFuntores}.\\
In addition, if $\mc{D}$ is a simplicial descent category and $x$
is an object in $\mc{D}$, the ``evaluation at'' $x$ functor,
$ev_x:I\mc{D}\rightarrow \mc{D}$ with $ev_x(f)=f(x)$, is a functor
of simplicial descent categories.
\end{ej}

\begin{obs} In the following lemma we will prove that the
composition of two functors of simplicial descent categories is
again a functor in this way. Hence, we have the category (in a
convenient universe) $\mathfrak{D}es\mc{S}imp$ of simplicial
descent categories together with the functors of simplicial
descent categories.
\end{obs}

\begin{lema} The composition of two functors of simplicial descent
categories is again a functor of simplicial descent categories.
\end{lema}

\begin{proof} Let $\psi:\mc{D}\rightarrow\mc{D}'$ and
$\psi':\mc{D}'\rightarrow\mc{D}''$ be functors of simplicial
descent categories. We shall study the commutativity of the
diagrams
$$\xymatrix@H=4pt@C=20pt@R=20pt{\simp\mc{D} \ar[r]^-{\simp\psi}\ar[d]_{\mbf{s}} & \simp\mc{D}' \ar[d]^{\mbf{s'}} \ar[r]^-{\simp\psi'}  & \simp\mc{D}'' \ar[d]^{\mbf{s''}}\\
                                \mc{D}\ar[r]^-{\psi}                            & \mc{D}'        \ar[r]^-{\psi'}                       & \mc{D}'' \ .}
$$
Assume that the zig-zag $\Theta$ and $\Phi$ associated with $\psi$
and $\psi'$ are given respectively by
$$\xymatrix@M=4pt@H=4pt@C=25pt@R=9pt{
  \psi\comp\mbf{s}=A^0 \ar@{-}[r]^-{\Theta^0} & A^1 \ar@{-}[r]^-{\Theta^1} & \cdots \cdots\ar@{-}[r]^-{\Theta^{r-1}} & A^r=\mbf{s}'\comp\simp\psi \\
  \psi'\comp\mbf{s'}=B^{0} \ar@{-}[r]^-{\Phi^0} & B^1 \ar@{-}[r]^-{\Phi^1} & \cdots \cdots\ar@{-}[r]^-{\Phi^{s-1}} & B^s=\mbf{s''}\comp\simp\psi'\ .}
$$
Then $\Psi=(\Phi\comp \simp\psi)\comp (\psi'\comp\Theta)$ is the
zig-zag relating $\psi'\comp\psi\comp\mbf{s}$ to
$\mbf{s''}\simp(\psi'\comp\psi)$, whose value at $X\in\simp D$ is
$$\xymatrix@M=4pt@H=4pt@C=29pt{
  \psi'\psi\mbf{s}X \ar@{-}[r]^-{\psi'\Theta^0_X} &  \cdots \cdots\ar@{-}[r]^-{\psi'\Theta^{r-1}_X} & \psi'\mbf{s}'\comp\simp\psi X  \ar@{-}[r]^-{\Phi^0_{\simp\psi X}} &
  \cdots \cdots\ar@{-}[r]^-{\Phi^{s-1}_{\simp\psi X}} & \mbf{s''}\simp(\psi'\psi) X }
$$
The compatibility of $\Psi$ with $\lambda$ and $\lambda''$ follows
from the one of $\Theta$ and $\Phi$ with $\lambda$, $\lambda'$ and
$\lambda''$. On the other hand, if $Z\in\simp\simp\mc{D}$, the
square (\ref{compatFDi}) can be drawn in this case as
\begin{equation}\label{CuadrAuxCompfunt}
\xymatrix@M=4pt@H=4pt@C=40pt@R=20pt{\mbf{s''}\mrm{D}(\psi'\psi Z)\ar[r]^-{\mu''_{\psi'\psi Z}}                         &   \mbf{s''} \simp\mbf{s''}(\psi'\psi Z)   \\
                                                                                                                          & \mbf{s''}\psi'\psi(\simp\mbf{s}Z) \ar[u]_{(\mbf{s}''\simp\Psi)_Z}\\
                                         \psi'\psi\mbf{s}(\mrm{D}Z)\ar[uu]^{\Psi_{\mrm{D}Z}}    \ar[r]^-{\psi'\psi\mu_{Z}}  & \psi'\psi\mbf{s}(\simp\mbf{s}Z) \ar[u]_{\Psi_{\simp\mbf{s}Z}}   \ .}\end{equation}
It is enough to check that $(\mbf{s}''\simp\Psi)_Z\comp
\Psi_{\simp\mbf{s}Z}$ agrees with
\begin{equation}\label{RistAuxiliar}\xymatrix@M=4pt@H=4pt@C=45pt@R=9pt{
  \psi'\psi\mbf{s}(\simp\mbf{s}Z) \ar[r]^-{\psi'\Theta_{\simp\mbf{s}Z}} & \psi'\mbf{s'}\psi(\simp\mbf{s}Z) \ar[r]^-{\psi'(\mbf{s}'\simp\Theta)_Z} & \psi'\mbf{s'} \simp\mbf{s'}(\psi Z)\\
  \psi'\mbf{s'} \simp\mbf{s'}(\psi Z) \ar[r]^-{\Phi_{\simp\mbf{s'}(\psi Z)}} &  \mbf{s''}\psi'\mbf{s'}(\psi Z) \ar[r]^-{(\mbf{s}''\simp\Phi)_{\psi Z}}  &  \mbf{s''} \simp\mbf{s''}(\psi'\psi Z)  }
\end{equation}
since in this case diagram (\ref{CuadrAuxCompfunt}) is the
composition of the following commutative squares
$$\xymatrix@M=4pt@H=4pt@C=40pt@R=20pt{
 \psi'\mbf{s'}\mrm{D}(\psi Z)\ar[r]^-{\psi'\mu'_{\psi Z}}                             & \psi'\mbf{s'} \simp\mbf{s'}(\psi Z)                                      & \mbf{s''}\mrm{D}(\psi'\psi Z)\ar[r]^-{\mu''_{\psi'\psi Z}}                                  &   \mbf{s''} \simp\mbf{s''}(\psi'\psi Z)                            \\
                                                                                       & \psi'\mbf{s'}\psi(\simp\mbf{s}Z)  \ar[u]_{\psi'(\mbf{s}'\simp\Theta)_Z}   &                                                                                              & \mbf{s''}\psi'\simp\mbf{s'}(\psi Z) \ar[u]_{(\mbf{s}''\simp\Phi)_{\psi Z}} \\
 \psi'\psi\mbf{s}(\mrm{D}Z)\ar[uu]^{\psi'\Theta_{\mrm{D}Z}}\ar[r]^-{\psi'\psi(\mu_{Z})} & \psi'\psi\mbf{s}(\simp\mbf{s}Z) \ar[u]_{\psi'\Theta_{\simp\mbf{s}Z}}    & \psi'\mbf{s'}\mrm{D}(\psi Z)\ar[r]^-{\psi'\mu'_{\psi(Z)}} \ar[uu]^{\Phi_{\psi\mrm{D}Z}} & \psi'\mbf{s'} \simp\mbf{s'}(\psi Z) \ar[u]_{\Phi_{\simp\mbf{s'}(\psi Z)}}     \ .}$$
By def{i}nition $(\mbf{s}''\simp\Psi)_Z\comp \Psi_{\simp\mbf{s}Z}$
is the zig-zag
$$\xymatrix@M=4pt@H=4pt@C=45pt@R=9pt{
  \psi'\psi\mbf{s}(\simp\mbf{s}Z) \ar[r]^-{\psi'\Theta_{\simp\mbf{s}Z}} & \psi'\mbf{s}'\psi (\simp\mbf{s}Z)  \ar[r]^-{\Phi_{\psi (\simp\mbf{s}Z)}} & \mbf{s''}\psi'\psi(\simp\mbf{s}Z) \\
  \mbf{s''}\psi'\psi(\simp\mbf{s}Z) \ar[r]^-{(\mbf{s''}\psi'\simp\Theta)_{Z}}  & \mbf{s''}\psi'\mbf{s}'\psi (\simp\mbf{s}Z)  \ar[r]^-{(\mbf{s''}\simp\Phi)_{Z}}  & \mbf{s''}\simp\mbf{s''}(\psi'\psi Z)  }
$$
and to see that it agrees with (\ref{RistAuxiliar}), it
suf{f}{i}ces to prove that the following diagram is commutative in
$Ho\mc{D''}$
$$\xymatrix@M=4pt@H=4pt@C=40pt@R=30pt{
  \psi'\mbf{s}'\psi (\simp\mbf{s}Z)  \ar[r]^-{\Phi_{\psi (\simp\mbf{s}Z)}} \ar[d]_{\psi'(\mbf{s}'\simp\Theta)_Z} & \mbf{s''}\psi'\psi(\simp\mbf{s}Z) \ar[d]^-{(\mbf{s''}\psi'\simp\Theta)_{Z}} \\
  \psi'\mbf{s'} \simp\mbf{s'}(\psi Z) \ar[r]^-{\Phi_{\simp\mbf{s'}(\psi Z)}}                                     &  \mbf{s''}\psi'\mbf{s'}(\psi Z) \ . }
$$
Expanding this diagram we get
$$\xymatrix@M=4pt@H=4pt@C=35pt@R=28pt{
  \psi'\mbf{s}'\psi (\simp\mbf{s}Z)  \ar@{-}[r]^-{\Phi^0_{\psi (\simp\mbf{s}Z)}}\ar@{-}[d]^{\psi'(\mbf{s}'\simp\Theta^0)_Z} & B^1(\psi (\simp\mbf{s}Z))\cdots B^{s-1}(\psi (\simp\mbf{s}Z))\ar@{-}[r]^-{\Phi^{s-1}_{\psi (\simp\mbf{s}Z)}}         & \mbf{s''}\psi'\psi(\simp\mbf{s}Z)\ar@{-}[d]_-{(\mbf{s''}\psi'\simp\Theta^0)_{Z}} \\
  \psi'\mbf{s}'(\simp A^1 Z) \ar@{.}[d]                                                   &                                                                                                                      & \mbf{s''}\psi'(\simp A^1 Z) \ar@{.}[d]\\
              \ar@{-}[d]^{\psi'(\mbf{s}'\simp\Theta^{r-1})_Z}                                                         &                                                                                                                      &  \ar@{-}[d]_-{(\mbf{s''}\psi'\simp\Theta^{r-1})_{Z}} \\
  \psi'\mbf{s'} \simp\mbf{s'}(\psi Z)\ar@{-}[r]_-{\Phi^0_{\simp\mbf{s'}(\psi Z)}}                                     & B^1(\simp\mbf{s'}(\psi Z))\cdots B^{s-1}(\simp\mbf{s'}(\psi Z))\ar@{-}[r]_-{\Phi^{s-1}_{\simp\mbf{s'}(\psi Z)}} &  \mbf{s''}\psi'\mbf{s}'\psi (\simp\mbf{s}Z)
 }
$$
Since $\Phi^j$ is a natural transformation relating $B^j$ and
$B^{j+1}$, we deduce that the two upper rows of the above square
can be completed to the following diagram, where each square
commutes
$$\xymatrix@M=4pt@H=4pt@C=35pt@R=30pt{
  \psi'\mbf{s}'\psi (\simp\mbf{s}Z)  \ar@{-}[r]^-{\Phi^0_{\psi (\simp\mbf{s}Z)}}\ar@{-}[d]_{\psi'(\mbf{s}'\simp\Theta^0)_Z} & B^1(\psi (\simp\mbf{s}Z)) \ar@{-}[d]_{B^1(\simp\Theta^0)_Z} \ar@{..}[r] & B^{s-1}(\psi (\simp\mbf{s}Z))\ar@{-}[r]^-{\Phi^{s-1}_{\psi (\simp\mbf{s}Z)}}\ar@{-}[d]_{B^{s-1}(\simp\Theta^0)_Z} & \mbf{s''}\psi'\psi(\simp\mbf{s}Z)\ar@{-}[d]_-{(\mbf{s''}\psi'\simp\Theta^0)_{Z}} \\
  \psi'\mbf{s}'(\simp A^1 Z)         \ar@{-}[r]_-{\Phi^0_{\psi (\simp A^1 Z)}}                                              & B^1(\simp A^1 Z)                                            \ar@{..}[r] & B^{s-1}(\simp A^1 Z)    \ar@{-}[r]_-{\Phi^{s-1}_{\psi (\simp A^1 Z)}}                                                    & \mbf{s''}\psi'(\simp A^1 Z) \\
 }
$$
and iterating this procedure we get that the required diagram
commutes.
\end{proof}

Next we introduce the ``transfer lemma''. To this end we need the
following remark about the commutativity between the simplicial
cylinder functor and the functor
$\simp\psi:\simp\mc{D}\rightarrow\simp\mc{D}'$ induced by a
quasi-strict monoidal functor $\psi:\mc{D}\rightarrow \mc{D}'$.

\begin{obs}\label{PsiConmCyl}
If $\psi:\mc{D}\rightarrow \mc{D}'$ satisf{i}es (FD 1), that is,
it is quasi-strict monoidal, then the diagram
$$\xymatrix@H=4pt@C=30pt@R=30pt{\Omega(\mc{D}) \ar[r]^{\simp\psi}\ar[d]_{Cyl} & \Omega(\mc{D'}) \ar[d]^{Cyl}\\
                                \simp\mc{D}\ar[r]^{\simp\psi}         &  \simp\mc{D'}}$$
commutes up to (degreewise) equivalence.\\
More concretely, there exists $\tau
:Cyl\,\simp\psi\rightarrow\simp\psi\,Cyl$ such that
$(\tau_D)_n\in\mrm{E'}$ $\forall D\in\Omega(\mc{D})$, $\forall n\geq 0$.\\
\indent Indeed, if $D\equiv \xymatrix@M=4pt@H=4pt{ X_{-1}\times\Dl
& X \ar[r]^{f} \ar[l]_-{\varepsilon} &  Y }$ then
$$\simp\psi(\mc{D})\equiv \xymatrix@M=4pt@H=4pt@C=35pt{ \psi(X_{-1})\times\Dl & \simp\psi(X)\ar[r]^-{\simp\psi(f)}\ar[l]_-{\simp\psi(\varepsilon)} & \simp\psi(Y) \ .}$$
The morphisms $(\tau_D)_n=\sigma_\psi :
\psi(Y_n)\sqcup\psi(X_{n-1})\sqcup\cdots\sqcup\psi(X_{-1})\rightarrow
\psi(Y_n \sqcup X_{n-1}\sqcup\cdots\sqcup {X_{-1}})$ are
equivalences, and $\tau_D=\{(\tau_D)_n\}_n$ is a morphism in
$\simp\mc{D'}$, due to the universal property of the coproduct.\\
Note also that by def{i}nition of
$i_{X_{-1}}:X_{-1}\times\Dl\rightarrow Cyl(D)$ and
$i_{\psi(X_{-1})}:\psi(X_{-1})\times\Dl\rightarrow
Cyl(\simp\psi(D))$, the diagram
\begin{equation}\label{PsiConmCylInclusiones}\xymatrix@M=4pt@H=4pt@C=40pt@R=40pt{\psi(X_{-1})\times\Dl \ar[r]^{\simp\psi(i_{X_{-1}})}\ar[dr]_{i_{\psi(X_{-1})}}  & \simp\psi(Cyl(D))\\
                                                                                                                         & Cyl(\simp\psi(D))\ar[u]_{\tau_D} }\end{equation}
commutes.
\end{obs}

\begin{thm}[Transfer lemma]\label{FDfuerte}\index{Index}{Transfer lemma}\mbox{}\\
Consider the data $(\mc{D},\mbf{s},\mu,\lambda)$ satisfying\\
$\mathbf{(SDC\; 1)}$ $\mc{D}$ is a category with f{i}nite
coproducts
and f{i}nal object $1$.\\
$\mathbf{(SDC\; 3)'}$ $\mbf{s}:\simp\mc{D}\rightarrow \mc{D}$ is
a functor.\\
$\mathbf{(SDC\; 4)'}$ $\mu:\mbf{s}\mrm{D} \rightarrow
\mbf{s}(\simp\mathbf{s})$ is a natural transformation.\\
$\mathbf{(SDC\; 5)'}$ $\lambda:\mathbf{s}(-\times \Dl)\rightarrow
Id_{\mc{D}}$ is a natural transformation compatible with $\mu$,
that is, the equalities in $($\ref{compatibLambdaMu}$)$ hold.\\
Assume that $(\mc{D'},\mrm{E'},\mbf{s'},\mu',\lambda')$ is a
simplicial descent category and $\psi:\mc{D}\rightarrow\mc{D'}$ a
functor such that $\mathbf{(FD\; 1)}$ and
$\mathbf{(FD\; 2)}$ hold.\\
Then, taking $\mrm{E}=\{f\; |\; \psi(f)\in \mrm{E'}\}$,
$(\mc{D},\mrm{E},\mbf{s},\mu,\lambda)$ is a simplicial descent
category.\\
In addition, $\psi:\mc{D}\rightarrow\mc{D}'$ is a functor of
simplicial descent categories.
\end{thm}
\begin{proof}
We must see that the data $(\mc{D},\mrm{E},\mbf{s},\mu,\lambda)$
satisfy the axioms of simplicial descent category.\\[0.2cm]
For clarity, assume that the functorial zig-zag
(\ref{RistraTheta}) in $\mc{D}'$ associated with
$\Theta:\gamma\comp\psi\comp\mbf{s} \rightarrow
\gamma\comp\mbf{s}'\comp\simp\psi$ is
$$\xymatrix@M=4pt@H=4pt@C=25pt{ \psi\comp\mbf{s}  & A \ar[r]^-{\Theta^1} \ar[l]_-{\Theta^0}& \mbf{s}'\comp\simp\psi\ .}$$
(SDC 2) By def{i}nition, $\mrm{E}$ is the class consisting of
those morphisms that are mapped by the composition
$$ \mc{D}  \stackrel{\psi}{\rightarrow} \mc{D}'\rightarrow Ho\mc{D}'$$
into isomorphisms. Hence $\mrm{E}$ is saturated. Let us check that $\mrm{E}$ is stable under coproducts.\\
Let $f_j:X_j\rightarrow Y_j$ be morphisms in $\mc{D}$ for $j=1,2$
such that $\psi(f_j)\in \mrm{E}'$. Since
$$\xymatrix@H=4pt@C=50pt@R=25pt{\psi(X_1\sqcup X_2)\ar[r]^{\psi(f_1\sqcup f_2)} & \psi(Y_1\sqcup Y_2)\\
                                            \psi(X_1)\sqcup\psi(X_2)\ar[r]^{\psi(f_1)\sqcup\psi(f_2)}\ar[u]^{\sigma_{\psi}} & \psi(Y_1)\sqcup\psi(Y_2)\ar[u]_{\sigma_{\psi}} }$$
commutes and $\sigma_{\psi}\in\mrm{E'}$, it follows from the 2 out
of 3 property that $\psi(f_1\sqcup
f_2)\in \mrm{E'}$, therefore $f_1\sqcup f_2 \in \mrm{E}$.\\[0.2cm]
(SDC 3) If $X,Y\in\simp\mc{D}$, we must see that
$\psi(\sigma_\mbf{s}):\psi(\mbf{s}(X)\sqcup\mbf{s}(Y))\rightarrow
\psi(\mbf{s}(X\sqcup Y))$ is in $\mrm{E}'$.\\
Consider the diagram
$$\xymatrix@M=4pt@H=4pt@C=45pt@R=45pt{ \psi(\mbf{s}X\sqcup\mbf{s}Y)\ar[r]^{\psi(\sigma_{\mbf{s}})}                                                         & \psi\mbf{s}(X\sqcup Y)\\
                                       \psi\mbf{s}X\sqcup\psi\mbf{s}Y\ar[r]^{\sigma_{\psi\mbf{s}}}  \ar[u]^{\sigma_{\psi}}                                 & \psi\mbf{s}(X\sqcup Y)\ar[u]_{Id} \\
                                       A(X)\sqcup A(Y) \ar[d]_{\Theta^1_{X}\sqcup\Theta^1_Y}\ar[r]^{\sigma_A}\ar[u]^{\Theta^0_{X}\sqcup\Theta^0_Y}         & A(X\sqcup Y)\ar[d]^{\Theta^1_{X\sqcup Y}} \ar[u]_{\Theta^0_{X\sqcup Y}}\\
                                       \mbf{s'}(\simp\psi X)\sqcup\mbf{s'}(\simp\psi Y)\ar[r]^-{\sigma_{\mbf{s'}\simp\psi}}\ar[d]_{\sigma_{\mbf{s'}}}      & \mbf{s'}\simp\psi(X\sqcup Y)\ar[d]^{Id} \\
                                       \mbf{s'}(\simp\psi X\sqcup\simp\psi Y)\ar[r]^-{\mbf{s'}(\sigma_{\simp\psi})}                                        & \mbf{s'}\simp\psi(X\sqcup Y)\ .}$$
The top and bottom squares commute in $\mc{D}'$ (because of the
universal property of the coproduct and the def{i}nition of
$\sigma$, (\ref{FuntorMonoidalCasiEstricto})). The two central
squares commute by the same reason, since every natural
transformation is monoidal with respect to the coproduct.

On the other and, the morphism $\sigma_{\simp\psi}$ is
$\sigma_{\psi}: \psi(X_n)\sqcup\psi(Y_n)\rightarrow\psi(X_n\sqcup
Y_n)$ in degree $n$, that is in $\mrm{E'}$ for all $n$. Therefore,
from the exactness of $\mbf{s'}$ we deduce that
$\mbf{s'}(\sigma_{\simp\psi})\in \mrm{E}'$, and by the 2 out of 3
property we get that $\sigma_{\mbf{s}'\simp\psi}\in\mrm{E}'$.
In addition $\Theta^i_X\sqcup\Theta^i_Y,\Theta^i_{X\sqcup
Y}\in\mrm{E}'$ for $i=0,1$, and hence
$\sigma_{\psi\mbf{s}}\in\mrm{E}'$.
Consequently $\psi(\sigma_{\mbf{s}})\in\mrm{E}'$.\\[0.2cm]
(SDC 4) Let $Z\in\simp\simp\mc{D}$. The square (\ref{compatFDi})
commutes in $Ho\mc{D}'$ and by def{i}nition the horizontal
morphisms
are isomorphisms, as well as $\mu'_{\psi(Z)}$.\\
Hence $\psi(\mu_Z)$ is an isomorphism in $Ho\mc{D}'$, so it follows that $\psi(\mu_Z)\in\mrm{E}'$, and $\mu_Z\in\mrm{E}$.\\[0.2cm]
(SDC 5) Analogously, given $X\in\mc{D}$, we deduce from the
commutativity of (\ref{compatFDii}) that
$\psi(\lambda_X)\in\mrm{E'}$, so $\lambda_X\in\mrm{E}$.\\[0.2cm]
(SDC 6) Let $f:X\rightarrow Y$ be a morphism in $\simp\mc{D}$ with
$f_n \in \mrm{E}$ $\forall n$. Then
$[\simp\psi(f)]_n=\psi(f_n)\in\mrm{E'}$ $\forall n$, and
consequently $\mbf{s'}(\simp\psi(f))\in\mrm{E'}$.\\
It follows from the naturality of $\Theta$ that
$\mathbf{s}'(\simp\psi(f))\comp\Theta_X=\Theta_Y\comp\psi(\mbf{s}(f))$
in $Ho\mc{D}'$, and hence $\psi(\mbf{s}(f))\in\mrm{E'}$, therefore
$\mbf{s}(f)\in\mrm{E}$.\\[0.2cm]
(SDC 7) Given a morphism $f:X\rightarrow Y$ in $\simp\mc{D}$ we
have to prove that $\psi\mbf{s}f\in\mrm{E'}$ if and only if
$\psi\mbf{s}(Cf)\rightarrow \psi(1)$ is so.
By (FD 2), we have that $\psi\mbf{s}f\in\mrm{E}'$ if and only if
$\mbf{s}'(\simp\psi f)\in\mrm{E}'$.
Let $\tau :Cyl\,\simp\psi\rightarrow\simp\psi\,Cyl$ be the
morphism def{i}ned as in \ref{PsiConmCyl}.\\
If $D\equiv \xymatrix@M=4pt@H=4pt{ 1\times\Dl & X \ar[r]^{f}
\ar[l]_-{\rho} &  Y }$, applying $\simp\psi$ we obtain
$$\simp\psi(\mc{D})\equiv \xymatrix@M=4pt@H=4pt@C=35pt{ \psi(1)\times\Dl & \simp\psi(X)\ar[r]^-{\simp\psi(f)}\ar[l]_-{\simp\psi(\rho)} & \simp\psi(Y) \ ,}$$
so $\tau_D:Cyl(\simp\psi f,\simp\psi\rho)\rightarrow
\simp\psi(Cyl(f,\rho))=\simp\psi (Cf)$ is a degreewise
equivalence, and we deduce that
$\mbf{s}'(\tau_D):\mbf{s}'Cyl(\simp\psi
f,\simp\psi\rho)\rightarrow \mbf{s}'\simp\psi( Cf)$ is in $\mrm{E}'$.\\
In addition, it follows from the commutativity of
(\ref{PsiConmCylInclusiones}) that
$$\mbf{s}'i_{\psi 1}\!:\!\mrm{R}\psi(1)\rightarrow\!\mbf{s}'Cyl(\simp\psi f,\simp\psi\rho)\!\in\!\mrm{E}'\mbox{ if and only if }\mbf{s}'(\simp\psi i_1)\!:\!\mrm{R}\psi(1)\rightarrow\! \mbf{s}'(\simp\psi Cf)\!\in\!\mrm{E}'$$
By (FD 2), $\mbf{s}'(\simp\psi i_1)\in\mrm{E}'$ if and only if
$\psi(\mbf{s}i_1):\psi(\mrm{R}1)\rightarrow \psi(\mbf{s}Cf)$ is so.\\
But we proved in (SDC 5) that $\psi(\mrm{R}1)\rightarrow
\psi(1)\in\mrm{E}'$, then
\begin{center}$\mbf{s}'i_{\psi 1}:\mrm{R}\psi(1)\rightarrow\mbf{s}'Cyl(\simp\psi f,\simp\psi\rho)\in\mrm{E}'$ if and only if $\psi(\mbf{s}Cf)\rightarrow \psi(1)\in\mrm{E}'$.\end{center}
On the other hand, if $\rho':Y\rightarrow 1\times\Dl$ is the
trivial morphism, we have that $\rho'\comp f=\rho$, and
$(\simp\psi f)\comp(\simp\psi\rho')=\simp\psi \rho$.
Hence, by \ref{AxCil+retracto}, it holds that
\begin{center}$\mbf{s}'(\simp\psi f)\in \mrm{E}'$ if and only if $\mbf{s}'i_{\psi(1)}:\mrm{R}\psi(1)\rightarrow \mbf{s}' Cyl(\simp\psi f,\simp\psi\rho)$ is in $\mrm{E}'$,\end{center}
so (SDC 7) is already proven.\\
(SDC 8) Given a morphism $f$ in $\simp\mc{D}$ then
$$ \mbf{s}f\in\mrm{E} \Leftrightarrow  \psi(\mbf{s}f)\in\mrm{E}' \Leftrightarrow \mbf{s}' (\simp\psi f)\in\mrm{E}' \Leftrightarrow
\mbf{s}'(\Upsilon\comp\simp\psi f)=\mbf{s}'(\simp \psi (\Upsilon
f))\in\mrm{E}' \Leftrightarrow $$
$$\Leftrightarrow \psi(\mbf{s}(\Upsilon f))\in\mrm{E}' \Leftrightarrow \mbf{s}(\Upsilon f)\in\mrm{E}\ .$$
\end{proof}

\begin{cor}
If $\mc{D}$ is a subcategory of a simplicial descent category
$(\mc{D'},\mrm{E'},\mbf{s'},\mu',\lambda')$ such that
\begin{enumerate}
 \item $\mc{D}$ is closed under coproducts
 \item $\mc{D}$ is closed under the simple functor, that is, if $X\in\simp
 \mc{D}$ then $\mbf{s'}(X)$ is in $\mc{D}$.
 \item The value of the natural transformation $\lambda'$ (resp. $\mu'$) at any object of $\mc{D}$  (resp. $\simp\simp\mc{D}$) is a morphism
 of $\mc{D}$. For instance, if $\mc{D}$ is full.
\end{enumerate}
Then
$(\mc{D},\mrm{E'}\cap\mc{D},\mbf{s'}|_{\mc{D}},\mu'|_{\mc{D}},\lambda'|_{\mc{D}})$
is a simplicial descent category.
\end{cor}

\begin{proof} It suf{f}{i}ces to take $\psi=i:\mc{D}\rightarrow \mc{D}'$ in the lemma transfer.
\end{proof}

\subsection{Associativity  of $\mu$}

If $\mc{D}$ is a (simplicial) descent category and the natural
transformation $\mu$ is ``associative'' then $\simp\mc{D}$ has a
second structure of descent category in addition to the one
introduced in \ref{SimpDcatDescGamma}, this time taking the
diagonal functor $\mrm{D}:\simp\simp\mc{D}\rightarrow \simp\mc{D}$
as simple functor.

The associativity property of $\mu$ will not be used in the
sequel. However, this is a relevant property and means that the
descent structure can be iterated in a ``suitable'' way in the
category of multisimplicial objects in $\mc{D}$.\\
Moreover, in this case the transformations $\lambda$ and $\mu$
give rise to a cotriple (cf. \cite{Dus}) in $\mc{D}$.

\begin{def{i}}[Associativity of $\mu$]\label{Def{i}asocitividadMu}\mbox{}\\
Let $\mrm{D}_{1,2}:\simp\simp\simp\mc{D}\rightarrow
\simp\simp\mc{D}$ (resp.
$\mrm{D}_{2,3}:\simp\simp\simp\mc{D}\rightarrow \simp\simp\mc{D}$)
be the functor that makes equal the two f{i}rst indexes (resp. the
two last indexes) of a trisimplicial object.\\
We will say that the natural transformation $\mu$ is
\textit{associative} if for every $T\in\simp\simp\simp\mc{D}$ the
following diagram commutes in $\mc{D}$
\begin{equation}\label{asociatividadMu}
\xymatrix@M=4pt@H=4pt@C=40pt{
 \mbf{s}\mrm{D}\mrm{D}_{1,2}T=\mbf{s}\mrm{D}\mrm{D}_{2,3}T\ar[r]^-{\mu_{\mrm{D}_{2,3}T}}\ar[d]_{\mu_{\mrm{D}_{1,2}T}} & \mbf{s}\simp\mbf{s}\mrm{D}_{2,3}T \ar[d]^{\mbf{s}(\simp\mu_T)} \\
 \mbf{s}\simp\mbf{s}\mrm{D}_{1,2}T=\mbf{s}\mrm{D}\simp\simp\mbf{s}T \ar[r]^-{\mu_{\simp\simp\mbf{s}T}}                                                 & \mbf{s}\comp\simp\mbf{s}\comp\simp\simp\mbf{s}(T) \ ,} \end{equation}
where $\simp\mu_T:\simp\mbf{s}\mrm{D}_{2,3}T\rightarrow
\simp\mbf{s}\comp\simp\simp\mbf{s}(T)$ is in degree $n$ the
morphism between simplicial objects
$$(\simp\mu_T)_n =\mu_{T_{n,\cdot,\cdot}}:\mbf{s}\mrm{D}T_{n,\cdot,\cdot}\rightarrow \mbf{s}(\simp\mbf{s}T_{n,\cdot,\cdot})\ .$$
\end{def{i}}

\begin{prop} Assume that $\mc{D}$ is a simplicial descent category with $\mu$
associative. Then the data
$$\begin{array}{l}
 \mrm{R}:\mc{D}\rightarrow\mc{D}\mbox{ where }\mrm{R}X=\mbf{s}(X\times\Dl)\\
  \lambda:\mrm{R}\rightarrow Id_{\mc{D}}\\
  \mu':\mrm{R}\rightarrow \mrm{R}^{2}\mbox{ where }\mu'_X=\mu_{X\times\Dl\times\Dl}
\end{array}$$
is a cotriple $(\mrm{R},\lambda,\mu')$ in the category $\mc{D}$.
\end{prop}

\begin{proof}
The statement is an immediate consequence of the associativity of
$\mu$ and the equations describing the compatibility between
$\lambda$ and $\mu$ (\ref{compatibLambdaMu}).\\
Indeed, assume given an object $Y$ in $\mc{D}$ and set
$X=Y\times\Dl$. Following the notations in
(\ref{compatibLambdaMu}), we  have that $\mu_{\Dl\times X\times
\Dl}=\mu_{X\times\Dl\times\Dl}=\mu'_{Y}$, whereas
$\lambda_{\mbf{s}X}=\lambda_{\mrm{R}Y}$ and
$\mbf{s}(\lambda_{X_n})=\mbf{s}(\lambda_Y\times\Dl)=\mrm{R}\lambda_Y$.
Hence, the compatibility equations between $\lambda$ and $\mu$
are in this case $\mu'_Y\comp \lambda_{\mrm{R}Y}=\mu'_Y \comp \mrm{R}\lambda_Y=Id_Y$.\\
On the other hand, if $T=Y\times\Dl\times\Dl\in\simp\simp\mc{D}$
then
$$ \mu_{\mrm{D}_{1,2}T}=\mu_{\mrm{D}_{2,3}T}=\mu'_Y \, ,\
 \mu_{\simp\simp\mbf{s}T}=\mu'_{\mrm{R}Y} \mbox{ and } \mbf{s}(\simp\mu_T)=\mbf{s}(\mu_{Y\times\Dl\times\Dl}\times\Dl)=\mrm{R}\mu'_Y$$
and the commutativity of diagram (\ref{asociatividadMu}) in this
case is just the equality $\mu'_{\mrm{R}Y}\comp
\mu'_Y=(\mrm{R}\mu'_Y)\comp\mu'_Y$.
\end{proof}

\begin{prop}\label{SimpDdeDescConDiag}
If $(\mc{D},\mrm{E},\mbf{s},\mu,\lambda)$ is a simplicial descent
category with $\mu$ associative then $\simp\mc{D}$ is also a
simplicial descent category, where the simple functor is the
diagonal functor $\mrm{D}:\simp\simp\mc{D}\rightarrow\simp\mc{D}$
and the class of equivalences is
$$\mrm{E}'_{\simp\mc{D}}=\{f\ |\ \mbf{s}f\in\mrm{E}\, \}\ .$$
In addition, $\mbf{s}:\simp\mc{D}\rightarrow \mc{D}$ is a functor
of descent categories.
\end{prop}

\begin{proof}
The result follows from the transfer lemma, setting
$\psi=\mbf{s}:\simp\mc{D}\rightarrow \mc{D}$.\\
Axioms (SDC 1) and (SDC 3$)'$ hold. The natural transformations
$\mu_{\simp\mc{D}}$ and $\lambda_{\simp\mc{D}}$
are both the identity natural transformation, therefore they satisfy trivially the equalities (\ref{compatibLambdaMu}).\\
On the other hand, (FD 1) is a consequence of the additivity
axiom, whereas by the normalization one, the natural
transformation $\Theta=\mu: \mbf{s}\mrm{D}\rightarrow
\mbf{s}\simp\mbf{s}$ is a
(pointwise) equivalence.\\
If $X$ is a simplicial object in $\mc{D}$, diagram
(\ref{compatFDii}) is just
$$\xymatrix@M=4pt@H=4pt@C=50pt@R=40pt{\mbf{s}\mrm{D}(X\times\Dl) \ar[r]^-{\mbf{s}(Id)=Id}   \ar[d]_-{\mu_{X\times\Dl}}  & \mbf{s}X\\
                                      \mbf{s}((\mbf{s}X)\times\Dl) \ar[ru]_-{\lambda_{\mbf{s}X}} &}$$
that commutes in $\mc{D}$ by the compatibility condition between $\lambda$ and $\mu$.\\
Consider now $T\in\simp\simp\simp\mc{D}$. Under this setting, the
usual diagonal functor $\mrm{D}:\simp\simp(\simp\mc{D})\rightarrow
\simp\mc{D}$ is $\mrm{D}_{1,2}$, whereas $\simp\mrm{D}:\simp\simp
\simp\mc{D}\rightarrow \simp\mc{D}$ is by def{i}nition $\mrm{D}_{2,3}$.\\
Hence, diagram (\ref{compatFDi}) can be written as
$$\xymatrix@M=4pt@H=4pt@C=40pt@R=20pt{ \mbf{s}\mrm{D}_{1,2}(\simp\simp\mbf{s} T)\ar[r]^-{\mu_{\simp\simp\mbf{s} T}} &   \mbf{s} \simp\mbf{s}(\simp\simp\mbf{s} T)   \\
                                                                                                         &   \mbf{s}\simp\mbf{s}\mrm{D}_{2,3}T  \ar[u]_{(\mbf{s}\simp\mu)_T}\\
                                       \mbf{s}\mrm{D}\mrm{D}_{1,2}T\ar[uu]^{\mu_{\mrm{D}_{1,2}T}}    \ar[r]^{\mbf{s}(Id)=Id}  & \mbf{s}\mrm{D}\mrm{D}_{2,3}T \ar[u]_{\mu_{\mrm{D}_{2,3}T}}   }
$$
so its commutativity is just the associativity condition satisfied
by $\mu$.
\end{proof}

%

\chapter{The homotopy category of a simplicial descent category}\label{CapituloLocalizada}

\section{Description of $Ho\mc{D}$}
This section is devoted to the study of the homotopy category
associated with a simplicial category $\mc{D}$, that is by
def{i}nition $\mc{D}[\mrm{E}^{-1}]$.\\
In general the class $\mrm{E}$ does not has calculus of fractions,
for instance when $\mc{D}=$chain complexes and $\mrm{E}=$morphisms
inducing isomorphism in homology. However, some of the properties
satisf{i}ed by the functor $cyl$ developed in the last chapter are
similar, but in a more general sense, to the left calculus of
fractions (or to the right calculus of fractions in the
cosimplicial case). This fact will allow us to exhibit a
``reasonable'' description of the morphisms in $Ho\mc{D}$. From
now on $\mc{D}$ will be a simplicial descent category.

\begin{num}
Let $X$ be an object in $\mc{D}$. We remind that $\mrm{R}:\mc{D}\rightarrow\mc{D}$ is $\mrm{R}X=\mbf{s}(X\times\Dl)$.\\
In addition, if $T=X\times\Dl\in\simp\simp\mc{D}$ then
$\mu_{T}:\mrm{R}X\rightarrow \mrm{R}^2 X$. Denote also by
$\mu:\mrm{R}\rightarrow \mrm{R}^{2}$ the natural transformation
obtained in this way, that is, $\mu_X$ means $\mu_{X\times\Dl\times\Dl}$.\\
From the compatibility between $\lambda$ and $\mu$
(\ref{compatibLambdaMu}) we deduce that the following compositions
must be the identity in $\mc{D}$
\begin{equation}\label{comptLambdaMuCtes}\xymatrix@M=4pt@H=4pt@R=8pt{\mrm{R}X \ar[r]^{\mu_X} &\mrm{R}^2 X \ar[r]^{\lambda_{\mrm{R}X}}& \mrm{R}X\\
                                          \mrm{R}X \ar[r]^{\mu_X} &\mrm{R}^2 X \ar[r]^{\mrm{R}\lambda_{X}}& \mrm{R}X \ .}
\end{equation}
Note also that from the naturality of $\lambda:\mrm{R}\rightarrow
Id_{\mc{D}}$ it follows that
$\lambda_X\comp\lambda_{\mrm{R}X}=\lambda_X\comp\mrm{R}\lambda_X$.
\end{num}

\begin{def{i}}\label{Def{i}HoD}\index{Index}{description of $Ho\mc{D}$}
Let $Ho\mc{D}$ be the category with the same objects as $\mc{D}$ and whose morphisms are described as follows.\\
Given objects $X$, $Y$ in $\mc{D}$ then
$$ Hom_{Ho\mc{D}(X,Y)}=   \begin{array}{cc}\mrm{T}(X,Y) &                                                           \\[-0.5cm]
                                                       &\hspace{-0.2cm} \mbox{\large{$\diagup$}} \!\! \sim \end{array}$$
where an element $F$ of
$\mrm{T}(X,Y)$\index{Symbols}{$\mrm{T}(X,Y)$} is a zig-zag
$$\xymatrix@M=4pt@H=4pt{X & \mrm{R}X \ar[l]_{\lambda_X}\ar[r]^{{f}} & T & \mrm{R}Y \ar[l]_w\ar[r]^{\lambda_{Y}}& Y} ,\  w\in\mrm{E} \ .$$
If $\xymatrix@M=4pt@H=4pt{X & \mrm{R}X
\ar[l]_{\lambda_X}\ar[r]^{{g}} & S & \mrm{R}Y
\ar[l]_u\ar[r]^{\lambda_{Y}}& Y}$ is another element
$G\in\mrm{T}(X,Y)$, then $F$ is related to $G$, $F\sim G$, if and
only there exists a `hammock' (that is, a commutative diagram in
$\mc{D}$)
\begin{equation}\label{hamaca}\xymatrix@M=4pt@H=4pt{
              &\mrm{R}^2 X \ar[ld]_{Id} \ar[r]^{\mrm{R}{f}}\ar[d]& \mrm{R}T\ar[d] & \mrm{R}^2 Y\ar[l]_{\mrm{R}w}\ar[rd]^{Id}\ar[d]  &         \\
\mrm{R}^2 X   & \widetilde{X}       \ar[l]\ar[r]^{h}      & U       &   \widetilde{Y}  \ar[l]\ar[r]         & \mrm{R}^2 Y \ ,\\
              & \mrm{R}^2 X\ar[lu]^{Id}\ar[r]^{\mrm{R}{g}}\ar[u]& \mrm{R}S\ar[u] & \mrm{R}^2Y\ar[l]_{\mrm{R}u}\ar[ru]_{Id}\ar[u]  & }\end{equation}
relating $F$ and $G$, such that all maps except $f$, $g$ and $h$
are equivalences.
\end{def{i}}

\begin{num} Given two composable morphisms in $Ho\mc{D}$
represented by zig-zags  $F$ and $G$ given respectively by
$$\xymatrix@M=4pt@H=4pt@R=8pt{X & \mrm{R}X \ar[l]_{\lambda_X}\ar[r]^{{f}} & T & \mrm{R}Y \ar[l]_u\ar[r]^{\lambda_{Y}}& Y\\
                              Y & \mrm{R}Y \ar[l]_{\lambda_Y}\ar[r]^{{g}} & S & \mrm{R}Z \ar[l]_v\ar[r]^{\lambda_{Z}}& Z\ ,}$$
then their composition is represented by the zig-zag $G\comp F$
def{i}ned as
$$\xymatrix@M=4pt@H=4pt{X& \mrm{R}X \ar[l]_{\lambda_X}\ar[r]^-{h} & cyl(u,g) & \mrm{R}Z \ar[l]_-{w} \ar[r]^{\lambda_{Z}}& Z \ ,}$$
where the morphisms $h:\mrm{R}X\rightarrow  cyl(u,g)$ and
$w:\mrm{R}Z \rightarrow cyl(u,g)\in\mrm{E}$ are the respective
compositions
$$\xymatrix@M=4pt@H=4pt@R=8pt{
 \mrm{R}X \ar[r]^-{\mu_X} & \mrm{R}^2 X \ar[r]^-{\mrm{R}f} & \mrm{R}T \ar[r]^-{I_T} & cyl(u,g)\\
 \mrm{R}Z \ar[r]^-{\mu_Z} & \mrm{R}^2 Z \ar[r]^-{\mrm{R}v} & \mrm{R}S \ar[r]^-{I_S} & cyl(u,g)\ .}$$
By \ref{axiomacylenD}, $I_S\in\mrm{E}$. Therefore $w\in\mrm{E}$
since it is the composition of two equivalences.
\end{num}

\begin{obs} Note that if we compose the hammock (\ref{hamaca}) with $\lambda$ we get the following hammock
where the upper zig-zag is $F$ and the lower one is $G$
\begin{equation}\label{diamante}
\xymatrix@M=4pt@H=4pt{
              &\mrm{R} X  \ar[r]^{{f}}   \ar[ldd]_{\lambda_X}                   & T                                  & \mrm{R}Y\ar[l]_{w}\ar[rdd]^{\lambda_Y}                                     &         \\
              &\mrm{R}^2 X \ar[ld]\ar[u]_{\lambda_{\mrm{R}X}} \ar[r]^{\mrm{R}{f}}\ar[d]& \mrm{R}T\ar[d] \ar[u]^{\lambda_{T}}& \mrm{R}^2 Y\ar[l]_{\mrm{R}w}\ar[rd]\ar[d] \ar[u]^{\lambda_{\mrm{R}Y}} &         \\
        X     & \widetilde{X}       \ar[l]\ar[r]^{h}                                        & U                                  & \widetilde{Y}  \ar[l]\ar[r]                                                &  Y \ .\\
              & \mrm{R}^2 X\ar[lu]\ar[r]^{\mrm{R}{g}}\ar[u]\ar[d]^{\lambda_{\mrm{R}X}} & \mrm{R} S\ar[u]\ar[d]_{\lambda_S}  & \mrm{R}^2Y\ar[l]_{\mrm{R}u}\ar[ru]\ar[u]\ar[d]_{\lambda_{\mrm{R}Y}}   &\\
              & \mrm{R} X\ar[r]^{{g}}\ar[luu]^{\lambda_X}                                   &  S                                 & \mrm{R}Y\ar[l]_{u}\ar[ruu]_{\lambda_Y}  &}
\end{equation}
\end{obs}

\begin{thm}\label{ProHoDcat}
$Ho\mc{D}$ is in fact a category. Moreover, the functor
$$\gamma: \mc{D} \rightarrow Ho\mc{D} $$
def{i}ned as the identity over objects, and over morphisms as
$$\gamma(X\stackrel{f}{\rightarrow}Y) = \xymatrix@M=4pt@H=4pt{X & \mrm{R}X\ar[l]_{\lambda_X} \ar[r]^{\mrm{R}f} & \mrm{Y} &\mrm{R}Y\ar[l]_{Id} \ar[r]^{\lambda_Y} & Y}$$
is a localization of $\mc{D}$ with respect to $\mrm{E}$.
\end{thm}

The proof of the above theorem is very similar to the analogue
proof in the calculus of fractions case, except that now a great
number of technical problems must be solved. This is why we decide
to divide this proof into the following lemmas and preliminary
results.

\begin{lema}\label{LemaHamacaconRk}
Two elements $F$ and $G$ in $\mrm{T}(X,Y)$ given by
$$\xymatrix@M=4pt@H=4pt@R=8pt{X & \mrm{R}X \ar[l]_{\lambda_X}\ar[r]^{{f}} & T & \mrm{R}Y \ar[l]_w\ar[r]^{\lambda_{Y}}& Y\\
                              X & \mrm{R}X \ar[l]_{\lambda_X}\ar[r]^{{g}} & S & \mrm{R}Y \ar[l]_u\ar[r]^{\lambda_{Y}}& Y}$$
are such that $F\sim G$ if and only if there exists $k\geq 0$ and
a hammock $\mc{H}^k$
$$\xymatrix@M=4pt@H=4pt{
              &\mrm{R}^{k+1}  X \ar[ld]_{Id} \ar[r]^{\mrm{R}^k{f}}\ar[d]& \mrm{R}^kT\ar[d] & \mrm{R}^{k+1} Y\ar[l]_{\mrm{R}^k w}\ar[rd]^{Id}\ar[d]  &         \\
\mrm{R}^{k+1} X   & \widetilde{X}       \ar[l]\ar[r]^{h}      & U       &   \widetilde{Y}  \ar[l]\ar[r]         & \mrm{R}^{k+1} Y \ ,\\
              & \mrm{R}^{k+1} X\ar[lu]^{Id}\ar[r]^{\mrm{R}^k{g}}\ar[u]& \mrm{R}^kS\ar[u] & \mrm{R}^{k+1}Y\ar[l]_{\mrm{R}^ku}\ar[ru]_{Id}\ar[u]  & }$$
where all maps except $f$, $g$ and $h$ are equivalences.
\end{lema}

\begin{proof} If $k=0$, it is enough to apply $\mrm{R}$ to the hammock $\mc{H}^0$, since $\mrm{R}(\mrm{E})\subseteq \mrm{E}$.\\
If $k> 1$, every hammock $H^k$ gives rise to a new one with $k=1$,
through the natural transformations
$\lambda^{k}:\mrm{R}^k\rightarrow \mrm{R}$ and
$\mu^k:\mrm{R}\rightarrow \mrm{R}^k$.\\
More specif{i}cally, let $\lambda^{n}$ be the natural
transformation def{i}ned as
$\lambda_X=\lambda_{\mrm{R}^{k-1}X}\comp
\lambda_{\mrm{R}^{k-1}X}\comp \cdots\comp
\lambda_{\mrm{R}X}:\mrm{R}^kX\rightarrow \mrm{R}X$, as well as
$\mu^{k}_X=\mu_{\mrm{R}^{k-2}X}\comp
\cdots\comp\mu_X:\mrm{R}\rightarrow \mrm{R}^k$.\\
Then, from (\ref{comptLambdaMuCtes}) we deduce that the
composition
$$\xymatrix@M=4pt@H=4pt{\mrm{R}^k X\ar[r]^{\mu^{k}_X} & \mrm{R} X \ar[r]^{\lambda^{k}_X} & \mrm{R}^k X }$$
is equal to the identity.\\
Let us do the computation for the upper half of $\mc{H}^k$, since
it can be argued similarly for the lower one. To this end, just
note that the naturality of $\mu^k$ together with the equality
$\lambda^{k}\comp\mu^k=Id$ imply that the following diagram is
commutative
$$\xymatrix@M=4pt@H=4pt{
             &                                                   & \mrm{R}^{2}  X \ar[lldd]_{Id} \ar[r]^{\mrm{R}{f}}\ar[d]_{\mu^k_{\mrm{R}X}}  & \mrm{R}T  \ar[d]_{\mu^k_T} & \mrm{R}^{2} Y\ar[l]_{\mrm{R} w}\ar[rrdd]^{Id}\ar[d]_{\mu^k_{\mrm{R}Y}}  &         \\
             &                                                   & \mrm{R}^{k+1}  X \ar[ld]_{Id} \ar[r]^{\mrm{R}^k{f}}\ar[d]                   & \mrm{R}^kT\ar[d]           & \mrm{R}^{k+1} Y\ar[l]_{\mrm{R}^k w}\ar[rd]^{Id}\ar[d]                 &         \\
 \mrm{R}^2 X & \mrm{R}^{k+1} X \ar[l]_{\lambda^{k}_{\mrm{R}X}} & \widetilde{X} \ar[l]\ar[r]^{h}                                              & U                          & \widetilde{Y}  \ar[l]\ar[r]                                           & \mrm{R}^{k+1} Y \ar[r]^{\lambda^{k}_{\mrm{R}Y}} & \mrm{R}^{2} Y \ .}$$
\end{proof}

\begin{lema}
The relation $\sim$ is an equivalence relation over $T(X,Y)$.
\end{lema}

\begin{proof}
The relation $\sim$ is symmetric by def{i}nition. The reflexivity
is
also clear, just take the vertical morphisms in the corresponding hammock as identities.\\
It remains to check that $\sim$ is transitive. Assume that
 $F\sim G$ and $G\sim L$ through the hammocks
$\mathcal{H}$ and $\mathcal{H}'$ given respectively by
$$\xymatrix@M=4pt@H=4pt@C=15pt{
              & \mrm{R}^2 X\ar[ld]_{Id} \ar[r]^{\mrm{R}f}\ar[d]              & \mrm{R}T\ar[d]       & \mrm{R}^2 Y\ar[l]_{\mrm{R}w}\ar[rd]^{Id}\ar[d]       &                                      &                          &\mrm{R}^2 X\ar[ld]_{Id} \ar[r]^{\mrm{R}{g}}\ar[d]_{q} & \mrm{R}S\ar[d]_{q''} & \mrm{R}^2 Y\ar[l]_{\mrm{R}w'}\ar[rd]^{Id}\ar[d]_{q'}  &\\
\mrm{R}^2 X   & \widetilde{X}       \ar[l]_{\alpha}\ar[r]^{h}                         & U                    &  \widetilde{Y}  \ar[l]_{u}\ar[r]^{\beta}    & \mrm{R}^2Y\!\!\!\! \ar@{}[r]|{;}     & \!\!\!\! \mrm{R}^2 X & \widehat{X}       \ar[l]_{\alpha'}\ar[r]^{h'}                  & W                    &   \widehat{Y}  \ar[l]_{u'}\ar[r]^{\beta'}    & \mrm{R}^2 Y \\
              & \mrm{R}^2 X\ar[lu]^{Id}\ar[r]^{\mrm{R}g}\ar[u]^{p}           & \mrm{R}S\ar[u]^{p''} & \mrm{R}^2 Y\ar[l]_{\mrm{R}w'}\ar[ru]_{Id}\ar[u]^{p'} &                                      &                          & \mrm{R}^2 X\ar[lu]^{Id}\ar[r]^{\mrm{R}{l}}\ar[u]     & \mrm{R}V\ar[u]       & \mrm{R}^2 Y\ar[l]_{\mrm{R}w''}\ar[ru]_{Id}\ar[u]  &
}$$
Applying by columns the functor $cyl$ we get
$$
\xymatrix@M=4pt@H=4pt{  \widetilde{X}  \ar[r]^{h}                       & U                                  & \widetilde{Y}  \ar[l]_u\\
                        \mrm{R}^2 X\ar[r]^{\mrm{R}{g}}\ar[d]_q \ar[u]^p & \mrm{R} S\ar[u]^{p''} \ar[d]_{q''} & \mrm{R}^2 Y\ar[l]_{\mrm{R}w'}\ar[u]^{p'}\ar[d]_{q'}  \\
                        \widehat{X}    \ar[r]^{h'}                      & W                                  & \widehat{Y} \ar[l]_{u'} }$$
and setting $\overline{X}=cyl(q,p)$, $M=cyl(q'',p'')$ and
$\overline{Y}=cyl(q',p')$ we obtain the commutative diagram in
$\mc{D}$
\begin{equation}\label{diagrinteriorhamaca}
\xymatrix@M=4pt@H=4pt{  \mrm{R}\widetilde{X}  \ar[r]^{\mrm{R}h} \ar[d]_{s}  & \mrm{R}U  \ar[d]     &  \mrm{R} \widetilde{Y} \ar[d]_{s'} \ar[l]_{\mrm{R}u}\\
                        \overline{X}          \ar[r]                    & M                    & \overline{Y}\ar[l]  \\
                        \mrm{R}\widehat{X}    \ar[r]^{\mrm{R}h'} \ar[u]^{t} & \mrm{R} W \ar[u]     &  \mrm{R} \widehat{Y}\ar[u]^{t'} \ar[l]_{\mrm{R}u'} }
\end{equation}
where all vertical arrows are equivalences by \ref{axiomacylenD},
as well as $\overline{Y}\rightarrow M$
by \ref{exactitudcyl}.\\
In addition, since $\alpha\comp p=\alpha'\comp q=Id_{\mrm{R}^2X}$
(resp. $\beta\comp p'=\beta'\comp q'=Id_{\mrm{R}^2Y}$), if follows
from \ref{cilaumentado} the existence of a morphism
$\rho:\overline{X}\rightarrow \mrm{R}^3 X$ (resp.
$\rho':\overline{Y}\rightarrow \mrm{R}^3 Y$) such that $\rho\comp
s=\mrm{R}\alpha$ and $\rho\comp t= \mrm{R}\alpha'$ (resp.
$\rho'\comp s'=\mrm{R}\beta$ and $\rho'\comp t'=
\mrm{R}\beta'$).\\
By the 2 out of 3 property we have that $\rho$ (resp. $\rho'$) is
an equivalence.\\
On the other hand, applying $\mrm{R}$ to the upper half of
$\mc{H}$ and to the lower half of $\mathcal{H}'$ we obtain
$$
\xymatrix@M=4pt@H=4pt@C=15pt{
               & \mrm{R}^3 X\ar[ld]_{Id} \ar[r]^{\mrm{R}^2f}\ar[d]            & \mrm{R}^2T\ar[d] & \mrm{R}^3 Y\ar[l]_{\mrm{R}^2w}\ar[rd]^{Id}\ar[d]            & \mrm{R}^3 X   & \mrm{R}\widehat{X}\ar[l]_{\mrm{R}\alpha'}\ar[r]^{\mrm{R}h'} & \mrm{R}W         & \mrm{R} \widehat{Y} \ar[l]_{\mrm{R}u'}\ar[r]^{\mrm{R}\beta'} & \mrm{R}^3 Y \\
 \mrm{R}^3 X   & \mrm{R}\widetilde{X} \ar[l]_{\mrm{R}\alpha}\ar[r]^{\mrm{R}h} & \mrm{R}U         &\mrm{R}\widetilde{Y} \ar[l]_{\mrm{R}u}\ar[r]^{\mrm{R}\beta} &\mrm{R}^3Y    & \mrm{R}^3 X\ar[lu]^{Id}\ar[r]^{\mrm{R}^2{l}}\ar[u]&\mrm{R}^2V\ar[u] &\mrm{R}^3Y\ar[l]_{\mrm{R}^2w''}\ar[ru]_{Id}\ar[u] & }$$
and after adjoining them to (\ref{diagrinteriorhamaca}) and
composing, the result is
$$
\xymatrix@M=4pt@H=4pt{
                & \mrm{R}^3 X\ar[ld]_{Id} \ar[r]^{\mrm{R}^2f}\ar[d] & \mrm{R}^2T\ar[d] & \mrm{R}^3 Y \ar[l]_{\mrm{R}^2w}\ar[rd]^{Id}\ar[d] & \\
 \mrm{R}^3 X    & \overline{X} \ar[l]_{\rho}\ar[r]                  & M                & \overline{Y} \ar[l] \ar[r]^{\rho'}                & \mrm{R}^3Y \\
                & \mrm{R}^3 X\ar[lu]^{Id}\ar[r]^{\mrm{R}^2{l}}\ar[u]&\mrm{R}^2V\ar[u]  & \mrm{R}^3Y\ar[l]_{\mrm{R}^2w''}\ar[ru]_{Id}\ar[u] &
}$$
Then, by the previous lemma, $F\sim L$.
\end{proof}

\begin{lema}
The composition of morphisms in $Ho\mc{D}$ is well def{i}ned, that
is, if $F\sim F'$ and $G\sim G'$ represent two composable
morphisms in $Ho\mc{D}$ then $G\comp F \sim G'\comp F'$.
\end{lema}

\begin{proof}
Assume that $F\sim F'$ and $G\sim G'$ through the hammocks
$\mathcal{H}$ and $\mathcal{H}'$ given respectively by
$$\xymatrix@M=4pt@H=4pt@C=15pt{
              & \mrm{R}^2 X\ar[ld]_{Id} \ar[r]^{\mrm{R}f}\ar[d]_{s}& \mrm{R}T\ar[d]_{s''} & \mrm{R}^2 Y\ar[l]_{\mrm{R}w}\ar[rd]^{Id}\ar[d]_{s'}  &                                      &                      &\mrm{R}^2 Y\ar[ld]_{Id} \ar[r]^{\mrm{R}{g}}\ar[d]_{q}  & \mrm{R}U\ar[d]_{q''}  & \mrm{R}^2 Z\ar[l]_{\mrm{R}v}\ar[rd]^{Id}\ar[d]_{q'}  &\\
\mrm{R}^2 X   & \widetilde{X}       \ar[l]_{\alpha}\ar[r]^{h}      & L                    &  \widetilde{Y}  \ar[l]_{u}\ar[r]^{\beta}             & \mrm{R}^2Y\!\!\!\! \ar@{}[r]|{;}     & \!\!\!\! \mrm{R}^2 Y & \widehat{Y}       \ar[l]_{\alpha'}\ar[r]^{h'}         & W                     &   \widehat{Z}  \ar[l]_{u'}\ar[r]^{\beta'}            & \mrm{R}^2 Z \\
              & \mrm{R}^2 X\ar[lu]^{Id}\ar[r]^{\mrm{R}f'}\ar[u]^{p}& \mrm{R}S\ar[u]^{p''} & \mrm{R}^2 Y\ar[l]_{\mrm{R}w'}\ar[ru]_{Id}\ar[u]^{p'} &                                      &                      & \mrm{R}^2 Y\ar[lu]^{Id}\ar[r]^{\mrm{R}{g'}}\ar[u]^{t} & \mrm{R}V\ar[u]^{t''}  & \mrm{R}^2 Z\ar[l]_{\mrm{R}v'}\ar[ru]_{Id}\ar[u]^{t'}&
}$$
We have to f{i}nd a hammock relating the zig-zags
$$\xymatrix@M=4pt@H=4pt@R=8pt@C=15pt{
 X & \mrm{R}X \ar[l]_{\lambda_X} \ar[r]^{\mu_X} & \mrm{R}^2 X \ar[r]^{\mrm{R}f} & \mrm{R}T \ar[r]^-{I_T} & cyl(w,g) & \mrm{R}U \ar[l]_-{I_U} &  \mrm{R}^2 Z \ar[l]_{\mrm{R}v} & \mrm{R}Z \ar[l]_{\mu_Z} \ar[r]^{\lambda_Z} & Z \\
 X & \mrm{R}X \ar[l]_{\lambda_X} \ar[r]^{\mu_X} & \mrm{R}^2 X \ar[r]^{\mrm{R}f'} & \mrm{R}S \ar[r]^-{I_S} & cyl(w',g') & \mrm{R}V \ar[l]_-{I_V} &  \mrm{R}^2 Z \ar[l]_{\mrm{R}v'} & \mrm{R}Z \ar[l]_{\mu_Z} \ar[r]^{\lambda_Z} & Z  \ .}$$
To this end, apply by rows the functor $cyl$ to the diagram
$$
\xymatrix@M=4pt@H=4pt{\mrm{R}T\ar[d]_{s''}   & \mrm{R}^2 Y\ar[l]_{\mrm{R}w}\ar[r]^{Id}\ar[d]           & \mrm{R}^2Y \ar[d]_{Id}  \\
                      L                    &  \widetilde{Y}  \ar[l]_{u}   \ar[r]^{\beta}             & \mrm{R}^2Y              \\
                      \mrm{R}S\ar[u]^{p''} & \mrm{R}^2 Y\ar[l]_{\mrm{R}w'}\ar[r]^{Id}\ar[u]^{p'}     & \mrm{R}^2Y \ar[u]^{Id} }$$
Then we obtain
$$
\xymatrix@M=4pt@H=4pt{  \mrm{R}^2T \ar[r]^-{J_T} \ar[d]_{\mrm{R}s''} & cyl(\mrm{R}w,Id)  \ar[d]     &  \mrm{R}^3Y \ar[d]_{Id}  \ar[l]_-{J_Y} \\
                        \mrm{R}L         \ar[r]                   & cyl(u,\beta)                 &  \mrm{R}^3Y  \ar[l]                \\
                        \mrm{R}^2S\ar[u]^{\mrm{R}p''}\ar[r]^-{J_S} & cyl(\mrm{R}w',Id)\ar[u]      &  \mrm{R}^3Y \ar[u]^{Id}  \ar[l]_-{I_Y} }
$$
where all arrows are equivalences, by properties \ref{axiomacylenD} and \ref{exactitudcyl} of $cyl$.\\
On the other hand, since $Id_{\mrm{R}T}\comp \mrm{R}w
=\mrm{R}w\comp Id_{\mrm{R}^2Y}$ (resp. $Id_{\mrm{R}S}\comp
\mrm{R}w' =\mrm{R}w'\comp Id_{\mrm{R}^2Y}$), it follows from
\ref{cilaumentado} the existence of a morphism
$\rho:cyl(\mrm{R}w,Id)\rightarrow \mrm{R}^2 T$ (resp.
$\rho':cyl(\mrm{R}w',Id)\rightarrow \mrm{R}^2 S$) such that
$\rho\comp J_T=Id_{\mrm{R}^2T}$ and $\rho\comp J_Y= \mrm{R}^2w$
(resp. $\rho'\comp
J_S=Id_{\mrm{R}^2S}$ and $\rho'\comp I_Y= \mrm{R}^2w'$).\\
From the 2 out of 3 property we deduce that $\rho$ and $\rho'$ are
equivalences. So we can construct the following diagram from these
two maps, together with the result of applying $\mrm{R}$ to some
parts of $\mc{H}$ and $\mc{H}'$, getting
$$
\xymatrix@M=4pt@H=4pt@C=14pt@R=16pt{
 \mrm{R}^3 X\ar[r]^{\mrm{R}^2f}                              & \mrm{R}^2T \ar[r]^-{Id}                               &  \mrm{R}^2T                           &  \mrm{R}^3Y   \ar[l]_-{{\mrm{R}^2w}}              & \mrm{R}^3 Y  \ar[r]^{\mrm{R}^2{g}}\ar[l]_{Id}                              & \mrm{R}^2U \\
 \mrm{R}^3 X\ar[u]^{Id} \ar[r]^{\mrm{R}^2f}\ar[d]_{\mrm{R}s} & \mrm{R}^2T \ar[r]^-{J_T} \ar[d]^{\mrm{R}s} \ar[u]_{Id}& cyl(\mrm{R}w,Id) \ar[u]^\rho \ar[d]   &  \mrm{R}^3Y \ar[d]_{Id}\ar[l]_-{J_Y} \ar[u]^{Id}  & \mrm{R}^3 Y \ar[u]^{Id} \ar[r]^{\mrm{R}^2{g}}\ar[d]_{\mrm{R}q}\ar[l]_{Id}  & \mrm{R}^2U\ar[d]_{\mrm{R}q''}\ar[u]^{Id}  \\
 \mrm{R}\widetilde{X} \ar[r]^{\mrm{R}h}                      & \mrm{R}L         \ar[r]                               & cyl(u,\beta)                          &  \mrm{R}^3Y  \ar[l]                               & \mrm{R}\widehat{Y} \ar[l]_{\mrm{R}\alpha'}\ar[r]^{\mrm{R}h'}               & \mrm{R}W                     \\
 \mrm{R}^3 X\ar[d]_{Id}\ar[r]^{\mrm{R}^2f'}\ar[u]^{\mrm{R}p} & \mrm{R}^2S\ar[u]_{\mrm{R}p''}\ar[r]^-{J_S} \ar[d]^{Id}& cyl(\mrm{R}w',Id)\ar[u] \ar[d]_{\rho'}&  \mrm{R}^3Y \ar[u]^{Id} \ar[d]_{Id} \ar[l]_-{I_Y} & \mrm{R}^3 Y \ar[d]_{Id}\ar[r]^{\mrm{R}^2{g'}}\ar[u]^{\mrm{R}t}\ar[l]_{Id}  & \mrm{R}^2V\ar[u]^{\mrm{R}t''} \ar[d]_{Id} \\
 \mrm{R}^3 X       \ar[r]^{\mrm{R}^2f'}                      &  \mrm{R}^2S\ar[r]^-{Id}                                 & \mrm{R}^2S                           &  \mrm{R}^3Y  \ar[l]_-{\mrm{R}^2w'}                & \mrm{R}^3 Y \ar[r]^{\mrm{R}^2{g'}} \ar[l]_{Id}                             & \mrm{R}^2V \ .}
$$
Composing arrows in the above diagram and attaching the remaining
part of $\mc{H}'$ we get
$$
\xymatrix@M=4pt@H=4pt@C=14pt@R=14pt{
 \mrm{R}^3 X\ar[r]^{\mrm{R}^2f}                  &  \mrm{R}^2T                           &   \mrm{R}^3 Y  \ar[r]^{\mrm{R}^2{g}}\ar[l]_{\mrm{R}^2w}                 & \mrm{R}^2U                               & \mrm{R}^3 Z\ar[l]_{\mrm{R}^2v}                               \\
 \mrm{R}^3 X\ar[u]^{Id} \ar[r] \ar[d]_{\mrm{R}s} & cyl(\mrm{R}w,Id) \ar[u]^\rho \ar[d]   &   \mrm{R}^3 Y \ar[u]^{Id} \ar[r]^{\mrm{R}^2{g}}\ar[d]_{\mrm{R}q} \ar[l] & \mrm{R}^2U\ar[d]_{\mrm{R}q''}\ar[u]^{Id} & \mrm{R}^3 Z\ar[l]_{\mrm{R}^2v}\ar[u]^{Id}\ar[d]_{\mrm{R}q'}  \\
 \mrm{R}\widetilde{X}  \ar[r]                    & cyl(u,\beta)                          &   \mrm{R}\widehat{Y} \ar[r]^{\mrm{R}h'}  \ar[l]                         & \mrm{R}W                                 & \mrm{R}\widehat{Z}  \ar[l]_{\mrm{R}u'} \\
 \mrm{R}^3 X\ar[d]_{Id}\ar[u]^{\mrm{R}p} \ar[r]  & cyl(\mrm{R}w',Id)\ar[u] \ar[d]_{\rho'}&   \mrm{R}^3 Y \ar[d]_{Id}\ar[r]^{\mrm{R}^2{g'}}\ar[u]^{\mrm{R}t}\ar[l]  & \mrm{R}^2V\ar[u]^{\mrm{R}t''} \ar[d]_{Id}& \mrm{R}^3 Z\ar[l]_{\mrm{R}^2v'}\ar[d]_{Id}\ar[u]^{\mrm{R}t'} \\
 \mrm{R}^3 X       \ar[r]^{\mrm{R}^2f'}          & \mrm{R}^2S                            &   \mrm{R}^3 Y \ar[r]^{\mrm{R}^2{g'}} \ar[l]_{\mrm{R}^2w'}               & \mrm{R}^2V                               & \mrm{R}^3 Z\ar[l]_{\mrm{R}^2v'}                                }
$$
where all vertical maps are equivalences, as well as the columns
that contains $\mrm{R}^2w$ and $\mrm{R}^2v$. Now compose with the
natural transformation $\lambda^2:\mrm{R}^2 \rightarrow Id$ to
obtain
$$
\xymatrix@M=4pt@H=4pt@C=14pt@R=14pt{
 \mrm{R} X\ar[r]^{f}                                             & T                                         & \mrm{R} Y  \ar[r]^{{g}}\ar[l]_{w}                                                       & U                                                   & \mrm{R} Z\ar[l]_{v}                               \\
 \mrm{R}^3 X\ar[u]^{\lambda^2_{\mrm{R}X}}\ar[r]\ar[d]_{\mrm{R}s} & cyl(\mrm{R}w,Id) \ar[u]^{\varrho} \ar[d]  & \mrm{R}^3 Y \ar[u]^{\lambda^2_{\mrm{R}Y}} \ar[r]^{\mrm{R}^2{g}}\ar[d]_{\mrm{R}q} \ar[l] & \mrm{R}^2U\ar[d]_{\mrm{R}q''}\ar[u]^{\lambda^2_{U}} & \mrm{R}^3 Z\ar[l]_{\mrm{R}^2v}\ar[u]^{\lambda^2_{\mrm{R} Z}}\ar[d]_{\mrm{R}q'}  \\
 \mrm{R}\widetilde{X}  \ar[r]                                    & cyl(u,\beta)                              & \mrm{R}\widehat{Y} \ar[r]^{\mrm{R}h'}  \ar[l]                                           & \mrm{R}W                                            & \mrm{R}\widehat{Z}  \ar[l]_{\mrm{R}u'} \\
 \mrm{R}^3 X\ar[d]_{\lambda^2_{\mrm{R}X}}\ar[u]^{\mrm{R}p} \ar[r]& cyl(\mrm{R}w',Id)\ar[u] \ar[d]_{\varrho'} & \mrm{R}^3 Y \ar[d]_{\lambda^2_{\mrm{R}Y}}\ar[r]^{\mrm{R}^2{g'}}\ar[u]^{\mrm{R}t}\ar[l]  & \mrm{R}^2V\ar[u]^{\mrm{R}t''} \ar[d]_{\lambda^2_{V}}& \mrm{R}^3 Z\ar[l]_{\mrm{R}^2v'}\ar[d]_{\lambda^2_{\mrm{R} Z}}\ar[u]^{\mrm{R}t'} \\
 \mrm{R} X       \ar[r]^{f'}                                     & S                                         & \mrm{R} Y \ar[r]^{{g'}} \ar[l]_{w'}                                                     & V                                                   & \mrm{R} Z \ar[l]_{{v}'}                                }
$$
where $\varrho=\lambda_T\comp\rho$ and
$\varrho'=\lambda_S\comp\rho'$. The result of applying $cyl$ to
the two middle columns in the above diagram is
$$
\xymatrix@M=4pt@H=4pt@C=14pt@R=16pt{
 \mrm{R}^2 X\ar[r]^{\mrm{R}f}                                             & \mrm{R}T                   \ar[r]^-{I_T}                     & cyl(w,g)         & \mrm{R}U   \ar[l]_-{I_U}                                           & \mrm{R}^2Z\ar[l]_{\mrm{R}v}                                                             \\
 \mrm{R}^4 X\ar[u]^{\mrm{R}\lambda^2_{\mrm{R}X}}\ar[r]\ar[d]_{\mrm{R}^2s} & \mrm{R}cyl(\mrm{R}w,Id) \ar[u]^{\mrm{R}\varrho} \ar[d]\ar[r] & M' \ar[u] \ar[d] & \mrm{R}^3U\ar[d]_{\mrm{R}^2q''}\ar[u]^{\mrm{R}\lambda^2_{U}}\ar[l] & \mrm{R}^4 Z\ar[l]_{\mrm{R}^3v}\ar[u]_{\mrm{R}\lambda^2_{\mrm{R}^2 Z}}\ar[d]^{\mrm{R}^2q'}  \\
 \mrm{R}^2\widetilde{X}  \ar[r]                                           & \mrm{R}cyl(u,\beta) \ar[r]                                   & M                & \mrm{R}^2W \ar[l]                                                  & \mrm{R}^2\widehat{Z}  \ar[l]_{\mrm{R}^2u'}                                          \\
 \mrm{R}^4 X\ar[d]_{\mrm{R}\lambda^2_{\mrm{R}X}}\ar[u]^{\mrm{R}^2p} \ar[r]& \mrm{R}cyl(\mrm{R}w',Id)\ar[u] \ar[d]_{\mrm{R}\varrho'}\ar[r]& M'' \ar[d] \ar[u]& \mrm{R}^3V\ar[u]^{\mrm{R}^2t''} \ar[d]_{\mrm{R}\lambda^2_{V}}\ar[l]& \mrm{R}^4 Z\ar[l]_{\mrm{R}^3v'}\ar[d]^{\mrm{R}\lambda^2_{\mrm{R} Z}}\ar[u]_{\mrm{R}^2t'} \\
 \mrm{R}^2 X       \ar[r]^{\mrm{R}f'}                                     & \mrm{R}S  \ar[r]^-{I_S}                                      & cyl(w',g')       & \mrm{R}V   \ar[l]_-{I_V}                                           & \mrm{R}^2 Z\ar[l]_{\mrm{R}v'}                                }
$$
Now apply $cyl$ to the three f{i}rst rows, getting
\begin{equation}\label{diagramaInfernalLema}
\xymatrix@M=4pt@H=4pt@C=15pt@R=16pt{
 \mrm{R}^3 X\ar[r]^{\mrm{R}^2f} \ar[d]_{}                                   & \mrm{R}^2T    \ar[d]         \ar[r]^-{\mrm{R}I_T}                &\mrm{R} cyl(w,g)\ar[d]    & \mrm{R}^2U   \ar[l]_-{\mrm{R}I_U}\ar[d]                              & \mrm{R}^3 Z\ar[l]_{\mrm{R}^2v}  \ar[d]^{}                                                  \\
 cyl(\mrm{R}\lambda^2_{\mrm{R}X},\mrm{R}^2s)\ar[r]                          & N' \ar[r]                                                        & N                       & N'' \ar[l]                                                            & cyl(\mrm{R}\lambda^2_{\mrm{R}^2 Z},\mrm{R}^2q')\ar[l]  \\
 \mrm{R}^3\widetilde{X}  \ar[r]\ar[u]                                       & \mrm{R}^2cyl(u,\beta) \ar[r]\ar[u]                               &\mrm{R} M  \ar[u]        & \mrm{R}^3W \ar[l] \ar[u]                                             & \mrm{R}^3\widehat{Z}  \ar[l]_{\mrm{R}^3u'} \ar[u]_{}                            \\
 \mrm{R}^5 X\ar[d]_{\mrm{R}^2\lambda^2_{\mrm{R}X}}\ar[u]^{\mrm{R}^3p} \ar[r]& \mrm{R}^2cyl(\mrm{R}w',Id)\ar[u] \ar[d]_{\mrm{R}^2\varrho'}\ar[r]&\mrm{R} M'' \ar[d] \ar[u]& \mrm{R}^4V\ar[u]^{\mrm{R}^3t''} \ar[d]_{\mrm{R}^2\lambda^2_{V}}\ar[l]& \mrm{R}^5 Z\ar[l]_{\mrm{R}^4v'}\ar[d]^{\mrm{R}^2\lambda^2_{\mrm{R} Z}}\ar[u]_{\mrm{R}^3t'} \\
 \mrm{R}^3 X       \ar[r]^{\mrm{R}^2f'}                                     & \mrm{R}^2S  \ar[r]^-{\mrm{R}I_S}                                 &\mrm{R} cyl(w',g')       & \mrm{R}^2V   \ar[l]_-{\mrm{R}I_V}                                    & \mrm{R}^3 Z\ar[l]_{\mrm{R}^2v'}                                }
\end{equation}
On the other hand, since $\mrm{R}^2\alpha \comp\mrm{R}^2s=
Id_{\mrm{R}^4X}$ and $\mrm{R}^2\beta' \comp\mrm{R}^2q'=
Id_{\mrm{R}^4Z}$, the equivalences
$\sigma=\mrm{R}\lambda^2_{\mrm{R}X}\comp \mrm{R}^2\alpha$ and
$\sigma'=\mrm{R}\lambda^2_{\mrm{R}Z}\comp \mrm{R}^2\beta'$ fit
into the commutative diagrams of $\mc{D}$
$$
\xymatrix@M=4pt@H=4pt@C=20pt@R=16pt{
 \mrm{R}^4 X \ar[r]^-{\mrm{R}\lambda^2_{\mrm{R}X}}\ar[d]_{\mrm{R}^2s} & \mrm{R}^2 Z \ar[d]_{Id} &  \mrm{R}^4 Z \ar[r]^-{\mrm{R}\lambda^2_{\mrm{R}Z}}\ar[d]_{\mrm{R}^2q'} & \mrm{R}^2 Z \ar[d]_{Id}\\
 \mrm{R}^2\widetilde{X}  \ar[r]^{\sigma}                              & \mrm{R}^2 Z             &  \mrm{R}^2\widehat{Z}  \ar[r]^{\sigma'}                                & \mrm{R}^2 Z \ .}
$$
It follows the existence of equivalences $\delta$ and $\delta'$
such that the diagrams
$$
\xymatrix@M=4pt@H=4pt@C=15pt@R=16pt{
                                & \mrm{R}^3 X                                                 &                                                         &                                & \mrm{R}^3 Z                                                 &                                                        \\
 \mrm{R}^3 X \ar[r] \ar[ru]^{Id}& cyl(\mrm{R}\lambda^2_{\mrm{R}X},\mrm{R}^2s) \ar[u]^{\delta} &  \mrm{R}^3 \widetilde{X} \ar[lu]_{\mrm{R}\sigma} \ar[l] & \mrm{R}^3 Z \ar[r] \ar[ru]^{Id}&cyl(\mrm{R}\lambda^2_{\mrm{R}^2 Z},\mrm{R}^2q') \ar[u]^{\delta'} &  \mrm{R}^3 \widehat{Z} \ar[lu]_{\mrm{R}\sigma'} \ar[l] \ .}
$$
are commutative. In addition, by def{i}nition we have that $\sigma
\comp\mrm{R}^2p=\mrm{R}\lambda^2_{\mrm{R}X}$ and
$\sigma'\comp\mrm{R}^2t'=\mrm{R}\lambda^2_{\mrm{R}Z}$.\\
If $A=cyl(\mrm{R}\lambda^2_{\mrm{R}X},\mrm{R}^2s)$ and
$A'=cyl(\mrm{R}\lambda^2_{\mrm{R}^2 Z},\mrm{R}^2q')$, attaching
these data to (\ref{diagramaInfernalLema}) and composing arrows we
obtain
$$
\xymatrix@M=4pt@H=4pt@C=17pt@R=15pt{
 \mrm{R}^3 X\ar[d]_{Id}           & \mrm{R}^3 X\ar[r]^{\mrm{R}^2f} \ar[d]\ar[l]_{Id}                                                    & \mrm{R}^2T    \ar[d]         \ar[r]^-{\mrm{R}I_T}                &\mrm{R} cyl(w,g)\ar[d]    & \mrm{R}^2U   \ar[l]_-{\mrm{R}I_U}\ar[d]                              & \mrm{R}^3 Z\ar[l]_{\mrm{R}^2v}  \ar[d]^{}   \ar[r]^{Id}                                                                           &\mrm{R}^3 Z \ar[d]^{Id} \\
 \mrm{R}^3 X                      & A\ar[l]_{\delta}  \ar[r]                                                                            & N' \ar[r]                                                        & N                       & N'' \ar[l]                                                            & A' \ar[l] \ar[r]^{\sigma'}                                                                                                      &\mrm{R}^3 Z\\
 \mrm{R}^3 X\ar[d]_{Id}\ar[u]^{Id}& \mrm{R}^5 X\ar[d]_{\mrm{R}^2\lambda^2_{\mrm{R}X}}\ar[u]\ar[r]\ar[l]_{\mrm{R}^2\lambda^2_{\mrm{R}X}} & \mrm{R}^2cyl(\mrm{R}w',Id)\ar[u] \ar[d]_{\mrm{R}^2\varrho'}\ar[r]&\mrm{R} M'' \ar[d] \ar[u]& \mrm{R}^4V\ar[u]^{\mrm{R}^3t''} \ar[d]_{\mrm{R}^2\lambda^2_{V}}\ar[l]& \mrm{R}^5 Z\ar[l]_{\mrm{R}^4v'}\ar[d]^{\mrm{R}^2\lambda^2_{\mrm{R} Z}}\ar[u]_{\mrm{R}^3t'}\ar[r]^{\mrm{R}^2\lambda^2_{\mrm{R} Z}}&\mrm{R}^3 Z \ar[u]_{Id}\ar[d]^{Id} \\
 \mrm{R}^3 X                      & \mrm{R}^3 X       \ar[r]^{\mrm{R}^2f'}    \ar[l]_{Id}                                               & \mrm{R}^2S  \ar[r]^-{\mrm{R}I_S}                                 &\mrm{R} cyl(w',g')       & \mrm{R}^2V   \ar[l]_-{\mrm{R}I_V}                                    & \mrm{R}^3 Z\ar[l]_{\mrm{R}^2v'}    \ar[r]^{Id}                                                                                   &\mrm{R}^3 Z \ .}
$$
F{i}nally, apply $cyl$ to the three lower rows to get
\begin{equation}\label{diagrAuxiliar}
\xymatrix@M=4pt@H=4pt@C=16pt@R=14pt{
 \mrm{R}^4 X\ar[d]_{Id}    & \mrm{R}^4 X\ar[r]^{\mrm{R}^3f} \ar[d]\ar[l]_{Id}   & \mrm{R}^3T \ar[d] \ar[r]^-{\mrm{R}^2I_T}  & \mrm{R}^2 cyl(w,g)\ar[d]    & \mrm{R}^3U   \ar[l]_-{\mrm{R}^2I_U}\ar[d]  & \mrm{R}^4 Z\ar[l]_{\mrm{R}^3v} \ar[d] \ar[r]^{Id} & \mrm{R}^4 Z \ar[d]^{Id} \\
 \mrm{R}^4 X \ar[d]        & \mrm{R}A\ar[l]_{\mrm{R}\delta}  \ar[r] \ar[d]      & \mrm{R}N' \ar[r] \ar[d]                   & \mrm{R}N  \ar[d]            & \mrm{R} N'' \ar[l] \ar[d]                  & \mrm{R}A' \ar[l] \ar[r]^{\mrm{R}\sigma'} \ar[d]   & \mrm{R}^4 Z \ar[d] \\
 cyl(\mrm{R}^3X)           & \widehat{B} \ar[r]\ar[l]                           & B' \ar[r]                               & B                           & B'' \ar[l]                                 & \widetilde{B} \ar[l]\ar[r]                        & cyl(\mrm{R}^3 Z) \\
 \mrm{R}^4 X  \ar[u]       & \mrm{R}^4 X \ar[u]\ar[r]^{\mrm{R}^3f'} \ar[l]_{Id} & \mrm{R}^3S \ar[r]^-{\mrm{R}^2 I_S} \ar[u] & \mrm{R}^2 cyl(w',g') \ar[u] & \mrm{R}^3V   \ar[l]_-{\mrm{R}^2I_V} \ar[u] & \mrm{R}^4 Z\ar[l]_{\mrm{R}^3v'}\ar[u]\ar[r]^{Id}  & \mrm{R}^4 Z  \ar[u] \ .}
\end{equation}
Let $\eta$ and $\eta'$ be the equivalences deduced from
\ref{cilaumentado}. Then the diagrams
$$
\xymatrix@M=4pt@H=4pt@C=15pt@R=16pt{
                                & \mrm{R}^4 X                   &                                    &                                & \mrm{R}^4 Z                    &                                  \\
 \mrm{R}^4 X \ar[r] \ar[ru]^{Id}& cyl(\mrm{R}^3X) \ar[u]^{\eta} &  \mrm{R}^4 {X} \ar[lu]_{Id} \ar[l] & \mrm{R}^4 Z \ar[r] \ar[ru]^{Id}& cyl(\mrm{R}^3Z) \ar[u]^{\eta'} &  \mrm{R}^4 {Z} \ar[lu]_{Id} \ar[l]}
$$
commute. If we adjoin $\eta$ and $\eta'$ to (\ref{diagrAuxiliar})
and compose arrows, the result is the hammock
$$
\xymatrix@M=4pt@H=4pt@C=16pt@R=14pt{
                      & \mrm{R}^4 X\ar[r]^{\mrm{R}^3f} \ar[d]_{m}\ar[ld]_{Id}    & \mrm{R}^3T \ar[d]_{l'} \ar[r]^-{\mrm{R}^2I_T}  & \mrm{R}^2 cyl(w,g)\ar[d]_{l}    & \mrm{R}^3U   \ar[l]_-{\mrm{R}^2I_U}\ar[d]_{l''}  & \mrm{R}^4 Z\ar[l]_{\mrm{R}^3v} \ar[d]_{k} \ar[rd]^{Id} & \\
 \mrm{R}^4X           & \widehat{B} \ar[r]\ar[l]_{\tau}                                 & B' \ar[r]                                      & B                               & B'' \ar[l]                                       & \widetilde{B} \ar[l]\ar[r]^{\tau'}              & \mrm{R}^4 Z \ , \\
                      & \mrm{R}^4 X \ar[u]^{m'}\ar[r]^{\mrm{R}^3f'} \ar[lu]^{Id} & \mrm{R}^3S \ar[r]^-{\mrm{R}^2 I_S} \ar[u]^{r'} & \mrm{R}^2 cyl(w',g') \ar[u]^{r} & \mrm{R}^3V   \ar[l]_-{\mrm{R}^2I_V} \ar[u]^{r''} & \mrm{R}^4 Z\ar[l]_{\mrm{R}^3v'}\ar[u]^{k'}\ar[ru]_{Id} & }
$$
that gives rise to
\begin{small}
$$
\xymatrix@M=4pt@H=4pt@C=16pt{
            & \mrm{R}^3X \ar[ld]_{Id}\ar[r]^{\mrm{R}^2\mu_X} \ar[d] & \mrm{R}^4 X\ar[r]^{\mrm{R}^3f} \ar[d]_{m}   & \mrm{R}^3T \ar[d]_{l'} \ar[r]^-{\mrm{R}^2I_T}  & \mrm{R}^2 cyl(w,g)\ar[d]_{l}    & \mrm{R}^3U   \ar[l]_-{\mrm{R}^2I_U}\ar[d]_{l''}  & \mrm{R}^4 Z\ar[l]_{\mrm{R}^3v} \ar[d]_{k} & \mrm{R}^3Z \ar[rd]^{Id}\ar[l]_{\mrm{R}^2\mu_Z} \ar[d] & \\
 \mrm{R}^3X & \widehat{B}\ar[l]\ar[r]^{Id}                          & \widehat{B} \ar[r]                          & B' \ar[r]                                      & B                               & B'' \ar[l]                                       & \widetilde{B} \ar[l]                      & \widetilde{B}\ar[r]\ar[l]_{Id}       & \mrm{R}^3 Z  \\
            & \mrm{R}^3X \ar[lu]^{Id}\ar[r]^{\mrm{R}^2\mu_X}\ar[u]  & \mrm{R}^4 X \ar[u]^{m'}\ar[r]^{\mrm{R}^3f'} & \mrm{R}^3S \ar[r]^-{\mrm{R}^2 I_S} \ar[u]^{r'} & \mrm{R}^2 cyl(w',g') \ar[u]^{r} & \mrm{R}^3V   \ar[l]_-{\mrm{R}^2I_V} \ar[u]^{r''} & \mrm{R}^4 Z\ar[l]_{\mrm{R}^3v'}\ar[u]^{k'}& \mrm{R}^3Z \ar[ru]_{Id}\ar[l]_{\mrm{R}^2\mu_Z}\ar[u] & }
$$
\end{small}
where the left triangle consists of
$$
\xymatrix@M=4pt@H=4pt@C=15pt@R=15pt{
                                                            &                        & \mrm{R}^3 X                                      &                           &                                                      \\
                                                            &                        & \mrm{R}^4 X \ar[u]^{\mrm{R}^2\lambda_{\mrm{R}X}} &                           &                                                    \\
 \mrm{R}^3 X \ar[r]^{\mrm{R}^2\mu_X} \ar@<1.4ex>[rruu]^{Id} & \mrm{R}^4 X \ar[r]^{m'}& \widehat{B} \ar[u]^{\tau}                        &  \mrm{R}^4 {X} \ar[l]_{m} & \mrm{R}^3 X  \ar@<-1.4ex>[lluu]_{Id}\ar[l]_{\mrm{R}^2\mu_X} \ .}
$$
whereas the right triangle is constructed analogously. Then by
\ref{LemaHamacaconRk} we have that $G\comp F\sim G'\comp F'$, that
f{i}nishes the proof.
\end{proof}

\begin{lema}\label{FrelRF}
The zig-zag $F$ given by
$$\xymatrix@M=4pt@H=4pt{ X &\mrm{R}X \ar[l]_{\lambda_X}\ar[r]^{f} & T &\mrm{R}Y \ar[r]^{\lambda_Y}\ar[l]_{w}& Y}$$
is related to the following zig-zag, that will be denote by
$\mrm{R}F$
$$\xymatrix@M=4pt@H=4pt@C=15pt{ X &\mrm{R}X \ar[l]_{\lambda_X}\ar[r]^{\mu_X} & \mrm{R}^2X \ar[r]^{\mrm{R}f}  & \mrm{R}T & \mrm{R}^2 Y \ar[l]_{\mrm{R}w} & \mrm{R}Y\ar[l]_{\mu_Y} \ar[r]^{\lambda_Y}& Y.}$$
\end{lema}

\begin{proof}
It is enough to consider the small hammock $\mc{H}^0$
$$\xymatrix@M=4pt@H=4pt@R=7pt{
              & \mrm{R}X \ar[ld]_{Id} \ar[r]^{{\mu_{X}}}\ar[dd]^{Id} &  \mrm{R}^2X  \ar[r]^-{\mrm{R}f} \ar[dd]_{\lambda_{\mrm{R}X}} & \mrm{R}T \ar[dd]_{\lambda_T} & \mrm{R}^2 Y \ar[l]_{\mrm{R}w} \ar[dd]_{\lambda_{\mrm{R}Y}} & \mrm{R}Y \ar[rd]^{Id}\ar[l]_{\mu_Y} \ar[dd]_{Id}& \\
 \mrm{R}X     &                                                      &                                                              &                              &                                                   &                                                 & \mrm{R}Y\\
              & \mrm{R}X \ar[lu]^{Id} \ar[r]^{Id}                    &  \mrm{R}X   \ar[r]^f                                         & T                            &  \mrm{R}Y   \ar[l]_{w}                            & \mrm{R}Y \ar[l]_{Id}\ar[ru]_{Id}                &  }$$
that is a particular case of a usual hammock. Then the statement
follows from \ref{LemaHamacaconRk}.
\end{proof}

\begin{lema}\label{FrelFs}
Given an element $F$ of $\mathrm{T}(X,Y)$ given by
$$\xymatrix@M=4pt@H=4pt{ X &\mrm{R}X \ar[l]_{\lambda_X}\ar[r]^{f} & T &\mrm{R}Y \ar[r]^{\lambda_Y}\ar[l]_{w}& Y}$$
and an equivalence $s:\mrm{R}X\rightarrow S$, then $F$ is related
to the following zig-zag, that will be denote by $F_S$
$$\xymatrix@M=4pt@H=4pt@C=15pt{ X &\mrm{R}X \ar[l]_{\lambda_X}\ar[r]^{\mu_X} & \mrm{R}^2X \ar[r]^{\mrm{R}s} & \mrm{R}S\ar[r]^-{I_S} & cyl(s,f) & \mrm{R}T \ar[l]_{I_T} & \mrm{R}^2 Y \ar[l]_{\mrm{R}w} & \mrm{R}Y\ar[l]_{\mu_Y} \ar[r]^{\lambda_Y}& Y\ .}$$
\end{lema}

\begin{proof}
By \ref{pdadesCylenD} we deduce the existence of
$H:cyl(\mrm{R}X)\rightarrow cyl(s,f)$ satisfying the following
property. If $I_{\mrm{R}X},J_{\mrm{R}X}$ are the canonical
inclusions of $\mrm{R}^2X$ into $cyl(\mrm{R}X)$ then $H\comp
I_{\mrm{R}X}=I_S\comp\mrm{R}s$ and
$H\comp J_{\mrm{R}X}=I_T\comp\mrm{R}f$.\\
Hence, we have the commutative diagram
$$\xymatrix@M=4pt@H=4pt@C=15pt@R=20pt{
  \mrm{R}^2X    \ar[rr]^{\mrm{R}f}\ar[d]^{J_{\mrm{R}X}}  &                     & \mrm{R}T \ar[d]_{I_T} &                        & \mrm{R}^2Y \ar[ll]_-{\mrm{R}w} \ar[d]_{Id} \\
  cyl(\mrm{R}X) \ar[rr]^{H}                              &                     & cyl(s,f)              &                        & \mrm{R}^2Y  \ar[ll]_{I_T\comp\mrm{R}w}      \\
  \mrm{R}^2X \ar[r]^{\mrm{R}s} \ar[u]_{I_{\mrm{R}X}}     &\mrm{R}S\ar[r]^{I_S} & cyl(s,f) \ar[u]^{Id}  &  \mrm{R}T \ar[l]_{I_T} & \mrm{R}^2Y\ar[u]^{Id} \ar[l]_{\mrm{R}w}   }$$
where all the vertical arrows are equivalences by the properties of $cyl$.\\
If $\rho:cyl(\mrm{R}X)\rightarrow \mrm{R}^2X$ is such that
$\rho\comp J_{\mrm{R}X}=\rho\comp J_{\mrm{R}X}=Id$, then the above
diagram can be completed to
$$\xymatrix@M=4pt@H=4pt@C=23pt@R=20pt{
            & \mrm{R} X \ar[ld]_{Id} \ar[r]^{\mu_{X}}\ar[d]^{J_{\mrm{R}X\comp\mu_{X}}} &  \mrm{R}^2X    \ar[r]^{\mrm{R}f}\ar[d]^{J_{\mrm{R}X}}        & \mrm{R}T \ar[d]_{I_T} & \mrm{R}^2Y \ar[l]_-{\mrm{R}w} \ar[d]_{Id}       & \mrm{R}Y \ar[l]_{\mu_Y} \ar[d]_{\mu_Y} \ar[rd]^{Id}  & \\
 \mrm{R}X   & cyl(\mrm{R}X) \ar[l]_{\lambda_{\mrm{R}X}\comp\rho}\ar[r]^{Id}                     &  cyl(\mrm{R}X) \ar[r]^{H}                                    & cyl(s,f)              & \mrm{R}^2Y  \ar[l]_{I_T\comp\mrm{R}w}           & \mrm{R}^2Y     \ar[l]_{Id}   \ar[r]^{\lambda_{\mrm{R}Y}}    & \mrm{R}Y ,   \\
            & \mrm{R} X\ar[lu]^{Id} \ar[r]^{\mu_{X}}\ar[u]_{I_{\mrm{R}X\comp\mu_{X}}}  &  \mrm{R}^2X \ar[r]^{I_S\comp\mrm{R}s} \ar[u]_{I_{\mrm{R}X}}  & cyl(s,f) \ar[u]^{Id}  & \mrm{R}^2Y\ar[u]^{Id} \ar[l]_{I_T\comp\mrm{R}w} & \mrm{R}Y  \ar[ru]_{Id}\ar[u]^{\mu_{Y}} \ar[l]_{\mu_Y}&  }$$
that implies that $F_S$ is related to the zig-zag $\mrm{R}F$ given
in \ref{FrelRF}. Then the result follows from the transitivity of
$\sim$.
\end{proof}

The next lemma can be proved in a similar way.

\begin{lema}\label{FrelFu}
Consider the element $F$ of $\mathrm{T}(X,Y)$
$$\xymatrix@M=4pt@H=4pt{ X &\mrm{R}X \ar[l]_{\lambda_X}\ar[r]^{f} & T &\mrm{R}Y \ar[r]^{\lambda_Y}\ar[l]_{w}& Y}$$
and an equivalence $u:\mrm{R}Y\rightarrow U$. Then $F$ is related
to the zig-zag $F^u$ given by
$$\xymatrix@M=4pt@H=4pt@C=15pt{ X &\mrm{R}X \ar[l]_{\lambda_X}\ar[r]^{\mu_X} & \mrm{R}^2X \ar[r]^{\mrm{R}f} & \mrm{R}T\ar[r]^-{I_T} & cyl(w,u) & \mrm{R}U \ar[l]_{I_U} & \mrm{R}^2 Y \ar[l]_{\mrm{R}u} & \mrm{R}Y\ar[l]_{\mu_Y} \ar[r]^{\lambda_Y}& Y\ .}$$
\end{lema}

\begin{lema}\label{EliminarFlechasComp}

The composition of morphisms in $Ho\mc{D}$ allows to ``delete
pairs of inverse arrows''. In other words, assume given composable
zig-zags $F$ and $G$ described respectively as
$$\xymatrix@M=4pt@H=4pt@R=8pt{
 X &\mrm{R}X \ar[l]_{\lambda_X}\ar[r]^{f} & U & \mrm{R}S\ar[l]_u & \mrm{R}Y \ar[r]^{\lambda_Y}\ar[l]_{s}& Y \\
 Y &\mrm{R}Y \ar[l]_{\lambda_Y} \ar[r]^s  &\mrm{R} S  \ar[r]^{g} & V &\mrm{R}Z \ar[r]^{\lambda_Z}\ar[l]_{v}& Z\ .}
$$
Then the zig-zag $G\comp F$ is equivalent to the composition of
the zig-zags $F'$ and $G'$ given by
$$\xymatrix@M=4pt@H=4pt@R=8pt{
 X &\mrm{R}X \ar[l]_{\lambda_X} \ar[r]^{f} & U & \mrm{R}S \ar[r]^{\lambda_S}\ar[l]_{u}& S\ar@{}[r]|{;} & S &\mrm{R}S \ar[l]_{\lambda_S} \ar[r]^{g} & V &\mrm{R}Z \ar[r]^{\lambda_Z}\ar[l]_{v}& Z\ .}
$$
\end{lema}

\begin{proof}
Since the following diagram commutes
$$\xymatrix@M=4pt@H=4pt@R=16pt@C=16pt{
 U \ar[d]_{Id} &\mrm{R}Y \ar[l]_{s\comp u} \ar[r]^{g\comp s} \ar[d]_s & V \ar[d]_{Id} \\
 U             &\mrm{R}S \ar[l]_{u} \ar[r]^g                          & V }
$$
we deduce from the functoriality of $cyl$ a morphism
$\alpha:cyl(s\comp u,g\comp s)\rightarrow cyl(u,g)$ such that the
diagram
$$\xymatrix@M=4pt@H=4pt@R=16pt@C=16pt{
 \mrm{R}U \ar[d]_{Id}\ar[r]^-{I_U} &cyl(s\comp u,g\comp s) \ar[d]_{\alpha} & \mrm{R}V \ar[l]_-{I_V}\ar[d]_{Id} \\
 \mrm{R}U   \ar[r]^-{J_U}          &cyl(u,g)                               & \mrm{R}V \ar[l]_-{J_V}}
$$
commutes. In addition, by \ref{exactitudcyl}, $\alpha\in\mrm{E}$.
Then the above diagram can be completed to the following small
hammock
$$\xymatrix@M=4pt@H=4pt@R=10pt@C=15pt{
          & \mrm{R}X \ar[ld]_{Id}\ar[dd]_{Id} \ar[r]^{\mu_X} & \mrm{R}^2X \ar[dd]_{Id} \ar[r]^{\mrm{R}f} & \mrm{R}U \ar[dd]_{Id}\ar[r]^-{I_U} & cyl(s\comp u,g\comp s) \ar[dd]_{\alpha} & \mrm{R}V \ar[l]_-{I_V}\ar[dd]_{Id} & \mrm{R}^2Z \ar[l]_{\mrm{R}v}\ar[dd]_{Id} & \mrm{R}Z \ar[dd]_{Id}\ar[l]_{\mu_X}\ar[dd]_{Id} \ar[rd]^{Id}& \\
 \mrm{R}X &                                                 &                                          &                                  &                                        &                                  &                             &                                                &\mrm{R}Z \\
          & \mrm{R}X \ar[lu]^{Id} \ar[r]^{\mu_X}            & \mrm{R}^2X \ar[r]^{\mrm{R}f}             & \mrm{R}U   \ar[r]^-{J_U}          & cyl(u,g)                               & \mrm{R}V \ar[l]_-{J_V}            & \mrm{R}^2Z \ar[l]_{\mrm{R}v}& \mrm{R}Z \ar[l]_{\mu_X}    \ar[ru]_{Id}        & }
$$
that relates $G\comp F$ to $G'\comp F'$.
\end{proof}

\begin{lema}\label{EliminarFlechasCompTotal} The composition of
the following zig-zags $F$ and $G$
$$\xymatrix@M=4pt@H=4pt@R=8pt{
 X &\mrm{R}X \ar[l]_{\lambda_X}\ar[r]^{f} & U                             & \mrm{R}Y \ar[r]^{\lambda_Y}\ar[l]_{u}& Y \\
 Y &\mrm{R}Y \ar[l]_{\lambda_Y} \ar[r]^u  & U \ar[r]^{g} & V              &\mrm{R}Z \ar[r]^{\lambda_Z}\ar[l]_{v}& Z}
$$
is equivalent to
$$\xymatrix@M=4pt@H=4pt@R=8pt{
 X &\mrm{R}X \ar[l]_{\lambda_X}\ar[r]^{f} & U \ar[r]^{g} & V   &\mrm{R}Z \ar[r]^{\lambda_Z}\ar[l]_{v}& Z\ .}
$$
\end{lema}

\begin{proof} By def{i}nition $G\comp F$ is
$$\xymatrix@M=4pt@H=4pt@R=10pt@C=15pt{
  X       & \mrm{R}X \ar[l]_{\lambda_X }\ar[r]^{\mu_X} & \mrm{R}^2X \ar[r]^{\mrm{R}f} & \mrm{R}U \ar[r]^-{I_U} & cyl(u,g\comp u)  & \mrm{R}V \ar[l]_-{I_V} & \mrm{R}^2Z \ar[l]_{\mrm{R}v}& \mrm{R}Z \ar[l]_{\mu_X} \ar[r]^{\lambda_X} & X \ .}
$$
Since $g\comp u=Id_V\comp (g\comp u)$, from \ref{cilaumentado} we
deduce the existence of $\beta:cyl(u,g\comp u)\rightarrow
\mrm{R}V$ such that $\beta\comp I_U=\mrm{R}g$ and $\beta\comp
I_V=Id_{\mrm{R}V}$. Then, $\beta\in\mrm{E}$ and we have the
hammock
$$\xymatrix@M=4pt@H=4pt@R=10pt@C=15pt{
          & \mrm{R}X \ar[ld]_{Id}\ar[dd]_{Id} \ar[r]^{\mu_X} & \mrm{R}^2X \ar[dd]_{Id} \ar[r]^{\mrm{R}f} & \mrm{R}U \ar[dd]_{Id}\ar[r]^-{I_U} & cyl(u,g\comp u) \ar[dd]_{\beta} & \mrm{R}V \ar[l]_-{I_V}\ar[dd]_{Id} & \mrm{R}^2Z \ar[l]_{\mrm{R}v}\ar[dd]_{Id} & \mrm{R}Z \ar[dd]_{Id}\ar[l]_{\mu_X}\ar[dd]_{Id} \ar[rd]^{Id}& \\
 \mrm{R}X &                                                  &                                           &                                    &                          &                                    &                                          &                                                             &\mrm{R}Z \ .\\
          & \mrm{R}X \ar[lu]^{Id} \ar[r]^{\mu_X}             & \mrm{R}^2X \ar[r]^{\mrm{R}f}              & \mrm{R}U   \ar[r]^-{\mrm{R}g}      & \mrm{R}V                 & \mrm{R}V \ar[l]_-{Id}              & \mrm{R}^2Z \ar[l]_{\mrm{R}v}             & \mrm{R}Z \ar[l]_{\mu_X}    \ar[ru]_{Id}                     & }
$$
To f{i}nish the proof, it suf{f}{i}ces to take into account
\ref{FrelRF}.
\end{proof}

\begin{proof}[\textbf{Proof of \ref{ProHoDcat}}]\mbox{}\\
To see that $Ho\mc{D}$ is a category, it remains to check that the
composition of morphisms is associative and that is has a unit.\\[0.2cm]
\textbf{Unit:} given $X$ in $\mc{D}$, the unit for the composition
in $\mathrm{Hom}_{Ho\mc{D}}(X,X)$ is $1_X=\gamma(Id_X)$, that is
the morphism represented by
$$\xymatrix@M=4pt@H=4pt{X &\mrm{R}X \ar[l]_{\lambda_X}\ar[r]^{Id} &\mrm{R}X &\mrm{R}X  \ar[r]^{\lambda_X}\ar[l]_{Id} & X\ .}$$
Indeed, if $\widehat{f}\in \mathrm{Hom}_{Ho\mc{D}}(X,Y)$ is the
class of $\xymatrix@M=4pt@H=4pt{X &\mrm{R}X
\ar[l]_{\lambda_X}\ar[r]^{f} & T &\mrm{R}Y
\ar[r]^{\lambda_Y}\ar[l]_{w}& Y}$
then by def{i}nition $\widehat{f}\comp 1_X$ is represented by
$$\xymatrix@M=4pt@H=4pt@C=15pt{
 X & \mrm{R}X \ar[l]_{\lambda_X} \ar[r]^{\mu_X}  & \mrm{R}^2 X \ar[r]^-{I_{\mrm{R}X}} & cyl(Id_{\mrm{R}X},f) & \mrm{R}T \ar[l]_-{I_{T}} &  \mrm{R}^2 Y \ar[l]_{\mrm{R}w} & \mrm{R}Y \ar[l]_{\mu_Y} \ar[r]^{\lambda_Y} & Y   .}$$
From the commutative diagram
$$\xymatrix@M=4pt@H=4pt{
 \mrm{R}X\ar[r]^{Id} \ar[d]_f & \mrm{R}X \ar[d]_{f} \\
 T \ar[r]^{Id}& T }$$
we deduce by \ref{cilaumentado} the existence of an equivalence
$\alpha:cyl(\mrm{R}X,f)\rightarrow \mrm{R}T$ such that
$\alpha\comp I_{\mrm{R}X}=\mrm{R}f$ and $\alpha\comp I_T=Id_T$, so
we have the hammock
$$\xymatrix@M=4pt@H=4pt{
              & \mrm{R}^2X \ar[ld]_{Id} \ar[r]^{\mrm{R}{\mu_{X}}}\ar[d]^{\mrm{R}{\mu_{X}}}  &  \mrm{R}^3X  \ar[r]^-{ \mrm{R}I_{\mrm{R}X}} & \mrm{R}cyl(Id,f) \ar[d]_{\mrm{R}\alpha} & & \mrm{R}^2Y \ar[ll]_-{\mrm{R}I_T\comp\mrm{R}^2w\comp\mrm{R}\mu_Y} \ar[d]_{\mrm{R}\mu_Y} \ar[rd]^{Id} & \\
 \mrm{R}^2X   & \mrm{R}^3X \ar[l]_{\mrm{R}\lambda_{\mrm{R}X}}\ar[rr]^{\mrm{R}^2{f}}         &                                             & \mrm{R}^2T                              & & \mrm{R}^3Y     \ar[ll]_{\mrm{R}^2w}   \ar[r]^{\mrm{R}\lambda_{\mrm{R}X}}                                             & \mrm{R}^2Y .   \\
              & \mrm{R}^2X\ar[lu]^{Id} \ar[rr]^{\mrm{R}{f}}\ar[u]_{\mu_{\mrm{R}X}}          &                                             & \mrm{R}T  \ar[u]^{\mu_T}                & & \mrm{R}^2Y  \ar[ru]_{Id}\ar[u]^{\mu_{\mrm{R}Y}} \ar[ll]_{\mrm{R}w} &  }$$
Therefore, $\widehat{f}\comp 1_X=\widehat{f}$, and it can be
proved analogously that $1_Y\comp\widehat{f}=\widehat{f}$.\\[0.2cm]
\textbf{Associativity:} Let $\widehat{f}$, $\widehat{g}$ and
$\widehat{h}$ three composable morphisms in $Ho\mc{D}$,
represented respectively by the following zig-zags $F$, $G$ and
$H$
$$\xymatrix@M=4pt@H=4pt@R=8pt{
 X &\mrm{R}X \ar[l]_{\lambda_X}\ar[r]^{f} & U &\mrm{R}Y \ar[r]^{\lambda_Y}\ar[l]_{u}& Y\\
 Y &\mrm{R}Y \ar[l]_{\lambda_Y}\ar[r]^{g} & V &\mrm{R}Z \ar[r]^{\lambda_Z}\ar[l]_{v}& Z\\
 Z &\mrm{R}Z \ar[l]_{\lambda_Z}\ar[r]^{h} & W &\mrm{R}S \ar[r]^{\lambda_S}\ar[l]_{w}& S}
$$
Let us see that $(\widehat{h}\comp \widehat{f})\comp \widehat{g}=
\widehat{h}\comp (\widehat{g}\comp \widehat{f})$. Consider the
diagram
$$\xymatrix@M=4pt@H=4pt@R=15pt{
 \mrm{R}^2U \ar[rd]_{\mrm{R}I_U} & \mrm{R}^3Y \ar[l]_{\mrm{R}^2u}\ar[r]^{\mrm{R}^2g} & \mrm{R}^2V\ar[ld]^{\mrm{R}I_V} \ar[rd]_{\mrm{R}J_V} & \mrm{R}^3Z \ar[l]_{\mrm{R}^2v}\ar[r]^{\mrm{R}^2h} & \mrm{R}^2W\ar[ld]^{\mrm{R}I_W} \\
                                 & \mrm{R}cyl(u,g)\ar[rd]_{I}                        &                                                     & \mrm{R} cyl(v,h)\ar[ld]^{J}                      &            \\
                                 &                                                   & cyl(I_V,J_V)                                            &                                        &}
$$
where $I_V$, $I_W$ and $J$ are equivalences.\\
Let us check that $(\widehat{h}\comp \widehat{g})\comp
\widehat{f}$ and $\widehat{h}\comp (\widehat{g}\comp \widehat{f})$
coincides with the morphism represented by the zig-zag $L$ given
by
$$\xymatrix@M=4pt@H=4pt@R=15pt@C=18pt{
 X &\mrm{R}X \ar[l]_{\lambda_X}\ar[r]^{\hat{\mu}_{X}} &\mrm{R}^3X \ar[r]^{\mrm{R}^2f} & \mrm{R}^2 U \ar[r]^-{I\comp \mrm{R} I_U}& cyl(I_V,J_V) &  \mrm{R}^2W \ar[l]_-{J\comp\mrm{R}I_W} & \mrm{R}^3S \ar[l]_{\mrm{R}^2w }& \mrm{R}S \ar[l]_{\hat{\mu}_S} \ar[r]^{\lambda_S}& S ,}
$$
where $\hat{\mu}:\mrm{R}\rightarrow\mrm{R}^3$ is
$\hat{\mu}_M=\mrm{R}\mu_M\comp\mu_M$, for every object $M$ in $\mc{D}$.\\
Having into account the equivalences $u:\mrm{R}Y\rightarrow U$ and
$v:\mrm{R}Z\rightarrow V$, by \ref{FrelFs} it follows that $G\sim
G_u$ and $H\sim H_v$, where $G_u$ and $H_v$ are respectively
$$\xymatrix@M=4pt@H=4pt@R=14pt@C=16pt{
 Y &\mrm{R}Y \ar[l]_{\lambda_Y}\ar[r]^{\mu_{Y}} &\mrm{R}^2Y \ar[r]^{\mrm{R}u} & \mrm{R} U \ar[r]^-{I_U} & cyl(u,g) &  \mrm{R}V \ar[l]_-{I_V} & \mrm{R}^2Z \ar[l]_{\mrm{R}v }& \mrm{R}Z \ar[l]_{\mu_Z} \ar[r]^{\lambda_Z}& Z \\
 Z &\mrm{R}Z \ar[l]_{\lambda_Z}\ar[r]^{\mu_{Z}} &\mrm{R}^2Z \ar[r]^{\mrm{R}v} & \mrm{R} V \ar[r]^-{J_V} & cyl(v,h) &  \mrm{R}W \ar[l]_-{I_W} & \mrm{R}^2S \ar[l]_{\mrm{R}w }& \mrm{R}S \ar[l]_{\mu_S} \ar[r]^{\lambda_S}& S \ .}
$$
Therefore $\widehat{h}\comp \widehat{g}$ is the class of $H_v\comp
G_u$, that by \ref{EliminarFlechasComp} coincides with the
composition of zig-zags
$$\xymatrix@M=4pt@H=4pt@R=11pt@C=16pt{
 Y & \mrm{R}Y \ar[l]_{\lambda_Y}\ar[r]^{\mu_{Y}} & \mrm{R}^2Y \ar[r]^{\mrm{R}u} & \mrm{R} U \ar[r]^-{I_U} & cyl(u,g) &  \mrm{R}V \ar[l]_-{I_V} \ar[r]^{\lambda_{V}} & V \\
 V & \mrm{R}V \ar[l]_{\lambda_V}\ar[r]^-{J_V} & cyl(v,h) &  \mrm{R}W \ar[l]_-{I_W} & \mrm{R}^2S \ar[l]_{\mrm{R}w }& \mrm{R}S \ar[l]_{\mu_S} \ar[r]^{\lambda_S}& S \ ,}
$$
and this is by def{i}nition the zig-zag $(H\comp G)'$
$$\xymatrix@M=4pt@H=4pt@R=14pt@C=16pt{
 Y & \mrm{R}Y \ar[l]_{\lambda_Y}\ar[r]^{\hat{\mu}_{Y}} & \mrm{R}^3Y \ar[r]^{\mrm{R}^2u} & \mrm{R}^2U \ar[r]^-{I\comp\mrm{R}I_U} & cyl(I_V,J_V) & \mrm{R}^2S \ar[l]_-{J\comp\mrm{R}t'} &\mrm{R}S \ar[l]_{\mu_W}\ar[r]^{\lambda_{S}} & S }
$$
where $t'=I_W\comp\mrm{R}w\comp\mu_S:\mrm{R}S\rightarrow cyl(v,h)$.\\
In addition, $F$ is related to the zig-zag $\mrm{R}^2F$ consisting
of
$$\xymatrix@M=4pt@H=4pt{
 X & \mrm{R}X\ar[l]_{\lambda_X}\ar[r]^{\hat{\mu}_{X}} & \mrm{R}^3X \ar[r]^{\mrm{R}^2f} &\mrm{R}^2U & \mrm{R}^3Y\ar[l]_-{\mrm{R}^2u} & \mrm{R}Y \ar[l]_{\hat{\mu}_Y}\ar[r]^{\lambda_{Y}} & Y }
$$
and by \ref{EliminarFlechasCompTotal} it follows that $
(H\comp G)' \comp\mrm{R}^2F \sim L$.\\[0.1cm]
\indent On the other hand, $G\comp F\sim G_u\comp \mrm{R}F$, that
becomes after deleting arrows in
$$\xymatrix@M=4pt@H=4pt@R=14pt@C=16pt{
 X &\mrm{R}X \ar[l]_{\lambda_X}\ar[r]^{\mu_{X}} &\mrm{R}^2X \ar[r]^{\mrm{R}f} & \mrm{R} U \ar[r]^-{I_U} & cyl(u,g) &  \mrm{R}V \ar[l]_-{I_V} & \mrm{R}^2Z \ar[l]_{\mrm{R}v }& \mrm{R}Z \ar[l]_{\mu_Z} \ar[r]^{\lambda_Z}& Z }
$$
Hence the result of composing with $H_v$ is related to the
composition of
$$\xymatrix@M=4pt@H=4pt@R=11pt@C=16pt{
 X &\mrm{R}X \ar[l]_{\lambda_X}\ar[r]^{\mu_{X}} &\mrm{R}^2X \ar[r]^{\mrm{R}f} & \mrm{R} U \ar[r]^-{I_U} & cyl(u,g) &  \mrm{R}V \ar[l]_-{I_V} \ar[r]^{\lambda_V} & V \\
 V & \mrm{R} V \ar[l]_{\lambda_V} \ar[r]^-{J_V} & cyl(v,h) &  \mrm{R}W \ar[l]_-{I_W} & \mrm{R}^2S \ar[l]_{\mrm{R}w }& \mrm{R}S \ar[l]_{\mu_S} \ar[r]^{\lambda_S}& S \ ,}
$$
that is $L$ by def{i}nition.\\[0.2cm]
\textbf{Functoriality of $\gamma:\mc{D}\rightarrow Ho\mc{D}$:}
we have by def{i}nition that the image under $\gamma$ of an
identity morphism in $\mc{D}$ is an identity morphism in
$Ho\mc{D}$. Given composable morphisms $f:X\rightarrow Y$ and
$g:Y\rightarrow Z$ in
$\mc{D}$, let us check that $\gamma(g\comp f)=\gamma(g)\comp\gamma(f)$.\\
By def{i}nition, $\gamma(g\comp f)$ is represented by
$$\xymatrix@M=4pt@H=4pt@R=11pt@C=16pt{
 X &\mrm{R}X \ar[l]_{\lambda_X}\ar[r]^{\mrm{R}f} & \mrm{R}Y \ar[r]^{\mrm{R}g}  & \mrm{R} Z  &  \mrm{R}Z \ar[l]_-{Id} \ar[r]^{\lambda_Z} & Z }
$$
In the same way, the composition of the zig-zags\\
$$\xymatrix@M=4pt@H=4pt@R=11pt@C=15pt{
 X &\mrm{R}X \ar[l]_{\lambda_X}\ar[r]^{\mrm{R}f} & \mrm{R}Y   &  \mrm{R}Y \ar[l]_-{Id} \ar[r]^{\lambda_Y} & Y \ar@{}[r]|{;} &  Y & \mrm{R}Y \ar[l]_{\lambda_Y}\ar[r]^{\mrm{R}g} & \mrm{R}Z  &  \mrm{R}Z \ar[l]_-{Id} \ar[r]^{\lambda_Z} & Z}$$
is $\xymatrix@M=4pt@H=4pt@R=16pt@C=15pt{
 X & \mrm{R}X \ar[l]_{\lambda_X}  \ar[r]^{\mu_X} &\mrm{R}^2X \ar[r]^{\mrm{R}^2 f} & \mrm{R} ^2Y \ar[r]^-{I_{\mrm{R}Y}}  & cyl(Id_{\mrm{R}Y},\mrm{R}g) & \mrm{R}^2Z \ar[l]_-{I_{\mrm{R}Z}}  & \mrm{R}Z  \ar[l]_{\mu_Z} \ar[r]^{\lambda_Z} &
 Z}$ by def{i}nition.
It follows from \ref{cilaumentado} the existence of
$\rho:cyl(Id,\mrm{R}g)\rightarrow \mrm{R}^2Z$ such that $\rho\comp
I_{\mrm{R}Y }=\mrm{R}^2g$ and $\rho\comp I_{\mrm{R}Z
}=Id_{\mrm{R}^2Z}$. Hence, $\rho$ is an equivalence and we have
the following small hammock
$$\xymatrix@M=4pt@H=4pt@R=10pt@C=15pt{
          & \mrm{R}X \ar[ld]_{Id}\ar[dd]_{Id} \ar[r]^{\mu_X} & \mrm{R}^2X \ar[dd]_{Id} \ar[r]^{\mrm{R}^2f} & \mrm{R}^2Y \ar[dd]_{Id}\ar[r]^-{I_{\mrm{R}Y}} & cyl(Id,\mrm{R}g) \ar[dd]_{\rho} & \mrm{R}^2Z \ar[dd]_{Id} \ar[l]_-{I_{\mrm{R}Z}} & \mrm{R}Z \ar[l]_{\mu_Z}\ar[dd]_{Id} \ar[rd]^{Id}&  \\
 \mrm{R}X &                                                  &                                             &                                               &                                 &                                                &                                                 & \mrm{R}Z \\
          & \mrm{R}X \ar[lu]^{Id} \ar[r]^{\mu_X}             & \mrm{R}^2X \ar[r]^{\mrm{R}^2f}              & \mrm{R}^2Y    \ar[r]^{\mrm{R}^2g}             & \mrm{R}^2Z                      & \mrm{R}^2Z \ar[l]_{Id}                         & \mrm{R}Z \ar[l]_{\mu_Z}\ar[ru]_{Id}             & }
$$
where the lower zig-zag represents the morphism $\gamma(g\comp
f)$ by \ref{FrelRF}, and then the required equality holds.\\[0.2cm]
$\textbf{Universal property}:$ First, $\gamma:\mc{D}\rightarrow
Ho\mc{D}$ is such that $\gamma(\mrm{E})\subseteq \{$isomorphisms
of $Ho\mc{D}\}$. Indeed, if $w:X\rightarrow Y$ is an equivalence,
then $\mrm{R}w$ is so, and $\gamma(w)^{-1}$ is the morphism in
$Ho\mc{D}$ given by
$$\xymatrix@M=4pt@H=4pt@R=11pt@C=16pt{
 Y & \mrm{R}Y \ar[l]_{\lambda_Y}\ar[r]^{Id} & \mrm{R}Y  &  \mrm{R}X\ar[l]_-{\mrm{R}w} \ar[r]^{\lambda_X} & X }
$$
that clearly is the inverse of $\gamma(w)$ (by \ref{EliminarFlechasCompTotal}).\\
\indent It remains to see that the pair $(Ho\mc{D},\gamma)$
satisf{i}es the universal property of the localized category $\mc{D}[\mrm{E}^{-1}]$.\\
Let $\mathcal{F}:\mc{D}\rightarrow \cont$ be a functor such that
$\mathcal{F}$ maps equivalences into isomorphisms.\\
We must prove that there exists an unique functor
$\mc{G}:Ho\mc{D}\rightarrow
\cont$ such that $\mc{G}\comp\gamma=\mc{F}$.\\
\textit{Existence:} def{i}ne $\mathcal{G}$ as
$\mathcal{G}X=\mathcal{F}X$ if $X$ is any object of $\mc{D}$. The
image under $\mathcal{G}$ of a morphism $\hat{f}$ in $Ho\mc{D}$
given by $\xymatrix@M=4pt@H=4pt@C=15pt{ X & \mrm{R}X
\ar[l]_{\lambda_X}\ar[r]^{f} & T  &
\mrm{R}Y\ar[l]_-{w}\ar[r]^{\lambda_Y} & Y}$ is
$$\mathcal{G}(\hat{f})=(\mc{F}\lambda_Y)\comp(\mc{F}w)^{-1}\comp(\mc{F}f)\comp(\mc{F}\lambda_X)^{-1} \ .$$
Note that $\mc{G}(\hat{f})$ does not depend on the zig-zag chosen.
Indeed, if two zig-zags $L$ and $L'$ represent $\hat{f}$, both
will be related by a hammock as in (\ref{diamante}). This hammock
becomes, after applying $\mc{F}$, a commutative diagram in
$\cont$. As the equivalences are now
isomorphisms, it follows that $\mc{F}(L)=\mc{F}(L')$.\\
\indent In addition, the equality $\mc{G}\comp\gamma=\mc{F}$
holds. To see this, let $f:X\rightarrow Y$ be a morphism in
$\mc{D}$. By def{i}nition $\mc{G}(\gamma f)$ is
$$\mc{F}\lambda_Y\comp(\mc{F}Id_{\mrm{R}Y})^{-1}\comp\mc{F}(\mrm{R}f)\comp(\mc{F}\lambda_X)^{-1}=
\mc{F}\lambda_Y\comp\mc{F}(\mrm{R}f)\comp(\mc{F}\lambda_X)^{-1}\
,$$ that agrees with $F(f)$ since $\lambda_Y\comp \mrm{R}f =
f\comp\lambda_X$ in $\mc{D}$. \\
\indent Next we check that $\mc{G}$ is functorial. It is clear
that $\mc{G}$ maps identities into identities, since $\mc{G}(Id_X)=\mc{G}(\gamma(Id_X))=Id_{F(X)}$.\\
On the other hand, let $\hat{f}$ and $\hat{g}$ be composable
morphisms in $Ho\mc{D}$ represented by the respective zig-zags
$$\xymatrix@M=4pt@H=4pt@R=10pt@C=15pt{
 X & \mrm{R}X \ar[l]_{\lambda_X}\ar[r]^{f} & T  & \mrm{R}Y\ar[l]_-{w}\ar[r]^{\lambda_Y} & Y\ar@{}[r]|{;} & Y & \mrm{R}Y \ar[l]_{\lambda_Y}\ar[r]^{g} & S  & \mrm{R}Z\ar[l]_-{v}\ar[r]^{\lambda_Z} & Z.}
$$
We must see that
$\mc{G}(\hat{g})\comp\mc{G}\hat{f})=\mc{G}(\hat{g}\comp\hat{f})$,
that is, the following composition of morphisms must coincide in
$\cont$
$$\xymatrix@M=4pt@H=4pt@R=10pt@C=16pt{
 \mc{F}(\mrm{R}X)  \ar[r]^{\mc{F}f} & \mc{F} T \ar[r]^-{(\mc{F}w)^{-1}} & \mc{F}\mrm{R}Y\ar[r]^{ \mc{F}g} &  \mc{F}S \ar[r]^-{(\mc{F}v)^{-1}} & \mc{F}(\mrm{R}Z)\\
 \mc{F}(\mrm{R} X) \ar[r]^{\mc{F}\mu_X} &\mc{F}(\mrm{R}^2X) \ar[r]^{\mc{F}(\mrm{R}f)} & \mc{F} (\mrm{R}T) \ar[r]^-{\mc{F}I_T} & \mc{F}cyl(w,g)\ar[r]^-{(\mc{F}I_S)^{-1}}  & \mc{F}\mrm{R}S\ar[r]^{ (\mc{F}\mrm{R}v)^{-1}} &  \mc{F}(\mrm{R}^2Z) \ar[r]^-{(\mc{F}\mu_Z)^{-1}} & \mc{F}(\mrm{R}Z) .}
$$
where we have already deleted in
$\mc{G}(\hat{g})\comp\mc{G}(\hat{f})$ and
$\mc{G}(\hat{g}\comp\hat{f})$ the isomorphisms
$(\mc{F}\lambda_X)^{-1}$ and $\mc{F}\lambda_Z$.
Then, it suf{f}{i}ces to prove the commutativity in $\cont$ of the
following diagrams (I), (II) and (III)
$$ \xymatrix@M=4pt@H=4pt@R=24pt@C=24pt{
 \mc{F}(\mrm{R}X)  \ar[r]^{\mc{F}f} \ar@{}[rd]|{(I)}\ar[d]_{\mc{F}\mu_X} & \mc{F} T \ar[r]^-{(\mc{F}w)^{-1}}                             & \mc{F}\mrm{R}Y\ar[r]^{ \mc{F}g}\ar@{}[d]|{(II)} & \mc{F}S \ar[r]^-{(\mc{F}v)^{-1}}  \ar@{}[rd]|{(III)}                    & \mc{F}(\mrm{R}Z) \ar[d]^{\mc{F}\mu_Z} \\
 \mc{F}(\mrm{R}^2X) \ar[r]^{\mc{F}(\mrm{R})f}                            & \mc{F}(\mrm{R}T) \ar[u]_{\mc{F}\lambda_T} \ar[r]^-{\mc{F}I_T} & \mc{F}cyl(w,g)\ar[r]^-{(\mc{F}I_S)^{-1}}        & \mc{F}\mrm{R}S\ar[r]^{ (\mc{F}(\mrm{R}v))^{-1}}\ar[u]^{\mc{F}\lambda_S} &  \mc{F}(\mrm{R}^2Z) \ .}
$$
To see (I) and (III), just note that in $\mc{D}$ we have the
commutativity of
$$
\xymatrix@M=4pt@H=4pt@R=24pt@C=24pt{
 \mrm{R}X\ar[r]^{f} \ar[d]_{\mu_X} & T    \ar@{}[rd]|{;}        & \mrm{R}Z\ar[r]^{v} \ar[d]_{\mu_Z} & S  \\
 \mrm{R}^2X \ar[r]^{\mrm{R}f}     & \mrm{R}T\ar[u]^{\lambda_T} & \mrm{R}^2Z \ar[r]^{\mrm{R}v}     & \mrm{R}S\ar[u]^{\lambda_S} }
$$
Indeed, the statement follows from the equalities
$\lambda_T\comp\mrm{R}f=f\comp\lambda_{\mrm{R}X}$ and
$\lambda_{\mrm{R}X}\comp\mu_X=Id_{\mrm{R}X}$, and analogously for the right diagram.\\
The commutativity of (II) follows from the commutativity of the
diagrams bellow in $\mc{D}$ and $\cont$ respectively
$$
\xymatrix@M=4pt@H=4pt@R=24pt@C=24pt{
 T                                         & \mrm{R}Y\ar[r]^{g}\ar[l]_-{w}                                               & S \ar@{}[rd]|{;}           & \mc{F}(\mrm{R}T)\ar[rd]_-{\mc{F}I_T} & \mc{F}(\mrm{R}^2Y) \ar[r]^{\mc{F}(\mrm{R}g)} \ar[l]_{\mc{F}(\mrm{R}w)}& \mc{F}(\mrm{R}S)\ar[ld]^-{\mc{F}I_S} \ .\\
 \mrm{R}T \ar[u]_{\lambda_T}               & \mrm{R}^2Y \ar[r]^-{\mrm{R}g}\ar[u]_{\lambda_{\mrm{R}Y}}\ar[l]_{\mrm{R}w}  & \mrm{R}S\ar[u]^{\lambda_S} &                                      & \mc{F}cyl(w,g)                                                        &}
$$
Let us see that
$\mc{F}I_T\comp\mc{F}(\mrm{R}w)=\mc{F}I_S\comp\mc{F}(\mrm{R}g)$ in
$\cont$. From the property \ref{cilaumentado} of $cyl$ we deduce
the following diagram involving
$cyl(\mrm{R}Y)=cyl(Id_{\mrm{R}Y},Id_{\mrm{R}Y})$
$$
\xymatrix@M=4pt@H=4pt@R=15pt@C=24pt{
 \mrm{R}^2Y\ar[rd]_{I_{\mrm{R}Y}} \ar@/^1pc/[rrd]^{Id} &                           & \\
                                                      & cyl(\mrm{R}Y) \ar[r]^\rho  & \mrm{R^2Y} \ . \\
 \mrm{R}^2Y\ar[ru]^{J_{\mrm{R}Y}} \ar@/_1pc/[rru]_{Id} &                           & }
$$
Hence, $\rho\in\mrm{E}$ and since $\mc{F}\rho \comp
\mc{F}I_{\mrm{R}Y} = \mc{F}\rho \comp \mc{F}J_{\mrm{R}Y} =
Id_{\mrm{R}^2Y}$ then $\mc{F}I_{\mrm{R}Y}=\mc{F}J_{\mrm{R}Y}$ in
$\cont$.
On the other hand, by \ref{pdadesCylenD} we have a morphism
$H:cyl(\mrm{R}Y)\rightarrow cyl(w,g)$ such that $H\comp
J_{\mrm{R}Y}=I_S\comp\mrm{R}g$ and $H\comp I_{\mrm{R}Y}=I_T\comp
\mrm{R}w$. Applying $\mc{F}$ we deduce that $\mc{F}(I_T\comp
\mrm{R}w)=\mc{F}(I_S\comp\mrm{R}g)$.\\[0.2cm]
\textit{Uniqueness:} Assume that $\mc{G}':Ho\mc{D}\rightarrow
\cont$ is such that $\mc{G}'\comp\gamma=\mc{F}$. The equality
$\mc{G}'=\mc{G}$ is deduced from
$\mc{G}'\comp\gamma=\mc{G}\comp\gamma$, together with the
following lemma.
\end{proof}

\begin{lema}\label{Mof{i}smoHoD}
The morphism $\hat{f}$ in $Ho\mc{D}$ given by
$\xymatrix@M=4pt@H=4pt@C=15pt{ X & \mrm{R}X
\ar[l]_{\lambda_X}\ar[r]^{f} & T  &
\mrm{R}Y\ar[l]_-{w}\ar[r]^{\lambda_Y} & Y}$ is equal to
$\gamma(\lambda_Y)\comp(\gamma(w))^{-1}\comp\gamma(f)\comp(\gamma(\lambda_X))^{-1}$.
\end{lema}

\begin{proof}
Firstly, note that if $S$ is any object in $\mc{D}$, then
$\gamma(\lambda_S:\mrm{R}S\rightarrow S)$ is represented by
$\xymatrix@M=4pt@H=4pt@C=15pt{
 \mrm{R}S & \mrm{R}^2S \ar[l]_{\lambda_{\mrm{R}S}}\ar[r]^{\lambda_{\mrm{R}S}} & \mrm{R}S & \mrm{R}S\ar[l]_-{Id} \ar[r]^{\lambda_{S}} & S}$.\\
Indeed, by def{i}nition $\gamma(\lambda_S)$ is the class of
$\xymatrix@M=4pt@H=4pt@C=15pt{
 \mrm{R}S & \mrm{R}^2S \ar[l]_{\mrm{R}\lambda_{S}}\ar[r]^{\lambda_{\mrm{R}S}} & \mrm{R}S & \mrm{R}S\ar[l]_-{Id} \ar[r]^{\lambda_{S}} &
 S}$, and from the naturality of $\lambda$ it follows that
$\lambda_{\mrm{R}S}\comp\lambda_{S}=\mrm{R}\lambda_S\comp\lambda_S=\alpha$.
Therefore we have the following hammock relating both zig-zags
$$\xymatrix@M=4pt@H=4pt@C=15pt{
            & \mrm{R}^2S \ar[ld]_{Id}\ar[d]_{Id}\ar[r]^{\lambda_{\mrm{R}S}} & \mrm{R}S\ar[d]_{\lambda_S} & \mrm{R}S\ar[l]_-{Id} \ar[rd]^{Id}\ar[d]_{Id} & \\
 \mrm{R}^2S & \mrm{R}^2S \ar[l]_{Id}\ar[r]^{\alpha}                         & S                          & \mrm{R}S\ar[l]_{\lambda_S}\ar[r]^{Id}        & \mrm{R}S\ .\\
            & \mrm{R}^2S \ar[lu]^{Id}\ar[u]^{Id}\ar[r]^{\mrm{R}\lambda_{S}} & \mrm{R}S\ar[u]^{\lambda_S} & \mrm{R}S\ar[l]_-{Id} \ar[ru]_{Id}\ar[u]^{Id} &}
$$
Hence, $\gamma(f)\comp(\gamma(\lambda_X))^{-1}$ is represented by
the composition of the zig-zags
$$\xymatrix@M=4pt@H=4pt@C=15pt{
 X & \mrm{R}X \ar[l]_{\lambda_{X}}\ar[r]^{Id} & \mrm{R}X & \mrm{R}^2X\ar[l]_-{\lambda_{\mrm{R}X}} \ar[r]^{\lambda_{\mrm{R}X}} & \mrm{R} X \ar@{}[r]|{;} & \mrm{R}X & \mrm{R}^2X \ar[l]_{\lambda_{\mrm{R}X}}\ar[r]^{\mrm{R}f} & \mrm{R}T & \mrm{R}T\ar[l]_-{Id} \ar[r]^{\lambda_{T}} & T
}$$
that is by def{i}nition $\xymatrix@M=4pt@H=4pt@C=15pt{
 X & \mrm{R}X \ar[l]_{\lambda_{X}}\ar[r]^{\mu_X} & \mrm{R}^2X \ar[r]^-{I_{\mrm{R}X}} & cyl(\lambda_{\mrm{R}X},\mrm{R}f) & \mrm{R}^2T  \ar[l]_-{I_{\mrm{R}T}} &\mrm{R}T \ar[r]^{\lambda_{T}} \ar[l]_{\mu_T} & T}$.
From the equality
$\lambda_T\comp\mrm{R}f=g\comp\lambda_{\mrm{R}X}$ and the property
\ref{cilaumentado} of $cyl$ we get the small hammock
$$\xymatrix@M=4pt@H=4pt@C=15pt@R=9pt{
         & \mrm{R}X \ar[dd]_{Id}\ar[r]^{\mu_X}\ar[ld]_{Id} & \mrm{R}^2X \ar[r]^-{I_{\mrm{R}X}}\ar[dd]_{Id} & cyl(\lambda_{\mrm{R}X},\mrm{R}f)\ar[dd] & \mrm{R}^2T \ar[l]_-{I_{\mrm{R}T}}\ar[dd]_{\lambda_T} &\mrm{R}T \ar[rd]^{Id} \ar[l]_{\mu_T} \ar[dd]_{Id}& \\
\mrm{R}X &                                                 &                                               &                                         &                                                     &                                                 & \mrm{R}T\ .\\
         & \mrm{R}X \ar[lu]^{Id}\ar[r]^{\mu_X}             & \mrm{R}^2X \ar[r]^{\mrm{R}f}       & \mrm{R}T                                & \mrm{R}T\ar[l]_{Id}             & \mrm{R}T\ar[l]_{Id}  \ar[ru]_{Id}               &\\
}$$
Then $\gamma(f)\comp(\gamma(\lambda_X))^{-1}$ is represented by
$\xymatrix@M=4pt@H=4pt@C=15pt{   X & \mrm{R}X \ar[l]_{\lambda_X}\ar[r]^{\mu_X} & \mrm{R}^2X \ar[r]^{\mrm{R}f} & \mrm{R}T & \mrm{R}T\ar[l]_{Id}\ar[r]^{\lambda_T}  & T}$.\\
Compose with $\xymatrix@M=4pt@H=4pt@C=15pt{ T & \mrm{R}T
\ar[l]_{\lambda_T}\ar[r]^{Id} & \mrm{R}T & \mrm{R}^2Y
\ar[l]_{\mrm{R}w} \ar[r]^{\lambda_{\mrm{R}Y}}& \mrm{R}Y}$ (that
represents the morphism $[\gamma(w)]^{-1}$) and delete arrows in a
suitable way to obtain that
$[\gamma(w)]^{-1}\comp\gamma(f)\comp(\gamma(\lambda_X))^{-1}$ is
given by
$\xymatrix@M=4pt@H=4pt@C=15pt{X & \mrm{R}X \ar[l]_{\lambda_X}\ar[r]^{\mu_X} & \mrm{R}^2X\ar[r]^{\mrm{R}f} & \mrm{R}T   & \mrm{R}^2Y\ar[l]_{\mrm{R}w} \ar[r]^{\lambda_{\mrm{R}Y}}  & \mrm{R}Y}$.\\
To f{i}nish it suf{f}{i}ces to compose this zig-zag with
$\xymatrix@M=4pt@H=4pt@C=15pt{
 \mrm{R}Y & \mrm{R}^2Y \ar[l]_{\lambda_{\mrm{R}Y}}\ar[r]^{\lambda_{\mrm{R}Y}} & \mrm{R}Y & \mrm{R}Y\ar[l]_-{Id} \ar[r]^{\lambda_{Y}} & Y}$.
Again $\lambda_T\comp\mrm{R}w=w\comp\lambda_{\mrm{R}Y}$, and as
before we get
$$\xymatrix@M=4pt@H=4pt@C=14pt@R=9pt{
         & \mrm{R}X \ar[dd]_{Id}\ar[r]^{\mu_X}\ar[ld]_{Id} & \mrm{R}^2X \ar[r]^-{\mrm{R}\mu_X}\ar[dd]_{\lambda_{\mrm{R}X}} & \mrm{R}^3X \ar[r]^-{\mrm{R}^f}\ar[dd]_{\lambda_{\mrm{R}^2X}}  & \mrm{R}^2T \ar[r]^-{I_{\mrm{R}X}}\ar[dd]_{\lambda_T} & cyl(\mrm{R}w,\lambda_{\mrm{R}Y})\ar[dd] & \mrm{R}^2Y \ar[l]_-{I_{\mrm{R}Y}}\ar[dd]_{Id} & \mrm{R}Y \ar[rd]^{Id} \ar[l]_{\mu_Y} \ar[dd]_{Id}&  \\
\mrm{R}X &                                                 &                                                               &                                                               &                                                      &                                         &                                               &                                                  & \mrm{R}Y \\
         & \mrm{R}X \ar[lu]^{Id}\ar[r]^{Id}                & \mrm{R}X \ar[r]^{\mu_X}                                       & \mrm{R}^2X \ar[r]^{\mrm{R}f}                                  & \mrm{R}T   \ar[r]^{Id}                               & \mrm{R}T                                & \mrm{R}^2Y\ar[l]_{\mrm{R}w}                   & \mrm{R}Y \ar[l]_{\mu_Y}\ar[ru]_{Id}              &
}$$
Hence
$\gamma(\lambda_Y)\comp[\gamma(w)]^{-1}\comp\gamma(f)\comp(\gamma(\lambda_X))^{-1}$
is the class of the hammock induced by the lower row of the above
hammock, and from \ref{FrelRF} the required equality follows.
\end{proof}

\begin{cor}\label{objetoF{i}nal} A f{i}nal object in $\mc{D}$ is
also a f{i}nal object in $Ho\mc{D}$.
\end{cor}

\begin{proof} Indeed, if  $\xymatrix@M=4pt@H=4pt@C=15pt{
 X & \mrm{R}X \ar[l]_{\lambda_{X}}\ar[r]^{f} & T & \mrm{R}1\ar[l]_-{w} \ar[r]^{\lambda_{1}} &
 1}$ represents the morphism $\hat{f}:X\rightarrow 1$ in $Ho\mc{D}$, then
 $\hat{f}=\gamma(\rho)$, where $\rho:X\rightarrow 1$ is the trivial morphism in $\mc{D}$.
To see this, just consider the hammock
$$\xymatrix@M=4pt@H=4pt@C=15pt{
          & \mrm{R}X \ar[ld]_{Id}\ar[d]_{Id}\ar[r]^{\mrm{R}\rho} & \mrm{R}1\ar[d]_{\lambda_1} & \mrm{R}1\ar[l]_-{\mrm{R}w} \ar[rd]^{Id}\ar[d]_{Id} & \\
 \mrm{R}X & \mrm{R}X \ar[l]_{Id}\ar[r]                           & 1                          & \mrm{R}1\ar[l]_{\lambda_S}\ar[r]^{Id}        & \mrm{R}1\ .\\
          & \mrm{R}X \ar[lu]^{Id}\ar[u]^{Id}\ar[r]^{\mrm{R}f}    & \mrm{R}1\ar[u]^{\lambda_1} & \mrm{R}1\ar[l]_-{Id} \ar[ru]_{Id}\ar[u]^{Id} &}
$$
\end{proof}

\begin{cor}\label{objetoInicial} Assume that the trivial morphism $\sigma_0:0\rightarrow \mrm{R}0$ is an isomorphism in $\mc{D}$,
where $0$ is an initial object in $\mc{D}$. Then $0$ is also an
initial object in $Ho\mc{D}$.
\end{cor}

\begin{obs}\mbox{}\\
i) Since $\mbf{s}$ commutes with coproducts up to equivalence, we
deduce that $\sigma_0$ is always an equivalence, and
$\lambda_0\comp\sigma_0=Id$ because $0$ is initial.\\
ii) In the examples considered in this work, the simple functor
$\mbf{s}$ always commutes with coproducts (that is, the
transformation $\sigma$ of \ref{FuntorMonoidalCasiEstricto} is an
isomorphism). In particular, the hypothesis in the previous
corollary holds.
\end{obs}

\begin{proof}
Let $F:$ $\xymatrix@M=4pt@H=4pt@C=15pt{
 0 & \mrm{R}0 \ar[l]_{\lambda_{0}}\ar[r]^{f} & T & \mrm{R}X\ar[l]_-{w} \ar[r]^{\lambda_{X}} &
 X}$ be a zig-zag representing the morphism $\hat{f}:0\rightarrow X$ in $Ho\mc{D}$.\\
By assumption, $\mrm{R}0$ is isomorphic to $0$, so $\mrm{R}0$ is
an initial object in $\mrm{D}$. Then we have the following
commutative diagram
$$
\xymatrix@M=4pt@H=4pt@C=16pt@R=18pt{
 \mrm{R}0 \ar[r]^{f}          & T                    & \mrm{R}X \ar[l]_{w}       \\
 \mrm{R}0 \ar[r] \ar[u]_{Id}  & \mrm{R}X \ar[u]_{w}  & \mrm{R}X \ar[u]_{Id} \ar[l]_{Id}\ ,}
$$
that gives rise to a hammock relating $F$ to $\gamma(0\rightarrow
X)$.
\end{proof}

%
%
%
%

\section{Descent categories with $\lambda$ quasi-invertible}

In the description of $Ho\mc{D}$ given in the previous section,
the zig-zags that represent the morphisms in $Ho\mc{D}$ consists
of four arrows instead of two (that is the situation in the
calculus of fractions case). The reason is that the cylinder of
two morphisms
$A\stackrel{f}{\leftarrow}B\stackrel{g}{\rightarrow}C$ gives rise
to $\mrm{R}A{\rightarrow}cyl(f,g){\leftarrow}\mrm{R}C$,
so we need to attach $\lambda$s in order to recover $A$ and $C$.\\
If $\lambda$ is ``quasi-invertible'' this problem disappears, and
the description of the morphisms in $Ho\mc{D}$ become easier.\\

Throughout this section, $(\mc{D},\mrm{E},\mbf{s},\mu,\lambda)$
denotes a simplicial descent category.

\begin{def{i}}
We will say that $\lambda:\mrm{R}\rightarrow Id_{\mc{D}}$ is
\textit{quasi-invertible}\index{Index}{quasi-invertible} if there
exists a natural transformation $\rho:Id_{\mc{D}}\rightarrow
\mrm{R}$ such that $\lambda\comp \rho:Id_{\mc{D}}\rightarrow
Id_{\mc{D}}$ is the identity natural transformation. That is,
$\lambda_X\comp \rho_X=Id_X$ for every object $X$ in $\mc{D}$.
\end{def{i}}

An example of such situation is the category of chain complexes.

\begin{prop}\label{CylLambdaCasiinvertible}
Assume that $\lambda$ is quasi-invertible. Then the cylinder
object of two morphisms
$A\stackrel{f}{\leftarrow}B\stackrel{g}{\rightarrow}C$
in $\mc{D}$ has the following properties\\
\textbf{1)} there exists morphisms in $\mc{D}$, functorial in
$(f,g)$
$$ \jmath_A:A\rightarrow cyl(f,g) \ \ \jmath_B:B\rightarrow cyl(f,g)$$
such that $\jmath_A$ (resp. $\jmath_C$) is in $\mrm{E}$ if and
only if $g$ (resp. $f$) is so.\\
\textbf{2)} If $f=g=Id_A$, there exists an equivalence
$P:cyl(A)\rightarrow A$ in $\mc{D}$ such that the composition of
$P$ with the inclusions $\jmath_A,\jmath'_A:A\rightarrow cyl(A)$
given in 1) is equal to the identity $
A$.\\
\textbf{3)} there exists $H:cyl(A)\rightarrow cyl(f,g)$ such that
$H$ composed with $\jmath_A$ and $\jmath'_A$ is equal to
$\jmath_A\comp f$ and $\jmath_C\comp g$ respectively.
\end{prop}

\begin{proof}
As usual, 3) follows from 1).\\
To see 1), $\jmath_A$ and $\jmath_B$ are def{i}ned as the
compositions
$$\xymatrix@H=4pt@M=4pt{A\ar[r]^{\rho_A} & \mrm{R}A  \ar[r]^-{I_A} & cyl(f,g)& B\ar[r]^{\rho_B} & \mrm{R}B  \ar[r]^-{I_B} & cyl(f,g) }\ .$$
Since $\lambda_A\comp\rho_A=Id_A$ and $\lambda_A\in\mrm{E}$, we
deduce that $\rho_A\in\mrm{E}$. Hence, 1) follows from the
properties of the functor $cyl$.\\
\indent F{i}nally, if $\mc{D}$ is any descent category, there
exists $P':cyl(A)\rightarrow\mrm{R}A$ such that $P\comp I_A=P\comp
J_A=Id_{\mrm{R}A}$, where $I_A,J_A$ denote the canonical inclusions of $\mrm{R}A$ into $cyl(A)$.\\
Therefore, $P=\lambda_A\comp P'$ satisf{i}es 3) trivially.
\end{proof}

We can replace the maps $I_{A}:\mrm{R}A\rightarrow cyl(f,g)$
(resp. $I_C$) by $\jmath_A$ (resp. $\jmath_C$) in the proofs of
the previous section. In this way we obtain the following
proposition.

\begin{prop}
If $\mc{D}$ is a simplicial descent category with $\lambda$
quasi-invertible, then the morphisms in $Ho\mc{D}$ can be described as follows\\
Given objects $X$, $Y$ in $\mc{D}$,
$$ Hom_{Ho\mc{D}(X,Y)}=   \begin{array}{cc}\mrm{T}'(X,Y) &                                                           \\[-0.5cm]
                                                       &\hspace{-0.2cm} \mbox{\large{$\diagup$}} \!\! \sim \end{array}$$
where an element $F$ of
$\mrm{T}'(X,Y)$\index{Symbols}{$\mrm{T}(X,Y)$} is a zig-zag
$$\xymatrix@M=4pt@H=4pt{X \ar[r]^{{f}} & T & \ar[l]_w Y} ,\  w\in\mrm{E} \ .$$
If $\xymatrix@M=4pt@H=4pt{X \ar[r]^{{g}} & S & \ar[l]_u Y}$ is
another element $G\in\mrm{T}'(X,Y)$, then $F\sim G$ if and only if
there exists a hammock
\begin{equation}\label{hamacaLambdaCasiinv}\xymatrix@M=4pt@H=4pt{
     &  X \ar[ld]_{Id} \ar[r]^{{f}}\ar[d]    &  T\ar[d] &  Y\ar[l]_{w}\ar[rd]^{Id}\ar[d]  &         \\
 X   & \widetilde{X}       \ar[l]\ar[r]^{h}  & U       &   \widetilde{Y}  \ar[l]\ar[r]         & Y \ ,\\
     &  X\ar[lu]^{Id}\ar[r]^{{g}}\ar[u]      & S\ar[u] & Y\ar[l]_{u}\ar[ru]_{Id}\ar[u]  & }\end{equation}
relating $F$ to $G$ and where all maps except $f$, $g$ and $h$ are
equivalences.
\end{prop}

%
%
%
%
%
\section{Additive descent categories}

\begin{def{i}}\label{def{i}CatDescAditiva}
An additive simplicial descent category\index{Index}{additive
descent cat.} is a simplicial descent category that is also an
additive category and such that the simple functor is additive.\\
A functor of additive simplicial descent
categories\index{Index}{functor!of additive descent cat.} is a
functor of simplicial descent categories which is also additive.
\end{def{i}}

\begin{prop}\label{aditividad}
If $\mc{D}$ is an additive simplicial descent category then
$Ho\mc{D}$ is an additive category, and the localization functor
$\gamma:\mc{D}\rightarrow Ho\mc{D}$ is additive.\\
In addition, every functor $F:\mc{D}\rightarrow\mc{D'}$ of
additive simplicial descent categories gives rise to an additive
functor $HoF:Ho\mc{D}\rightarrow Ho\mc{D'}$.
\end{prop}

Throughout this section, we will assume that $\mc{D}$ is an
additive simplicial descent category.

\begin{lema}\label{denominadorComun}
If $\hat{f},\hat{g}:X\rightarrow Y$ are morphisms in $Ho\mc{D}$
then there exists zig-zags representing $\hat{f}$ and $\hat{g}$
with ``common denominator''
$$\xymatrix@M=4pt@H=4pt@R=11pt{
 \rho_{\hat{f}} : X & \mrm{R}X \ar[l]_-{\lambda_X}\ar[r]^{{f}} & L & \mrm{R}Y \ar[l]_w\ar[r]^{\lambda_{Y}}& Y\\
 \rho_{\hat{g}} : X & \mrm{R}X \ar[l]_-{\lambda_X}\ar[r]^{{g}} & L &\mrm{R}Y \ar[l]_w\ar[r]^{\lambda_{Y}}& Y.}$$
\end{lema}

\begin{proof}
Let $\hat{f}$ and $\hat{g}$ be the respective classes of the
following zig-zags $F$ and $G$
$$\xymatrix@M=4pt@H=4pt@R=15pt{
 X & \mrm{R}X \ar[l]_{\lambda_X}\ar[r]^{{f}'}  & T & \mrm{R}Y \ar[l]_-{u}\ar[r]^{\lambda_{Y}}& Y \ar@{}[r]|{;} & X & \mrm{R}X \ar[l]_{\lambda_X}\ar[r]^{{g}'} &
S &\mrm{R}Y \ar[l]_-{v}\ar[r]^{\lambda_{Y}}& Y\ .}$$
If $L=cyl(u,v)$, then $I_T:\mrm{R}T\rightarrow L$ and
$I_S:\mrm{R}S\rightarrow L$ are equivalences by
\ref{axiomacylenD}. Let $w:\mrm{R}Y\rightarrow L$ be the
equivalence def{i}ned as the composition
$$\xymatrix@M=4pt@H=4pt@R=15pt{\mrm{R}Y \ar[r]^{\mu_Y} & \mrm{R}^2Y \ar[r]^{\mrm{R}v} & \mrm{R}S \ar[r]^{I_S} & L \ .}$$
On one hand, by \ref{FrelFu} we have that $F$ is related to the
zig-zag $F^v$ given by $\xymatrix@M=4pt@H=4pt@R=15pt{
 X & \mrm{R}X \ar[l]_{\lambda_X} \ar[r]^{\mu_X} & \mrm{R}^2X \ar[r]^{\mrm{R}{f}'}  &  \mrm{R}T \ar[r]^-{I_T} & cyl(u,v) &\mrm{R}Y \ar[l]_-{w}\ar[r]^{\lambda_{Y}}&
 Y}$.\\
On the other hand, it is clear that the zig-zag $\mrm{R}G$ (see
\ref{FrelRF}) is related to
$$\xymatrix@M=4pt@H=4pt@R=15pt{ X & \mrm{R}X \ar[l]_{\lambda_X} \ar[r]^{\mu_X} & \mrm{R}^2X \ar[r]^{\mrm{R}g'}  &  \mrm{R}S \ar[r]^-{I_S} & L & \mrm{R}S \ar[l]_-{I_S} & \mrm{R}^2 Y \ar[l]_-{\mrm{R}v}  &\mrm{R}Y \ar[l]_{\mu_Y}\ar[r]^{\lambda_{Y}}& Y}\ .$$
Indeed, we have the hammock
$$\xymatrix@M=4pt@H=4pt@C=13pt@R=9pt{
         & \mrm{R}X \ar[dd]_{Id}\ar[r]^{\mu_X}\ar[ld]_{Id} & \mrm{R}^2X \ar[rr]^-{\mrm{R}g'}\ar[dd]_{Id} &                       & \mrm{R}S  \ar[dd]_{I_S}  &  \mrm{R}^2Y \ar[l]_-{\mrm{R}v}  & \mrm{R}Y \ar[rd]^{Id} \ar[l]_{\mu_Y} \ar[dd]_{Id} &  \\
\mrm{R}X &                                                 &                                             &                       &                          &                                 &                                                   & \mrm{R}Y\ . \\
         & \mrm{R}X \ar[lu]^{Id}\ar[r]^{\mu_X}             & \mrm{R}^2X \ar[r]^{\mrm{R}g'}               & \mrm{R}S\ar[r]^-{I_S} &  L                       &                                 & \mrm{R}Y   \ar[ll]_{w}  \ar[ru]_{Id}              &
}$$
Hence, the proof is f{i}nished.
\end{proof}

\begin{def{i}}[\textbf{Def{i}nition of sum in} $Ho\mc{D}$]\label{def{i}Suma}\index{Index}{sum in $Ho\mc{D}$}\mbox{}\\
Let $\hat{f},\hat{g}:X\rightarrow Y$ be morphisms in $Ho\mc{D}$
and $\rho_{\hat{f}}$, $\rho_{\hat{g}}$ be zig-zags representing
them as in \ref{denominadorComun}. We def{i}ne
$\hat{f}+\hat{g}:X\rightarrow Y$ as the class of
$$\xymatrix@M=4pt@H=4pt{\rho_{\hat{f}}+\rho_{\hat{g}}:&  X & \mrm{R}X \ar[l]_-{\lambda_X}\ar[r]^-{{f}+{g}} & L & \mrm{R}Y\ar[l]_w\ar[r]^-{\lambda_{Y}}& Y\ .}$$
\end{def{i}}

\begin{lema}
The sum of $f$ and $g$ in $Ho\mc{D}$ is well def{i}ned.
Equivalently, $\hat{f}+\hat{g}$ does not depend on the zig-zags
${\rho}_{\hat{f}}$ and ${\rho}_{\hat{g}}$ with ``common
denominator'' chosen for $f$ and $g$.
\end{lema}
\begin{proof}
Consider two dif{f}erent zig-zags $\rho_{\hat{f}}$,
$\rho_{\hat{f}'}$ representing $\hat{f}$, as well as
$\rho_{\hat{g}}$, $\rho_{\hat{g}'}$ representing $\hat{g}$
$$\xymatrix@M=4pt@H=4pt@R=10pt@C=15pt{
 \rho_{\hat{f}} : X & \mrm{R}X \ar[l]_-{\lambda_X}\ar[r]^{{f}} & S & \mrm{R}Y \ar[l]_u\ar[r]^{\lambda_{Y}}& Y \ar@{}[r]|{;} & \rho'_{\hat{f}} : X & \mrm{R}X \ar[l]_{\lambda_X}\ar[r]^{{f'}} & T & \mrm{R}Y \ar[l]_v\ar[r]^{\lambda_{Y}}& Y \\
 \rho_{\hat{g}} : X & \mrm{R}X \ar[l]_-{\lambda_X}\ar[r]^{{g}} & S &\mrm{R}Y \ar[l]_u\ar[r]^{\lambda_{Y}} & Y \ar@{}[r]|{;} & \rho'_{\hat{g}} : X & \mrm{R}X \ar[l]_{\lambda_X}\ar[r]^{{g'}} & T &\mrm{R}Y \ar[l]_v\ar[r]^{\lambda_{Y}} & Y .}$$
We must check that the following zig-zags are related
$$\xymatrix@M=4pt@H=4pt@R=10pt{
  \rho_{\hat{f}}+\rho_{\hat{g}}:  X & \mrm{R}X \ar[l]_-{\lambda_X}\ar[r]^-{{f}+{g}} & S & \mrm{R}Y\ar[l]_u\ar[r]^-{\lambda_{Y}}& Y\\
 \rho'_{\hat{f}}+\rho'_{\hat{g}}:  X & \mrm{R}X \ar[l]_-{\lambda_X}\ar[r]^-{{f}'+{g}'} & T & \mrm{R}Y\ar[l]_v\ar[r]^-{\lambda_{Y}}& Y\ .}
$$
We will see that both are related to a zig-zag
$\rho'_{\hat{f}}+\rho_{\hat{g}}$.\\
The equivalences $\mrm{R}u:\mrm{R}^2Y\rightarrow\mrm{R}S$ and
$\mrm{R}v:\mrm{R}^2Y\rightarrow\mrm{R}T$ induce equivalences
$$I_{\mrm{R}S}:\mrm{R}^2S\rightarrow cyl(\mrm{R}v,\mrm{R}u),\ I_{\mrm{R}T}:\mrm{R}^2T\rightarrow cyl(\mrm{R}v,\mrm{R}u)\ .$$
Consider the morphisms
$\widetilde{f}'=I_{\mrm{R}T}\comp\mrm{R}^2f'$,
$\widetilde{g}=I_{\mrm{R}S}\comp\mrm{R}^2g:\mrm{R}^3X\rightarrow
cyl(\mrm{R}v,\mrm{R}u)$. Then $\rho'_{\hat{f}}+\rho_{\hat{g}}$ is
def{i}ned by the zig-zag
$$\xymatrix@M=4pt@H=4pt@R=14pt@C=16pt{
 X & \mrm{R}X \ar[l]_-{\lambda_X} \ar[r]^{\hat{\mu}_X} & \mrm{R}^3X \ar[r]^-{{\widetilde{f}'}+\widetilde{g}} & cyl(\mrm{R}v,\mrm{R}u) & \mrm{R}^2S\ar[l]_-{I_{\mrm{R}S}} & \mrm{R}^3Y \ar[l]_{\mrm{R}^2u} &\mrm{R}Y \ar[l]_{\hat{\mu}_Y} \ar[r]^{\lambda_{Y}} & Y  }$$
where $\hat{\mu}_X=\mu_X\comp\mrm{R}\mu_X$, and
$\hat{\mu}_Y=\mu_Y\comp\mrm{R}\mu_Y$.\\
Consider the hammocks relating $\rho_{\hat{f}}$ to
$\rho_{\hat{f}'}$, and $\rho_{\hat{g}}$ to $\rho_{\hat{g}'}$
$$\xymatrix@M=4pt@H=4pt@C=14pt@R=14pt{
            & \mrm{R}^2X \ar[ld]_{Id}\ar[d]_{p}\ar[r]^{\mrm{R}f} & \mrm{R}S\ar[d]_{p'} & \mrm{R}^2Y\ar[l]_-{\mrm{R}u} \ar[rd]^{Id}\ar[d]_{p''} &            &            & \mrm{R}^2X \ar[ld]_{Id}\ar[d]_{s}\ar[r]^{\mrm{R}g}   & \mrm{R}S\ar[d]_{s'} & \mrm{R}^2Y\ar[l]_-{\mrm{R}u} \ar[rd]^{Id}\ar[d]_{s''} &           \\
 \mrm{R}^2X & \overline{X} \ar[l]_-{\alpha}\ar[r]^{h}            & M                   & \overline{Y}\ar[l]_{w}\ar[r]^-{\beta}                 & \mrm{R}^2Y & \mrm{R}^2X & \widetilde{X} \ar[l]_-{\alpha'}\ar[r]^{h'}           & N                   & \widetilde{Y}\ar[l]_{w'}\ar[r]^-{\beta'}              & \mrm{R}^2Y\ .\\
            & \mrm{R}^2X \ar[lu]^{Id}\ar[u]^{q}\ar[r]^{\mrm{R}f'}& \mrm{R}T\ar[u]^{q'} & \mrm{R}^2Y\ar[l]_-{\mrm{R}v} \ar[ru]_{Id}\ar[u]^{q''} &            &            & \mrm{R}^2X \ar[lu]^{Id}\ar[u]^{t}\ar[r]^{\mrm{R}g'}  & \mrm{R}T\ar[u]^{t'} & \mrm{R}^2Y\ar[l]_-{\mrm{R}v} \ar[ru]_{Id}\ar[u]^{t''} &            }
$$
We will denote them by $\mc{H}$ and $\mc{H'}$ respectively.\\
\textit{F{i}rst step:} $ \rho'_{\hat{f}}+\rho_{\hat{g}}$ is related to $ \rho_{\hat{f}}+\rho_{\hat{g}}$.\\
The hammocks $\mc{H}$ and $\mc{H}'$ give rise to the commutative
diagram
$$\xymatrix@M=4pt@H=4pt@C=14pt@R=14pt{
 \mrm{R}^2X \ar[d]_{p}\ar[r]^{\mrm{R}f} & \mrm{R}S\ar[d]_{p'} & \mrm{R}^2Y\ar[l]_-{\mrm{R}u} \ar[d]_{p''}\ar[r]^{\mrm{R}u} & \mrm{R}S\ar[d]_{Id} & \mrm{R}^2X\ar[l]_-{\mrm{R}g}\ar[d]_{p}   \\
 \overline{X} \ar[r]^{h}                & M                   & \overline{Y}\ar[l]_{w}                 \ar[r]^{{u'}}       & \mrm{R}S            & \overline{X}\ar[l]_-{l}   \\
 \mrm{R}^2X \ar[u]^{q}\ar[r]^{\mrm{R}f'}& \mrm{R}T\ar[u]^{q'} & \mrm{R}^2Y\ar[l]_-{\mrm{R}v} \ar[u]^{q''}\ar[r]^{\mrm{R}u} & \mrm{R}S\ar[u]^{Id} & \mrm{R}^2X\ar[l]_-{\mrm{R}g} \ar[u]^{q}  }
$$
where $u'$ is the composition
$\xymatrix@M=4pt@H=4pt@C=14pt{\overline{Y}\ar[r]^{\beta} &
\mrm{R}^2Y \ar[r]^{\mrm{R}u}& \mrm{R}S}$ and $l$ is
$\xymatrix@M=4pt@H=4pt@C=14pt{\overline{X}\ar[r]^{\alpha} &
\mrm{R}^2X \ar[r]^{\mrm{R}g}& \mrm{R}S}$. Then, applying the
functor $cyl$  to the two rows in the middle of the above diagram
we get
$$\xymatrix@M=4pt@H=4pt@C=14pt@R=14pt{
 \mrm{R}^3X \ar[d]_{\mrm{R}p}\ar[r]^{\mrm{R}^2f} & \mrm{R}^2S\ar[d]_{\mrm{R}p'} \ar[r]^-{J_{\mrm{R}S}} & cyl(\mrm{R}u,\mrm{R}u) \ar[d]_{r} & \mrm{R}^2S\ar[d]_{Id}\ar[l]_-{K_{\mrm{R}S}} & \mrm{R}^3X\ar[l]_-{\mrm{R}^2g}\ar[d]_{\mrm{R}p}   \\
 \mrm{R}\overline{X} \ar[r]^{\mrm{R}h}           & \mrm{R}M                     \ar[r]^-{I_{M}}        & cyl(u',w)                         & \mrm{R}^2S           \ar[l]_-{L_{\mrm{R}S}} & \mrm{R}\overline{X}\ar[l]_-{\mrm{R}l}   \\
 \mrm{R}^3X \ar[u]^{\mrm{R}q}\ar[r]^{\mrm{R}^2f'}& \mrm{R}^2T\ar[u]^{\mrm{R}q'} \ar[r]^-{I_{\mrm{R}T}} & cyl(\mrm{R}v,\mrm{R}u) \ar[u]^{r'}& \mrm{R}^2S\ar[u]^{Id}\ar[l]_-{I_{\mrm{R}S}} & \mrm{R}^3X\ar[l]_-{\mrm{R}^2g} \ar[u]^{\mrm{R}q}  }
$$
that becomes, after composing arrows, in
$$\xymatrix@M=4pt@H=4pt@C=14pt@R=14pt{
 \mrm{R}^3X \ar[d]_{\mrm{R}p}\ar[r]^-{\bar{f}}       & cyl(\mrm{R}u,\mrm{R}u) \ar[d]_{r}  & \mrm{R}^3X\ar[l]_-{\bar{g}}\ar[d]_{\mrm{R}p}   \\
 \mrm{R}\overline{X} \ar[r]^-{\bar{h}}               & cyl(u',w)                          & \mrm{R}\overline{X}\ar[l]_-{\bar{l}}   \\
 \mrm{R}^3X \ar[u]^{\mrm{R}q}\ar[r]^-{\widetilde{f}'}& cyl(\mrm{R}v,\mrm{R}u) \ar[u]^{r'} & \mrm{R}^3X\ar[l]_-{\widetilde{g}} \ar[u]^{\mrm{R}q}  .}
$$
Then, it holds that
$r\comp(\bar{f}+\bar{g})=(\bar{h}+\bar{l})\comp\mrm{R}p$ and
$r'\comp(\widetilde{f}'+\widetilde{g})=(\bar{h}+\bar{l})\comp\mrm{R}q$.
Therefore, the following diagram commutes
\begin{equation}\label{diagAuxiliar2}
\xymatrix@M=4pt@H=4pt@C=18pt@R=16pt{
 \mrm{R}^3X \ar[d]_{\mrm{R}p}\ar[r]^-{\bar{f}+\bar{g}}                & cyl(\mrm{R}u,\mrm{R}u) \ar[d]_{r} & \mrm{R}^2S\ar[d]_{Id}\ar[l]_-{K_{\mrm{R}S}} & \mrm{R}^3Y\ar[l]_-{\mrm{R}^2u} \ar[d]_{Id}            \\
 \mrm{R}\overline{X} \ar[r]^-{\bar{h}+\bar{l}}                        & cyl(u',w)                         & \mrm{R}^2S           \ar[l]_-{L_{\mrm{R}S}} & \mrm{R}^3Y\ar[l]_-{\mrm{R}^2u}             \\
 \mrm{R}^3X \ar[u]^{\mrm{R}q}\ar[r]^-{\widetilde{f}'+\widetilde{g}}   & cyl(\mrm{R}v,\mrm{R}u) \ar[u]^{r'}& \mrm{R}^2S\ar[u]^{Id}\ar[l]_-{I_{\mrm{R}S}} & \mrm{R}^3Y\ar[l]_-{\mrm{R}^2u} \ar[u]^{Id} .}
\end{equation}
As usual, we complete this diagram to a hammock through the
natural transformations  $\lambda$ and $\mu$. Indeed, consider the
diagrams
\begin{equation}\label{diagAuxiliar3}
\xymatrix@M=4pt@H=4pt@C=18pt@R=18pt{
             & \mrm{R}X \ar[r]^{\mu_X} \ar[d]_{\tau} \ar[ld]_{Id} & \mrm{R}^2X \ar[r]^{\mrm{R}\mu_X} \ar[d]_{\eta}       & \mrm{R}^3X \ar[d]_{\mrm{R}p}& & \mrm{R}^3Y \ar[d]_{Id}& \mrm{R}^2Y \ar[d]_{Id}\ar[l]_-{\mrm{R}\mu_Y}& \mrm{R}Y\ar[l]_-{\mu_Y}\ar[rd]^{Id}\ar[d]_{Id}  &             \\
 \mrm{R}X    & \mrm{R}\overline{X} \ar[l]_{\varrho}  \ar[r]^{Id}  & \mrm{R}\overline{X}              \ar[r]^{Id}         & \mrm{R}\overline{X}         & & \mrm{R}^3Y            & \mrm{R}^2Y            \ar[l]_-{\mrm{R}\mu_Y}& \mrm{R}Y\ar[l]_-{\mu_Y} \ar[r]^-{Id}            & \mrm{R}Y  \\
             & \mrm{R}X \ar[u]^{\tau'}\ar[lu]^{Id} \ar[r]^{\mu_X} & \mrm{R}^2X \ar[u]^{\eta'} \ar[r]^{\mrm{R}\mu_X}      & \mrm{R}^3X \ar[u]^{\mrm{R}q}& & \mrm{R}^3Y \ar[u]^{Id}& \mrm{R}^2Y \ar[u]^{Id}\ar[l]_-{\mrm{R}\mu_Y}& \mrm{R}Y\ar[l]_-{\mu_Y}\ar[ru]_{Id} \ar[u]^{Id} &              }
\end{equation}
where $\eta\!\!=\!\!\mrm{R}p\comp\mrm{R}\mu_X$,
$\eta'\!\!=\!\!\mrm{R}q\comp\mrm{R}\mu_X$,
$\tau\!\!=\!\!\mrm{R}p\comp\mrm{R}\mu_X\comp\mu_X$,
$\tau'\!\!=\!\!\mrm{R}q\comp\mrm{R}\mu_X\comp\mu_X$ and $\varrho\!\!=\!\!\lambda_{\mrm{R}X}\comp\mrm{R}\lambda_{\mrm{R}X}\comp\mrm{R}\alpha$.\\
Then, the squares in the left diagram commute by def{i}nition of
the arrows involved in them. In addition,
$\varrho\comp\tau=\varrho\comp\tau'=Id_{\mrm{R}X}$, since
$$\varrho\comp\tau=\lambda_{\mrm{R}X}\comp\mrm{R}\lambda_{\mrm{R}X}\comp(\mrm{R}\alpha\comp\mrm{R}p)\comp\mrm{R}\mu_X\comp\mu_X=
\lambda_{\mrm{R}X}\comp(\mrm{R}\lambda_{\mrm{R}X}\comp\mrm{R}\mu_X)\comp\mu_X=\lambda_{\mrm{R}X}\comp\mu_X=Id_{\mrm{R}X}.$$
The equality $\varrho\comp\tau'=Id_{\mrm{R}X}$ is checked
analogously. Therefore, the diagrams in (\ref{diagAuxiliar3}) are
commutative. Attaching them to (\ref{diagAuxiliar2}) we obtain the
hammock
$$
\xymatrix@M=4pt@H=4pt@C=16pt@R=16pt{
             & \mrm{R}X \ar[r]^{\hat{\mu}_X} \ar[d] \ar[ld]_{Id}  & \mrm{R}^3X \ar[d]\ar[r]^-{\bar{f}+\bar{g}}                & cyl(\mrm{R}u,\mrm{R}u) \ar[d] & \mrm{R}^2S\ar[d]\ar[l]_-{K_{\mrm{R}S}} & \mrm{R}^3Y\ar[l]_-{\mrm{R}^2u} \ar[d]  & \mrm{R}Y\ar[l]_-{\hat{\mu}_Y}\ar[rd]^{Id}\ar[d]  &             \\
 \mrm{R}X  & \mrm{R}\overline{X} \ar[l]  \ar[r]                 & \mrm{R}\overline{X} \ar[r]                                & cyl(u',w)                     & \mrm{R}^2S           \ar[l]            & \mrm{R}^3Y\ar[l]                       & \mrm{R}Y\ar[l]_-{\hat{\mu}_Y} \ar[r]                    & \mrm{R}Y  \\
             & \mrm{R}X \ar[u]\ar[lu]^{Id} \ar[r]^{\hat{\mu}_X}   & \mrm{R}^3X \ar[u]\ar[r]^-{\widetilde{f}'+\widetilde{g}}   & cyl(\mrm{R}v,\mrm{R}u) \ar[u] & \mrm{R}^2S\ar[u]\ar[l]_-{I_{\mrm{R}S}} & \mrm{R}^3Y\ar[l]_-{\mrm{R}^2u} \ar[u]  & \mrm{R}Y\ar[l]_-{\hat{\mu}_Y}\ar[ru]_{Id} \ar[u] &              }
$$
relating $\rho'_{\hat{f}}+\rho_{\hat{g}}$ to the zig-zag
$\widetilde{\rho}$
$$
\xymatrix@M=4pt@H=4pt@C=16pt@R=16pt{
   X          & \mrm{R}X \ar[r]^{\hat{\mu}_X} \ar[l]_{\lambda_X}  & \mrm{R}^3X \ar[r]^-{\bar{f}+\bar{g}}    & cyl(\mrm{R}u,\mrm{R}u)& \mrm{R}^2S\ar[l]_-{K_{\mrm{R}S}} & \mrm{R}^3Y\ar[l]_-{\mrm{R}^2u}  & \mrm{R}Y\ar[l]_-{\hat{\mu}_Y}\ar[r]^{\lambda_Y} & Y .}
$$
On the other hand, applying twice \ref{FrelRF}, it follows that
$\rho_{\hat{f}}+\rho_{\hat{g}}$ is related to the zig-zag
$\mrm{R}^2(\rho_{\hat{f}}+\rho_{\hat{g}})$, given by
$$\xymatrix@M=4pt@H=4pt@C=18pt{
  X & \mrm{R}X \ar[l]_-{\lambda_X} \ar[r]^{\hat{\mu}_X} & \mrm{R}^3X \ar[r]^-{\mrm{R}^2({f}+{g})} &  \mrm{R}^2S & \mrm{R}^3Y\ar[l]_{ \mrm{R}^2u} &  \mrm{R}Y \ar[l]_{\hat{\mu}_Y} \ar[r]^-{\lambda_{Y}} & Y }
$$
and $\mrm{R}^2(f+g)=\mrm{R}^2f + \mrm{R}^2g$ since $\mrm{R}$ is
additive. In addition, by the properties of the cylinder functor
we have an equivalence
$\theta:cyl(\mrm{R}u,\mrm{R}u)\rightarrow\mrm{R}^2S$ such that
$\theta\comp K_{\mrm{R}S}=\theta\comp
J_{\mrm{R}S}=Id_{\mrm{R}^2S}$.\\
Hence
$\theta\comp(\bar{f}+\bar{g})=\theta\comp(J_{\mrm{R}S}\comp\mrm{R}^2f+K_{\mrm{R}S}\comp\mrm{R}^2g)=\mrm{R}^2f+\mrm{R}^2g$.
Therefore, we get the following hammock relating
$\widetilde{\rho}$ to $\mrm{R}^2(\rho_{\hat{f}}+\rho_{\hat{g}})$
$$
\xymatrix@M=4pt@H=4pt@C=16pt@R=10pt{
             & \mrm{R}X \ar[r]^{\hat{\mu}_X} \ar[dd]_{Id} \ar[ld]_{Id}  & \mrm{R}^3X \ar[r]^-{\bar{f}+\bar{g}}\ar[dd]_{Id} & cyl(\mrm{R}u,\mrm{R}u) \ar[dd]_{\theta} & \mrm{R}^2S\ar[dd]_{Id}\ar[l]_-{K_{\mrm{R}S}} & \mrm{R}^3Y\ar[l]_-{\mrm{R}^2u} \ar[dd]_{Id} & \mrm{R}Y\ar[l]_-{\hat{\mu}_Y}\ar[rd]^{Id}\ar[dd]_{Id}  &             \\
 \mrm{R}^2X  &                                                          &                                                  &                                         &                                              &                                             &                                                        & \mrm{R}Y  \\
             & \mrm{R}X \ar[lu]^-{Id} \ar[r]^{\hat{\mu}_X}              & \mrm{R}^3X \ar[r]^-{\mrm{R}^2({f}+{g})}          &  \mrm{R}^2S                             & \mrm{R}^2S           \ar[l]_{Id}             & \mrm{R}^3Y\ar[l]_{ \mrm{R}^2u}              & \mrm{R}Y \ar[l]_{\hat{\mu}_Y} \ar[ru]_-{Id}    &  }
$$
that f{i}nishes the proof of the f{i}rst step.\\[0.2cm]
\textit{Second step:} $\rho'_{\hat{f}}+\rho'_{\hat{g}}$ is related to $\rho'_{\hat{f}}+\rho_{\hat{g}}$.\\
Consider this time the following commutative diagram induced by
the hammocks $\mc{H}$ and $\mc{H}'$
$$\xymatrix@M=4pt@H=4pt@C=14pt@R=14pt{
 \mrm{R}^2X \ar[d]_{s}\ar[r]^{\mrm{R}g} & \mrm{R}S\ar[d]_{s'} & \mrm{R}^2Y\ar[l]_-{\mrm{R}u} \ar[d]_{s''}\ar[r]^{\mrm{R}v} & \mrm{R}T\ar[d]_{Id} & \mrm{R}^2X\ar[l]_-{\mrm{R}f'}\ar[d]_{s}   \\
 \widetilde{X} \ar[r]^{h'}                & N                   & \widetilde{Y}\ar[l]_{w'}                 \ar[r]^{{v'}}       & \mrm{R}T            & \widetilde{X}\ar[l]_-{l'}   \\
 \mrm{R}^2X \ar[u]^{t}\ar[r]^{\mrm{R}g'}& \mrm{R}T\ar[u]^{t'} & \mrm{R}^2Y\ar[l]_-{\mrm{R}v} \ar[u]^{t''}\ar[r]^{\mrm{R}v} & \mrm{R}T\ar[u]^{Id} & \mrm{R}^2X\ar[l]_-{\mrm{R}f'} \ar[u]^{t}  }
$$
where $v'$ is the composition
$\xymatrix@M=4pt@H=4pt@C=14pt{\widetilde{Y}\ar[r]^-{\beta'} &
\mrm{R}^2Y \ar[r]^{\mrm{R}v}& \mrm{R}T}$ and $l'$ is
$\xymatrix@M=4pt@H=4pt@C=14pt{\widetilde{X}\ar[r]^-{\alpha'} &
\mrm{R}^2X \ar[r]^{\mrm{R}f'}& \mrm{R}T}$. Again, applying
 $cyl$ (this time changing the order of the arrows) we obtain
$$\xymatrix@M=4pt@H=4pt@C=14pt@R=14pt{
 \mrm{R}^3X \ar[d]_{\mrm{R}s}\ar[r]^{\mrm{R}^2g} & \mrm{R}^2S\ar[d]_{\mrm{R}s'} \ar[r]^-{I_{\mrm{R}S}} & cyl(\mrm{R}v,\mrm{R}u) \ar[d]_{k} & \mrm{R}^2T\ar[d]_{Id}\ar[l]_-{I_{\mrm{R}T}} & \mrm{R}^3X\ar[l]_-{\mrm{R}^2f'}\ar[d]_{\mrm{R}s}   \\
 \mrm{R}\widetilde{X} \ar[r]^{\mrm{R}h'}           & \mrm{R}N                     \ar[r]^-{I_{N}}      & cyl(v',w')                        & \mrm{R}^2T           \ar[l]_-{L_{\mrm{R}T}} & \mrm{R}\widetilde{X}\ar[l]_-{\mrm{R}l'}   \\
 \mrm{R}^3X \ar[u]^{\mrm{R}t}\ar[r]^{\mrm{R}^2g'}& \mrm{R}^2T\ar[u]^{\mrm{R}t'} \ar[r]^-{K_{\mrm{R}T}} & cyl(\mrm{R}v,\mrm{R}v) \ar[u]^{k'}& \mrm{R}^2T\ar[u]^{Id}\ar[l]_-{J_{\mrm{R}T}} & \mrm{R}^3X\ar[l]_-{\mrm{R}^2f'} \ar[u]^{\mrm{R}t}  }
$$
that becomes, after composing arrows, in
$$\xymatrix@M=4pt@H=4pt@C=14pt@R=14pt{
 \mrm{R}^3X \ar[d]_{\mrm{R}s}\ar[r]^-{\widetilde{g}}   & cyl(\mrm{R}v,\mrm{R}u) \ar[d]_{k}   & \mrm{R}^3X\ar[l]_-{\widetilde{f}'}\ar[d]_{\mrm{R}s}   \\
 \mrm{R}\widetilde{X} \ar[r]^-{\check{h}}              & cyl(v',w')                          & \mrm{R}\widetilde{X}\ar[l]_-{\check{l}}   \\
 \mrm{R}^3X \ar[u]^{\mrm{R}t}\ar[r]^-{\check{g}'}       &  cyl(\mrm{R}v,\mrm{R}v) \ar[u]^{k'}& \mrm{R}^3X\ar[l]_-{\check{f}'} \ar[u]^{\mrm{R}t}  ,}
$$
and again we deduce the commutative diagram
$$
\xymatrix@M=4pt@H=4pt@C=18pt@R=16pt{
 \mrm{R}^3X \ar[d]_{\mrm{R}s}\ar[r]^-{\widetilde{f}'+\widetilde{g}}& cyl(\mrm{R}v,\mrm{R}u) \ar[d]_{k} & \mrm{R}^2T\ar[d]_{Id}\ar[l]_-{I_{\mrm{R}T}} & \mrm{R}^3Y\ar[l]_-{\mrm{R}^2v} \ar[d]_{Id}            \\
 \mrm{R}\widetilde{X} \ar[r]^-{\check{h}+\check{l}}               & cyl(v',w')                        & \mrm{R}^2T           \ar[l]_-{L_{\mrm{R}T}} & \mrm{R}^3Y\ar[l]_-{\mrm{R}^2v}             \\
 \mrm{R}^3X \ar[u]^{\mrm{R}t}\ar[r]^-{\check{f}'+\check{g}'}      & cyl(\mrm{R}v,\mrm{R}v) \ar[u]^{k'}& \mrm{R}^2T\ar[u]^{Id}\ar[l]_-{J_{\mrm{R}T}} & \mrm{R}^3Y\ar[l]_-{\mrm{R}^2v} \ar[u]^{Id} .}
$$
As in the previous step, it is possible to complete this diagram
to the hammock
$$
\xymatrix@M=4pt@H=4pt@C=16pt@R=16pt{
             & \mrm{R}X \ar[r]^{\hat{\mu}_X} \ar[d] \ar[ld]_{Id}   & \mrm{R}^3X \ar[d]_{\mrm{R}s}\ar[r]^-{\widetilde{f}'+\widetilde{g}}& cyl(\mrm{R}v,\mrm{R}u) \ar[d]_{k} & \mrm{R}^2T\ar[d]_{Id}\ar[l]_-{I_{\mrm{R}T}} & \mrm{R}^3Y\ar[l]_-{\mrm{R}^2v} \ar[d]_{Id} & \mrm{R}Y\ar[l]_-{\hat{\mu}_Y}\ar[rd]^{Id}\ar[d]  &             \\
 \mrm{R}X    & \mrm{R}\widetilde{X} \ar[l]  \ar[r]                 & \mrm{R}\widetilde{X} \ar[r]^-{\check{h}+\check{l}}               & cyl(v',w')                        & \mrm{R}^2T           \ar[l]_-{L_{\mrm{R}T}} & \mrm{R}^3Y\ar[l]_-{\mrm{R}^2v}             & \mrm{R}Y\ar[l]_-{\hat{\mu}_Y} \ar[r]                    & \mrm{R}Y  \\
             & \mrm{R}X \ar[u]\ar[lu]^{Id} \ar[r]^{\hat{\mu}_X}    & \mrm{R}^3X \ar[u]^{\mrm{R}t}\ar[r]^-{\check{f}'+\check{g}'}      & cyl(\mrm{R}v,\mrm{R}v) \ar[u]^{k'}& \mrm{R}^2T\ar[u]^{Id}\ar[l]_-{J_{\mrm{R}T}} & \mrm{R}^3Y\ar[l]_-{\mrm{R}^2v} \ar[u]^{Id} & \mrm{R}Y\ar[l]_-{\hat{\mu}_Y}\ar[ru]_{Id} \ar[u] &              }
$$
relating the zig-zags $\widehat{\rho}$ and $\overline{\rho}$
consisting respectively of the top and bottom rows of this
hammock.\\
Note that $\widehat{\rho}\sim \rho'_{\hat{f}}+\rho_{\hat{g}}$.
Indeed, we have the hammock
$$
\xymatrix@M=4pt@H=4pt@C=16pt@R=16pt{
             & \mrm{R}X \ar[r]^{\hat{\mu}_X} \ar[d] \ar[ld]_{Id}   & \mrm{R}^3X \ar[d]_{Id}\ar[r]^-{\widetilde{f}'+\widetilde{g}}  & cyl(\mrm{R}v,\mrm{R}u) \ar[d]_{Id} & \mrm{R}^2T\ar[l]_-{I_{\mrm{R}T}} & \mrm{R}^3Y\ar[l]_-{\mrm{R}^2v} \ar[d]_{I_{\mrm{R}^2Y}} & \mrm{R}Y\ar[l]_-{\hat{\mu}_Y}\ar[rd]^{Id}\ar[d]  &             \\
 \mrm{R}X    & \mrm{R}\widetilde{X} \ar[l]  \ar[r]                 & \mrm{R}^3X \ar[r]^-{\widetilde{f}'+\widetilde{g}}             & cyl(\mrm{R}v,\mrm{R}u)             &                                  & cyl(\mrm{R}^2Y)\ar[ll]_-{H}                            & cyl(\mrm{R}^2Y)\ar[l]_-{Id} \ar[r]               & \mrm{R}Y  \\
             & \mrm{R}X \ar[u]\ar[lu]^{Id} \ar[r]^{\hat{\mu}_X}    & \mrm{R}^3X \ar[u]^{Id}\ar[r]^-{\widetilde{f}'+\widetilde{g}}  & cyl(\mrm{R}v,\mrm{R}u) \ar[u]^{Id} & \mrm{R}^2S\ar[l]_-{I_{\mrm{R}S}} & \mrm{R}^3Y\ar[l]_-{\mrm{R}^2u} \ar[u]^{J_{\mrm{R}^2Y}} & \mrm{R}Y\ar[l]_-{\hat{\mu}_Y}\ar[ru]_{Id} \ar[u] &              }
$$
where $H$ is the morphism provided by \ref{pdadesCylenD}. If
$\varphi:cyl(\mrm{R}^2Y)\rightarrow \mrm{R}^3Y$ is such that
$\varphi\comp I_{\mrm{R}^2Y}=\varphi\comp
J_{\mrm{R}^2Y}=Id_{\mrm{R}^3Y}$ (see \ref{cilaumentado}), and
$\hat{\lambda}_X=\mrm{R}\lambda_{\mrm{R}X}\comp\lambda_{\mrm{R}X}$,
then the triangle in the right hand side of the previous hammock
consists of
$$
\xymatrix@M=4pt@H=4pt@C=15pt@R=15pt{
                                                      &                                    & \mrm{R} X                               &                                        &                                                         \\
                                                      &                                    & \mrm{R}^3 X \ar[u]^{\hat{\lambda}_{X}}  &                                        &                                                          \\
 \mrm{R}Y \ar[r]^{\hat{\mu}_X} \ar@<1.4ex>[rruu]^{Id} & \mrm{R}^3 X \ar[r]^{I_{\mrm{R}^2Y}}& cyl(\mrm{R}^2Y) \ar[u]^{\varphi}        &  \mrm{R}^3 {X} \ar[l]_{J_{\mrm{R}^2Y}} & \mrm{R} X  \ar@<-1.4ex>[lluu]_{Id}\ar[l]_{\hat{\mu}_X} \ .}
$$
Then $\widehat{\rho}\sim \rho'_{\hat{f}}+\rho_{\hat{g}}$, and the
fact $\overline{\rho}\sim \rho'_{\hat{f}}+\rho'_{\hat{g}}$ can be
proved as in the previous step, so we are done.
\end{proof}
\begin{proof}[\textbf{Proof of \ref{aditividad}}]
Let us prove that the axioms of additive category are hold in
$Ho\mc{D}$. We follow here the presentation given for additive categories in \cite{GM}.\\[0.2cm]
\textsl{Axiom $($A1$)$:} The sum in $Ho\mc{D}$ clearly makes each
$\mrm{Hom}_{Ho\mc{D}}(X,Y)$ into an abelian group.\\
Let us check that the composition $\mrm{Hom}_{Ho\mc{D}}(Z,X)\times
\mrm{Hom}_{Ho\mc{D}}(X,Y)\!\rightarrow\!
\mrm{Hom}_{Ho\mc{D}}(Z,Y)$
is bilineal, that is, $\hat{h}\comp(\hat{f}+\hat{g})=\hat{h}\comp\hat{f}+\hat{h}\comp\hat{g}$; $(\hat{f}+\hat{g})\comp\hat{h}=\hat{f}\comp\hat{h}+\hat{g}\comp\hat{h}$.\\
We will see that
$(\hat{f}+\hat{g})\comp\hat{h}=\hat{f}\comp\hat{h}+\hat{g}\comp\hat{h}$
(the other equality can be proved similarly).\\
Consider the morphisms $\hat{f},\hat{g}:X\rightarrow Y$ and
$\hat{h}:Z\rightarrow X$ in $Ho\mc{D}$. As we saw before, we can
assume that these morphisms are represented respectively by
$$\xymatrix@M=4pt@H=4pt@R=10pt{
 \rho_{\hat{f}}:  X & \mrm{R}X \ar[l]_-{\lambda_X}\ar[r]^-{{f}} & T & \mrm{R}Y \ar[l]_{w}\ar[r]^-{\lambda_{Y}}& Y\\
 \rho_{\hat{g}}:  X & \mrm{R}X \ar[l]_-{\lambda_X}\ar[r]^-{{g}} & T & \mrm{R}Y \ar[l]_{w}\ar[r]^-{\lambda_{Y}}& Y\\
 \rho_{\hat{h}}:  Z & \mrm{R}Z \ar[l]_-{\lambda_Z}\ar[r]^-{{h}} & L & \mrm{R}X \ar[l]_{t}\ar[r]^-{\lambda_{X}}& X.}$$
The functor $cyl$ provides the diagrams
$$\xymatrix@M=4pt@H=4pt{\mrm{R}^{2}X \ar[r]^{\mrm{R}{f}}\ar[d]_{\mrm{R}t} & \mrm{R}T \ar[d]^{I_T} & \ar@{}[d]|{;}&   \mrm{R}^{2}X \ar[r]^{\mrm{R}{g}}\ar[d]_{\mrm{R}t} & \mrm{R}T \ar[d]^{J_T}\\
                        \mrm{R}L \ar[r]^{I_L}                             & C                     &              &   \mrm{R}L \ar[r]^{J_L}                             & C' \ ,}$$
where the arrows $i_T$, $j_T\in\mrm{E}$. Again, we have the
commutative diagram in $Ho\mc{D}$
$$\xymatrix@M=4pt@H=4pt{\mrm{R}^{2}T \ar[r]^{\mrm{R}I_T}\ar[d]_{\mrm{R}J_T} & \mrm{R}C\ar[d]^{I_C} \\
                        \mrm{R}C' \ar[r]^{I_{C'}}                           & N\ ,}$$
where all maps are in $\mrm{E}$.\\
Let $u:\mrm{R}^{3}Y\rightarrow N$ be the composition
$\xymatrix@M=4pt@H=4pt{\mrm{R}^{3}Y\ar[r]^{\mrm{R}^2w} &
\mrm{R}^2T \ar[r]^{\mrm{R}I_T} & \mrm{R}C \ar[r]^{I_C}& N}$. Set
${f}'=I_C\comp\mrm{R}I_L$ and ${g}'=I_{C'}\comp\mrm{R}J_L$.
Assume that the following zig-zags represent $\hat{f}$ and
$\hat{g}$ respectively
$$\xymatrix@M=4pt@H=4pt@C=20pt@R=10pt{
 \rho'_{\hat{f}}: X & \mrm{R}X\ar[r]^{\mu_X} \ar[l]_-{\lambda_X} & \mrm{R}^2X\ar[r]^{\mrm{R}\mu_X} & \mrm{R}^{3}X \ar[r]^-{{f}'\comp\mrm{R}^2t} & N & \mrm{R}^{3}Y \ar[l]_{u} & \mrm{R}^2Y\ar[l]_{\mrm{R}\mu_Y} & \mrm{R}Y\ar[l]_{\mu_Y} \ar[r]^-{\lambda_{Y}}& Y\\
 \rho'_{\hat{g}}: X & \mrm{R}X\ar[r]^{\mu_X} \ar[l]_-{\lambda_X} & \mrm{R}^2X\ar[r]^{\mrm{R}\mu_X} & \mrm{R}^{3}X \ar[r]^-{{g}'\comp\mrm{R}^2t} & N & \mrm{R}^{3}Y \ar[l]_{u} & \mrm{R}^2Y\ar[l]_{\mrm{R}\mu_Y} & \mrm{R}Y\ar[l]_{\mu_Y} \ar[r]^-{\lambda_{Y}}& Y \ .}$$
In this case, $\hat{f}+\hat{g}$ is given by the zig-zag
$$\xymatrix@M=4pt@H=4pt@C=17pt{
 X & \mrm{R}X\ar[r]^{\mu_X} \ar[l]_-{\lambda_X} & \mrm{R}^2X\ar[r]^{\mrm{R}\mu_X} & \mrm{R}^{3}X \ar[r]^{\mrm{R}^2t} &\mrm{R}^2L \ar[r]^-{{f}'+g'} & N & \mrm{R}^{3}Y \ar[l]_{u} & \mrm{R}^2Y\ar[l]_{\mrm{R}\mu_Y} & \mrm{R}Y\ar[l]_{\mu_Y} \ar[r]^-{\lambda_{Y}}& Y}$$
and if we compose it with the zig-zag $\mrm{R}^2\rho_{\hat{h}}$
representing $\hat{h}$
$$\xymatrix@M=4pt@H=4pt@C=19pt{
 Z & \mrm{R}Z\ar[r]^{\mu_Z} \ar[l]_-{\lambda_Z} & \mrm{R}^2Z\ar[r]^{\mrm{R}\mu_Z} & \mrm{R}^{3}Z \ar[r]^{\mrm{R}^2h} &\mrm{R}^2L  & \mrm{R}^{3}X \ar[l]_{\mrm{R}^2t} & \mrm{R}^2X\ar[l]_{\mrm{R}\mu_X} & \mrm{R}X\ar[l]_{\mu_X} \ar[r]^-{\lambda_{X}}& X}
$$
then by \ref{EliminarFlechasComp} we can delete arrows getting
that $(\hat{f}+\hat{g})\comp\hat{h}$ is given by
$$\xymatrix@M=4pt@H=4pt@C=17pt{
 Z & \mrm{R}Z\ar[r]^{\mu_Z} \ar[l]_-{\lambda_Z} & \mrm{R}^2Z\ar[r]^{\mrm{R}\mu_Z} & \mrm{R}^{3}Z \ar[r]^{\mrm{R}^2h} &\mrm{R}^2L \ar[r]^-{{f}'+g'} & N & \mrm{R}^{3}Y \ar[l]_{u} & \mrm{R}^2Y\ar[l]_{\mrm{R}\mu_Y} & \mrm{R}Y\ar[l]_{\mu_Y} \ar[r]^-{\lambda_{Y}}& Y}
$$
On the other hand, we can use again $\mrm{R}^2\rho_{\hat{h}}$ to
compute $\hat{f}\comp \hat{h}$ and $\hat{g}\comp \hat{h}$. After
deleting arrows, they are given by
$$\xymatrix@M=4pt@H=4pt@C=17pt@R=15pt{
 Z & \mrm{R}Z\ar[r]^{\mu_Z} \ar[l]_-{\lambda_Z} & \mrm{R}^2Z\ar[r]^{\mrm{R}\mu_Z} & \mrm{R}^{3}Z \ar[r]^{\mrm{R}^2h} &\mrm{R}^2L \ar[r]^-{{f}'} & N & \mrm{R}^{3}Y \ar[l]_{u} & \mrm{R}^2Y\ar[l]_{\mrm{R}\mu_Y} & \mrm{R}Y\ar[l]_{\mu_Y} \ar[r]^-{\lambda_{Y}}& Y \\
 Z & \mrm{R}Z\ar[r]^{\mu_Z} \ar[l]_-{\lambda_Z} & \mrm{R}^2Z\ar[r]^{\mrm{R}\mu_Z} & \mrm{R}^{3}Z \ar[r]^{\mrm{R}^2h} &\mrm{R}^2L \ar[r]^-{g'} & N & \mrm{R}^{3}Y \ar[l]_{u} & \mrm{R}^2Y\ar[l]_{\mrm{R}\mu_Y} & \mrm{R}Y\ar[l]_{\mu_Y} \ar[r]^-{\lambda_{Y}}& Y}$$
and therefore their sum is $(\hat{f}+\hat{g})\comp\hat{h}$, so in this case (A1) holds.\\
It remains to see that $\rho'_{\hat{f}}$ and $\rho'_{\hat{g}}$
represent $\hat{f}$ and $\hat{g}$ respectively.
Considering the zig-zags $\rho_{\hat{f}}$ and $\rho_{\hat{g}}$,
and applying \ref{FrelFs} to the equivalence
$t:\mrm{R}X\rightarrow L$, it follows that $\hat{f}$ and $\hat{g}$
are given by
$$\xymatrix@M=4pt@H=4pt@C=23pt@R=15pt{
 X & \mrm{R}X\ar[r]^{\mu_X} \ar[l]_-{\lambda_X} & \mrm{R}^2X \ar[r]^-{I_L\comp\mrm{R}t} & C  &  \mrm{R}^2Y\ar[l]_{I_T\comp\mrm{R}w} & \mrm{R}Y\ar[l]_{\mu_Y} \ar[r]^-{\lambda_{Y}}& Y\\
 X & \mrm{R}X\ar[r]^{\mu_X} \ar[l]_-{\lambda_X} & \mrm{R}^2X \ar[r]^-{J_L\comp\mrm{R}t} & C' &  \mrm{R}^2Y\ar[l]_{J_T\comp\mrm{R}w} & \mrm{R}Y\ar[l]_{\mu_Y} \ar[r]^-{\lambda_{Y}}& Y  .}$$
Set $\tau=\mrm{R}I_T\comp\mrm{R}^2w$ and
$\tau'=\mrm{R}J_T\comp\mrm{R}^2w$. By \ref{FrelRF} these zig-zags
are related respectively to
$$\xymatrix@M=4pt@H=4pt@C=17pt@R=15pt{
 X & \mrm{R}X\ar[r]^{\mu_X} \ar[l]_-{\lambda_X} & \mrm{R}^2X\ar[r]^{\mrm{R}\mu_X} & \mrm{R}^{3}X \ar[r]^-{\mrm{R}I_L\comp\mrm{R}^2t} & \mrm{R}C  & \mrm{R}^{3}Y \ar[l]_{\tau} & \mrm{R}^2Y\ar[l]_{\mrm{R}\mu_Y} & \mrm{R}Y\ar[l]_{\mu_Y} \ar[r]^-{\lambda_{Y}}& Y\\
 X & \mrm{R}X\ar[r]^{\mu_X} \ar[l]_-{\lambda_X} & \mrm{R}^2X\ar[r]^{\mrm{R}\mu_X} & \mrm{R}^{3}X \ar[r]^-{\mrm{R}J_L\comp\mrm{R}^2t} & \mrm{R}C' & \mrm{R}^{3}Y \ar[l]_{\tau'} & \mrm{R}^2Y\ar[l]_{\mrm{R}\mu_Y} & \mrm{R}Y\ar[l]_{\mu_Y} \ar[r]^-{\lambda_{Y}}& Y .}
$$
In the f{i}rst case, it suf{f}{i}ces to replace $\mrm{R}C$ by
$\mrm{R}C \stackrel{I_C}{\rightarrow} N \stackrel{I_C}{\leftarrow}
\mrm{R}C$, obtaining in this way $\rho'_{\hat{f}}$.\\
In the second case, replacing $\mrm{R}C'$ by $\mrm{R}C'
\stackrel{I_C'}{\rightarrow} N \stackrel{I_{C'}}{\leftarrow}
\mrm{R}C'$ we deduce that $\rho_{\hat{g}}$ is related to
$$\xymatrix@M=4pt@H=4pt@C=17pt@R=15pt{
 \rho''_{\hat{g}}: X & \mrm{R}X\ar[r]^{\mu_X} \ar[l]_-{\lambda_X} & \mrm{R}^2X\ar[r]^{\mrm{R}\mu_X} & \mrm{R}^{3}X \ar[r]^-{g'\comp\mrm{R}^2t} & N & \mrm{R}^{3}Y \ar[l]_{u'} & \mrm{R}^2Y\ar[l]_{\mrm{R}\mu_Y} & \mrm{R}Y\ar[l]_{\mu_Y} \ar[r]^-{\lambda_{Y}}& Y .}
$$
where $u'=I_{C'}\comp\mrm{R}J_T\comp\mrm{R}^2w$. But
$I_{C'}\comp\mrm{R}J_T$ is ``homotopic'' to
$I_{C}\comp\mrm{R}I_T$, that is, \ref{pdadesCylenD} provides
$H:cyl(\mrm{R}T)\rightarrow N$ such that composing with the
inclusions of $\mrm{R}^2T$ into $cyl(\mrm{R}T)$ we obtain just
$I_{C'}\comp\mrm{R}J_T$ and $I_{C}\comp\mrm{R}I_T$.\\
Then, a hammock relating $\rho''_{\hat{g}}$ to $\rho'_{\hat{g}}$
can be constructed in the usual way using $H$. So (A1) is already proved.\\[0.2cm]
\textsl{Axiom $($A2$)$:} We must show that $Ho\mc{D}$ has a zero
object, that is, an object $0$ such that
$\mrm{Hom}_{Ho\mc{D}}(0,0)=\{\ast\}=$trivial group. Since $\mc{D}$
has a zero object $0_{\mc{D}}$, that is at the same time an
initial and f{i}nal object, the from \ref{objetoF{i}nal} (or
\ref{objetoInicial}) (A2) follows.\\[0.2cm]
\textsl{Axiom $($A3$)$:} Given objects $X$, $Y$ in $Ho\mc{D}$, we
must show the existence of an object $Z$ and morphisms
$$\xymatrix@M=4pt@H=4pt@C=20pt@R=20pt{ X \ar@<-0.7ex>[r]_{i_1} & Z \ar@<0.7ex>[r]^{p_2} \ar@<-0.7ex>[l]_{p_1} & Y \ar@<0.7ex>[l]^{i_2}}$$
such that $p_1\comp i_1=Id_X$; $p_2\comp i_2=Id_Y$; $i_1\comp
p_1+i_2\comp p_2=Id_Z$ and
$p_2\comp i_1=p_1\comp i_2=0$ in $Ho\mc{D}$.\\
But, since $X$ and $Y$ are objects in $\mc{D}$, the data
$Z,i_1,i_2,p_1$ and $p_2$ so exists in $\mc{D}$, and it is enough
to
take the image under $\gamma:\mc{D}\rightarrow Ho\mc{D}$ of these morphisms.\\
Hence, since $\gamma$ is functorial it is clear that
$\gamma(p_1)\comp
\gamma(i_1)=Id_X$ and $\gamma(p_2)\gamma( i_2)=Id_Y$.\\
On the other hand, it follows from the def{i}nitions of sum in
$Ho\mc{D}$ and of $\gamma$ that $\gamma(f)+\gamma(g)=\gamma(f+g)$,
then $\gamma(i_1)\comp
\gamma(p_1)+\gamma(i_2)\comp\gamma(p_2)=Id_Z$.
To f{i}nish, the morphism $0$ in $\mc{D}$ is def{i}ned as the
unique morphism that factors through the zero object $0_{\mc{D}}$,
and since $\gamma(0_{\mc{D}})=0_{Ho\mc{D}}$, we deduce that
$\gamma$ maps the morphism $0$ into the morphism $0$, so the
equality $\gamma(p_2)\comp\gamma( i_1)=\gamma(p_1)\comp\gamma(
i_2)=0$ holds, and (A3) is proven.\\[0.2cm]
\indent In addition, as mentioned before,
$\gamma:\mrm{Hom}_{\mc{D}}(X,Y)\rightarrow\mrm{Hom}_{Ho\mc{D}}(X,Y)$
is lineal, so $\gamma$ is additive.\\[0.2cm]
\indent F{i}nally, let $F:\mc{D}\rightarrow\mc{D'}$ be a functor of additive simplicial descent categories.\\
Assume given $\hat{f},\hat{g}:X\rightarrow Y$ in $Ho\mc{D}$, and
zig-zags representing them as in \ref{denominadorComun}. By
def{i}nition, we have that $\hat{f}+\hat{g}$ is represented by
$$\xymatrix@M=4pt@H=4pt{\rho_{\hat{f}}+\rho_{\hat{g}}:&  X & \mrm{R}X \ar[l]_-{\lambda_X}\ar[r]^-{{f}+{g}} & L & \mrm{R}Y\ar[l]_w\ar[r]^-{\lambda_{Y}}& Y\ .}$$
Then, by \ref{Mof{i}smoHoD} it holds that
$\hat{f}+\hat{g}=\gamma(\lambda_Y)\comp(\gamma(w))^{-1}\comp\gamma(f+g)\comp(\gamma(\lambda_X))^{-1}$,
and since $\gamma$ and $F$ preserve sums, then
$HoF(\hat{f}+\hat{g})=HoF(\hat{f})+HoF(\hat{g})$ follows from the
equality $HoF\comp \gamma=\gamma\comp F$.
\end{proof}


\chapter{Relationship with triangulated categories}\label{CapitDescTriang}

 \setcounter{section}{1}
 \setcounter{subsection}{1}
\setcounter{thm}{0}

The aim of this chapter is to describe the ``left unstable''
triangulated structure existing on the homotopy category
$Ho\mc{D}$ associated with any simplicial descent category
$\mc{D}$.

The distinguished triangles will be those def{i}ned through the
cone functor $c:Maps(\mc{D})\rightarrow\mc{D}$.

With no extra assumption, neither additivity, this class of
triangles will satisfy all axioms of triangulated category but
TR3, that is the one involving the shift of distinguished
triangles. When $\mc{D}$ is an additive category then $Ho\mc{D}$
is ``right triangulated'' or ``suspended'' \cite{KV}, so if
in addition the shift functor is an equivalence of categories then $Ho\mc{D}$ is a triangulated category.\\
For triangulated categories we will follow the notations of \cite{GM}.\\

In order to simplify the notations we will write $f:X\rightarrow
Y$ to denote the morphism $\gamma(f)$ in $Ho\mc{D}$, for any
morphism $f$ in $\mc{D}$.

\begin{def{i}}[shift functor]\index{Index}{shift}\mbox{}\\
The \textit{shift functor}
$\mrm{T}:\mc{D}\rightarrow\mc{D}$\index{Symbols}{$\mrm{T}$} is
def{i}ned as
\begin{center}$\mrm{T}(X)=cyl(1\leftarrow X\rightarrow 1)=c(X\rightarrow 1)$.\end{center}
As usual, by $X[n]$ we mean $\mrm{T}^n(X)$.
\end{def{i}}
\begin{obs}\label{translacionLocalizado}
It follows from \ref{exactitudcyl} that $\mrm{T}$ preserve
equivalences, that is,
$\mrm{T}(\mrm{E})\subseteq \mrm{E}$ (in fact we will see that $\mrm{T}^{-1}(\mrm{E})= \mrm{E}$).\\
Then $\mrm{T}$ induces a functor between the localized categories,
that we will denote also by $\mrm{T}:Ho\mc{D}\rightarrow
Ho\mc{D}$. If $f$ is a morphism in $Ho\mc{D}$ we will also write
$f[1]$ instead of $\mrm{T}f$.\\
Note that our shift functors $\mrm{T}:\mc{D}\rightarrow\mc{D}$,
$\mrm{T}:Ho\mc{D}\rightarrow Ho\mc{D}$ may not be equivalences of
categories.
\end{obs}

\begin{num} If $X$ is an object of $\mc{D}$, there exists an isomorphism
$\theta_X:\mrm{R}(X[1])\rightarrow (\mrm{R}X)[1]$ in $Ho\mc{D}$
functorial in $X$, given by
$$\xymatrix@M=4pt@H=4pt@C=30pt{\mrm{R}(X[1])\ar[r]^{\lambda_{X[1]}} & X[1] & (\mrm{R}X)[1]\ar[l]_{\lambda_X[1]} \ .}$$
\end{num}

The next proposition and its corollary will not be used in this
work, we introduce them just for completeness.

\begin{prop}
If $f$ is a morphism in $\mc{D}$, there exists an isomorphism in
$Ho\mc{D}$
$$\theta_f:c(f[1])\rightarrow(c(f))[1]\ ,$$
functorial on $f$ and such that the following diagram commutes
$($in $Ho\mc{D})$
$$\xymatrix@M=4pt@H=4pt@R=20pt{
 \mrm{R}(Y[1]) \ar[d]_{\theta_Y} \ar[r]^-{I_{Y[1]}} & c(f[1])\ar[d]_{\theta_f}\\
 (\mrm{R}Y)[1] \ar[r]^-{I_Y[1]}                     &(c(f))[1]\ .}$$
\end{prop}

\begin{proof}%
Let $f:X\rightarrow Y$ be a morphism in $\mc{D}$. To see the
existence of the isomorphism $\theta_f$ it suf{f}{i}ces to apply
the factorization property of the cone,
 \ref{coherenciaCono}, to the square
$$\xymatrix@M=4pt@H=4pt{X \ar[r]^f  \ar[d]_{Id} &  Y \ar[d]\\
                        X \ar[r]  &  1   }$$
and substitute $c(Id_1)$ by $1$ when needed.\\
Note that $\theta_Y$ is just $\theta_{0\rightarrow Y}$, and by
def{i}nition $\theta_f$ is functorial on $f$. Hence the equality
$\theta_f\comp I_{Y[1]}=I_Y[1]\comp \theta_Y$ holds.
\end{proof}

\begin{cor}
A morphism $f$ of $\mc{D}$ is in $\mrm{E}$ if and only if
$f[1]$ is so.\\
Therefore, a morphism $f$ of $Ho\mc{D}$ is an isomorphism if and
only $f[1]$ is so.
\end{cor}

\begin{proof} Consider a morphism $f:X\rightarrow Y$ of
$\mc{D}$ such that $f[1]\in\mrm{E}$. From the acyclicity axiom (or
its corollary \ref{axiomaconoenD}) we deduce that $c(f[1])\simeq
(c(f))[1]\rightarrow 1$ is in $\mrm{E}$. But
$(c(f))[1]=c(c(f)\rightarrow 1)$, and again from (SDC 7) follows
that
$c(f)\rightarrow 1$ is an equivalence, hence $f\in\mrm{E}$.\\
The last statement is a formal consequence of the f{i}rst one
since $\mrm{E}$ is saturated.
\end{proof}

\begin{def{i}}[triangles in $Ho\mc{D}$]\index{Index}{triangles!}\mbox{}\\
A triangle in $Ho\mc{D}$ is a sequence of morphisms in $Ho\mc{D}$
of the form
\begin{center}$X\rightarrow Y\rightarrow Z\rightarrow X[1]$.\end{center}
A morphism between the triangles $X\rightarrow Y\rightarrow
Z\rightarrow X[1]$ and $X'\rightarrow Y'\rightarrow Z'\rightarrow
X'[1]$ is a commutative diagram in $Ho\mc{D}$
$$\xymatrix@M=4pt@H=4pt@R=20pt{X \ar[r] \ar[d]_{\alpha} & Y \ar[r]\ar[d]_\beta & Z \ar[r]\ar[d]_\delta & X[1]\ar[d]_{\alpha[1]}\\
                                  X' \ar[r]  & Y' \ar[r]& Z' \ar[r]& X'[1].}$$
\end{def{i}}

\begin{num}\label{construccionTriangDist}
Given a morphism $f:X\rightarrow Y$ in $\mc{D}$, if we apply $cyl$
to the diagram
$$\xymatrix@M=4pt@H=4pt{ 1 \ar[d] & X \ar[r]^f \ar[l] \ar[d]_{Id} &  Y \ar[d]\\
                         1 & X \ar[r] \ar[l] &  1   }$$
we will obtain, by \ref{cubConmenD}, the diagram
\begin{equation}\label{CuboConmTriang}\xymatrix@R=33pt@C=33pt@H=4pt@M=4pt@!0{                                                                &  \mrm{R}X\ar[rr]^{\mrm{R}f}\ar'[d][dd]\ar[ld]&                            & \mrm{R}Y \ar[dd]\ar[ld]_{I_Y}\\
                                         \mrm{R}1 \ar[rr]^{\;\;\;\;\;\;I_1}\ar[dd] &                                                                                       & c(f)  \ar[dd]^(0.4){p}\\
                                                                                                         &  \mrm{R}X \ar'[r][rr]\ar[ld]                                &                            & \mrm{R}1\, ,\ar[ld]^{J_1}\\
                                             \mrm{R}1  \ar[rr]^{K_{1}}                                     &                                                                                       &  X[1]}\end{equation}
where all faces commute in $\mc{D}$ except the top and bottom ones, that commutes in $Ho\mc{D}$.\\
Therefore, $f$ gives rise to the following sequence of morphisms
of $\mc{D}$
$$\xymatrix@M=4pt@H=4pt{\mrm{R}X \ar[r]^{\mrm{R}f} &  \mrm{R}Y
\ar[r]^{I_Y} & c(f) \ar[r]^{p_f} &  X[1].}$$
It holds that the compositions $P\comp I_Y$ and $I_Y\comp f$ are
trivial in $Ho\mc{D}$, that is, they factor through the f{i}nal
object $1$ (since $1\simeq \mrm{R}1$ in $Ho\mc{D}$).
\end{num}

\begin{def{i}}[distinguished triangles in $Ho\mc{D}$]\label{}\index{Index}{triangles!distinguished}\mbox{}\\
Def{i}ne \textit{distinguished triangles} in $Ho\mc{D}$ as those
triangles isomorphic, for some morphism $f$ of $\mc{D}$, to
\begin{equation}\label{def{i}triangDisting}
\xymatrix@M=4pt@H=4pt{
 X\ar[r]^f &  Y \ar[r]^{\iota_Y} & c(f) \ar[r]^{p} & X[1]}
\end{equation}
where $\iota_Y$ is the composition in $Ho\mc{D}$ of
$\xymatrix@M=4pt@H=4pt{Y \ar[r]^{\lambda_Y^{-1}}&\mrm{R}Y
\ar[r]^{I_Y} & c(f)}$.
\end{def{i}}
\begin{obs}
Therefore, we have automatically that the composition of two
consecutive maps in a distinguished triangle is trivial, that is,
it factors through the object 1 in $Ho\mc{D}$.
\end{obs}

Next we will see that most of the axioms of triangulated
categories hold for this class of distinguished triangles.

\begin{prop}[\textbf{TR 1}]\label{TR1}\mbox{}\\
i) the triangle $X\stackrel{Id}{\rightarrow}X\rightarrow
1\rightarrow X[1]$ is distinguished, where the map $1\rightarrow
X[1]$ is the composition
$\xymatrix@M=4pt@H=4pt{1\ar[r]^{\lambda_1^{-1}} & \mrm{R}1
\ar[r]^{{K_1}} & X[1] }$ $($see $($\ref{CuboConmTriang}$))$.\\
ii) Every triangle isomorphic to a distinguished triangle is also distinguished.\\
iii) Given $f:X\rightarrow Y$ in $Ho\mc{D}$, there exists a
distinguished triangle of the form
\begin{center}$\xymatrix@M=4pt@H=4pt{X\ar[r]^f
&Y\ar[r]&Z\ar[r]&X[1]}$.\end{center}
\end{prop}

\begin{proof}\mbox{}\\
To see \textit{i}), consider (\ref{CuboConmTriang}) for $f=Id_X$.
Set $\rho=K_1\comp \lambda_1^{-1}: 1\rightarrow X[1]$. The map
$I_{1}:\mrm{R}1\rightarrow c(Id)$ is in $\mrm{E}$, so
$c(Id)\stackrel{\eta}{\rightarrow} 1$ is so (since $\eta\comp
I_1=\lambda_1\in\mrm{E}$). The diagram
$$\xymatrix@M=4pt@H=4pt@R=20pt{X \ar[r] \ar[d]_{Id} & X \ar[r]\ar[d]_{Id}& c(Id) \ar[r]^{p_{Id}}\ar[d]_{\eta} & X[1]\ar[d]_{Id}\\
                                  X \ar[r]  & X \ar[r]& 1 \ar[r]^{\rho}& X[1].}$$
provides a morphism of triangle. Indeed, to see the commutativity
in $Ho\mc{D}$ of the right square just note that by
(\ref{CuboConmTriang}), $K_1=p\comp I_1$, and then $\rho\comp
\eta=K_1\comp (\lambda_1^{-1}\comp\eta)=K_1\comp I_1^{-1}=p$.\\
\textit{ii}) holds by def{i}nition of distinguished triangle.\\
Then it remains to prove \textit{iii}). Given a morphism
$f:X\rightarrow Y$ in $Ho\mc{D}$, from \ref{ProHoDcat} we have
that $f$ is represented by a zig-zag of the form
$\xymatrix@M=4pt@H=4pt{X & \mrm{R}X
\ar[l]_{\lambda_X}\ar[r]^{\overline{f}} &  T & \mrm{R}Y
\ar[l]_w\ar[r]^{\lambda_{Y}}& Y }$. If $Z=c(\overline{f})$, we
consider the distinguished triangle
$$\xymatrix@M=4pt@H=4pt{\mrm{R}X\ar[r]^{\overline{f}} &  T \ar[r]^{\iota_{T}} & Z\ar[r]^-{p} & (\mrm{R}X)[1]}\ .$$
Let $g:Y\rightarrow Z$ and $h:Z\rightarrow X[1]$ be the
compositions given respectively by
$$\xymatrix@M=4pt@H=4pt@R=10pt{
 Y \ar[r]^{\lambda^{-1}_{Y}} & \mrm{R}Y \ar[r]^w & T \ar[r]^{\iota_{T}} & c(\overline{f})=Z \\
 Z \ar[r]^-{p} & (\mrm{R}X)[1] \ar[r]^-{\lambda_X[1]} & X[1] \ .}$$
Then, setting $\alpha=\lambda_Y\comp w^{-1}:T\rightarrow Y$ we
deduce the following commutative diagram
$$\xymatrix@M=4pt@H=4pt@R=20pt{\mrm{R}X\ar[r]^{\overline{f}}\ar[d]_{\lambda_X} &  T\ar[d]_{\alpha} \ar[r]^{\iota_{T}} & Z\ar[d]_{Id} \ar[r]^-{p} &(\mrm{R}X)[1]\ar[d]^{(\lambda_X)[1]}\\
                               X \ar[r]^f                                      & Y\ar[r]^g                       & Z\ar[r]^h                & X[1]}$$
that is in fact an isomorphism of triangles, so the bottom
triangle is distinguished.
\end{proof}
\begin{prop}[\textbf{TR 3}]\label{AxiomaTR3}\mbox{}\\
If $X\rightarrow Y\rightarrow Z\rightarrow X[1]$ and
$X'\rightarrow Y'\rightarrow Z'\rightarrow X'[1]$ are
distinguished triangles and
$$\xymatrix@M=4pt@H=4pt{  X \ar[r]^-{f} \ar[d]_{\alpha} &  Y \ar[d]^{\beta}\\
                                                                              X' \ar[r]^-{g}        &  Y'   ,}$$
commutes in $Ho\mc{D}$, then there exists an isomorphism of
triangles
$$\xymatrix@M=4pt@H=4pt{X\ar[r]^f \ar[d]_{\alpha} & Y \ar[r]\ar[d]_{\beta}    & Z \ar[r]\ar[d]_h & X[1]\ar[d]_{\alpha[1]}\\
                         X'\ar[r]^{g}                 & Y'\ar[r]   & Z' \ar[r]& X'[1].}$$
In addition, if $\alpha$ and $\beta$ are isomorphisms of
$Ho\mc{D}$ then so is $h$.
\end{prop}

\begin{proof}\mbox{}\\
By def{i}nition of distinguished triangle we can assume that $f$
and $g$ are morphisms of $\mc{D}$ and the triangles $X\rightarrow
Y\rightarrow Z\rightarrow X[1]$, $X'\rightarrow Y'\rightarrow
Z'\rightarrow X'[1]$ are those obtained from $f$ and $g$
respectively as in (\ref{def{i}triangDisting}).\\[0.2cm]
\textit{Case 1}: $\alpha$ and $\beta$ are morphisms in $\mc{D}$
and
$\beta\comp f = g\comp\alpha$ in $\mc{D}$.\\
In this case it follows directly from the functoriality of the
cone the existence of $h: c(f)=Z\rightarrow c(g)=Z'$ in $\mc{D}$
such that the required diagram is commutative. If $\alpha,\beta$
are isomorphisms in $Ho\mc{D}$ then they are in $\mrm{E}$, since
$\mrm{E}$ is saturated. Hence, we deduce from corollary
$\ref{exactitudcono}$ that $h\in\mrm{E}$.\\[0.2cm]
\textit{Case 2}: $\alpha$ and $\beta$ are morphisms of $\mc{D}$
and
$\beta f=g\alpha$ in $Ho\mc{D}$.\\
In this case the zig-zags
$$\xymatrix@M=4pt@H=4pt{
 X & \mrm{R}X  \ar[r]^{ \mrm{R}(\beta\comp f)}\ar[l]_{\lambda_X} &  \mrm{R}Y' &   \mrm{R}Y'\ar[l]_{Id} \ar[r]^{\lambda_{Y'}}& Y'\\
 X &  \mrm{R}X \ar[r]^{ \mrm{R}(g\comp\alpha)}\ar[l]_{\lambda_X}  &  \mrm{R}Y' &  \mrm{R}Y'\ar[l]_{Id} \ar[r]^{\lambda_{Y'}}& Y'}$$
def{i}ne the same morphism of $Ho\mc{D}$, and by \ref{ProHoDcat}
we have a hammock in $\mc{D}$ in the form
\begin{equation}\label{diagramaAuxiliarC4}
\xymatrix@M=4pt@H=4pt@C=40pt{
               & \mrm{R}^{2}X\ar[ld]_{Id}\ar[r]^{\mrm{R}^2(\beta\comp f)}\ar[d]_{w}   & \mrm{R}^2Y'\ar[d]_{t}    & \mrm{R}^{2}Y'\ar[l]_{Id}\ar[rd]^{Id}\ar[d]  &\\
 \mrm{R}^{2}X  & \overline{X} \ar[l]\ar[r]^H                                                & L                        & \overline{Y}     \ar[l]\ar[r]               &\mrm{R}^{2}Y' \\
               & \mrm{R}^{2}X\ar[lu]^{Id}\ar[r]^{\mrm{R}^2(g\comp\alpha)}\ar[u]^{w'} & \mrm{R}^2Y'\ar[u]^{t'}   & \mrm{R}^{2}Y',\ar[l]_{Id}\ar[ru]_{Id}\ar[u] &}
\end{equation}
where all maps are in $\mrm{E}$ except $\mrm{R}^2(\beta\comp f)$,
$\mrm{R}^2(g\comp\alpha)$ and $H$.
Hence, if we denote by $\lambda^2:\mrm{R}^2\rightarrow \mrm{R}$
the natural transformation with
$\lambda_S^2=\lambda_S\comp\lambda_{\mrm{R}S}$, we have the
following diagram consisting of the commutative squares in
$\mc{D}$
$$
\xymatrix@M=4pt@H=4pt@C=30pt{
 X \ar[d]^f & \mrm{R}^{2}X \ar[d]^{\mrm{R}^{2}f}\ar[l]_{\lambda^2_X}\ar[r]^{Id} & \mrm{R}^{2}X \ar[d]^{\mrm{R}^2 (\beta\comp f)}\ar[r]^{w} & T\ar[d]^H   & \mrm{R}^{2}X\ar[l]_{w'}\ar[d]^{\mrm{R}^2(g\comp\alpha)}\ar[r]^{\mrm{R}^2\alpha} & \mrm{R}^{2}X'\ar[d]^{\mrm{R}^2 g}\ar[r]^{\lambda^2_{X'}} & X'\ar[d]^{g}\\
 Y          & \mrm{R}^{2}Y \ar[r]^{\mrm{R}^{2}\beta}\ar[l]_{\lambda^2_Y}        & \mrm{R}^2Y'\ar[r]^{t}                                    & S           & \mrm{R}^2Y'\ar[l]_{t'}\ar[r]^{Id}                                               & \mrm{R}^{2}Y'\ar[r]^{\lambda^2_{Y'}}               & Y' \ .}
$$
By the f{i}rst case there exists morphisms $\widetilde{\lambda}$,
$\widetilde{\beta}$, $u$, $u'$, $\widetilde{\alpha}$ and
$\widehat{\lambda}$ such that the diagram
$$\xymatrix@M=4pt@H=4pt@C=30pt{  X \ar[r]^f                                                                       & Y   \ar[r]                                                        & Z=c(f) \ar[r]                                                               & X[1]\\
                                \mrm{R}^{2}X \ar[r]^{\mrm{R}^{2}f}\ar[u]^{\lambda^2_X}\ar[d]_{Id}                 & \mrm{R}^{2}Y \ar[d]_{\mrm{R}^{2}\beta}\ar[r]\ar[u]^{\lambda^2_Y}  & c(\mrm{R}^{2}f)\ar[u]^{\widetilde{\lambda}}\ar[d]_{\widetilde{\beta}}\ar[r] & (\mrm{R}^{2}X)[1]\ar[d]_{Id}\ar[u]^{\lambda^2_X[1]} \\
                                \mrm{R}^{2}X \ar[r]^{\mrm{R}^2 (\beta\comp f)}\ar[d]_{w}                          & \mrm{R}^2Y'\ar[d]_{t} \ar[r]                                      & c(\mrm{R}^{2}(\beta f))\ar[d]_{u}\ar[r]                                     & (\mrm{R}^{2}X)[1]\ar[d]_{w[1]}\\
                                T\ar[r]^H                                                                         & S\ar[r]                                                           & c(H)\ar[r]                                                                  & T[1]\\
                                \mrm{R}^{2}X\ar[u]^{w'}\ar[r]^{\mrm{R}^2 (g\comp\alpha)}\ar[d]_{\mrm{R}^2\alpha}  & \mrm{R}^2Y'\ar[u]^{t'}\ar[r]\ar[d]_{Id}                           & c(\mrm{R}^2 (g\alpha))\ar[r] \ar[d]_{\widetilde{\alpha}}\ar[u]^{u'}         & (\mrm{R}^{2}X)[1]\ar[u]^{w'[1]}\ar[d]_{(\mrm{R}^{2}\alpha)[1]}\\
                                \mrm{R}^{2}X'\ar[r]^{\mrm{R}^2 g}\ar[d]_{\lambda^2_{X'}}                          & \mrm{R}^{2}Y'\ar[r]\ar[d]_{\lambda^2_{Y'}}                        & c(\mrm{R}^2 g)\ar[r]\ar[d]_{\widehat{\lambda}}                              & (\mrm{R}^{2}X')[1]\ar[d]_{\lambda^2_{X'}[1]}  \\
                                X'\ar[r]^{g}                                                                      & Y'\ar[r]                                                          & Z'=c(g) \ar[r] &X'[1]}
$$
commutes in $Ho\mc{D}$.
On the other hand, the morphisms $\widetilde{\lambda}$, $u$, $u'$
and $\widehat{\lambda}$ are in $\mrm{E}$. Observe that the
composition in $Ho\mc{D}$
of the morphisms in the f{i}rst column is just $\alpha$.\\
Indeed, from (\ref{diagramaAuxiliarC4}) it follows that
$w'=w^{-1}$, and it is enough to have into account the equality
$\lambda^2_{X'}\comp\mrm{R}^2\alpha=\alpha\comp\lambda^2_{X}$,
that holds since $\lambda^2$ is a natural transformation. In the
same way, the second column
is the morphism $\beta$, whereas the fourth one is $\alpha[1]$.\\
Summing all up, we get a morphism
$h=\widehat{\lambda}\comp\widetilde{\alpha}\comp (u')^{-1}\comp
u\comp \widetilde{\beta}\comp\widetilde{\lambda}^{-1}$ such that the requested diagram commutes.\\[0.2cm]
F{i}nally, if $\alpha$ and $\beta$ are in $\mrm{E}$, then
$\mrm{R}^2\alpha$ and $\mrm{R}^2\beta$ are also equivalences, and
by the previous case the same holds for $\widetilde{\alpha}$ and
$\widetilde{\beta}$. Therefore $h$ an isomorphism in $Ho\mc{D}$.\\[0.2cm]
\textit{Case 3}: General case: $\alpha$ and $\beta$ are morphism
in
$Ho\mc{D}$.\\
Let $A$ and $B$ be zig-zags representing $\alpha$ and $\beta$
respectively, given by
$$\xymatrix@M=4pt@H=4pt@R=10pt{
 X & \mrm{R}X \ar[r]^{{\alpha}'} \ar[l]_{\lambda_X}   & S & \mrm{R}X' \ar[l]_{u} \ar[r]^{\lambda_{X'}}& X'\\
 Y & \mrm{R}Y \ar[r]^{{\beta}'}  \ar[l]_{\lambda_{Y}} & T & \mrm{R}Y' \ar[l]_{v} \ar[r]^{\lambda_{Y'}}& Y'\ .}
 $$
Consider the diagram
$$\xymatrix@M=4pt@H=4pt{
 X\ar[d]_{f} &  \mrm{R}X \ar[r]^{{\alpha}'}\ar[d]_{\mrm{R}f}\ar[l]_{\lambda_X}  & S & \mrm{R}X' \ar[d]_{\mrm{R}g} \ar[l]^{u} \ar[r]^{\lambda_{X'}}& X'\ar[d]_{g} \\
 Y           &  \mrm{R}Y \ar[r]^{{\beta}'}                  \ar[l]_{\lambda_{Y}}& T & \mrm{R}Y'                   \ar[l]^{v} \ar[r]^{\lambda_{Y'}}& Y' \ .}
$$
If there exists $t:S\rightarrow T$ such that
$t\comp\alpha'=\beta'\comp\mrm{R}f$ and $t\comp u=v\comp\mrm{R}g$
in $Ho\mc{D}$, then it suf{f}{i}ces to apply the case 2 to the
squares in the above diagram.
Let us check that we can always choose zig-zags representing
$\alpha$ and $\beta$ satisfying this property. That is to say, it
is enough see that there exists zig-zags $A'$ and $B'$
$$\xymatrix@M=4pt@H=4pt@R=10pt{
 X & \mrm{R}X \ar[r]^{{\alpha}''} \ar[l]_{\lambda_X}   & S' & \mrm{R}X' \ar[l]_{u'} \ar[r]^{\lambda_{X'}}& X'\\
 Y & \mrm{R}Y \ar[r]^{{\beta}''}  \ar[l]_{\lambda_{Y}} & T' & \mrm{R}Y' \ar[l]_{v'} \ar[r]^{\lambda_{Y'}}& Y'}
 $$
representing $\alpha$ and $\beta$ and such that there exists
$s:S\rightarrow T$ that makes the following diagram commute in
$Ho\mc{D}$
$$\xymatrix@M=4pt@H=4pt{
 X\ar[d]_{f} &  \mrm{R}X \ar[r]^{{\alpha}''}\ar[d]_{\mrm{R}f}\ar[l]_{\lambda_X}  & S'\ar[d]_s & \mrm{R}X' \ar[d]_{\mrm{R}g} \ar[l]_{u'} \ar[r]^{\lambda_{X'}}& X'\ar[d]_{g} \\
 Y           &  \mrm{R}Y \ar[r]^{{\beta}''}                  \ar[l]_{\lambda_{Y}}& T'         & \mrm{R}Y'                   \ar[l]_{v'} \ar[r]^{\lambda_{Y'}}& Y' \ .}
$$
By \ref{FrelRF}, the zig-zags $\mrm{R}^2A$ and $\mrm{R}^2B$ given
by
$$\xymatrix@M=4pt@H=4pt@R=10pt{
 X & \mrm{R}X\ar[r]^{{\mu}_X}\ar[l]_{\lambda_X} & \mrm{R}^2 X  \ar[r]^{\mrm{R}\alpha'}  & \mrm{R}S  & \mrm{R}^2 X' \ar[l]_{\mrm{R}u} & \mrm{R}X'\ar[l]_{{\mu}_{X'}}\ar[r]^{\lambda_{X'}} & X'\\
 Y & \mrm{R}Y\ar[r]^{{\mu}_Y}\ar[l]_{\lambda_Y} & \mrm{R}^2 Y  \ar[r]^{\mrm{R}{\beta}'} & \mrm{R}T  & \mrm{R}^2 Y' \ar[l]_{\mrm{R}v} & \mrm{R}Y'\ar[l]_{{\mu}_{Y'}}\ar[r]^{\lambda_{Y'}} & Y' }
$$
represent also the morphisms $\alpha$ and $\beta$ respectively.\\
imply that the square
$$\xymatrix@M=4pt@H=4pt@C=20pt@R=20pt{
  \mrm{R}X'\ar[d]_{\mrm{R}u}\ar[r]^{\mrm{R}g} & \mrm{R}Y'\ar[d]^{I_{Y'}}\\
  \mrm{R}S  \ar[r]^-{I_S}               & cyl(g,u)}$$
commutes in $Ho\mc{D}$, and $I_{Y'}$ is an equivalence. Moreover,
in the same way as before we can build the square
$$\xymatrix@M=4pt@H=4pt@C=30pt@R=23pt{
  \mrm{R}^2Y'\ar[d]_{\mrm{R}I_{Y'}}\ar[r]^{\mrm{R}v} & \mrm{R}T\ar[d]_{I_{T}}\\
  \mrm{R}cyl(g,u)  \ar[r]^-{I_{cyl(g,u)}}             & cyl(\mrm{R}v,I_{Y'})}$$
that commutes in $Ho\mc{D}$ and such that all maps are in
$\mrm{E}$. Set $T'=cyl(\mrm{R}v,I_{Y'})$. Since $I_T:
\mrm{R}T\rightarrow T'\in\mrm{E}$, it is clear that $\mrm{R}^2A$
is related to
$$\xymatrix@M=4pt@H=4pt@C=20pt{
 Y & \mrm{R}Y\ar[r]^{{\mu}_Y}\ar[l]_{\lambda_Y} & \mrm{R}^2 Y  \ar[r]^{\mrm{R}{\beta}'} & \mrm{R}T\ar[r]^{I_T} & T'& \mrm{R}T \ar[l]_{I_T} & \mrm{R}^2 Y' \ar[l]_{\mrm{R}v} & \mrm{R}Y'\ar[l]_{{\mu}_{Y'}}\ar[r]^{\lambda_{Y'}} & Y' \ .}
$$
Consequently it suf{f}{i}ces to check that the morphism
$s=I_{cyl(g,u)}\comp \mrm{R}I_S:\mrm{R}T\rightarrow S'$ is such
that
$$\xymatrix@M=4pt@H=4pt{
 X \ar[d]_{f} & \mrm{R}X\ar[r]^{{\mu}_X}\ar[l]_{\lambda_X}\ar[d]_{\mrm{R}f} & \mrm{R}^2 X \ar[d]_{\mrm{R}^2f} \ar[rr]^{\mrm{R}\alpha'} & \ar@{}[d]|{\mathrm{I}}        &\mrm{R}S\ar[d]_s  & \ar@{}[d]|{\mathrm{II}}        & \mrm{R}^2 X' \ar[ll]_{\mrm{R}u}\ar[d]_{\mrm{R}^2g} & \mrm{R}X'\ar[d]_{\mrm{R}f}\ar[l]_{{\mu}_{X'}}\ar[r]^{\lambda_{X'}} & X'\ar[d]_{g}\\
 Y            & \mrm{R}Y\ar[r]^{{\mu}_Y}\ar[l]_{\lambda_Y}                  & \mrm{R}^2 Y \ar[r]^{\mrm{R}{\beta}'}                     & \mrm{R}T\ar[r]^{I_T} & T'               & \mrm{R}T \ar[l]_{I_T} & \mrm{R}^2 Y' \ar[l]_{\mrm{R}v}                     & \mrm{R}Y'\ar[l]_{{\mu}_{Y'}}\ar[r]^{\lambda_{Y'}}                  & Y' }
$$
commutes in $Ho\mc{D}$. In order to see that, it is clear that the square II commutes.\\
To see I, since $\beta\comp f=g\comp\alpha$ and $\lambda_W$ is an
isomorphism in $Ho\mc{D}$ for any $W$ in $\mc{D}$, we deduce the
commutativity in $Ho\mc{D}$ of the diagram
$$\xymatrix@M=4pt@H=4pt{
  \mrm{R}X \ar[r]^{{\alpha}'}\ar[d]_{\mrm{R}f}& S  \ar[r]^{u^{-1}} & \mrm{R}X' \ar[d]_{\mrm{R}g} \\
  \mrm{R}Y \ar[r]^{{\beta}'}                  & T  \ar[r]^{v^{-1}} & \mrm{R}Y'                   \ .}
$$
Then we have that
$$I_T\comp\mrm{R}\beta'\comp\mrm{R}^2f= I_T\comp\mrm{R}v\comp(\mrm{R}v^{-1}\mrm{R}\beta'\comp\mrm{R}^2f)
=(I_T\comp\mrm{R}v)\comp\mrm{R}^2g\comp\mrm{R}u^{-1}\comp\mrm{R}\alpha'=$$
$$I_{cyl(g,u)}\comp(\mrm{R}I_{Y'}\comp\mrm{R}^2g)\comp\mrm{R}u^{-1}\comp\mrm{R}\alpha'=
 I_{cyl(g,u)}\comp\mrm{R}I_S\mrm{R}^2u\comp\mrm{R}u^{-1}\comp\mrm{R}\alpha'=s\comp\mrm{R}\alpha'\ .$$
\end{proof}

Now we will begin the proof of the octahedron axiom.

\begin{num}\label{def{i}Octaedro}
Two composable morphisms
$X\stackrel{u}{\rightarrow}Y\stackrel{v}{\rightarrow}Z$ in
$\mc{D}$ gives rise in a natural way to the triangle
$$c(u)\stackrel{\alpha}{\rightarrow} c(v\comp u)\stackrel{\beta}{\rightarrow} c(v)\stackrel{\gamma}{\rightarrow} c(u)[1]\ .$$
Indeed, applying the cone functor to the following squares
$$\xymatrix@M=4pt@H=4pt{ X \ar[d]_{Id}\ar[r]^u  \ar@{}[rd]|{A}        & Y \ar[d]^{v}     && X\ar[r]^{v\comp u} \ar[d]_{u} & Z\ar[d]^{Id} \\
                         X \ar[r]^{v\comp u}                          & Z\ar@{}[rru]|{;} && Y\ar[r]^{v}    \ar@{}[ru]|{B} & Z           .}$$
we obtain $c(u)\stackrel{\alpha}{\rightarrow} c(vu)$ and
$c(vu)\stackrel{\beta}{\rightarrow} c(v)$ respectively.\\
On the other hand $c(v)\stackrel{\gamma}{\rightarrow}c(u)[1]$ is
def{i}ned as the composition $c(v)\stackrel{p}{\rightarrow} Y[1]
\stackrel{\iota_Y[1]}{\rightarrow} c(u)[1]$.
\end{num}

\begin{prop}\label{ConoDeComposicion}
Under the notations given above, the triangle
$c(u)\stackrel{\alpha}{\rightarrow} c(v\comp
u)\stackrel{\beta}{\rightarrow} c(v)\stackrel{\gamma}{\rightarrow}
c(u)[1]$ is distinguished in $Ho\mc{D}$.
\end{prop}
\begin{proof}
Let us see that the above triangle is isomorphic to the one
induced by $\alpha$, that is
$$\xymatrix@M=4pt@H=4pt@C=30pt{
 c(u)\ar[r]^{\alpha} &  c(v\comp u) \ar[r]^{\iota} & c(\alpha) \ar[r]^{p'} & c(u)[1]\ .}
$$
We will apply the factorization property of the cone
\ref{coherenciaCono} to the square
$$\xymatrix@M=4pt@H=4pt{
 X \ar[d]_{u}\ar[r]^{Id}  & X \ar[d]^{v\comp u}\\
 Y \ar[r]^{v}             & Z \ .}$$
To that end we introduce some notations.\\
Let $\widehat{u}:c(Id_X)\rightarrow c(v)$ be the morphism obtained
by applying the cone functor by rows to the previous square, as
well as $\psi':c(\mrm{R}(v\comp u))\rightarrow c(\widehat{u})$,
$\psi:c(\mrm{R}v)\rightarrow c(\alpha)$,
$\widehat{\lambda}:c(\mrm{R}(v\comp u))\rightarrow c(v\comp u)$
and $\widetilde{\lambda}:c(\mrm{R}v)\rightarrow c(v)$ the
morphisms obtained in the same way from the squares
$$\xymatrix@M=4pt@H=4pt{\mrm{R}X\ar[r]^{\mrm{R}(v\comp u)}\ar[d]_{I} & \mrm{R}Z\ar[d]_{I}\ar@{}[rd]|{;} & \mrm{R}Y\ar[d]_{I}\ar[r]^{\mrm{R}v} & \mrm{R}Z\ar[d]_{I}\ar@{}[rd]|{;} & \mrm{R}X\ar[r]^{\mrm{R}(v\comp u)}\ar[d]_{\lambda_{X}} & \mrm{R}Z\ar[d]_{\lambda_{Z}}\ar@{}[rd]|{;}& \mrm{R}Y\ar[r]^{\mrm{R}v}\ar[d]_{\lambda_{Y}} & \mrm{R}Z\ar[d]_{\lambda_{Z}}\\
                                   c(Id_X)          \ar[r]^{\widehat{u}}   & c(v)                             & c(u) \ar[r]^{\alpha}                      & c(v\comp u)                             & X\ar[r]^{v\comp u}                                    & Z                                          & Y\ar[r]^{v}                                    & Z\ ,}$$
where each $I$ denotes the corresponding canonical inclusion.\\
Denote by $\widetilde{T}\in\simp\mc{D}$ the image under
$\widetilde{Cyl}$ of $1\times\Dl\leftarrow
C(u)\stackrel{\widetilde{\alpha}}{\rightarrow} C(v\comp u)$. Take
isomorphisms $\Phi:\mbf{s}(\widetilde{T})\rightarrow c(\alpha)$
and $\Psi:\mbf{s}(\widetilde{T})\rightarrow c(\widehat{u})$ such
that the diagram
\begin{equation}\label{diagrama2AuxiliarC4}
\xymatrix@R=50pt@C=50pt@H=4pt@M=4pt@!0{                                                  & \mrm{R}^2 Z \ar'[d][dd]_{I}\ar[ld]_{\mrm{R}I}\ar[rr]^{\lambda_{\mrm{R}Z}}&                                                   &  \mrm{R}Z\ar'[d][dd]_{I}\ar[ld]_{I}  &                                                                              & \mrm{R}^2 Z \ar[dd]^{\mrm{R}I}\ar[ld]_{I} \ar[ll]_{\lambda_{\mrm{R}Z}}  \\
                                             \mrm{R}c(v)\ar[dd]_{I}\ar[rr]^(0.6){\lambda}&                                                                          & c(v)\ar[dd]_(0.7){\eta}                           &                                      &  c(\mrm{R}v) \ar[dd]\ar[ll]_(0.4){\widetilde{\lambda}}\ar[dd]^(0.6){\psi}  &                                               \\
                                                                                         & c(\mrm{R}(v\comp u))\ar[ld]_{\psi'}\ar'[r][rr]^(0.2){\widehat{\lambda}}  &                                                   & c(v\comp u)\ar[ld]_{\eta'}           &                                                                              & \mrm{R}c(v\comp u)\, ,\ar[ld]^{I}\ar'[l][ll]_(0.2){\lambda}\\
                                             c(\widehat{u})                              &                                                                          & \mbf{s}\widetilde{T}  \ar[rr]^{\Phi}\ar[ll]_{\Psi}&                                      &  c(\alpha)                                                                   &  }
\end{equation}
commutes, where $\eta$ is the image under $\mbf{s}$ of the
canonical inclusion of $C(v)$ into
$\widetilde{T}$, whereas $\eta'$ the image under $\mbf{s}$ of the morphism induced by the canonical inclusions of $X$ and $Z$ into $C(Id_X)$ and $C(v)$ respectively.\\
Since $c(Id_X)\rightarrow 1$ is an equivalence, by
\ref{axiomacylenD} we have that $I:\mrm{R}c(v)\rightarrow
c(\widehat{u})$ is in $\mrm{E}$. Hence, we deduce from the
commutativity of the front face of the above diagram that
$\eta,\psi\in\mrm{E}$.\\
Set $\tau=\psi\comp (\widetilde{\lambda})^{-1}=\Phi\comp
\eta:c(v)\rightarrow c(\alpha)$. It is enough to see that the
diagram
$$\xymatrix@M=4pt@H=4pt@C=30pt{
 c(u)\ar[r]^{\alpha}\ar[d]_{Id} &  c(v\comp u) \ar[r]^{\beta}\ar[d]_{Id}      & c(v) \ar[r]^{\gamma}\ar[d]_{\tau}      & c(u)[1]\ar[d]_{Id}\\
 c(u)\ar[r]^{\alpha}            &  c(v\comp u) \ar[r]^{\iota}\ar@{}[ru]|{(1)} & c(\alpha) \ar[r]^{p'} \ar@{}[ru]|{(2)} & c(u)[1]
 }
$$
is a morphism of triangles. In other words, we must prove that (1)
and
(2) commute in $Ho\mc{D}$.\\
Let us see f{i}rst the commutativity of (2), that is
$$\xymatrix@M=4pt@H=4pt@C=20pt@R=15pt{
  c(v) \ar[r]^{p}                                      & Y[1]    & (\mrm{R}Y)[1] \ar[l]_{\lambda_Y[1]}\ar[r]^{I_Y[1]} & c(u)[1]\ar[dd]_{Id}\\
  c(\mrm{R}v) \ar[d]_{\psi}\ar[u]^{\widetilde{\lambda}} &         &                                                    &                   \\
  c(\alpha) \ar[rrr]^{p'}                              &         &                                                    &            c(u)[1] \ .}
$$
Let $p'':c(\mrm{R}v)\rightarrow (\mrm{R}Y)[1]$ be the morphism
induced by $\mrm{R}v$ (see \ref{def{i}triangDisting}). Then
$\lambda_Y[1]\comp p''=p \comp \widetilde{\lambda}$ in $\mc{D}$,
since both morphisms agree with the result of applying the cone
functor to the following compositions
$$\xymatrix@M=4pt@H=4pt{
 \mrm{R}Y \ar[r]^{\mrm{R}v}\ar[d]_{\lambda_Y} & \mrm{R}Z\ar[d]_{\lambda_Z} & \mrm{R}Y  \ar[r]^{\mrm{R}v} \ar[d]_{Id}  & \mrm{R}Z\ar[d] \\
 Y \ar[r]^v   \ar[d]_{Id}                     & Z \ar[d]                   & \mrm{R}Y  \ar[r] \ar[d]_{\lambda_Y}      & 1\ar[d] \\
 Y \ar[r]                                     & 1                          & Y  \ar[r]                                & 1 \ .}$$
Hence, it remains to see that $p'\comp \psi=I_Y[1]\comp p''$, but
again this equality holds in $\mc{D}$ because both morphisms are
equal to the image under the cone functor of the compositions
$$\xymatrix@M=4pt@H=4pt{
 \mrm{R}Y  \ar[r]^{\mrm{R}v} \ar[d]_{I_Y} & \mrm{R}Z\ar[d]_{I_Z} & \mrm{R}Y  \ar[r]^{\mrm{R}v} \ar[d]_{Id} & \mrm{R}Z\ar[d]_{I_Z}   \\
 c(u) \ar[r]^{\alpha}\ar[d]_{Id}          & c(v\comp u)\ar[d]    & \mrm{R}Y  \ar[r] \ar[d]_{I_Y}           & 1 \ar[d] \\
 c(u) \ar[r]                              & 1                    & c(u)  \ar[r]                            & 1 \ .}$$
Now we study the square (1), that consist of
$$\xymatrix@M=4pt@H=4pt@C=20pt@R=15pt{
  c(v\comp u) \ar[rr]^{\beta} \ar[dd]_{Id} &                                                                           & c(v)\ar[d]_{\eta}\\
                                           &                                                                           & \mbf{s}\widetilde{T}\ar[d]_{\Phi}\\
  c(v\comp u)                              & \mrm{R}c(v\comp u)\ar[r]^{I_{c(v\comp u)}}\ar[l]_{\lambda_{c(v\comp u)}}  & c(\alpha)\ .}
$$
The strategy will be the following. We will def{i}ne a simplicial
morphism $\Theta:\widetilde{T}\rightarrow C(v)$
such that\\[0.1cm]
\indent a) If $i_{C(v\comp u)}:C(v\comp u)\rightarrow
\widetilde{T}$ is the canonical inclusion and
$\widetilde{\beta}:C(v\comp u)\rightarrow C(v)$ the simplicial
morphism def{i}ned through the diagram B of
(\ref{def{i}Octaedro}),
then $\Theta\comp i_{C(v\comp u)} =\widetilde{\beta}$.\\
\indent b) If $\rho:C(v)\rightarrow \widetilde{T}$ is the map induced by the canonical inclusions of $Y$ and $Z$ into $C(u)$ and $C(v\comp u)$, then $\Theta\comp\rho=Id_{C(v)}$.\\[0.2cm]
Assume that a) and b) are satisf{i}ed. Since
$\mbf{s}(\rho)=\eta:c(v)\rightarrow \mbf{s}\widetilde{T}$ is an
equivalence, $\theta=\mbf{s}\Theta=(\eta)^{-1}$ in $Ho\mc{D}$. On
the other hand, $\mbf{s}{i_{C(v\comp u)}}=\eta'$ and
$\mbf{s}(\widetilde{\beta})=\beta:c(v\comp u)\rightarrow c(v)$.
Hence we deduce from a) that $\theta\comp\eta'=\beta$ in $\mc{D}$.\\
Therefore, (1) commutes, because on one hand $\eta\comp
\beta=\theta^{-1}\comp \beta =\eta'$, and on the other hand, by
(\ref{diagrama2AuxiliarC4}) $\eta'=\Phi^{-1}\comp I_{c(v\comp
u)}\comp(\lambda_{c(v\comp u)})^{-1}$.\\
Hence, it remains to prove a) and b). Def{i}ne
$\Theta:\widetilde{T}\rightarrow C(v)$ as follows.\\
Recall that $C(u)$ is def{i}ned in degree $n$ as $Y\sqcup
\coprod^{n} X \sqcup 1$. Following the notations in
\ref{def{i}CylenTerminosdeSigma} it can be described as
$$C(u)_n=Y^{u_1}\sqcup \coprod_{\sigma\in\Lambda_n} X^{\sigma}\sqcup 1^{u_0}\ .$$
Similarly, by \ref{def{i}CilTildeAbstracto},
$\widetilde{T}_n=C(v\comp u)_n^{u_1}\sqcup
\coprod_{\sigma\in\Lambda_n} C(u)_n^{\sigma}\sqcup 1^{u_0}$, that
is
$$\widetilde{T}_n=(Z^{u_1,u_1}\sqcup \coprod_{\rho\in\Lambda_n} X^{\rho,u_1}\sqcup 1^{u_0,u_1}) \sqcup\coprod_{\sigma\in\Lambda_n} (Y^{u_1,\sigma}\sqcup \coprod_{\rho\in\Lambda_n} X^{\rho,\sigma}\sqcup 1^{u_0,\sigma})\sqcup 1^{u_0,u_0}$$
where the superscripts are mute, and are just used as labels for
indexing the coproduct. Def{i}ne
$$ \Theta_n:(Z^{u_1,u_1}\sqcup \coprod_{\rho\in\Lambda_n} X^{\rho,u_1}\sqcup 1^{u_0,u_1}) \sqcup\coprod_{\sigma\in\Lambda_n} (Y^{u_1,\sigma}\sqcup \coprod_{\rho\in\Lambda_n} X^{\rho,\sigma}\sqcup 1^{u_0,\sigma})\sqcup 1^{u_0,u_0}
 \longrightarrow Z^{u_1}\sqcup \coprod_{\sigma\in\Lambda_n} Y^{\sigma}\sqcup 1^{u_0}$$
as the morphism whose restriction to the component $\rho,\sigma$
is
$$
\Theta_n|_{\rho,\sigma}=
\begin{cases}
 Id:Z^{u_1,u_1}\rightarrow Z^{u_0} & \mbox{ if } \rho=\sigma=u_1\\
 Id:Y^{u_1,\sigma}\rightarrow Y^{\sigma} & \mbox{ if } \sigma\in\Lambda_n,\,\rho=u_1\\
 Id:1^{u_0,\sigma}\rightarrow 1^{u_0} & \mbox{ if } \rho=u_0\\
 u:X^{\rho,\sigma}\rightarrow Y^{\sigma} & \mbox{ if } \rho,\sigma\in\Lambda_n,\,\sigma^{-1}(1)\subseteq \rho^{-1}(1)\\
 u:X^{\rho,\sigma}\rightarrow Y^{\rho} & \mbox{ if } \rho\in\Lambda_n,\,\sigma\neq u_0,\,\rho^{-1}(1)\subseteq \sigma^{-1}(1)\ .
\end{cases}
$$
Provided that $\Theta$ is an isomorphism of simplicial objects, it is clear that a) and b) hold.\\
Therefore, it remains to see that $\Theta$ is in fact a morphism between simplicial objets.\\
Given an order preserving map $\nu:[m]\rightarrow[n]$, we must
check that $\Theta_m\comp
\widetilde{T}(\nu)=[C(v)](\nu)\comp\Theta_n:\widetilde{T}_n\rightarrow C(v)_m$ in $\mc{D}$.\\
Recall that $[C(v)](\nu):Z^{u_1}\sqcup
\coprod_{\sigma\in\Lambda_n} Y^{\sigma}\sqcup 1^{u_0}\rightarrow
Z^{u_1}\sqcup \coprod_{\sigma\in\Lambda_m} Y^{\sigma}\sqcup
1^{u_0}$ is given by (see \ref{def{i}CylenTerminosdeSigma})
$$[C(v)](\nu)|_{\sigma}=\begin{cases}
 Id:Y^{\sigma}\rightarrow Y^{\sigma\nu}& \mbox{if } \sigma\nu\in\Lambda_m\\
 Id:Z^{u_1}\rightarrow Z^{u_1}& \mbox{if } \sigma=u_1\\
 v:Y^{\sigma}\rightarrow Z^{u_1}& \mbox{if } \sigma\in\Lambda_n\mbox{ and }\sigma\nu=u_1\\
 Id:1^{u_0}\rightarrow 1_{u_0}& \mbox{if } \sigma=u_0\\
 Y^{\sigma}\rightarrow 1^{u_0}& \mbox{if }\sigma\in\Lambda_n\mbox{ and }\sigma\nu=u_0\ .
\end{cases}$$
On the other hand, $\widetilde{T}(\nu):\widetilde{T}_n\rightarrow
\widetilde{T}_m$ is (see \ref{def{i}CilTildeAbstracto})
$$
\widetilde{T}(\nu)|_{\sigma}=
\begin{cases}
   [C(u)](\nu):C(u)_n^\sigma\rightarrow C(u)_m^{\sigma\nu} & \mbox{ if } \sigma\nu\in\Lambda\\
   \widetilde{\alpha}_m \comp [C(u)](\nu):C(u)_n^\sigma\rightarrow C(v\comp u)_m^{u_1} &\mbox{ if } \sigma\in\Lambda, \sigma\nu=u_1\\
    C(u)_n^\sigma\rightarrow 1^{u_1} &\mbox{ if } \sigma\in\Lambda, \sigma\nu=u_0\\
    [C(v\comp u)](\theta):C(v\comp u)_n^{u_1}\rightarrow C(v\comp u)_m^{u_1} &\mbox{ if } \sigma=u_1\\
    Id:1^{u_1}\rightarrow 1^{u_1} &\mbox{ if } \sigma=u_0\ .
\end{cases}
$$
that, by def{i}nition of the cone functor, is equal to
$$\widetilde{T}(\nu)|_{\rho,\sigma}=
\begin{cases}
     Id:X^{\rho,\sigma}\rightarrow X^{\rho\nu,\sigma\nu} &\mbox{ if } \rho\nu\in\Lambda_n,\,\sigma\nu\neq u_0\\
     Id:Y^{u_1,\sigma}\rightarrow Y^{u_1,\sigma\nu} &\mbox{ if } \sigma\nu\in\Lambda_m, \rho=u_1\\
     Id:1^{u_0,\sigma}\rightarrow 1^{u_0,u_0} &\mbox{ if }  \rho=u_0\\
     X^{\rho,\sigma}\rightarrow 1^{u_0,\sigma\nu} &\mbox{ if } \sigma\nu\neq u_0,\,\rho\nu=u_0,\, \rho\in\Lambda_n\\
     u:X^{\rho,\sigma}\rightarrow Y^{u_1,\sigma\nu} &\mbox{ if } \sigma\nu\in\Lambda_m, \rho\nu=u_1,\,\rho\in\Lambda_n\\
     v:Y^{u_1,\sigma}\rightarrow Z^{u_1,u_1} &\mbox{ if } \sigma\nu=u_1,\sigma\in\Lambda_n, \rho=u_1\\
     v\comp u:X^{\rho,\sigma}\rightarrow Z^{u_1,u_1} &\mbox{ if } \sigma\nu=u_1,\rho\in\Lambda_n, \rho\nu=u_1\\
     Id:Z^{u_1,u_1}\rightarrow Z^{u_1,u_1} &\mbox{ if } \sigma=\rho=u_1\\
     X^{\rho,\sigma}\rightarrow 1^{u_0,u_0} &\mbox{ if } \sigma\nu=u_0,\sigma\neq u_0,\,\rho\in\Lambda_n\\
     Y^{u_1,\sigma}\rightarrow 1^{u_0,u_0} &\mbox{ if } \sigma\nu=u_0,\sigma\in\Lambda_n, \rho=u_1\ .
\end{cases}$$
Hence, the equality $\Theta_m\comp
\widetilde{T}(\nu)=[C(v)](\nu)\comp\Theta_n$ is clearly
satisf{i}ed over the components $Z^{u_1,u_1}$, $Y^{u_1,\sigma}$
and $1^{u_0,\sigma}$ of
$\widetilde{T}_n$. Let us check it over the components of the form $X^{\rho,\sigma}$, with $\rho\in\Lambda_n$ and $\sigma\neq u_0$.\\[0.2cm]
\indent Case $\sigma^{-1}(1)\subseteq\rho^{-1}(1)$.\\
In this case $\sigma\neq u_1$ (otherwise $\rho=u_1$), so
$\sigma\in\Lambda_n$ and we have that
$$[C(v)](\nu)\comp\Theta_n|_{X^{\rho,\sigma}}=
\begin{cases}
     u:X^{\rho,\sigma}\rightarrow Y^{\sigma\nu} &\mbox{ if } \sigma\nu\in\Lambda_m\\
     v\comp u:X^{\rho,\sigma}\rightarrow Z^{u_1} &\mbox{ if } \sigma\nu=u_1\\
     X^{\rho,\sigma}\rightarrow 1^{u_0} &\mbox{ if } \sigma\nu=u_0\ .
\end{cases}$$
On the other hand, since $\sigma^{-1}(1)\subseteq\rho^{-1}(1)$
then $\nu^{-1}\sigma^{-1}(1)\subseteq\nu^{-1}\rho^{-1}(1)$, that is, $(\sigma\nu)^{-1}(1)\subseteq(\rho\nu)^{-1}(1)$.\\
Therefore, if $\sigma\nu\in\Lambda_m$ in particular
$(\sigma\nu)^{-1}(1)\neq\emptyset$ and consequently $\rho\nu\neq
u_0$.
If $\rho\nu=u_1$, by def{i}nition
$\tilde{T}(\nu)|_{X^{\rho,\sigma}}=u:X^{\rho,\sigma}\rightarrow
Y^{u_1,\sigma\nu}$ and $\Theta_m\comp
\widetilde{T}(\nu)|_{X^{\rho,\sigma}}=u:X^{\rho,\sigma}\rightarrow Y^{\sigma\nu}$.\\
Otherwise, we have that $\rho\nu\in\Lambda_m$ and then
$\tilde{T}(\nu)|_{X^{\rho,\sigma}}=Id:X^{\rho,\sigma}\rightarrow
X^{\rho\nu,\sigma\nu}$. As
$(\sigma\nu)^{-1}(1)\subseteq(\rho\nu)^{-1}(1)$ then
$\Theta_m\comp
\widetilde{T}(\nu)|_{X^{\rho,\sigma}}=u:X^{\rho,\sigma}\rightarrow Y^{\sigma\nu}$.\\
Now assume that $\sigma\nu=u_1$. Then
$(\sigma\nu)^{-1}(1)=[m]\subseteq (\rho\nu)^{-1}(1)$ and
$\rho\nu=u_1$, so $\widetilde{T}(\nu)|_{X^{\rho,\sigma}}=
     v\comp u:X^{\rho,\sigma}\rightarrow Z^{u_1,u_1}$, and $\Theta_m\comp
\widetilde{T}(\nu)|_{X^{\rho,\sigma}}=v\comp u:X^{\rho,\sigma}\rightarrow Z^{u_1}$.\\
On the other hand, if $\sigma\nu=u_0$, it is clear that $\Theta_m\comp\widetilde{T}(\nu)|_{X^{\rho,\sigma}}:X^{\rho,\sigma}\rightarrow 1^{u_0}$.\\[0.2cm]
\indent Case $\rho^{-1}(1)\subseteq\sigma^{-1}(1)$. Again by
def{i}nition
$$[C(v)](\nu)\comp\Theta_n|_{X^{\rho,\sigma}}=
\begin{cases}
     u:X^{\rho,\sigma}\rightarrow Y^{\rho\nu} &\mbox{ if } \rho\nu\in\Lambda_m\\
     v\comp u:X^{\rho,\sigma}\rightarrow Z^{u_1} &\mbox{ if } \rho\nu=u_1\\
     X^{\rho,\sigma}\rightarrow 1^{u_0} &\mbox{ if } \rho\nu=u_0\ .
\end{cases}$$
Note that $(\rho\nu)^{-1}(1)\subseteq(\sigma\nu)^{-1}(1)$.\\
If $\rho\nu\in\Lambda_m$, it follows that
$(\rho\nu)^{-1}(1)\neq\emptyset$,
then $\sigma\nu\neq u_0$.\\
Hence,
$\widetilde{T}(\nu)|_{X^{\rho,\sigma}}=Id:X^{\rho,\sigma}\rightarrow
X^{\rho\nu,\sigma\nu}$
and $\Theta_m\comp\widetilde{T}(\nu)|_{X^{\rho,\sigma}}=u:X^{\rho,\sigma}\rightarrow Y^{\rho\nu}$.\\
If $\rho\nu=u_1$, we have that $\sigma\nu=u_1$ and $\Theta_m\comp\widetilde{T}(\nu)|_{X^{\rho,\sigma}}=v\comp u:X^{\rho,\sigma}\rightarrow Z^{u_1}$.\\
F{i}nally, if $\rho\nu=u_0$, by def{i}nition
$\widetilde{T}(\nu)|_{X^{\rho,\sigma}}:X^{\rho,\sigma}\rightarrow
1^{u_0,\sigma\nu}$ and
$\Theta_m\comp\widetilde{T}(\nu)|_{X^{\rho,\sigma}}:X^{\rho,\sigma}\rightarrow
1^{u_0}$, that f{i}nish the proof.
\end{proof}

In order to prove the octahedron axiom in the general case, we
will need the following notations.

\begin{num}\label{MitadSupOctaedro}
Denote by $f:X\stackrel{[1]}{\longrightarrow} Y$ a morphism of the
form $f:X\rightarrow Y[1]$ of $Ho\mc{D}$. Then the distinguished
triangle $X\rightarrow Y\rightarrow Z\rightarrow X[1]$ can be
written as
$$\xymatrix@M=4pt@H=4pt@R=10pt{
 X \ar[rr] &               & Y \ar[ld] \\
           & Z\ar[lu]^{[1]}&}$$
We will call ``\textit{octahedron upper half}'' a diagram in
$Ho\mc{D}$ as in the following picture
\begin{equation}\label{MitadSupOctaedro}\xymatrix@M=4pt@H=4pt{
 X \ar[rd]^{u}\ar[rr]^{w} &                         & Z\ar[dd]^s \\
 \ar@{}[r]|{\ast}         & Y\ar[ru]^{v}\ar[ld]^{q} & \ar@{}[l]|{\ast} \\
 M\ar[uu]^{[1]}_t         &                         & N\ar[lu]^{[1]}_r\ar[ll]^{[1]}_p}
\end{equation}
where the triangles labelled with the symbol $\ast$ are
distinguished and the two remaining commute (in $Ho\mc{D}$).
\end{num}

\begin{prop}[\textbf{TR 4, Octahedron axiom}]\index{Index}{octahedron, axiom} \mbox{}\\
Every octahedron upper half can be completed to an octahedron.
More precisely, given an octahedron upper half as
\ref{MitadSupOctaedro}, there exists a diagram
$$\xymatrix@M=4pt@H=4pt{
 X \ar[rr]^{w}               & \ar@{}[d]|{\ast}                 & Z\ar[dd]^{s} \ar[ld]_{v'}\\
                             & Y'\ar[lu]^{[1]}_{u'}\ar[rd]^{r'} & \\
 M\ar[uu]^{[1]}_t\ar[ru]^{q'}& \ar@{}[u]|{\ast}                 & N\ar[ll]^{[1]}_p}
$$
where, again, the triangles labelled with $\ast$ are distinguished
and the others commute $($in $Ho\mc{D})$. Moreover, the following
diagrams commute in $Ho\mc{D}$
$$\xymatrix@M=4pt@H=4pt@R=8pt{
                       & Z\ar[rd]^{v'} &    &                             & X[1]\ar[rd]^{u[1]} &     \\
 Y \ar[ru]^v\ar[rd]_q  &               & Y' & Y' \ar[ru]^{u'}\ar[rd]_{r'} &                    & Y[1] \ . \\
                       & M\ar[ru]_{q'} &    &                             & N\ar[ru]_{r}       &     }
$$
\end{prop}

\begin{proof} First, suppose that $u$ and $v$ are
morphisms of $\mc{D}$. In this case, using the notations given in
\ref{def{i}triangDisting} and \ref{def{i}Octaedro}, it follows
from TR3 that the given octahedron upper half is isomorphic to
$$\xymatrix@M=4pt@H=4pt{
 {X} \ar[rd]^{{u}}\ar[rr]^{{v}\comp{u}}  &                                       & {Z}\ar[dd]^{\iota_{{Z}}} \\
 \ar@{}[r]|{\ast}                        & {Y}\ar[ru]^{{v}}\ar[ld]^{\iota_{{Y}}} & \ar@{}[l]|{\ast}                  \\
 c({u})\ar[uu]^{[1]}_{p_{{u}}}           &                                       & c({v})\ar[lu]^{[1]}_{p_{{v}}}\ar[ll]^{[1]}_{\gamma}\ .}
$$
Hence, it suf{f}{i}ces to prove that this octahedron upper half
can be completed into a whole octahedron.
Consider the distinguished triangle obtained from $v\comp u$
$$\xymatrix@M=4pt@H=4pt@R=15pt{
 {X}\ar[r]^{v\comp{u}}  & Z\ar[r]^{\jmath_{Z}} & c(v\comp{u})\ar[r]^{p_{v\comp{u}}}   & {X}[1]\ .}$$
Following the notations given in \ref{def{i}Octaedro}, we consider
the diagram
$$\xymatrix@M=4pt@H=4pt{
 X \ar[rr]^{{v}\comp{u}}                       & \ar@{}[d]|{\ast}                                       & Z\ar[dd]^{\iota_{Z}} \ar[ld]_{\jmath_{Z}}\\
                                               & c(v\comp u)\ar[lu]^{[1]}_{p_{v\comp u}}\ar[rd]^{\beta} & \\
 c({u})\ar[uu]^{[1]}_{p_{{u}}} \ar[ru]^{\alpha}& \ar@{}[u]|{\ast}                                       & c(v)\ar[ll]^{[1]}_{\gamma}\ .}
$$
I claim that the triangles labelled with $\ast$ are distinguished.
The upper triangle is clearly distinguished, whereas the lower one is so because of proposition \ref{ConoDeComposicion}.\\
Since $\alpha$ and $\beta$ are the morphisms obtained as the image
under the cone functor of the squares A and B in
\ref{def{i}Octaedro}, it follows that the above triangles not
labelled with $\ast$ are commutative, as well as
$$\xymatrix@M=4pt@H=4pt@R=8pt{
                               & Z\ar[rd]^{\jmath_Z}   &             &                                                   & X[1]\ar[rd]^{u[1]} &     \\
 Y \ar[ru]^v\ar[rd]_{\iota_Y}  &                       & c(v\comp u) & c(v\comp u) \ar[ru]^{p_{v\comp u}}\ar[rd]_{\beta} &                    & Y[1] \ , \\
                               & c(u)\ar[ru]_{\alpha}  &             &                                                   & c(v)\ar[ru]_{p_v}  &     }
$$
This f{i}nish the proof of TR 4 when $u,v$ are morphisms of $\mc{D}$.\\[0.1cm]
\indent To see the general case, when $u$ and $v$ are morphisms of $Ho\mc{D}$, let us check that each octahedron upper half (\ref{MitadSupOctaedro}) is isomorphic to an octahedron upper half where $u$ and $v$ are in $\mc{D}$.\\
Since the triangle $\xymatrix@M=4pt@H=4pt@R=10pt{
 X\ar[r]^u & Y\ar[r]^q & M\ar[r]^t & X[1]}$
is distinguished, by def{i}nition there exists a morphism
$\bar{u}:\bar{X}\rightarrow \bar{Y}$ in $\mc{D}$ and an
isomorphism of triangles
$$\xymatrix@M=4pt@H=4pt@R=15pt{
 X\ar[r]^u   \ar[d]_{\tau}     & Y\ar[r]^q \ar[d]_{\tau'}        & M\ar[r]^t\ar[d]_{\tau''}         & X[1]\ar[d]_{\tau[1]}\\
 \bar{X}\ar[r]^{\bar{u}}       & \bar{Y}\ar[r]^{\iota_{\bar{Y}}} & c(\bar{u})\ar[r]^{p_{\bar{u}}}   & \bar{X}[1]\ .}
$$
Hence, the isomorphisms $\tau,\tau',\tau''$ provide an isomorphism
between the given octahedron upper half and the following one
$$\xymatrix@M=4pt@H=4pt{
 \bar{X} \ar[rd]^{\bar{u}}\ar[rr]^{\bar{v}\comp\bar{u}} &                                                   & {Z}\ar[dd]^{s} \\
 \ar@{}[r]|{\ast}                                       & \bar{Y}\ar[ru]^{\bar{v}}\ar[ld]^{\iota_{\bar{Y}}} & \ar@{}[l]|{\ast}                  \\
 c(\bar{u})\ar[uu]^{[1]}_{p_{\bar{u}}}                  &                                                   & N\ar[lu]^{[1]}_{\bar{r}}\ar[ll]^{[1]}_{\bar{p}}}
$$
where $\bar{v}=v\comp(\tau')^{-1}$, $\bar{r}=\tau'\comp r$ and $\bar{p}=\iota_{\bar{Y}}\comp\bar{r}=(\tau'')^{-1}\comp p$.\\
Therefore we can assume that the morphism $u$ in our octahedron
upper half
$$\xymatrix@M=4pt@H=4pt{
 X \ar[rd]^{u}\ar[rr]^{w} &                         & Z\ar[dd]^s \\
 \ar@{}[r]|{\ast}         & Y\ar[ru]^{v}\ar[ld]^{q} & \ar@{}[l]|{\ast} \\
 M\ar[uu]^{[1]}_t         &                         & N\ar[lu]^{[1]}_r\ar[ll]^{[1]}_p}$$
is a morphism in $\mc{D}$.\\
On the other hand, we deduce from theorem \ref{ProHoDcat} that
$v:Y\rightarrow Z$ is represented by a zig-zag of morphisms of
$\mc{D}$ in the form
$$\xymatrix@M=4pt@H=4pt{Y & \mrm{R}Y \ar[l]_{\lambda_Y}\ar[r]^{\bar{v}} & T & \mrm{R}Z \ar[l]_l\ar[r]^{\lambda_{Z}}& Z},\  l\in\mrm{E}\ .$$
F{i}nally, let us see that the original octahedron upper half is
isomorphic to
\begin{equation}\label{DiagAuxOctaedro}\xymatrix@M=4pt@H=4pt{
 \mrm{R}X \ar[rd]^{\mrm{R}u}\ar[rr]^{\bar{v}\comp\mrm{R}u} &                                       & T\ar[dd]^{\bar{s}} \\
 \ar@{}[r]|{\ast}                                    & \mrm{R}Y\ar[ru]^{\bar{v}}\ar[ld]^{\mrm{R}q} & \ar@{}[l]|{\ast} \\
 \mrm{R}M\ar[uu]^{[1]}_{\bar{t}}                    &                                       & \mrm{R}N\ar[lu]^{[1]}_{\hat{r}}\ar[ll]^{[1]}_{\mrm{R}p}\ .}\end{equation}
To this end, consider for any $A$ in $\mc{D}$ the isomorphism
$\theta_A$ of $Ho\mc{D}$ def{i}ned as the composition
$\xymatrix@M=4pt@H=4pt{\mrm{R}(A[1])\ar[r]^{\lambda_{A[1]}} & A[1]\ar[r]^{(\lambda_A[1])^{-1}} & (\mrm{R}A)[1]}$.\\
Set $\bar{t}=\theta_X\comp\mrm{R}t:\mrm{R}M\rightarrow
(\mrm{R}X)[1]$. Then the following diagram commutes
$$\xymatrix@M=4pt@H=4pt@R=15pt{
 \mrm{R}X\ar[r]^{\mrm{R}u}\ar[d]_{\lambda_X}  & \mrm{R}Y\ar[r]^{\mrm{R}q} \ar[d]_{\lambda_Y}   & \mrm{R} M\ar[r]^{\bar{t}}\ar[d]_{\lambda_M}  & (\mrm{R}X)[1]\ar[d]_{\lambda_X[1]}\\
 {X}\ar[r]^{{u}}                              & {Y}\ar[r]^{q}                                  & M\ar[r]^{t}                                  & {X}[1]\ .}
$$
In the same way, $\bar{s}$ and $\hat{r}$ are the respective
compositions
$$\xymatrix@M=4pt@H=4pt{
 T\ar[r]^{l^{-1}}           & \mrm{R}Z \ar[r]^{\mrm{R}s}      & \mrm{R}N \\
 \mrm{R}N \ar[r]^{\mrm{R}r} & \mrm{R}(Y[1]) \ar[r]^{\theta_Y} & (\mrm{R}Y)[1]}$$
that give rise to the isomorphism of triangles
$$\xymatrix@M=4pt@H=4pt@R=25pt@C=30pt{
 \mrm{R}Y\ar[r]^{\bar{v}}\ar[d]_{\lambda_Y}   & T\ar[r]^{\mrm{R}q} \ar[d]_{\lambda_Z\comp l^{-1}}   & \mrm{R} N\ar[r]^{\bar{t}}\ar[d]_{\lambda_N}  & (\mrm{R}X)[1]\ar[d]_{\lambda_X[1]}\\
 {Y}\ar[r]^{w}                                & {Z}\ar[r]^{q}                                       & N\ar[r]^{\hat{r}}                            & {X}[1]\ .}
$$
Therefore it is clear that (\ref{DiagAuxOctaedro}) is an
octahedron upper half isomorphic to the original, that f{i}nish
the proof.
\end{proof}

In order to study the remaining axiom TR 2 of triangulated
category, we need $Ho\mc{D}$ to be additive, since TR 2 involves a
``minus'' sign. Recall that by \ref{aditividad}, if we assume that
$\mc{D}$ is an additive simplicial descent category (def{i}nition
\ref{def{i}CatDescAditiva}) then so is $Ho\mc{D}$.

\begin{prop}[\textbf{TR 2}]\label{TR2}\mbox{}\\
i) Suppose that $\mc{D}$ is an \emph{additive simplicial descent category}.\\
If the triangle $X\stackrel{u}{\rightarrow} Y\stackrel{v}{\rightarrow} Z\stackrel{w}{\rightarrow} X[1]$ is distinguished in $Ho\mc{D}$, then so is $Y\stackrel{v}{\rightarrow} Z\stackrel{w}{\rightarrow} X[1]\stackrel{-u[1]}{\rightarrow}Y[1]$.\\
ii) If in addition the shift functor $T:Ho\mc{D}\rightarrow
Ho\mc{D}$ is fully faithful, so the converse statement also holds.
\end{prop}

\begin{proof}\mbox{}\\
Proof of \textit{i)} By def{i}nition of distinguished triangle we
can assume that $u$ is a morphism of $\mc{D}$, and that
$X\stackrel{u}{\rightarrow} Y\stackrel{v}{\rightarrow}
Z\stackrel{w}{\rightarrow} X[1]$ is the triangle obtained from
$u$, that is
$$\xymatrix@M=4pt@H=4pt@R=15pt{
 {X}\ar[r]^{{u}}  & Y\ar[r]^{\iota_Y} & c({u})\ar[r]^{p_{{u}}}   & {X}[1]\ .}$$
We must prove that the triangle
$$\xymatrix@M=4pt@H=4pt@R=15pt{
  Y\ar[r]^{\iota_Y} & c({u})\ar[r]^{p_{{u}}}   & {X}[1] \ar[r]^{-{u}[1]} & Y[1]  \ .}$$
is distinguished. Def{i}ne the morphism $I_Y:\mrm{R}Y\rightarrow
c(u)$ as in (\ref{inclusionescylenD}). We will see that there
exists an isomorphism of triangles
\begin{equation}\label{diagramaAuxiliar2C4}\xymatrix@M=4pt@H=4pt@R=15pt@C=25pt{
  Y\ar[r]^{\iota_Y}                       & c({u})\ar[r]^{p_{{u}}}             & {X}[1] \ar[r]^{-{u}[1]}                   & Y[1]\\
  \mrm{R}Y\ar[r]^{I_Y}\ar[u]_{\lambda_Y}  & c({u})\ar[r]^{\iota_{c({u})}} \ar[u]_{Id} & c(I_Y) \ar[r]^{p_{I_Y}}\ar[u]_{\theta}& (\mrm{R}Y)[1] \ar[u]_{\lambda_Y[1]}\ .}
\end{equation}
Let $0$ be a zero object of $\mc{D}$, that is at the same time
initial and f{i}nal object. Given a simplicial object $S$ in
$\mc{D}$ it holds that, at the simplicial level, the simplicial
cone of $0\rightarrow S$ is by def{i}nition the simplicial
cylinder of $0\leftarrow 0\rightarrow S$, that coincides with $S$
(since $0$ is the unit for the coproduct).
In particular, if $S=T\times\Dl$ for some object $T$ of $\mc{D}$, it follows that $c(0\rightarrow T)=\mbf{s}(T\times\Dl)=\mrm{R}T$.\\
By assumption $\mbf{s}$ is additive, so $\mrm{R}0=0$.
Consequently, we can consider the morphism $f:c(I_Y)\rightarrow
\mrm{R}(X[1])$ in $\mc{D}$ def{i}ned as the image under the cone
functor of the square
$$\xymatrix@M=4pt@H=4pt@R=15pt{
  \mrm{R}Y\ar[r]^{I_Y} \ar[d] & c(u)\ar[d]^{p_u}\\
  0\ar[r]                     & X[1]  }
$$
obtained from diagram (\ref{CuboConmTriang}) in
\ref{construccionTriangDist}.\\
Let us see that $f\in\mrm{E}$. Consider the following commutative
cube in $\mc{D}$
\begin{equation}\label{diagramaAuxiliar3C4}
\xymatrix@R=30pt@C=30pt@H=4pt@M=4pt@!0{
                                & 0 \ar'[d][dd]\ar[ld]\ar[rr]    &          &  X\ar[dd]_{Id}\ar[ld]_{u} \\
 Y\ar[dd]_{I}\ar[rr]^(0.6){Id}  &                                & Y\ar[dd] &                                \\
                                & 0\ar[ld]\ar'[r][rr]            &          &  X\ar[ld]                   \\
 0      \ar[rr]                 &                                & 0        &                              }
\end{equation}
We will apply the factorization property of the cone,
\ref{coherenciaCono}, to the upper and lower faces of this cube.
Begin with the lower one
$$\xymatrix@M=4pt@H=4pt@R=15pt{
  0\ar[r] \ar[d] & 0\ar[d]\\
  0\ar[r]        & X  \, .}
$$
Applying the cone functor by rows and columns we obtain the
morphisms
$\mrm{R}X\rightarrow 0$ and $0\rightarrow X[1]$ respectively.\\
Denote by $X\{1\}$ the simplicial cone object associated with the
morphism $X\times\Dl\rightarrow 0\times\Dl$. If $\widehat{T}$ is
the image under $\widetilde{Cyl}$ of $0\times\Dl\leftarrow
0\times\Dl{\rightarrow} X<1>$, then $\widehat{T}=X\{1\}$.
Note that $\mbf{s}(X\{1\})=X[1]$ by def{i}nition of $X[1]$.\\
The natural isomorphisms $\Psi':X[1]\rightarrow
 (\mrm{R}X)[1]$ and $\Phi':X[1]\rightarrow
\mrm{R}(X[1])$ obtained from \ref{coherenciaCono} are such that
the diagram
$$\xymatrix@R=50pt@C=50pt@H=4pt@M=4pt@!0{
                 & (\mrm{R}X)[1]\ar[ld]_{Id}\ar[rr]^-{\lambda_X[1]}    &                                     & X[1]\ar[ld]^{Id}           &                                           & \mrm{R}(X[1])\, ,\ar[ld]^{Id}\ar[ll]_-{\lambda_{X[1]}}\\
 (\mrm{R}X)[1]   &                                                             & X[1]  \ar[rr]^{\Phi'}\ar[ll]_{\Psi'}&                                     &  \mrm{R}(X[1])                            &  }$$
Therefore
$\Psi'=(\lambda_X[1])^{-1}$ and $\Phi'=\lambda_{X[1]}^{-1}$.\\
On the other hand, consider now
$$\xymatrix@M=4pt@H=4pt@R=15pt{
  0\ar[r] \ar[d] & X\ar[d]^{u}\\
  Y\ar[r]^{Id}   & Y  \, .}
$$
Let ${g}:\mrm{R}X\rightarrow c(Id_Y)$ be the result of applying
the cone functor by rows to the above square, and
$\widetilde{T}\in\simp\mc{D}$ the image under $\widetilde{Cyl}$ of
the diagram $0\times\Dl\leftarrow
Y\times\Dl\stackrel{i_Y}{\rightarrow}
C(u)$.\\
It follows from \ref{coherenciaCono} the existence of natural
isomorphisms $\Psi:\mbf{s}\widetilde{T}\rightarrow
 c({g})$ and $\Phi:\mbf{s}\widetilde{T}\rightarrow
c(I_Y)$ in $Ho\mc{D}$.\\
Diagram (\ref{diagramaAuxiliar3C4}) provides the morphisms
$$ f:c(I_Y)\rightarrow (\mrm{R}X)[1]\ ; \ f':\mbf{s}\widetilde{T}\rightarrow X[1]\ ; \ f'':c(g)\rightarrow (\mrm{R}X)[1]$$. We deduce from the functoriality of the isomorphisms
$\Psi$, $\Phi$, $\Psi'$ and $\Phi'$ the commutativity in
$Ho\mc{D}$
 of the following diagram
\begin{equation}\label{diagramaAuxiliar4C4}
\xymatrix@M=4pt@H=4pt{
  c(g)  \ar[d]_{f''}                  &  \mbf{s}\widetilde{T} \ar[d]_{f'} \ar[l]_{\Phi}\ar[r]^{\Psi} &  c(I_Y) \ar[d]_{f} \\
  (\mrm{R}X)[1] \ar[r]^{\lambda_X[1]} & X[1]                                                         & \mrm{R}(X[1])\ar[l]_{\lambda_{X[1]}} \ .}
\end{equation}
On the other hand, the morphism $f'':c(g)\rightarrow
(\mrm{R}X)[1]$ is obtained as the image under the cone functor of
the square
$$\xymatrix@M=4pt@H=4pt@R=15pt{
  \mrm{R}X \ar[r]^{g} \ar[d]_{Id} & c(Id_Y)\ar[d]\\
  \mrm{R}X \ar[r]                 & 0  \, .}
$$
Since $Id_{Y}\in\mrm{E}$, from
\ref{axiomaconoenD} we deduce that $c(Id_Y)\rightarrow 0$ is an equivalence. Hence it follows from corollary \ref{exactitudcono} that $f''$ is in $\mrm{E}$.\\
Then by (\ref{diagramaAuxiliar4C4})we f{i}nd that $f\in\mrm{E}$.
Take $\theta=\lambda_{X[1]}\comp f=f''\comp
\Psi^{-1}:c(I_Y)\rightarrow
X[1]$, that is isomorphism in $Ho\mc{D}$ by def{i}nition.\\
We must check that diagram (\ref{diagramaAuxiliar2C4}) is in fact
a morphism of triangles. In other words, we
must see that the squares appearing in this diagram are commutative in $Ho\mc{D}$.\\
The equality $\iota_Y\comp \lambda_Y=I_Y$, follows from the def{i}nition of $\iota_Y$.\\
Let us prove that $\theta\comp\iota_{c(u)}=p_u:c(u)\rightarrow
X[1]$ in $Ho\mc{D}$. By def{i}nition $f$ comes from the
commutative square
$$\xymatrix@M=4pt@H=4pt@R=15pt{
  \mrm{R}Y \ar[r]^{I_Y} \ar[d] & c(u)\ar[d]_{p_u}\\
  0 \ar[r]                     & X[1]  \ ,}
$$
and hence $f\comp I_{c(u)}=\mrm{R}p_u$, so
$\theta\comp\iota_{c(u)}=\lambda_X[1]\comp f\comp I_{c(u)}\comp
\lambda_{c(u)}^{-1}=\lambda_X[1]\comp \mrm{R}p_u\comp
\lambda_{c(u)}^{-1}$. But by the naturality of
$\lambda$, this is just $p_u$.\\
Therefore, it remains to check the equality $-u[1]\comp
\theta=\lambda_Y[1]\comp p_{I_Y}$ in $Ho\mc{D}$. To this end, it is enough to def{i}ne a simplicial morphism $H:X\{1\}\rightarrow \widetilde{T}$ such that\\
a) $f'\comp\mbf{s}H=Id_{X[1]}$, hence $\mbf{s}H=(f')^{-1}$ and $\theta^{-1}=\Psi\comp \mbf{s}H$.\\
b) $\lambda_{Y}[1]\comp p_{I_Y}\comp \Psi\comp \mbf{s}H=-u[1]$.\\
By def{i}nition $X\{1\}_n$ is the coproduct (that is, the direct
sum) of $n$ copies of $X$. We will index this sum over the set
$\{\sigma:[n]\rightarrow [1]\, \sigma\neq u_0,u_1\}=\Lambda_n$
(see \ref{def{i}CylenTerminosdeSigma}). Then
$$X\{1\}_n= \bigoplus_{\sigma\in\Lambda_n} X^{\sigma}\ .$$
On the other hand,
$\widetilde{T}=\widetilde{Cyl}(0\times\Dl\leftarrow
Y\times\Dl\stackrel{i_Y}{\rightarrow} C(u))$, that in degree $n$
is (see \ref{def{i}CilTildeAbstracto})
$\widetilde{T}_n=C(u)_n^{u_1}\sqcup \bigoplus_{\sigma\in\Lambda_n}
Y^{\sigma}$. Again, by def{i}nition of the simplicial cone
functor, $\widetilde{T}$ can be described as
$$\widetilde{T}_n=(Y^{u_1,u_1}\oplus \bigoplus_{\rho\in\Lambda_n} X^{\rho,u_1}) \oplus\bigoplus_{\sigma\in\Lambda_n} Y^{u_1,\sigma}\ .$$
Def{i}ne the restriction of $H_n:X\{1\}_n\rightarrow
\widetilde{T}_n$ to the component $\sigma$ of $X\{1\}_n$ as the
map
$$H_n|_{\sigma}=(Id,-u):X^{\sigma}\rightarrow X^{\sigma,u_1}\oplus Y^{u_1,\sigma}\ .$$
\indent Now we are ready to check a) and b).\\
F{i}rstly, let us prove that $f'\comp \mbf{s}H=Id_{X[1]}$. Let
$Q:C(u)\rightarrow X\{1\}$ be the morphism obtained from the
square
$$\xymatrix@M=4pt@H=4pt@R=15pt{
  X \ar[r]^{u} \ar[d]_{Id} & Y\ar[d]\\
  X \ar[r]                 & 0  \ .}
$$
By def{i}nition $f'=\mbf{s}F'$, where $F':\widetilde{T}\rightarrow
X\{1\}$ is the morphism obtained from applying $\widetilde{Cyl}$
to the diagram
$$\xymatrix@M=4pt@H=4pt@R=15pt{
 0\ar[d] &   Y\ar[l] \ar[r]^{i_Y} \ar[d] & C(u)\ar[d]_{Q}\\
 0       &   0 \ar[l]\ar[r]              & X\{1\}  \ .}
$$
Then $F'_n:(Y^{u_1,u_1}\oplus \bigoplus_{\rho\in\Lambda_n}
X^{\rho,u_1}) \oplus\bigoplus_{\sigma\in\Lambda_n}
Y^{u_1,\sigma}\longrightarrow \bigoplus_{\sigma\in\Lambda_n}
X^{\sigma}$ is given by
$$F'|_{\rho,\sigma}=
\begin{cases}
     0:Y^{u_1,\sigma}\rightarrow X\{1\}_n   &\mbox{ if } \sigma\neq u_0,\, \rho=u_1\\
     Id:X^{\rho,u_1}\rightarrow X^{\rho} & \mbox{if } \rho\in\Lambda_n,\, \sigma=u_1\ .
\end{cases}$$
Therefore it is clear that $F'\comp H=Id_{X\{1\}}$, so $f'\comp \mbf{s}H=Id_{X[1]}$.\\
\indent Proof of b). We may check the commutativity of
$$\xymatrix@M=4pt@H=4pt@R=15pt{
 X[1]  \ar[rr]^{-u[1]} \ar[d]_{\mbf{s}H}   &                         & Y[1]         \\
 \mbf{s}\widetilde{T} \ar[r]^{\Psi}        & c(I_Y) \ar[r]^{p_{I_Y}} & (\mrm{R}Y)[1]\ar[u]_{\lambda_Y[1]} \ .}
$$
Denote by $P:\widetilde{T}\rightarrow Y\{1\}$ the image under
$\widetilde{Cyl}$ of
$$\xymatrix@M=4pt@H=4pt@R=15pt{
 0\ar[d] &   Y\ar[l] \ar[r]^{i_Y} \ar[d]_{Id} & C(u)\ar[d]\\
 0 &   Y \ar[l]\ar[r]                   & 0  \, ,}
$$
Consider now the cube
$$\xymatrix@R=30pt@C=30pt@H=4pt@M=4pt@!0{
                                & 0 \ar'[d][dd]\ar[ld]\ar[rr]    &               &  X\ar[dd]\ar[ld]_{u} \\
 Y\ar[dd]_{Id}\ar[rr]^(0.6){Id} &                                & Y\ar[dd]      &                                \\
                                & 0\ar[ld]\ar'[r][rr]            &               &  0\ar[ld]                   \\
 Y      \ar[rr]                 &                                & 0             &                              }
$$
The isomorphisms provided by \ref{coherenciaCono} are natural, so
the above cube gives rise as before to the following commutative
diagram in $Ho\mc{D}$
$$\xymatrix@M=4pt@H=4pt@R=20pt@C=25pt{
    \mbf{s}\widetilde{T} \ar[r]^{\Psi} \ar[d]_{\mbf{s}P} & c(I_Y)\ar[d]_{p_{I_Y}}\\
    Y[1]                            & (\mrm{R}Y)[1]\ar[l]_{\lambda_Y[1]}  \ ,}$$
consequently $\lambda_Y[1]\comp p_{I_Y}\comp \Psi=\mbf{s}P$ in
$Ho\mc{D}$.\\
Moreover, we have trivially the equality of simplicial morphisms
$P\comp H=-u\{1\}:X\{1\}\rightarrow Y\{1\}$, that is just the
morphism induced by
$$\xymatrix@M=4pt@H=4pt@R=20pt@C=25pt{
    X \ar[r] \ar[d]_{-u} & 0\ar[d]\\
    Y    \ar[r]          & 0  \, .}$$
Then $\lambda_Y[1]\comp p_{I_Y}\comp \Psi\comp
\mbf{s}H=\mbf{s}P\comp \mbf{s}H=-u[1]$, so b) is proven.\\
To f{i}nish the proof it remains to see that $H$
is a morphism of simplicial objects.\\
Recall that $\widetilde{T}_n=(Y^{u_1,u_1}\oplus
\bigoplus_{\rho\in\Lambda_n} X^{\rho,u_1})
\oplus\bigoplus_{\sigma\in\Lambda_n} Y^{u_1,\sigma}$ and if
$\alpha:[m]\rightarrow [n]$ is a morphism of $\Dl$, then
$\widetilde{T}(\alpha):\widetilde{T}_n\rightarrow\widetilde{T}_m$
is
$$
\widetilde{T}(\alpha)|_{\sigma}=\begin{cases}
     Id:Y^\sigma\rightarrow Y^{\sigma\alpha}   &\mbox{ if } \sigma\alpha\in\Lambda_m\\
     (i_Y)_m :Y^\sigma\rightarrow C(u)^{u_1}_m &\mbox{ if } \sigma\in\Lambda_n, \sigma\alpha=u_1\\
     0:Y^\sigma\rightarrow \widetilde{T}_m     &\mbox{ if } \sigma\in\Lambda_n, \sigma\alpha=u_0\\
     [C(u)](\alpha): C(u)_n^{u_1}\rightarrow C(u)_m^{u_1} &\mbox{ if } \sigma=u_1\ .
\end{cases}
$$
That is
$$\widetilde{T}(\alpha)|_{\rho,\sigma}=
\begin{cases}
     Id :Y^{u_1,\sigma}\rightarrow Y^{u_1,\sigma\alpha} &\mbox{ if } \sigma\alpha\neq u_0,\,  \rho=u_1\\
     0:Y^{u_1,\sigma}\rightarrow \widetilde{T}_m     &\mbox{ if } \sigma\in\Lambda_n,\, \sigma\alpha=u_0,\, \rho=u_1\\
 Id:X^{\rho,u_1}\rightarrow X^{\rho\alpha,u_1} & \mbox{if } \rho\alpha\in\Lambda_m,\, \sigma=u_1\\
 u:X^{\rho,u_1}\rightarrow Y^{u_1,u_1} & \mbox{if } \rho\in\Lambda_n\mbox{ and }\rho\alpha=u_1,\,\sigma=u_1\\
 0:X^{\rho,u_1}\rightarrow \widetilde{T}_m& \mbox{if }\rho\in\Lambda_n,\,\rho\alpha=u_0,\,\alpha=u_1 \ .
\end{cases}$$
On the other hand, $X\{1\}_n= \bigoplus_{\sigma\in\Lambda_n}
X^{\sigma}\ $, and the restriction of
$(X\{1\})(\alpha):X\{1\}_n\rightarrow X\{1\}_m$ to the component
$\sigma$ is def{i}ned as
$$(X\{1\})(\alpha)|_{\sigma}=
\begin{cases}
 Id:X^{\sigma}\rightarrow X^{\sigma\alpha} & \mbox{if } \sigma\alpha\in\Lambda_m\\
 0:X^{\sigma}\rightarrow X\{1\}_m & \mbox{if } \sigma\alpha\not\in \Lambda_m\ .
\end{cases}$$
In addition, $H_n:X\{1\}_n\rightarrow \widetilde{T}_n$ is
$H_n|_{\sigma}=(Id,-u):X^{\sigma}\rightarrow
X^{\sigma,u_1}\oplus Y^{u_1,\sigma}$.\\
Let us see that $H_m\comp
X\{1\}(\alpha)=\widetilde{T}(\alpha)\comp H_n$. We have that
$$H_m\comp X\{1\}(\alpha)|_{\sigma}=
\begin{cases}
 (Id,-u):X^{\sigma}\rightarrow X^{\sigma\alpha,u_1}\oplus  Y^{u_1,\sigma\alpha}\mbox{ if }\sigma\alpha\in\Lambda_m\\
 0:X^{\sigma}\rightarrow \widetilde{T}_m \mbox{ if }\sigma\alpha \not\in  \Lambda_m \ .
\end{cases}$$
If $\sigma\alpha\neq u_0,u_1$, it follows from the def{i}nitions
that $\widetilde{T}(\alpha)\comp
H_n|_{\sigma}=(Id,-u):X^{\sigma}\rightarrow
X^{\sigma\alpha,u_1}\oplus  Y^{u_1,\sigma\alpha}$.\\
If $\sigma\alpha=u_0$ then
$\widetilde{T}(\alpha)|_{X^{\sigma,u_1}\oplus Y^{u_1,\sigma}}$ is
the trivial morphism,
so the equality is satisf{i}ed.\\
F{i}nally, if $\sigma\alpha=u_1$ then $\widetilde{T}(\alpha)\comp
H_n|_{\sigma}$ is the composition
$$\xymatrix@M=4pt@H=4pt{X^{\sigma}\ar[r]^-{(Id,-u)} & X^{\sigma,u_1}\oplus Y^{u_1,\sigma}\ar[r]^-{u+Id_Y} & Y^{u_1,u_1}}$$
Then $\widetilde{T}(\alpha)\comp H_n|_{\sigma}=0$, and the proof of \textit{i)} is f{i}nished.\\[0.4cm]
Proof of \textit{ii)}.
Assume that $Y\stackrel{v}{\rightarrow} Z\stackrel{w}{\rightarrow}
X[1]\stackrel{-u[1]}{\rightarrow}Y[1]$ is distinguished. Applying
ii) twice we obtain that $X[1]\stackrel{-u[1]}{\longrightarrow}
Y[1]\stackrel{-v[1]}{\longrightarrow}
Z[1]\stackrel{-w[1]}{\rightarrow}X[2]$ is also distinguished. If
we take the trivial isomorphism of triangles consisting of $\pm
Id$, we deduce that $X[1]\stackrel{u[1]}{\longrightarrow}
Y[1]\stackrel{v[1]}{\longrightarrow}
Z[1]\stackrel{w[1]}{\rightarrow}X[2]$ is distinguished.
By the axiom TR1, the morphism $u:X\rightarrow Y$ can be inserted
into a distinguished triangle
$$X\stackrel{u}{\longrightarrow}
Y\stackrel{v'}{\longrightarrow} Z'\stackrel{w'}{\rightarrow}X[1]$$
If we apply three times i) to it, we obtain the distinguished
triangle $X[1]\stackrel{u[1]}{\longrightarrow}
Y[1]\stackrel{v'[1]}{\longrightarrow}
Z'[1]\stackrel{w'[1]}{\rightarrow}X[2]$.\\
Then it follows from TR3 the existence of an isomorphism
$\Theta:Z[1]\rightarrow Z'[1]$, such that the diagram
$$\xymatrix@M=4pt@H=4pt@C=25pt{
 X[1]\ar[r]^{u[1]}\ar[d]_{Id} & Y[1] \ar[r]^{v[1]}\ar[d]_{Id}   &  Z[1]\ar[d]_{\Theta}    \ar[r]^{w[1]}   &  X[2]\ar[d]_{Id}\\
 X[1]\ar[r]^{u[1]} & Y[1] \ar[r]^{v'[1]}   &  Z'[1]    \ar[r]^{w'[1]}   &  X[2] \ .}
$$
commutes. Since $\mrm{T}$ is fully faithful, there exists an
isomorphism $\Theta':Z\rightarrow Z'$ of $Ho\mc{D}$ such that
$\Theta=\Theta'[1]$, and the diagram
$$\xymatrix@M=4pt@H=4pt@C=25pt{
 X\ar[r]^{u}\ar[d]_{Id} & Y \ar[r]^{v}\ar[d]_{Id}   &  Z\ar[d]_{\Theta'}    \ar[r]^{w}   &  X[1]\ar[d]_{Id}\\
 X\ar[r]^{u} & Y \ar[r]^{v'}   &  Z'    \ar[r]^{w'}   &  X[1] \ .}
$$
is commutative, so the upper triangle is distinguished.
\end{proof}

Summing all up, we have the following

\begin{thm}\label{TmaHoDtriangulada} If $\mc{D}$ is a simplicial descent category then the axioms $TR1$, $TR3$ and $TR4$ of triangulated category hold.\\
If moreover $\mc{D}$ is an additive simplicial descent category then $Ho\mc{D}$ is a ``suspended'' (or right triangulated, \cite{KV}).\\
A functor $F:\mc{D}\rightarrow \mc{D}'$ of additive simplicial
descent categories \ref{def{i}CatDescAditiva} induces a functor of
suspended categories $F:Ho\mc{D}\rightarrow Ho\mc{D}'$.
\end{thm}

\begin{proof}
Except the last part, the theorem is already proven.\\
Let $F:\mc{D}\rightarrow \mc{D}'$ be a functor of additive simplicial descent categories , and $\Theta:\mbf{s}'\comp \simp F\rightarrow F\comp \mbf{s}$ a natural transformation as in def{i}nition \ref{Def{i}FuntorcatDesc}.\\
If $f:X\rightarrow Y$ is a morphism of $\mc{D}$, let us see that
$F(c(f))$ is isomorphic to $c'(F(f))$ in $Ho\mc{D}'$, through a functorial isomorphism $\theta$.\\
Since $F$ is additive and the simplicial cone is def{i}ned
degreewise using direct sums, it follows that the canonical
morphism $\sigma_F:\simp F(C(f\times\Dl))\simeq C(F(f)\times\Dl)$
in $\simp\mc{D}'$. On the other hand $\sigma_F$ commutes with the
canonical inclusions $F(i_Y):F(Y)\times\Dl\rightarrow
\simp F(C(f\times\Dl))$ and $i_{FY}:F(Y)\times\Dl\rightarrow C(F(f)\times\Dl)$.\\
Hence,
$c'(F(f))=\mbf{s}'C(F(f)\times\Dl)\stackrel{\mbf{s}'\sigma_F}{\simeq}
\mbf{s}'\simp F(C(f\times\Dl))\stackrel{\Theta}{\simeq}
F(\mbf{s}(f\times\Dl))=F(c(f))$ gives rise to the natural
isomorphism
$\theta$ between $c'(F(f))$ and $F(c(f))$.\\
In particular $\theta:F(X[1])=F(c(X\rightarrow 0))\simeq
(FX)[1]=c'(FX\rightarrow 0)$ in $Ho\mc{D}'$, and from the
functoriality of $\Theta$ and $\sigma_F$ follows that the diagram
$$\xymatrix@M=4pt@H=4pt@R=15pt{
 F{X}\ar[r]^-{F(f)}\ar[d]_{Id}  & FY \ar[r]^{F(\iota_Y)}\ar[d]_{Id} & c'(F(f))\ar[d]_{\theta}\ar[r]^-{F(p_{{f}})}   &  F({X}[1])\ar[d]_{Id}\\
 F{X}\ar[r]^-{F(f)}             & FY \ar[r]^{F(\iota_Y)}            & F(c({f}))\ar[r]^-{F(p_{{f}})}                 &  F({X}[1])\stackrel{\theta}{\simeq} (FX)[1] \ .}$$
is commutative, providing an isomorphism of triangles.
\end{proof}

\begin{cor}
If $\mc{D}$ is an additive simplicial descent category such that
$\mrm{T}:Ho\mc{D}\rightarrow Ho\mc{D}$ is an automorphism of categories, then $Ho\mc{D}$ is a triangulated category.\\
A functor $F:\mc{D}\rightarrow \mc{D}'$ of additive simplicial
descent categories induces a functor of triangulated categories
$F:Ho\mc{D}\rightarrow Ho\mc{D}'$.
\end{cor}

\begin{obs} In fact, the axioms TR1,$\ldots$,TR4
hold just assuming that $\mc{D}$ is an additive simplicial descent
category and $\mrm{T}$ is fully faithful.
\end{obs}



\chapter{Examples of Simplicial Descent Categories}

In the previous chapters we developed the notion and properties of
simplicial descent categories. Now we introduce examples of such
categories. The f{i}rst one consists of the category of chain
complexes. The axioms of simplicial descent categories will be
checked by hand in this case, whereas in the remaining examples we
will use the transfer lemma to this end.

\section{Chain complexes and homotopy equivalences}\label{SeccionCdcHomt}

\subsection{Preliminaries}

\begin{num} Let $\mc{A}$ be an additive category. Denote by $\cdc$\index{Symbols}{$\cdc$} the category of chain complexes\index{Index}{chain complexes!}
in $\mc{A}$. We will assume that $A$ has numerable sums, that is,
if $\{A_k\}_{k\in\mathbb{Z}}$ is a family of objects of
$\mc{A}$, then $\bigoplus_{k\in\mathbb{Z}}A_k$ exists in $\mc{A}$.\\
If we consider the category of uniformly bounded bellow chain
complexes (see \ref{Def{i}ComplAcotUnif}), for instance positive
chain complexes, the assumption of the existence of numerable sums
in $\mc{A}$ can be dropped, so in this case $\mc{A}$ is just an
additive category $\mc{A}$.
\end{num}

\begin{def{i}}[Double and triple complexes; total functor]\label{totales}\mbox{}\\
$\bullet$ Let $Ch_\ast \cdc$ be the category of double chain
complexes\index{Symbols}{$Ch_\ast\cdc$}\index{Index}{chain
complexes!double}, also called ``naif'' (cf. \cite{Del}). An
object $A$ of $Ch_\ast \cdc$ consists of the data
$$A=\{A_{i_1,i_2}\, ;\ d^1:A_{i_1,i_2}\rightarrow A_{i_1-1,i_2} \, ,\ d^2:A_{i_1,i_2}\rightarrow A_{i_1,i_2-1}\}$$
$$\xymatrix@H=4pt@M=4pt@C=20pt@R=20pt{
 {}     & \vdots  \ar[d]            & \vdots  \ar[d]                   &  \vdots  \ar[d]                    & {}           \\
 \ldots & A_{0,2} \ar[d]^{d^2}\ar[l]& A_{1,2}\ar[l]_{d^1} \ar[d]^{d^2} &  A_{2,2}\ar[l]_{d^1} \ar[d]^{d^2}  & \ldots \ar[l]\\
 \ldots & A_{0,1} \ar[d]^{d^2}\ar[l]& A_{1,1}\ar[l]_{d^1} \ar[d]^{d^2} &  A_{2,1}\ar[l]_{d^1} \ar[d]^{d^2}  & \ldots \ar[l]\\
 \ldots & A_{0,0}  \ar[d]     \ar[l]& A_{1,0}\ar[l]_{d^1} \ar[d]       &  A_{2,0}\ar[l]_{d^1} \ar[d]        & \ldots \ar[l]\\
 {}     & \vdots              & \vdots                     &  \vdots                       & {}
 }$$
\indent The functor $Tot: Ch_\ast \cdc\rightarrow
\cdc$\index{Symbols}{$Tot$}\index{Index}{total functor! of double complexes} is def{i}ned as follows.\\
If $A=\{A_{i_1,i_2} ; d^1:A_{i_1,i_2}\rightarrow A_{i_1 -1,i_2} ,
d^2:A_{i_1,i_2}\rightarrow A_{i_1,i_2-1}\}$ is a double complex,
$Tot A$ is the (single) chain complex given by
$$ (Tot A )_n = \ds\bigoplus_{i_1+i_2=n} A_{i_1,i_2} ; \;\;\; d= \oplus (-1)^{i_2}d^{1}+d^2\ .$$
\noindent $\bullet$ Consider now the category $Ch_\ast
Ch_\ast\cdc=3-\cdc$ of triple chain complexes.\\
Given an object $A=\{A_{i_1,i_2,i_3} ; d^1 , d^2 , d^3 \}$ of
$3-\cdc$, where $d^j$ is the boundary map corresponding to the
index $i_j$, set
$$Tot^{1,2}(A)_{p,q}=\ds\bigoplus_{i_1 + i_2 = p} A_{i_1,i_2,q} ;\;\;\; d^p=\oplus (-1)^{i_2}d^{1}+d^2,\, d^q=\oplus d^{i_3} ;$$
$$Tot^{2,3}(A)_{p,q}=\ds\bigoplus_{i_2 + i_3 = p} A_{q,i_2,i_3} ;\;\;\; d^p=\oplus (-1)^{i_3}d^{2}+d^3,\, d^q=\oplus d^{i_1} ;$$
$$Tot^{1,3}(A)_{p,q}=\ds\bigoplus_{i_1 + i_3 = p} A_{i_1,q,i_3} ;\;\;\; d^p=\oplus (-1)^{i_3}d^{1}+d^3,\, d^q=\oplus d^{i_2} .$$
In this way we obtain the functors $Tot^{1,2},Tot^{2,3},Tot^{1,3}:
3-\cdc\rightarrow
Ch_\ast\cdc$\index{Symbols}{$Tot^{1,2},Tot^{2,3},Tot^{1,3}$}.
\end{def{i}}

\begin{obs}\label{TotordenOpuesto}
Following \cite{Del}, the functor $Tot: Ch_\ast \cdc\rightarrow
\cdc$ is the total functor corresponding  to the order $i_1
> i_2$ and the respective total functor corresponding to
$i_2>i_1$ is canonically isomorphic to the one used here.\\
In other words, if $\Gamma:Ch_\ast \cdc\rightarrow Ch_\ast\cdc$ is
the functor which swaps the indexes of a double complex, then the
following diagram commutes (up to canonical isomorphism)
$$\xymatrix@M=4pt@H=4pt{Ch_\ast \cdc \ar[rr]^{\Gamma}\ar[rd]_{Tot} &       & Ch_\ast \cdc \ar[ld]^{Tot} \\
                                                                & \cdc  & .}$$
\end{obs}

\begin{lema}\label{coherenciatotales}
The functors $Tot: 3-\cdc\rightarrow \cdc$ obtain by composing
$Tot:Ch_\ast\cdc\rightarrow\cdc$ with
$Tot^{1,2},Tot^{2,3},Tot^{1,3}$ are canonically isomorphic.
\end{lema}

\begin{proof}
Following the notations in \cite{Del}, $Tot\small{\comp}
Tot^{1,2}$ is the total functor given by the order $i_3<i_2<i_1$,
whereas $Tot\small{\comp} Tot^{2,3}$ corresponds to $i_1<i_3<i_2$
and $Tot\small{\comp} Tot^{1,3}$ to $i_2<i_3<i_1$. So it follows
from loc. cit. that these three compositions are canonically
isomorphic.
\end{proof}
Moreover, the total functor ``commutes'' with cones in this case.
F{i}rstly, we will remind the classical construction of cone
functor in the chain complex case, that will be also denoted by
$c:Fl(\cdc)\rightarrow \cdc$. In fact, the functor $c$ is obtained
as a particular case of the total functor $Tot$.

\begin{def{i}}\label{def{i}ConoCDC}
If $\mc{B}$ is an additive category, the cone
functor\index{Index}{cone functor!of chain
complexes}\index{Symbols}{$c$}
$$c:Fl(Ch_\ast \mc{B})\rightarrow Ch_\ast(\mc{B})$$
assigns to the morphism $f:X\rightarrow Y$ of chain complexes the
chain complex
$$ c(f)_n = Y_n \oplus X_{n-1}\;\;\;\;\;
d^{cf}=\left(\begin{array}{cc}
               d^Y  & 0\\
               f    & -d^X\end{array}\right) \ .$$
Equivalently, if $\mathcal{J}:Fl(Ch_\ast \mc{B})\rightarrow
Ch_\ast Ch_\ast \mc{B}$ is the functor with $\mc{J}f$ equal to the
double complex
$$\xymatrix@H=4pt@M=4pt{
 {}     &                                  &                                   &                                  & {}           \\
 \ldots &  0 \ar[d]\ar@{}[u]|{\vdots}\ar[l]& 0\ar[l]  \ar[d]\ar@{}[u]|{\vdots} &  0\ar[l] \ar[d]\ar@{}[u]|{\vdots}& \ldots \ar[l]\\
 \ldots &  X_{0} \ar[d]^{f_0}        \ar[l]& X_{1}\ar[l]_{d^X} \ar[d]^{f_1}    &  X_{2}\ar[l]_{d^X} \ar[d]^{f_2}  & \ldots \ar[l]\\
 \ldots &  Y_{0} \ar@{}[d]|{\vdots}\ar[l]  & Y_{1}\ar[l]_{d^Y}\ar@{}[d]|{\vdots} &  Y_{2}\ar[l]_{d^Y} \ar@{}[d]|{\vdots}& \ldots \ar[l] \\
 {}     &                                  &                                   &                                  & {}
}$$ then $c:Fl(Ch_\ast \mc{B})\rightarrow Ch_\ast(\mc{B})$ is the
composition
$$\xymatrix@H=4pt@M=4pt{Fl(Ch_\ast \mc{B})\ar[r]^-{\mc{J}} & Ch_{\ast}Ch_\ast \mc{B} \ar[r]^-{Tot} & Ch_\ast\mc{B} \ .}$$
\end{def{i}}

\begin{lema}\label{totcono}\mbox{}\\
\textbf{i)} If $\mathcal{B}=\cdc$ and $c^\mathcal{A}:
Fl(\cdc)\longrightarrow \cdc$,
$c^\mathcal{B}:Fl(Ch_\ast(\mathcal{B}))\longrightarrow
Ch_\ast(\mathcal{B})$ are the respective cone functors, the
following diagram commutes (up to canonical isomorphism)
$$\xymatrix@M=4pt@H=4pt{ Fl(Ch_\ast(\mathcal{B})) \ar[r]^-{c^\mathcal{B}} \ar[d]^{Tot} & Ch_\ast(\mathcal{B})\ar[d]^{Tot} \\
                         Fl(\cdc) \ar[r]^-{c^\mathcal{A}}                                & \cdc .}$$
\textbf{ii)} The functor $Tot$ preserve homotopies. That is, if
$f,g: X\rightarrow Y$ are homotopic morphisms of chain complexes,
then the induced morphism in the total complex are also homotopic.
\end{lema}

\begin{proof}\mbox{}\\
i) It suf{f}{i}ces to check the commutativity (up to isomorphism)
of the diagrams (I) and (II) bellow
$$\xymatrix@M=4pt@H=4pt{ Fl(Ch_\ast\cdc) \ar[r]^-{\mc{J}^\mathcal{B}} \ar[d]^{Tot} & 3-Ch_\ast\mc{A} \ar[r]^-{Tot_{1,2}}\ar[d]^{Tot_{1,3}}  & Ch_\ast\cdc \ar[d]^{Tot} \\
                         Fl(\cdc) \ar[r]^-{\mc{J}^\mathcal{A}}\ar@{}[ru]|{(I)}    & Ch_\ast\cdc    \ar[r]^-{Tot} \ar@{}[ru]|{(II)}         & \cdc \ .}$$
The commutativity of (II) follows from lemma
\ref{coherenciatotales}, whereas the commutativity of (I)
is an easy computation.\\
Indeed, if $f:B_{n,m}\rightarrow C_{n,m}$ if a morphism of
$Ch_\ast\cdc=Ch_\ast\mc{B}$ then $Tot(f):Tot(B)\rightarrow
Tot(C)$, and $\mc{J}(Tot(f))$ is the double complex
$\{A_{i_1,i_2}; d^1,d^2\}$ given by
$$\begin{cases}
 0            & \mbox{ if }i_2> 1\\
 Tot(B)_{i_1} & \mbox{ if }i_2= 1\\
 Tot(C)_{i_1} & \mbox{ if }i_2=0
 \end{cases} $$
with boundary map $d^1$ is equal to either $d^{Tot(B)}$ or $d^{Tot(C)}$ depending on the case, and $d^2=Tot(f)$.\\
On the other hand, $\mc{J}^\mc{B}(f)$ is the triple complex
$D_{i_1,i_2,i_3}$ given by
$$\begin{cases}
 0         & \mbox{ if }i_2> 1\\
 B_{i_1,i_3} & \mbox{ if }i_2= 1\\
 C_{i_1,i_3} & \mbox{ if }i_2=0
 \end{cases} $$
The boundary map $d^2$ is equal to $f$, whereas $d^1$ is equal to
the boundary map of either $B_{n,m}$ or $C_{n,m}$ (depending on
the case) respect to the index $n$,  and analogously $d^3$ is the
one respective to the index $m$. Thus
$$Tot_{1,3}(D)_{p,q}=\bigoplus_{n+m=p} D_{n,q,m}=\begin{cases}
 0         & \mbox{ if }q> 1\\
 \bigoplus_{n+m=p}B_{n,m}  & \mbox{ if }q= 1\\
 \bigoplus_{n+m=p}C_{n,m} & \mbox{ if }q=0
 \end{cases}\ \ = \mc{J}(Tot(f))_{p,q}\ .$$
Moreover, the boundary maps coincide, since
$d^p:Tot_{1,3}(D)_{p,q}\rightarrow Tot_{1,3}(D)_{p-1,q}$ is by
def{i}nition
$$d^p=\begin{cases}
 0         & \mbox{ if }q> 1\\
 \oplus (-1)^m d_B^n + d_B^m  & \mbox{ if }q= 1\\
 \oplus(-1)^m d_C^n + d_C^m   & \mbox{ if }q=0
 \end{cases}= \begin{cases}
 0         & \mbox{ if }q> 1\\
 d:Tot(B)_p\rightarrow Tot(B)_{p-1}  & \mbox{ if }q= 1\\
 d:Tot(C)_p\rightarrow Tot(C)_{p-1}   & \mbox{ if }q=0
 \end{cases}\ .$$
F{i}nally, $d^q:Tot_{1,3}(D)_{p,q}\rightarrow
Tot_{1,3}(D)_{p,q-1}$
is $\oplus d_D^{i_2}=\oplus f =Tot(f)$.\\[0.2cm]
ii) can be deduced from i) having in mind that a morphism
$p:X\rightarrow Y$ of $Ch_\ast(\mc{B})$ is homotopic to $0$ if and
only if $p$ can be extended to the cone of $X$. Hence, if $f$ and
$g$ are homotopic in $Ch_\ast(\mc{B})$ then $\exists
H:c^{\mc{B}}(X)\rightarrow Y $ such that the diagram
$$\xymatrix@M=4pt@H=4pt@C=15pt@R=17pt{X \ar[rr]^{f-g}\ar[rd] &                          & Y \\
                                                      & c^{\mc{B}}(X)\ar[ru]_{H} &    .}$$
is commutative. Applying $Tot$ to the previous diagram, it follows
from i) that $Totf-Totg$ can be extended to $c(Tot X)$, so the
statement is proven.
\end{proof}

\begin{obs}\label{pdadscilindro}

Recall that the cone functor $c:Fl(\cdc)\rightarrow \cdc$ satisf{i}es the following properties\\
i) $f$ is a homotopy equivalence if and only if $c(f)$ is
contractible, that is, is and only if the morphism
$c(f)\rightarrow 0$ is a homotopy equivalence.\\
ii) if $\mc{A}$ is abelian, $f$ induces an isomorphism in homology
if and only if the homology of $c(f)$ is equal to $0$.
\end{obs}

\begin{def{i}}\label{estructuraDescensoCDC}\mbox{}
\textbf{Simple functor:} The simple functor $\mathbf{s}:
\simp\cdc\longrightarrow \cdc$ is def{i}ned as the composition
$$ \xymatrix@M=4pt{ \simp\cdc \ar[r]^-{K} & Ch_\ast\cdc \ar[r]^-{Tot} & \cdc}$$
where $K(\{X,d_i,s_j\})=\{X,\sum (-1)^i d_i\}$.

More explicitly, let $X=\{X_n , d_i , s_j \}$ be a simplicial
chain complex. Then, each $X_n$ is a chain complex, that will be
referred to as $\{ X_{n,p} , d_{X_n} \}_{p\in\mathbb{Z}}$.

Note that $X$ induces the double complex (\ref{cdc}), with
vertical boundary map $d_{X_n}: X_{n,p}\rightarrow X_{n,p-1}$ and
horizontal boundary map $\partial: X_{n,p}\rightarrow
X_{n-1,p}$ def{i}ned as $\partial=\sum_{i=0}^{n}(-1)^{i} d_i \,$.\\
\begin{equation}\label{cdc}
\xymatrix{
{}     & \vdots  \ar[d]                      & \vdots                                         \ar[d]                       & \vdots        \ar[d]                                              & {}    \\
\ldots & X_{n-1,p+1} \ar[d]^{d_{X_{n-1}}} \ar[l]         & X_{n,p+1}\ar[l]_{\partial} \ar[d]^{d_{X_{n}}}           & X_{n+1,p+1} \ar[l]_{\partial} \ar[d]^{d_{X_{n+1}}} & \ldots \ar[l]\\
\ldots & X_{n-1,p}  \ar[d]^{d_{X_{n-1}}} \ar[l]         &  X_{n,p}\ar[l]_{\partial}   \ar[d]^{d_{X_{n}}}           & X_{n+1,p} \ar[l]_{\partial} \ar[d]^{d_{X_{n+1}}}     & \ldots \ar[l]\\
\ldots & X_{n-1,p-1}   \ar[d]            \ar[l]         &  X_{n,p-1}\ar[l]_{\partial} \ar[d]                       & X_{n+1,p-1} \ar[l]_{\partial}  \ar[d]     & \ldots \ar[l]\\
   {}  & \vdots                                             & \vdots                                                                 &  \vdots & {}
   }\end{equation}

Thus, the image under the simple functor of $X$ is the chain
complex $\mathbf{s}X$ given by

$$\begin{array}{llr} (\mathbf{s}X)_q=\ds\bigoplus_{p+n=q} X_{n,p} & &
 d= \oplus (-1)^p\partial + d_{X_n}:\ds\bigoplus_{p+n=q} X_{n,p}\longrightarrow\ds\bigoplus_{p+n=q-1} X_{n,p} \, .\end{array} $$
\textbf{Weak equivalences:} Def{i}ne $\mrm{E}$ as the class of
homotopy equivalences.\\[0.2cm]
\textbf{Transformation $\mathbf{\lambda}$:} Given $A\in\cdc$,
$\mbf{s}(A\times\Dl)$ in degree $n$ is $\bigoplus_{k\leq n}A_k$,
in such a way that $A$ is canonically a direct summand of
$\mbf{s}(A\times\Dl)$. Then,
$\lambda_A:\mbf{s}(A\times\Dl)\rightarrow A$
is just the projection.\\[0.2cm]
\textbf{Transformation $\mathbf{\mu}$:} Given
$Z\in\simp\simp\cdc$, $\mu_Z:\mbf{s}\mrm{D}(Z)\rightarrow
\mbf{s}\simp\mbf{s}(Z)$ is obtained from the
Alexander-Whitney map \ref{Alex-Whitney-Map}.\\
In degree $n$, the restriction of
$(\mu_Z)_n:\ds\oplus_{p+q=n}Z_{p,p,q}\rightarrow
\ds\oplus_{i+j+q=n}Z_{i,j,q}$ to the component $Z_{p,p,q}$ is
$\oplus_{i+j=p}{\mu_Z}_{i,j,q}:Z_{p,p,q}\rightarrow Z_{i,j,q}$,
where
$${\mu_Z}_{i,j,q}=Z(d^0\stackrel{j)}{\cdots}d^0,d^pd^{p-1}\cdots d^{j+1})\ .$$
\end{def{i}}

\begin{prop}\label{cdcaditivCDS}
Let $\mc{A}$ be an additive category with numerable sums. Then
$\cdc\mc{A}$ together with the homotopy equivalences, and together
with the simple functor, $\lambda$ and $\mu$ def{i}ned above is an
additive simplicial descent category. In addition, $\mu$ is
associative and $\lambda$ is quasi-invertible.
\end{prop}

\begin{proof}[\begin{large}\textbf{Proof of \ref{cdcaditivCDS}}\end{large}]\mbox{}\\
The f{i}rst two axioms are well known, whereas (SDC 3) follows
from the additivity of $\mbf{s}$ (since it is the composition of
additive functors).\\[0.3cm]
\textsl{Proof of axiom $\mathrm{(SDC\; 4)}$:}
We must check that the diagram bellow commutes up to natural
homotopy equivalence

$$\xymatrix@M=4pt@H=4pt@C=30pt{ \simp\simp\cdc \ar[r]^{\simp K}\ar[dd]_{\mrm{D}} & \simp Ch_\ast \cdc \ar[r]^{\simp Tot} & \simp\cdc \ar[d]_K \\
                                                                  &                                   &  Ch_\ast\cdc \ar[d]_{Tot}  \\
                                \simp\cdc \ar[r]^K                & Ch_\ast\cdc \ar[r]^{Tot}             & \cdc.}$$
Following the notations given in \ref{totales}, we can split our
diagram into
\begin{equation}\label{factorizacioncdc}\xymatrix@M=4pt@H=4pt@C=30pt{ \simp\simp\cdc \ar[r]^{\simp K}\ar[dd]_{\mrm{D}} & \simp Ch_\ast \cdc \ar[r]^{\simp Tot} \ar[d]_K         & \simp\cdc \ar[d]_K \\
                                                                         & 3-\cdc \ar[d]_{Tot^{1,2}} &  Ch_\ast\cdc \ar[d]_{Tot}  \\
                                \simp\cdc \ar[r]^K                       & Ch_\ast\cdc \ar[r]^{Tot}  \ar@{}[uul]|{(\mrm{I})}        & \cdc  \ar@{}[uul]|{\#}.}\end{equation}
The right hand side of (\ref{factorizacioncdc}) is commutative,
whereas (I) commutes up to homotopy, by the
Eilenberg-Zilber-Cartier's theorem (see \ref{E-Z-C}), taking
$\mc{U}=\cdc$.\\
Then, given $Z\in\simp\simp\cdc$, consider
$\mu_{E-Z}(Z):K\mrm{D}(Z)\rightarrow Tot^{1,2} K \simp K (Z)$ as
in \ref{E-Z-C}. Hence $\mu=Tot\comp
\mu_{E-Z}:\mbf{s}\mrm{D}\rightarrow \mbf{s}\simp\mbf{s}$ is a
homotopy equivalence because $Tot$ preserve homotopies.\\[0.3cm]
\textsl{Proof if axiom $\mathrm{(SDC\; 5)}$:}\\
F{i}rstly, consider an additive category $\mc{B}$ and the functor
$\mc{I}:\mc{B}\rightarrow Ch_\ast(\mc{B})$ that maps $A$ into the
complex $\mc{I}(A)_0=A$, $\mc{I}(A)_n=0$ if
$n>0$.\\
Let us see that there exists a functor $G:\mc{B}\rightarrow
Ch_\ast \mc{B}$ such that $G(A)$ is contractible for every $A$.
Recall that contractible means that the identity over $G(A)$ is
homotopic to the zero morphism. In addition, we will check that
$K(-\times\Dl) = \mc{I}\oplus G$.\\[0.1cm]
Indeed, if $A\in\mc{B}$ then $K(A\times\Dl)$ is the chain complex
$$\xymatrix@H=4pt@M=4pt{ \cdots & 0 \ar[l]&  A\ar[l] & \ar[l]_{0} A & A \ar[l]_{Id} & A \ar[l]_0 & A \ar[l]_{Id} & \ar[l] \cdots.}$$
Then, def{i}ne $G(A)$ as
$$\xymatrix@H=4pt@M=4pt{ \cdots & 0 \ar[l]& 0\ar[l] & \ar[l] A & A \ar[l]_{Id} & A \ar[l]_0 & A \ar[l]_{Id} & \ar[l] \cdots.}$$
Clearly, $K(A\times\Dl)=\mc{I}(A)\oplus G(A)$. To see that $G(A)$
is contractible, just take the homotopy $H_n
=Id_A :G(A)_n \rightarrow G(A)_{n+1}$, $n>0$.\\[0.2cm]
Set now $\mc{B}=Ch_\ast\mc{A}$. If $A$ is a chain complex of
$\mc{B}$, it holds in $Ch_\ast(\mc{B})$ that
$K(A\times\Dl)=\mc{I}(A)\oplus G(A)$ and then
$$\mbf{s}(A\times\Dl)=Tot K(A\times\Dl)=Tot\mc{I}(A)\,\oplus\, Tot\,G(A)=A\oplus TotG(A).$$
By def{i}nition $\lambda_A$ is the projection $A\oplus
Tot\,G(A)\rightarrow A$. Since $Tot$ preserve homotopies and is
additive, then $Tot\,G(A)$ is contractible. Hence
$\lambda_A$ is a homotopy equivalence.\\[0.3cm]
\textsl{Proof of axiom $\mathrm{(SDC\; 6)}$:}\\
The case $Y=0\times\Dl$ can be found, for instance, in \cite{B}.\\
The general case follows from proposition \ref{AxiomaExactitudPrima}.\\[0.3cm]
\textsl{Proof of axiom $\mathrm{(SDC\; 7)}$:}\\
Let $\mc{B}$ an additive category.
The functor $K:\simp\mc{B}\rightarrow Ch_\ast\mc{B}$ induces in a
natural way a functor
$Fl(\simp\mc{B})\rightarrow Fl(Ch_\ast\mc{B})$, that will be denoted also by $K$.\\
Assume proven the commutativity up to homotopy equivalence of the
diagram bellow\footnote{ This is a well known fact and due to Dold
and Puppe. An analogous proof for the cosimplicial case appears in
\cite{H} p.21}
\begin{equation}\label{cylcompatibilidad}
\xymatrix@H=4pt@M=4pt{Fl(\simp\mc{B})\ar[r]^{C} \ar[d]_{{K}} & \simp\mc{B} \ar[d]^{K} \\
                        Fl(Ch_\ast\mc{B}) \ar[r]^{c^{\mc{B}}}                     & Ch_\ast\mc{B} \, ,}\end{equation}
where $C$ denotes the simplicial cone functor
\ref{def{i}ConoSimplicial}, and $c:Fl(Ch_\ast\mc{B})\rightarrow
Ch_\ast\mc{B}$  is the cone of chain complexes given in
\ref{def{i}ConoCDC}.\\
In this case, if $\mc{B}=\cdc$, since $Tot$ preserves homotopies
(see \ref{totcono}) we obtain that $Tot\comp K\comp C=\mbf{s}\comp
C$ is homotopic to $Tot\comp c^{\mc{B}}\comp {K}$. Again by lemma
\ref{totcono}, $Tot\comp c^{\mc{B}}$ is isomorphic to
$c^{\mc{A}}\comp Tot$. Therefore $\mbf{s}\comp C$ is homotopic to
$c^{\mc{A}}\comp Tot\comp {K}=c^{\mc{A}}\comp\mbf{s}$.
Equivalently, the following diagram commutes up to (natural)
homotopy equivalence
\begin{equation}\label{compatibilidadConoconSimple}
\xymatrix@H=4pt@M=4pt{Fl(\simp\cdc)\ar[r]^{C} \ar[d]_{\mbf{s}} & \simp\cdc \ar[d]^{\mbf{s}} \\
                      Fl(Ch_\ast\mc{A}) \ar[r]^{c^{\mc{A}}}                     & Ch_\ast\mc{A} \, .}\end{equation}
Hence, given $f:A\rightarrow B$ in $\simp\cdc$, it follows that
$\mbf{s}Cf\rightarrow 0$ is a homotopy equivalence if and only if
$c^{\mc{A}}(\mbf{s}f)\rightarrow 0$ is so, and this happens if and
only if $\mbf{s}f$ is a homotopy equivalence, by the classical
properties satisf{i}ed by the cone functor
$c^{\mc{A}}:Fl(\cdc)\rightarrow \cdc$. Consequently (SDC 7) would be proven.\\
Hence it remains to prove the commutativity up to equivalence of
diagram (\ref{cylcompatibilidad}). Indeed, let $f:A\rightarrow B$
be a morphism of simplicial chain complexes. By def{i}nition
$C(f)$ is the total simplicial object associated with the
following biaugmented bisimplicial object (see \ref{conoSimpl})
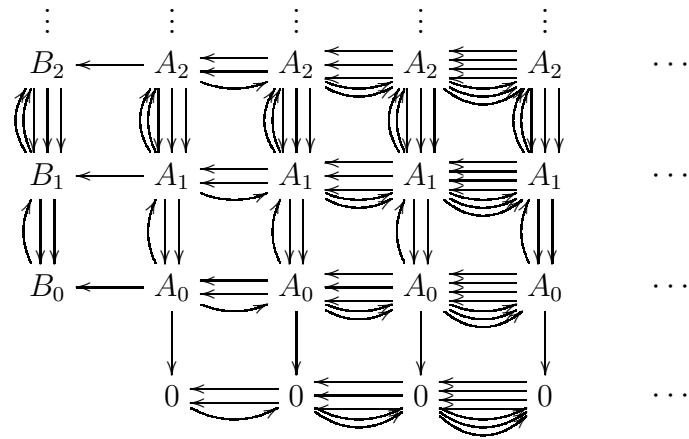
$$\xymatrix@M=4pt@H=4pt@C=25pt{                                                                          &                                                                                             &                                                                                                                                              &                                                                                                                                                          &                                                                                                                             & \\
                                      B_{2}\ar@{}[u]|{\vdots}\ar@<0ex>[d]\ar@<1ex>[d] \ar@<-1ex>[d]      & {A_{2}}\ar[l]\ar@<0ex>[d]\ar@<1ex>[d] \ar@<-1ex>[d]\ar@{}[u]|{\vdots}\ar@/_0.75pc/[r]     & {A_{2}} \ar@{}[u]|{\vdots}\ar@<0.5ex>[l] \ar@<-0.5ex>[l]\ar@<0ex>[d]\ar@<1ex>[d] \ar@<-1ex>[d]\ar@/_1pc/[r]\ar@/_0.75pc/[r]                &A_{2} \ar@{}[u]|{\vdots}\ar@<0ex>[l]\ar@<1ex>[l] \ar@<-1ex>[l]\ar@<0ex>[d]\ar@<1ex>[d] \ar@<-1ex>[d]\ar@/_1pc/[r]\ar@/_0.75pc/[r]\ar@{-}@/_1.25pc/[r]   & A_{2}\ar@{}[u]|{\vdots}\ar@<0.33ex>[l]\ar@<-0.33ex>[l]\ar@<1ex>[l]\ar@<-1ex>[l]\ar@<0ex>[d]\ar@<1ex>[d] \ar@<-1ex>[d]     & \cdots \\
                                      B_{1}\ar@<0.5ex>[d] \ar@<-0.5ex>[d] \ar@/^1pc/[u]\ar@/^0.75pc/[u]  & {A_{1}}\ar[l]\ar@<0.5ex>[d] \ar@<-0.5ex>[d] \ar@/^1pc/[u]\ar@/^0.75pc/[u]\ar@/_0.75pc/[r] &  A_{1} \ar@/^1pc/[u]\ar@/^0.75pc/[u]\ar@<0.5ex>[l]\ar@<-0.5ex>[l]\ar@<0.5ex>[d] \ar@<-0.5ex>[d]\ar@/_1pc/[r]\ar@/_0.75pc/[r]               &A_{1} \ar@/^1pc/[u]\ar@/^0.75pc/[u] \ar@<0ex>[l]\ar@<1ex>[l] \ar@<-1ex>[l] \ar@<0.5ex>[d] \ar@<-0.5ex>[d]\ar@/_1pc/[r]\ar@/_0.75pc/[r]\ar@/_1.25pc/[r]  & A_{1}\ar@<0.33ex>[l]\ar@<-0.33ex>[l]\ar@<1ex>[l]\ar@<-1ex>[l] \ar@<0.5ex>[d] \ar@<-0.5ex>[d]\ar@/^1pc/[u]\ar@/^0.75pc/[u] & \cdots \\
                                      B_{0}  \ar@/^0.75pc/[u]                                            &  A_{0}\ar[d] \ar[l]\ar@/_0.75pc/[r] \ar@/^0.75pc/[u]                                      & {A_{0}} \ar[d] \ar@<0.5ex>[l] \ar@<-0.5ex>[l]  \ar@/^0.75pc/[u]  \ar@/_1pc/[r]\ar@/_0.75pc/[r]                                             &A_{0} \ar[d] \ar@<0ex>[l]\ar@<1ex>[l] \ar@<-1ex>[l]\ar@/^0.75pc/[u]\ar@/_1pc/[r]\ar@/_0.75pc/[r]\ar@/_1.25pc/[r]                                        & A_{0}\ar[d] \ar@<0.33ex>[l]\ar@<-0.33ex>[l]\ar@<1ex>[l]\ar@<-1ex>[l] \ar@/^0.75pc/[u]                                     & \cdots \\
                                                                                                         &  0 \ar@/_0.75pc/[r]                                                                       & 0 \ar@<0.5ex>[l]\ar@<-0.5ex>[l]\ar@/_1pc/[r]\ar@/_0.75pc/[r]                                                                               &0 \ar@<0ex>[l]\ar@<1ex>[l] \ar@<-1ex>[l]\ar@/_1pc/[r]\ar@/_0.75pc/[r]\ar@/_1.25pc/[r]                                                                   & 0\ar@<0.33ex>[l]\ar@<-0.33ex>[l]\ar@<1ex>[l]\ar@<-1ex>[l]                                                                 & \cdots
                                                                                                         }$$
The image under $K$ of $C(f)$ is the same as the total chain
complex of the double complex obtained by applying $K$ to
$$
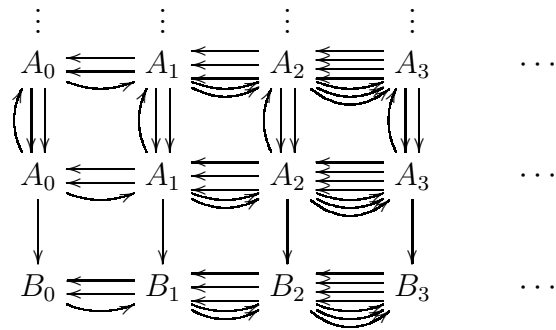
\xymatrix@M=4pt@H=4pt@C=25pt{                                                                               &                                                                                                                                 &                                                                                                                                             &                                                                                                                & \\
                                     {A_{0}}\ar@<0.5ex>[d] \ar@<-0.5ex>[d] \ar@{}[u]|{\vdots}\ar@/_0.75pc/[r] &  A_{1} \ar@{}[u]|{\vdots}\ar@<0.5ex>[l]\ar@<-0.5ex>[l]\ar@<0.5ex>[d] \ar@<-0.5ex>[d]\ar@/_1pc/[r]\ar@/_0.75pc/[r]               &A_{2} \ar@{}[u]|{\vdots} \ar@<0ex>[l]\ar@<1ex>[l] \ar@<-1ex>[l] \ar@<0.5ex>[d] \ar@<-0.5ex>[d]\ar@/_1pc/[r]\ar@/_0.75pc/[r]\ar@/_1.25pc/[r]  & A_{3}\ar@<0.33ex>[l]\ar@<-0.33ex>[l]\ar@<1ex>[l]\ar@<-1ex>[l] \ar@<0.5ex>[d] \ar@<-0.5ex>[d]\ar@{}[u]|{\vdots} & \cdots \\
                                      A_{0}\ar[d] \ar@/_0.75pc/[r] \ar@/^0.75pc/[u]                           & {A_{1}} \ar[d] \ar@<0.5ex>[l] \ar@<-0.5ex>[l]  \ar@/^0.75pc/[u]  \ar@/_1pc/[r]\ar@/_0.75pc/[r]                                  &A_{2} \ar[d] \ar@<0ex>[l]\ar@<1ex>[l] \ar@<-1ex>[l]\ar@/^0.75pc/[u]\ar@/_1pc/[r]\ar@/_0.75pc/[r]\ar@/_1.25pc/[r]                             & A_{3}\ar[d] \ar@<0.33ex>[l]\ar@<-0.33ex>[l]\ar@<1ex>[l]\ar@<-1ex>[l] \ar@/^0.75pc/[u]                          & \cdots \\
                                      B_{0} \ar@/_0.75pc/[r]                                                  & B_{1} \ar@<0.5ex>[l]\ar@<-0.5ex>[l]\ar@/_1pc/[r]\ar@/_0.75pc/[r]                                                                &B_{2} \ar@<0ex>[l]\ar@<1ex>[l] \ar@<-1ex>[l]\ar@/_1pc/[r]\ar@/_0.75pc/[r]\ar@/_1.25pc/[r]                                                    & B_{3}\ar@<0.33ex>[l]\ar@<-0.33ex>[l]\ar@<1ex>[l]\ar@<-1ex>[l]                                                  & \cdots}$$
This double complex is homotopic by columns to
$$\xymatrix@H=4pt@M=4pt{
                             &                                   &                                     & {}    \\
 0 \ar[d] \ar@{}[u]|{\vdots} & 0\ar[l]  \ar[d] \ar@{}[u]|{\vdots}&  0\ar[l] \ar[d] \ar@{}[u]|{\vdots}  & \ldots \ar[l]\\
 A_{0} \ar[d]^{f_0}          & A_{1}\ar[l]_{d^{KA}} \ar[d]^{f_1} &  A_{2}\ar[l]_{d^{KA}} \ar[d]^{f_2}  & \ldots \ar[l]\\
 B_{0}                       & B_{1}\ar[l]_{d^{KB}}              &  B_{2}\ar[l]_{d^{KB}}               & \ldots \ar[l]}$$
whose associated  complex is just $c(Kf)$.\\
We will give explicitly a homotopy equivalence between $K(C(f))\in
Ch_\ast\mc{B}$ and $F=c({K}(f))$. We have that $F$ is the chain
complex given by
$$ F_n=c({K}(f))_n=B_n\oplus A_{n-1}\,  \; ; \ \ d^F(b,a)=(d^{KB}(b)+f(a),-d^{KA}(a)). $$
By def{i}nition, it holds that
$$ K(C(f))_n = B_n \oplus A_{n-1}\oplus \cdots \oplus A_0 \ \mbox{ and }\ d^{K(C(f))}\mbox{ is }$$
$$d^{KB_n}+ f_{n-1} + \ds\sum_{k=1}^{n}(-1)^k d^{KA_{n-k}} +\!\!\!\! \ds\sum_{1\leq k \leq n/2}\!\!\!\! Id_{A_{n-2k-1}}:K(C(f))_n\rightarrow K(C(f))_{n-1}.$$
Consider the chain complex $\widetilde{A}$ def{i}ned as
$$\widetilde{A}_n=A_{n-2}\oplus\cdots\oplus A_{0} \mbox{ if }n\geq 2 \mbox{ and } \widetilde{A}_0=\widetilde{A}_1=0 $$
$$d^{\widetilde{A}}= \ds\sum_{k=2}^{n}(-1)^k d^{KA_{n-k}} + \ds\sum_{1\leq k \leq n/2} Id_{A_{n-2k-1}}:\widetilde{A}_n \rightarrow \widetilde{A}_{n-1}.$$
Then $K(C(f))= F \oplus \widetilde{A}$, since it holds that
$d^F\oplus d^{\widetilde{A}}=d^{K(C(f))}$. In addition, $\widetilde{A}$ is contractible.\\
To see that, let $h_n =\ds\sum_{k=2}^{n+1} Id_{A_{n-k}}:
A_{n-2}\oplus \cdots\oplus A_{0} \rightarrow A_{n-1}\oplus
A_{n-2}\oplus\cdots\oplus A_{0} $ if $n\geq 2$ and $h_0=h_1=0$. If
$n\geq 2$ then
$$h_{n-1}d^{\widetilde{A}}+d^{\widetilde{A}}h_n\!=\!\!\ds\sum_{k=2}^{n}(-1)^k d^{A_{n-k}} + \!\!\!\!\!\!\ds\sum_{1\leq k \leq n/2}\!\!\!\!\!\! Id_{A_{n-2k-1}} + \!\ds\sum_{k=3}^{n+1}(-1)^k d^{A_{n+1-k}} + \!\!\!\!\!\!\!\!\!\ds\sum_{1\leq k \leq (n+1)/2}\!\!\!\!\! Id_{A_{n-2k}}\!=\!Id_{\widetilde{A}_n}$$
Thus, the projection $K(C(f))\rightarrow c({K}(f))$ is a homotopy
equivalence.\\[0.3cm]
\textsl{Proof of axiom $\mathrm{(SDC\; 8)}$:}\\
The functor $\Upsilon:\simp\cdc\rightarrow\simp\cdc$ assigns to
the simplicial chain complex $X$ the chain complex $\Upsilon X$
whose face morphisms are $d^{\Upsilon
X}_i=d^X_{n-i}:X_n\rightarrow X_{n-1}$. Then $K(\Upsilon X)$ has
as horizontal boundary map the morphism $d^{K(\Upsilon
X)}=\sum_{i=0}^n (-1)^i d^X_{n-i}=(-1)^n d^{KX}$, whereas the
vertical boundary
maps coincide in both cases.\\
The double complexes $KX$ and $K(\Upsilon X)$ are then canonically
isomorphic. Composing with $Tot$, we deduce that the functors
$\mbf{s},\mbf{s}\comp\Upsilon:\simp\cdc\rightarrow\cdc$  are
also canonically isomorphic, so (SDC 8) is proven.\\[0.3cm]
\textsl{Compatibility between $\lambda$ and $\mu$:}\\
Given $X$ in $\simp\cdc$, we must check that the following
composition of morphisms are equal to the identity
$$\xymatrix@M=4pt@H=4pt@R=10pt{ \mbf{s}X \ar[r]^-{\mu_{\Dl\times X}} &  \mbf{s}((\mbf{s}X)\times\Dl)\ar[r]^-{\lambda_{\mbf{s}X}} & \mbf{s}X\\
                                                                      \mbf{s}X \ar[r]^-{\mu_{X\times \Dl}} &  \mbf{s}(n\rightarrow \mbf{s}(X_n \times\Dl))\ar[r]^-{\mbf{s}(\lambda_{X_n})} & \mbf{s}X\ .}$$
F{i}rstly, consider the bisimplicial chain complex $Z=\Dl\times
X$. By def{i}nition
$$\mbf{s}((\mbf{s}X)\times\Dl)_n=\bigoplus_{p+q=n}\mbf{s}(X)_q=\bigoplus_{p+i+j=n}X_{i,j}$$
and
$({\lambda_{\mbf{s}X}})_n:\bigoplus_{p+q=n}(\mbf{s}X)_p\rightarrow
(\mbf{s}X)_n$ is the projection of $\bigoplus_{p+i+j=n}X_{i,j}$ onto $\bigoplus_{i+j=n}X_{i,j}$.\\
On the other hand, the restriction of
$({\mu_Z})_n:\ds\bigoplus_{l+k=n}X_{l,k}\rightarrow
\ds\bigoplus_{s+t+k=n}X_{t,k}$ to $X_{l,k}$ is
$$\bigoplus_{s+t=l} X(d^ld^{l-1}\cdots d^{t+1}):X_{l,k}\rightarrow X_{t,k}\ .$$
To compose with $(\lambda_{\mbf{s}Z})_n$ is the same as to project
over the components with $s=0$, that is, $t=l$ in the above
equation. But the restriction of $(\mu_Z)_n$ to these components is the identity.\\

Therefore $\lambda_{\mbf{s}X}\comp \mu_{\Dl\times X}=Id$. The case
$Z=X\times\Dl$ is completely similar.\\
F{i}nally, the associativity of $\mu$ follows from the
associativity of $\mu_{E-Z}$ (proposition \ref{AsociatMuE-Z}),
whereas the quasi-inverse of $\lambda$ is just the inclusion of
$A$ as direct summand of
$\mbf{s}(A\times\Dl)$.\\
\end{proof}

\begin{obs}
As in the proof of the commutativity up to homotopy of diagram
(\ref{compatibilidadConoconSimple}) in (SDC 7), it holds that
given any chain complex $A$, there exists a natural homotopy
equivalence between $cyl(A)$ and the classical cylinder of $A$. In
addition this equivalence is compatible with the respective
inclusions of $A$ into both cylinders.\\
Hence, the properties deduced for $\cdc$ from sections
\ref{pdadesDescenso}, \ref{SeccionFactCyl} and
\ref{SeccionCritAcCyl}, and from chapters \ref{CapituloLocalizada}
and \ref{CapituloLocalizada}, recover the classical treatment of
the homotopy category associated with $\cdc$.\\
Another consequence is that each class $\mrm{E}$ making $\cdc$
into a simplicial descent category must contain the homotopy
equivalences.
\end{obs}

Again by (\ref{compatibilidadConoconSimple}), the shift functor
induced by the descent structure coincide up to homotopy
equivalence with the usual one. Thus it is an automorphism of
categories $Ho\cdc\rightarrow Ho\cdc$. Therefore, we obtain in
this way the usual triangulated structure on $Ho\cdc$ by theorem
\ref{TmaHoDtriangulada}.


\section{Chain complexes and quasi-isomorphisms}\label{SeccionCdcQuis}

Now, let $\mc{A}$ be an abelian category with numerable sums (or
just an abelian category if we work in the uniformly
bounded-bellow case, \ref{Def{i}ComplAcotUnif}).\\

As in the additive case, consider the simple functor
$\mbf{s}=Tot\comp K:\simp\cdc\rightarrow\cdc$, and the natural
transformations $\lambda$ and $\mu$ given in def{i}nition
\ref{estructuraDescensoCDC}.

As usual, a quasi-isomorphism\index{Index}{quasi-isomorphism} is a
morphism of $\cdc$ that induces isomorphism in homology.

\begin{prop}\label{cdcabCDS}
Let $\mc{A}$ be an abelian category with numerable sums. Then the
category of chain complexes over $\mc{A}$, together with the
quasi-isomorphisms as equivalences are an additive simplicial
descent category. In addition $\lambda$ is quasi-invertible and
$\mu$ is associative.
\end{prop}

\begin{proof}
Again, the two f{i}rst axioms are well known properties. The
axioms (SDC 3), (SDC 4), (SDC 5) and (SDC 8) follow
directly from the additive case.\\
The axiom (SDC 7) is again a consequence of \ref{pdadscilindro}
and of the commutativity up to homotopy equivalence of diagram
\ref{compatibilidadConoconSimple} in the proof of (SDC 7) in proposition \ref{cdcaditivCDS}.\\
F{i}nally, to see (SDC 6) it is enough to proof that if
$X\in\simp\cdc$ is such that $X_n$ is acyclic for all $n\geq 0$,
then $\mbf{s}X$ is so, by \ref{AxiomaExactitudPrima}.\\
If $X$ is such a simplicial chain complex, then $KX$ is a double
complex located in the right-half of the plane, and whose columns
are acyclic. The following lemma state that $Tot(KX)=\mbf{s}X$ is
acyclic in this case.
\end{proof}

\begin{lema}[\cite{B}, p. 98, exercise 1.]\mbox{}\\
Let $\{X_{p,q}\ ;\ d^1:X_{p,q}\rightarrow X_{p-1,q}\ ;\
d^2:X_{p,q}\rightarrow X_{p,q-1}
\}$ be a double chain complex such that\\
\indent \textbf{1.} $X_{p,q}=0$ if $p< 0$.\\
\indent \textbf{2.} Given $p\geq 0$, the complex $X_p=\{X_{p,q}\
;\
d^2:X_{p,q}\rightarrow X_{p,q-1}\}$ is acyclic.\\
Then $Tot(X)$ is acyclic.
\end{lema}

\begin{proof}
The proof given here is an adaptation of the same fact for
``f{i}rst
quadrant'' double complexes with acyclic columns.\\
Given $q\geq 0$, consider the subcomplex $F^q\subseteq Tot(X)$
def{i}ned as
$$ (F^l)_n = \ds\bigoplus_{p+q =n\, ;\  p\leq l} X_{p,q}$$
that is in fact a subcomplex of $Tot(X)$ since
$d^{Tot(X)}(F^q)\subseteq F^q$.\\
In this way we obtain the increasing chain of subcomplexes of
$Tot(X)$
$$ 0\subseteq X_0 =F^0\subseteq \cdots \subseteq F^q \subseteq
F^{q+1}\subseteq \cdots \subseteq Tot(X).$$
Moreover, we have the short exact sequence
$$\xymatrix@H=4pt@M=4pt{0\ar[r] & F^{l}  \ar[r] & F^{l+1} \ar[r] & F^{l+1}/ F^{l} \ar[r] & 0\ .}$$
Since $F^{l+1}/ F^{l}\simeq X_{l+1}$ and $F^0=X_0$ are acyclic, it
follows by induction that $F^q$ is acyclic for every
$q\geq 0$.\\
If $[x]\in H^r(Tot(X))$ is the class of $x\in
\ker\{d_{r}:Tot(X)_r\rightarrow Tot(X)_{r-1}\}$, in particular
$x\in\bigoplus_{p+q=r}X_{p,q}$. Then, by def{i}nition of direct
sum, there exists a f{i}nite set of indexes $I$ such that $x
\in\bigoplus_{(p,q)\in I}X_{p,q}$. Therefore, we can f{i}nd $l$
such that $x\in F^l_r$. But $d^{F_l}x=d^{Tot(X)}_rx=0$ and $F^l$
is acyclic, so $x=d^{F^l}x'=d^{Tot(X)}x'$. Hence $[x]=0$ and
$Tot(X)$ is acyclic.\\
In other words, $Tot(X)$ is just the colimit of $F^0\subseteq
\cdots \subseteq F^q \subseteq F^{q+1}\subseteq \cdots$, where
each $F^q$ is acyclic, so $Tot(X)$ is also acyclic since homology
commutes with f{i}ltered colimits.
\end{proof}

As in the additive case, the shift functor induced by the descent
structure coincide up to homotopy equivalence with the usual one.
So it is an automorphism of $Ho\cdc$. Therefore,
\ref{TmaHoDtriangulada} recover the usual triangulated structure
on the derived category of $\mc{A}$.\\

\begin{obs}
In the abelian case, apart from the usual simple we can also
consider the ``normalized simple'' that is def{i}ned using the
normalized version of $K$
$$K_N:\simp\mc{B}\rightarrow Ch_\ast\mc{B}$$
instead $K$. Given a simplicial object $B$ in $\mc{B}$, $K_N(B)$
is just the quotient of  $K(B)$ over the degenerate part of $B$
\cite{May}.
Then, $K_N$ provides the normalized simple functor
$$\mbf{s}_N:\simp\cdc\rightarrow\cdc\ ,$$
that also gives rise to a simplicial descent category structure on
$\cdc$. This time $\lambda=Id$, whereas the transformation $\mu$
of the non-normalized structure pass to quotient, inducing
$\mu_N:\mbf{s}_N\comp \mrm{D}\rightarrow
\mbf{s}_N\simp\mbf{s}_N$.\\
This new descent structure is of course equivalent to the
non-normalized one, since the projection $K\rightarrow K_N$ is a
homotopy equivalence \cite{May}, and applying $Tot$ we obtain a
homotopy equivalence relating $\mbf{s}$ and $\mbf{s}_N$.\\
Consequently, the identity functor $\cdc\rightarrow\cdc$ is an
equivalence of descent categories.\\
In the normalized case, the corresponding diagram
(\ref{compatibilidadConoconSimple}) appearing in the proof of (SDC
7) is commutative. Similarly, if we compute $cyl(A)$ using this
normalized simple we obtain the usual cylinder associated with
$A$, for each $A$ in $\cdc$. In particular, the shift functor
coincides with the usual one in this case.
\end{obs}

\begin{obs}\label{Def{i}ComplAcotUnif}
In sections \ref{SeccionCdcQuis} and \ref{SeccionCdcHomt} we have
considered non bounded chain complexes, but all the properties
contained in this section remain valid for uniformly
bounded-bellow chain complexes.\\

Denote by $Ch_{q}\mathcal{A}$ the category of chain complexes
$\{A_n,d\}$ with $A_n=0$ for all $n$ smaller than the f{i}xed
bound $q\in\mathbb{Z}$. In this case, we don't need to impose the
existence of numerable sums to def{i}ne the simple functor, since
we deal now with f{i}rst-quadrant double complexes.
\end{obs}

\begin{prop}\label{ComplPositivosCDS} Let $\mc{A}$ be an additive $($resp. abelian$)$ category.
Then $Ch_{q}\mathcal{A}$ with the homotopy equivalences $($resp.
quasi-isomorphisms$)$ as equivalences and the transformations
$\lambda$ and $\mu$ given in \ref{estructuraDescensoCDC}, are an
additive simplicial descent category. In addition $\lambda$  is
quasi-invertible and $\mu$ is associative.
\end{prop}

In this case the shift functor is not an automorphism of $Ho
Ch_{q}\mathcal{A}$, so $Ho Ch_{q}\mathcal{A}$ is a suspended
category (that is, right triangulated).

\begin{obs} The case $Ch_{\mathbf{b}}\mc{A}$ of (non-uniformly) bounded-bellow chain complexes cannot
be considered directly as an example of simplicial descent
category, since the simple functor doest not preserve bounded-bellow chain complexes.\\
However, we can use the previous proposition and argue as follows
in order to give a proof based in these techniques of the well
known triangulated structure on the derived category of
$Ch_b\mc{A}$ (that is, the bounded-bellow derived category
associated with $\mc{A}$).
\end{obs}

\begin{cor}\label{ComplAcotadoesTriang} Let $Ch_{\mbf{b}}\mc{A}$ be the category of bounded-bellow chain
complexes. Then the localized category $Ho Ch_{\mbf{b}}\mc{A}$ of
$Ch_{\mbf{b}}\mc{A}$ with respect to the quasi-isomorphisms (resp.
the homotopy equivalences) is a triangulated category.
\end{cor}

\begin{proof}
Let us prove the case $Ho Ch_{\mbf{b}}\mc{A}=D_{\mbf{b}}\mc{A}$,
the localized category of $Ch_{\mbf{b}}\mc{A}$ with respect to the
quasi-isomorphisms. The other case is completely similar.\\
The idea of the proof is just to induce in
$D_{\mbf{b}}\mc{A}``=''\bigcup_{k\in\mathbb{Z}}Ho Ch_k\mc{A}$ the
suspended category structure coming from each $Ho Ch_k\mc{A}$, and
since in $D_{\mbf{b}}\mc{A}$ the shift functor is an automorphism,
it follows that $D_{\mbf{b}}\mc{A}$
is triangulated.\\
Before localizing, we have the chain of inclusions of categories
$$\cdots \subset Ch_k\mc{A}\subset Ch_{k+1}\mc{A}\subset \cdots \subset Ch_{\mbf{b}}\mc{A}$$
and $Ch_{\mbf{b}}\mc{A}=\bigcup_{k\in\mathbb{Z}}Ch_k\mc{A}$.\\
For any $k\in\mathbb{Z}$, the category $Ch_k\mc{A}$ is an additive
simplicial descent category. In particular, we have the cone
functor $c^k:Fl(Ch_k\mc{A})\rightarrow Ch_k\mc{A}$, that is
compatible with the inclusions $Ch_k\mc{A}\subset Ch_{k+1}\mc{A}$.
This compatibility holds because the simple functor does not
depend on $k$. Then, the family $\{c^k\}$ induces the cone functor
$c:Fl(Ch_{\mbf{b}}\mc{A})\rightarrow Ch_{\mbf{b}}\mc{A}$, that is
well def{i}ned.
Therefore, the shift functor $[1]:Ch_{\mbf{b}}\mc{A}\rightarrow
Ch_{\mbf{b}}\mc{A}$ is also def{i}ned. Moreover, it preserves
quasi-isomorphisms,
so it passes to the derived categories.\\
\indent Given a morphism $f:X\rightarrow Y$ in
$Ch_{\mbf{b}}\mc{A}$, there exists $K\in\mathbb{Z}$ such that $f$
is in $Ch_K\mc{A}$, so $f$ gives rise to the triangle
\begin{equation}\label{triangDistCFb}
 \xymatrix@M=4pt@H=4pt{  X  \ar[r]^f & Y \ar[r] &  c(f)\ar[r] & X[1] }
\end{equation}
where all arrows are in $Ch_K\mc{A}$ because $\lambda$ is quasi-invertible.\\
\indent Def{i}ne the class of distinguished triangles of
$D_{\mbf{b}}\mc{A}$ as those isomorphic (in $D_{\mbf{b}}\mc{A}$)
to some triangle in the form (\ref{triangDistCFb}).\\
This class of distinguished triangles is by def{i}nition closed by
isomorphism. Each distinguished triangle of $D_{\mbf{b}}\mc{A}$ is
isomorphic to some distinguished triangle of $D_K\mc{A}$, for some
$K$. Since $D_K\mc{A}$ is suspended (theorem
\ref{TmaHoDtriangulada}). Thus the distinguished triangles of
$D_{\mbf{b}}\mc{A}$
satisfy all axioms of suspended categories.\\
\indent Consequently, $D_{\mbf{b}}\mc{A}$ is suspended, and since
$[1]$ is an automorphism, it is triangulated.
\end{proof}

\begin{obs} The above proof can be generalized to the case in
which a category $\mc{D}$ is the inductive limit of a family of
simplicial descent categories, because the argument used above
just means that theorem \ref{TmaHoDtriangulada} is preserved by
``inductive limits''.
\end{obs}


%
%
\section{Simplicial objects in additive or abelian categories}

The Eilenberg-Zilber-Cartier theorem \ref{E-Z-C} admits an
interesting interpretation in the context of descent categories.
Under our setting, this theorem means that the simplicial objects
in an additive category are a simplicial descent category, taking
as simple the diagonal functor, and as equivalences the morphisms
that are mapped into homotopy equivalences in $Ch_\ast\mc{A}$.

Let $Ch_+\mc{A}$ be the category of positive chain complexes of
$\mc{A}$ (that is, $Ch_+\mc{A}=Ch_0\mc{A}$ in
\ref{Def{i}ComplAcotUnif}).

Remind that the functor $K:\simp\mathcal{A}\rightarrow Ch_+\mc{A}$
given in \ref{estructuraDescensoCDC} is def{i}ned by taking as
boundary map the alternate sum of the face maps.\\

Next def{i}nition describes the descent structure on
$\simp\mc{A}$.

\begin{def{i}}\mbox{}\\
\textbf{Simple functor:} The simple functor is the diagonal
functor $\mathrm{D}:\simp\simp\mathcal{A}\rightarrow \simp\mathcal{A}$.\\[0.2cm]
\textbf{Equivalences:} The equivalences are the class
$$\mrm{E}=\{f\in\simp\mathcal{A} \;|\; K(f)\mbox{ is a homotopy
equivalence}\}\ .$$
\textbf{Transformations $\mathbf{\lambda}$ and $\mu$:} The natural
transformations $\lambda$ and $\mu$ are def{i}ned as the
corresponding identity natural transformation.
\end{def{i}}

\begin{prop}\label{simpad}
Under the previous notations,
$(\simp\mc{A},\mrm{E},\mrm{D},\mu,\lambda)$ is an additive
simplicial descent category, in which $\lambda$ is
quasi-invertible and $\mu$ associative.\\
In addition, $K:\simp\mc{A}\rightarrow Ch_+\mc{A}$ is a functor of
additive simplicial descent categories, where we consider in
$Ch_+\mc{A}$ the descent structure given in proposition
\ref{ComplPositivosCDS}.
\end{prop}

\begin{proof}
Let us check that $K:\simp\mathcal{A}\rightarrow Ch_+\mc{A}$
satisfies the hypothesis of transfer lemma \ref{FDfuerte}.\\
The axioms (SDC 1) and (SDC 3)$'$ are clear. Let us see (SDC 4)$'$.\\
Since $\mbf{s}=\mathrm{D}$, the compositions $\mbf{s}\simp\mbf{s}$
and $\mbf{s}\mathrm{D}$ are equal, so it suf{f}{i}ces to take $\mu=Id:\mbf{s}\mathrm{D}\rightarrow \mbf{s}\simp\mbf{s}$.\\
To see (SDC 5)$'$, consider $X\in\simp\simp\mathcal{A}$. Then
$\mathrm{D}(X\times\Dl)=X$, and we can set $\lambda=Id$.\\
The compatibility between  $\lambda$ and $\mu$ holds trivially.
In addition, $K$ is an additive functor, in particular it is
quasi-strict monoidal, so (FD 1) holds.\\[0.1cm]
It remains to see (FD 2). Denote by
$(Ch_+\mc{A},\mrm{E}',\mbf{s}',\mu',\lambda')$ the descent
structure on $Ch_+\mc{A}$ given in \ref{ComplPositivosCDS}.
By theorem \ref{E-Z-C} a) for $\mc{U}=\mc{A}$, we deduce that
$\Theta=\mu_{E-Z}:K\comp \mrm{D}\rightarrow \mbf{s}'\comp \simp K$
is a homotopy equivalence when we evaluate it at each bisimplicial
object of $\mc{A}$.\\
Then, if $X\in\simp\simp\mathcal{A}$, the morphism
$\Theta_X=\mu_{E-Z}:K \mrm{D}(X)\rightarrow \mbf{s}'\simp K(X)$ is
``universal'' and such that
$(\Theta_X)_0=Id_{X_{0,0}}:X_{0,0}\rightarrow X_{0,0}$.\\
Let us check the compatibility between $\lambda$, $\mu$,
$\Theta$, $\lambda'$ and $\mu'$.\\
Given $X\in\simp\mathcal{A}$, denote by $\widetilde{X}$ the
associated constant simplicial object, that is
$\widetilde{X}_{n,m}=X_m$ for all $n$, $m$. We must see that
$\lambda'_{KX}\comp \mu_{E-Z}(\widetilde{X})=Id_{KX}$ in
$ Ch_+\mc{A}$.\\
By def{i}nition, $(\lambda'_{KX})_n:X_{n}\oplus\cdots\oplus X_0
\rightarrow X_n$ is the projection, whereas
$\mu_{E-Z}(\widetilde{X}):X_n\rightarrow X_{n}\oplus\cdots\oplus
X_0$ is $\mu_{E_Z}(\widetilde{X})=(Id,d_n,d_{n-1}\comp
d_n,\ldots,d_1\comp\cdots\comp d_n)$. Therefore, when we project
over the component $X_n$, we obtain the identity.\\[0.1cm]
The compatibility between $\Theta=\mu_{E_Z}$ and $\mu'$ (also
obtained from $\mu_{E-Z}$) is consequence of the associativity of
this transformation \ref{AsociatMuE-Z}.\\
F{i}nally, $\simp\mc{A}$ is additive since $\mc{A}$ is, and the
diagonal functor $\mrm{D}:\simp\simp\mc{A}\rightarrow\mc{A}$ is
additive, so $\simp\mc{A}$ is an additive simplicial descent
category.
\end{proof}
%


Assume now that $\mc{A}$ is an abelian category, and
$K:\simp\mc{A}\rightarrow Ch_+\mc{A}$ the usual functor. We have
the additional descent structure on $\simp\mc{A}$.
\begin{def{i}}\mbox{}\\
\textbf{Simple functor:} Again, the simple functor is the diagonal
functor $\mathrm{D}:\simp\simp\mathcal{A}\rightarrow
\simp\mathcal{A}$.\\[0.2cm]
\textbf{Equivalences:} The class of equivalences is
$$\mrm{E}'=\{f\in\simp\mathcal{A} \;|\; K(f)
\mbox{ is a quasi-isomorphism in } Ch_+\mc{A}\}\ .$$
\textbf{Transformations $\mathbf{\lambda}$ and $\mu$:} The natural
transformations $\lambda$ and $\mu$ are def{i}ned as the identity
natural transformation.
\end{def{i}}

\begin{prop}\label{simpab}
Under the above notations,
$(\simp\mc{A},\mrm{E}',\mrm{D},\mu,\lambda)$ is an additive
simplicial descent category such that $\lambda$ is
quasi-invertible and $\mu$ is associative.\\
In addition, $K:\simp\mc{A}\rightarrow Ch_+\mc{A}$ is a functor of
additive simplicial descent categories, where the descent
structure on $ Ch_+\mc{A}$ is the one given in proposition
\ref{ComplPositivosCDS}.
\end{prop}

\begin{proof} From the proof of proposition
\ref{simpad} we deduce that $K:\simp\mathcal{A}\rightarrow
Ch_+\mc{A}$ satisf{i}es the conditions of the transfer lemma
\ref{FDfuerte}, now taking the descent structure on $Ch_+\mc{A}$
in with the equivalences are the quasi-isomorphisms
\ref{ComplPositivosCDS}.
\end{proof}
%

%


%
%
\section{Simplicial Sets}

Denote by $Set$\index{Symbols}{$Set$} the category of sets, and by
$Ab$\index{Symbols}{$Ab$} the category of abelian groups.

In this section we will give a descent structure to $\simp
Set$\index{Index}{simplicial sets}, in which the equivalences will
be the quasi-isomorphisms.

\begin{def{i}}\label{HomologiaSimpSet}
If $L:Set\rightarrow Ab$ is the functor that maps a set $T$ to the
free group with base $T$, then the homology of a simplicial set
$W$ is the homology of the chain complex $K\comp\simp L (W)$, that
is the image of $W$ under the composition of functors
$$\xymatrix@M=4pt@H=4pt@C=35pt{ \simp Set \ar[r]^-{\simp L} & \simp Ab \ar[r]^-{K}& Ch_+(Ab) .}$$
\end{def{i}}
\begin{def{i}}\mbox{}\\
\textbf{Simple functor:} Again, the simple functor is the diagonal
functor $\mathrm{D}:\simp\simp {Set}\rightarrow \simp
{Set}$.\\[0.2cm]
\textbf{Equivalences:} The class $\mrm{E}$ of equivalences
consists of those morphisms that induce isomorphism in homology.\\[0.2cm]
\textbf{Transformations $\mathbf{\lambda}$ and $\mu$:} The natural
transformations $\lambda$ and $\mu$ are def{i}ned as the identity
natural transformation.
\end{def{i}}

\begin{prop}\label{sSet}
Under the above notations,
$(\simp{Set},\mrm{E},\mrm{D},\mu,\lambda)$ is a simplicial descent
category such that $\lambda$ is
quasi-invertible and $\mu$ is associative.\\
In addition, $\simp {L}:\simp{Set}\rightarrow \simp Ab$ is a
functor of simplicial descent categories, where the descent
structure on $\simp Ab$ is the one given in proposition
\ref{simpab}.
\end{prop}

\begin{proof} The compatibility between $\lambda$ and
$\mu$ is clear. The functor $\simp L:\simp Set\rightarrow \simp
Ab$ satisf{i}es the hypothesis of the transfer lemma
\ref{FDfuerte} trivially, where $\Theta$ is again the identity
natural transformation.
\end{proof}

\begin{obs} The homology of a simplicial set $W$ coincides with the singular homology of its geometric realization $|W|$ (see \cite{May} 16.2 ii)).
Then
$$\mrm{E}=\{f:W\rightarrow W'\ |\ |f|:|W|\rightarrow |L|\mbox{ induces isomorphism in singular homology}\}$$
\end{obs}

%



%
\section{Topological Spaces}

Consider the category $Top$\index{Symbols}{$Top$} of topological
spaces and continuous maps.
We will endow the category $Top$ with a descent structure in which
the equivalences are the quasi-isomorphisms (that is, morphisms
inducing isomorphism in singular homology).\\
\indent The usual geometric realization $|\cdot|:\simp
Top\rightarrow Top$ present some disadvantages for our purposes.
For instance, we need to impose some extra conditions to a map
$f:X\rightarrow Y$ such that $f_n$ is an equivalence for all $n$
in order to have that $|f|$ is again an equivalence (see, for
instance, \cite{M} 11.13). In other words, the exactness axiom of
simplicial descent categories is not satisf{i}ed by $|\cdot|$
under the generality needed here.\\
\indent This is the reason why we consider as simple the so called
``fat'' geometric realization, def{i}ned in a similar way as
$|\cdot|$, except that now we do not identify those terms related
through the degeneracy maps (we only identify those terms related
through the face maps).\\
\indent The natural transformation $\mbf{s}\rightarrow |\cdot| $
is a homotopy equivalence when evaluated at those $X\in\simp Top$
such that the degeneracy maps are closed cof{i}brations (see
\cite{S}, appendix A), for instance when evaluated at simplicial
sets.

\begin{def{i}}
Let $\bigtriangleup : \Dl\longrightarrow
Top$\index{Symbols}{$\bigtriangleup$} be the functor which maps
the ordinal $[m]$ in $\Dl$ to the standard $m$-dimensional simplex
$\bigtriangleup^m \subset \mathbb{R}^{m+1}$ given by
$$\bigtriangleup^m = \{(t_0,\ldots,t_{m})\in \mathbb{R}^{m+1}\;\mid\;\ds\sum_{k=0}^m t_k=1\mbox{ and } t_k \geq 0 \}.$$
If $f:[n]\rightarrow [m]$ is a morphism of $\Dl$, then $f$ induces
a continuous map $\bigtriangleup(f):
\bigtriangleup^n\rightarrow\bigtriangleup^m$.\\
Setting $J_i=f^{-1}(\{i\})$, then
$\bigtriangleup(f)(t_0,\ldots,t_{n})=(r_1,\ldots ,r_m)$ where
$r_i=\sum_{j\in J_i} t_j$ if $J_i$ is not empty, and $r_i=0$ otherwise.\\
Recall that the singular homology of a topological space is by
def{i}nition the homology of the simplicial set obtained through
the ``singular chains'' functor\index{Index}{singular chains}.
This functor $S:Top\rightarrow \simp Set$ assigns to a topological
space $X$ the simplicial set
$$SX=\{Hom_{Top}(\bigtriangleup^n , X)\}_n\ .$$
Then, the singular homology of $X$ is just the homology of the
chain complex $K\simp L (SX)$ given in def{i}nition
\ref{HomologiaSimpSet}.
\end{def{i}}

\begin{def{i}}\mbox{}\\
\textbf{Simple functor:} the simple functor $\mbf{s}:\simp
Top\rightarrow Top$ is the ``fat'' geometric
realization\index{Index}{fat geometric realization}. Given a
simplicial topological space
$$
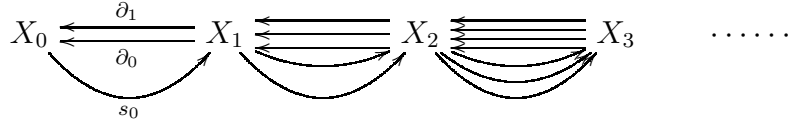
\xymatrix@1{ X_0 \; \ar@/_2pc/[rr]_{s_0} && {\;} X_1 \;
\ar@<0.5ex>[ll]^-{\partial_0} \ar@<-0.5ex>[ll]_-{\partial_1}
\ar@/_2pc/[rr] \ar@/_1pc/[rr] && \; X_2 \; \ar@<0ex>[ll]
\ar@<1ex>[ll] \ar@<-1ex>[ll] \ar@/_1.5pc/[rr] \ar@/_2pc/[rr]
\ar@/_1pc/[rr] &&   {\;} X_3
\ar@<0.33ex>[ll]\ar@<-0.33ex>[ll]\ar@<1ex>[ll]\ar@<-1ex>[ll]&
\cdots\cdots }$$
consider the bifunctor $\xymatrix@M=4pt@R=1pt{
    \Dl_e^\comp \times \Delta_e \ar[rr] && Top \\
                  ([n], [m]) \ar[rr]&& X_n \times \bigtriangleup^m }$.\\
The fat geometric realization of $X$ is def{i}ned as the cof{i}nal
of this bifunctor (cf. \cite{ML}):
$$ \mathbf{s}X= \int^n X_n \times \bigtriangleup^n $$
more specif{i}cally,\\
$$
\mathbf{s}X=  \begin{array}{cc}\ds{\coprod_{n\geq 0} X_n \times \bigtriangleup^n} &\\[-0.7cm]
                                      & \hspace{-0.5cm} \mbox{\huge{$\diagup$}} \!\! \sim \end{array}$$
where $\sim$ is the equivalence relation generated by
$$ (\partial_i(x),u)\sim (x,\bigtriangleup(d_i)(u))\mbox{ if }d_i:[n-1]\rightarrow [n],\mbox{ and } (x,u)\in X_n\times\bigtriangleup^{n-1}.$$
We will write $[x,t]$ for the equivalence class of an element
$(x,t)\in\coprod_{n\geq 0} X_n \times \bigtriangleup^n$.\\[0.2cm]
\noindent\textbf{Equivalences:} Consider the class $\mrm{E}$
consisting of those continuous maps that induce isomorphism in
singular homology. This kind of maps will be called
quasi-isomorphisms as well.\\[0.2cm]
\noindent\textbf{Transformation $\mathbf{\lambda}$:} Given $X\in
Top$, the projections $p_n:X\times\bigtriangleup^n\rightarrow X$
induce by the universal property of cof{i}nals a continuous map
$\lambda_X:\mbf{s}(X\times\Dl)\rightarrow X$, natural in
$X$, with $\lambda_X[x,t]=x$.\\[0.2cm]
\noindent\textbf{Transformation $\mathbf{\mu}$:} If
$Z\in\simp\simp Top$, since the simple functor is def{i}ned as a
cof{i}nal, the Fubini theorem holds (\cite{ML},IX.8), and we
deduce that $\mbf{s}\,\simp\mbf{s}Z$ is the quotient of
$\coprod_{n,m\geq
0}Z_{n,m}\times\bigtriangleup^n\times\bigtriangleup^m$ over the
obvious identif{i}cations.\\
Then, the maps $(\mu_Z)_n:Z_{n,n}\times\bigtriangleup^n\rightarrow
\mbf{s}\,\simp\mbf{s}Z$ with $(\mu_Z)_n[z,t]=[z,t,t]$ provides a
continuous map $\mu_Z:\mbf{s}\mathrm{D}Z\rightarrow
\mbf{s}\,\simp\mbf{s}Z$ such that
$\mu_Z([z_{nn},t_n])=[z_{n,n},t_n,t_n]$.
\end{def{i}}

\begin{prop}\label{esptopologic}
Under the previous notations, $(Top,\mrm{E},\mbf{s},\mu,\lambda)$
is a simplicial descent category, such that $\mu$ is associative
and $\lambda$ is quasi-invertible.
In addition, the singular functor $S:Top\rightarrow \simp Set$ is
a functor of simplicial descent categories.
\end{prop}

Now we will begin with the proof of this proposition. To this end
we need some preliminary results.

\begin{lema}\label{exactitudgruesa}
If $f:X\rightarrow Y$ is a morphism in $\simp Top$ such that for
all $n$, $f_n$ induces isomorphism in singular homology $($resp.
$f_n$ is a homotopy equivalence$)$ then the same holds for
$\mbf{s}f:\mbf{s}X\rightarrow \mbf{s}Y$.
\end{lema}

The proof for $f_n$ quasi-isomorphism can be found in \cite{Dup}
5.16, and the case $f_n$ homotopy equivalence appears in \cite{S}
A.1.

\begin{obs}
The previous lemma justif{i}es the choice of $\mbf{s}$ as simple
functor instead of the usual geometric realization $|\cdot|$. On
the other hand, one of the advantages of $|\cdot|$ is the
existence of the adjoint pair (\cite{May}, $\S$16)
\begin{equation}\label{adjuncionRealGeom}\xymatrix@1@C=25pt@H=4pt@M=4pt{ \simp Set
\ar@<0.5ex>[r]^-{|\cdot|} &  Top  \ar@<0.5ex>[l]^-{S} \
.}\end{equation}
Our simple functor $\mbf{s}$ is def{i}ned by forgetting the
degeneracy maps of a simplicial topological space. A consequence
of this fact is that the above adjunction does not hold at the
level of simplicial sets between $\mbf{s}$ and $S$, but it holds
at the
level of strict simplicial sets. That is to say, there is an adjunction\\
$$\xymatrix@1@C=25pt@H=4pt@M=4pt{ \simp_e Set \ar@<0.5ex>[r]^-{\mbf{s}} &  Top  \ar@<0.5ex>[l]^-{S} \ .}$$
Due to this fact we will have to solve some technical
dif{f}{i}culties in the proof of proposition \ref{esptopologic}.
\end{obs}

\begin{obs}
Note that we can consider as well the homology of a strict
simplicial set $W\in\simp_e Set$, because the functor $K$
does not use the degeneracy maps of a simplicial set, that is, $K:\simp_e Ab\rightarrow Ch_+(Ab)$.\\
Then, quasi-isomorphisms between strict simplicial sets can be
def{i}ned in the same way as those morphisms $f:W\rightarrow W'$
in $\simp_e Set$ that induce isomorphisms in homology.
\end{obs}

\begin{def{i}} The geometric realization $|\cdot |:\simp Top\rightarrow
Top$ is def{i}ned as
$$
|X|=  \begin{array}{cc}\ds{\coprod_{n\geq 0} X_n \times \bigtriangleup^n} &\\[-0.7cm]
                                      & \hspace{-0.5cm} \mbox{\huge{$\diagup$}} \!\! \sim \end{array}$$
where $\sim$ is the equivalence relation generated by
$$ (\theta(x),u)\sim (x,\bigtriangleup(\theta)(u))\mbox{ if
}\theta:[m]\rightarrow [n]\mbox{ is a morphisms of }\Dl,\mbox{ and
} (x,u)\in X_n\times\bigtriangleup^{m}.$$
We will write $|x,t|$ for the equivalence class of the element
$(x,t)\in\coprod_{n\geq 0} X_n \times \bigtriangleup^n$.
\end{def{i}}

\begin{def{i}} A simplicial topological space is called ``good'' if its degeneracy maps are closed cof{i}brations.\\
Simplicial sets (with the discrete topology) are always good in
this sense.\\
Moreover, if $T\in\simp\simp Set$ is a bisimplicial set, then the
simplicial topological space $X$ with $X_n=|m\rightarrow
T_{m,n}|$, obtained by applying the geometric realization to $T$
with respect to one of its indexes is also an example of ``good''
simplicial topological space.\\
This is consequence of the fact that any degeneracy map
$s_j:T_{\cdot,n}\rightarrow T_{\cdot,n+1}$ is an inclusion of
simplicial sets (because of the simplicial identities), and the
geometric realization of any inclusion of simplicial sets is
always a closed cof{i}bration.
\end{def{i}}

Next we recall the following connection between $|\cdot|$ and
$\mbf{s}$, given in \cite{S}, appendix A.

\begin{lema}\label{RelsimplRealGeom}
If $X$ is a good simplicial topological space then the morphism
$$ \tau_X: \mbf{s}X \rightarrow |X|\ \ \ [x,t]\rightarrow |x,t|$$
is a homotopy equivalence.
\end{lema}

The fat geometric realization satisf{i}es as well the classical
Eilenberg-Zilber property.

\begin{lema}[Eilenberg-Zilber]\label{E-Z}
Given $W\in \simp\simp Set$, the map
$\eta(W):|\mrm{D}(W)|\rightarrow |\simp|W||$, with
$\eta(W)([w_{nn},t_n])=|w_{n,n},t_n,t_n|$ is an homeomorphism.\\
The map $\mu(W):\mbf{s}(\mrm{D}(W))\rightarrow
\mbf{s}\simp\mbf{s}(W)$, with
$\mu(W)([w_{nn},t_n])=[w_{n,n},t_n,t_n]$ is a homotopy
equivalence. In addition, the following diagram commutes
\begin{equation}\label{DiagE-Z}\xymatrix@M=4pt@H=4pt{ \mbf{s}\mrm{D}(W)\ar[d]^{\tau} \ar[r]^-{\mu(W)} & \mbf{s}\simp\mbf{s}(W)\ar[d]^{P} \\
                           |\mrm{D}(W)| \ar[r]^-{\eta(W)} &  |\simp|W|| \ ,}\end{equation}
where $P:\mbf{s}\simp\mbf{s}(W)\rightarrow |\simp|W||$,
$[x,p,q]\rightarrow |x,p,q|$.
\end{lema}

\begin{proof}
F{i}rstly, the usual Eilenberg-Zilber theorem (\cite{GM}. I.3.7)
states that $\eta(W):|\mrm{D}(W)|\rightarrow |\simp|W||$ is a
homoeomorphism. The second part of the lemma is a consequence of
the commutativity of diagram (\ref{DiagE-Z}), since both $\tau$
(lemma \ref{RelsimplRealGeom}) and $P$
are homotopy equivalences.\\
F{i}x $n$ as the f{i}rst index of $W$, obtaining
$W_{n,\cdot}\in\simp Set$. The projections $p_n=\tau_{W_n} :
\mbf{s}(W_{n,\cdot})\rightarrow |W_{n,\cdot}|$ are homotopy equivalences for all $n$ by \ref{RelsimplRealGeom}.\\
Hence, $p_\cdot=\{p_n\}_n:\simp\mbf{s}(W)\rightarrow\simp|W|$ is
such that $\mbf{s}(p_\cdot):\mbf{s}\simp\mbf{s}(W)\rightarrow
\mbf{s}\simp|W|$ is a homotopy equivalence by \ref{exactitudgruesa}.\\
Since $\simp|W|$ is a good simplicial topological space then the
projection $\mbf{s}\simp|W|\rightarrow |\simp|W||$ is also a
homotopy equivalence, and composing both maps we deduce that the
projection $P:\mbf{s}\simp\mbf{s}(W)\rightarrow |\simp|W||$ is a
homology equivalence as required.
\end{proof}

The following statement is similar to the classical result
satisf{i}ed by the usual geometric realization $|\cdot|$.

\begin{lema}\label{adjuncion}
Consider the natural transformations $\Phi:\mbf{s}\,\mrm{S}
\Longrightarrow Id_{Top}$ and $\Psi: Id_{\simp_e Set}
\Longrightarrow \mrm{S}\, \mbf{s}$ def{i}ned as
$\Phi(Y)([\lambda_n ,t_n])=\lambda_n (t_n)$, $(\Psi(W)_n(w_n))
(t_n)=[w_n,t_n]$. Then $\Phi$ and $\Psi$ induce isomorphism in
homology for all $Y\in Top$ and for all $W\in \simp_e Set$.
\end{lema}
\begin{proof}
As stated in the proof of (\cite{Dup}, 5.15), given $W\in \simp_e
Set$, then the fat geometric realization of $W$, $\mbf{s}W$, is a
CW-complex whose $n$-cell are the set $W_n$. Then, the group of
n-cellular chains of $\mbf{s}W$ is just $L(W_n)$ and the cellular
boundary map is $\sum_i (-1)^i d_i$.\\
Since the cellular homology and the singular homology of a
CW-complex coincide, then the morphism $\Psi(W)$ of $\simp_e Set$
induces isomorphism in homology.

Consider now $Y\in Top$. By the f{i}rst part of the lemma,
$\Psi(\mathrm{S}Y):\mathrm{S}Y\rightarrow
\mathrm{S}\mbf{s}(\mathrm{S}Y)$ induces isomorphism in homology.
But $\mathrm{S}\Phi(Y)\comp \Psi \mrm{S}(Y)=Id:
\mathrm{S}Y\rightarrow \mathrm{S}Y$, so $\mathrm{S}\Phi(Y)$ is
also a quasi-isomorphism. Hence $\Phi(Y)$ induces isomorphism in
singular homology.
\end{proof}

We need use the next technical result.

\begin{lema}\label{LemaF{i}el}
Let $Ho \simp Set$ $($resp. $Ho\simp_e Set)$ be the localized
category of simplicial sets $($resp. strict simplicial sets$)$
with respect to the quasi-isomorphisms. Then, the forgetful
functor $\mrm{U}:\simp Set\rightarrow \simp_e Set$ preserves
quasi-isomorphisms, giving rise to the functor
$$\mrm{U}: Ho \simp Set\rightarrow Ho \simp_e Set\ .$$
This is a faithful functor, that is, the map $ Hom_{Ho \simp
Set}(W,L) \rightarrow Hom_{Ho \simp_e Set}(\mrm{U}W,\mrm{U}L)$ is
injective.
\end{lema}

\noindent We will delay the proof of this lemma to give before the
one of proposition \ref{esptopologic}.

\begin{proof}[\begin{large}\textbf{Proof of \ref{esptopologic}}\end{large}]\mbox{}\\
Let us check that the hypothesis of proposition \ref{FDfuerte} are
satisf{i}ed by the singular chain complex functor $S:Top\rightarrow \simp Set$.\\
First, note that the transformations $\lambda$ and $\mu$ are
compatible. To see this, let $X$ be a simplicial topological
space, and $[x,t]$ an element of $\mbf{s}X$
representing the pair $(x,t)\in X_n \times \bigtriangleup^n$.\\
Then $\lambda_{\mbf{s}X}\comp \mu_{\Dl\times
X}([x,t])=\lambda_{\mbf{s}X}([[x,t],t])=[x,t]$, so
$\lambda_{\mbf{s}X}\comp \mu_{\Dl\times X}=Id$, and similarly
$\mbf{s}(\lambda_{X_n})\comp \mu_{X\times \Dl}=Id$.\\
\indent The disjoint union is the coproduct in $Top$, and the
singular chain functor commutes with coproducts, so (FD 1) holds.\\
\indent %
In order to relax the notations, we will write also $\psi$ for the
induced functors $\simp\psi:\simp\mc{D}\rightarrow \simp\mc{D}'$
and
$\simp\simp\psi:\simp\simp\mc{D}\rightarrow \simp\simp\mc{D}'$.\\
To see (FD 2), we must study the commutativity of the diagram
$$\xymatrix@M=4pt@H=4pt{ \simp Top \ar[r]^-{S}\ar[d]^{\mbf{s}} &  \simp\simp Set \ar[d]^{\mrm{D}}   \\
                               Top \ar[r]^-{S}                 &  \simp Set     \ .}$$
Def{i}ne the isomorphism $\Theta_X:S(\mbf{s}X)\rightarrow
\mrm{D}(SX)$ of $Ho\simp Set$ as the one coming from the zig-zag
in $\simp Set$
$$\xymatrix@M=4pt@H=4pt{S(\mbf{s}X) & S(\mbf{s}\simp\mbf{s}(SX)) \ar[l]_-{\Theta^0_X} \ar[r]^-{\Theta^1_X} & S|\simp|SX||  & \mrm{D}(SX)\ar[l]_-{\Theta^2_X}}$$
def{i}ned as follows.\\[0.1cm]
\indent The transformation
$\Theta^0_X:S(\mbf{s}\simp\mbf{s}(SX))\rightarrow S(\mbf{s}X)$ is
just the image under $S\mbf{s}$ of the morphism
$\phi:\simp\mbf{s}(SX)\rightarrow X$ of $\simp Top$, which in
degree $n$ is given by $\Phi(X_n):\mbf{s}SX_n\rightarrow X_n$ (see
\ref{adjuncion}).
Therefore, if $\alpha:\triangle^n\rightarrow
\mbf{s}\simp\mbf{s}(SX)$ is the morphism that assigns to
$t\in\triangle^n$ the class of $(\beta,p,q)\in
S_kX_m\times\triangle^k\times\triangle^m$, it follows that
$$\Theta^0_X(\alpha)(t)= \mbf{s}\phi \comp\alpha (t)= [\beta(p),q]\in \mbf{s}X\ .$$
Then, as each $\Phi(X_n)$ is a quasi-isomorphism, we deduce from
\ref{exactitudgruesa} that $\mbf{s}\phi$ (and hence
$S\mbf{s}\phi=\Theta^0_X$) is so.\\[0.1cm]
\indent The transformation
$\Theta^1_X:S(\mbf{s}\simp\mbf{s}(SX))\rightarrow S|\simp|SX||$ is
the image under $S$ of the projection
$P:\mbf{s}\simp\mbf{s}(SX)\rightarrow S|\simp|SX||$,
$P([x,t,r])=|x,t,r|$. In \ref{E-Z} we checked that $P$ is a
homotopy equivalence for any bisimplicial set $W$, in particular it is so for $W=SX$.\\[0.1cm]
\indent Secondly, $\Theta^2_X:\mrm{D}(SX)\rightarrow S|\simp|SX||$
assigns to $\alpha:\triangle^n\rightarrow X_n$ the morphism
$$\Theta^2_X(\alpha):\triangle^n\rightarrow |\simp|SX||\mbox{ given by }t\rightarrow |\alpha,t,t|, \mbox{ where }(\alpha,t,t)\in S_nX_n\times\triangle^n\times\triangle^n\ .$$
Note that $\Theta^2_X$ induces isomorphism in homology, since it
is just the composition
$$\xymatrix@M=4pt@H=4pt@C=50pt@R=5pt{ \mrm{D}(SX) \ar[r]^{\Psi'(SX)}                                   & S|\mrm{D}(SX)| \ar[r]^{S(\eta(SX))}                                                & S|\simp|SX||\\
                        \mbox{\begin{small}$\alpha:\bigtriangleup^n\rightarrow X_n$\end{small}} \ar[r] & \mbox{\begin{small}$t\in\bigtriangleup^n\rightarrow |\alpha,t|$\end{small}} \ar[r] & \mbox{\begin{small}$t\in\bigtriangleup^n\rightarrow |\alpha,t,t|$\end{small}} } $$
The morphism $\Psi'(SX)$ comes from the adjunction
(\ref{adjuncionRealGeom}) and it is a quasi-isomorphism (cf.
\cite{May} 16.2), whereas $\eta(SX)$ is a homeomorphism by
\ref{E-Z}, so $S(\eta(SX))$ is an isomorphism of simplicial sets.\\[0.1cm]
\indent It remains to prove that $\Theta$ is compatible with the
natural transformations $\lambda$ and $\mu$ of $Top$ and $\simp
Set$.
To see this, by lemma \ref{LemaF{i}el} it suf{f}{i}ces to check
that the corresponding diagrams commute in $Ho\simp_e Set$. The
advantage of wording in $\simp_e Set$ is that $\mrm{U}\Theta$
coincides in $Ho\simp_e Set$ with the morphism of strict
simplicial sets
$$\xymatrix@M=4pt@H=4pt{S(\mbf{s}X) & \mrm{D}(SX)\ar[l]_-{\theta(X)}}$$
such that $\theta(X)_n : S_n X_n= \{\gamma
:\bigtriangleup^n\rightarrow X_n\} \rightarrow S
\mbf{s}(X)=\{\zeta:\bigtriangleup^n\rightarrow \mbf{s}(X)\}$ is
$\theta(X)_n (\gamma)(t)= [\gamma(t),t]\in \mbf{s}(X)$.\\
Indeed, the morphism $\theta'(X)$ of $\simp_e Set$ def{i}ned as
$$\xymatrix@M=4pt@H=4pt@C=50pt@R=5pt{ \mrm{D}(SX) \ar[r]^-{\theta'(X)}                                   & S(\mbf{s}\simp\mbf{s}(SX)) \\
    \mbox{\begin{small}$\alpha:\bigtriangleup^n\rightarrow X_n$\end{small}} \ar[r] & \mbox{\begin{small}$t\in\bigtriangleup^n\rightarrow [\alpha,t,t]$\end{small}}  } $$
f{i}ts into the commutative diagram
$$\xymatrix@M=4pt@H=4pt@R=15pt@C=15pt{
 S(\mbf{s}\simp\mbf{s}(SX)) \ar[rd]_{\mrm{U}\Theta^1_X} &              & \mrm{D}(SX) \ar[ld]^-{\mrm{U}\Theta^2_X}\ar[ll]_-{\theta'(X)}\\
                                                        & S|\simp|SX|| &\ .}$$
Hence $\mrm{U}\Theta=(\mrm{U}\Theta^1_X \comp \theta'(X))^{-1}=
\theta(X)^{-1}:S(\mbf{s}X)\rightarrow \mrm{D}(SX)$.\\[0.1cm]
F{i}rstly, given $Y\in Top$, diagram (\ref{compatFDii}) in
$Ho\simp_e Set$ is now
$$\xymatrix@M=4pt@H=4pt@C=40pt@R=12pt{   {S(\mbf{s}(Y\times\Dl))} \ar[rd]^-{S(\lambda_Y)}          & \\
                                                                                                    &     SY\ .\\
                                 \mrm{D}((SY)\times\Dl) \ar[ru]_-{Id} \ar[uu]^-{\theta(Y\times\Dl)} &              }$$
whose commutativity in $\simp_e Set$ follows directly from the def{i}nitions.\\
On the other hand, given $Z\in\simp\simp Top$, (\ref{compatFDii})
is now the following diagram of $Ho\simp_e Set$
$$\xymatrix@M=4pt@H=4pt@C=40pt@R=20pt{
 \mrm{D}\mrm{D}(S Z)\ar[dd]_{\theta(\mrm{D}Z)}   \ar[r]^-{Id} & \mrm{D}\mrm{D}(S Z) \ar[d]^{\mrm{D}\theta(Z)}  \\
                                                              & \mrm{D}S(\simp\mbf{s}Z) \ar[d]^{\theta(\simp \mbf{s}Z)} \\
 S\mbf{s}(\mrm{D}Z)  \ar[r]^{S\mu_{Z}}                        & S\mbf{s}(\simp\mbf{s}Z)  }
$$
whose commutativity is again a direct consequence of def{i}nitions.\\
Therefore, the transfer lemma is proven for $S:Top\rightarrow \simp Set$.\\
To f{i}nish the proof, $\mu$ is clearly associative, whereas the
quasi-inverse of $\lambda$, $\lambda':Id_{Top}\rightarrow
\mbf{s}(-\times\Dl)$, can be def{i}ned as follows. If $X\in Top$
and $x\in X$, then $\lambda'_X(x)$ is the equivalence class in
$\mbf{s}(X\times\Dl)$ of the pair $(x,\ast)\in
X\times\bigtriangleup^0$.
\end{proof}

\vspace{0.3cm}

\begin{proof}[\textbf{Proof of lemma \ref{LemaF{i}el}}]
The proof is based in the properties satisf{i}ed by the adjoint
pair
$$\xymatrix@1@C=25pt@H=4pt@M=4pt{ \simp_e Set \ar@<0.5ex>[r]^-{\pi} &  \simp Set  \ar@<0.5ex>[l]^-{\mrm{U}} \ .}$$
where $\pi:\simp _e Set\rightarrow \simp Set$ is the
Dold-Puppe transform (see \ref{Dold-Puppe}).\\[0.2cm]
\textit{Step 1}: Denote also by $\pi:\simp _e Ab\rightarrow \simp
Ab$ to the Dold-Puppe transform in the category of abelian groups.
We will use the functors $K:\simp_e Ab\rightarrow Ch_+ Ab$,
$K:\simp Ab\rightarrow Ch_+ Ab$ and $K^{N}:\simp Ab\rightarrow
Ch_+ Ab$, where $K$ is as usual the functor that takes the
alternate sums of face maps as boundary map, whereas
$$K^N(A)=K(A)/ D(A),\mbox{ where }D(A)_n=\bigcup_{i=0}^{n-1}s_j A_{n-1}\ .$$
Given any $W\in\simp_e Ab$, we will see that
$$ K^N(\pi A) = K(A) \ .$$

Indeed, if $n\geq 0$, $K(\pi
A)_n=\coprod_{\theta:[n]\twoheadrightarrow [m]} A_{m}^{\theta}$
and it is enough to check that
$$D(\pi A)_n=\coprod_{\theta:[n]\twoheadrightarrow [m], \theta\neq Id} A_{m}^{\theta}\ .$$
By def{i}nition, the restriction of $s_i:(\pi A)_{n-1}\rightarrow
(\pi A)_n$ to the component $A^\sigma_{m}$ corresponding to
$\sigma:[n-1]\twoheadrightarrow [m]$ is just
$$s_i|_{A^\sigma_{m}}=Id:A^{\sigma}_{m}\rightarrow A_m^{\sigma\comp s_i}\ .$$
Then, we deduce that $D(\pi A)_n\subseteq
\coprod_{\theta:[n]\twoheadrightarrow [m], \theta\neq Id} A_{m}^{\theta}$.\\
The other inclusion is also clear since if
$\theta:[n]\twoheadrightarrow [m]$ is a non-identity surjection,
then there exists $0\leq i\leq n-1$ and
$\widetilde{\theta}:[n-1]\twoheadrightarrow [m]$ such that
$\theta=\widetilde{\theta}\comp s_i$. To see this, just take $i$
such that $\theta(i)=\theta(i+1)$ and def{i}ne $\theta'$ in the
natural way.\\[0.2cm]
\textit{Step 2}: It holds that $\pi:\simp _e Set\rightarrow \simp Set$ preserves quasi-isomorphisms.\\
Indeed, if $f:W\rightarrow W'$ is a quasi-isomorphism, this means
that $K(L f)$ is so in $Ch_+ Ab$, where $L:\simp_e Set\rightarrow
\simp_e Ab$ (or $L:\simp Set\rightarrow \simp
Ab$) is def{i}ned just by taking free groups. Note that $\pi\comp L =L\comp \pi$.\\
From the previous step it follows that $K^N(L \pi f)=K^N(\pi L
f)=K(Lf)$ is also a quasi-isomorphism. Since $K^N$ and $K$ are
homotopic functors (\cite{May} 22.3), we deduce that $\pi f$ is a quasi-isomorphism as well.\\[0.2cm]
\textit{Step 3}: Given $X\in \simp Set$, the morphism
$a_X:X\rightarrow \pi\mrm{U}X$ coming from the adjunction
$(\pi,\mrm{U})$,
is a quasi-isomorphism.\\
The morphism $a_X$ is in degree $n$ the inclusion $X_n\rightarrow
X_n^{Id}\sqcup\coprod_{\theta:[n]\twoheadrightarrow [m],
\theta\neq Id} X_{m}^{\theta}$. Then by step 2 we get that $K^N(L
a_X):K^N(LX)\rightarrow K^N(\pi LX)=K(LX)$ coincides with the
inclusion $K^N(LX)_n\rightarrow K(LX)_n= K^N(LX)_n\oplus D(LX)_n$,
that again by loc. cit. is a homotopy equivalence.\\[0.2cm]
\textit{Step 4}: The functor $\mrm{U}: Ho \simp Set\rightarrow Ho \simp_e Set$ is faithful.\\
Let $f,g:X\rightarrow Y$ be morphisms in $Ho\simp Set$ such that
$\mrm{U}f=\mrm{U}g$ in $Ho\simp_e Set$. By step 2, $\pi$ pass to
the localized categories, $\pi: Ho \simp_e Set\rightarrow
Ho \simp Set$.\\
Then $\pi \mrm{U}f=\pi\mrm{U}g$ in $Ho\simp Set$. On the other
hand, it follows from the functoriality of $a$ that $\pi
\mrm{U}f\comp a_X=a_Y\comp f$, and $\pi \mrm{U}g\comp a_X=a_Y\comp
g$. From step 3 we deduce that $a_Y$ is an isomorphism in $Ho\simp
Set$, so $f=g$ in $Ho\simp Set$.
\end{proof}

To f{i}nish this section we give the following consequence of the
properties previously developed of functors $|\cdot|$ and
$\mbf{s}$.

\begin{cor} The geometric realization $|\cdot|:\simp Set\rightarrow Top$ and the fat geometric realization $\mbf{s}:\simp Set\rightarrow
Top$ are functors of simplicial descent categories.
\end{cor}

\begin{proof}
F{i}rstly, let us begin with $\mbf{s}:\simp Set\rightarrow Top$.
The Eilenberg-Zilber theorem \ref{E-Z} provides the
quasi-isomorphism $\Theta=\mu_W:\mbf{s}\mrm{D}W\rightarrow
\mbf{s}\simp\mbf{s}W$ for any bisimplicial set $W$. The
compatibility between $\Theta=\mu$ and the transformations
$\lambda$ of $\simp Set$ and $Top$ follows from the compatibility
between those transformations $\lambda$ and
$\mu$ of $Top$.\\
On the other hand, the compatibility between $\Theta$ and the
respective transformations $\mu$ can be deduced from the
associativity
property of the transformation $\mu$ of $Top$.\\
To see that $|\cdot|:\simp Set\rightarrow Top$ is a functor of
simplicial descent categories, we can take this time $\Theta$ as
the zig-zag whose value at a bisimplicial set $W$ is
$$\xymatrix@M=4pt@H=4pt{|\mrm{D}W|\ar[r]^-{\eta_W} & |\simp|W|| & \mbf{s}(\simp|W|) \ar[l]_-{\tau_{\simp|W|}} \ ,}$$
where $\eta_W$ is the homeomorphism given in \ref{E-Z} and
$\tau_{\simp|W|}$ is the homotopy equivalence from
\ref{RelsimplRealGeom} (indeed $\simp|W|$ is a good simplicial
topological space since it is obtained from bisimplicial set).\\
Consider now a simplicial set $T$. The compatibility between
$\Theta_T$ and the respective transformations $\lambda_T$ is just
the commutativity of the diagram
$$\xymatrix@M=4pt@H=4pt@C=40pt{
   |\mrm{D}(\Dl\times T)|\ar[r]^-{\eta_{\Dl\times T}}\ar[rd]_{Id} & ||T|\times\Dl| \ar[d]_{\overline{\lambda}_T} & \mbf{s}(|T|\times\Dl) \ar[l]_-{\tau_{|T|\times\Dl}} \ar[ld]^{\lambda_{|T|}} \\
                                                                  & | T |                                              & }$$
where, if $(t,p,q)\in T_n\times \triangle^n\times\triangle^m$ then $\overline{\lambda}_T(|t,p,q|)=|t,p|$.\\
On the other hand, given $Z\in\simp\simp(\simp Set)$, under the
notations of def{i}nition \ref{Def{i}asocitividadMu}, we must
check the commutativity  in $Ho Top$ of the diagram
$$\xymatrix@M=4pt@H=4pt@C=27pt{
  \mbf{s}\mrm{D}_{1,2}\simp\simp|Z|\ar[rr]^-{\tau_{\simp|\mrm{D}_{1,2}Z |}} \ar[d]^{\mu_{\simp\simp|Z|}} &                             & |\simp|\mrm{D}_{1,2}Z||                                               &       & |\mrm{D}\mrm{D}_{1,2}Z |  \ar[d]_{Id} \ar[ll]_-{\eta_{\mrm{D}_{1,2}Z}} \\
  \mbf{s}\simp\mbf{s}\simp\simp|Z| \ar[r]^-{\mbf{s}\simp\tau}                                         & \mbf{s}\simp|\simp\simp|Z|| & \mbf{s}\simp|\mrm{D}_{2,3}Z| \ar[l]_-{\mbf{s}\simp\eta}\ar[r]^-{\tau} & |\simp|\mrm{D}_{2,3}Z|| & |\mrm{D}\mrm{D}_{2,3}Z| \ar[l]_-{\eta} }
$$
By the commutativity of the square
$$\xymatrix@M=4pt@H=4pt@C=40pt{
 \mbf{s}\simp|\mrm{D}_{2,3}Z| \ar[r]^-{\mbf{s}\simp\eta}\ar[d]_{\tau} & \mbf{s}\simp|\simp\simp|Z|| \ar[d]_{\tau}\\
 |\simp|\mrm{D}_{2,3}Z|| \ar[r]^-{|\simp \eta|}                       & |\simp|\simp\simp|Z|||\ ,}
$$
we can just see the commutativity of
$$\xymatrix@M=4pt@H=4pt@C=30pt{
  \mbf{s}\mrm{D}_{1,2}\simp\simp|Z|\ar[rr]^-{\tau_{\simp|\mrm{D}_{1,2}Z |}} \ar[d]^{\mu_{\simp\simp|Z|}} &  & |\simp|\mrm{D}_{1,2}Z|| &  & |\mrm{D}\mrm{D}_{1,2}Z |  \ar[d]_{Id} \ar[ll]_-{\eta_{\mrm{D}_{1,2}Z}} \\
  \mbf{s}\simp\mbf{s}\simp\simp|Z| \ar[rr]^-{\tau\comp \mbf{s}\simp\tau}                                 &  & |\simp|\simp\simp|Z|||  &  & |\mrm{D}\mrm{D}_{2,3}Z| \ar[ll]_-{|\simp \eta|\comp\eta} \ ,}
$$
One can see this commutativity in $Ho Top$ dividing it into two
squares through the map $|\simp\eta| :
|\simp|\mrm{D}_{1,2}Z||\longrightarrow |\simp|\simp\simp|Z|||$,
and it follows from the def{i}nitions that these two squares
commute in $Top$.
\end{proof}

%
%
%
%

\chapter{Examples of Cosimplicial Descent Categories}

\section{Cochain complexes}\label{Cocadenas}

If $\mc{A}$ is an additive category, the category
$Ch^\ast\mc{A}$\index{Symbols}{$Ch^\ast\mc{A}$} of cochain
complexes can be identif{i}ed with\index{Index}{cochain
complexes} $(Ch_\ast(\mc{A}^\comp))^\comp$.\\
Assume moreover that $\mc{A}$ has numerable products, that is,
given a family $\{A_k\}_{k\in\mathbb{Z}}$ of objects of
$\mc{A}$, then $\prod_{k\in\mathbb{Z}}A_k$ exists in $\mc{A}$.\\
In this case $\mc{A}^\comp$ is an additive category with numerable
products, so $Ch_\ast(\mc{A}^\comp)$ is a simplicial descent
category with respect to the homotopy equivalences.\\
We can argue analogously if $\mc{A}$ is abelian with respect to
the quasi-isomorphisms.

Again, we can drop the condition of existence of numerable
products
in the case of uniformly bounded-bellow cochain complexes.\\

Next we introduce the dual constructions to those given in
\ref{estructuraDescensoCDC}.

\begin{def{i}}\label{estructuraDescensoCDCos}\mbox{}\\
\textbf{Simple functor:} The simple functor $\mathbf{s}:
\Dl\cdco\longrightarrow \cdco$ is def{i}ned as the composition
$$ \xymatrix@M=4pt{ \Dl\cdco \ar[r]^-{K} & Ch^\ast\cdco \ar[r]^-{Tot} & \cdco}$$
where $K(\{X,d^i,s^j\})=\{X,\sum (-1)^i d^i\}$.

More concretely, let $X=\{X^n , d^i , s^j \}$ be a cosimplicial
cochain complex. Each $X^n$ is an object of $\cdco$, that will be
written as $\{ X^{n,p} , d^{X^n} \}_{p\in\mathbb{Z}}$.

Then $X$ induces the double cochain complex (\ref{cdcos}), with
vertical boundary map $d^{X^n}: X^{n,p}\rightarrow X^{n,p+1}$ and
horizontal boundary map $\partial: X^{n,p}\rightarrow X^{n+1,p}$,
$\partial=\sum_{i=0}^{n+1}(-1)^{i} d^i$.
\begin{equation}\label{cdcos}
\xymatrix{
{}            & \vdots                                            & \vdots                                         & \vdots       & {}    \\
\ldots \ar[r] & X^{n-1,p+1}\ar[u]               \ar[r]^{\partial} & X^{n,p+1}  \ar[u]             \ar[r]^{\partial}& X^{n+1,p+1} \ar[u]              \ar[r]  & \ldots \\
\ldots \ar[r] & X^{n-1,p}  \ar[u]^{d^{X^{n-1}}} \ar[r]^{\partial} & X^{n,p}    \ar[u]^{d^{X^{n}}} \ar[r]^{\partial}& X^{n+1,p}   \ar[u]^{d^{X^{n+1}}}\ar[r]  & \ldots \\
\ldots \ar[r] & X^{n-1,p-1}\ar[u]^{d^{X^{n-1}}} \ar[r]^{\partial} & X^{n,p-1}  \ar[u]^{d^{X^{n}}} \ar[r]^{\partial}& X^{n+1,p-1} \ar[u]^{d^{X^{n+1}}}\ar[r]  & \ldots \\
   {}         & \vdots     \ar[u]                                 & \vdots     \ar[u]                              &  \vdots     \ar[u]                  & {}
   }\end{equation}

Hence, the image of $X$ under the simple functor is the cochain
complex $\mathbf{s}X$ given by
$$\begin{array}{llr} (\mathbf{s}X)^q=\ds\prod_{p+n=q} X^{n,p} & &
 d= \prod(-1)^p\partial + d^{X^n}:\ds\prod_{p+n=q} X^{n,p}\longrightarrow\ds\prod_{p+n=q+1} X^{n,p} \, .\end{array} $$
\textbf{Transformation $\mathbf{\lambda}$:} If $A\in\cdco$,
$\mbf{s}(A\times\Dl)$ is in degree $n$ the product $\prod_{k\leq
n}A^k$, in such a way that the inclusion $\lambda_A:A\rightarrow
\mbf{s}(A\times\Dl)$ is a morphism of cochain complexes.\\[0.2cm]
\textbf{Transformation $\mathbf{\mu}$:} If
$Z\in\Delta\Delta\cdco$,
$\mu_Z:\mbf{s}\simp\mbf{s}(Z)\rightarrow\mbf{s}\mrm{D}(Z)$ is in
degree $n$
$$(\mu_Z)^n=\prod(\mu_Z)^{p,q}:\ds\prod_{i+j+q=n}Z^{i,j,q}\rightarrow\ds\prod_{p+q=n}Z^{p,p,q}$$
where, given $p,q$ with $p+q=n$, $(\mu_Z)^{p,q}$ is
$$(\mu_Z)^{p,q}=\sum_{i+j=p} Z(d^0\stackrel{j)}{\cdots}d^0,d^pd^{p-1}\cdots d^{j+1}):Z^{i,j,q}\rightarrow Z^{p,p,q}\ .$$
If $p$ is a f{i}xed integer, note that the sum
$\sum_{i+j=p}Z(d^0\stackrel{j)}{\cdots}d^0,d^pd^{p-1}\cdots
d^{j+1})$ is f{i}nite since the indexes $i,j$ in $Z^{i,j,q}$ are
positive (because $Z\in\Dl\Dl\cdco$).
\end{def{i}}

Thus, the next proposition follows directly from the def{i}nition
of cosimplicial descent category.

\begin{prop}
Under the above notations, if $\mc{A}$ is an additive category
with numerable products, then
$(Ch^\ast\mc{A},\mbf{s},\mu,\lambda)$ is a cosimplicial descent
category with respect to the homotopy
equivalences.\\[0.1cm]
\indent If moreover $\mc{A}$ is abelian then
$(Ch^\ast\mc{A},\mbf{s},\mu,\lambda)$ is a cosimplicial descent
category, where $\mrm{E}$ is the class of quasi-isomorphisms.\\
In both descent structures $\lambda$ is quasi-invertible and $\mu$
is associative.
\end{prop}

This time we will give the specif{i}c consequences of the previous
proposition, because while the cone in $Ch^\ast\mc{A}$ is widely
known, does not occur the same to its dual construction, the path
object.

\begin{prop}\label{pddesCaminoenComplejos}
Given morphisms
$A\stackrel{f}{\rightarrow}B\stackrel{g}{\leftarrow}C$ of cochain
complexes, the \textit{path object} associated with $f$ and $g$ is
a cochain complex $path(f,g)$, functorial in the pair $(f,g)$,
which satisf{i}es the following properties\\[0.1cm]
\textbf{1)} there exists maps in $Ch^\ast\mc{A}$, functorial in
$(f,g)$
$$ \jmath_A:path(f,g)\rightarrow  A\ \ \jmath_B:path(f,g)\rightarrow B$$
such that $\jmath_A$ is a quasi-isomorphism $($resp. homotopy
equivalence$)$ if and only if $g$ is so. Similarly, $\jmath_C$ is
a quasi-isomorphism $($resp. homotopy
equivalence$)$ if and only if $f$ is so.\\[0.1cm]
\textbf{2)} If $f=g=Id_A$, there exists a homotopy equivalence
$P:A\rightarrow path(A)$ in $Ch^\ast\mc{A}$ such that the
composition of $P$ with the projections
$\jmath_A,\jmath'_A:path(A)\rightarrow
A$ given in 1) is equal to the identity.\\[0.1cm]
\textbf{3)} The following square commutes up to homotopy
equivalence
$$\xymatrix@M=4pt@H=4pt{
 B  & A  \ar[l]_-f\\
 C\ar[u]_g & path(f,g)\ar[u]_{\jmath_A}\ar[l]_-{\jmath_C}\ .}$$
\end{prop}

\begin{obs}
When $C=0$, $Path(f,0)$ is (up to homotopy equivalence) the
cochain complex $c(f)[-1]$, where $c(f)$ is the classical cone of
$f$.\\
The proof is just the dual of the one where we checked the
commutativity of diagram (\ref{compatibilidadConoconSimple}) in
proposition \ref{cdcaditivCDS} up to homotopy, and can be found in \cite{H} 2.2.11.\\
Then, under the classical approach of the homotopy theory of
$Ch^\ast\mc{A}$, $f$ gives rise to the distinguished canonical
triangle
$$A\stackrel{f}{\rightarrow}B\stackrel{i}{\rightarrow}c(f)\stackrel{p}{\rightarrow}A[1]$$
and, under our settings, the morphism
$\jmath_A:Path(f,0)\rightarrow A$ corresponds to the projection
$c(f)[-1]\stackrel{p[-1]}{\rightarrow}A$.
Therefore the induced triangulated structure on $HoCh^\ast\mc{A}$
coincides with the usual one.
\end{obs}

\begin{lema}\label{DescripcionPath}
Given two morphisms
$A\stackrel{f}{\rightarrow}B\stackrel{g}{\leftarrow}C$ in $\cdco$,
there exists a natural homotopy equivalence between  $path(f,g)$
and the cochain complex $path_r(f,g)$ given by
$$path_r(f,g)^n=A^n\oplus B^{n-1}\oplus C^n \mbox{ ; }
d=\left( \begin{array}{ccc}
               d_A & 0    & 0\\
               f   & -d_B & g\\
               0   & 0    & 0 \end{array}\right)\ .$$
In addition, the morphisms $\jmath_A$ and $\jmath_C$ correspond to
the respective projections $A^n\oplus B^{n-1}\oplus C^n\rightarrow
A^n$ and $A^n\oplus B^{n-1}\oplus C^n\rightarrow C^n$.
\end{lema}

\begin{proof}[Idea of the proof] By def{i}nition
$path(f,g)=\mbf{s}(Path(f\times\Dl,g\times\Dl))$. The cosimplicial
object $Path(f\times\Dl,g\times\Dl)$ is the image under the total
functor (dual of def{i}nition \ref{def{i}TotCombinatorio}) of $T$,
consisting of the biaugmented bisimplicial cochain complex
$$\xymatrix@M=4pt@H=4pt@C=25pt{                                                                         &                                                                                            &                                                                                                                 &                                                                                                                                              &                                                                                                            & \\
                                      A \ar@/^1pc/[d]\ar@/^0.75pc/[d]\ar@{}[u]|{\vdots}    \ar[r]^f     &  B \ar@/^1pc/[d]\ar@/^0.75pc/[d] \ar@{}[u]|{\vdots}        \ar@<0.5ex>[r] \ar@<-0.5ex>[r]  & B \ar@/^1pc/[d]\ar@/^0.75pc/[d] \ar@{}[u]|{\vdots}      \ar@/_0.75pc/[l]\ar@<0ex>[r]\ar@<1ex>[r] \ar@<-1ex>[r]  &B \ar@/^1pc/[d]\ar@/^0.75pc/[d] \ar@{}[u]|{\vdots}     \ar@/_1pc/[l]\ar@/_0.75pc/[l] \ar@<0.33ex>[r]\ar@<-0.33ex>[r]\ar@<1ex>[r]\ar@<-1ex>[r] & B \ar@/^1pc/[d]\ar@/^0.75pc/[d] \ar@{}[u]|{\vdots}      \ar@/_1pc/[l]\ar@/_0.75pc/[l]\ar@{-}@/_1.25pc/[l]  & \cdots \\
                                      A\ar@/^0.75pc/[d]  \ar[r]^f \ar@<0ex>[u]\ar@<1ex>[u] \ar@<-1ex>[u]&  B \ar@/^0.75pc/[d]\ar@<0ex>[u]\ar@<1ex>[u] \ar@<-1ex>[u]  \ar@<0.5ex>[r]\ar@<-0.5ex>[r]   & B \ar@/^0.75pc/[d]\ar@<0ex>[u]\ar@<1ex>[u]\ar@<-1ex>[u] \ar@/_0.75pc/[l]\ar@<0ex>[r]\ar@<1ex>[r] \ar@<-1ex>[r]  &B \ar@/^0.75pc/[d]\ar@<0ex>[u]\ar@<1ex>[u]\ar@<-1ex>[u]\ar@/_1pc/[l]\ar@/_0.75pc/[l] \ar@<0.33ex>[r]\ar@<-0.33ex>[r]\ar@<1ex>[r]\ar@<-1ex>[r] & B \ar@/^0.75pc/[d]\ar@<0ex>[u]\ar@<1ex>[u]\ar@<-1ex>[u] \ar@/_1pc/[l]\ar@/_0.75pc/[l]\ar@/_1.25pc/[l]      & \cdots \\
                                      A  \ar@<0.5ex>[u] \ar@<-0.5ex>[u]                      \ar[r]^f   &  B \ar@<0.5ex>[u] \ar@<-0.5ex>[u]                          \ar@<0.5ex>[r] \ar@<-0.5ex>[r]  & B \ar@<0.5ex>[u] \ar@<-0.5ex>[u]                        \ar@/_0.75pc/[l]\ar@<0ex>[r]\ar@<1ex>[r] \ar@<-1ex>[r]  &B \ar@<0.5ex>[u] \ar@<-0.5ex>[u]                       \ar@/_1pc/[l]\ar@/_0.75pc/[l] \ar@<0.33ex>[r]\ar@<-0.33ex>[r]\ar@<1ex>[r]\ar@<-1ex>[r] & B \ar@<0.5ex>[u] \ar@<-0.5ex>[u]                        \ar@/_1pc/[l]\ar@/_0.75pc/[l]\ar@/_1.25pc/[l]      & \cdots \\
                                                                                                        &  C  \ar[u]_g                                               \ar@<0.5ex>[r]\ar@<-0.5ex>[r]   & C  \ar[u]_g                                             \ar@/_0.75pc/[l]\ar@<0ex>[r]\ar@<1ex>[r] \ar@<-1ex>[r]  &C  \ar[u]_g                                            \ar@/_1pc/[l]\ar@/_0.75pc/[l] \ar@<0.33ex>[r]\ar@<-0.33ex>[r]\ar@<1ex>[r]\ar@<-1ex>[r] & C  \ar[u]_g                                             \ar@/_1pc/[l]\ar@/_0.75pc/[l]\ar@/_1.25pc/[l]      & \cdots
                                                                                                         }$$
where all maps without label are identities. Then, it follows from
the def{i}nitions that $ K (Tot( T))$ coincides with $Tot\comp
K(T)$. Therefore $path(f,g)=Tot(Tot(K(T)))$ is the total cochain
complex associated with the following triple complex
$$\xymatrix@M=4pt@H=4pt@C=25pt{                                    &                                                        &                                                                              &                                                                                                           &                                                                         & \\
                                      A \ar@{}[u]|{\vdots}\ar[r]^f &  B \ar@{}[u]|{\vdots} \ar[r]_0 & B \ar@{}[u]|{\vdots} \ar[r]_{Id}  &B \ar@{}[u]|{\vdots}\ar[r]_0 & B \ar@{}[u]|{\vdots}  & \cdots \\
                                      A \ar[r]^f \ar[u]_{Id}       &  B \ar[u]_{Id}        \ar[r]_0 & B \ar[u]_{Id}        \ar[r]_{Id}  &B \ar[u]_{Id}       \ar[r]_0 & B \ar[u]_{Id}         & \cdots \\
                                      A \ar[u]_0 \ar[r]^f          &  B \ar[u]_0           \ar[r]_0 & B \ar[u]_0           \ar[r]_{Id}  &B \ar[u]_0          \ar[r]_0 & B \ar[u]_0            & \cdots \\
                                                                   &  C  \ar[u]_g          \ar[r]_0 & C  \ar[u]_g          \ar[r]_{Id}  &C  \ar[u]_g         \ar[r]_0 & C  \ar[u]_g           & \cdots
                                                                                                         }$$
which is homotopic by columns to
$$\xymatrix@M=4pt@H=4pt@C=25pt{                                    &                                &                                   &                             &                       & \\
                                      0 \ar@{}[u]|{\vdots}\ar[r]   &  0 \ar@{}[u]|{\vdots} \ar[r]   & 0 \ar@{}[u]|{\vdots} \ar[r]       &0 \ar@{}[u]|{\vdots}\ar[r]   & 0 \ar@{}[u]|{\vdots}  & \cdots \\
                                      A \ar[u]   \ar[r]^f          &  B \ar[u]             \ar[r]_0 & B \ar[u]             \ar[r]_{Id}  &B \ar[u]            \ar[r]_0 & B \ar[u]              & \cdots \\
                                                                   &  C  \ar[u]_g          \ar[r]_0 & C  \ar[u]_g          \ar[r]_{Id}  &C  \ar[u]_g         \ar[r]_0 & C  \ar[u]_g           & \cdots
                                                                                                         }$$
that is homotopic by rows to
\begin{equation}\label{diagAuxiliarPath}\xymatrix@M=4pt@H=4pt@C=25pt{                                    &                              &                           &                       & \\
                                      0 \ar@{}[u]|{\vdots}\ar[r]   &  0 \ar@{}[u]|{\vdots} \ar[r] &0 \ar@{}[u]|{\vdots}\ar[r] & 0 \ar@{}[u]|{\vdots}  & \cdots \\
                                      A \ar[u]   \ar[r]^f          &  B \ar[u]             \ar[r] &0 \ar[u]            \ar[r] & 0 \ar[u]              & \cdots \\
                                                                   &  C  \ar[u]_g          \ar[r] &0  \ar[u]           \ar[r] & 0  \ar[u]           & \cdots
}\end{equation}
Therefore, the projection of $path(f,g)$ onto the total complex of
(\ref{diagAuxiliarPath}) is a homotopy equivalence. But this total
complex coincides with $path_r(f,g)$, so $\jmath_A$ and $\jmath_C$
correspond to the projections of $path_r(f,g)$ onto $A$ and $C$,
respectively.
\end{proof}

\begin{proof}[Proof of \ref{pddesCaminoenComplejos}]
The statement follows from proposition
\ref{CylLambdaCasiinvertible}. Let us see 3). By loc. cit. there
exists $H:path(f,g)\rightarrow path(B)$ such that $\jmath'_B\comp
H=f\comp \jmath_A$ and $\jmath''_B\comp H= g\comp \jmath_B$, where
$\jmath'_B,\jmath''_B:path(A)\rightarrow A$ are the canonical projections.\\
Then, it suf{f}{i}ces to see that $\jmath'_B$ and $\jmath''_B$ are
homotopic morphisms of $Ch^\ast\mc{A}$, or equivalently, that the
projections $P,Q:path_r(Id_B,Id_B)\rightarrow B$ are homotopic.\\
To this end, consider the homotopy $H:path_r(Id_B,Id_B)^n
\rightarrow B^{n-1}$ def{i}ned as the projection onto the second
summand $B^n\oplus B^{n-1}\oplus B^n\rightarrow B^{n-1}$.
\end{proof}

\begin{obs}\label{simpleNormalizadoCosimpl} Similarly to the chain complexes case,
in the abelian we can consider as well the normalized version of
the simple functor, $\mbf{s}_N$.\\
The functor $path$ obtained using the normalized simple functor is
equal to $path_r$.
\end{obs}

It also holds the dual of propositions \ref{ComplPositivosCDS} and
corollary \ref{ComplAcotadoesTriang}.

\begin{prop} Let $\mc{A}$ be an additive category $($resp.
abelian$)$. Then the category $Ch^{q}\mathcal{A}$ of uniformly
bounded-bellow cochain complexes, together with the homotopy
equivalences $($resp. the quasi-isomorphisms$)$ as equivalences
and the data $\lambda$, $\mu$ given in
\ref{estructuraDescensoCDCos}, is an additive cosimplicial descent
category. In addition, $\lambda$ is quasi-invertible and $\mu$ is
associative.
\end{prop}

We obtain in this way the usual ``cosuspended'' (or left
triangulated) category structure on $Ho Ch^{q}\mathcal{A}$.

\begin{cor} Let $Ch^{\mbf{b}}\mc{A}$ be the category of bounded-bellow cochain complexes.
Then the localized category $Ho Ch^{\mbf{b}}\mc{A}$ of
$Ch^{\mbf{b}}\mc{A}$ with respect to the quasi-isomorphisms
$($resp. homotopy equivalences$)$ is a triangulated category.
\end{cor}

%
%
%
%
%
%
\section{Commutative dif{f}erential graded algebras}\label{SecADGC}

The Thom-Whitney functor and its properties were developed in
\cite{N}. This simple functor gives rise to a cosimplicial descent
category structure on the category of commutative dif{f}erential
graded algebras over a f{i}eld of characteristic 0.

\begin{def{i}}\mbox{}\\
Let $\mbf{Cdga}(k)$\index{Symbols}{$\mbf{Cdga}(k)$} be the
category of commutative dif{f}erential graded algebras
\index{Index}{commutative dif{f}erential graded algebras} (or cdg
algebras) over a f{i}eld $k$ of characteristic 0.\\
The product of two cdg algebras $A$ and $B$ has as underlying
cochain complex the product (i.e. the direct sum) of those
underlying complexes of $A$ and $B$.
\end{def{i}}

\begin{def{i}}[Descent structure on $\mbf{Cdga}(k)$]\mbox{}\\
\noindent\textbf{Simple functor:}
The Thom-Whitney simple functor $\mbf{s}_{TW} :
\Dl\mbf{Cdga}(k)\rightarrow \mbf{Cdga}(k)$, introduced in \cite{N}, is def{i}ned as follows.\\
F{i}rstly, let $\mrm{L}\in\simp\mbf{Cdga}(k)$ \cite{BG} be the
simplicial object of $\mbf{Cdga}(k)$ such that
$$\mrm{L}_n = \ds\frac{\Lambda(x_0,\ldots ,x_n,dx_0,\ldots,dx_n)}{\left( \sum x_i -1 ,\, \sum dx_i\right)},$$
where $\Lambda_n=\Lambda(x_0,\ldots ,x_n,dx_0,\ldots,dx_n)$ is the
free cdg algebra in which $x_k$ has degree 0 and $dx_k$ has degree
1, $0\leq k\leq n$. The boundary map is the unique derivation in
$\Lambda_n$ such that $d(x_k)=dx_k$, $d(dx_k)=0$.\\
The face maps $d_i:\mathrm{L}_{n+1}\rightarrow\mathrm{ L}_{n}$ and
the degeneracy maps $s_j:\mathrm{L}_n\rightarrow \mathrm{L}_{n+1}$
are def{i}ned as
$d_i(x_k)=\left\{\begin{array}{ll} x_k, & k<i\\
                                    0, & k=i\\
                                    x_{k-1}, & k>i\end{array}\right.$
                                    and
$s_j(x_k)=\left\{\begin{array}{ll} x_k, & k<j\\
                                    x_k + x_{k+1}, & k=j\\
                                    x_{k+1}, &
                                    k>j\end{array}\right.$.\\[0.3cm]
Given $A\in \Dl\mbf{Cdga}(k)$, denote by $T_A :
\simp_e\times\Dl_e\rightarrow \mbf{Cdga}(k)$ the bifunctor
obtained from $\mathrm{L}\otimes A :\simp\times\Dl\rightarrow
\mbf{Cdga}(k)$ by forgetting the degeneracy maps. Then the
Thom-Whitney\index{Index}{ Thom-Whitney
simple}\index{Symbols}{$\mbf{s}_{TW}$} simple is the f{i}nal
$$\mbf{s}_{TW}(A)=\int_n T_A (n,n).$$
\noindent \textbf{Equivalences:} The class
$\mrm{E}$ consists of the quasi-isomorphisms, that is to say, those morphisms of cdg algebras which induce isomorphism in cohomology.\\[0.2cm]
\noindent \textbf{Transformation $\mathbf{\lambda}$:} If
$A\in\mbf{Cdga}(k)$, the morphisms $A \rightarrow A \otimes
\mathrm{L}_n$; $a\rightarrow a\otimes 1$ def{i}ne the morphism
$\lambda(A):A\rightarrow\mbf{s}_{TW}(A\times\Dl)$.\\[0.2cm]
\noindent \textbf{Transformation $\mathbf{\mu}$:} If
$Z\in\Dl\Dl\mbf{Cdga}(k)$,
$\mu_{TW}(Z):\mbf{s}_{TW}\simp\mbf{s}_{TW}Z\rightarrow
\mbf{s}_{TW}\mrm{D}Z=\int_p Z^{p,p}\otimes \mrm{L}_p$ is given by
the morphisms
$$\xymatrix@M=4pt@H=4pt{\mbf{s}_{TW}\simp\mbf{s}_{TW}Z \ar[r]^-{\pi} & Z^{p,p}\otimes \mrm{L}_p\otimes\mrm{L}_p \ar[r]^{Id\otimes\tau_p}& Z^{p,p}\otimes \mrm{L}_p}$$
where $\pi$ is the iterated projection and
$\tau_p:\mrm{L}_p\otimes\mrm{L}_p\rightarrow \mrm{L}_p$ is the
structural morphisms of the cdg algebra $\mrm{L}_p$, that are
morphisms of cdg algebras since $L_p$ is commutative.
\end{def{i}}

\begin{prop}
The category $\mbf{Cdga}(k)$ together with the quasi-isomorphisms
and the Thom-Whitney simple is a cosimplicial descent category. In
addition, the forgetful functor $\mrm{U}:\mbf{Cdga}(k)\rightarrow
Ch^\ast k$ is a functor of cosimplicial descent categories.
\end{prop}

We will use the following notations in the proof of the previous
proposition.

\begin{num}\label{stomwitni2*}
%
Let $k[p]$ be the cochain complex that is equal to 0 in all
degrees except $p$, and $(k[p])^p=k$. As in \cite{N} (2.2), we
denote by $\int_{\triangle^p}\; : L_p\rightarrow
k[p]$\index{Symbols}{$\int_{\triangle^p}$} the map of cochain
complexes which in degree $p$ is the ``integral over the simplex
$\triangle^p$'', $L_p^p\rightarrow k$.
\end{num}

\begin{proof}
Let $Ch^\ast k$ be the category of cochain complexes of $k$-vector
spaces. Let us check that the forgetful functor
$\mrm{U}:\mbf{Cdga}(k)\rightarrow Ch^\ast k$ satisf{i}es the dual
proposition to \ref{FDfuerte}, where $Ch^\ast k$ is considered
with the usual descent structure given in def{i}nition
\ref{estructuraDescensoCDCos}.\\
The axiom $\mathrm{(FD 1)}^\comp$ hold, and it is straightforward
to check that $\lambda$ and $\mu$ are compatible
natural transformations.\\
Let us see $\mathrm{(FD 2)}^\comp$. Denote by $\mbf{s}_{TW}:\Dl
Ch^{\ast} k\rightarrow Ch^{\ast} k$ the functor def{i}ned as the
f{i}nal of $\mrm{U}(\mrm{L})\otimes
A:\simp_e\times\Dl_e\rightarrow Ch^{\ast} k$. Then we have the
commutative diagram
$$\xymatrix@M=4pt@H=4pt{ \Dl\mbf{Cdga}(k) \ar[r]^{\Dl\mrm{U}} \ar[d]_{\mbf{s}_{TW}} & \Dl Ch^{\ast}k \ar[d]^{\mbf{s}_{TW}}\\
                           \mbf{Cdga}(k) \ar[r]^{\mrm{U}}                             &   Ch^{\ast}k   .}$$
By (\cite{N}, 2.15), there exists a natural transformation
$I:\mbf{s}_{TW}\rightarrow\mbf{s}:\Dl Ch^{\ast}k \rightarrow
Ch^{\ast}k$ such that $I(A)$ is a homotopy equivalence for each
$A\in Ch^{\ast} k$. Set $\Theta=I\,\Dl \mrm{U} :
\mrm{U}\mbf{s}_{TW}\rightarrow \mbf{s}\Dl\mrm{U}:\Dl
\mbf{Cdga}(k)\rightarrow Ch^{\ast} k$.
More concretely, given $A\in\Dl \mbf{Cdga}(k)$, the morphism
$\Theta_A:\mrm{U}\mbf{s}_{TW}(A)\rightarrow \mbf{s}\Dl\mrm{U}(A)$
in degree $n$ is the map $(\mbf{s}_{TW}(A))^n\rightarrow
(\mbf{s}\Dl\mrm{U}(A))^n=\prod_{p+q=n}A^{p,q}$ whose projection
onto the component $A^{p,q}$ is given by the composition
$$(\mbf{s}_{TW}(A))^n=\int_m (A^m\otimes \mrm{L}_m)^n\stackrel{\pi}{\longrightarrow} (A^p\otimes \mrm{L}_p)^n \stackrel{Id\otimes \int_{\triangle^p}}{\longrightarrow} A^{p,q}\ ,$$
where $\pi$ denotes the projection and $\int_{\triangle^p}$ is the
morphism introduced in \ref{stomwitni2*}.\\
The compatibility between
$\lambda:Id_{\mbf{Cdga}(k)}\rightarrow\mbf{s}_{TW}(-\times\Dl)$
and $\lambda':Id_{Ch^{\ast} k}\rightarrow \mbf{s}(-\times\Dl)$
means the commutativity up to quasi-isomorphism of the diagram
$$\xymatrix@M=4pt@H=4pt@C=50pt@R=40pt{ \mrm{U}(\mbf{s}_{TW}(A\times\Dl))   \ar[d]_-{\Theta(A\times\Dl)}      & \mrm{U}(A) \ar[l]_-{\mrm{U}(\lambda(A))} \ar[ld]^-{\lambda'(\mrm{U}(X))} .\\
                                      \mbf{s}(\mrm{U}(A)\times\Dl) &                    }$$
that in degree $n$ becomes
$$\xymatrix@M=4pt@H=4pt@C=50pt@R=40pt{  \mbf{s}_{TW}^n (A\times\Dl)  \ar[d]_-{\Theta_n (A\times\Dl)}      & A^n \ar[l]_-{\lambda_n(A)} \ar@{(->}[ld]^-{i_n} .\\
                                        \prod_{i\leq n}A^i &                    }$$
If $a\in A^n$, $\lambda_n(A)(a)=\{(1_m\otimes a,0,\ldots)\}_{m\geq
0}\in \mbf{s}_{TW}^n (A\times\Dl) \subset
\ds\prod_{k+l=n}\ds\prod_m  \mathrm{L}_m^k \otimes A^l$, where $1_m\in\mrm{L_m^0}$ is the unit of $\mathrm{L}_m$. Therefore\\
$$\Theta_n (A\times\Dl)(\lambda_n(A)(a))=((\int_{\triangle^0}1)a,0,\ldots)=(a,0,\ldots)=i_n(a).$$
It remains to see that $\mu_{TW}: \mbf{s}_{TW}\Dl\mbf{s}_{TW}
\rightarrow \mbf{s}_{TW} \mrm{D}$, $\mu':\mbf{s}\,\Dl\mbf{s}
\rightarrow \mbf{s} \mrm{D}$ and $\Theta$ are compatible. To see
this, it is enough to prove the commutativity up to homotopy of
the diagram
\begin{equation}\label{compat}\xymatrix@M=4pt@H=4pt@C=40pt@R=40pt{\mrm{U}\mbf{s}_{TW}\Dl\mbf{s}_{TW}(A) \ar[d]_{\widetilde{\Theta}} \ar[r]^-{\mrm{U}(\mu_{TW})} & \mrm{U}\mbf{s}_{TW}(\mrm{D}A))  \ar[d]^{\Theta_{\mrm{D}A}}   \\
                                                      \mbf{s}( \Dl\mbf{s}\,\Dl\Dl\mrm{U}(A))                   \ar[r]^{\mu'}                                   &  \mbf{s}\mrm{D}(\Dl\Dl\mrm{U}(A)),}\end{equation}
where $\widetilde{\Theta}$ is the composition
$$\xymatrix@M=4pt@H=4pt@C=50pt{\mrm{U}\, \mbf{s}_{TW}\,\Dl\mbf{s}_{TW}(A) \ar[r]^-{\Theta_{\Dl\mbf{s}_{TW}(A)}} & \mbf{s}(\Dl \mrm{U}\mbf{s}_{TW}A) \ar[r]^-{\mbf{s}(\Dl\Theta_A)} & \mbf{s}( \Dl\mbf{s}\,\Dl\Dl\mrm{U}(A))}.$$
Using the homotopy inverse $E$ of $\Theta$ given in \cite{N}, one
can obtain a homotopy inverse of $\widetilde{\Theta}$, that will
be referred to as $\widetilde{E}$. Then, it holds that
$\Theta_{\mrm{D}A}\comp \mrm{U}(\mu_{TW}) \comp \widetilde{E}$
coincides with $\mu'$ up to homotopy equivalence.\\
Indeed, (\ref{compat}) is homotopy equivalent to the image under
the total complex of a square of double cochain complexes of
$k$-vector spaces of the form
$$\xymatrix@M=4pt@H=4pt@C=40pt@R=40pt{\mrm{U}\mbf{s}_{TW}^{**}(A)  \ar[r]^-{\mrm{U}(\mu_{TW}^{**})} & \mrm{U}\, \mbf{s}_{TW}^{**}(\mrm{D}A)  \ar[d]^{I^{**}}   \\
                                      \mbf{s}^{**}\mrm{U}(A)  \ar[u]^{{E}^{**}}                 \ar[r]^-{\mu_{E-Z}}                                   & K\,\mrm{D}\,\mrm{U}(A).}$$
where $\mu_{E-Z}$ is the Alexander-Whitney map (see
\ref{Alex-Whitney-Map}). By the Eilenberg-Zilber-Cartier theorem
\ref{E-Z-C}, it is enough to check that the above square commutes
in degree $0$, but this is a straightforward calculation, totally
similar to \cite{N}, (3.4).
\end{proof}


\begin{center} \end{center}
\subsection{Comments on the non-commutative case}

Consider now a commutative (associative and unitary) ring $R$, and
let $\mbf{Dga}$\index{Symbols}{$\mbf{Dga}$} be the category of
dif{f}erential graded $R$-algebras\index{Index}{dif{f}erential
graded algebras} (not necessarily commutative). Denote by $Ch^\ast R$ the category of cochain complexes of $R$-modules.\\
In this case, as we will explain later, the simple functor in
$Ch^\ast R$ (def{i}nition \ref{estructuraDescensoCDCos}) induces
in a natural way the so called Alexander-Whitney simple
$\mbf{s}_{AW}:\Dl\mbf{Dga}\rightarrow \mbf{Dga}$ \cite{N}.\\
Then, a natural question is if this Alexander-Whitney simple
provides a descent structure on $\mbf{Dga}$, together with the
quasi-isomorphisms. The answer is that all axioms of cosimplicial
descent category are satisf{i}ed by this simple, except the
factorization one.\\
The reason is that the transformation $\mu$ of the cosimplicial
descent structure on $Ch^\ast R$ does not induce a natural
transformation in  $\mbf{Dga}$, because it is not compatible with
the multiplicative structures involved.\\
Hence, the descent structure on $Ch^\ast R$ doest not induce a
descent structure on $\mbf{Dga}$. We will give in this subsection
an explicit counterexample of this fact. That is, we will exhibit
a bicosimplicial graded algebra $Z$ such that the morphism of
cochain complexes $\mu_Z:\mbf{s}\Dl\mbf{s}Z\rightarrow
\mbf{s}\mrm{D}$ does not preserve the multiplicative structure, so
it is not a morphism of the category $\mbf{Dga}$.

\begin{num}
By def{i}nition \ref{estructuraDescensoCDCos}, the simple functor
$\mbf{s}:\Dl Ch^\ast R\rightarrow Ch^\ast R$ is the composition of
functors $K$ and $Tot$. Since both functors are monoidal with
respect to the tensorial product of $R$-modules, so is $\mbf{s}$.
Hence, given cosimplicial cochain complexes $X$ and $Y$, we have
the K{\"u}nneth morphism
$$ k: \mbf{s}X\otimes\mbf{s}Y\longrightarrow \mbf{s}(X\otimes Y) $$
that is obtained from the Alexander-Whitney map (see
\ref{Alex-Whitney-Map}) as follows. In degree $n$,
$$k^n : \ds\bigoplus_{p+q=n} (\prod_{i+j=p}X^{i,j}\otimes \prod_{s+t=q}Y^{s,t}) \longrightarrow \prod_{l+m=n}( X^{l}\otimes Y^{l})^m$$
where, as usual, $X^d\in Ch^\ast R$ is denoted by
$\{X^{d,r}\}_{r\in\mathbb{Z}}$ for any $d\geq 0$. Then $k^n$ is
determined by the morphisms
$$\sum_{u+v=l}k^{u,v}: \prod_{i+j=p}X^{i,j}\otimes \prod_{s+t=q}Y^{s,t}\longrightarrow ( X^{l}\otimes Y^{l})^m \ ,$$
If $u,v\geq 0$ with $u+v=l$, $k^{u,v}$ is given by the composition
$$\xymatrix@M=4pt@H=4pt{\prod_{i+j=p}X^{i,j}\otimes \prod_{s+t=q}Y^{s,t} \ar[r]^-{p\otimes p} &  X^{u,p-u}\otimes Y^{v,q-v} \ar[r]^-{A-W} & X^{u+v,p-u}\otimes Y^{u+v,q-v} }$$
where each $p$ denotes the corresponding projection, whereas
${A-W}=(-1)^{uq}X(d^0\stackrel{v)}{\cdots}d^0)\otimes
Y(d^ld^{l-1}\cdots d^{v+1})$. The sign $(-1)^{uq}$ appearing in
the last equation comes from the K{\"u}nneth morphism of the functor
$Tot$.
\end{num}

\begin{def{i}}[Alexander-Whitney Simple] \cite{N}\mbox{}\\
Let $A\in\Dl \mbf{Dga}$ and $\mrm{U}A\in\Dl Ch^\ast R$ be the
cosimplicial cochain complex obtained by forgetting the
multiplicative structure.\\
We have that $\mbf{s}(\mrm{U}A)$ is a dif{f}erential graded
algebra through the morphism
$\tau_{A}:\mbf{s}(\mrm{U}A)\otimes{}\mbf{s}(\mrm{U}A)\rightarrow
{}\mbf{s}(\mrm{U}A)$ def{i}ned as the composition
\begin{equation}\label{productoAW}\mbf{s}(\mrm{U}A)\otimes
\mbf{s}(\mrm{U}A)\stackrel{k}{\longrightarrow}
\mbf{s}(\mrm{U}A\otimes \mrm{U}A)
\stackrel{\mbf{s}\tau}{\longrightarrow} \mbf{s}(A) \end{equation}
where $k$ is the K{\"u}nneth morphism and $\tau^n:\mrm{U}A^n\otimes
\mrm{U}A^n\rightarrow \mrm{U}A^n$ is the structural morphism of
the dif{f}erential graded algebra $A^n$.\\
The Alexander-Whitney simple (\cite{N}, 3.1) is the functor
$\mbf{s}_{AW}:\Dl\mbf{Dga}\rightarrow
\mbf{Dga}$\index{Symbols}{$s_{AW}$} obtained in this way.
\end{def{i}}

\begin{obs} Consider now $B\in\mbf{Dga}$ and $Z\in\Dl\Dl\mbf{Dga}$. We have the following morphisms in $Ch^\ast R$
$$\lambda_{\mrm{U}B}: \mrm{U}B\rightarrow \mbf{s}(B\times\Dl)\ \ \mu_{\mrm{U}Z}:\mbf{s}\mrm{D}\mrm{U}Z\rightarrow \mbf{s}\Dl\mbf{s}\mrm{U}Z \ .$$
It can be checked easily that $\lambda_{\mrm{U}B}$ is compatible
with the respective multiplicative structures of $B$ and
$\mbf{s}_{AW}(B\times\Dl)$, giving rise to a natural
transformation
$\lambda_{AW}:Id_{\mbf{Dga}}\rightarrow \mbf{s}_{AW}(-\times\Dl)$.\\
However, it does not happen the same with $\mu$, as we will see in
the following counterexample.
\end{obs}

\begin{ej} Consider $\mathbb{Z}$ as a differential graded algebra
concentrated in degree 0 and let $\mathbb{Z}\times\Dl\in
\Dl\mbf{Dga}$ be the associated constant cosimplicial object.
Denote by $B\in\Dl\mbf{Dga}$ the path object associated with
$\mathbb{Z}\times\Dl$. In other words, $B$ is equal to
$Path(\mathbb{Z}\leftarrow \mathbb{Z}\times\Dl\rightarrow
\mathbb{Z}\times\Dl)$ (definition \ref{DefiPath}), and can be
visualized as
$$\xymatrix@M=4pt@H=4pt@C=25pt{ \mathbb{Z}\times \mathbb{Z}
\ar@<0.5ex>[r]^-{d^0} \ar@<-0.5ex>[r]_-{d^1}  & \mathbb{Z}\times
\mathbb{Z}\times \mathbb{Z}
\ar@/^1.5pc/[l]\ar@<0ex>[r]\ar@<1ex>[r] \ar@<-1ex>[r] &
\mathbb{Z}\times \mathbb{Z}\times \mathbb{Z}\times \mathbb{Z}
\ar@/^1.25pc/[l]\ar@/^1.5pc/[l] & \cdots }$$
where $d^0(a,b)=(a,a,b)$ and $d^1(a,b)=(a,b,b)$, for any integers $a$ and $b$.\\
Set $Z=B\times \Dl\in\Dl\Dl\mbf{Dga}$, that is,
$Z^{n,m}=B^n=\mathbb{Z}\times\stackrel{n+2}{\cdots}\times
\mathbb{Z}$. In this case the morphism of cochain complexes
$\mu_{Z}:\mbf{s}_{AW}\mrm{D}Z\rightarrow
\mbf{s}_{AW}\Dl\mbf{s}_{AW}Z$ is not a morphism of algebras.\\
Indeed, if we consider $x=(1,0)\in Z^{0,1,0}\subset
(\mbf{s}_{AW}\Dl\mbf{s}_{AW}Z)^1$ and $y=(1,2)\in Z^{0,0,0}\subset
(\mbf{s}_{AW}\Dl\mbf{s}_{AW}Z)^0$, it holds that $\mu(x)\cdot
\mu(y)\neq \mu(x\cdot y)$, (where each $\cdot$ denotes the
corresponding product in $\mbf{s}_{AW}\mrm{D}Z$ and $\mbf{s}_{AW}\Dl\mbf{s}_{AW}Z$).\\
Let us compute f{i}rst $\mu(x)\cdot \mu(y)$.\\
By def{i}nition (see \ref{estructuraDescensoCDCos}) we have that
$\mu(x)=Z(d^0,Id)x=(1,1,0)\in Z^{1,1,0}\subset
(\mbf{s}_{AW}\mrm{D}Z)^1$ and $\mu(y)=Z(Id,d^1)y=y\in
Z^{0,0,0}\subset (\mbf{s}_{AW}\mrm{D}Z)^0$. The product of these
two elements is, following (\ref{productoAW}), equal to the
product in $Z$ of $(1,1,0)$ and $Z(d^1,d^1)y=(1,2,2)$, so
$\mu(x)\cdot
\mu(y)=(1,1,0)\cdot(1,2,2)=(1,2,0) \in Z^{1,1,0}\subset (\mbf{s}_{AW}\mrm{D}Z)^1$.\\
Secondly, $x\cdot y$ is the product in $Z$ of $x\in Z^{0,1,0}$ and
$Z(Id,d^1)y=y\in Z^{0,1,0}$, that
is, $x\cdot y =(1,0) \in Z^{0,1,0}$.\\
Therefore $\mu(x\cdot y)=Z(d^0,Id)(1,0)=(1,1,0)\in Z^{1,1,0}$.\\
Consequently $\mu(x)\cdot \mu(y)=(1,2,0) \neq (1,1,0)=\mu(x\cdot
y)$.
\end{ej}

\begin{obs} Given $Z\in \Dl\Dl\mbf{Dga}$, let $T$ be the 4-simplicial object in $\mbf{Dga}$
given by $T^{i,j,k,l}=Z^{i,j}\otimes Z^{k,l}$. Under the notations
of \ref{Def{i}asocitividadMu}, $Z\otimes
Z=\mrm{D}^{1,3}\mrm{D}^{2,4}T$. In addition
$\mbf{s}\Dl\mbf{s}Z\otimes \mbf{s}\Dl\mbf{s}Z$ is obtained by applying four times $\mbf{s}$ to $T$.\\
In this way, the K{\"u}nneth morphism $k:\mbf{s}\Dl\mbf{s}Z\otimes
\mbf{s}\Dl\mbf{s}Z\rightarrow \mbf{s}\Dl\mbf{s}(Z\otimes Z)$ is
just an iteration of $\mu$. Analogously, under this point of view,
$k:\mbf{s}\mrm{D}Z\otimes \mbf{s}\mrm{D}Z\rightarrow
\mbf{s}\mrm{D}(Z\otimes Z)$ is just the image under $\mu$ of
$\mrm{D}^{1,2}\mrm{D}^{3,4}T$.\\
The preservation of the multiplicative structure by $\mu$ means
the commutativity of the diagram
$$\xymatrix@M=4pt@H=4pt@C=30pt{ \mbf{s}\Dl\mbf{s}Z\otimes \mbf{s}\Dl\mbf{s}Z \ar[r]^-k \ar[d]_{\mu\otimes\mu} & \mbf{s}\Dl\mbf{s}(Z\otimes Z)\ar[r]^-{\mbf{s}\Dl\mbf{s}\tau}\ar[d]_{\mu_{Z\otimes Z}} & \mbf{s}\Dl\mbf{s}Z\ar[d]_{\mu_Z}\\
                                 \mbf{s}\mrm{D}Z\otimes \mbf{s}\mrm{D}Z\ar[r]^-k & \mbf{s}\mrm{D}(Z\otimes Z) \ar[r]^-{\mbf{s}\mrm{D}\tau} & \mbf{s}\mrm{D}Z}$$
The right hand side commutes by the naturality of $\mu$, but the
left hand side commutes provided that $\mu$ is associative and
commutative. This is not the case because the Alexander-Whitney
map fails to be commutative, then $\mu$ is not a
morphism of algebras.\\
In the commutative case we have used the Thom-Whitney simple. The
corresponding natural transformation $\mu$ comes from the product
in the cdg algebra $\mrm{L}$, so this time $\mu$ is actually an
associative and commutative natural transformation.
\end{obs}


%
%
%

\section{DG-modules over a DG-category}\label{SecDGmod}

In this section we study the category of DG-modules over a f{i}xed
DG-category \cite{K} as an example of cosimplicial descent
category. We begin by recalling the def{i}nition of DG-category.

Along this section $R$ will be a f{i}xed commutative ring, and the
tensorial product $\otimes_\mrm{R}$ over $\mrm{R}$ will be written
as $\otimes$.

\begin{def{i}}
If $\mc{A}$ is a category and $A,B$ are objects of $\mc{A}$,
denote by $\mc{A}(A,B)$ the set of morphisms of from $A$ to
$B$ in the category $\mc{A}$.\\
A DG-category $\mc{A}$ (or a \textit{dif{f}erential graded
category})\index{Index}{DG-category} is a category such that given
objects $A,B$ of $\mc{A}$ then
$$\mc{A}(A,B)=\{\mc{A}(A,B)^r\}_{r\in\mathbb{Z}}$$
where each $\mc{A}(A,B)^r$ is an $R$-module.\\
Moreover, $\mc{A}(A,B)$ has a boundary map
$d:\mc{A}(A,B)^r\rightarrow\mc{A}(A,B)^{r+1}$ satisfying the following properties\\
\indent \textbf{0.} the composition of morphisms of $\mc{A}$ is a
homogeneous map of degree 0
$$\mc{A}(A,B)\otimes\mc{A}(B,C)\rightarrow\mc{A}(A,C)$$
\indent \textbf{1.} $d^2=0$, that is, $\mc{A}(A,B)$ is a cochain complex of $R$-modules.\\
\indent \textbf{2.} if $f$ and $g$ are composable morphisms of
$\mc{A}$ and $f$ is homogeneous of degree $p$ then
$$d(f\comp g)=(df)\comp g + (-1)^p f\comp (dg)\ .$$
\end{def{i}}

\begin{ej} The DG-category $\mrm{Dif }\,R$\index{Symbols}{$\mrm{Dif }\,R$} has as objects the cochain complexes of $R$-modules.
If $V$ and $W$ are such cochain complexes, set $\mrm{Dif
}\,R(V,W)=\{\mrm{Dif }\,R(V,W)^p\}_{p\in\mathbb{Z}}$
where
$$\mrm{Dif }\,R(V,W)^p=\{f:V\rightarrow W\mbox{
morphism of }R-\mbox{modules } |\ f(V^k)\subseteq W^{p+k}\}\ .$$
Denote by $f^k=f|_{V^k}:V^k\rightarrow W^{p+k}$.\\
The image under the boundary map $d:\mrm{Dif
}\,R(V,W)^p\rightarrow \mrm{Dif }\,R(V,W)^{p+1}$ of
$f=\{f^k\}_{k\in\mathbb{Z}}$ is
$$\{d^W\comp f^k-(-1)^p f^{k+1}\comp d^V\}_{k\in\mathbb{Z}}\ .$$
\end{ej}

\begin{obs} Note that a morphism $f\in\mrm{Dif }\,R(V,W)^p$ does not commute in general with the boundary maps of $V$ and $W$. Actually,
the cochain complex
$$\cdots {\longrightarrow} \mrm{Dif }\,R(V,W)^{-1} \stackrel{d^{-1}}{\longrightarrow} \mrm{Dif }\,R(V,W)^{0} \stackrel{d^0}{\longrightarrow}\mrm{Dif }\,R(V,W)^{1}\longrightarrow \cdots  $$
is such that $Ker{d^1}$ consists of the morphisms of cochain
complexes between $V$ and $W$, whereas $H^0(\mrm{Dif }\,R(V,W))$
consists of the morphisms between them in the homotopy category
$K(R-mod)$, that is, are equivalence classes of morphisms of
complexes modulo homotopy.
\end{obs}

From now until the end of this section $\mc{A}$ will denote a
f{i}xed DG-category, that we will assume to be a small category.

\begin{def{i}} The category $\mc{CA}$\index{Symbols}{$\mc{CA}$} of \textit{dif{f}erential graded modules}\index{Index}{modules over a DG-category}
over $\mc{A}$ has as objects the functors of DG-categories
$$M:\mc{A}^\comp \rightarrow \mrm{Dif }\,R \ .$$
More concretely, given objects $A$ and $B$ of $\mc{A}$,
$M:\mc{A}(A,B)\rightarrow \mrm{Dif }\,R(MB,MA)$ is a morphism
of cochain complexes (it is $R$-lineal, homogeneous of degree 0 and commutes with the dif{f}erentials).\\
A morphism of $\mc{CA}$ between $M$ and $N$ is a natural
transformation $\tau:M\rightarrow N$ such that for each object $A$
of $\mc{A}$, $\tau_A:MA\rightarrow NA$ is a morphism of cochain
complexes.\\
The category $\mc{CA}$ is an additive category. This additive
structure is induced in a natural way from the additivity of
$\mrm{Dif }\,R $ and $Ch^\ast R$. Actually,  $\mc{CA}$ is an exact
category (cf. \cite{K} 2.2).
\end{def{i}}

Denote by $(Ch^\ast R, {}_{R}\mbf{s},{}_{R}\mu,{}_{R}\lambda)$ the
descent structure on the category of cochain complexes of
$R$-modules $Ch^\ast R$ given in section \ref{Cocadenas}.

\begin{num}\label{objcosCA} Let $M=\{M^n,\, d^i,s^j \}$ be a cosimplicial object of $\mc{CA}$. Then, for each $n\geq 0$, the functor
$M^n:\mc{A}^\comp\rightarrow\mrm{Dif }\,R$ is a functor of
DG-categories. In particular, for a f{i}xed $A\in\mc{CA}$, $M^nA$
is
a cochain complex  that will be written as $\{(M^nA)^q,d^{N^nA}\}_{q\in\mathbb{Z}}$.\\
\indent On the other hand, the face and degeneracy maps of $M$ are
natural transformations $d^i:M^n\rightarrow M^{n+1}$,
$s^j:M^n\rightarrow M^{n-1}$ satisfying the simplicial identities,
and such that their value at each object $A$ of $\mc{A}$ is a
morphism of cochain complexes
$$d^i_Z:M^nA\rightarrow M^{n+1}A\, ,\ \ s^j_A:M^nA\rightarrow M^{n-1}A\ .$$
Therefore, f{i}xed $A\in\mc{A}$, it follows that $MA=\{M^nA,\,
d^i_A,s^j_A\}$ is a cosimplicial cochain complex of $R$-modules.
\end{num}

\begin{def{i}} [Descent structure on $\mc{CA}$]\mbox{}\\
\noindent\textbf{Simple functor:} Given $M\in\Dl\mc{CA}$, the
image under $\mbf{s}M:\mc{A}^\comp\rightarrow\mrm{Dif }\,R$ of an
object $A$ of $\mc{A}$ is def{i}ned as the cochain complex
$$ (\mbf{s}M)A:= {}_{R}\mbf{s}(MA)\mbox{ which is in degree }m((\mbf{s}M)A)^m=\prod_{p+q=m}(M^pA)^q\ .$$
If $f\in \mc{A}(A,B)^r$ and $n\geq 0$ then $M^nf\in \mrm{Dif
}\,R(M^nB,M^n A)^r$ is the morphism $M^nf\!=\! \{ \! M^nf^k\! :\!
(M^nA)^k\rightarrow (M^nA)^{k+r}\}_{k\in\mathbb{Z}}$, given by
$$ (\mbf{s}M)f=\{((\mbf{s}M)f)^k\}_{k\in\mathbb{Z}}\mbox{ where }((\mbf{s}M)f)^k=\!\!\!\!\prod_{p+q=k+r}\!\!\!\! M^p f^{q-r} :((\mbf{s}M)B)^k\rightarrow ((\mbf{s}M)A)^{k+r} \ .$$
On the other hand, if $\tau:M\rightarrow N$ is a morphism of
$\simp\mc{CA}$ then
$$(\mbf{s}\tau)_A ={}_{R}\mbf{s}(\tau_A):(\mbf{s}M)A\rightarrow (\mbf{s}N)A \ .$$
\noindent \textbf{Equivalences:} The class $\mrm{E}$ of
equivalences consists of those morphisms $\rho:M\rightarrow N$
such that $\rho_A:MA\rightarrow NA$ is a quasi-isomorphism in
$Ch^\ast R$ for all $A$ in $\mc{A}$.\\
\noindent \textbf{Transformations $\mathbf{\lambda}$ and $\mu$:}
Given $M\in\mc{CA}$ and $Z\in\simp\simp\mc{CA}$, the
transformations $\lambda(M):M \rightarrow \mbf{s}(M\times\Dl)$ and
$\mu(Z):\mbf{s}\simp\mbf{s}Z\rightarrow\mbf{s}\mrm{D}Z$ are
def{i}ne respectively as
$$\lambda(M)_A = {}_{R}\lambda_{MA}\ \ \mbox{ and } \ \ \mu(Z)_A = {}_{R}\mu_{ZA}$$
for each object $A$ of $\mc{A}$.
\end{def{i}}

\begin{lema} If $M\in\Dl\mc{CA}$, the mapping $A\rightarrow
(\mbf{s}M)A={}_R\mbf{s}(MA)$ def{i}nes a functor
$$\mbf{s}:\Dl\mc{CA}\rightarrow \mc{CA}\ .$$
Following the above notations, in addition
$\lambda:Id_{\mc{CA}}\rightarrow \mbf{s}(-\times\Dl)$ and
$\mu:\mbf{s}\Dl\mbf{s}\rightarrow\mbf{s}\mrm{D}$ are in fact
natural transformations.
\end{lema}

\begin{proof} By (\ref{objcosCA}) $MA\in \Dl Ch^\ast R$, so
${}_R\mbf{s}(MA)\in Ch^\ast R$ and it is and object of $\mrm{Dif
}\, R$. If $f:A\rightarrow B$ is a homogeneous morphism of
$\mc{A}$ of degree $r$, then each $M^nf:M^nB\rightarrow M^nA$ is
an $R$-lineal morphism, and homogeneous of degree $r$. Hence, it
is clear that
$(\mbf{s}M)f:(\mbf{s}M)B\rightarrow(\mbf{s}M)A$ is a morphism of $\mrm{Dif }\, R$.\\
Therefore, $\mbf{s}M:\mc{A}^\comp\rightarrow\mrm{Dif }\,R$ is a
functor, and to see that $\mbf{s}M$ is an object of $\mc{CA}$ it
remains to see that
$$\mbf{s}M:\mc{A}(A,B)\longrightarrow \mrm{Dif }\,R ((\mbf{s}M)B,(\mbf{s}M)A)$$
commutes with the respective boundary maps.\\
Given $n\geq 0$, it holds that $M^n:\mc{A}(A,B)\rightarrow
\mrm{Dif }\,R (M^nf,M^nA)$ commutes with the boundary maps, that
is, if $f\in\mc{A}(A,B)^r$ and $M^nf=\{M^nf^k:M^nB^k\rightarrow
M^nA^{k+r}\}_{k\in\mathbb{Z}}$ then
$$ M^{n}(df)\!=\! d(M^{n}f)\!=\! \{d^{M^nA}\comp M^nf^k -(-1)^r M^nf^{k+1}\comp d^{M^nB}\}_{k\in\mathbb{Z}}\in\! \mrm{Dif }\,R (M^nB,M^nA)^{r+1}$$
and $M^n(df)^k=d^{M^nA}\comp M^nf^k -(-1)^r M^nf^{k+1}\comp
d^{M^nB}$. It follows that
$$((\mbf{s}M)(df))^k=\!\!\!\!\!\!\!\!\prod_{p+q=k+r+1}\!\!\!\!\!\!\!\! M^p(df)^{q-r-1}=\!\!\!\!\!\!\!\!\prod_{p+q=k+r+1}\!\!\!\!\!\!( d^{M^pA}\comp M^pf^{q-r-1} -(-1)^r M^pf^{q-r}\comp d^{M^pB})\ .$$
On the other hand,
$(\mbf{s}M)f=\{((\mbf{s}M)f)^k\}_{k\in\mathbb{Z}}$ with
$((\mbf{s}M)f)^k=\prod_{p+q=k+r}M^p f^{q-r}$, so
$$d((\mbf{s}M)f)^k=d^{(\mbf{s}M)A}\comp \left(\prod_{p+q=k+r}\!\!\!\!\! M^p f^{q-r}\right) -(-1)^r \left(\prod_{s+t=k+r+1}\!\!\!\!\!\! M^s f^{t-r-1}\right)\comp d^{(\mbf{s}M)B}\ .$$
Set
$\partial^{(M^pB)^q}=\sum_{i=0}^p(-1)^{i}d^i_{(M^pB)^{q}}:(M^pB)^q\rightarrow
(M^{p+1}B)^q$, and denote by $d^{(M^pB)^q}:(M^pB)^q\rightarrow
(M^pB)^{q+1}$ the boundary maps of the double complex that is
induced by $MB$, and similarly for $MA$. Note that since
$d^i:M^p\rightarrow M^{p+1}$ is a natural transformation, then
$M^{p+1} f^{q}\comp d^i_{(M^pB)^q} =d^i_{(M^pA)^{q+r}}\comp M^p
f^q$, and we deduce that
$M^{p} f^{q-r}\comp \partial^{(M^{p-1}B)^{q-r}} =\partial^{(M^{p-1}A)^{q}}\comp M^{p-1} f^{q-r}$.\\
By def{i}nition, $d^{(\mbf{s}M)A}:((\mbf{s}M)A)^{k+r}\rightarrow
((\mbf{s}M)A)^{k+r+1}$ and
$d^{(\mbf{s}M)B}:((\mbf{s}M)B)^{k}\rightarrow ((\mbf{s}M)B)^{k+1}$
are
$$d^{(\mbf{s}M)A}=\!\!\!\!\!\!\!\!\prod_{p+q=k+r+1}\!\!\!\!\!\!\!\! d^{(M^pA)^{q-1}}+(-1)^{q}\partial^{(M^{p-1}A)^{q}}\ ; \ d^{(\mbf{s}M)B}\!\!\!\!=\prod_{s+t=k+1}\!\!\!\!\! d^{(M^sB)^{t-1}}+(-1)^{t}\partial^{(M^{s-1}A)^{t}}$$
Therefore
$$d((\mbf{s}M)f)^k=\!\!\!\!\!\!\!\!\prod_{p+q=k+r+1}\!\!\!\!\left( d^{(M^pA)^{q-1}}\comp M^p f^{q-r-1} + (-1)^{q}\partial^{(M^{p-1}A)^{q}}\comp M^{p-1} f^{q-r} +\right. $$
$$\left. -(-1)^r(M^p f^{q-r}\comp d^{(M^pB)^{q-r-1}}+ (-1)^{q-r} M^{p} f^{q-r}\comp \partial^{(M^{p-1}B)^{q-r}})\right)=$$
$$=\prod_{p+q=k+r+1} d^{(M^pA)^{q-1}}\comp M^p f^{q-r-1} -(-1)^r(M^p f^{q-r}\comp d^{(M^pB)^{q-r-1}}=((\mbf{s}M)(df))^k \ .$$
Consequently $\mbf{s}M$ is indeed an object of $\mc{CA}$. The
functoriality of $\mbf{s}$ with respect to the morphisms of
$\mc{CA}$ is clear, so
$\mbf{s}:\Dl\mc{CA}\rightarrow\mc{CA}$ is a functor as requested.\\
The naturality of $\lambda$ and $\mu$ is a straightforward
computation, left to the reader.
\end{proof}

\begin{prop}
Under the above notations, $(\mc{CA},\mbf{s},\mrm{E},\lambda,\mu)$
is an additive cosimplicial descent category. In addition, the
natural transformation $\mu$ is associative and $\lambda$ is
quasi-invertible.
\end{prop}

\begin{proof}
We will see that $(\mc{CA},\mbf{s},\mrm{E},\lambda,\mu)$ satisf{i}es the axioms of the notion of cosimplicial descent category.\\
(CDC 1) is clear. Let us check that the class $\mrm{E}$ is
saturated. Let $A$ be an object of $\mc{A}$. By def{i}nition of
$\mc{CA}$, the evaluation on $A$ provides a functor
$$ev_A:\mc{CA}\rightarrow Ch^\ast R \rightarrow HoCh^\ast R\ .$$
Moreover, if $\rho$ is an equivalence in $\mc{CA}$ then
$ev_A(\rho)=\rho_A$ is an isomorphism of $HoCh^\ast R$. Hence, the
evaluation functor induces $ev_A:Ho\mc{CA}\rightarrow HoCh^\ast
R$, which f{i}ts into the commutative diagram of functors
$$\xymatrix@M=4pt@H=4pt{\mc{CA}\ar[d]^{ev_A} \ar[r]^-{\gamma} & Ho\mc{CA}\ar[d]^{ev_A}\\
                         Ch^\ast R \ar[r]^-{\gamma} & HoCh^\ast R \ .}$$
Therefore, if $\gamma(\rho)$  is an isomorphism of $Ho\mc{CA}$
then $ev_A(\gamma(\rho))=\gamma(\rho_A)$ is an isomorphism of
$HoCh^\ast R $ for each $A\in\mc{CA}$. Thus $\rho_A$ is a
quasi-isomorphism for each $A$, and this means that
$\rho\in\mrm{E}$. It is also clear that $\mrm{E}$ is closed by
products, so (CDC 2) holds, as well as (CDC 3).\\
Given $M\in\mc{A}$ and $Z\in\simp\simp\mc{A}$, the transformations
$\lambda(M):\mbf{s}(M\times\Dl)\rightarrow M$ and
$\mu(Z):\mbf{s}\simp\mbf{s}Z\rightarrow\mbf{s}\mrm{D}Z$ are
equivalences since ${}_{R}\lambda_{MA}$ and ${}_{R}\mu_{ZA}$ are
quasi-isomorphisms for each object $A$ of $\mc{A}$, because
$Ch^\ast R$ is a cosimplicial descent category.\\
Thus (CDC 4), (CDC 5) hold, and we can argue similarly for (CDC 6).\\
Given a morphism $\rho:M\rightarrow N$ of $\Dl\mc{CA}$, to see
(CDC 7) it is enough to note that $[C(\rho)](A)=C(\rho_A)$ in
$\simp Ch^\ast R$.\\
F{i}nally, (CDC 8) follows from the equality
$(\mbf{s}\Upsilon\rho)(A))={}_{R}\mbf{s}(\Upsilon\rho_A)$.
\end{proof}

\begin{obs}
Given $M\in\mc{CA}$, we denote by $[-1]:\mc{CA}\rightarrow\mc{CA}$
the usual shift functor $[-1]:Ch^\ast R\rightarrow Ch^\ast R$.\\
On the other hand, we have the shift functor induced by the
descent structure on $\mc{CA}$, that is
$\mrm{S}M=\mbf{s}Path(0\rightarrow M \leftarrow 0)$.\\
If $A\in\mc{A}$ then $(\mrm{S}M)A={}_R\mbf{s}Path(0\rightarrow MA
\leftarrow 0)$ and the inclusion $(MA)[-1]=path_r(0\rightarrow MA
\leftarrow 0)={}_R\mbf{s}_N Path(0\rightarrow MA \leftarrow 0)$
into $(\mrm{S}M)A$ is a natural homotopy equivalence (see
\ref{DescripcionPath} and \ref{simpleNormalizadoCosimpl}).\\
Then the functors $\mrm{S}$, $[-1]:Ho\mc{CA}\rightarrow Ho\mc{CA}$
are isomorphic, so $\mrm{S}:Ho\mc{CA}\rightarrow Ho\mc{CA}$ is an
isomorphism of categories. Hence, we obtain from the dual of
theorem \ref{TmaHoDtriangulada} the (well-known) triangulated
structure on $Ho\mc{CA}$.
\end{obs}

%
%
%
%
%

\section{F{i}ltered cochain complexes}\label{CF}

Given an abelian category $\mc{A}$, let $\mrm{CF}^+\mc{A}$ be the
category of f{i}ltered positive cochain complexes, f{i}ltered by a
biregular f{i}ltration. In this section we will endow
$\mrm{CF}^+\mc{A}$ with two dif{f}erent descent structures, whose
equivalences will be the f{i}ltered quasi-isomorphisms on one
hand, and ${E}_2$-isomorphism on the other hand. Both structures
will be
related through the ``decalage'' functor of a f{i}ltered complex \cite{DeII}.\\

The category of positive cochain complexes will be written as
$Ch^+\mc{A}$, whose objects are complexes $\{X^n,d\}$ such that
$X^n=0$ if $n< 0$.

\begin{def{i}}\label{DefCatObjF{i}lt} A (decreasing) \textit{f{i}ltration} $\mrm{F}$ of an object $K$ of $\mc{A}$ is a family $\{\mrm{F}^k K\}_{k\in\mathbb{Z}}$
of subobjects of  $K$ such that $\mrm{F}^k K\subseteq \mrm{F}^{l} K$ if $l\leq k$.\\
Denote by $\mrm{F}\mc{A}$\index{Symbols}{$\mrm{F}\mc{A}$} the
additive category whose objects are pairs $(K,\mrm{F})$ consisting
of an object $K$ of $\mc{A}$ together with a f{i}ltration
$\mrm{F}$ of $K$, and whose morphisms are those morphisms of
$\mc{A}$
compatible with the f{i}ltrations.\\
A f{i}ltration $\mrm{F}$ is said to be
\textit{f{i}nite}\index{Index}{f{i}ltration!f{i}nite} if there
exists integers $n,m\in\mathbb{Z}$ such that $\mrm{F}^nK=K$ and $\mrm{F}^mK=0$.\\
The full subcategory of $\mrm{F}\mc{A}$ whose objects are the
complexes f{i}ltered by a f{i}nite f{i}ltration will be denoted by
$\mrm{F}_f\mc{A}$\index{Symbols}{$\mrm{F}_f\mc{A}$}. Of course,
this is an additive category as well.
\end{def{i}}

\begin{obs}
We can consider similarly increasing f{i}ltrations instead of
decreasing ones. If $\mrm{F}$ is a decreasing f{i}ltration, the
convention $\mrm{F}_kK=\mrm{F}^{-k}K$ allows us to reduce our
study to the case of decreasing f{i}ltrations.

\end{obs}

\begin{def{i}}\label{f{i}ltracionBirregular}\index{Index}{f{i}ltration!biregular}
Let
$\mrm{CF}^+\mc{A}$\index{Symbols}{$\mrm{CF}^+\mc{A}$}\index{Index}{f{i}ltered
complexes} be the additive category of pairs $(A,\mrm{F})$, where
$A$ is a f{i}ltered  positive cochain complex and $\mrm{F}$ is a
biregular decreasing f{i}ltration of $A$. In other words,
$\mrm{F}=\{\mrm{F}^k A\}_{k\in\mathbb{Z}}$ is such that
\begin{itemize}
 \item[1.] $\mrm{F}^k A$ is a subcomplex of $A$ $\forall\, k$ and $A=\cup_k \mrm{F}^kA$.
 \item[2.]$\mrm{F}^{k+1}A \subseteq \mrm{F}^{k}A$ $\forall\, k$.
 \item[3.] Given $q\geq 0$, the f{i}ltration $\{(\mrm{F}^k A)^q\}_k$ of $A^q$ is f{i}nite. Then, there exists integers $a$ and $b$
 such that $(\mrm{F}^aA)^q=A^q$ and $(\mrm{F}^bA)^q=0$.
\end{itemize}
a morphism $f:(A,\mrm{F})\rightarrow (B,\mrm{G})$ of $\filt$ is a
morphism $f:A\rightarrow B$ of cochain complexes such that
$f(\mrm{F}^kA)\subseteq \mrm{G}^kB$, for all $k$.
\end{def{i}}

\begin{obs} Equivalently, $\mrm{CF}^+\mc{A}$ is the category of positive cochain complexes of the additive category
$\mrm{F}_f\mc{A}$.
\end{obs}

\subsection{F{i}ltered quasi-isomorphisms}

\begin{def{i}}\label{Def{i}QuisF{i}lt}
For each $k\in\mathbb{Z}$, the \textit{graded} functor
$\mbf{Gr}_k:\filt\rightarrow
Ch^+\mc{A}$\index{Symbols}{$\mbf{Gr}$}\index{Index}{graded
complex} is def{i}ned as
$${}_{\mrm{F}}\mbf{Gr}_k A = \begin{array}{cc} \mrm{F}^k A\hspace{-0.2cm}&\\[-0.5cm]
                                                                         & \hspace{-0.2cm} \mbox{\Large{$\diagup$}} \! \mrm{F}^{k+1}A \end{array} \ .$$
for a f{i}ltered cochain complex $(A,\mrm{F})$.
A morphism $f$ of $\filt$ is a f{i}ltered
quasi-isomorphism\index{Index}{f{i}ltered quasi-isomorphism} if
$\mbf{Gr}_k(f)$ is a quasi-isomorphism for all $k\in\mathbb{Z}$.

Let
$(Ch^+\mc{A})^\mathbb{Z}$\index{Symbols}{$(Ch^+\mc{A})^\mathbb{Z}$}
be the category of graded cochain complexes, whose objects are
families indexed over $\mathbb{Z}$ of positive cochain complexes.\\
The \textit{graded} functor $\mbf{Gr}:\filt\rightarrow
(Ch^+\mc{A})^\mathbb{Z}$ applied to $(A,\mrm{F})$ is in degree $k$
the complex ${}_{\mrm{F}}\mbf{Gr}_k A$.
\end{def{i}}

\begin{def{i}}[Descent structure on $\filt$]\label{EstructuraDescensoCatF{i}lt}\mbox{}
\begin{itemize}
 \item[$\bullet$] Let $\mbf{s}:\Dl Ch^+\mc{A}\rightarrow Ch^+\mc{A}$ be the simple introduced in \ref{Cocadenas}.
Given $(A,F)\in\Dl\filt$ denote by $\mbf{s}(F)$ the f{i}ltration
of $\mbf{s}(A)$ def{i}ned as $(\mbf{s}(\mrm{F}))^k (\mbf{s}A) =
\mbf{s}(\mrm{F}^kA)$. The simple functor
$(\mbf{s},\mbf{s}):\Dl\filt\rightarrow\filt$ is given by
$(\mbf{s},\mbf{s})(A,\mrm{F})=(\mbf{s}(A),\mbf{s}(\mrm{F}))$.
 \item[$\bullet$] The class $\mrm{E}$ consists of the f{i}ltered quasi-isomorphisms.
 \item[$\bullet$] The natural transformations $\lambda$ and $\mu$ in $\filt$
are the same as in the cochain complex case.
\end{itemize}
\end{def{i}}

As well as in the cubical case, \cite{GN} 1.7.5, it holds the
following proposition.

\begin{prop}\label{f{i}ltrados-quis}
Under the notations introduced in
\ref{EstructuraDescensoCatF{i}lt},
$(\filt,(\mbf{s},\mbf{s}),\mrm{E},\lambda,\mu)$ is an additive
cosimplicial descent category. In addition, $\mu$ is associative,
and $\lambda$ is quasi-invertible.
\end{prop}

\begin{proof}
F{i}rstly, note that if $\mathbb{Z}$ is the discrete category
whose objects are the integers and whose morphisms are the
identities, then $(Ch^+\mc{A})^\mathbb{Z}$ is just the category of
functors from $\mathbb{Z}$ with values in $Ch^+\mc{A}$. By
\ref{CategFuntores}, $(Ch^+\mc{A})^\mathbb{Z}$ is a cosimplicial
descent category, with the simple functor induced degreewise, and
with equivalences those morphisms that are degreewise
quasi-isomorphisms.\\
Let us see that proposition \ref{FDfuerte}${}^{\mrm{op}}$ holds
for $\mbf{Gr}:\filt\rightarrow (Ch^+\mc{A})^\mathbb{Z}$.

(SDC 1)$^{\mrm{op}}$ is clear since $\filt$ is additive. To see
(SDC 3)'$^{\mrm{op}}$, we will check that for each $(A,\mrm{F})$
in $\Dl\filt$, $\mbf{s}(\mrm{F})$ is biregular. Denote
$A=\{A^{nm}\}_{n,m}$ where $n$ is the cosimplicial degree and $m$ is the degree relative to $Ch^+\mc{A}$.\\
F{i}xed $k\in\mathbb{Z}$, the complex $(\mbf{s}(\mrm{F}))^k
(\mbf{s}(A))$ is in degree $q$
$\mbf{s}(\mrm{F}^kA)^q=\bigoplus_{i+j=q}\mrm{F}^k A^{i,j} $. By
assumption $\mrm{F}$ is biregular on each $A^n$, so given $p\geq
0$ there exists $a=a(n,p)$ and $b=b(n,p)$ with $\mrm{F}^a
A^{n,p}=A^{n,p}$ and $\mrm{F}^b A^{n,p}=0$. Let
$\alpha=\alpha(q)=\min\{a(i,j)\; | \; i+j=q\; ; \;  i,j\geq 0\}$
and $\beta=\beta(q)=\max\{b(i,j)\; | \; i+j=q\; ; \; i,j\geq 0\}$.
Then $\mbf{s}(\mrm{F}^\alpha A)^q=\mbf{s}(A)^q$ and
$\mbf{s}(\mrm{F}^\beta A)^q=0$, so $\mbf{s}(\mrm{F})$ is
biregular.

To see (SDC 4$)'^{\mrm{op}}$ and (SDC 5$)'^{\mrm{op}}$, let us
check that  the transformations $\mu$ and $\lambda$ of
\ref{Cocadenas} are indeed morphism in $\filt$.\\
If $(A,\mrm{F})\in\filt$, then $\mbf{s}(A\times\Dl)^n=A^n\oplus
A^{n-1}\oplus\cdots\oplus A^0$ and $\lambda(A):A\rightarrow
\mbf{s}(A\times\Dl)$ is the inclusion. Therefore
$$(\lambda(A))(\mrm{F}^k A^{n})=\mrm{F}^k A^{n}\subseteq (\mbf{s}(\mrm{F}))^k
(\mbf{s}(A\times\Dl)^n)=\mrm{F}^k A^{n}\oplus \mrm{F}^k
A^{n-1}\oplus\cdots\oplus \mrm{F}^k A^{0}\ .$$
On the other hand, if $(Z,\mrm{F})\in\Dl\Dl\filt$, the restriction
of
$\mu(Z):\ds\oplus_{i+j+q=n}Z^{i,j,q}\rightarrow\ds\oplus_{p+q=n}Z^{p,p,q}$
to $Z^{i,j,q}$ is $Z(d^0\stackrel{j
)}{\cdots}d^0,d^pd^{p-1}\cdots d^{j+1})$, where $p=i+j$.\\
Moreover
$(\mbf{s}\Dl\mbf{s}\mrm{F})^k(\mbf{s}\Dl\mbf{s}Z)^n=\ds\oplus_{i+j+q=n}\mrm{F}^k
Z^{i,j,q}$ and
$(\mbf{s}\mrm{D}\mrm{F})^k(\mbf{s}\mrm{D}Z)^n=\ds\oplus_{p+q=n}\mrm{F}^k
Z^{p,p,q}$. Thus, $\mu(Z)$ is morphism in $\filt$ since
$Z(d^0\stackrel{j )}{\cdots}d^0,d^pd^{p-1}\cdots d^{j+1})$
preserves the f{i}ltration $\mrm{F}$ for each $i,j,q$.

Secondly, (FD 1)$^{\mrm{op}}$ holds because $\mbf{Gr}$ is
additive. It remains to see (FD 2)$^{\mrm{op}}$. We have that the
diagram
$$\xymatrix@H=4pt@C=30pt@R=20pt{\Dl\filt \ar[r]^-{\Dl\mbf{Gr}}\ar[d]_{(\mbf{s},\mbf{s})} & \Dl (Ch^+\mc{A})^\mathbb{Z} \ar[d]^{\mbf{s}}\\
                                \filt\ar[r]^{\mbf{Gr}}         & (Ch^+\mc{A})^\mathbb{Z}}$$
commutes up to canonical isomorphism. Indeed, given
$(A,\mrm{F})\in\filt$, $k\in\mathbb{Z}$ and $n\geq 0$ it holds
that
$${}_{\mbf{s}(\mrm{F})}\mbf{Gr}_k
(\mbf{s}(A)^n)=\ds\frac{\mbf{s}(\mrm{F})^k(\mbf{s}(A)^n)}{\mbf{s}(\mrm{F})^{k+1}(\mbf{s}(A)^n)}=\ds\frac{\ds\bigoplus_{i+j=n}\mrm{F}^k
A^{i,j}}{\ds\bigoplus_{i+j=n}\mrm{F}^{k+1} A^{i,j}}\simeq
\ds\bigoplus_{i+j=n}\ds\frac{\mrm{F}^{k} A^{i,j}}{\mrm{F}^{k+1}
A^{i,j}}=\mbf{s}({}_{F}\mbf{Gr}_k (A))^n$$
and the boundary maps coincide.\\
\indent If $(A,\mrm{F})\in\filt$, we have that
$\mbf{Gr}(\lambda(A))$ corresponds to the inclusion
${}_{\mrm{F}}\mbf{Gr}_k (A)\rightarrow
\mbf{s}(({}_\mrm{F}\mbf{Gr}_k A)\times\Dl)\simeq
{}_{\mbf{s}(\mrm{F})}\mbf{Gr}_k
(\mbf{s}(A\times\Dl))$.\\
F{i}nally, if $(Z,\mrm{F})\in\Dl\Dl\filt$, it is clear that
$\mbf{Gr}_k(\mu(Z))$ corresponds to $\mu({}_{\mrm{F}}\mbf{Gr}_kZ)$
through the isomorphisms
$${}_{\mbf{s}\Dl\mbf{s}(\mrm{F})}\mbf{Gr}_k(\mbf{s}\Dl\mbf{s}(Z))\simeq \mbf{s}\Dl\mbf{s}({}_{\mrm{F}}\mbf{Gr}_kZ))\mbox{ and } {}_{\mbf{s}\mrm{D}(\mrm{F})}\mbf{Gr}_k(\mbf{s}\mrm{D}(Z))\simeq \mbf{s}\mrm{D}({}_{\mrm{F}}\mbf{Gr}_kZ))\ .$$
\end{proof}

\begin{obs}\label{CFk}
By simplicity, we have considered the category $\mrm{CF}^+\mc{A}$
of uniformly bounded-bellow cochain complexes, with uniform bound
equal to 0, but the arguments remain valid for any f{i}xed value
of the bound. So, if $k$ is a f{i}xed integer and
$\mrm{CF}^k\mc{A}$ is the category of f{i}ltered (by a biregular
f{i}ltration) cochain complexes $(A,\mrm{F})$ such that $A^n=0$
when $n< k$, then
$$(\mrm{CF}^k\mc{A},(\mbf{s},\mbf{s}),\mrm{E},\lambda,\mu)$$
is an additive cosimplicial descent category, where $\lambda$ is
quasi-invertible and $\mu$ associative.
\end{obs}

\begin{def{i}}[F{i}ltered homotopies]\label{Def{i}HomotF{i}lt}\mbox{}\\
Since $\mrm{CF}^+\mc{A}=Ch^+(\mrm{F}_f\mc{A})$, the homotopy
theory of $\mrm{CF}^+\mc{A}$ is just the one coming from $Ch^+\mrm{F}_f\mc{A}$.\\
Then, a \textit{f{i}ltered homotopy}\index{Index}{f{i}ltered
homotopy} between the morphisms $f,g:(A,\mrm{F})\rightarrow
(B,\mrm{G})$ in $\mrm{CF}^+\mc{A}$ is a homotopy
$h:A^{i+1}\rightarrow B^i$ that preserve the f{i}ltrations (that
is, $h(\mrm{F}^k(A^{i+1}))\subseteq \mrm{G}^k(B^i)$) and such that
it is a usual homotopy  between $f$ and $g$ (that is, $d^B\comp
h+h\comp d^A=f-g$). In this case we will say that $f$ is homotopic
to $g$ in $\filt$.
\end{def{i}}

\begin{cor}\label{propiedadesObjCaminoenCF}
Given morphisms
$A\stackrel{f}{\rightarrow}B\stackrel{g}{\leftarrow}C$ of
f{i}ltered cochain complexes, the \textit{path object} associated
with $f$ and
$g$ is a cochain complex $path(f,g)$, functorial in $(f,g)$, which satisf{i}es the following properties\\[0.1cm]
\textbf{1)} there exists functorial maps in $\filt$
$$ \jmath_A:path(f,g)\rightarrow  A\ \ \jmath_B:path(f,g)\rightarrow B$$
such that $\jmath_A$ $($resp. $\jmath_C)$ is a f{i}ltered
quasi-isomorphism if and only if $g$ $($resp. $f)$ is so.\\[0.1cm]
\textbf{2)} If $f=g=Id_A$, then there exists a f{i}ltered
quasi-isomorphism $P:A\rightarrow path(A)$ of $\filt$ such that
the composition of $P$ with the projections
$\jmath_A,\jmath'_A:path(A)\rightarrow A$ given in 1) is equal to the identity on $A$.\\[0.1cm]
\textbf{3)} The following square commutes up to f{i}ltered
homotopy equivalence
$$\xymatrix@M=4pt@H=4pt{
 B  & A  \ar[l]_-f\\
 C\ar[u]_g & path(f,g)\ar[u]_{\jmath_A}\ar[l]_-{\jmath_C}\ .}$$
\end{cor}

\begin{obs}
If $C=0$,  then $Path(f,0)$ is (up to natural f{i}ltered homotopy
equivalence) the cochain complex $c(f)[-1]$, where $c(f)$ denotes
the classical cone, f{i}ltered by the induced f{i}ltration by those of $A$ and $B$.\\
Then, following the classical homotopy theory of
$\mrm{CF}^+\mc{A}$, $f$ gives rise to the distinguished triangle
$$A\stackrel{f}{\rightarrow}B\stackrel{i}{\rightarrow}c(f)\stackrel{p}{\rightarrow}A[1]\ .$$
In our setting the morphism $\jmath_A:Path(f,0)\rightarrow A$
corresponds to $c(f)[-1]\stackrel{p[-1]}{\rightarrow}A$.\\
On the other hand, the object $path(f,g)$ is homotopic to the
complex given in \ref{DescripcionPath}, f{i}ltered by the
f{i}ltration which is induced in a natural way by those of $A$,
$B$ and $C$.
\end{obs}

The category $Ho\filt$ is a subcategory of the usual f{i}ltered
derived category associated with $\mc{A}$,
$D\mrm{F}\mc{A}=\mrm{CF}\mc{A}[\mrm{E}^{-1}]$ (where the cochain complexes does not need to be bounded).\\
It is known that the class of equivalences $\mrm{E}$ has calculus
of fractions in $KF\mc{A}$, and the description of the f{i}ltered
derived category deduced of this fact is similar to the one given
in the following corollary, obtained using our descent techniques.

\begin{cor}\label{descripcionCatDerivF{i}ltr}
The category $Ho\filt$ is additive. A morphism $F:X\rightarrow Y$
of $Ho\filt$ can be represented by a zig-zag in the form
$$\xymatrix@M=4pt@H=4pt{X  & T \ar[r]^{{f}} \ar[l]_w &  Y} ,\  w\mbox{ is a f{i}ltered quasi-isomorphism} \ .$$
Another zig-zag $\xymatrix@M=4pt@H=4pt{X  & S \ar[l]_u
\ar[r]^{{g}} & Y}$ represents $F$ if and only if there exists a
hammock $($commutative in $\filt)$, relating both zig-zags, in the
form
$$
\xymatrix@M=4pt@H=4pt{
                                   &  X                         & T\ar[r]^{{f}}\ar[l]_{w}                 &  Y                        &          \\
 X \ar[ru]^{Id}\ar[rd]_{Id} \ar[r] & \widetilde{X}\ar[d]\ar[u]  & U\ar[d]\ar[u]  \ar[r]^{h} \ar[l]\ar[r]  & \widetilde{Y}\ar[u]\ar[d] & Y\ar[lu]_{Id}\ar[ld]^{Id}\ar[l] \ ,\\
                                   &  X                         & S \ar[r]^{{g}}\ar[l]_{u}                & Y                         &                             }$$
where all maps except $f$, $g$ and $h$ are f{i}ltered
quasi-isomorphisms.
\end{cor}

One can proceed similarly as in proposition
\ref{ComplAcotadoesTriang} to deduce the following corollary.

\begin{cor} Let $\mrm{CF}^\mathbf{b}\mc{A}$ be the category of (non-uniformly) bounded-bellow cochain complexes that are f{i}ltered
by a biregular f{i}ltration. Then the localized  category
$D\mrm{F}^\mathbf{b}\mc{A}$ of $\mrm{CF}^\mathbf{b}\mc{A}$ with
respect to the f{i}ltered quasi-isomorphisms is a triangulated
category.
\end{cor}

The well-known triangulated structure on
$D\mrm{F}^\mathbf{b}\mc{A}$ is usually obtained in the literature
as a consequence of the exact structure on $\mrm{CF}A$. However,
this triangulated structure can be obtained directly (see
\cite{IlI} p. 271), and this is the approach recover here.\\

\begin{obs} In the case of unbounded cochain complexes, the
simple functor of a biregular f{i}ltration is not in general a
biregular f{i}ltration. This is why we have reduced ourselves to the uniformly bounded-bellow case.\\
However, we can also apply these techniques in the case of not
necessarily regular f{i}ltrations, dropping the boundness
condition. The same happens with biregular f{i}ltrations that are
zero outside uniform upper and lower bounds. In other words, let
$\mrm{CF}_f\mc{A}$ be the category whose objects are cochain
complexes together with a f{i}ltration $\mrm{F}$ such that
$$0=\mrm{F}^{M}A \subseteq \mrm{F}^{1}A \subseteq \cdots \subseteq \mrm{F}^0A=A$$
where $M$ is a f{i}xed integer. In this case, it can be proved
similarly that $(\mrm{CF}_f\mc{A},(\mbf{s},\mbf{s}),\mu,\lambda)$
is an additive cosimplicial descent category.
\end{obs}

%
%
%
%
%
%

\subsection{\textit{E}$_\mbf{2}$-isomorphisms}

\noindent We recall the def{i}nition of the spectral sequence
associated with a f{i}ltered cochain complex.

\begin{def{i}}\label{SucEspectral}\index{Index}{spectral sequence}
Let $(A,\mrm{F})\in\filt$, and $r\geq 0$ $p,q\in\mathbb{Z}$.
Def{i}ne $Z_r^{p,q}$, $B_r^{p,q}$ and
$E_r^{p,q}$\index{Symbols}{$E_r^{p,q}$}, where $E_r^{p,q}$ is
called the spectral sequence associated with the f{i}ltration
$\mrm{F}$, as follows
$$\begin{array}{ccrl}
 Z_r^{p,q} &=& \ker &\left\{d:\mrm{F}^pA^{p+q}\rightarrow\ds\frac{A^{p+q}}{\mrm{F}^{p+r}A^{p+q}}\right\}\\[0.1cm]
 \ds\frac{A^{p+q}}{B_r^{p,q}}&=&\mathrm{coker}&\left\{d:\mrm{F}^{p-r+1}A^{p+q-1}\rightarrow \ds\frac{A^{p+q}}{\mrm{F}^{p+1}A^{p+q}}\right\}\\[0.1cm]
 E_r^{p,q}&=&\mathrm{Im}&\left\{Z_r^{p,q}\rightarrow\ds\frac{A^{p+q}}{B_r^{p,q}}\right\}=\ds\frac{Z_r^{p,q}}{B_r^{p,q}\cap Z_r^{p,q}}\end{array}\ .$$

The boundary map $d_r:E_r^{p,q}\rightarrow E_r^{p+r,q-r+1}$ is
induced by the one of $A$. The equality $d_rd_r=0$ holds and
\begin{equation}\label{E_r+1}E_{r+1}^{p,q}=\mrm{H}(\xymatrix@H=4pt@M=4pt{E_r^{p-r,q+r-1}\ar[r]^-{d_r}&
E_r^{p,q}\ar[r]^-{d_r}&E_r^{p+r,q-r+1}})\ .\end{equation}
\end{def{i}}

\begin{obs}\label{pdadesSucEsp}\mbox{}\\[-0.3cm]
\begin{itemize}
 \item[a)] For $r=0$, it holds that $E_0^{p,q}=
{}_\mrm{F}\mbf{Gr}_p(A^{p+q})$ and $d=d_0:E_0^{p,q}\rightarrow
E_0^{p,q+1}$.\\
Then $E_0: \filt\rightarrow(Ch\mc{A})^{\mathbb{Z}}$, so in
$E_0^{p,q}$ we have that $q$ is the degree corresponding to
$Ch\mc{A}$ and $p$ the one corresponding to $\mathbb{Z}$.
 \item[b)] For $r=1$,
 $E_1^{p,q}\stackrel{\star}{=}\ds\frac{d^{-1}(\mrm{F}^{p+1}A^{p+q+1})\cap\mrm{F}^{p}A^{p+q}}{d(\mrm{F}^pA^{p+q-1})+\mrm{F}^{p+1}A^{p+q}}$,
 and $d=d_1:E_1^{p,q}\rightarrow E_1^{p+1,q}$.\\[0.2cm]
Therefore $E_1: \filt\rightarrow(Ch\mc{A})^{\mathbb{Z}}$, in such
a way that in $E_1^{p,q}$, $p$ is the degree corresponding to
$Ch\mc{A}$ whereas $q$ corresponds to $\mathbb{Z}$.
 \item[c)] By (\ref{E_r+1}), a morphism $f$ of $\filt$ is a f{i}ltered quasi-isomorphism if and only if $E_1(f)$ is an isomorphism.
 \item[d)] Similarly, $E_2^{p,q}(f)$ is an isomorphism for all $p,q$ if and only if $E_1(f)$ is a quasi-isomorphism in
 $(Ch\mc{A})^{\mathbb{Z}}$ (that is, it is a quasi-isomorphism in $Ch\mc{A}$ degreewise).
We will refer to such morphisms as $E_2$-isomorphisms.
\end{itemize}
\end{obs}

\begin{def{i}}[Second descent structure on $\filt$]\label{SegundaEstructuraDescensoCatF{i}lt}\mbox{}\\[-0.3cm]
\begin{itemize}
\item[$\bullet$] The simple functor
$(\mbf{s},\delta):\Dl\filt\rightarrow\filt$ is def{i}ned as
follows. If $(A,\mrm{F})\in\filt$, then
$(\mbf{s},\delta)(A,\mrm{F})=(\mbf{s}(A),\delta\mrm{F})$, where
$\mbf{s}(A)$ is the usual simple of cochain complexes. On the
other hand, $\delta\mrm{F}$ is the \textit{diagonal
f{i}ltration}\footnote{The f{i}ltration $\delta F$ is the diagonal
f{i}ltration of $\mbf{s}\mrm{F}$ and of the natural f{i}ltration
$\mrm{G}$ of $\mbf{s}A$ given by $\mrm{G}^q=\oplus_{p\leq q}
A^{p,\cdot}$} over $\mbf{s}(A)$, given by
$$ (\delta\mrm{F})^k(\mbf{s}(A)^n) = \ds\bigoplus_{i+j=n}\mrm{F}^{k-i}A^{i,j} \ .$$
\item[$\bullet$] The class of equivalences $\mrm{E}_2$ consists of
the $E_2$-isomorphisms.
\item[$\bullet$] The transformations $\lambda$ and $\mu$ are the
same as in \ref{EstructuraDescensoCatF{i}lt}, that is, the same as
in the cochain complexes case.
\end{itemize}
\end{def{i}}

\begin{def{i}}[``decalage'' functor]\mbox{}\\
Let $Dec:\filt\rightarrow\filt$\index{Index}{decalage!of a
f{i}ltration}\index{Symbols}{$Dec$} be the functor that maps the
f{i}ltered complex $(A,\mrm{F})$ into the f{i}ltered complex
$(A,{Dec (F)})$, where $Dec(\mrm{F})$ is the ``decalage''
f{i}ltration of $\mrm{F}$ (\cite{DeII} I.3.3), that is also
biregular. This f{i}ltration is def{i}ned as
$$ Dec(\mrm{F})^pA^n=\mrm{Z}_1^{p+n,-p}=\ker\left\{d:\mrm{F}^{p+n}A^{n}\rightarrow\ds\frac{A^{n+1}}{\mrm{F}^{p+n+1}A^{n+1}}\right\}\ .$$
\end{def{i}}

We recall the following result (\cite{DeII} I.3.4).

\begin{lema}\label{PdadesDecalado}\mbox{}\\
Consider $(A,\mrm{F})\in\filt$.
\begin{itemize}
\item[i)] The sequence of inclusions $Z_1^{p+n+1,-p-1}\subseteq
\mrm{F}^{p+n+1}A^n\subseteq B_1^{p+n,-p}\subseteq Z_1^{p+n,-p}$
induces a natural morphism
$$u^{n,p}:E_0^{p,n-p}(Dec(\mrm{F}))=\ds\frac{Z_1^{p+n,-p}}{Z_1^{p+n+1,-p-1}}\longrightarrow
E_1^{p+n,-p}(F)=\frac{Z_1^{p+n,-p}}{B_1^{p+n,-p}}\ .$$
\item[ii)] Given $p$, the morphisms $u^{\ast,p}$ gives rise to a
morphism of $Ch\mc{A}$, natural in $(A,\mrm{F})$
$$u(A,\mrm{F}):E_0^{p,\ast-p}(Dec(\mrm{F}))\rightarrow E_1^{p+\ast,-p}\ .$$
\item[iii)] The morphism $u(A,\mrm{F})$ is a quasi-isomorphism,
that it, it induces isomorphism in cohomology.
\item[iv)] Using the equation $($\ref{E_r+1}$)$ in def{i}nition
\ref{SucEspectral}, we have for all $r\geq 1$ that $u$ induces an
isomorphism of graded complexes
$$E_r(Dec(\mrm{F}))\stackrel{\sim}{\longrightarrow} E_{r+1}(\mrm{F})\ .$$
\end{itemize}
\end{lema}

By \ref{PdadesDecalado} \textit{iv)} for $r=1$ and
\ref{pdadesSucEsp} c) we deduce the following corollary.

\begin{cor}\label{E_2VERSUSE_1}
A morphism $f$ of $\filt$ is an $E_2$-isomorphism if and only if
$Dec(f)$  is a f{i}ltered quasi-isomorphism. In other words
$$\mrm{E}_2=\{f\in\filt \; |\; Dec(f)\in\mrm{E}\}\ .$$
\end{cor}

\begin{prop}\label{CF+E2esDescenso}
Under the notations given in
\ref{SegundaEstructuraDescensoCatF{i}lt},
$(\filt,(\mbf{s},\delta),\mrm{E}_2,\lambda,\mu)$ is an additive
cosimplicial descent category. In addition, $\mu$ is associative
and $\lambda$ is quasi-invertible.
\end{prop}

\begin{proof}
Having into account \ref{E_2VERSUSE_1}, it suf{f}{i}ces to prove
the transfer lemma \ref{FDfuerte}${}^{\mrm{op}}$ for the functor
$Dec:\filt\rightarrow\filt$.\\
Again $\filt$ is additive, so (SDC 1)${}^{\mrm{op}}$ holds. Let us
see (SDC 3)'${}^{\mrm{op}}$, or equivalently, that given
$(A,\mrm{F})\in\filt$, the f{i}ltration
$\delta\mrm{F}$ of $\mbf{s}(A)$ is biregular.\\
By assumption $\mrm{F}$ is a biregular of $A^{i,\ast}$ for all
$i\geq 0$. So f{i}xed $j\geq 0$, there exists $a(i,j),\,
b(i,j)\in\mathbb{Z}$ such that $\mrm{F}^kA^{i,j}=A^{i,j}$ $\forall
k\leq a(i,j)$ and $\mrm{F}^kA^{i,j}=0$ $\forall k\geq b(i,j)$.
Then setting  $\alpha=\min\{a(i,j)+i\; | \; i+j=q\; ; \; i,j\geq
0\}$ and $\beta=\max\{b(i,j)+i\; | \; i+j=q\; ; \; i,j\geq 0\}$ we
have that
$$(\delta\mrm{F})^\alpha(\mbf{s}(A)^n) = \ds\bigoplus_{i+j=n}\mrm{F}^{\alpha-i}A^{i,j}=\mbf{s}(A)^n\;\mbox{ and }\;(\delta\mrm{F})^\beta(\mbf{s}(A)^n) = \ds\bigoplus_{i+j=n}\mrm{F}^{\beta-i}A^{i,j}=0\ .$$

Let us prove (SDC 4$)'{}^{\mrm{op}}$. If
$(Z,\mrm{F})\in\Dl\Dl\filt$, in degree $n$
$\mu(Z):\mbf{s}\Dl\mbf{s}(Z)\rightarrow\mbf{s}\mrm{D}Z$ is the sum
of the morphisms $\mu(Z)_{i,j,q}=Z(d^0\stackrel{j
)}{\cdots}d^0,d^pd^{p-1}\cdots d^{j+1}):Z^{i,j,q}\rightarrow
Z^{p,p,q}$, where $p=i+j$ and $p+q=n$.\\
The f{i}ltration $\delta\Dl\delta\mrm{F}$ of
$\mbf{s}\Dl\mbf{s}(Z)$ if
$(\delta\Dl\delta\mrm{F})^k(\mbf{s}\Dl\mbf{s}(Z))^n=\bigoplus_{i+j+q=n}\mrm{F}^{k-i,k-j}Z^{i,j,q}$,
whereas $(\delta\mrm{D}\mrm{F})^k(\mbf{s}\mrm{D}(Z))^n=\bigoplus_{p+q=n}\mrm{F}^{k-p,k-p}Z^{p,p,q}$.\\
Since $\mrm{F}$ is a decreasing f{i}ltration,
$\mu(Z)_{i,j,q}(\mrm{F}^{k-i,k-j}Z^{i,j,q})\subseteq
\mrm{F}^{k-i,k-j}Z^{p,p,q}\subseteq \mrm{F}^{k-p,k-p}Z^{p,p,q}$,
so $\mu(Z)$ preserves the f{i}ltrations.

Let us see now (SDC 5)'${}^{\mrm{op}}$. If $(A,\mrm{F})\in\filt$,
then $\lambda(A)^n:A^n\rightarrow \mbf{s}(A\times\Dl)^n=A^n\oplus
A^{n-1}\oplus\cdots A^0$ is the inclusion, and
$(\delta(\mrm{F}\times\Dl))^k(\mbf{s}(A\times\Dl))^n=\mrm{F}^kA^n\oplus
\mrm{F}^{k-1}A^{n-1}\oplus\cdots \oplus\mrm{F}^0A^0$, then
$\lambda(A)$ preserves the f{i}ltrations as well.

It is clear that $Dec$ is an additive functor, so (FD
1)${}^{\mrm{op}}$ holds. To f{i}nish the proof it remains to see
(FD 2)${}^{\mrm{op}}$. Let us check the commutativity of the
following diagram (see \cite{DeIII} 8.I.16)
$$\xymatrix@H=4pt@C=30pt@R=20pt{\Dl\filt \ar[r]^-{\Dl Dec}\ar[d]_{(\mbf{s},\delta)} & \Dl\filt \ar[d]^{(\mbf{s},\mbf{s})}\\
                                \filt\ar[r]^{Dec}         & \filt}$$
Let $(A,\mrm{F})\in\filt$. By def{i}nition
$$\begin{array}{l}Dec(\delta\mrm{F})^p(\mbf{s}(A))^n={}_{\delta\mrm{F}}\mrm{Z}_1^{p+n,-p}=\ker\left\{d:(\delta\mrm{F})^{p+n}\mbf{s}(A)^{n}\rightarrow\ds\frac{\mbf{s}(A)^{n+1}}{(\delta\mrm{F})^{p+n+1}\mbf{s}(A)^{n+1}}\right\}=\\
                  =\ker\left\{d:\ds\bigoplus_{i+j=n}\mrm{F}^{p+n-i}A^{i,j}\rightarrow\ds\frac{\ds\bigoplus_{l+s=n+1}A^{l,s}}{\ds\bigoplus_{l+s=n+1}\mrm{F}^{p+n+1-i}A^{l,s}}\right\}\end{array}$$
The restriction of the boundary map $d$ to $A^{i,j}$ is
$d_{A^i}+\sum_k (-1)^{k+j}\partial^k$, where
$d_{A^i}:A^{i,j}\rightarrow A^{i,j+1}$ is the boundary map of the
complex $A^i$, and $\partial^k:A^{i,j}\rightarrow A^{i+1,j}$ is
the $k$-th face map of $A$.\\
Since $\partial^k(\mrm{F}^{p+n-i}A^{i,j})\subseteq
\mrm{F}^{p+n-i}A^{i+1,j}= \mrm{F}^{p+n-(i+1)+1}A^{i+1,j}\subseteq
(\delta\mrm{F})^{p+n+1}\mbf{s}(A)^{n+1}$, then the restriction of
$d$ to $\mrm{F}^{p+n-i}A^{i,j}$ and modulo
$(\delta\mrm{F})^{p+n+1}\mbf{s}(A)^{n+1}$ coincides with
$\bigoplus_i d_{A^i}$. Thus
$$Dec(\delta\mrm{F})^p(\mbf{s}(A))^n=\ds\bigoplus_{i+j=n}\ker\left\{d_{A^i}:\mrm{F}^{p+j}A^{i,j}\rightarrow \ds\frac{A^{i,j+1}}{\mrm{F}^{p+j+1}A^{i,j+1}}\right\}\!=\!\ds\bigoplus_{i+j=n} Dec(\mrm{F})^p(A^{i,j})\, .$$
Therefore
$Dec(\delta\mrm{F})^p(\mbf{s}(A))^n=\mbf{s}(Dec(\mrm{F}))^p\mbf{s}(A)^n$
and $(\mbf{s},\mbf{s})\Dl Dec=Dec(\mbf{s},\delta)$.

F{i}nally, since the image under $Dec$ of a morphism of $\filt$ is
the same morphism between the underlying cochain complexes, it is
clear that $Dec(\lambda(A,\mrm{F}))=\lambda(A,Dec(\mrm{F}))$ if
$(A,\mrm{F})\in\filt$, and
$Dec(\mu(Z,\mrm{F}))=\mu(Z,Dec(\mrm{F}))$ if
$(Z,\mrm{F})\in\Dl\Dl\filt$.
\end{proof}

\begin{cor}
If we denote by ${}_1\filt$ the category $\filt$ with the descent
structure given in \ref{EstructuraDescensoCatF{i}lt} and by
${}_2\filt$ the category $\filt$ with the one given in
\ref{SegundaEstructuraDescensoCatF{i}lt}, then
$Dec:{}_2\filt\rightarrow{}_1\filt$ is a functor of additive
descent categories.
\end{cor}

\begin{cor}\label{propiedadesObjCaminoenCF}
Given morphisms
$A\stackrel{f}{\rightarrow}B\stackrel{g}{\leftarrow}C$ of
f{i}ltered cochain complexes, the \textit{path object} associated
with $f$ and $g$ is a f{i}ltered cochain complex $path(f,g)$,
which is functorial
in $(f,g)$ and such that satisf{i}es the following properties\\[0.1cm]
\textbf{1)} there exists functorial maps in $\filt$
$$ \jmath_A:path(f,g)\rightarrow  A\ \ \jmath_B:path(f,g)\rightarrow B$$
such that $\jmath_A$ $($resp. $\jmath_C)$ is an $E_2$-isomorphism
if and only if $g$ $($resp. $f)$ is so.\\[0.1cm]
\textbf{2)} If $f=g=Id_A$, there exists an $E_2$-isomorphism
$P:A\rightarrow path(A)$ of $\filt$ such that the composition of
$P$ with the projections $\jmath_A,\jmath'_A:path(A)\rightarrow A$
given in 1) is equal to the identity on $A$.\\[0.1cm]
\textbf{3)} The following square commutes up to $E_2$-isomorphism
$$\xymatrix@M=4pt@H=4pt{
 B  & A  \ar[l]_-f\\
 C\ar[u]_g & path(f,g)\ar[u]_{\jmath_A}\ar[l]_-{\jmath_C}\ .}$$
\end{cor}

\begin{obs}
The underlying cochain complex of $path(f,g)$ coincides with one
of the path object given in proposition
\ref{propiedadesObjCaminoenCF}, but they are not the same object
of $\filt$, since now the f{i}ltration of
$\mbf{s}(Path(f\times\Dl,g\times\Dl))$ is the diagonal
f{i}ltration.
\end{obs}

\begin{cor}\label{descripcionCatDerivE2}
The category $Ho_2\filt=\filt[\mrm{E}_2^{-1}]$ is additive. A
morphism $F:X\rightarrow Y$ of $Ho\filt$ is represented by a
zig-zag in the form
$$\xymatrix@M=4pt@H=4pt{X  & T \ar[r]^{{f}} \ar[l]_w &  Y} ,\  w\mbox{ is an $E_2$-isomorphism} \ .$$
Another zig-zag $\xymatrix@M=4pt@H=4pt{X  & S \ar[l]_u
\ar[r]^{{g}} & Y}$ represents $F$ if and only if there exists a
hammock $($commuting in $\filt)$ relating both zig-zags, in the
form
$$
\xymatrix@M=4pt@H=4pt{
                                   &  X                         & T\ar[r]^{{f}}\ar[l]_{w}                 &  Y                        &          \\
 X \ar[ru]^{Id}\ar[rd]_{Id} \ar[r] & \widetilde{X}\ar[d]\ar[u]  & U\ar[d]\ar[u]  \ar[r]^{h} \ar[l]\ar[r]  & \widetilde{Y}\ar[u]\ar[d] & Y\ar[lu]_{Id}\ar[ld]^{Id}\ar[l] \ ,\\
                                   &  X                         & S \ar[r]^{{g}}\ar[l]_{u}                & Y                         &                             }$$
where all maps except $f$, $g$ and $h$ are $E_2$-isomorphisms.
\end{cor}

Again, one can proceed as in proposition
\ref{ComplAcotadoesTriang} to deduce the following

\begin{cor} Let $\mrm{CF}^b\mc{A}$ be the category of bounded-bellow cochain complexes, f{i}ltered by a biregular
f{i}ltration. Then the category localized category
$\mrm{CF}^b\mc{A}[\mrm{E}_2^{-1}]$ of $\mrm{CF}^b\mc{A}$ with
respect to the $E_2$-isomorphisms is a triangulated category.\\
In addition, the ``decalage'' functor induces a functor of
triangulated categories
$$Dec: D^b\mrm{F}\mc{A} \rightarrow \mrm{CF}^b\mc{A}[\mrm{E}_2^{-1}]\ .$$
\end{cor}

All the results given in this section are satisf{i}ed in the case
of decreasing f{i}ltrations instead of decreasing ones. In
particular the following proposition holds.

\begin{prop}\label{F{i}lt+crecient}
If $\mrm{CF}^+\mc{A}$ denotes the category of cochain complexes
f{i}ltered by a biregular increasing f{i}ltration, then
${}_2\filt=(\filt,(\mbf{s},\delta),\mrm{E}_2,\mu,\lambda)$ is an
additive cosimplicial descent category. The diagonal f{i}ltration
$\delta$ of the simple of a cosimplicial f{i}ltered cochain
complex
 $(A,\mrm{W})$ is def{i}ned this time as
$$(\delta\mrm{W})_k(\mbf{s}(A)^n) = \ds\bigoplus_{i+j=n}\mrm{W}_{k+i}A^{i,j}\ .$$
In addition, the functor $Dec:{}_2\filt\rightarrow{}_1\filt$ is a
functor of additive cosimplicial descent categories.
\end{prop}

%

\section{Mixed Hodge Complexes}

In \cite{DeIII} it is introduced the notion of mixed Hodge
complex. The morphisms between such complexes are not given
explicitly, but it can be understood that they are those morphisms
living in the respective (bi)f{i}ltered derived categories that
are compatible with the structural morphisms of the mixed Hodge
complexes
involved.\\
We will def{i}ne in this section a category of mixed Hodge
complexes, and we will endow it with a structure of cosimplicial
descent category using the simple functor developed in
\cite{DeIII}. The homotopy category associated with this
cosimplicial descent category is mapped into the ``category''
appearing in loc. cit..

From now on $\mc{A}$ will denote an abelian category. Before
giving the notion of mixed Hodge complex, we need to introduce the
following preliminaries.

\begin{def{i}}
Given a f{i}ltered complex $(A,\mrm{W})$ of $\filt$, the boundary
map $d:A^i\rightarrow A^{i+1}$ is just a morphism of
$\mrm{F}_f\mc{A}$ (see \ref{DefCatObjF{i}lt}). Then, $d$ is said
to be \textit{strictly compatible with the
f{i}ltration}\index{Index}{strictly compatible with the
f{i}ltration} $\mrm{W}$ if the morphism
$$ A^i / \ker(d) \longrightarrow \mathrm{Im}(d)$$
induced by $d$ is an isomorphism in $\mrm{F}_f\mc{A}$, where $A^i
/ \ker(d)$ and $\mathrm{Im}(d)$ are endowed with the f{i}ltrations
induced by $\mrm{W}$ (cf. \cite{DeII}[I.I]).
\end{def{i}}

\begin{obs} The boundary map of a f{i}ltered complex $(A,\mrm{W})$
is compatible with the f{i}ltration if and only if the spectral
sequence associated with $\mrm{W}$ degenerates at $E_1$
\cite{DeII}[I.3.2].
\end{obs}

\begin{def{i}}[bif{i}ltered complexes]\mbox{}\\
Denote by
$\mrm{CF}^+_2\mc{A}$\index{Symbols}{$\mrm{CF}^+_2\mc{A}$}\index{Index}{complexes!bif{i}ltered}
the category whose objects are triples $(K,\mrm{W},\mrm{F})$,
where

1.- $K$ is a positive cochain complex.

2.- $\mrm{W}$ is an increasing biregular f{i}ltration of $K$ (see
\ref{f{i}ltracionBirregular}).

3.- $\mrm{F}$ is a decreasing biregular f{i}ltration.\\
A morphism $f:(K,\mrm{W},\mrm{F})\rightarrow
(K',\mrm{W}',\mrm{F}')$ of $\mrm{CF}^+_2\mc{A}$ is a morphism of
cochain complexes $f:K\rightarrow K'$ such that
$f:(K,\mrm{W})\rightarrow (K',\mrm{W}')$ and
$f:(K,\mrm{F})\rightarrow (K,\mrm{F}')$ are morphisms of
f{i}ltered complexes.
\end{def{i}}

\begin{num}
Let $k$ be a f{i}eld and $\mc{A}$ be the category of $k$-vector spaces.\\
In order to relax the notations, we will write $Ch^+k$ instead of
$Ch^+\mc{A}$, $\mrm{CF}^+k$ instead of $\filt$ and $\mrm{CF}^+_2k$
instead of $\mrm{CF}^+_2\mc{A}$.
\end{num}

\begin{def{i}}[Hodge complex of weight $n$]\label{ComplHodge}\mbox{}\\
A \textit{Hodge complex of weight} $n$\index{Index}{Hodge
complex!of weight $n$} is the data
$(K_{\mathbb{Q}},(K_{\mathbb{C}},\mrm{F}),\alpha)$, consisting of
\begin{itemize}
\item[a)] A complex $K_{\mathbb{Q}}$ of $Ch^+\mathbb{Q}$ such that
its cohomology $\mathrm{H}^k K_{\mathbb{Q}}$ has f{i}nite
dimension over $\mathbb{Q}$, for all $k$.
\item[b)] A f{i}ltered complex $(K_{\mathbb{C}},\mrm{F})$ in
$\mrm{CF}^+\mathbb{C}$.
\item[c)] $\alpha$ is $(\alpha_0,\alpha_1,\widetilde{K})$, where
$\widetilde{K}$ is in $Ch^+\mathbb{C}$ and $\alpha_i$, $i=0,1$,
are quasi-isomorphisms
$$\xymatrix@M=4pt@H=4pt@C=35pt{ K_{\mathbb{C}} & \widetilde{K} \ar[l]_-{\alpha_0} \ar[r]^-{\alpha_1} & K_{\mathbb{Q}}\otimes \mathbb{C}}\ .$$
\end{itemize}
In addition, the following properties must be satisf{i}ed
\begin{itemize}
\item[(HCI)] The boundary map of $K_\mathbb{C}$ is strictly
compatible with $\mrm{F}$.
\item[(HCII)] For any $k$, the f{i}ltration over $H^k(K_C) =
H^k(K_Q) \otimes C$ induced by $\mrm{F}$ def{i}nes a Hodge
structure on $H^k(K_Q)$ of weight $n+k$.
\end{itemize}
\end{def{i}}

\begin{obs}
(HCII) means that the f{i}ltration $\mrm{F}$ of
$\mrm{H}^k(K_\mathbb{C})$ is $n+k$-opposite to its conjugate
$\overline{\mrm{F}}$, that is, $\mrm{H}^k K_\mathbb{C}$ admits the
following decomposition into a direct sum
$$\mrm{H}^k K_\mathbb{C}=\!\!\!\ds\bigoplus_{p+q=n+k}\!\!\! H^{p,q}\ \mbox{ where }\ \mrm{F}^m(\mrm{H}^k K_\mathbb{C})=\!\ds\bigoplus_{p\geq m}\! H^{p,q}\ \mbox{ and } \ \overline{\mrm{F}}^m(\mrm{H}^k K_\mathbb{C})=\!\ds\bigoplus_{q\geq m}\! H^{p,q}$$
or equivalently, \cite{DeII}[I.2.5]
$${}_{\mrm{F}}\mbf{Gr}_p \left( {}_{\overline{\mrm{F}}}\mbf{Gr}_q (\mrm{H}^k K_\mathbb{C})\right)=0\ \mbox{ if }\ p+q\neq n+k\ .$$
\end{obs}

\noindent In the previous def{i}nition the zig-zag $\alpha$ is in
the form $\cdot\leftarrow\cdot\rightarrow\cdot$, but we can
consider as well any other kind of zig-zag relating
$K_{\mathbb{C}}$ and $K_{\mathbb{Q}}\otimes \mathbb{C}$.

\begin{def{i}} Let $\mc{A}$ be a category endowed with a class of morphisms $\mc{W}$.
A $\mc{W}$-\textit{zig-zag}\index{Index}{$\mc{W}$-zig-zag} of
$\mc{A}$ (or just zig-zag, if $\mc{W}$ is understood) is a pair
$(\underline{A},\underline{w})$ consisting of a family
$\underline{A}=\{A_0,\ldots,A_r\}$ of objects of $\mc{A}$ together
with morphisms $\underline{w}=\{w_0,\ldots ,w_{r-1}\}$ of
$\mc{W}$, such that each $w_i$ is a morphism between $A_i$ and
$A_{i+1}$ (that is, either $w_i:A_i\rightarrow A_{i+1}$ or $w_i:A_{i+1}\rightarrow A_{i}$).\\
In addition, we will assume that two consecutive arrows have
opposite senses, that is, either
$A_{i-1}\stackrel{w_{i-1}}{\rightarrow}A_i
\stackrel{w_i}{\leftarrow} A_{i+1}$ or
$A_{i-1}\stackrel{w_{i-1}}{\leftarrow}A_i
\stackrel{w_i}{\rightarrow}
A_{i+1}$.\\[0.1cm]
\indent Two $\mc{W}$-zig-zags $(\underline{A},\underline{w})$,
$(\underline{B},\underline{v})$ are said to be of the \textit{same
kind} if their respective families of objects have the same
cardinality $r$ and if each $w_i$ has the same sense as $v_i$, for $i=0,\ldots,r-1$.\\[0.1cm]
\indent A morphism between two $\mc{W}$-zig-zags
$(\underline{A},\underline{w})$, $(\underline{B},\underline{v})$
of the same kind is a family of morphisms
$\underline{f}=\{f_i:A_i\rightarrow B_i\}_i$ such that each
diagram involving the maps $f_\cdot$, $w_\cdot$ and $v_\cdot$ is commutative.\\[0.1cm]
\indent A $\mc{W}$-zig-zag between $A$ and $B$ is just a
$\mc{W}$-zig-zag $(\underline{A},\underline{w})$ such that $A_0=A$
and $A_r=B$.\\[0.1cm]
\indent Given objects $A$ and $B$ of $\mc{A}$, the
$\mc{W}$-zig-zags between $A$ and $B$ are the objects of a
category, that will be denoted by
$\mc{R}ist^{\mc{W}}(A,B)$\index{Symbols}{$\mc{R}ist^{\mc{W}}(A,B)$}.
The $\mc{W}$-zig-zags between $A$ and $B$ of the same kind as
$\cdot\leftarrow\cdot\rightarrow\cdot$ gives rise to the full
subcategory of $\mc{R}ist^{\mc{W}}(A,B)$, that will be denoted by
$\mc{R}ist_{red}^{\mc{W}}(A,B)$\index{Symbols}{$\mc{R}ist_{red}^{\mc{W}}(A,B)$}.\\
In addition, if we invert the morphisms of $\mc{W}$ in $A$, then
each $\mc{W}$-zig-zag becomes a morphism, so we have the functors
$$\xymatrix@M=4pt@H=4pt@R=9pt{
 \mc{R}ist_{red}(A,B)\ar[r]^-{\gamma} & \mathrm{Hom}_{\mc{A}[\mc{W}^{-1}]}(A,B) & \mc{R}ist^{\mc{W}}(A,B)\ar[r]^-{\gamma} & \mathrm{Hom}_{\mc{A}[\mc{W}^{-1}]}(A,B)\\
 \mc{R}ist_{red}\ar[r]^-{\gamma}      & Fl(\mc{A}[\mc{W}^{-1}])                 & \mc{R}ist^{\mc{W}}\ar[r]^-{\gamma}      & Fl(\mc{A}[\mc{W}^{-1}])\ .} $$
\end{def{i}}

\begin{num} Let $\mathbf{Quis}$ be the class of quasi-isomorphisms of $Ch^+k$, and $\mathbf{QuisF}$
be the class of f{i}ltered quasi-isomorphisms of $\mrm{CF}^+k$, where $k=\mathbb{Q},\mathbb{C}$.\\
Then, the data $\alpha$ of a Hodge complex of weight $n$ is just
an object of the category
$\mc{R}ist_{red}^{\mathbf{Quis}}(K_{\mathbb{C}},K_{\mathbb{Q}}\otimes
\mathbb{C})$.
\end{num}

\begin{def{i}}
A \textit{generalized} \textit{Hodge complex of weight} $n$
\index{Index}{Hodge complex!of weight $n$ generalized} consists of
$(K_{\mathbb{Q}},(K_{\mathbb{C}},\mrm{F}),\alpha)$, where
$K_{\mathbb{Q}}$ and $(K_{\mathbb{C}},\mrm{F})$ satisf{i}es
conditions a), b) (HCI) and (HCII) in def{i}nition
\ref{ComplHodge}, whereas $\alpha$ is a $\mathbf{Quis}$-zig-zag
between $K_{\mathbb{C}}$ and $K_{\mathbb{Q}}\otimes \mathbb{C}$ of
$Ch^+\mathbb{C}$.
\end{def{i}}

Due to properties 1) and 3) in proposition
\ref{pddesCaminoenComplejos}, we can associate a Hodge complex of
weight $n$ to a generalized Hodge complex of weight $n$ in a
functorial way.

\begin{lema}\label{RistrRedQuis} If $\mc{A}$ is an abelian category, there exists a functor\index{Symbols}{$\underline{red}$}
$$\underline{red}:\mc{R}ist^{\mathbf{Quis}}\longrightarrow \mc{R}ist_{red}^{\mathbf{Quis}}$$
such that the composition $\xymatrix@M=4pt@H=4pt@C=30pt{
 \mc{R}ist^{\mathbf{Quis}}\ar[r]^{\underline{red}}& \mc{R}ist_{red}^{\mathbf{Quis}} \ar[r]^-{\gamma} & Fl(Ho \,Ch^+\mc{A})}$
is just $\gamma:\mc{R}ist^{\mathbf{Quis}}\longrightarrow
Fl(Ho \,Ch^+\mc{A})$.\\
In addition, if $A$, $B$ are objects of $Ch^+\mc{A}$, the functor
$\underline{red}$ restrict to
$$\underline{red}:\mc{R}ist^{\mathbf{Quis}}(A,B)\longrightarrow \mc{R}ist_{red}^{\mathbf{Quis}}(A,B)$$
We will say that $\underline{red}(\underline{A},\underline{w})$ is
the \textit{reduced zig-zag}\index{Index}{reduced zig-zag}
associated with $(\underline{A},\underline{w})$.
\end{lema}

\begin{proof}
Given a $\mathbf{Quis}$-zig-zag $R=(\underline{A},\underline{w})$
between $A$ and $B$ in the form
$$\xymatrix@M=4pt@H=4pt{A=A_0\ar@{-}[r]^{w_0} & A_1 \ar@{-}[r]^{w_1} & \cdots \ar@{-}[r]^{w_{r-2}} & A_{r-1}\ar@{-}[r]^{w_{r-1}} & A_r}$$
its associated reduced zig-zag is obtained through the following
procedure.\\
We take the f{i}rst pair of consecutive arrows in $R$ of the form
$A_{i-1}\stackrel{w_{i-1}}{\rightarrow}A_i
\stackrel{w_i}{\leftarrow} A_{i+1}$. If there is no such pair of
consecutive arrows in $R$, then this zig-zag is already a zig-zag
in $\mc{R}ist_{red}^{\mathbf{Quis}}(A,B)$, and we def{i}ne
$\underline{red}(R)=R$.\\
If there exists such $w_i,w_{i-1}$, properties 1) and 2) of
proposition \ref{pddesCaminoenComplejos} provide the following
square, commutative up to homotopy,
$$\xymatrix@M=4pt@H=4pt{
 A_i  & A_{i+1}  \ar[l]_-{w_{i}}\\
 A_{i-1}\ar[u]_{w_{i-1}} & path(w_i,w_{i-1})\ar[u]_{\jmath_{A_{i+1}}}\ar[l]_-{\jmath_{A_{i-1}}}\ ,}$$
where $\jmath_{A_{i+1}}$ and $\jmath_{A_{i-1}}$ are
quasi-isomorphisms, that is, morphisms of $\mathbf{Quis}$. Hence,
replacing in $R$ the maps
$A_{i-1}\stackrel{w_{i-1}}{\rightarrow}A_i
\stackrel{w_i}{\leftarrow} A_{i+1}$ with
$A_{i-1}\stackrel{\jmath_{A_{i-1}}}{\leftarrow} path(w_i,w_{i-1})
\stackrel{\jmath_{A_{i+1}}}{\rightarrow} A_{i+1}$ and composing
maps we obtain a new $\mathbf{Quis}$-zig-zag $\widetilde{R}$
between $A$ and $B$ of length strictly smaller than the length of $R$.\\
Since $w_i\comp \jmath_{A_{i+1}}$ is homotopic to $w_{i-1}\comp
\jmath_{A_{i-1}}$, then $\gamma(R)$ and $\gamma(\widetilde{R})$
coincides in $Ho \,Ch^+\mc{A}$.\\
Moreover, the mapping $R\rightarrow \widetilde{R}$ def{i}nes a
functor $\mc{R}ist^{\mathbf{Quis}}(A,B)\rightarrow
\mc{R}ist^{\mathbf{Quis}}(A,B)$. Indeed, if we have a commutative
diagram in $Ch^+\mc{A}$
$$\xymatrix@M=4pt@H=4pt{
 A_{i-1}\ar[r]^-{w_{i-1}}\ar[d]_{f_{i-1}} & A_i \ar[d]_{f_{i}} & A_{i+1}  \ar[l]_-{w_{i}}\ar[d]_{f_{i+1}}\\
 B_{i-1}\ar[r]^-{v_{i-1}} & B_i  & B_{i+1}  \ar[l]_-{v_{i}}}$$
then from the functoriality of $path$ it follows the existence of
a morphism $\widetilde{f}$ that f{i}ts into the commutative
diagram
$$\xymatrix@M=4pt@H=4pt@C=30pt{
 A_{i-1}\ar[d]_{f_{i-1}} & path(w_i,w_{i-1}) \ar[r]^-{\jmath_{A_{i+1}}}\ar[d]_{\widetilde{f}}\ar[l]_-{\jmath_{A_{i+1}}} & A_{i+1}  \ar[d]_{f_{i+1}}\\
 B_{i-1}                 & path(v_i,v_{i-1}) \ar[r]^-{\jmath_{B_{i+1}}}\ar[l]_-{\jmath_{B_{i+1}}}                       & B_{i+1}  \ .}$$
Therefore, it suf{f}{i}ces to iterate this procedure until we get
the desired zig-zag
$\underline{red}(\underline{A},\underline{w})$.
\end{proof}

\begin{obs} Note that the reduced zig-zag associated with $(\underline{A},\underline{w})$
not only preserves the morphism in $Ho\,\cdc$ represented by this
zig-zag. In addition, the original and reduced zig-zags are in
some sense ``homotopic''.
\end{obs}

\begin{cor} Each generalized Hodge complex of weight $n$ gives rise to a  Hodge complex of weight $n$, just by replacing $\alpha$
with $\underline{red}(\alpha)$.
\end{cor}

Next we recall the notion of mixed Hodge complex, and introduce a
category consisting of these complexes.

\begin{num}\label{GrConmutaconTensor} The tensor product over $\mathbb{C}$,
$-\otimes\mathbb{C}:\mathbb{Q}$-vector spaces
$\rightarrow\mathbb{C}$-vector spaces, is an exact functor. Thus
it induces
$$ -\otimes\mathbb{C} : \mrm{CF}^+\mathbb{Q}\rightarrow \mrm{CF}^+\mathbb{C}$$
in such a way that the functor $\mbf{Gr}_n$ commutes with
$-\otimes\mathbb{C}$.
\end{num}

\begin{def{i}}[\textbf{Mixed Hodge Complex}]\label{Def{i}ComplHodgeMixto}\mbox{}\\
A \textit{mixed Hodge complex}\index{Index}{Hodge complex!mixed}
consists of the data
$((K_\mathbb{Q},\mrm{W}),(K_\mathbb{C},\mrm{W},\mrm{F}),\alpha)$,
where
\begin{itemize}
\item[a)] $(K_\mathbb{Q},\mrm{W})$ is a cochain complex of
$\mathbb{Q}$-vector spaces, f{i}ltered by the increasing
f{i}ltration $\mrm{W}$. In other words, $(K_\mathbb{Q},\mrm{W})$
is an object of $\mrm{CF}^+\mathbb{Q}$. In addition,
$\mrm{H}^kK_\mathbb{Q}$ has f{i}nite dimension over $\mathbb{Q}$
for all $k$.
\item[b)] $(K_\mathbb{C},\mrm{W},\mrm{F})$ is an object of
$\mrm{CF}^+_2\mathbb{C}$.
\item[c)] $\alpha$ is the data
$(\alpha_0,\alpha_1,(\widetilde{K},\widetilde{\mrm{W}}))$, where
$(\widetilde{K},\widetilde{\mrm{W}})$ is an object of
$\mrm{CF}^+\mathbb{C}$ and $\alpha_i$, $i=0,1$, is a f{i}ltered
quasi-isomorphism (see \ref{Def{i}QuisF{i}lt})
$$\xymatrix@M=4pt@H=4pt@C=35pt{ (K_\mathbb{C},\mrm{W}) & (\widetilde{K},\widetilde{\mrm{W}}) \ar[l]_-{\alpha_0} \ar[r]^-{\alpha_1} & (K_\mathbb{Q},\mrm{W})\otimes\mathbb{C} }\ .$$
\end{itemize}
The following axiom must be satisf{i}ed\\
(MHC) For each $n$,
$({}_\mrm{W}\mbf{Gr}_nK_\mathbb{Q},({}_\mrm{W}\mbf{Gr}_nK_\mathbb{C},\mrm{F}),\mbf{Gr}_n(\alpha))$
is a Hodge complex of weight $n$, where $\mbf{Gr}_n(\alpha)$
denotes the zig-zag
$$\xymatrix@M=4pt@H=4pt@C=35pt{ {}_\mrm{W}\mbf{Gr}_nK_\mathbb{C} & {}_{\widetilde{\mrm{W}}}\mbf{Gr}_n\widetilde{K} \ar[l]_-{\mbf{Gr}_n\alpha_0} \ar[r]^-{\mbf{Gr}_n\alpha_1} &  {}_{\mrm{W}\otimes\mathbb{C}}\mbf{Gr}_n\left(K_\mathbb{Q}\otimes \mathbb{C}\right) \stackrel{(\ref{GrConmutaconTensor})}{\simeq}   \left({}_\mrm{W}\mbf{Gr}_nK_\mathbb{Q}\right)\otimes \mathbb{C}}\ .$$
\end{def{i}}

\begin{obs}
Again, the data $\alpha$ of a mixed Hodge complex is just an
object of the category
$\mc{R}ist_{red}^{\mathbf{Quis}}((K_{\mathbb{C}},\mrm{W}),(K_{\mathbb{Q}},\mrm{W})\otimes
\mathbb{C})$, where $\mathbf{QuisF}$ is the class of the
f{i}ltered quasi-isomorphisms of $\mrm{CF}^+\mathbb{C}$.
\end{obs}

Analogously to the case of Hodge complexes of weight $n$, we can
consider any kind of zig-zag to def{i}ne the data $\alpha$ of a
mixed Hodge complex.

\begin{def{i}}
A \textit{generalized mixed Hodge complex}\index{Index}{Hodge
complex!mixed and generalized} consists of\\
$((K_\mathbb{Q},\mrm{W}),(K_\mathbb{C},\mrm{W},\mrm{F}),\alpha)$,
where $(K_\mathbb{Q},\mrm{W})$ and
$(K_\mathbb{C},\mrm{W},\mrm{F})$ satisf{i}es conditions a) and b)
of
mixed Hodge complex.\\
As before, $\alpha$ is a $\mathbf{QuisF}$-zig-zag between
$(K_{\mathbb{C}},\mrm{W})$ and $(K_{\mathbb{Q}},\mrm{W})\otimes
\mathbb{C}$ in $\mrm{CF}^+\mathbb{C}$.
In addition, the following axiom must be satisf{i}ed\\[0.1cm]
(MHC) For each $n$,
$({}_\mrm{W}\mbf{Gr}_nK_\mathbb{Q},({}_\mrm{W}\mbf{Gr}_nK_\mathbb{C},\mrm{F}),\mbf{Gr}_n(\alpha))$
is a generalized Hodge complex of weight $n$ (where
$\mbf{Gr}_n(\alpha)$ is def{i}ned analogously).
\end{def{i}}

Similarly to \ref{RistrRedQuis}, properties 1) and 3) of the
functor $path$ given in proposition \ref{propiedadesObjCaminoenCF}
allows us to associate a mixed Hodge complex to any generalized
mixed Hodge complex in a functorial way.

\begin{lema} If $\mc{A}$ is an abelian category, there exists a functor\index{Symbols}{$\underline{red}$}
$$\underline{red}:\mc{R}ist^{\mathbf{QuisF}}\longrightarrow \mc{R}ist_{red}^{\mathbf{QuisF}}$$
such that the composition $\xymatrix@M=4pt@H=4pt@C=30pt{
 \mc{R}ist^{\mathbf{QuisF}}\ar[r]^{\underline{red}}& \mc{R}ist_{red}^{\mathbf{QuisF}} \ar[r]^-{\gamma} & Fl(Ho \filt)}$
is just $\gamma:\mc{R}ist^{\mathbf{QuisF}}\longrightarrow
Fl(Ho \filt)$.\\
In addition, if $(A,\mrm{F})$, $(B,\mrm{G})$ are objects of
$\filt$, functor $\underline{red}$ restricts to
$$\underline{red}:\mc{R}ist^{\mathbf{QuisF}}((A,\mrm{F}),(B,\mrm{G}))\longrightarrow \mc{R}ist_{red}^{\mathbf{QuisF}}((A,\mrm{F}),(B,\mrm{G}))$$
We will refer to
$\underline{red}((\underline{A},\underline{\mrm{F}}),\underline{w})$
as the \textit{reduced zig-zag}\index{Index}{reduced zig-zag}
associated with
$((\underline{A},\underline{\mrm{F}}),\underline{w})$.
\end{lema}

\begin{obs} The square appearing in property 3) of proposition
\ref{propiedadesObjCaminoenCF} commutes up to f{i}ltered homotopy,
so the reduced zig-zag associated with
$((\underline{A},\underline{\mrm{F}}),\underline{w})$ is
``homotopy equivalent'' (in some sense) to $((\underline{A},\underline{\mrm{F}}),\underline{w})$.\\
In addition, this is a constructive procedure, consisting just in
iterate functor $path$ (in an ordered way).
\end{obs}

\begin{cor} Each generalized mixed Hodge complex gives rise to a mixed Hodge complex in a functorial way by replacing $\alpha$
by $\underline{red}(\alpha)$.
\end{cor}

\begin{ej} \cite{DeIII}, 8.I.8
Let $j:U\rightarrow X$ be an open immersion  of smooth varieties
(here variety means a separated, reduced and of f{i}nite type
$\mathbb{C}$-scheme). Assume that $X$ is proper and $Y=X\backslash
U$ is a normal crossing divisor.\\
Let $(Rj_{\ast}\mathbb{Q},\mrm{W})$ be the f{i}ltered complex of
sheaves of $\mathbb{Q}$-vector spaces on $X$, where
$\mrm{W}=\tau_{\leq}$ is the ``canonical'' f{i}ltration. That is
to say, $\tau_{\leq p}Rj_{\ast}\mathbb{Q}$ is given in degree $n$
by $Rj_{\ast}\mathbb{Q}$
if $n<p$, $\mrm{Ker}d$ if $p=n$ and 0 otherwise.\\
Let $(\Omega_X\langle Y\rangle,\mrm{W},\mrm{F})$ be the
logarithmic De Rham complex of $X$ along $Y$ \cite{DeII} 3.I. The
f{i}ltration $\mrm{W}$ is the so-called ``weight f{i}ltration'',
consisting in f{i}ltering by the order of poles in
$\Omega_X\langle Y\rangle$. The f{i}ltration $\mrm{F}$, called
``Hodge f{i}ltration'', is just the f{i}ltration ``b{\^e}te''
associated with $\Omega_X\langle Y\rangle$, that is,
$\mrm{F}^n\Omega^p_X\langle Y\rangle=\Omega^p_X\langle Y\rangle$
if $p\geq n$ and $0$ otherwise.\\
Then, there exists a zig-zag $\alpha$ of f{i}ltered
quasi-isomorphisms such that
$$(R\Gamma (j_{\ast}\mathbb{Q},\mrm{W}),R\Gamma
(\Omega_X\langle Y\rangle,\mrm{W},\mrm{F}),\alpha)$$ is a mixed Hodge complex.\\
The zig-zag $\alpha$ involves the result \cite{DeII} 3.I.8 that
relates $\Omega_X\langle Y\rangle$ to $j_\ast \Omega_U$, together
with Poincar{\'e} lemma (that is, $\Omega_U$ is a resolution of the
constant sheaf $\mathbb{C}$), and together with Godement resolutions.\\
The zig-zag $\alpha=\alpha(U,X)$ is natural in $(U,X)$. Moreover,
since Godement resolutions are functorial, the natural
transformations that gives rise to $\alpha(U,X)$ given in
\cite{DeII} has values in the category of f{i}ltered complexes
instead of in the f{i}ltered derived category (cf.
\cite{Be}, 4 or \cite{H}, {8.2}).\\
It should be pointed out that the zig-zag $\alpha$ is also
considered as a zig-zag of length 2 in \cite{H}, using a dual
procedure to the one given here, that is called ``quasi-pushout''
in loc. cit..
\end{ej}

\begin{def{i}}[\textbf{Category of mixed Hodge complexes}]\mbox{}\\
Let $\mc{H}dg$ be the category whose objects are the mixed Hodge
complexes, and whose morphisms are def{i}ned as follows.\\
A morphism $f\!=\!(f_\mathbb{Q},f_\mathbb{C},\widetilde{f})\!
:\!((K_\mathbb{Q},\mrm{W}),(K_\mathbb{C},\mrm{W},\mrm{F}),\alpha)\rightarrow
((K'_\mathbb{Q},\mrm{W}'),(K'_\mathbb{C},\mrm{W}',\mrm{F}'),\alpha')$
consists of
\begin{itemize}
\item[1)] A morphism
$f_\mathbb{Q}:(K_\mathbb{Q},\mrm{W})\rightarrow
(K'_\mathbb{Q},\mrm{W}')$ of $\mrm{CF}^+\mathbb{Q}$.
\item[2)] A morphism
$f_\mathbb{C}:(K_\mathbb{C},\mrm{W},\mrm{F})\rightarrow
(K'_\mathbb{C},\mrm{W}',\mrm{F}')$ of $\mrm{CF}_2^+\mathbb{C}$.
\item[3)] If $\alpha$ and $\alpha'$ are the respective zig-zags
$$\xymatrix@M=4pt@H=4pt@C=32pt@R=9pt{ (K_\mathbb{C},\mrm{W}) & (\widetilde{K},\widetilde{\mrm{W}}) \ar[l]_-{\alpha_0} \ar[r]^-{\alpha_1} & (K_\mathbb{Q},\mrm{W})\otimes\mathbb{C} \\
 (K'_\mathbb{C},\mrm{W}') & (\widetilde{K}',\widetilde{\mrm{W}}') \ar[l]_-{\alpha'_0} \ar[r]^-{\alpha'_1} & (K'_\mathbb{Q},\mrm{W}')\otimes\mathbb{C}}$$
then $\widetilde{f}:(\widetilde{K},\widetilde{\mrm{W}})\rightarrow
(\widetilde{K}',\widetilde{\mrm{W}}')$ is a morphism of
$\mrm{CF}_2^+\mathbb{C}$ such that the squares I and II of diagram
\begin{equation}\label{MorfismHodge}\xymatrix@H=4pt@C=30pt@R=35pt{(K_\mathbb{C},\mrm{W}) \ar[d]_{f_{\mathbb{C}}} & (\widetilde{K},\widetilde{\mrm{W}}) \ar[l]_-{\alpha_0} \ar[r]^-{\alpha_1}  \ar[d]_{\widetilde{f}}   & (K_\mathbb{Q},\mrm{W})\otimes\mathbb{C}  \ar[d]^{f_{\mathbb{Q}}\otimes\mathbb{C}} \\
                                (K'_\mathbb{C},\mrm{W}')\ar@{}[ru]|{\mrm{I}}   & (\widetilde{K}',\widetilde{\mrm{W}}') \ar[l]_-{\alpha'_0} \ar[r]^-{\alpha'_1} \ar@{}[ru]|{\mrm{II}} &
                                (K'_\mathbb{Q},\mrm{W}')\otimes\mathbb{C}}\end{equation}
\textbf{commutes} in $\mrm{CF}_2^+\mathbb{C}$. For a similar
definition see \cite{Be}, 3.
\end{itemize}
\end{def{i}}

\begin{def{i}} The category of generalized mixed Hodge complexes,
$\mc{H}dg_{\mc{G}}$, is def{i}ned analogously.\\
A morphism $f\!=\!(f_\mathbb{Q},f_\mathbb{C},\widetilde{f})\!
:\!((K_\mathbb{Q},\mrm{W}),(K_\mathbb{C},\mrm{W},\mrm{F}),\alpha)\rightarrow
((K'_\mathbb{Q},\mrm{W}'),(K'_\mathbb{C},\mrm{W}',\mrm{F}'),\alpha')$
between two generalized mixed Hodge complexes such that $\alpha$
and $\alpha'$ are of the same kind, consists of
\begin{itemize}
\item[1)] A morphism
$f_\mathbb{Q}:(K_\mathbb{Q},\mrm{W})\rightarrow
(K'_\mathbb{Q},\mrm{W}')$ of $\mrm{CF}^+\mathbb{Q}$.
\item[2)] A morphism
$f_\mathbb{C}:(K_\mathbb{C},\mrm{W},\mrm{F})\rightarrow
(K'_\mathbb{C},\mrm{W}',\mrm{F}')$ of $\mrm{CF}_2^+\mathbb{C}$.
\item[3)] A morphism $\widetilde{f}$ between $\alpha$ and
$\alpha'$ in the category $\mc{R}ist^{\mathbf{QuisF}}$ of zig-zags
of f{i}ltered quasi-isomorphisms in $\mrm{CF}^+\mathbb{C}$.
\end{itemize}
\end{def{i}}

\begin{cor} The functor ``reduced zig-zag'' gives rise to a functor
$$\underline{red}:\mc{H}dg_{\mc{G}}\longrightarrow \mc{H}dg $$
that maps the generalized mixed Hodge complex
$((K_\mathbb{Q},\mrm{W}),(K_\mathbb{C},\mrm{W},\mrm{F}),\alpha)$
into the mixed Hodge complex
$((K_\mathbb{Q},\mrm{W}),(K_\mathbb{C},\mrm{W},\mrm{F}),\underline{red}\,\alpha)$.\\
Moreover, the zig-zags $\alpha$ and $\underline{red}\,\alpha$
def{i}ne the same morphism of $Ho\mrm{CF}^+\mathbb{Q}$ and, in
addition, they are ``homotopy equivalent''.
\end{cor}

\begin{proof}
The functoriality of
$\underline{red}:\mc{H}dg_{\mc{G}}\longrightarrow \mc{H}dg$ is clear.\\
If $f\!=\!(f_\mathbb{Q},f_\mathbb{C},\widetilde{f})\!
:\!((K_\mathbb{Q},\mrm{W}),(K_\mathbb{C},\mrm{W},\mrm{F}),\alpha)\rightarrow
((K'_\mathbb{Q},\mrm{W}'),(K'_\mathbb{C},\mrm{W}',\mrm{F}'),\alpha')$,
is a morphism in $\mc{H}dg_{\mc{G}}$ then
$(f_\mathbb{Q},f_\mathbb{C},\underline{red}\widetilde{f})$ is a
morphism in $\mc{H}dg$, because of the functoriality of
$\underline{red}:\mc{R}ist^{\mathbf{QuisF}}\longrightarrow
\mc{R}ist_{red}^{\mathbf{QuisF}}$.
\end{proof}

\begin{obs} Assume given
$$f\!=\!(f_\mathbb{Q},f_\mathbb{C},\widetilde{f})\!
:\!((K_\mathbb{Q},\mrm{W}),(K_\mathbb{C},\mrm{W},\mrm{F}),\alpha)\rightarrow
((K'_\mathbb{Q},\mrm{W}'),(K'_\mathbb{C},\mrm{W}',\mrm{F}'),\alpha')$$
such that $f_\mathbb{Q}:(K_\mathbb{Q},\mrm{W})\rightarrow
(K'_\mathbb{Q},\mrm{W}')$ and
 $f_\mathbb{C}:(K_\mathbb{C},\mrm{W},\mrm{F})\rightarrow
(K'_\mathbb{C},\mrm{W}',\mrm{F}')$.
Assume also that $\alpha$ and $\alpha'$ consists of the respective
zig-zags
$$\xymatrix@M=4pt@H=4pt@C=32pt@R=9pt{ (K_\mathbb{C},\mrm{W}) & (\widetilde{K},\widetilde{\mrm{W}}) \ar[l]_-{\alpha_0} \ar[r]^-{\alpha_1} & (K_\mathbb{Q},\mrm{W})\otimes\mathbb{C} \\
 (K'_\mathbb{C},\mrm{W}') & (\widetilde{K}',\widetilde{\mrm{W}}') \ar[l]_-{\alpha'_0} \ar[r]^-{\alpha'_1} & (K'_\mathbb{Q},\mrm{W}')\otimes\mathbb{C}}$$
and that
$\widetilde{f}:(\widetilde{K},\widetilde{\mrm{W}})\rightarrow
(\widetilde{K}',\widetilde{\mrm{W}}')$ is a morphism of
$\mrm{CF}_2^+\mathbb{C}$ such that the squares I and II of diagram
$$\xymatrix@H=4pt@C=30pt@R=35pt{(K_\mathbb{C},\mrm{W}) \ar[d]_{f_{\mathbb{C}}} & (\widetilde{K},\widetilde{\mrm{W}}) \ar[l]_-{\alpha_0} \ar[r]^-{\alpha_1}  \ar[d]_{\widetilde{f}}   & (K_\mathbb{Q},\mrm{W})\otimes\mathbb{C}  \ar[d]^{f_{\mathbb{Q}}\otimes\mathbb{C}} \\
                                (K'_\mathbb{C},\mrm{W}')\ar@{}[ru]|{\mrm{I}}   & (\widetilde{K}',\widetilde{\mrm{W}}') \ar[l]_-{\alpha'_0} \ar[r]^-{\alpha'_1} \ar@{}[ru]|{\mrm{II}} & (K'_\mathbb{Q},\mrm{W}')\otimes\mathbb{C}}$$
commutes \textbf{up to f{i}ltered homotopy} in $\mrm{CF}_2^+\mathbb{C}$.\\
Let $cyl(\widetilde{K},\widetilde{\mrm{W}})$ be the ``classical''
cylinder object in the category
$\mrm{CF}_2^+\mathbb{C}=C(\mrm{F}_f\mathbb{C})$ (see
\ref{Def{i}HomotF{i}lt}), and
$i,j:(\widetilde{K},\widetilde{\mrm{W}})\rightarrow
cyl(\widetilde{K},\widetilde{\mrm{W}})$ be the canonical inclusions.\\
Recall that $f_{\mathbb{C}}\comp\alpha_0$ is homotopic to
$\alpha'_0\comp \widetilde{f}$ in $\mrm{CF}_2^+\mathbb{C}$ if and
only if there exists a homotopy
$H:cyl(\widetilde{K},\widetilde{\mrm{W}})\rightarrow
(\widetilde{K},\widetilde{\mrm{W}})$ giving rise to the following
morphism of $\mathbf{QuisF}$-zig-zags in $\mrm{CF}_2^+\mathbb{C}$
$$\xymatrix@H=4pt@M=4pt@C=30pt{
 (K_\mathbb{C},\mrm{W}) \ar[d]_{f_{\mathbb{C}}} & (\widetilde{K},\widetilde{\mrm{W}})\ar[l]_-{\alpha_0}\ar[d]_{f_{\mathbb{C}}\comp\alpha_0}\ar[r]^-{i} & cyl(\widetilde{K},\widetilde{\mrm{W}})\ar[d]_H  & (\widetilde{K},\widetilde{\mrm{W}})\ar[l]_-{j}\ar[d]_{\widetilde{f}} \\
  (K'_\mathbb{C},\mrm{W}')                      & (K'_\mathbb{C},\mrm{W}') \ar[l]_-{Id}    \ar[r]^-{Id}                                                & (K'_\mathbb{C},\mrm{W}')                        & (\widetilde{K}',\widetilde{\mrm{W}}')\ar[l]_-{\alpha'_0}}$$
One can argue in a similar way with square II, obtaining a
morphism of $\mc{H}dg_{\mc{G}}$. Therefore, we obtain in this way
a morphism in  $\mc{H}dg$ between the corresponding mixed
Hodge complexes.\\
This mapping in not functorial at all in the data
$(f_\mathbb{Q},f_\mathbb{C},\widetilde{f})$, since it depends on
chosen the homotopy  for the squares I and II.\\
Two dif{f}erent choices of homotopies for I and II provides two
morphisms of $\mc{H}dg$, that are no related in general.\\[0.2cm]
\end{obs}

\begin{obs} In \cite{PS} another def{i}nition of category of mixed Hodge
complexes is considered, in which a morphism is such that the
corresponding diagram (\ref{MorfismHodge}) commutes up to
homotopy. In this case some pathologies appear, for instance the
non-functoriality of the cone associated with a morphism of mixed
Hodge complexes (see loc. cit. 3.23).
\end{obs}

Next we endow $\mc{H}dg$ with a structure of cosimplicial descent
category, in which the simple functor
${\mbf{s}}_{\mc{H}dg}=(\mbf{s},\delta,\mbf{s}):\Dl\mc{H}dg\rightarrow\mc{H}dg$
is the one given in \cite{DeIII} 8.I.15.

\begin{obs}\label{SimpleConmutaTensor} Note that the simple functor
$(\mbf{s},\delta):\Dl\mrm{CF}^{+}\mathbb{Q}\rightarrow\mrm{CF}^{+}\mathbb{Q}$
(see \ref{F{i}lt+crecient}) commutes with $-\otimes\mathbb{C}$,
since the tensor product with $\mathbb{C}$ commutes with f{i}nite
sums.
\end{obs}

\begin{def{i}}[Descent structure on $\mc{H}dg$]\mbox{}\\
\textbf{Simple functor:}
Given a cosimplicial mixed Hodge complex
${K}=((K_\mathbb{Q},\mrm{W}),(K_\mathbb{C},\mrm{W},\mrm{F}),\alpha)$,
let ${\mbf{s}}_{\mc{H}dg}{K}$ be the mixed Hodge complex
$((\mbf{s} K_\mathbb{Q} ,\delta\mrm{W}),(\mbf{s}
K_\mathbb{C},\delta\mrm{W},\mbf{s}\mrm{F}),\mbf{s}\alpha)$, where
$\mbf{s}$ denotes the usual simple of cochain complexes and
$\delta\mrm{W}$ is def{i}ned as in \ref{F{i}lt+crecient}. More
concretely
$$\begin{array}{cl}
 \mbf{s}(K_{-})^n=\ds\bigoplus_{p+q=n} K_{-}^{p,q} \ ;\ \ (\delta\mrm{W})_k(\mbf{s}(K_{-})^n) = \ds\bigoplus_{i+j=n}\mrm{W}_{k+i}K_{-}^{i,j} \ , & \mbox{ if }-\mbox{ is }\mathbb{Q}\mbox{ or }\mathbb{C}\\
 (\mbf{s}(\mrm{F}))^k (\mbf{s}K_{\mathbb{C}})^n = \ds\bigoplus_{p+q=n} \mrm{F}^k K_\mathbb{C}^{p,q} \ .& {}
\end{array}$$
F{i}nally, if
$\alpha=(\alpha_0,\alpha_1,(\widetilde{K},\widetilde{\mrm{W}}))$
then $\mbf{s}\alpha$ denotes the zig-zag
\begin{equation}\label{Def{i}SimplAlpha}\xymatrix@M=4pt@H=4pt@C=35pt{
 (\mbf{s} K_\mathbb{C},\delta\mrm{W}) & (\mbf{s}\widetilde{K},\delta\widetilde{\mrm{W}}) \ar[l]_-{\mbf{s}\alpha_0} \ar[r]^-{\mbf{s}\alpha_1} & (\mbf{s} (K_\mathbb{Q}\otimes\mathbb{C}),\delta(\mrm{W}\otimes\mathbb{C}))\stackrel{\ref{SimpleConmutaTensor}}{\simeq }(\mbf{s} K_\mathbb{Q},\delta\mrm{W})\otimes \mathbb{C}}\ .
\end{equation}
\textbf{Equivalences:} the class of equivalences is def{i}ned as
$${\mrm{E}}_{\mc{H}dg}=\{(f_\mathbb{Q},f_\mathbb{C},\widetilde{f})\; |\; f_\mathbb{Q}\mbox{ is a quasi-isomorphism in }{Ch}^+\mathbb{Q}\}\ .$$
\textbf{Transformation $\mathbf{\lambda}$:}
$\lambda^{\mc{H}dg}:Id_{\mc{H}dg}\rightarrow
{\mbf{s}}_{\mc{H}dg}(-\times\Dl)$ is
${\lambda}^{\mc{H}dg}_K=(\lambda^{\mathbb{Q}}_{K_{\mathbb{Q}}},\lambda^{\mathbb{C}}_{K_{\mathbb{C}}},\lambda^{\mathbb{C}}_{\widetilde{K}})$
induced by the transformations $\lambda^\mathbb{Q}$ and $\lambda^{\mathbb{C}}$ of $Ch^+\mathbb{Q}$ and $Ch^+\mathbb{C}$ respectively.\\[0.2cm]
\textbf{Transformation $\mathbf{\mu}$:}
similarly, the transformation $\mu^{\mc{H}dg}_K:\mbf{s}_{\mc{H}dg}
\Dl \mbf{s}_{\mc{H}dg}\rightarrow \mbf{s}_{\mc{H}dg}\mrm{D}$ is
${\mu}^{\mc{H}dg}=(\mu^{\mathbb{Q}}_{K_{\mathbb{Q}}},\mu^{\mathbb{C}}_{K_{\mathbb{C}}},\mu^{\mathbb{C}}_{\widetilde{K}})$
where $\mu_{\mathbb{Q}}$, $\mu_{\mathbb{C}}$ and $\widetilde{\mu}$
are the usual natural transformations of $Ch^+\mathbb{Q}$ and
$Ch^+\mathbb{C}$ respectively.
\end{def{i}}

\begin{thm}\label{Hodge+E_2}
The category
$(\mc{H}dg,{\mbf{s}_{\mc{H}dg}},{\mrm{E}_{\mc{H}dg}},\mu_{\mc{H}dg},\lambda_{\mc{H}dg})$
is an additive cosimplicial descent category.\\
In addition, the forgetful functor $\mrm{U}:\mc{H}dg\rightarrow
Ch^+\mathbb{Q}$ given by
$\mrm{U}((K_\mathbb{Q},\mrm{W}),(K_\mathbb{C},\mrm{W},\mrm{F}),\alpha)=K_\mathbb{Q}$
is a functor of additive cosimplicial descent categories.
\end{thm}

By proposition \ref{F{i}lt+crecient}, the simple functor
$(\mbf{s},\delta):\Dl\mrm{CF}^+\mc{A}\rightarrow \mrm{CF}^+\mc{A}$
preserves $E_2$-isomorphisms. In addition, it also preserves
f{i}ltered quasi-isomorphisms, \cite{DeII}{7.I.6.2}.

\begin{lema}\label{SimplDeltaConsQuisF{i}tr} If $\mc{A}$ is an abelian category,
the functor $(\mbf{s},\delta):\Dl\mrm{CF}^+\mc{A}\rightarrow
\mrm{CF}^+\mc{A}$
preserves f{i}ltered quasi-isomorphisms.\\
That is to say, if $f:(A,\mrm{W})\rightarrow(B,\mrm{V})$ is a
morphism of cosimplicial f{i}ltered cochain complexes such that +
$f^m:(A^m,\mrm{W})\rightarrow(B^m,\mrm{V})$ is a f{i}ltered
quasi-isomorphism for each $m$, then $\mbf{s}f:(\mbf{s}A,\delta
\mrm{W})\rightarrow (\mbf{s}B,\delta\mrm{V})$ is also a f{i}ltered
quasi-isomorphism.
\end{lema}

\begin{proof}[\textbf{Proof of \ref{Hodge+E_2}}]\mbox{}\\
We will apply the transfer lemma \ref{FDfuerte}${}^{\mrm{op}}$ to
$\mrm{U}:\mc{H}dg\rightarrow Ch^+\mathbb{Q}$.\\
(SDC 1$)^{op}$ holds since $\mc{H}dg$ is additive. Let us see (SDC
3$)'^{op}$, that is, let us check that
${\mbf{s}}_{\mc{H}dg}=(\mbf{s},\delta,\mbf{s}):\Dl\mc{H}dg\rightarrow\mc{H}dg$
is indeed a functor.\\
Given ${K}\in\Dl\mc{H}dg$, then ${\mbf{s}}_{\mc{H}dg}{K}$ is a
Hodge complex by \cite{DeIII} 8.I.15 i).\\
Hence, ${\mbf{s}}_{\mc{H}dg}({K})$ satisf{i}es conditions a) and
b) of def{i}nition \ref{Def{i}ComplHodgeMixto}. Indeed, they are
consequences of the functoriality of $(\mbf{s},\mbf{s})$ and
$(\mbf{s},\delta)$, and it can be proven that
$\mrm{H}^k(\mbf{s}K_{\mathbb{Q}})$ is a f{i}nite dimensional
vector space using the standard argument of the proof of (SDC 6)
in proposition \ref{cdcabCDS} (or equivalently, using the spectral
sequence associated with $\mbf{s}K_{\mathbb{Q}}$).\\
On the other hand, by assumption
$\alpha=(\alpha_0,\alpha_1,(\widetilde{K},\widetilde{\mrm{W}}))$
is such that $\alpha_i$ is a degreewise f{i}ltered
quasi-isomorphism for $i=0,1$. Then, from
\ref{SimplDeltaConsQuisF{i}tr} we deduce that
$(\mbf{s},\delta)\alpha_i$ is so, for $i=0,1$. Therefore, the
zig-zag $\mbf{s}\alpha$ given by formula (\ref{Def{i}SimplAlpha})
satisf{i}es condition c) of the def{i}nition of mixed Hodge complex. Thus, it remains to see (MHC).\\
Given an integer $n$,
$({}_{\delta\mrm{W}}\mbf{Gr}_n(\mbf{s}K_\mathbb{Q}),({}_{\delta\mrm{W}}\mbf{Gr}_n(\mbf{s}K_\mathbb{C}),\mbf{s}\mrm{F}),\mbf{Gr}_n(\mbf{s}\alpha))$
satisf{i}es the hypothesis of def{i}nition \ref{ComplHodge} of
Hodge complex of weight $n$ by loc. cit., except condition c)
which is trivially satisf{i}ed because each
$\mbf{s}\alpha_i$ is a f{i}ltered quasi-isomorphism.\\[0.1cm]
Let us check now the functoriality of ${\mbf{s}}_{\mc{H}dg}$ with respect to the morphisms of $\Dl\mc{H}dg$.\\
A morphism
$f=(f_\mathbb{Q},f_\mathbb{C},\widetilde{f}):((K_\mathbb{Q},\mrm{W}),(K_\mathbb{C},\mrm{W},\mrm{F}),\alpha)\rightarrow
((K'_\mathbb{Q},\mrm{W}'),(K'_\mathbb{C},\mrm{W}',\mrm{F}'),\alpha')$
in $\Dl\mc{H}dg$ gives rise to the following commutative diagram
of $\Dl\mrm{CF}^+\mathbb{C}$
$$\xymatrix@H=4pt@C=30pt@R=35pt{(K_\mathbb{C},\mrm{W}) \ar[d]_{f_{\mathbb{C}}} & (\widetilde{K},\widetilde{\mrm{W}}) \ar[l]_-{\alpha_0} \ar[r]^-{\alpha_1}  \ar[d]_{\widetilde{f}}   & (K_\mathbb{Q},\mrm{W})\otimes\mathbb{C}  \ar[d]^{f_{\mathbb{Q}}\otimes\mathbb{C}} \\
                                (K'_\mathbb{C},\mrm{W}')                       & (\widetilde{K}',\widetilde{\mrm{W}}') \ar[l]_-{\alpha'_0} \ar[r]^-{\alpha'_1}                       & (K'_\mathbb{Q},\mrm{W}')\otimes\mathbb{C}\ .}$$
Therefore, applying
$(\mbf{s},\delta):\Dl\mrm{CF}^+\mathbb{C}\rightarrow
\mrm{CF}^+\mathbb{C}$ we get a commutative diagram in
$\mrm{CF}^+\mathbb{C}$, that gives rise to
$$\xymatrix@H=4pt@C=30pt@R=35pt{(\mbf{s}K_\mathbb{C},\delta\mrm{W}) \ar[d]_{\mbf{s}f_{\mathbb{C}}} & (\mbf{s}\widetilde{K},\delta\widetilde{\mrm{W}}) \ar[l]_-{\mbf{s}\alpha_0} \ar[r]^-{\mbf{s}\alpha_1}  \ar[d]_{\widetilde{f}}   & (\mbf{s}(K_\mathbb{Q}\otimes\mathbb{C}),\delta(\mrm{W}\otimes\mathbb{C}))  \ar[d]^{\mbf{s}(f_{\mathbb{Q}}\otimes\mathbb{C})}\ar[r]^-{\sim} & (\mbf{s}K_\mathbb{Q},\delta\mrm{W})\otimes\mathbb{C}  \ar[d]^{(\mbf{s}f_{\mathbb{Q}})\otimes\mathbb{C}} \\
                                (\mbf{s}K'_\mathbb{C},\delta\mrm{W}')                              & (\mbf{s}\widetilde{K}',\delta\widetilde{\mrm{W}}') \ar[l]_-{\mbf{s}\alpha'_0} \ar[r]^-{\mbf{s}\alpha'_1}                       & (\mbf{s}(K'_\mathbb{Q}\otimes\mathbb{C}),\delta(\mrm{W}'\otimes\mathbb{C}))                                                  \ar[r]^-{\sim}& (\mbf{s}K'_\mathbb{Q},\delta\mrm{W}')\otimes\mathbb{C}                                                  \ .}$$
Therefore,
${\mbf{s}}_{\mc{H}dg}{f}=(\mbf{s}f_\mathbb{Q},\mbf{s}f_\mathbb{C},\mbf{s}\widetilde{f})$
is a morphism in $\mc{H}dg$.\\[0.1cm]
Now we will prove (SDC 4$)'{}^{\mrm{op}}$ and (SDC
5)'${}^{\mrm{op}}$.
Denote by $\lambda^{\mathbb{Q}}$, $\lambda^{\mathbb{C}}$ the
natural transformations relative the descent categories
$\mrm{CF}^+\mathbb{Q}$ and $\mrm{CF}^+\mathbb{C}$ (with the
structure given in \ref{F{i}lt+crecient}). These transformations
coincide at the level of cochain complexes with the usual
transformation $\lambda$ of \ref{Cocadenas}.\\
If
$((K_\mathbb{Q},\mrm{W}),(K_\mathbb{C},\mrm{W},\mrm{F}),\alpha)$
is a mixed Hodge complex, from \ref{F{i}lt+crecient} and
\ref{f{i}ltrados-quis}, it follows that
$\lambda^\mathbb{Q}_{K_\mathbb{Q}}$,
$\lambda^\mathbb{C}_{K_\mathbb{C}}$ and
$\lambda_{\widetilde{K}}^{\mathbb{C}}$ preserve the f{i}ltrations.
Set $L=L\times\Dl$. We state that the following diagram commutes
in $\mrm{CF}^+\mathbb{C}$
$$\xymatrix@H=4pt@C=23pt@R=35pt{(K_\mathbb{C},\mrm{W}) \ar[d]_{\lambda^{\mathbb{C}}_{K_\mathbb{C}}}   & (\widetilde{K},\widetilde{\mrm{W}}) \ar[l]_-{\alpha_0} \ar[r]^-{\alpha_1}  \ar[d]_{\lambda^{\mathbb{C}}_{\widetilde{K}}}                           & (K_\mathbb{Q},\mrm{W})\otimes\mathbb{C}  \ar[d]^{\lambda^{\mathbb{C}}_{K_\mathbb{Q}\otimes\mathbb{C}}} \ar[rd]^{\lambda^{\mathbb{Q}}_{K_\mathbb{Q}}\otimes\mathbb{C}}      &  \\
                                (\mbf{s}(K_\mathbb{C}),\delta(\mrm{W}))             & (\mbf{s}(\widetilde{K}'),\delta(\widetilde{\mrm{W}})) \ar[l]_-{\mbf{s}(\alpha_0)} \ar[r]^-{\mbf{s}(\alpha_1)}  & (\mbf{s}(K'_\mathbb{Q}\otimes\mathbb{C}),\delta(\mrm{W}'\otimes\mathbb{C})) \ar[r]^-{\sim}   & (\mbf{s}(K'_\mathbb{Q}),\delta(\mrm{W}'))\otimes\mathbb{C}  \ .}$$
Indeed, the squares commutes by the functoriality of
$\lambda^{\mathbb{C}}$, as well as the right triangle since
$\lambda_L$ is just the inclusion of $L$ as direct summand of $\mbf{s}(L\times\Dl)$.\\
Consequently
$\lambda_{\mc{H}dg}=(\lambda^{\mathbb{Q}}_{K_{\mathbb{Q}}},\lambda^{\mathbb{C}}_{K_{\mathbb{C}}},\lambda^{\mathbb{C}}_{\widetilde{K}})$
is a morphism in $\mc{H}dg$. It can be argued similarly with
${\mu}^{\mc{H}dg}_K=(\mu^{\mathbb{Q}}_{K_{\mathbb{Q}}},\mu^{\mathbb{C}}_{K_{\mathbb{C}}},\mu^{\mathbb{C}}_{\widetilde{K}})$.\\[0.1cm]
(FD 1)${}^{\mrm{op}}$ is trivial since $\mrm{U}$ is additive, and
the diagram
$$\xymatrix@H=4pt@C=30pt@R=20pt{\Dl\mc{H}dg \ar[r]^-{\Dl \mrm{U}}\ar[d]_{(\mbf{s},\delta,\mbf{s})} & \Dl Ch^+\mathbb{Q} \ar[d]^{\mbf{s}}\\
                                \mc{H}dg \ar[r]^{\mrm{U}}         & Ch^+\mathbb{Q}}$$
commutes. F{i}nally, it is clear that the transformations
$\lambda$ of $\mc{H}dg$ and $Ch^+\mathbb{Q}$ are compatible, and
the same happens with $\mu$, so \ref{FDfuerte}${}^{\mrm{op}}$
holds.
\end{proof}

\begin{obs} The f{i}xed length of the zig-zag
$\alpha$ of a mixed Hodge complex has no relevance in the previous
proof. Consequently, the category
$$(\mc{H}dg_{\mc{G}},{\mbf{s}_{\mc{H}dg}},{\mrm{E}_{\mc{H}dg}},\mu_{\mc{H}dg},\lambda_{\mc{H}dg})$$
def{i}ned similarly, is an additive cosimplicial descent category,
and the forgetful functor $\mrm{U}:\mc{H}dg\rightarrow
Ch^+\mathbb{Q}$ is again a functor of additive cosimplicial descent categories.\\
\end{obs}

%


%
\appendix
\chapter{Eilenberg-Zilber-Cartier Theorem}
 \setcounter{section}{1}
 \setcounter{subsection}{1}
\setcounter{thm}{0}

We will need the following theorem, known as the
Eilenberg-Zilber-Cartier (\cite{DP},2.9).
\begin{thm}[Eilenberg-Zilber-Cartier]\label{E-Z-C}\index{Index}{Eilenberg-Zilber-Cartier thm.}\mbox{}\\
Consider an additive category $\mc{U}$, and the square
$$\xymatrix@M=4pt@H=4pt@C=30pt{\simp\simp\mc{U} \ar[r]^{\simp K}\ar[dd]_{\mrm{D}} & \simp Ch_+ \mc{U} \ar[d]_K \\
                                                                                        & Ch_+ Ch_+ \mc{U} \ar[d]_{Tot}  \\
                                             \simp\mc{U} \ar[r]^K                       & Ch_+\mc{U}.}$$
a) If $V\in\simp\simp\mc{U}$, then the morphisms $Id_{V_{0,0}}
:[Tot K \simp K (V)]_0 \rightarrow [K\mrm{D}(V)]_{0}$  and
$Id_{V_{0,0}} :[K\mrm{D}(V)]_{0} \rightarrow [Tot K \simp K
(V)]_0$ can be extended to universal morphisms
$$\left\{\begin{array}{l}\eta_{E-Z} (V): Tot K \simp K (V)\rightarrow K\mrm{D}(V)\\
                  \mu_{E-Z}(V): K\mrm{D}(V)\rightarrow Tot K \simp K (V)\end{array}\right. ,$$
that are homotopy inverse.\\
b) Each universal morphism $F: Tot K \simp K (V)\rightarrow
K\mrm{D}(V)$ with $F_0=Id$ is $($universally$)$ homotopic to
$\eta_{E-Z} (V)$. Analogously, if $G:K\mrm{D}(V) \rightarrow Tot K
\simp K (V))$ is universal and $G_0=Id$, then $G$ is
$($universally$)$ homotopic to $\mu_{E-Z}(V)$.
\end{thm}
\begin{obs}\mbox{}\\
\textsc{i)} Given $V\in\simp\simp\mc{U}$, a morphism $\oplus
V_{p,q}\rightarrow \oplus V_{p',q'}$ in $\mc{U}$ is universal if
each component $V_{p,q}\rightarrow V_{p',q'}$ is of the form
$\sum_{\alpha,\beta}n_{\alpha,\beta} V(\alpha,\beta)$, where
$n_{\alpha,\beta}\in\mathbb{Z}$, and $\alpha:[p']\rightarrow [p]$
and $\beta:[q']\rightarrow [q]$ are morphisms of $\Dl$.\\
Similarly, a morphism $F$ in $Ch_+(\mc{U})$ between $Tot K \simp
K (V)$ and $K\mrm{D}(V)$ is universal if each $F_n$ is so.\\[0.1cm]
\textsc{ii)} Since $\eta_{E-Z}$ and  $\mu_{E-Z}$ are universal, it
follows that they are functorial in $V$, so they def{i}ne natural
transformations between $K\mrm{D}$ and $Tot K \simp K$.
\end{obs}
\begin{proof}[\textbf{\textit{Proof}}]
Given $p,q,r,s \geq 0$, let $M(p,q;r,s)$ be the free abelian group
generated by the pairs $(\alpha,\beta)$, where
$\alpha:[r]\rightarrow
[p]$ and $\beta:[s]\rightarrow [q]$ are morphisms in $\Dl$.\\
Consider the category $\mc{M}$ whose objects are the symbols
$M_{p,q}$, $p,q\geq 0$. A morphism from $M_{p,q}$ to $M_{r,s}$ is
just an element of $M(p,q;r,s)$. Composition in $\mc{M}$ is
inherited from composition in $\Dl$. Hence
$(\alpha',\beta')(\alpha,\beta)=(\alpha\,\alpha',\beta\,\beta')$
if $(\alpha,\beta)\in M(p,q;r,s)$ and $(\alpha',\beta')\in
M(r,s;t,u)$. Consequently, $\simp\times\simp \subseteq \mc{M }$.

Attach objects and morphisms to $\mc{M}$ in such a way that we can
consider f{i}nite direct sums in $\mc{M}$. In this way we get the
additive category $\widetilde{\mc{M}}$, and again
$\simp\times\simp \subseteq \widetilde{\mc{M }}$. Therefore,
restricting the identity
$\widetilde{\mc{M}}\rightarrow\widetilde{\mc{M}}$ we obtain a simplicial object $M\in\simp\simp\widetilde{\mc{M }}$.\\
Consider $V\in\simp\simp\mc{U}$. Then $V$ can be extended in a
unique way to an additive functor $\mc{M}\rightarrow \mc{U}$, that
gives rise to $V': Ch_+(\widetilde{\mc{M}})\rightarrow
Ch_+(\mc{U})$.

Thus, a morphism in $Ch_+(\mc{U})$ between $Tot K \simp K (V)$ and
$K\mrm{D}(V)$ is universal if and only if it is the image under
$V'$ of a morphism of $Ch_+(\widetilde{\mc{M}})$ between $Tot K
\simp K (M)$ and $K\mrm{D}(M)$. Moreover, since
$V:\widetilde{\mc{M}}\rightarrow\mc{U}$ is additive, a homotopy
between two morphisms $F$ and $G$ from $Tot K \simp K (M)$ to
$K\mrm{D}(M)$ (or viceversa) is mapped by $V'$ into a (universal)
homotopy between $V'(F)$ and $V'(G)$, that are morphisms from $Tot
K \simp K (V)$ to $K\mrm{D}(V)$ (or viceversa).
Hence, we can restrict ourselves to $\mc{U}=\widetilde{\mc{M}}$
and $V=M$.

Let $Ab$ be the category of abelian groups. For each $l\geq 0$
denote by $K(l)\in\simp Ab$ the simplicial abelian group such that
$K(l)_{p}$ is the free abelian group generated by the morphisms
$\alpha:[p]\rightarrow [l]$. In other words, $K(l)$ is obtained
from the ``standard'' simplicial object $\triangle[l]$, by taking
free groups. Consider $K(l,m)\in\simp\simp Ab$
given by $K(l,m)_{p,q}=K(l)_p\otimes K(m)_q$.\\
Consider $R_{p,q}:\simp\simp Ab\rightarrow Ab$ with
$R_{p,q}(W)=W_{p,q}$ if $W\in\simp\simp Ab$, and the natural
transformations $\tau:R_{p,q}\rightarrow R_{p',q'}$ (also called
$FD$-operators). Denote by $N(p,q;p',q')$ the group consisting of
all of them. Examples of such $\tau$ are the basic transformations
$(\alpha,\beta)^*(W)=W_{\alpha,\beta}$
if $\alpha:[p']\rightarrow [p]$ and $\beta:[q']\rightarrow [q]$.\\
By \cite{EM} 3.1 we have that $N(p,q;p',q')$ is a free group
generated by the basic transformations. In addition, $\tau\in
N(p,q;p',q')$ is characterized by its value at the bisimplicial
abelian groups $K(l,m)$, $\forall l,m$.\\
Clearly, the mapping $(\alpha,\beta)\rightarrow (\alpha,\beta)^*$
is injective, since given any $\alpha:[p']\rightarrow [p]$ and
$\beta:[q']\rightarrow [q]$ then
$(\alpha,\beta)^*(K(p,q))(Id_{[p]},Id_{[q]})=(\alpha,\beta)\in
K(p,q)_{p',q'}$.
It follows that $N(p,q;p',q')\simeq M(p,q;p',q')$.

Therefore, $\mc{M}$ can be replaced by the category $\mc{N}$ whose
objects are symbols $N_{p,q}$ and whose morphisms between
$N_{p,q}$ and $N_{p',q'}$ are just $N(p,q;p',q')$. Similarly,
consider the restrictions of the identity functor
$\widetilde{\mc{N}}$ and $N\in\simp\simp\widetilde{\mc{N}}$. The
morphism $Id:N_{0,0}\rightarrow N_{0,0}$ in $N(0,0;0,0)$ is the
natural transformation def{i}ned by
$Id:K(l,m)_{0,0}\rightarrow K(l,m)_{0,0}$ for all $l,m$.\\
By the classical Eilenberg-Zilber theorem \cite{May} 29.3 it
follows the existence of natural transformations def{i}ned from
$\eta(l,m):\oplus_{p+q=n}K(l)_p\otimes K_{q}\rightarrow
K(l)_n\otimes K(m)_{n}$ and $\mu(l,m):K(l)_n\otimes
K(m)_{n}\rightarrow \oplus_{p+q=n}K(l)_p\otimes K_{q}$ for all
$l,m$, such that they are the identity in degree 0, and such that
they are homotopy inverse. In addition, any two such $\eta$ are
homotopic in a natural way, and the same holds for $\mu$.

By def{i}nition of $\widetilde{\mc{N}}$, the morphisms
$\{\eta(l,m)\}_{l,m}$ correspond to morphisms
$\eta\in\oplus_{p+q=n}N(p,q;n,n)$. In other words $\eta:Tot K
\simp K (N)\rightarrow K\mrm{D}(N)$, and similarly
$\{\mu(l,m)\}_{l,m}$ correspond to $\mu:K\mrm{D}(N)\rightarrow Tot
K \simp K (N)$, in such a way that (natural) homotopies are
preserved by this correspondence. Hence the proof is f{i}nished.
\end{proof}
\begin{obs}[Description of \textbf{$\eta_{E-Z}$} and \textbf{$\mu_{E-Z}$}]\label{Alex-Whitney-Map}\mbox{}\\
Given $V\in Ch_+ Ch_+ \mc{U}$, the
``shuf{f}le''\index{Index}{``shuf{f}le'' map} map $\eta_{E-Z} (V):
Tot K \simp K (V)\rightarrow
K\mrm{D}(V)$\index{Symbols}{$\eta_{E-Z}$} and the
Alexander-Whitney map\index{Index}{Alexander-Whitney map}
$\mu_{E-Z}(V): K\mrm{D}(V)\rightarrow Tot K \simp K
(V)$\index{Symbols}{$\mu_{E-Z}$} are (universal) inverse homotopy
equivalence (\cite{DP} 2.15).
They are given by $\eta_{E-Z}(V)=\Sigma \eta_{i,j}(V):
\ds\oplus_{i+j=k}V_{ij}\rightarrow V_{k,k}$
$$\eta_{i,j}(V)=\Sigma_{(\alpha,\beta)}sign(\alpha,\beta)V(\sigma^{\alpha_{j}}\sigma^{\alpha_{j-1}}\cdots\sigma^{\alpha_{1}},\sigma^{\beta_i}\sigma^{\beta_{i-1}}\cdots\sigma^{\beta_1})$$
where the sum is indexed over the $(i,j)$-``shuf{f}les''
$(\alpha,\beta)$ and $sign(\alpha,\beta)$ denotes the sign of
$(\alpha,\beta)$ (see \cite{EM}).\\
On the other hand
$\mu_{E-Z}(V)=\Sigma\mu_{i,j}(V):V_{k,k}\rightarrow
\ds\oplus_{i+j=k}V_{i,j}$, where
$$\mu_{i,j}(V)=V(d^0\stackrel{k-i)}{\cdots}d^0,d^kd^{k-1}\cdots d^{j+1})\ .$$
\end{obs}

\begin{obs}
The ``shuf{f}le'' $\widetilde{\eta}_{E-Z}$ and Alexander-Whitney
$\widetilde{\mu}_{E-Z}$ maps given in
\cite{DP} are not exactly those used in these notes.\\
The reason is that the total functor used in \cite{DP},
$\widetilde{Tot}:Ch_+ Ch_+ \mc{U} \rightarrow Ch_+\mc{U} $, is
isomorphic but not the same used here (see \ref{TotordenOpuesto}).
Indeed, given $\{V^{i,j}\ ;\ \ d^i ,\, d^j\}\in Ch_+ Ch_+ \mc{U}$,
then $\widetilde{Tot}(V)$ has as boundary map $\sum d^i+(-1)^id^j$.\\
Therefore, $\widetilde{\eta}_{E-Z}:\widetilde{Tot} K \simp K
\rightarrow K\mrm{D}$ and
$\widetilde{\mu}_{E-Z}:K\mrm{D}\rightarrow
\widetilde{Tot} K \simp K $.\\
Denote by $\Gamma:Ch_+ Ch_+ \mc{U}\rightarrow Ch_+ Ch_+ \mc{U}$
and $\Gamma:\simp\simp\mc{U}\rightarrow\simp\mc{U}$ the functors
that interchange the indexes in a double complex and in a
bisimplicial object respectively. Then
$$\eta_{E-Z}(V)=\widetilde{\eta}_{E-Z}(\Gamma V) \mbox{ and } \mu_{E-Z}(V)=\widetilde{\mu}_{E-Z}(\Gamma V)\ .$$
Note that $\widetilde{Tot}\Gamma=Tot$, $\mrm{D}\Gamma=\mrm{D}$ and
$K\simp K \Gamma= \Gamma K\simp K$. Hence,
$\widetilde{\mu}_{E-Z}(\Gamma V):K\mrm{D}(\Gamma
V)=K\mrm{D}(V)\rightarrow \widetilde{Tot} K \simp K (\Gamma V)=Tot
K\simp K$, and similarly for $\eta_{E-Z}$.
\end{obs}

We will use in these notes the Alexander-Whitney map in order to
proof the factorization axiom in the (co)chain complexes case. We
will need as well the following property of $\mu_{E-Z}$.

\begin{prop}\label{AsociatMuE-Z} The natural transformation $\mu_{E-Z}$ is associative. More concretely,
given a trisimplicial object $T$ in $\mc{U}$, the morphisms
$T_{n,n,n}\rightarrow \bigoplus_{r+s+t=n}T_{r,s,t}$ obtained by
applying twice $\mu_{E-Z}$
$$ \xymatrix@H=4pt@M=4pt@R=9pt{
                                                &  \bigoplus_{i+j=n}T_{i,j,n} \ar[rd]^{\mu_{EZ}} & \\
 T_{n,n,n}\ar[ru]^{\mu_{EZ}} \ar[rd]_{\mu_{EZ}} &                                                & \bigoplus_{r+s+t=n}T_{r,s,t} \\
                                                & \bigoplus_{p+q=n} T_{n,p,q} \ar[ru]_{\mu_{EZ}}  }$$
coincide.
\end{prop}

This is a well-known property (see, for instance, \cite{H}
14.2.3.). It holds that $\eta_{E-Z}$ is also associative, and in
addition it is symmetric (that is, it is invariant by swapping the
indexes). On the other hand, $\mu_{E-Z}$ is not symmetric, but it
is symmetric up to homotopy equivalence.

\begin{obs} Usually the transformations
$\mu_{E-Z}$ and $\eta_{E-Z}$ are used when $\mc{U}$ is a category
of $R$-modules, and the bisimplicial object considered is of the
form $\{X_n\otimes Y_m\}_{n,m}$ (for instance in the study of the
relationship between the homology of the cartesian product of
topological spaces and the tensorial product of their homologies).
\end{obs}

%
%
%


 \twocolumn{\chapter*{Symbols Index}}\addcontentsline{toc}{chapter}{Symbols Index}
  \input{Symbols.ind}

  \twocolumn{\chapter*{Index}}\addcontentsline{toc}{chapter}{Index}
  \input{Index.ind}

\end{document}